\documentclass[a4paper,11pt,leqno]{book}
\usepackage[T1]{fontenc}	
\usepackage[utf8]{inputenc}
\usepackage[french, english]{babel}
\usepackage{amssymb}
\usepackage{graphicx}
\usepackage{enumerate}
\usepackage{amsthm}
\usepackage[toc,page]{appendix}
\usepackage{amsmath}

\usepackage{imakeidx}
\makeindex[intoc]

\theoremstyle{plain}

\newtheorem{thm}{Theorem}

\newtheorem{cor}[thm]{Corollary}
\newtheorem{lem}[thm]{Lemma}
\newtheorem{prop}[thm]{Proposition}
\theoremstyle{definition}
\newtheorem{defn}[thm]{Definition}

\newtheorem{remark}[thm]{Remark}
\theoremstyle{remark}
\newtheorem{notation}{Notation}
{\left\lbrace\begin{array}{@{}l@{}}}%
{\end{array}\right.}
\usepackage{geometry} 
\usepackage{fullpage}
\usepackage{eucal}
\usepackage[babel]{csquotes}

\usepackage{pdfpages}
\usepackage{textcomp}
\usepackage{macros}
\usepackage{pdfpages}
\usepackage{hyperref}
\hypersetup{
    colorlinks=false,
    linkcolor=blue,
    filecolor=magenta,      
    urlcolor=cyan,
}
\numberwithin{equation}{section}
\numberwithin{thm}{section}
\renewcommand{\thesection}{\thechapter.\arabic{section}}

\usepackage{scalerel,stackengine}
\stackMath
\newcommand\reallywidehat[1]{%
\savestack{\tmpbox}{\stretchto{%
  \scaleto{%
    \scalerel*[\widthof{\ensuremath{#1}}]{\kern-.6pt\bigwedge\kern-.6pt}%
    {\rule[-\textheight/2]{1ex}{\textheight}}
  }{\textheight}%
}{0.5ex}}%
\stackon[1pt]{#1}{\tmpbox}%
}
 \usepackage{dsfont}

\newcommand\numberthis{\addtocounter{equation}{1}\tag{\theequation}}

\usepackage{color}

\usepackage[hang,flushmargin]{footmisc} 

\usepackage{setspace}

\usepackage[textsize=small]{todonotes}

\usepackage{hyperref}
\hypersetup{
    colorlinks=false,
    linkcolor=blue,
    filecolor=magenta,      
    urlcolor=cyan,
}

\usepackage{titling}
\usepackage{blindtext}

\newsavebox{\abstractbox}
\newenvironment{abstract}
 {%
  \global\setbox\abstractbox=\vtop\bgroup
  \begin{center}\bfseries\abstractname\end{center}%
 }
 {\par\egroup}

\usepackage[backend=biber,style=trad-plain, giveninits=true]{biblatex}

\usepackage{tikz}

\title{ Global existence of small amplitude solutions for a model quadratic quasi-linear coupled wave-Klein-Gordon system in two space dimension, with mildly decaying Cauchy data}

\author{A. Stingo
 \thanks{The author is supported by a PhD fellowship funded by the FSMP and the Labex MME-DII, and by For Women in Science fellowship funded by Fondation L'Or\'eal-UNESCO.\newline
   }
\\
Universit{\'e} Paris 13,\\
Sorbonne Paris Cit\'e, LAGA, CNRS (UMR 7539),\\
99, Avenue J.-B. Cl{\'e}ment,\\
F-93430 Villetaneuse}

\date{2018}

\begin{document}

\begin{abstract}
The aim of this paper is to study the global existence of solutions to a coupled wave-Klein-Gordon system in space dimension two when initial data are small, smooth and mildly decaying at infinity.
Some physical models strictly related to general relativity have shown the importance of studying such systems, but very few results are know at present in low space dimension. We study here a model two-dimensional system, in which the non-linearity writes in terms of ``null forms'', and show the global existence of small solutions.
Our goal is to prove some energy estimates on the solution when a certain number of Klainerman vector fields is acting on it, and some optimal uniform estimates.
The former ones are obtained using systematically quasi-linear normal forms, in their para-differential version; the latter ones are recovered by deducing a new coupled system of a transport equation and an ordinary differential equation from the starting PDE system, by means of a semi-classical micro-local analysis of the problem. 
We expect the strategy developed here to be robust enough to enable us, in the future, to treat the case of the most general non-linearities.
\end{abstract}

\maketitle 

\tableofcontents

\chapter*{Introduction}

\addcontentsline{toc}{chapter}{Introduction}

\renewcommand\thethm{\arabic{thm}}
\def\theequation{\arabic{equation}}
\renewcommand\thesection{\arabic{section}}
\setcounter{equation}{0}
\setcounter{thm}{0}

The result we present in this paper concerns the global existence of solutions to a quadratic quasi-linear coupled system of a wave equation and a Klein-Gordon equation in space dimension two, when initial data are small smooth and mildly decaying at infinity.
We prove this result for a model non-linearity with the aim of extend it, in the future, to the most general case.
Keeping this long term objective in mind, we shall try to develop a fairly general approach in spite of the fact that we are treating here a simple model.
The Cauchy problem we consider is the following
\begin{equation} \label{eq:WKG_intro}
\begin{cases}
& (\partial^2_t  - \Delta_x) u(t,x) = Q_0(v, \partial_1 v)\,, \\
& (\partial^2_t  - \Delta_x  + 1)v(t,x) = Q_0(v, \partial_1 u)\,,
\end{cases} \qquad (t,x)\in ]1,+\infty[\times \mathbb{R}^2
\end{equation}
with initial conditions
\begin{equation}\label{eq:data_intro}
\begin{cases}
& (u,v)(1,x) = \varepsilon(u_0(x),v_0(x))\,,\\
& (\partial_t u,\partial_t v)(1,x) = \varepsilon (u_1(x), v_1(x))\,,
\end{cases}
\end{equation}
where $\varepsilon>0$ is a small parameter, and $Q_0$ is the null form: 
\begin{equation*}
Q_0(v, w) = (\partial_t v)(\partial_t w) - (\nabla_x v)\cdot (\nabla_x w)\,.
\end{equation*}
We also suppose that, for some $n\in\mathbb{N}$ sufficiently large, $(\nabla_x u_0, u_1)$ is in the unit ball of $H^n(\mathbb{R}^2,\mathbb{R})\times H^n(\mathbb{R}^2,\mathbb{R})$, $(v_0,v_1)$ in the unit ball of $H^{n+1}(\mathbb{R}^2,\mathbb{R})\times H^n(\mathbb{R}^2,\mathbb{R})$, and that
\begin{equation}\label{eq:condition_data_intro}
 \sum_{1\leq\abs{\alpha}\leq 3}\bigl(\norm{x^\alpha\nabla_xu_0}_{H^{\abs{\alpha}}} + \norm{x^\alpha v_0}_{H^{\abs{\alpha}+1}}
  + \norm{x^\alpha u_1}_{H^{\abs{\alpha}}}  + \norm{x^\alpha v_1}_{H^{\abs{\alpha}}}\bigr)\leq 1.
\end{equation}

\smallskip

Some physical models, especially related to general relativity, have shown the importance of studying such systems to which several recent works have been dedicated.
Most of the results known at present concern wave-Klein-Gordon systems in space dimension 3.
One of the first ones goes back to Georgiev \cite{Georgiev:system}. He observed that the vector fields' method developed by Klainerman was not well adapted to handle at the same time massless and massive wave equations because of the fact that the scaling vector field $S = t\partial_t + x\cdot\nabla_x$ is not a Killing vector field for the Klein-Gordon equation. To overcome this difficulty he adapted Klainerman's techniques, introducing a \textit{strong null condition} to be satisfied by semi-linear nonlinearities that ensures global existence.
In 2012 Katayama \cite{Katayama:coupled_systems} showed the global existence of small amplitude solutions to coupled systems of wave and Klein-Gordon equations under certain suitable conditions on the non-linearity that include the \textit{null condition} of Klainerman (\cite{klainerman:null_condition}) on self-interactions between wave components, and are weaker than the \textit{strong null condition} of Georgiev. Consequently, the result he obtains applies also to certain other physical systems such as Dirac-Klein-Gordon equations, Dirac-Proca equations and Klein-Gordon-Zakharov equations. 
Later, this problem was also studied by LeFloch, Ma \cite{LeFloch_Ma} and Wang \cite{Q.Wang} as a model for the full Einstein-Klein-Gordon system (E-KG)
\begin{equation*}
\begin{cases}
& Ric_{\alpha\beta}=\textbf{D}_\alpha \psi \textbf{D}_\beta \psi + \frac{1}{2}\psi^2 g_{\alpha\beta} \\
& \Box_g \psi = \psi
\end{cases}
\end{equation*}
The authors prove global existence of solutions to wave-Klein-Gordon systems with quasi-linear quadratic non-linearities satisfying suitable conditions, when initial data are small, smooth and compactly supported, using the so-called \textit{hyperboloidal foliation method} introduced by Le Foch, Ma in \cite{LeFloch_Ma}.
Global stability for the full (E-KG) has been then proved by LeFloch-Ma \cite{LeFloch-Ma:global_nl_stability, lefloch-ma:global_stability_2} in the case of small smooth perturbations that agree with a Scharzschild solution outside a compact set (see also Wang \cite{Q.Wang:E-KG}).
In a recent paper \cite{Ionescu_Pausader:WKG} Ionescu and Pausader prove global regularity and modified scattering in the case of small smooth initial data that decay at suitable rates at infinity, but not necessarily compactly supported.
The quadratic quasi-linear problem they deal with is the following
\begin{equation*}
\begin{cases}
& -\Box u = A^{\alpha\beta}\partial_\alpha v\partial_\beta v + D v^2 \\
& -(\Box + 1) v = u B^{\alpha\beta}\partial_\alpha \partial_b v
\end{cases}
\end{equation*}
where $A^{\alpha\beta}, B^{\alpha\beta}, D$ are real constants. The system keeps the same linear structure as (E-KG) in harmonic gauge, but only keeps quadratic non-linearities that involve the massive scalar field $v$ (semilinear in the wave equation, quasi-linear in the Klein-Gordon equation). Moreover, the non-linearity they consider does not present a null structure but shows a particular resonant pattern. Their result relies, on the one hand, on a combination of energy estimates to control high Sobolev norms and weighted norms involving the admissible vector fields; on the other hand, on a Fourier analysis, in connection with normal forms and analysis of resonant sets, to prove dispersive estimates and decay in suitable lower regularity norms. 
The only results we know about global existence of small amplitude solutions in lower space dimension are due to Ma.
In space dimension 2 he considers the case of compactly supported Cauchy data and adapts the hyperboloidal foliation method mentioned above to $2+1$ spacetime wave-Klein-Gordon systems (see \cite{ma:2D_tools}). In particular, in \cite{ma:2D_quasilinear} he combines this method with a normal form argument to treat some quasi-linear quadratic non-linearities, while in \cite{ma:2D_semilinear} he studies the case of some semi-linear quadratic interactions. In a very recent paper \cite{ma:1D_semilinear} he also tackles the one-dimensional problem, studying a model semi-linear cubic wave-Klein-Gordon system. In this work he finally overcomes the restriction on the support of initial data and generalizes the hyperboloidal foliation method, combining the hyperboloidal foliation of the translated light cone with the standard time-constant foliation outside of it. The analysis of the problem and the deduction of the estimates of interest is then made separately inside and outside the mentioned light cone.

The result we prove in this paper is the following:

\begin{thm} \label{thm:main_intro}
There exists $\varepsilon_0>0$ such that for any $\varepsilon\in ]0,\varepsilon_0[$, system \eqref{eq:WKG_intro} with initial data satisfying \eqref{eq:data_intro}, \eqref{eq:condition_data_intro} admits a unique global solution defined on $[1,+\infty[$, with $\partial_{t,x}u\in C^0([1,+\infty[; H^n(\mathbb{R}^2))$ and $(v, \partial_tv)\in C^0([1,+\infty[; H^{n+1}(\mathbb{R}^2)\times H^n(\mathbb{R}^2))$.
\end{thm}

We describe below the strategy of the theorem's proof.
First of all, we rewrite system \eqref{eq:WKG_intro} in terms of unkowns
\begin{equation}\label{eq:u_pm, v_pm_intro}
u_\pm = \left(D_t \pm |D_x|\right)u, \quad v_\pm = \left(D_t \pm \langle D_x\rangle\right)v,
\end{equation}
where $D_{t,x}=-i\partial_{t,x}$, and introduce the admissible Klainerman vector fields for this problem, i.e.
\begin{equation*}
\Omega = x_1\partial_2 -x_2\partial_1, \quad Z_j = x_j\partial_t + t\partial_j, \quad j=1,2.
\end{equation*}
We also denote by $\mathcal{Z}=\{\Gamma_1,\dots, \Gamma_5\}$ the family made by above vector fields together with the two spatial derivatives, and if $I=(i_1,\dots,i_p)$ is an element of $\{1,\dots,5\}^p$, $\Gamma^I w$ is the function obtained letting $\Gamma_{i_1},\dots,\Gamma_{i_p}$ act successively on $w$.
We then set
\begin{equation}\label{eq:uI_pm, vI_pm_intro}
u^I_\pm = \left(D_t \pm |D_x|\right)\Gamma^I u, \quad v^I_\pm = \left(D_t \pm \langle D_x\rangle\right)\Gamma^I v,
\end{equation}
and introduce the following energies:
\[E_0(t;u_\pm,v_\pm) = \int_{\R^2}\bigl(\abs{u_+(t,x)}^2 + \abs{u_-(t,x)}^2 + \abs{v_+(t,x)}^2 +
\abs{v_-(t,x)}^2\bigr)\,dx,\]
then for $n\ge 3$,
\begin{equation*}
  E_n(t;u_\pm,v_\pm) = \sum_{\abs{\alpha}\leq n}E_0(t;D_x^{\alpha}u_\pm, D_x^{\alpha}v_\pm), 
\end{equation*}
which controls the $H^n$ regularity of $u_\pm, v_\pm$, and finally, for any integer $k$ between 0 and 2,
\begin{equation*}
  E_3^k(t;u_\pm,v_\pm) = \sum_{\substack{\abs{\alpha}+\abs{I}\leq 3\\ \abs{I}\le 3-k}} E_0(t;D_x^\alpha u_\pm^I, D_x^\alpha v_\pm^I)
\end{equation*}
that takes into account the decay in space of $u_\pm, v_\pm$ and of at most three of their spatial derivatives.
By a local existence argument, an a-priori estimate on $E_n$ on a certain time interval will be enough to ensure the extension of the solution to that interval. For this reason, we are led to prove a result as the following one, in which $\mathrm{R}=(\mathrm{R}_1,\mathrm{R}_2)$ denotes the Riesz transform:

\begin{thm} \label{thm:boostrap_intro}
Let $K_1,K_2$ two constants strictly bigger than 1. There exist two integers $n\gg \rho\gg 1$, $\varepsilon_0\in ]0,1[$ small enough, some small real $0<\delta\ll \delta_2\ll \delta_1\ll \delta_0\ll 1$ and two constants $A,B$ sufficiently large such that, if functions $u_\pm, v_\pm$, defined by \eqref{eq:u_pm, v_pm_intro} from a solution to \eqref{eq:WKG_intro}, satisfy
\begin{equation} \label{eq:boostrap_1_intro}
    \begin{split}
      &\norm{\absj{D_x}^{\rho+1} u_\pm(t,\cdot)}_{L^\infty} + \norm{\absj{D_x}^{\rho+1} \mathrm{R}u_\pm(t,\cdot)}_{L^\infty} \leq A\varepsilon t^{-\frac{1}{2}}\\
&\norm{\absj{D_x}^{\rho} v_\pm}_{L^\infty} \leq A\varepsilon t^{-1}\\
&E_n(t;u_\pm,v_\pm) \leq B^2\varepsilon^2 t^{2\delta}\\
&E_3^k(t; u_\pm,v_\pm) \leq B^2\varepsilon^2 t^{2\delta_{k}},\ 0\leq k\leq 2,
    \end{split}
  \end{equation}
  for every $t\in [1,T]$, then on the same interval $[1,T]$ we have
  \begin{equation}\label{eq:bootstrap_2_intro}
    \begin{split}
      &\norm{\absj{D_x}^{\rho+1} u_\pm(t,\cdot)}_{L^\infty} + \norm{\absj{D_x}^{\rho+1} \mathrm{R}u_\pm(t,\cdot)}_{L^\infty} \leq \frac{A}{K_1}\varepsilon t^{-\frac{1}{2}}\\
&\norm{\absj{D_x}^{\rho} v_\pm}_{L^\infty} \leq \frac{A}{K_1}\varepsilon t^{-1}\\
&E_n(t;u_\pm,v_\pm) \leq \frac{B^2}{K_2^2}\varepsilon^2 t^{2\delta}\\
&E_3^k(t; u_\pm,v_\pm) \leq \frac{B^2}{K_2^2}\varepsilon^2 t^{2\delta_{k}},\ 0\leq k\leq 2.
    \end{split}
  \end{equation}
\end{thm} 
The proof of the theorem consists, on the one hand, to prove that \eqref{eq:boostrap_1_intro} implies the latter two estimates in \eqref{eq:bootstrap_2_intro} by means of an energy inequality. On the other hand, by reduction of the starting problem to a coupled system of an ordinary differential equation and a transport equation, we prove that \eqref{eq:boostrap_1_intro} implies the first two estimates in \eqref{eq:bootstrap_2_intro}.

\smallskip
In order to recover the mentioned energy inequality that allows us to propagate the a-priori energy estimates, we let family $\Gamma^I$ of vector fields act on \eqref{eq:WKG_intro} and then pass to unknowns \eqref{eq:uI_pm, vI_pm_intro}.
We obtain a new system of the form
\begin{equation*}
  \begin{split}
    (D_t\mp \abs{D_x})u^I_\pm &= \textit{NL}_\mathrm{w}(v_\pm^I,v_\pm^I)\\
(D_t\mp \abs{D_x})v^I_\pm &= \textit{NL}_\mathrm{kg}(v_\pm^I,u_\pm^I)
  \end{split}
\end{equation*}
where the non-linearities (whose explicit expression may be found in the right hand side of \eqref{system for uI+-, vI+-}) are bilinear quantities of their arguments.
Because of the quasi-linear nature of our problem, the first step towards the derivation of the mentioned inequality is to highlight the very quasi-linear contribution to above non-linearities and make sure that it does not lead to a loss of derivatives. For this reason, we write the above system in a vectorial fashion by introducing vectors
\[U^I =
  \begin{bmatrix}
    u_+^I\\0\\u_-^I\\0
  \end{bmatrix}, \ 
V^I = \begin{bmatrix}
    0\\v_+^I\\0\\v_-^I
  \end{bmatrix},\ W^I = U^I + V^I,\]
and successively \emph{para-linearize} the vectorial equation satisfied by $W^I$ (using the tools introduced in subsection \ref{Subsection: Paradifferential Calculus}) to stress out the quasi-linear contribution to the non-linearity. Finally, we \emph{symmetrize} it (in the sense of subsection \ref{Subs: Symmetrization}) by introducing some new unknown $W^I_s$ comparable to $W^I$.
What we would need to show in order to prove the last two inequalities in \eqref{eq:bootstrap_2_intro} is that, using the estimates in \eqref{eq:boostrap_1_intro}, the derivative in time of the $L^2$ norm to the square of $W^I_s$ is bounded by $\frac{C\varepsilon}{t}\|W^I\|_{L^2}$. By analysing the semi-linear contributions in the symmetrized equation satisfied by $W^I_s$, we find out that the $L^2$ norm of some of those ones can only be estimated making appear the $L^\infty$ norm on the wave factor and the $L^2$ norm on the Klein-Gordon one.
Because of the very slow decay in time of the wave solution (the decay rate being $t^{-1/2}$, as assumed in the first inequality of \eqref{eq:boostrap_1_intro}), we are hence very far away from the wished estimate. Consequently, the second step for the derivation of the right energy inequality consists in performing a normal form argument to get rid of those quadratic terms and replace them with cubic ones. For that, we first use a Shatah' normal form adapted to quasi-linear equations (see subsection \ref{sub: a first normal form transformation}) as already used by several authors (we cite \cite{ozawaTT:remarks, D1, D2, D3} for quasi-linear Klein-Gordon equations, and \cite{HITW, HIT, IP1, AD2, IT1} for quasi-linear equations arising in fluids mechanics), but also a semi-linear normal form argument to treat some other terms on which we are allowed to lose some derivatives (see subsection \ref{sub: second normal form}).
These two normal forms' steps lead us to define some new energies $\widetilde{E}_n(t;u_\pm, v_\pm), \widetilde{E}^k_3(t;u_\pm, v_\pm)$, equivalent to the starting ones $E_n(t;u_\pm, v_\pm)$, $\widetilde{E}^k_3(t;u_\pm, v_\pm)$, that we are able to propagate. That concludes the first part of the proof.

\smallskip
The last thing that remains to prove is that \eqref{eq:boostrap_1_intro} implies the first two estimates in \eqref{eq:bootstrap_2_intro}.
The strategy we employ is very similar to the one developed in \cite{stingo:1D_KG}: we deduce from the starting system \eqref{eq:WKG_intro} a new coupled one of an ordinary differential equation, coming from the Klein-Gordon equation, and of a transport equation, derived from the wave one. The study of this system will provide us with the wished $L^\infty$ estimates.
We start our analysis by another normal form in order to replace almost all quadratic non-linear terms in the equations satisfied by $u_\pm, v_\pm$ with cubic ones. The only contributions that cannot be eliminated are those depending on $(v_+, v_{-})$ which are resonant and should be suitably treated.
We do not use directly the normal forms obtained in the previous step. In fact, our aim is basically to obtain an $L^\infty$ estimate for at most $\rho$ derivatives of the solution, having a control on their $H^s$ norm for $s\gg\rho$. This permits us to lose some derivatives in the normal form reduction, so the fact that the system is quasi-linear is no longer important.

We define two new unknowns $u^{NF}, v^{NF}$ by adding some quadratic perturbations to $u_{-}, v_{-}$, in such a way that they are solution to
\begin{equation}
  (D_t + \abs{D_x})\unf = q_w +c_w + \rnfw,\ (D_t + \abs{D_x})\vnf = \rnfkg,
\end{equation}
where $\rnfw, c_w, \rnfkg$ are cubic terms, whereas $q_w$ is the mentioned bilinear expression in $v_+, v_{-}$ that cannot be eliminated by normal forms but whose structure will successively provide us with remainder terms.
Then, if we define 
\begin{equation}
  \ut(t,x) = t\unf(t,tx),\ \vt(t,x) = t\vnf(t,tx),
\end{equation}
and introduce $h:=t^{-1}$ the \emph{semi-classical parameter}, we obtain that $\ut, \vt$ verify
\begin{equation} \label{eq:ut_vt_intro}
  \begin{split}
    (D_t - \oph(x\cdot\xi - \abs{\xi}))\ut &= h^{-1}\left[ q_w(t,tx) + c_w(t,tx)+ \rnfw(t,tx)\right]\\
(D_t - \oph(x\cdot\xi - \absj{\xi}))\vt &= h^{-1}\rnfkg(t,tx)
  \end{split}
\end{equation}
where $\oph$ is the Weyl quantization introduced, along with the semi-classical pseudo-differential calculus, in subsection \ref{Subsection: Paradifferential Calculus}.
We also consider the following operators
\begin{equation*} 
  \begin{split}
 \Mcal_j = \frac{1}{h}\Bigl(x_j\abs{\xi}  -\xi_j\Bigr),\  \Lcal_j = \frac{1}{h}\Bigl(x_j  -\frac{\xi_j}{\absj{\xi}}\Bigr),
\end{split}
\end{equation*}
whose symbols are given respectively (up to the multiplication by $|\xi|$ for the former case) by the derivative with respect to $\xi$ of symbols $x\cdot\xi - |\xi|$ and $x\cdot\xi - \langle\xi\rangle$ in \eqref{eq:ut_vt_intro}.
Using the equation satisfied by $\unf$ (resp. $\vnf$), we can express $\mathcal{M}_j\ut$ (resp. $\Lcal_j\vt$) in terms of $Z_j\unf$ (resp. $Z_j \vnf$) and of $q_w, c_w, \rnfw$ (resp. $\rnfkg$). 
As done in \cite{stingo:1D_KG}, we first introduce the lagrangian
\begin{equation*}
  \Lkg = \Bigl\{(x,\xi): x - \frac{\xi}{\absj{\xi}} = 0\Bigr\}
\end{equation*}
which is the graph of $\xi =-d\phi(x)$, with $\phi(x)=\sqrt{1-|x|^2}$,
and decompose $\vt$ into the sum of a contribution micro-localised on a neighbourhood of size $\sqrt{h}$ of $\Lkg$, and another one micro-localised out of that neighbourhood (in the spirit of \cite{ifrim_tataru:global_bounds}).
The second contribution can be basically estimated in $L^\infty$ by $h^{\frac{1}{2}-0}$ times the $L^2$ norm of some iterates of operator $\Lcal$ acting on $\vt$ (which are controlled by the $L^2$ hypothesis in theorem \ref{thm:boostrap_intro}).
The main contribution to $\vt$ is then represented by $\vt_\Lkg$, which appears to be solution to
\[[D_t-\oph(x\cdot\xi - \absj{\xi})]\vt_\Lkg =\textrm{ controlled terms}.\]
Developing the symbol in the above left hand side on $\Lkg$ we finally obtain the wished ODE, which combined with the a-priori estimate of the ``controlled terms'' allows us to deduce from \eqref{eq:boostrap_1_intro} the second estimate in \eqref{eq:bootstrap_2_intro} (with $\rho=0$, the general case being treated in the same way up to few more technicalities).

The same strategy is employed to obtain some uniform estimates on $\ut$.
We introduce the lagrangian
\begin{equation*}
  \Lw = \Bigl\{(x,\xi): x - \frac{\xi}{\abs{\xi}} = 0\Bigr\}
\end{equation*}
which, differently from $\Lkg$, is not a graph but projects on the basis as an hypersurface. For this reason, the classical problem associated to the first equation in \eqref{eq:ut_vt_intro} is rather a transport equation than an ordinary differential equation.
It is obtained in a similar way by decomposing $\ut$ into two contributions: one denoted by $\ut_\Lw$ and micro-localised in a neighbourhood of size $h^{\frac{1}{2}-\sigma}$ (for some small $\sigma>0$) of $\Lw$; another one micro-localised away from this neighbourhood. As for the Klein-Gordon component, this latter contribution can be easily controlled thanks to the $L^2$ estimates that the last two inequalities in \eqref{eq:boostrap_1_intro} infer on the iterates of $\Mcal_j$ acting on $\ut$.
By micro-localisation we derive that $\ut_\Lw$ satisfies
\[[D_t-\oph(x\cdot\xi-\abs{\xi})]\ut_\Lw = \textrm{ controlled terms},\]
and by developing symbol $x\cdot\xi-\abs{\xi}$ on $\Lw$ we obtain the wished transport equation. Integrating this equation by the method of characteristics, we finally recover the first estimate in \eqref{eq:boostrap_1_intro} and conclude the proof of theorem \ref{thm:boostrap_intro}.

\numberwithin{equation}{section}
\numberwithin{thm}{section}
\renewcommand{\thesection}{\thechapter.\arabic{section}}

\chapter{Main Theorem and Preliminary Results}

\section{Statement of the main theorem} \label{sec: statement of the main results}
\textsc{Notations:}
We warn the reader that, throughout the paper, we will often denote $\partial_t$ (resp. $\partial_{x_j}$, $j=1,2$) by $\partial_0$ (resp. $\partial_j$, $j=1,2$), while symbol $\partial$ without any subscript will stand for one of the three derivatives $\partial_a$, $a=0,1,2$. $\nabla_x f$ is the classical spatial gradient of $f$, $D:=-i\partial$ and $\mathrm{R}_j$ denotes the Riesz operator $D_j|D_x|^{-1}$, for $j=1,2$. We will also employ notation $\|\partial_{t,x}w\|$ with the meaning $\|\partial_tw\| + \|\partial_xw\|$ and $\|\mathrm{R}w\|=\sum_j\|\mathrm{R}_jw\|$.
\bigskip

We consider the following quadratic, quasi-linear, coupled wave-Klein-Gordon system 
\begin{equation} \label{wave KG system}
\begin{cases}
& (\partial^2_t  - \Delta_x) u(t,x) = Q_0(v, \partial_1 v)\,, \\
& (\partial^2_t  - \Delta_x  + 1)v(t,x) = Q_0(v, \partial_1 u)\,,
\end{cases} \qquad (t,x)\in ]1,+\infty[\times \mathbb{R}^2
\end{equation}
with initial conditions
\begin{equation} \label{initial data}
\begin{cases}
& (u,v)(1,x) = \varepsilon(u_0(x),v_0(x))\,,\\
& (\partial_t u,\partial_t v)(1,x) = \varepsilon (u_1(x), v_1(x))\,,
\end{cases}
\end{equation}
where $\varepsilon>0$ is a small parameter, and $Q_0$ is the null form:  \index{Q0@$Q_0(v, w)$, null form}
\begin{equation} \label{null form Q0}
Q_0(v, w) = (\partial_t v)(\partial_t w) - (\nabla_x v)\cdot (\nabla_x w)\,.
\end{equation}

Our aim is to prove that there is a unique solution to Cauchy problem \eqref{wave KG system}-\eqref{initial data} provided that $\varepsilon$ is sufficiently small and $u_0,v_0,u_1,v_1$ decay rapidly enough at infinity.
The theorem we are going to demonstrate is the following:

\begin{thm}[Main Theorem]\label{Thm: Main theorem}
There exist an integer $n$ sufficiently large and $\varepsilon_0\in]0,1[$ sufficiently small such that, for any $\varepsilon \in ]0,\varepsilon_0[$, any real valued $u_0,v_0,u_1,v_1$ satisfying:
\begin{equation}\label{condition_initial_data}
\begin{gathered}
\|\nabla_x u_0\|_{H^n}+ \|v_0\|_{H^{n+1}}+\|u_1\|_{H^n}+\|v_1\|_{H^n}\le 1,\\
\sum_{|\alpha|=1}^2 \left(\|x^\alpha \nabla_x u_0\|_{H^{|\alpha|}}+ \|x^\alpha v_0\|_{H^{|\alpha|+1}}+ \|x^\alpha u_1\|_{H^{|\alpha|}}+ \|x^\alpha v_1\|_{H^{|\alpha|}}\right)\le 1,
\end{gathered}
\end{equation}
system \eqref{wave KG system}-\eqref{initial data} admits a unique global solution $(u,v)$ with $\partial_{t,x}u\in C^0\left([1,\infty[;H^n(\mathbb{R}^2)\right)$, $v \in C^0\left([1,\infty[;H^{n+1}(\mathbb{R}^2)\right)\cap C^1\left([1,\infty[;H^n(\mathbb{R}^2)\right)$.
\end{thm}

The proof of the main theorem is based on the introduction of four new functions $u_+,u_{-},v_+,v_{-}$, defined in terms of $u,v$ as follows: \index{upm@$u_\pm$, wave components}\index{vpm@$v_\pm$ Klein-Gordon components}
\begin{equation} \label{def u+- v+-}
\begin{cases}
& u_+ : = (D_t + |D_x|)u\,, \\
& u_{-} : = (D_t - |D_x|)u \,,
\end{cases}\qquad
\begin{cases}
& v_+ : = (D_t +\langle D_x\rangle)v\,, \\
& v_{-} : = (D_t -\langle D_x\rangle)v\,,
\end{cases}
\end{equation}
and on the propagation of some a-priori estimates made on them in some interval $[1,T]$, for a fixed $T>1$.
In order to state this result we consider the admissible Klainerman vector fields for the wave-Klein-Gordon system: \index{Klainerman vector fields} \index{Omega@$\Omega$, Euclidean rotation} \index{Zj@$Z_j$, Lorentzian boost}
\begin{equation} \label{Omega, Zj}
\Omega:=x_1\partial_2 - x_2 \partial_1\,, \quad Z_j := x_j\partial_t + t\partial_j\,, j=1,2 
\end{equation}
and denote by $\Gamma$ a generic vector field in $\mathcal{Z}=\{\Omega, Z_j, \partial_j, j=1,2\}$. \index{Z@$\mathcal{Z}$, family of admissible vector fields}
If $\mathcal{Z}$ is assumed ordered, i.e.
\begin{equation} \label{order_Z}
\begin{gathered}
\mathcal{Z}=\{\Gamma_1,\dots, \Gamma_5\} \\
\text{with }\quad \Gamma_1=\Omega, \quad \Gamma_j= Z_{j-1} \ \text{for } j=2,3, \quad \Gamma_j=\partial_{j-3} \ \text{for } j=4,5,
\end{gathered}
\end{equation}
then for a multi-index $I=(i_1,\dots, i_n)$, $i_j\in\{1,\dots,5\}$ for $j=1,\dots,n$, we define the \textit{length} of $I$ as $|I|:=n$, and $\Gamma^I:=\Gamma_{i_1}\cdots\Gamma_{i_n}$ the product of vector fields $\Gamma_{i_j}\in \mathcal{Z}$, $j=1,\dots,n$. \index{GammaI@$\Gamma^I$, product of admissible vector fields}

Vector fields $\Gamma$ have two relevant properties: they act like derivations on non-linear terms; they exactly commute with the linear part of both wave and Klein-Gordon equation.
This is the reason why we exclude of our consideration the scaling vector field $S=t\partial_t+\sum_j x_j \partial_j$, which is always considered in the so-called \textit{Klainerman vector fields' method} for the wave equation, as it does not commute with the Klein-Gordon operator.

We also introduce the energy of $(u_+,u_{-},v_+,v_{-})$ at time $t\ge 1$ as \index{E0@$E_0(t; u_\pm, v_\pm)$, energy}
\begin{equation}
E_0(t; u_\pm, v_\pm) :=\int \left(|u_+(t,x)|^2 + |u_{-}(t,x)|^2+ |v_+(t,x)|^2+ |u_{-}(t,x)|^2\right) dx,
\end{equation}
together with the generalized energies \index{En@$E_n(t;u_\pm, v_\pm)$, generalized energy} \index{Ek@$E^k_3(t; u_\pm, v_\pm)$, generalized energy}
\begin{subequations} \label{def_generalized_energy}
\begin{equation}
E_n(t;u_\pm, v_\pm):= \sum_{|\alpha|\le n} E_0(t; D^\alpha_x u_\pm, D^\alpha_xv _\pm), \quad \forall n\in\mathbb{N}, n\ge 3,
\end{equation}
and
\begin{equation}
E^k_3(t; u_\pm, v_\pm) :=  \sum_{\substack{|\alpha|+|I|\le 3 \\  |I|\le 3-k}} E_0(t; D^\alpha_x u^I_\pm; D^\alpha_x v^I_\pm), \quad 0\le k\le 2,
\end{equation}
\end{subequations}
where, for any multi-index $I$, \index{uIpm@$u^I_\pm$, wave components with admissible vector fields} \index{vIpm@$v^I_\pm$, Klein-Gordon components with admissible vector fields}
\begin{equation}  \label{def uIpm vIpm}
u^I_\pm := (D_t \pm |D_x|)\Gamma^Iu, \quad v^I_\pm := (D_t\pm \langle D_x\rangle)\Gamma^I v.
\end{equation}
Energy $E_n(t;u_\pm, v_\pm)$, for $n\ge 3$, is introduced with the aim of controlling the Sobolev norm $H^n$ of $u_\pm, v_\pm$ for large values of $n$.
The reason of dealing with $E^k_3(t;u_\pm, v_\pm)$ is, instead, to control the $L^2$ norm of $\Gamma^I u_\pm, \Gamma^I v_\pm$, for any general $\Gamma\in\mathcal{Z}$ and $|I|\le 3$. 
In particular, superscript $k$ indicates that we are considering only products $\Gamma^I$ containing at most $3-k$ vector fields in $\{\Omega, Z_m,m=1,2\}$. For instance, the $L^2$ norms of $\Omega^3 u_\pm, \Omega Z^2_1 v_\pm$ are bounded by $E^0_3(t;u_\pm, v_\pm)$ but not by $E^1_3(t;u_\pm, v_\pm)$, while the $L^2$ norms of $Z^2_1u_\pm, \partial_2\Omega Z_2v_\pm$ are controlled by both $E^1_3(t;u_\pm, v_\pm), E^0_3(t;u_\pm, v_\pm)$, etc.
The interest of distinguishing between $k=0,1,2$, is to take into account the different growth in time of the $L^2$ norm of such terms depending on the number of vector fields $\Omega, Z_m$ acting on $u_\pm, v_\pm$, as emerges from a-priori estimate \eqref{est: bootstrap E02}.

\begin{thm}[Bootstrap Argument] \label{Thm: bootstrap argument}
Let $K_1,K_2>1$ and $H^{\rho,\infty}$ be the space defined in \ref{def Sobolev spaces-NEW} $(iii)$.
There exist two integers $n\gg \rho$ sufficiently large, some $0<\delta\ll \delta_2 \ll \delta_1\ll \delta_0\ll 1$ small, two constants $A,B>1$ sufficiently large and $\varepsilon_0\in ]0,(2A+B)^{-1}[$ such that, for any $0<\varepsilon<\varepsilon_0$, if $(u,v)$ is solution to \eqref{wave KG system}-\eqref{initial data} on some interval $[1,T]$, for a fixed $T>1$, and $u_\pm, v_\pm$ defined in \eqref{def u+- v+-} satisfy:
\begin{subequations}\label{est: bootstrap argument a-priori est}
\begin{gather}
\|u_\pm (t,\cdot)\|_{H^{\rho+1,\infty}}+ \|\mathrm{R} u_\pm(t,\cdot)\|_{H^{\rho+1,\infty}}\le A\varepsilon t^{-\frac{1}{2}}, \label{est: bootstrap upm} \\ 
\|v_\pm (t,\cdot)\|_{H^{\rho,\infty}}\le A\varepsilon t^{-1},\label{est: boostrap vpm} \\
E_n(t; u_\pm, v_\pm)^\frac{1}{2} \le B\varepsilon t^\frac{\delta}{2}, \label{est: bootstrap Enn}\\ 
E^k_3(t;u_\pm, v_\pm)^\frac{1}{2}\le B\varepsilon t^{\frac{\delta_k}{2}}, \quad \forall\ 0\le k\le 2, \label{est: bootstrap E02}
\end{gather}
\end{subequations}
for every $t\in [1,T]$, then in the same interval they verify also
\begin{subequations} 
\begin{gather}
\|u_\pm (t,\cdot)\|_{H^{\rho+1,\infty}}+ \|\mathrm{R} u_\pm(t,\cdot)\|_{H^{\rho+1,\infty}}\le \frac{A}{K_1}\varepsilon t^{-\frac{1}{2}},\label{est:bootstrap enhanced upm} \\
\|v_\pm (t,\cdot)\|_{H^{\rho,\infty}}\le \frac{A}{K_1}\varepsilon t^{-1}, \label{est:bootstrap enhanced vpm}\\
E_n(t; u_\pm, v_\pm)^\frac{1}{2} \le \frac{B}{K_2}\varepsilon t^\frac{\delta}{2}, \label{est: bootstrap enhanced Enn}\\
E^k_3(t;u_\pm, v_\pm)^\frac{1}{2}\le \frac{B}{K_2}\varepsilon t^{\frac{\delta_k}{2}}, \quad \forall \ 0\le k\le 2.\label{est: boostrap enhanced E02}
\end{gather}
\end{subequations}
\end{thm}

The a-priori estimates on the uniform norm of $u_\pm, \mathrm{R}u_\pm, v_\pm$ made in the above theorem translate in terms of $u_\pm, v_\pm$ the sharp decay in time we expect for the solution $(u,v)$ to starting problem \eqref{wave KG system}. Indeed, from definitions \eqref{def u+- v+-} it appears that
\begin{gather*}
D_tu = \frac{u_+ + u_{-}}{2}, \quad D_x u= \mathrm{R}\left(\frac{u_+-u_{-}}{2}\right), \\
D_t v= \frac{v_+ + v_{-}}{2}, \quad v= \langle D_x\rangle^{-1}\left(\frac{v_+-v_{-}}{2}\right), 
\end{gather*}
so \eqref{est: bootstrap upm}, \eqref{est: boostrap vpm} imply
\begin{equation*}
\|\partial_{t,x}u(t,\cdot)\|_{H^{\rho,\infty}}\le A\varepsilon t^{-\frac{1}{2}}, \quad \|\partial_t v(t,\cdot)\|_{H^{\rho,\infty}}+\|v(t,\cdot)\|_{H^{\rho+1,\infty}}\le A\varepsilon t^{-1}.
\end{equation*}
Furthermore, the following quantity
\begin{equation*}
\|\partial_t u(t,\cdot)\|_{H^n}+\|\nabla_x u(t,\cdot)\|_{H^n} + \|\partial_t v(t,\cdot)\|_{H^n}+\|\nabla_x v(t,\cdot)\|_{H^n}+\|v(t,\cdot)\|_{H^n}
\end{equation*}
is equivalent to the square root of $E_n(t;u_\pm, v_\pm)$, which implies that the propagation of a-priori energy estimate \eqref{est: bootstrap Enn} is equivalent to the propagation of a certain estimate on the above Sobolev norms. For this reason, the propagation of the a-priori estimate on $E_n(t;u_\pm, v_\pm)$ and a local existence argument will imply theorem \ref{Thm: Main theorem}.

Before ending this section and going into the core of the subject, we briefly remind the general definition of \textit{null condition} for a multilinear form on $\mathbb{R}^{1+n}$ and a result by H\"{o}rmander (see \cite{hormander:the_analysis_III}).

\begin{defn}
A $k$-linear form $G$ on $\mathbb{R}^{1+n}$ is said to satisfy the \textit{null condition} if and only if, for all $\xi\in\mathbb{R}^n, \xi=(\xi_0, \dots, \xi_n)$ such that $\xi_0^2 - \sum_{j=1}^n \xi_j^2 =0$, 
\begin{equation} \label{wave null condition}
G(\underbrace{\xi,\dots,\xi}_{k}) = 0.
\end{equation}
\end{defn}

\textsc{Example:} The trilinear form $\xi_0^2 \xi_a - \displaystyle\sum_{j=1,2} \xi_j^2 \xi_a$ associated to $Q_0(v,\partial_a w)$ satisfies the null condition \eqref{wave null condition}, for any $a=0,1,2$. This is the most common example of null form.

\begin{lem}[H\"{o}rmander \cite{hormander:non-linear}, Lemma 6.6.5.]  \label{Lemma : Vector Field on a null form}
Let $G$ be a $k$-linear form on $\mathbb{R}^{1+n}$, $k=k_1 + \dots + k_r$, with $k_j$ positive integers, and $\Gamma\in \mathcal{Z}$. For all $u_j\in C^{k+1}(\mathbb{R}^{1+n})$, all $\alpha_j\in \mathbb{N}^{1+n}$, $|\alpha_j|=k_j$, and $u_j^{(k_j)} := \partial^{\alpha_j} u_j$, 
\begin{equation}
\begin{split}
\Gamma G(u_1^{(k_1)},\dots, u_r^{(k_r)}) & = G((\Gamma u_1)^{(k_1)},\dots, u_r^{(k_r)}) + \dots \\
& + G(u_1^{(k_1)},\dots, (\Gamma u_r)^{(k_r)}) + G_1(u_1^{(k_1)},\dots, u_r^{(k_r)})\,,
\end{split}
\end{equation}
where $G_1$ satisfies the \textit{null condition}.
\end{lem}

\begin{remark} \label{Remark:Vector_field_on_null_structure}
Previous lemma simplifies when the multi-linear form $G$ satisfying the null condition is $Q_0(v, \partial_a w)$, for any $a=0,1,2$. Indeed, the structure of the null form is not modified by the action of vector field $\Gamma$ in the sense that 
 \begin{equation}\label{Gamma_nonlinearity}
\Gamma Q_0(v,\partial_a w) = Q_0(\Gamma v, \partial_a w) + Q_0(v, \partial_a \Gamma w) + G_1(v,\partial w)\,,
\end{equation} 
where $G_1(v,\partial w)=0$ if $\Gamma=\partial_m$, $m=1,2$, and
\begin{equation}\label{def_G1}
G_1(v,\partial w) =
\begin{cases}
-Q_0(v, \partial_m w), \ &\text{if } a=0, \Gamma=Z_m, m\in \{1,2\},\\
0, \ &\text{if } a=0, \Gamma =\Omega,\\
-Q_0(v, \partial_t w), \ &\text{if } a\ne 0, \Gamma =Z_a, \\
0, \ &\text{if } a\ne0, \Gamma=Z_m, m\in \{1,2\}\setminus\{a\},\\
(-1)^a Q_0(v, \partial_m w), \ &\text{with } m\in \{1,2\}\setminus \{a\}, \text{ if } a\ne0, \Gamma =\Omega.
\end{cases}
\end{equation}
If $\Gamma^I$ contains at least $k\le |I|$ space derivatives then
\begin{equation} \label{Gamma_I_nonlinearity}
\Gamma^I Q_0(v,\partial_1 w) = \sum_{|I_1|+|I_2|=|I|} Q_0(\Gamma^{I_1} v, \partial_1\Gamma^{I_2}  w) + \sum_{k\le |I_1|+|I_2|<|I|}c_{I_1,I_2} Q_0(\Gamma^{I_1} v, \partial\Gamma^{I_2}  w),
\end{equation}
with $c_{I_1,I_2}\in \{-1,0,1\}$.
In the above equality we should think of multi-index $I_1$ (resp. $I_2$) as obtained by extraction of a $|I_1|$-tuple (resp. $|I_2|$-tuple) from $I=(i_1,\dots i_n)$, in such a way that each $i_j$ appearing in $I$ and corresponding to a spatial derivative (e.g. $\Gamma_{i_j}=D_m$, for $m\in\{1,2\}$), appears either in $I_1$ or in $I_2$, but not in both. For further references, we define \index{I(I)@$\mathcal{I}(I)$, set of multi-indices}
\begin{equation} \label{def mathcal(I)}
\mathcal{I}(I):=\left\{(I_1,I_2)| I_1, I_2\text{ multi-indices obtained as described above}\right\}.
\end{equation}
\end{remark}

\section{Preliminary Results}
The aim of this section is to introduce most of the technical tools that will be used throughout the paper. 
In particular, subsections \ref{Subsection: Paradifferential Calculus} and \ref{Subsection: Semiclassical Pseudodifferential Calculus} are devoted to recall some definitions and results about paradifferential and pseudo-differential calculus respectively; subsection \ref{Subsection: Some Technical Estimates I} and \ref{Subsection: Some Technical Estimates II} are dedicated to the introduction of some special operators that we will frequently use when dealing with the wave and the Klein-Gordon component. Subsections \ref{Subsection: Paradifferential Calculus}, \ref{Subsection: Semiclassical Pseudodifferential Calculus} barely contain proofs (we refer for that to \cite{bony:calcul_symbo}, \cite{metivier:paradifferential}, \cite{dimassi:spectral}, \cite{zworski:semiclassical}), whereas subsections \ref{Subsection: Some Technical Estimates I}, \ref{Subsection: Some Technical Estimates II} are much longer and richer in proofs and technicalities.


\subsection{Paradifferential calculus} \label{Subsection: Paradifferential Calculus}
In the current subsection we recall some definitions and properties that will be useful in chapter \ref{Chap:Energy estimates}. We first recall the definition of some spaces (Sobolev, Lipschitz and H\"older spaces) in dimension $d\ge 1$ and afterwards some results concerning symbolic calculus and the action of paradifferential operators on Sobolev spaces (see for instance \cite{metivier:paradifferential}).
We warn the reader that we will use both notations $\hat{w}(\xi)$ and $\mathcal{F}_{x\mapsto \xi}w$ for the Fourier transform of a function $w=w(x)$.

\begin{defn}[Spaces] \label{def Sobolev spaces-NEW}
\begin{itemize}
\item[(i)] Let $s \in \mathbb{R}$. $H^s(\mathbb{R}^d)$ denotes the space of tempered distributions $w\in \mathcal{S}'(\mathbb{R}^d)$ such that $\hat{w}\in L^2_\textit{loc}(\mathbb{R}^d)$ and \index{Hs@$H^s(\mathbb{R}^d)$, space}
\begin{equation*}
\|w\|^2_{H^s(\mathbb{R}^d)} := \frac{1}{(2\pi)^d}\int (1+|\xi|^2)^s |\hat{w}(\xi)|^2d\xi < +\infty;
\end{equation*}
\item[(ii)] For $\rho\in\mathbb{N}$, $W^{\rho,\infty}(\mathbb{R}^d)$ denotes the space of distributions $w\in\mathcal{D}'(\mathbb{R}^d)$ such that $\partial^\alpha_x w\in L^\infty(\mathbb{R}^d)$, for any $\alpha\in\mathbb{N}^d$ with $|\alpha|\le \rho$, endowed with the norm $$\|w\|_{W^{\rho,\infty}} := \sum_{|\alpha|\le\rho}\|\partial^\alpha_x w\|_{L^\infty};$$ \index{Wrho@$W^{\rho,\infty}(\mathbb{R}^d)$, space}
\item[(iii)] For $\rho\in\mathbb{N}$, we also introduce $H^{\rho,\infty}(\mathbb{R}^d)$ as the space of tempered distributions $w\in\mathcal{S}'(\mathbb{R}^d)$ such that \index{Hrho@$H^{\rho,\infty}(\mathbb{R}^d)$, space}
\begin{equation*}
\|w\|_{H^{\rho,\infty}} :=\|\langle D_x\rangle^\rho w\|_{L^\infty} < +\infty .
\end{equation*}
\end{itemize}
\end{defn}

\smallskip

\begin{defn}
An operator $T$ is said of order $\le m\in \mathbb{R}$ if it is a bounded operator from $H^{s+m}(\mathbb{R}^d)$ to $H^s(\mathbb{R}^d)$ for all $s\in\mathbb{R}$.
\end{defn}

\begin{defn}[Smooth symbols]
Let $m\in \mathbb{R}$.
\begin{itemize}
\item[$(i)$] $S^m_0(\mathbb{R}^d)$\index{Sm0@$S^m_0(\mathbb{R}^d)$, class of symbols} denotes the space of functions $a(x,\eta)$ on $\mathbb{R}^d\times \mathbb{R}^d$ which are $C^\infty$ with respect to $\eta$ and such that, for all $\alpha\in\mathbb{N}^d$, there exists a constant $C_\alpha>0$ and
\begin{equation*}
\|\partial^\alpha_\eta a(\cdot,\eta)\|_{L^\infty} \le C_\alpha (1+|\eta|)^{m-|\alpha|}, \qquad \forall \eta \in\mathbb{R}^d.
\end{equation*}
$\Sigma^m_0(\mathbb{R}^d)$ denotes the subclass of symbols $a\in S^m_0(\mathbb{R}^d)$ satisfying \index{Sigmam0@$\Sigma^m_0(\mathbb{R}^d)$, class of symbols}
\begin{equation} \label{spectral condition class Sigma^m}
\exists \varepsilon<1\ : \mathcal{F}_{x\rightarrow \xi}a(\xi,\eta)=0 \quad \text{for}\ |\xi|>\varepsilon (1+|\eta|).
\end{equation}
$S^m_0$ is equipped with seminorm $M^m_0(a;n)$ given by \index{Mm0@$M^m_0(a;n)$, seminorm}
\begin{equation} \label{def: seminorm Mmo}
M^m_0(a ;n) = \sup_{|\beta|\le n}\sup_{\eta\in\mathbb{R}^2}\big\|(1+|\eta|)^{|\beta|-m}\partial^\beta_\eta a(\cdot,\eta)\big\|_{L^{\infty}}.
\end{equation}
\item[$(ii)$] For $r\in\mathbb{N}$, $S^m_r(\mathbb{R}^d)$ denotes more generally the space of symbols $a\in S^m_0(\mathbb{R}^d)$ such that, for all $\alpha\in \mathbb{N}^d$ and all $\eta \in \mathbb{R}^d$, function $x\rightarrow \partial^\alpha_\eta a(x,\eta)$ belongs to $W^{r,\infty}(\mathbb{R}^d)$ and there exists a constant $C_\alpha>0$ such that\index{Smr@$S^m_r(\mathbb{R}^d)$, class of symbols}
\begin{equation*}
\|\partial^\alpha_\eta a(\cdot,\eta)\|_{W^{r,\infty}}\le C_\alpha (1+|\eta|)^{m-|\alpha|}, \qquad \forall\eta\in\mathbb{R}^d.
\end{equation*}
$\Sigma^m_r(\mathbb{R}^d)$ denotes the subclass of symbols $a\in S^m_r(\mathbb{R}^d)$ satisfying the spectral condition \eqref{spectral condition class Sigma^m}. \index{Sigmamr@$\Sigma^m_r(\mathbb{R}^d)$, class of symbols}
$S^m_r$ is equipped with seminorm $M^m_r(a;n)$, given by
\begin{equation} \label{def: seminorm Mmr}
M^m_r(a ;n) = \sup_{|\beta|\le n}\sup_{\eta\in\mathbb{R}^2}\big\|(1+|\eta|)^{|\beta|-m}\partial^\beta_\eta a(\cdot,\eta)\big\|_{W^{r,\infty}}.
\end{equation}
\end{itemize}
These definitions extend to matrix valued symbols $a\in S^m_r$ ($a\in \Sigma^m_r$), $m\in\mathbb{R}$, $r\in\mathbb{N}$. If $a\in S^m_r$ (resp. $a\in \Sigma^m_r$), it is said \textit{of order $m$}.
\end{defn}

\begin{defn} \label{Def: admissible cut off function}
An admissible cut-off function $\psi(\xi, \eta)$ is a $C^\infty$ function on $\mathbb{R}^d\times\mathbb{R}^d$ such that \index{Admissible cut-off function}
\begin{itemize}
\item[(i)] there are $0<\varepsilon_1<\varepsilon_2<1$ and
\begin{equation} \label{adm_cut-off: property 1 }
\begin{cases}
& \psi(\xi,\eta) =1, \quad \text{for} \ |\xi|\le \varepsilon_1(1+|\eta|)\\
& \psi(\xi,\eta)=0, \quad \text{for} \ |\xi|\ge \varepsilon_2(1+|\eta|);
\end{cases}
\end{equation} \label{adm_cut-off: property 2}
\item[(ii)] for all $(\alpha,\beta)\in\mathbb{N}^d\times\mathbb{N}^d$ there is a constant $C_{\alpha,\beta}>0$ such that
\begin{equation}
 |\partial^\alpha_\xi \partial^\beta_\eta \psi(\xi,\eta)|\le C_{\alpha,\beta} (1+|\eta|)^{-|\alpha| - |\beta|}, \quad \forall (\xi, \eta).
\end{equation}
\end{itemize}
\end{defn}

\textsc{Example}: If $\chi$ is a smooth cut-off function such that $\chi(z)=1$ for $|z|\le \varepsilon_1$ and is supported in the open ball $B_{\varepsilon_2}(0)$, with $0<\varepsilon_1<\varepsilon_2<1$, function $\psi(\xi,\eta):= \chi\big(\frac{\xi}{\langle\eta\rangle}\big)$ is an admissible cut-off function in the sense of definition \ref{Def: admissible cut off function}. We will only consider this type of admissible cut-off functions for the rest of the paper and refer (abusively) to $\chi$ itself as an admissible cut-off.

\begin{defn} \label{Def: Paradiff_operator}
Let $\chi$ be an admissible cut-off function and $a(x,\eta)\in S^m_r$, $m\in \mathbb{R}, r\in\mathbb{N}$. The \textit{Bony quantization} (or \textit{paradifferential quantization}) $Op^B(a(x,\eta))$ associated to symbol $a$ and acting on a test function $w$ is defined as \index{OpB@$Op^B$, para-differential operator}
\begin{equation*}
\begin{split}
& Op^B(a(x,\eta))w(x)  :=\frac{1}{(2\pi)^d} \int_{\mathbb{R}^d} e^{i x\cdot\eta}\sigma^\chi_a(x,\eta) \hat{w}(\eta) d\eta\,, \\
& \text{with}\ \sigma^\chi_a(x,\eta)  := \frac{1}{(2\pi)^d}\int_{\mathbb{R}^d} e^{i(x-y)\cdot\zeta}\chi\left(\frac{\zeta}{\langle\eta\rangle}\right)a(y,\eta) dyd\zeta\,.
\end{split}
\end{equation*}
\end{defn}

The operator defined above depends on the choice of the admissible cut-off function $\chi$. However, if $a\in S^m_r$ for some $m\in\mathbb{R}, r\in\mathbb{N}$, a change of $\chi$ modifies $Op^B(a)$ only by the addition of a $r$-smoothing operator (i.e. an operator which is bounded from $H^s$ to $H^{s+r}$, see \cite{bony:calcul_symbo}), so the choice of $\chi$ will be substantially irrelevant as long as we can neglect $r$-smoothing operators. For this reason, we will not indicate explicitly the dependence of $Op^B$ (resp. of $\sigma^\chi_a$) on $\chi$ to keep notations as light as possible.
Let us also observe that, with such a definition, the Fourier transform of $Op^B(a)w$ has the following simple expression
\begin{equation} \label{fourier transform of paradiff op}
\mathcal{F}_{x\rightarrow \xi}\Big(Op^B(a(x,\eta))w(x)\Big)(\xi) = \frac{1}{(2\pi)^d}\int \chi\left(\frac{\xi - \eta}{\langle\eta\rangle}\right)\hat{a}_y(\xi - \eta, \eta)\hat{w}(\eta) d\eta\,,
\end{equation}
where $\hat{a}_y(\xi , \eta):=\mathcal{F}_{y\rightarrow\xi}\big(a(y,\eta)\big)$, and the product of two functions $u,v$ can be developed as
\begin{equation} \label{dev in paraproduct}
uv = Op^B(u)v + Op^B(v)u + R(u,v),
\end{equation}
where remainder $R(u,v)$ writes on the Fourier side as
\begin{equation}  \label{fourier transform of R(v,w)}
\reallywidehat{R(u,v)}(\xi) = \frac{1}{(2\pi)^d} \int \left(1 - \chi\left(\frac{\xi-\eta}{\langle\eta\rangle}\right) - \chi\left(\frac{\eta}{\langle\xi - \eta\rangle}\right)\right) \widehat{u}(\xi - \eta) \widehat{v}(\eta)d\eta.
\end{equation}
We remark that frequencies $\eta$ and $\xi - \eta$  in the above integral are either bounded or equivalent, and $R(u,v) = R(v,u)$.
With the aim of having uniform notations, we introduce the operator $Op^B_R$\index{OpBR@$Op^B_R$, remainder para-differential operator} associated to a symbol $a(x,\eta)$ and acting on a function $w$ as 
\begin{equation} \label{operator OpB_R}
\begin{split}
& Op^B_R(a(x,\eta))w(x):=\frac{1}{(2\pi)^d}\int e^{ix\cdot\eta}\delta^{\chi}_a(x,\eta)\hat{w}(\eta)d\eta\,, \\
& \text{with}\ \delta^\chi_a(x,\eta):= \frac{1}{(2\pi)^d}\int e^{i(x-y)\cdot\zeta}\left(1 - \chi\left(\frac{\zeta}{\langle\eta\rangle}\right) - \chi\left(\frac{\eta}{\langle\zeta\rangle}\right)\right)a(y,\eta) dyd\zeta\,.
\end{split}
\end{equation}

For future references, we recall the definition of the Littlewood-Paley decomposition of a function $w$.

\begin{defn}[Littlewood-Paley decomposition] \index{Littlewood Paley decomposition} \label{Def: Littlewood Paley decomposition}
Let $\chi :\mathbb{R}^2\rightarrow [0,1]$ be a smooth decaying radial function, supported for $|x|\le 2-\frac{1}{10}$ and identically equal to 1 for $|x|\le 1+\frac{1}{10}$. Let also $\varphi(\xi) : =\chi(\xi) -\chi(2\xi)\in C^\infty_0(\mathbb{R}^2 \setminus \{0\})$, supported for $\frac{1}{2}<|\xi|<2$, and $\varphi_k(\xi):=\varphi(2^{-k}\xi)$ for all $k\in \mathbb{N}^*$, with the convention that $\varphi_0 : = \chi$.
Then $\sum_{k\in\mathbb{N}} \varphi(2^{-k}\xi) = 1$, and for any $w\in \mathcal{S}'(\mathbb{R}^d)$
\begin{equation}
w = \sum_{k\in \mathbb{N}} \varphi_k(D_x) w
\end{equation}
is the Littlewood-Paley decomposition of $w$.
\end{defn}

The following proposition is a classical result about the action of para-differential operators on Sobolev spaces (see \cite{bony:calcul_symbo} for further details).
Proposition \ref{Prop: Paradiff action with non smooth symbols and R(u,v)} shows, instead, that some results of continuity over $L^2$ hold also for operators whose symbol $a(x,\eta)$ is not a smooth function of $\eta$, and that map $(u,v)\mapsto R(u,v)$ is continuous from $H^{4,\infty}\times L^2$ to $L^2$.

\begin{prop}[Action]  \label{Prop : Paradiff action on Sobolev spaces-NEW} 
Let $m\in \mathbb{R}$.
For all $s\in\mathbb{R}$ and $a \in S^m_0$, $Op^B(a)$ is a bounded operator from $H^{s+m}(\mathbb{R}^d)$ to $H^s(\mathbb{R}^d)$.
In particular,
\begin{equation}
\|Op^B(a)w\|_{H^s} \lesssim M^m_0\Big(a ; \Big[\frac{d}{2}\Big]+1\Big) \|w\|_{H^{s+m}}.
\end{equation}
\end{prop}

\begin{prop}\label{Prop: Paradiff action with non smooth symbols and R(u,v)}
\begin{itemize}
\item[(i)]Let $a(x, \eta) = a_1(x)b(\eta)$, with $a_1\in L^\infty(\mathbb{R}^2)$ and $b(\eta)$ bounded, supported in some ball centred in the origin and such that $|\partial^\alpha b(\eta)|\lesssim_\alpha |\eta|^{-|\alpha|+1}$ for any $\alpha\in\mathbb{N}^2$ with $|\alpha|\ge 1$. Then $Op^B(a(x,\eta)): L^2 \rightarrow L^2$ is bounded and for any $w\in L^2(\mathbb{R}^2)$
\begin{equation*}
\|Op^B(a(x,\eta))w\|_{L^2}\lesssim \|a_1\|_{L^\infty} \|w\|_{L^2}.
\end{equation*}
The same result is true for $Op^B_R(a(x,\eta))$;
\item[(ii)] Map $(u,v)\in H^{4,\infty}\times L^2 \mapsto R(u,v)\in L^2$ is well defined and continuous.
\end{itemize}
\proof
As concerns $(i)$ we have that
\begin{equation*}
Op^B(a(x, \eta))w(x) = \int K(x-z, x-y) a_1(y)w(z) dy dz
\end{equation*}
with
\begin{equation*}
K(x,y) :=  \frac{1}{(2\pi)^4}\int e^{i x\cdot\eta + i y\cdot\zeta} \chi\Big(\frac{\zeta}{\langle\eta\rangle}\Big) b(\eta)d\eta  d\zeta
\end{equation*}
and $\chi$ is an admissible cut-off function.
After the hypothesis on $b$ we have that for every $\alpha,\beta\in\mathbb{N}^2$, 
\begin{gather*}
\Big|\partial^\beta_\zeta \Big[\chi\Big(\frac{\zeta}{\langle\eta\rangle}\Big)b(\eta)\Big]\Big|\lesssim \mathds{1}_{\{|\eta|\lesssim 1\}} |g_\beta(\zeta)|, \\
\Big|\partial^\alpha_\eta\partial^\beta_\zeta \Big[\chi\Big(\frac{\zeta}{\langle\eta\rangle}\Big)b(\eta)\Big]\Big|\lesssim \mathds{1}_{\{|\eta|\lesssim 1\}} |\eta|^{1-|\alpha|}|g_\beta(\zeta)|, \qquad |\alpha|\ge 1,
\end{gather*}
for some bounded and compactly supported functions $g_\beta$.
Lemma \ref{Lem_appendix: Kernel with 1 function} $(i)$ and corollary \ref{Cor_appendix: decay of integral operators} $(i)$ of appendix \ref{Appendix A} imply that $|K(x,y)|\lesssim |x|^{-1}\langle x\rangle^{-2}\langle y\rangle^{-3}$ for any $(x,y)$, and statement $(i)$ follows by
an inequality such as \eqref{ineq: norm L(dx) kernel} with $L=L^2$.

In order to prove assertion $(ii)$ we consider a cut-off function $\psi\in C^\infty_0(\mathbb{R}^2)$ equal to 1 in some closed ball $\overline{B_C(0)}$, for a $C\gg 1$, and decompose $R(u,v)$ as follows, using \eqref{fourier transform of R(v,w)}:
\begin{equation*}
R(u,v) = \int K_0(x-y, y-z) u(y)v(z)  dy dz 
 + \int K_1(x-y,y-z)[\langle D_x\rangle^4 u](y)v(z) dy dz,
\end{equation*} 
with
\begin{align*}
K_0(x,y) &= \frac{1}{(2\pi)^2}\int e^{ix\cdot \xi + i y\cdot \eta } \left( 1-\chi\left(\frac{\xi -\eta}{\langle\eta\rangle}\right) - \chi\left(\frac{\eta}{\langle\xi-\eta\rangle}\right)\right)\psi(\eta) d\xi d\eta, \\
K_1(x,y)&= \frac{1}{(2\pi)^2}\int e^{ix\cdot \xi + i y\cdot \eta } \left( 1-\chi\left(\frac{\xi -\eta}{\langle\eta\rangle}\right) - \chi\left(\frac{\eta}{\langle\xi-\eta\rangle}\right)\right)(1-\psi)(\eta)\langle\xi-\eta\rangle^{-4} d\xi d\eta.
\end{align*}

Since frequencies $\xi, \eta$ are both bounded on the support of $\left( 1-\chi\left(\frac{\xi -\eta}{\langle\eta\rangle}\right) - \chi\left(\frac{\eta}{\langle\xi-\eta\rangle}\right)\right)\psi(\eta)$, one can show through some integration by parts that $|K_0(x,y)|\lesssim \langle x\rangle^{-3}\langle y \rangle^{-3}$ for any $(x,y)$, to then deduce that
\begin{equation*}
\left\| \int K_0(x-y, y-z) u(y)v(z)  dy dz \right\|_{L^2(dx)}\lesssim \|u\|_{L^\infty}\|v\|_{L^2}.
\end{equation*} 

Kernel $K_1(x,y)$ can be split using a Littlewood-Paley decomposition as follows \small
\begin{equation*}
K_1(x,y) =\sum_{k\ge 1}\frac{1}{(2\pi)^2} \underbrace{\int e^{ix\cdot \xi + i y\cdot \eta } \left( 1-\chi\left(\frac{\xi -\eta}{\langle\eta\rangle}\right) - \chi\left(\frac{\eta}{\langle\xi-\eta\rangle}\right)\right)(1-\psi)(\eta) \varphi(2^{-k}\eta)\langle\xi-\eta\rangle^{-4} d\xi d\eta}_{K_{1,k}(x,y)},
\end{equation*}\normalsize
for a suitable $\varphi \in C^\infty_0(\mathbb{R}^2\setminus\{0\})$. 
On the support of $ \left( 1-\chi\left(\frac{\xi -\eta}{\langle\eta\rangle}\right) - \chi\left(\frac{\eta}{\langle\xi-\eta\rangle}\right)\right)(1-\psi)(\eta)\varphi(2^{-k}\eta)$, frequencies $\eta, \xi-\eta$ are either bounded or equivalent and of size $2^k$ (which implies in particular that $\langle\xi-\eta\rangle^{-4}\lesssim \langle\xi\rangle^{-3}\langle\eta\rangle^{-1}$).
After a change of coordinates and some integration by parts one can show that $|K_{1,k}(x,y)|\lesssim 2^k \langle x\rangle^{-3}\langle 2^ky\rangle^{-3}$, for any $k\ge 1$, and therefore that
\begin{align*}
&\left\| \int e^{i(x-y)\cdot \xi + i(y-z)\cdot \eta } K_1(x-y,y-z)[\langle D_x\rangle^4 u](y)v(z)  dy dz \right\|_{L^2(dx)}\\
\lesssim & \sum_{k\ge 1}2^k \left\| \int \langle x-y\rangle^{-3}\langle 2^k(y-z)\rangle^{-3}|\langle D_x\rangle^4 u(y)||w(z)| dydz \right\|_{L^2(dx)} \\
\lesssim & \sum_{k\ge 1} 2^k \int \langle y\rangle^{-3} \langle 2^k z\rangle^{-3} \|[\langle D_x\rangle^4u](\cdot-y) w(\cdot-y-z)\|_{L^2(dx)} dydz \lesssim \|u\|_{H^{4,\infty}}\|w\|_{L^2},
\end{align*}
which concludes the proof of statement $(ii)$.
\endproof
\end{prop}

The last results of this subsection are stated without proofs. All the details can be found in chapter 6 of \cite{metivier:paradifferential} (see theorems 6.1.1, 6.1.4, 6.2.1, 6.2.4).

\begin{prop}[Composition] \label{Prop: paradifferential symbolic calculus}
Consider $a\in S^m_r$, $b\in S^{m'}_r$, $r\in\mathbb{N}^*$, $m, m'\in \mathbb{R}$.

$(i)$ Symbol $a \sharp b := \displaystyle\sum_{|\alpha|<r }\frac{1}{\alpha !}\partial^\alpha_\xi a(x,\xi) D^\alpha_x b(x,\xi)$
is well defined in $\sum_{j<r}S^{m+m'-j}_{r-j}$;

$(ii)$ $Op^B(a)Op^B(b) - Op^B(a \sharp b)$ is an operator of order $\le m + m' -r$, and for all $s\in\mathbb{R}$, there exists a constant $C>0$ such that, for all $a\in S^m_r(\mathbb{R}^d)$, $b\in S^{m'}_r(\mathbb{R}^d)$, and $w\in H^{s+m+m'-r}(\mathbb{R}^d)$,
\begin{multline*}
\|Op^B(a)Op^B(b)w - Op^B(a\sharp b)w\|_{H^s} \\
\le C \big(M^m_r(a;n)M^{m'}_0(b;n_0) + M^m_0(a;n)M^{m'}_r(b;n_0)\big) \|w\|_{H^{s+m+m'-r}},
\end{multline*}
where $n_0=\left[\frac{d}{2}\right]+1$, $n=n_0 + r$.
Moreover, $Op^B(a)Op^B(b)-Op^B(a\sharp b)=\widetilde{\sigma}_r(x,D_x)$ with 
\begin{multline*}
\widetilde{\sigma}_r(x,\xi) = (\sigma_a\sharp \sigma_b)(x,\xi) -\sigma_{a\sharp b}(x,\xi)\\
+ \sum_{|\alpha| = r}\frac{1}{r! (2\pi)^d}\int e^{i(x-y)\cdot\zeta}\left(\int_0^1 \partial^\alpha_\xi \sigma_a(x,\xi+t\zeta)(1-t)^{r-1}dt\right) \theta(\zeta,\xi) D^\alpha_x\sigma_b(y,\xi) dyd\zeta 
\end{multline*}
with $\theta\equiv 1$ in a neighbourhood of the support of $\mathcal{F}_{y\mapsto \eta}\sigma_b(\eta,\xi)$.

These results extend to matrix valued symbols and operators.
\end{prop}
 
\begin{remark} \label{Remark:on_symbol_Bony_calculus}
If symbol $a(x,\xi)$ only depends on $\xi$ then $\sigma_a \sharp \sigma_b - \sigma_{a\sharp b}=0$ and $\widetilde{\sigma}_r$ reduces to the only integral term. Moreover,
\begin{equation} \label{Fourier transform composition symbol}
\mathcal{F}_{x\mapsto \eta}\widetilde{\sigma}_r(\eta,\xi) = \sum_{|\alpha|=r}\frac{1}{r!} \left(\int_0^1 \partial^\alpha_\xi a(\xi+t\eta)(1-t)^{[r]_+-1}dt\right) \chi\Big(\frac{\eta}{\langle\xi\rangle}\Big)\eta^\alpha \hat{b}_y(\eta,\xi),
\end{equation}
where $\chi\big(\frac{\eta}{\langle\xi\rangle}\big)$ is the admissible cut-off function defining $\sigma_b$.
\end{remark}

\begin{cor} \label{Cor : paradiff symbolic calculus at order 1}
For $d=2$ and all $s\in\mathbb{R}$, there exists a constant $C>0$ such that, for $a\in S^m_r, b\in S^{m'}_r$, $r\in\mathbb{N}^*$, and $w\in H^{s+m+m'-1}$,
\begin{multline*}
\| Op^B(a)Op^B(b)w - Op(ab)w\|_{H^s} \\
\le C \big(M^m_1(a;3)M^{m'}_0(b;2) + M^m_0(a;3)M^{m'}_1(b;2)\big) \|w\|_{H^{s+m+m'-1}}.
\end{multline*}
\end{cor}

\begin{prop}[Adjoint]\label{Paradiff Ajoint}
Consider $a\in S^m_r(\mathbb{R}^d)$, denote by $Op^B(a)^*$ the adjoint operator of $Op^B(a)$ and by $a^*(x,\xi)=\overline{a}(x,\xi)$ the complex conjugate of $a(x,\xi)$. 

$(i)$ Symbol $b(x,\xi): = \displaystyle\sum_{|\alpha|<r}\frac{1}{\alpha !}D^\alpha_x \partial^\alpha_\xi a^*(x,\xi)$
is well defined in $ \sum_{j<r}S^{m-j}_{r-j}$;

$(ii)$ Operator $Op^B(a)^* - Op^B(b)$ is of order $\le m-r$. 
Precisely, for all $s\in \mathbb{R}$ there is a constant $C>0$ such that, for all $a\in S^m_r(\mathbb{R}^d)$ and $w\in H^{s + m -r}(\mathbb{R}^d)$,
\begin{equation*}
\left\| Op^B(a)^*w - Op^B(b)w\right\|_{H^s} \le C M^m_r(a;n)\|w\|_{H^{s+m -r}},
\end{equation*}
with $n_0 = \big[\frac{d}{2}\big]+1$, $n = n_0 + r$.

These results extend to matrix valued symbols $a$, with $a^*(x,\xi)$ denoting the adjoint of matrix $a(x,\xi)$.
\end{prop}

\begin{cor} \label{Cor : paradiff ajoint at order 1}
For $d=2$ and all $s\in\mathbb{R}$, there exists a constant $C>0$ such that, for $a\in S^m_r$, $r\in\mathbb{N}^*$ and $w\in H^{s+m-1}$,
\begin{equation*}
\| Op^B(a)^*w - Op(a^*)w\|_{H^s} \le C M^m_1(a;3)\|w\|_{H^{s+m-1}}.
\end{equation*}
\end{cor}

\subsection{Semi-classical pseudodifferential calculus} \label{Subsection: Semiclassical Pseudodifferential Calculus}

In this subsection we recall some definitions and results about semi-classical symbolic calculus in general space dimension $d\ge 1$ which will be used in section \ref{Sec: development of the PDE system}. We refer the reader to \cite{dimassi:spectral} and \cite{zworski:semiclassical} for more details. 

\begin{defn}
An order function on $\mathbb{R}^d\times\mathbb{R}^d$ is a smooth map from $\mathbb{R}^d\times\mathbb{R}^d$ to $\mathbb{R}_+$ : $(x,\xi)\rightarrow M(x,\xi)$ such that there exist $N_0\in \mathbb{N}$, $C>0$ and for any $(x,\xi), (y,\eta)\in \mathbb{R}^d\times\mathbb{R}^d$ \index{Order function}
\begin{equation} \label{def ineq order function} 
M(y,\eta) \le C \langle x-y\rangle^{N_0} \langle\xi-\eta\rangle^{N_0} M(x, \xi) \, ,
\end{equation}
where $\langle x\rangle=\sqrt{1+|x|^2}$.
\end{defn}

\begin{defn}
Let \emph{M} be an order function on $\mathbb{R}^d\times\mathbb{R}^d$, $\delta,\sigma \ge 0$. One denotes by $S_{\delta, \sigma}(M)$ the space of smooth functions \index{Sdelta@$S_{\delta, \sigma}(M)$, class of symbols}
\begin{align*}
(x,\xi, h)  & \rightarrow a(x,\xi, h) \\
\mathbb{R}^d\times\mathbb{R}^d\times ]0,1] & \rightarrow \mathbb{C}
\end{align*}
satisfying for any $\alpha_1, \alpha_2 \in \mathbb{N}^d, k, N \in \mathbb{N}$ 
\begin{equation} \label{symbol in S delta beta M}
|\partial_x^{\alpha_1}\partial_{\xi}^{\alpha_2}(h\partial_h)^k a(x,\xi, h)| \lesssim M(x,\xi)\, h^{-\delta(|\alpha_1|+|\alpha_2|)}(1+\sigma h^{\sigma}|\xi|)^{-N}\, .
\end{equation}
\end{defn}

A key role in this paper will be played by symbols $a$ verifying \eqref{symbol in S delta beta M} with $M(x,\xi)= \langle \frac{x+f(\xi)}{\sqrt{h}} \rangle^{-N}$, for $N \in \mathbb{N}$ and a certain smooth function $f(\xi)$. This function $M$ is no longer an order function because of the term $h^{-\frac{1}{2}}$, but nevertheless we keep writing $a \in S_{\delta,\sigma}(\langle \frac{x+f(\xi)}{\sqrt{h}} \rangle^{-N})$.

\begin{defn}
In the semi-classical setting we say that $a(x,\xi, h)$ is a symbol \textit{of order $r$} if $a \in S_{\delta,\sigma}(\langle \xi \rangle^r)$, for some $\delta,\sigma \ge 0$.
\end{defn}

Let us observe that when $\sigma>0$ the symbol decays rapidly in $h^{\sigma}|\xi|$, which implies the following inclusion for $r\ge 0$:
\begin{equation*}
S_{\delta,\sigma}(\langle\xi\rangle^r) \subset h^{-\sigma r}S_{\delta, \sigma}(1).
\end{equation*}
This means that, up to a small loss in $h$, this type of symbols can be always considered as symbols of order zero.
In the rest of the paper we will not indicate explicitly the dependence of symbols on $h$, referring to $a(x,\xi,h)$ simply as $a(x,\xi)$.

\begin{defn} \label{Def: Weyl and standard quantization}
Let $a\in S_{\delta, \sigma}(M)$ for some order function $M$, some $\delta,\sigma\ge 0$.
\begin{enumerate}[(i)]
\item We can define the \emph{Weyl quantization} of $a$ to be the operator $\oph(a)=a^{w}(x, hD)$ acting on $u \in \mathcal{S}(\mathbb{R}^d)$ by the following formula: \index{Opwh@$\oph$, semi-classical Weyl quantization}
\begin{equation*}
\oph(a(x,\xi))u(x) = \frac{1}{(2\pi h)^d}\int_{\mathbb{R}^d}\int_{\mathbb{R}^d} e^{\frac{i}{h}(x-y)\cdot\xi} a\Big(\frac{x+y}{2}, \xi\Big)\, u(y)\; dy d\xi \, ;
\end{equation*} 
\item We define also the \emph{standard quantization} of $a$: \index{Oph@$Op_h$, standard semi-classical quantization}
\begin{equation*}
Op_h(a(x,\xi))u(x) = \frac{1}{(2 \pi h)^d} \int_{\mathbb{R}^d}\int_{\mathbb{R}^d} e^{\frac{i}{h}(x-y)\cdot\xi} a(x, \xi)\, u(y)\; dy d\xi \, .
\end{equation*}
\end{enumerate}
\end{defn}

It is clear from the definition that the two quantizations coincide when the symbol does not depend on $x$.
We also introduce a semi-classical version of Sobolev spaces on which the above operators act naturally. \index{Sobolev injection, semi-classical}

\begin{defn} \label{def of h-Sobolev spaces}
\begin{enumerate}[(i)] 
\item Let $\rho \in \mathbb{N}$.
We define the semi-classical Sobolev space $H^{\rho, \infty}_h(\mathbb{R}^d)$ as the space of tempered distributions $w$ such that $\langle hD \rangle^\rho w := Op_h(\langle\xi\rangle^\rho)w \in L^\infty$, endowed with norm \index{Hrhoh@$H^{\rho, \infty}_h(\mathbb{R}^d)$, space}
\begin{equation*}
\|w\|_{H^{\rho,\infty}_h}=\|\langle hD \rangle^\rho w\|_{L^\infty};
\end{equation*}
\item \label{def of H^s_h} 
Let $s\in\mathbb{R}$. We define the semi-classical Sobolev space $H^s_h(\mathbb{R}^d)$ as the space of tempered distributions $w$ such that $\langle hD \rangle^s w := Op_h(\langle\xi\rangle^s)w\in L^2$, endowed with norm \index{Hsh@$H^s_h(\mathbb{R}^d)$, space}
\begin{equation*}
\|w\|_{H^s}=\|\langle hD \rangle^s  w\|_{L^2}.
\end{equation*}
\end{enumerate}
\end{defn}

For future references, we write down the semi-classical Sobolev injection in space dimension 2: \index{Sobolev injection, semi-classical}
\begin{equation} \label{semi-classical Sobolev injection}
\|v_h\|_{H_h^{\rho, \infty}(\mathbb{R}^2)} \lesssim_{\sigma}h^{-1} \|v_h\|_{H^{\rho + 1+\sigma}_h(\mathbb{R}^2)}\, , \qquad \forall \sigma> 0\, .
\end{equation}

The following two propositions are stated without proof. They concern the adjoint and the composition of pseudo-differential operators. All related details are provided in chapter 7 of \cite{dimassi:spectral} or in chapter 4 of \cite{zworski:semiclassical}.

\begin{prop}[Self-Adjointness]
If $a(x,\xi)$ is a real symbol its Weyl quantization is self-adjoint, i.e. 
\begin{equation*}
\big(\oph(a)\big)^*=\oph(a)\, .
\end{equation*}
\end{prop}

\begin{prop}[Composition for Weyl quantization] \label{Prop : Composition for Weyl quantization}
Let $a, b \in \mathcal{S}(\mathbb{R}^d)$. Then
\begin{equation*}
\oph(a)\circ \oph(b) = \oph (a\sharp b) \, ,
\end{equation*}
where
\begin{equation} \label{a sharp b integral formula}
a \sharp b \,(x,\xi) := \frac{1}{(\pi h)^{2d}}\int_{\mathbb{R}^d}\int_{\mathbb{R}^d}\int_{\mathbb{R}^d}\int_{\mathbb{R}^d} e^{\frac{2i}{h}\sigma (y, \eta; \, z, \zeta)} a(x+z, \xi + \zeta) b(x+y, \xi +\eta) \; dy d\eta dz d\zeta ,
\end{equation}
and
\begin{equation*}
\sigma (y, \eta; \, z, \zeta) =  \eta\cdot z - y\cdot \zeta\, .
\end{equation*}
\end{prop}

It is often useful to derive an asymptotic expansion for $a \sharp b$, as it allows easier computations than integral formula \eqref{a sharp b integral formula}.
This expansion is usually obtained by applying the stationary phase argument when $a, b \in S_{\delta, \sigma}(M)$, $\delta \in [0,\frac{1}{2}[$ (as shown in \cite{zworski:semiclassical}).
Here we provide an expansion at any order even when one of two symbols belongs to $S_{\frac{1}{2},\sigma_1}(M)$ (still having the other in $S_{\delta,\sigma_2}(M)$ for $\delta<\frac{1}{2}$, and $\sigma_1,\sigma_2$ either equal or, if not, one of them equal to zero), whose proof is based on the Taylor development of symbols $a, b$, and can be found in the appendix of \cite{stingo:1D_KG} (for $d=1$).

\begin{prop} \label{Prop: a sharp b}
Let $M_1,M_2$ be two order functions and $a\in S_{\delta_1, \sigma_1}(M_1)$, $b\in S_{\delta_2, \sigma_2}(M_2)$, $\delta_1, \delta_2 \in [0,\frac{1}{2}]$, $\delta_1 + \delta_2< 1$, $\sigma_1, \sigma_2 \ge 0$ such that
\begin{equation} \label{beta in symbolic calculus}
\sigma_1=\sigma_2\ge 0 \qquad \mbox{or} \qquad \big[\sigma_1\ne\sigma_2  \,\mbox{and } \, \sigma_i=0\,,\sigma_j>0\,, i\ne j\in\{1,2\} \big]\,.
\end{equation}
Then $a \sharp b \in S_{\delta, \sigma}(M_1 M_2)$, where $\delta = \max\{\delta_1, \delta_2 \}$, $\sigma=\max\{\sigma_1,\sigma_2\}$.
Moreover,
\begin{equation} \label{a sharp b asymptotic formula}
a \sharp b =  \sum_{\substack{\alpha=(\alpha_1,\alpha_2) \\ |\alpha|=0, \dots, N-1}} \frac{(-1)^{|\alpha_1|}}{\alpha !}\Big(\frac{h}{2i}\Big)^{|\alpha|}\partial^{\alpha_1}_x\partial^{\alpha_2}_\xi a(x,\xi) \ \partial^{\alpha_2}_x\partial^{\alpha_1}_\xi b(x,\xi)+ r_N \,,
\end{equation}
where $r_N \in h^{N(1-(\delta_1 + \delta_2))}S_{\delta, \sigma}(M_1 M_2)$ and
\begin{multline} \label{r_N 1}
r_N(x,\xi) = \, \left(\frac{h}{2i}\right)^N \frac{N}{(\pi h)^{2d}}   \sum_{\substack{\alpha= (\alpha_1, \alpha_2)\\ |\alpha|= N}} \frac{(-1)^{|\alpha_1|}}{\alpha!}
\int_{\mathbb{R}^4}e^{\frac{2i}{h}(\eta\cdot z - y\cdot\zeta)}   \\
\times \Big(\int_0^1  \partial_x^{\alpha_1}\partial_{\xi}^{\alpha_2}a(x+tz, \xi +t\zeta)(1-t)^{N-1} dt\Big)
\partial_x^{\alpha_2}\partial_{\xi}^{\alpha_1}b(x+y, \xi + \eta)\, dy d\eta dz d\zeta \,,
\end{multline}
or
\begin{multline} \label{r_N 2}
r_N(x,\xi) = \, \left(\frac{h}{2i}\right)^N \frac{N}{(\pi h)^{2d}}   \sum_{\substack{\alpha= (\alpha_1, \alpha_2)\\ |\alpha|= N}} \frac{(-1)^{|\alpha_1|}}{\alpha!}
\int_{\mathbb{R}^4}e^{\frac{2i}{h}(\eta\cdot z - y\cdot\zeta)}   \partial_x^{\alpha_1}\partial_{\xi}^{\alpha_2}a(x+z, \xi +\zeta) \\
\times \Big(\int_0^1 \partial_x^{\alpha_2}\partial_{\xi}^{\alpha_1}b(x+ty, \xi + t\eta)(1-t)^{N-1} dt\Big)\, dy d\eta dz d\zeta \,.
\end{multline}
More generally, if $h^{N\delta_1}\partial^{\alpha}a \in S_{\delta_1,\sigma_1}(M^{N}_1)$, $h^{N\delta_2}\partial^{\alpha}b \in S_{\delta_2,\sigma_2}(M^{N}_2)$, for $|\alpha|=N$ and some order functions $M_1^N, M_2^N$, then $r_N \in h^{N(1-(\delta_1+\delta_2))}S_{\delta,\sigma}(M^N_1M^N_2)$.
\end{prop}

\begin{remark}\label{Remark:symbols_with_null_support_intersection}
From the previous proposition it follows that, if symbols $a\in S_{\delta_1,\sigma_1}(M_1)$, $b\in S_{\delta_2,\sigma_2}(M_2)$ are such that $\text{supp}a \cap \text{supp}b = \emptyset$, then $a\sharp b = O(h^\infty)$, meaning that, for every $N\in \mathbb{N}$, $a\sharp b =r_N$ with $r_N\in h^{N(1-(\delta_1+\delta_2))}S_{\delta,\sigma}(M_1M_2)$.
\end{remark}
\begin{remark} 
We draw the reader's attention to the fact that symbol $\sharp$ is used simultaneously in Bony calculus (see proposition \ref{Prop: paradifferential symbolic calculus}) and in Weyl semi-classical calculus (as in \eqref{a sharp b asymptotic formula}) with two different meaning. However, we avoid to introduce different notations as it will be clear by the context if we are dealing with the former or the latter one.
\end{remark}

The result of proposition \ref{Prop: a sharp b} and remark \ref{Remark:symbols_with_null_support_intersection} are still true even when one of the two order functions, or both, has the form $\langle\frac{x+f(\xi)}{\sqrt{h}}\rangle^{-1}$, for a smooth function $f(\xi)$, $\nabla f(\xi)$ bounded, as stated below (see the appendix of \cite{stingo:1D_KG}).

\begin{lem} \label{Lem : a sharp b}
Let $f(\xi):\mathbb{R}^d \rightarrow \mathbb{R}$ be a smooth function, with $|\nabla f(\xi)|$ bounded.
Consider $a \in S_{\delta_1,\sigma_1}(\langle \frac{x+f(\xi)}{\sqrt{h}}\rangle^{-m})$, $m\in\mathbb{N}$, and $b\in S_{\delta_2,\sigma_2}(M)$, for $M$ order function or $M(x,\xi)=\langle\frac{x+f(\xi)}{\sqrt{h}}\rangle^{-n}$, $n \in \mathbb{N}$, some $\delta_1 \in [0,\frac{1}{2}]$, $\delta_2 \in [0,\frac{1}{2}[$, $\sigma_1, \sigma_2\ge 0$ as in \eqref{beta in symbolic calculus}.
Then $a\sharp b \in S_{\delta,\sigma}(\langle \frac{x+f(\xi)}{\sqrt{h}}\rangle^{-m}M)$, where $\delta = \max\{\delta_1,\delta_2\}$, $\sigma = \max\{\sigma_1,\sigma_2\}$,  and the asymptotic expansion \eqref{a sharp b asymptotic formula} holds, with $r_N \in h^{N(1-(\delta_1+\delta_2))}S_{\delta,\sigma}(\langle \frac{x+f(\xi)}{\sqrt{h}}\rangle^{-m}M)$ given by \eqref{r_N 1} (or equivalently \eqref{r_N 2}). \\
More generally, if $h^{N\delta_1}\partial^{\alpha}a \in S_{\delta_1,\sigma_1}(\langle \frac{x+f(\xi)}{\sqrt{h}}\rangle^{-m'})$ and $h^{N\delta_2}\partial^{\alpha}b \in S_{\delta_2,\sigma_2}(M^N)$, $|\alpha|=N$, $M^N$ order function or $M^N(x,\xi)=\langle \frac{x+f(\xi)}{\sqrt{h}}\rangle^{-n'}$, for some $m',n' \in \mathbb{N}$, then remainder $r_N$ belongs to $h^{N(1-(\delta_1+\delta_2))}S_{\delta,\sigma}(\langle \frac{x+f(\xi)}{\sqrt{h}}\rangle^{-m'}M^N)$.
\end{lem}

\subsection{Semi-classical Operators for the Wave Solution: Some Estimates} \label{Subsection: Some Technical Estimates I}

From now on we place ourselves in space dimension $d=2$.
This technical subsection focuses on the introduction and the analysis of some particular operators that we will use when dealing with the wave component in the semi-classical framework (subsection \ref{Subsection : The Derivation of the Transport Equation}). 
More precisely, lemma \ref{Lemma on inequalities for Op(A)} will be often recalled to prove that some operator belongs to $\mathcal{L}(L^2 ; L^\infty)$ and compute its norm; propositions  \ref{Prop : continuity Op(gamma) L2 to L2}, \ref{Prop : continuity of Op(gamma1):X to L2} concern the continuity of some important operators like $\Gamma^{w,k}$ defined in \eqref{def of Gamma_wk}, while propositions \ref{Prop: L2 est of integral remainders}, \ref{Prop : Linfty est of integral remainders} are devoted to prove the continuity of some other operators often arising when considering the quantization of symbolic integral remainders.
Finally, lemmas \ref{Lemma : symbolic product development} and \ref{Lemma : on the enhanced symbolic product} deal with the development of some special symbolic products. While \ref{Lemma : symbolic product development} will be used several times throughout the paper, lemma \ref{Lemma : on the enhanced symbolic product} is stated explicitly on purpose to prove lemma \ref{Lem: preliminary on Op(e)}.

\begin{lem} \label{Lemma on inequalities for Op(A)}
There exists a constant $C>0$ such that, for any function $A(x,\xi)$ with $\partial^\alpha_x \partial^\beta_\xi A \in L^2(\mathbb{R}^2\times\mathbb{R}^2)$ for $|\alpha|, |\beta|\le 3$, and any function $w\in L^2(\mathbb{R}^2)$, 
\begin{equation} \label{inequalities |Op(A)|}
\left| \oph\big(A(x,\xi)\big)w(x) \right|\le C \|w\|_{L^2} \int_{\mathbb{R}^2} \langle x-y \rangle^{-3} \sum_{|\alpha|,|\beta|\le 3} \Big\|\partial_y^{\alpha}\partial_\xi^{\beta}\Big[A\Big(\frac{x+y}{2},h\xi\Big)\Big]\Big\|_{L^2(d\xi)} dy.
\end{equation}
Moreover, if $A(x,\xi)$ is compactly supported in $x$ there exists a smooth function, supported in a neighbourhood of $supp A$, such that
\begin{equation} \label{inequality Op(A), A compactly supported}
\left| \oph\big(A(x,\xi)\big)w(x) \right|\le C \|w\|_{L^2} \int_{\mathbb{R}^2} \Big|\theta'\Big(\frac{x+y}{2}\Big)\Big| \sum_{|\alpha|\le 3} \Big\|\partial_y^{\alpha}\Big[A\Big(\frac{x+y}{2},h\xi\Big)\Big]\Big\|_{L^2(d\xi)} dy.
\end{equation}
\proof
Let us prove the statement for $A\in \mathcal{S}(\mathbb{R}^2\times \mathbb{R}^2)$ and $w\in\mathcal{S}(\mathbb{R}^2)$.
The density of $\mathcal{S}(\mathbb{R}^2\times\mathbb{R}^2)$ into $\{A\in L^2(\mathbb{R}^2\times \mathbb{R}^2)| \partial^\alpha_x\partial^\beta_\xi A \in L^2(\mathbb{R}^2 \times \mathbb{R}^2), |\alpha|, |\beta|\le 3\}$ and of $\mathcal{S}(\mathbb{R}^2)$ into $L^2(\mathbb{R}^2)$ will then justify the definition of $\oph(A(x,\xi))w$ for $A$ and $w$ as in the statement, together with inequalities \eqref{inequalities |Op(A)|}, \eqref{inequality Op(A), A compactly supported}.

Using integration by parts, Cauchy-Schwarz inequality, and Young's inequality for convolutions, we can write the following:
\begin{align*} 
& |\oph(A(x,\xi))w(x)|  = \frac{1}{(2\pi)^2}\left|\int_{\mathbb{R}^4} e^{i(x-y)\cdot\xi} A\Big(\frac{x+y}{2},h\xi\Big)w(y)\ dy d\xi \right| \\
& = \frac{1}{(2\pi)^4} \left| \int_{\mathbb{R}^2} \widehat{w}(\eta)  \int_{\mathbb{R}^2}\int_{\mathbb{R}^2} e^{i(x-y)\cdot\xi + i y\cdot \eta} A\Big(\frac{x+y}{2},h\xi\Big) \ dyd\xi\ d\eta\right| \\
& = \frac{1}{(2\pi)^4}\left| \int_{\mathbb{R}^2} \widehat{w}(\eta)  \int_{\mathbb{R}^2}\int_{\mathbb{R}^2} \left(\frac{1-i(x-y)\cdot\partial_\xi}{1+|x-y|^2}\right)^3 \left(\frac{1+i(\xi - \eta)\cdot\partial_y}{1+|\xi - \eta|^2}\right)^3 e^{i(x-y)\cdot\xi + i y\cdot \eta} \right.\\
& \left. \hspace{10cm} \times A\Big(\frac{x+y}{2},h\xi\Big) \ dyd\xi\ d\eta\right| \\
& \lesssim  \int_{\mathbb{R}^2} \left| \hat{w}(\eta)\right|  \int_{\mathbb{R}^2}\int_{\mathbb{R}^2}\langle x- y\rangle^{-3} \langle \xi - \eta\rangle^{-3} \sum_{|\alpha|,|\beta|\le 3} \Big|\partial_y^{\alpha}\partial_\xi^{\beta}\Big[A\Big(\frac{x+y}{2},h\xi\Big)\Big]\Big| dy d\xi \ d\eta \\
& \lesssim \|\hat{w}\|_{L^2(d\eta)} \|\langle\eta\rangle^{-3}\|_{L^1(d\eta)}\int_{\mathbb{R}^2} \langle x-y \rangle^{-3} \sum_{|\alpha|,|\beta|\le 3} \Big\|\partial_y^{\alpha}\partial_\xi^{\beta}\Big[A\Big(\frac{x+y}{2},h\xi\Big)\Big]\Big\|_{L^2(d\xi)} dy \\
& \lesssim \|w\|_{L^2} \int_{\mathbb{R}^2} \langle x-y \rangle^{-3} \sum_{|\alpha|,|\beta|\le 3} \Big\|\partial_y^{\alpha}\partial_\xi^{\beta}\Big[A\Big(\frac{x+y}{2},h\xi\Big)\Big]\Big\|_{L^2(d\xi)} dy\,. 
\end{align*}
If symbol $A(x,\xi)$ is compactly supported in $x$ we can consider a smooth function $\theta'\in C^\infty_0(\mathbb{R})$, identically equal to 1 on the support of $A(x,\xi)$, and write
\begin{align*} 
& |\oph(A(x,\xi))w(x)|  = \frac{1}{(2\pi)^2}\left| \int_{\mathbb{R}^2} \widehat{w}(\eta) d\eta \int_{\mathbb{R}^2}\int_{\mathbb{R}^2} \left(\frac{1+i(\xi - \eta)\cdot\partial_y}{1+|\xi - \eta|^2}\right)^3 e^{i(x-y)\cdot\xi + i y\cdot \eta} \right.\\
& \left. \hspace{12cm} \times A\Big(\frac{x+y}{2},h\xi\Big) \ dyd\xi\right| \\
& \lesssim \int_{\mathbb{R}^2} \left| \hat{w}(\eta)\right| d\eta \int_{\mathbb{R}^2}\int_{\mathbb{R}^2}\Big|\theta'\Big(\frac{x+y}{2}\Big)\Big|\langle \xi - \eta\rangle^{-3} \sum_{|\alpha|\le 3} \big|\partial_y^{\alpha} \Big[A\Big(\frac{x+y}{2},h\xi\Big)\Big]\Big| dy d\xi  \\
& \lesssim \|w\|_{L^2} \int_{\mathbb{R}^2}\Big|\theta'\Big(\frac{x+y}{2}\Big)\Big| \sum_{|\alpha| \le 3} \big\|\partial_y^{\alpha}\Big[A\Big(\frac{x+y}{2},h\xi\Big)\Big]\Big\|_{L^2(d\xi)} dy\,. 
\end{align*}
\endproof
\end{lem}

A very important role in this subsection and in subsection \ref{Subsection : The Derivation of the Transport Equation} will be played by functions of the form $\gamma\big(\frac{x|\xi| -\xi}{h^{1/2 - \sigma}}\big)\psi(2^{-k}\xi)$, where $\gamma \in C^\infty(\mathbb{R}^2)$ is such that $|\partial^\alpha\gamma(z)|\lesssim \langle z\rangle^{-|\alpha|}$, $\psi\in C^\infty_0(\mathbb{R}^2-\{0\})$, $\sigma>0$ is a small fixed constant and $k$ is an integer belonging to set $K$, with \index{K@$K$, set of integers}
\begin{equation}\label{set_frequencies_K}
K :=\{k\in\mathbb{Z}\ :\ h\lesssim 2^k \lesssim h^{-\sigma}\}.
\end{equation}
In various results, such as proposition \ref{Prop : continuity of Op(gamma1):X to L2}, we will need a more decaying smooth function $\gamma_1$ verifying that $|\partial^\alpha\gamma_1(z)|\lesssim \langle z \rangle^{-(1+|\alpha|)}$.
We introduce here some notations we will keep throughout the whole paper:
\begin{notation}
For any $n\in \mathbb{N}$, $\gamma_n$ denotes a smooth function in $\mathbb{R}^2$ such that $|\partial^\alpha \gamma_n(z)|\lesssim_\alpha \langle z\rangle^{-(n+|\alpha|)}$, for any $\alpha\in\mathbb{N}^2$. We use the simplest notation $\gamma$ for $\gamma_0$; \index{gamman@$\gamma_n$, function}
\end{notation}
\begin{notation}\index{bm@$b_m(\xi)$, function}
For any integer $m\in \mathbb{Z}$, $b_m(\xi)$ will denote any function satisfying $|\partial^\beta b_m(\xi)|\lesssim_\beta |\xi|^{m-|\beta|}$, for any $\xi$ in its domain, any $\beta\in\mathbb{N}^2$. 
\end{notation}

The following lemma is a useful reference when we need to deal with some derivatives of $\gamma\big(\frac{x|\xi| - \xi}{h^{1/2-\sigma}}\big)$.

\begin{lem} \label{Lem : est on gamma for wave}
Let $\sigma\in\mathbb{R}$ and $n\in\mathbb{N}$.
For any multi-indices $\alpha, \beta\in \mathbb{N}^2$ we have that
\begin{equation} \label{derivatives of gamma_n 1}
\partial^\alpha_x \partial^\beta_\xi \Big[\gamma_n\Big(\frac{x|\xi| - \xi}{h^{1/2-\sigma}}\Big)\Big] = \sum_{k=0}^{|\beta|} h^{-(|\alpha| + k)(\frac{1}{2}-\sigma)}\gamma_{n+|\alpha|+k}\Big(\frac{x|\xi| - \xi}{h^{1/2-\sigma}}\Big) b_{|\alpha| +k -|\beta|}(\xi).
\end{equation}
Furthermore, if $\theta=\theta(x)\in C^\infty_0(\mathbb{R}^2)$, there exists a set $\{\theta_k(x)\}_{1\le k\le |\beta|}$ of smooth compactly supported functions such that
\begin{equation} \label{derivatives of gamma_n 2}
\theta(x)\partial^\alpha_x \partial^\beta_\xi \Big[\gamma_n\Big(\frac{x|\xi| - \xi}{h^{1/2-\sigma}}\Big)\Big] = \sum_{k=1}^{|\beta|} h^{-(|\alpha| + k)(\frac{1}{2}-\sigma)}\gamma_{n+|\alpha|+k}\Big(\frac{x|\xi| - \xi}{h^{1/2-\sigma}}\Big)\theta_k(x) b_{|\alpha| +k -|\beta|}(\xi).
\end{equation}
\proof
Let $\delta_{ij}$ be the Kronecker delta and ${\sum}'$ be a concise notation to indicate a linear combination.
For $i=1,2$,
\begin{equation} \label{first derivative gamma_n}
\begin{split}
& \partial_{\xi_i}\Big[\gamma_n\Big(\frac{x|\xi| -\xi}{h^{1/2-\sigma}}\Big)\Big]  = h^{-(\frac{1}{2}-\sigma)}\sum_{j=1}^2 (\partial_j\gamma_n)\Big(\frac{x|\xi| -\xi}{h^{1/2-\sigma}}\Big) (x_j\xi_i |\xi|^{-1} - \delta_{ij}) \\
& = \sum_{j=1}^2 (\partial_j\gamma_n)\Big(\frac{x|\xi| -\xi}{h^{1/2-\sigma}}\Big) \Big(\frac{x_j|\xi| - \xi_j}{h^{1/2-\sigma}}\Big)\ \xi_i|\xi|^{-2}  + \sum_{j=1}^2 h^{-(\frac{1}{2}-\sigma)} (\partial_j\gamma_n)\Big(\frac{x|\xi| -\xi}{h^{1/2-\sigma}}\Big) [\xi_i\xi_j|\xi|^{-2} - \delta_{ij}],
\end{split}
\end{equation}
which can be summarized saying that
\begin{equation*}
\partial_{\xi_i}\Big[\gamma_n\Big(\frac{x|\xi| -\xi}{h^{1/2-\sigma}}\Big)\Big] = {\sum}' \gamma_n\Big(\frac{x|\xi| -\xi}{h^{1/2-\sigma}}\Big) b_{-1}(\xi) + h^{-(\frac{1}{2}-\sigma)}\gamma_{n+1}\Big(\frac{x|\xi| -\xi}{h^{1/2-\sigma}}\Big) b_{0}(\xi), 
\end{equation*}
for some new functions $\gamma_n, \gamma_{n+1}, b_0, b_{-1}$. Iterating this argument one finds that, for all $\beta\in\mathbb{N}^2$,
\begin{equation*}
\partial^\beta_\xi \Big[\gamma_n\Big(\frac{x|\xi| - \xi}{h^{1/2-\sigma}}\Big)\Big] = {\sum_{k=0, \dots, |\beta|}}' h^{-k(\frac{1}{2}-\sigma)} \gamma_{n+k}\Big(\frac{x|\xi| - \xi}{h^{1/2-\sigma}}\Big) b_{k-|\beta|}(\xi),
\end{equation*}
and obtains \eqref{derivatives of gamma_n 1} using that, for any $m\in\mathbb{N}, \alpha\in\mathbb{N}^2$, 
\begin{equation} \label{derivatives in x gamma_n}
\partial^\alpha_x \Big[ \gamma_m\Big(\frac{x|\xi| - \xi}{h^{1/2-\sigma}}\Big)\Big] = h^{-|\alpha|(\frac{1}{2}-\sigma)} (\partial^\alpha \gamma_m)\Big(\frac{x|\xi| - \xi}{h^{1/2-\sigma}}\Big) |\xi|^{|\alpha|}.
\end{equation}
Equality \eqref{derivatives of gamma_n 2} is obtained replacing \eqref{first derivative gamma_n} with
\begin{equation*}
\begin{split}
\theta(x)\partial_{\xi_i}\Big[\gamma_n\Big(\frac{x|\xi| -\xi}{h^{1/2-\sigma}}\Big)\Big] & = h^{-(\frac{1}{2}-\sigma)}\sum_{j=1}^2 (\partial_j\gamma_n)\Big(\frac{x|\xi| -\xi}{h^{1/2-\sigma}}\Big) (\theta(x) x_j\xi_i |\xi|^{-1} - \theta(x)\delta_{ij})\\
& = {\sum}' h^{-(\frac{1}{2}-\sigma)}\gamma_{n+1}\Big(\frac{x|\xi| -\xi}{h^{1/2-\sigma}}\Big) \theta_1(x)b_{0}(\xi),
\end{split} 
\end{equation*}
where $\theta_1(x)$ is a new compactly supported function. By iteration one finds that, for any $\beta\in\mathbb{N}^2$, there is a set of $|\beta|$ compactly supported functions $\theta_k(x)$, $1\le k\le |\beta|$, such that
\begin{equation*}
\begin{split}
\theta(x) \partial^\beta_\xi \Big[\gamma_n\Big(\frac{x|\xi| - \xi}{h^{1/2-\sigma}}\Big)\Big] & =\sum_{k=1}^{|\beta|} h^{-k(\frac{1}{2}-\sigma)}\gamma_{n+k}\Big(\frac{x|\xi| - \xi}{h^{1/2-\sigma}}\Big)\theta_k(x) b_{k - |\beta|}(\xi),
\end{split}
\end{equation*}
which combined with \eqref{derivatives in x gamma_n} gives \eqref{derivatives of gamma_n 2}.
\endproof
\end{lem}

In some of the following results we denote by $\Theta_h$ the operator of change of coordinates
\[\Theta_hw(x) = w(\sqrt{h}x),\]
for any $h\in ]0,1]$, and use that for any symbol $a(x,\xi)$,
\begin{equation}\label{unitary_transf_Thetah}
\oph\big(a(x,\xi)\big) = \Theta_h \oph\big(\widetilde{a}(x,\xi)\big)\Theta_h^{-1},
\end{equation}
with $\widetilde{a}(x,\xi)=a\Big(\frac{x}{\sqrt{h}},\sqrt{h}\xi\Big)$.

\begin{prop}[Continuity on $L^2$] \label{Prop : continuity Op(gamma) L2 to L2}
Let $\sigma>0$ be sufficiently small, $K$ be the set defined in \eqref{set_frequencies_K}, $k\in K$ and $p\in\mathbb{Z}$. Let also $\psi\in C^\infty_0(\mathbb{R}^2\setminus\{0\})$ and $a(x)$ be a smooth function, bounded together with all its derivatives.
Then $\oph\big(\gamma\big(\frac{x|\xi|-\xi}{h^{1/2-\sigma}}\big)\psi(2^{-k}\xi)a(x)b_p(\xi)\big):L^2 \rightarrow L^2$ is bounded and 
\begin{equation} \label{est L2-L2 Op(gamma1)}
\Big\| \oph\Big(\gamma\Big(\frac{x|\xi|-\xi}{h^{1/2-\sigma}}\Big)\psi(2^{-k}\xi)a(x)b_p(\xi)\Big)\Big\|_{\mathcal{L}(L^2)} \lesssim 2^{kp}.
\end{equation}
\proof
Let $A(x,\xi) = \gamma\big(\frac{x|\xi|-\xi}{h^{1/2-\sigma}}\big)\psi(2^{-k}\xi)a(x)b_p(\xi)$.
For indices $k\in K$ such that $h^{1/2-\sigma}\lesssim 2^k \lesssim h^{-\sigma}$ the statement follows from the fact that $A(x,\xi)\in 2^{kp}S_{\frac{1}{2},0}(1)$ and by theorem 7.11 of \cite{dimassi:spectral}.
For $k\in K$ such that $h\le 2^k\le h^{1/2-\sigma}$, $\widetilde{A}(x,\xi) := A(\frac{x}{\sqrt{h}},\sqrt{h}\xi)\in  2^{kp}S_{\frac{1}{2},0}(1)$ and the result follows by theorem 7.11 of \cite{dimassi:spectral} and equality \eqref{unitary_transf_Thetah}.
\endproof
\end{prop}

\begin{prop} \label{Prop: L2 est of integral remainders}
Let $\sigma,k,p$ be as in the previous proposition. Let also $q\in \mathbb{Z}$, $\widetilde{\psi}\in C^\infty_0(\mathbb{R}^2\setminus\{0\})$, $a'(x)$ be a smooth function, bounded together with all its derivatives, and $f\in C(\mathbb{R})$.
Define
\begin{multline} \label{integral Ik}
I^k_{p,q}(x,\xi):=\frac{1}{(\pi h)^4} \int e^{\frac{2i}{h}(\eta\cdot z - y \cdot\zeta)}
\left[ \int_0^1 \Big(\gamma\Big(\frac{x|\xi|-\xi}{h^{1/2-\sigma}}\Big)\psi(2^{-k}\xi)a(x)b_p(\xi)\Big)|_{(x+tz,\xi +t\zeta)} f(t)dt \right. \\
\left. \times \widetilde{\psi}(2^{-k}(\xi + \eta))a'(x+y)b_q(\xi+\eta)\right]\ dydzd\eta d\zeta
\end{multline}
and
\begin{multline} \label{integral Jk}
J^k_{p,q}(x,\xi):=\frac{1}{(\pi h)^4} \int e^{\frac{2i}{h}(\eta\cdot z - y \cdot\zeta)}
\left[ \int_0^1 \widetilde{\psi}(2^{-k}(\xi + t\zeta))a'(x+tz)b_q(\xi+t\zeta) f(t)dt \right. \\ 
\left. \times \Big(\gamma\Big(\frac{x|\xi|-\xi}{h^{1/2-\sigma}}\Big)\psi(2^{-k}\xi)a(x)b_p(\xi)\Big)|_{(x+y, \xi + \eta)}\right]\ dydzd\eta d\zeta.
\end{multline}
Then $\oph(I^k_{p,q}(x,\xi))$ and $\oph(J^k_{p,q}(x,\xi))$ are bounded operators on $L^2$ and
\[\left\| \oph(I^k_{p,q}(x,\xi))\right\|_{\mathcal{L}(L^2)} +  \left\| \oph(J^k_{p,q}(x,\xi))\right\|_{\mathcal{L}(L^2)} \lesssim 2^{k(p+q)}.\]
The same results holds also if $q=0$ and $\widetilde{\psi}(2^{-k}\xi)b_q(\xi)\equiv 1$.

\proof
We show the result for $\oph(I^k_{p,q})$, leaving the reader to check that a similar argument can be used for $\oph(J^k_{p,q})$.

We distinguish between two ranges of frequencies.
For indices $k\in K$ such that $h^{1/2-\sigma}\le 2^k \lesssim h^{-\sigma}$ we observe that $I^k_{p,q}(x,\xi)\in 2^{k(p+q)}S_{\frac{1}{2}, 0}(1)$. 
Indeed, $\gamma\big(\frac{x|\xi|-\xi}{h^{1/2-\sigma}}\big)\psi(2^{-k}\xi)a(x)b_p(\xi)\in 2^{kp}S_{\frac{1}{2}, \sigma}(1)$ by lemma \ref{Lem : est on gamma for wave} while $\widetilde{\psi}(2^{-k}\xi)a'(x)b_q(\xi)\in 2^{kq}S_{\frac{1}{2}-\sigma,\sigma}(1)$. Hence performing a change of variables $y\mapsto \sqrt{h}y$, $z\mapsto \sqrt{h}z$, $\eta \mapsto\sqrt{h}\eta$, $\zeta \mapsto \sqrt{h}\zeta$, writing 
\begin{equation} \label{complex exponential}
e^{2i (\eta\cdot z - y\cdot\zeta)}= \left(\frac{1+2iy\cdot \partial_\zeta}{1+4|y|^2}\right)^3 \left(\frac{1-2i z\cdot \partial_\eta}{1+4|z|^2}\right)^3 \left(\frac{1-2i\eta\cdot\partial_z}{1+4|\eta|^2}\right)^3 \left(\frac{1+2i\zeta\cdot \partial_y}{1+4|\zeta|^2}\right)^3 e^{2i (\eta\cdot z - y\cdot\zeta)},
\end{equation}
and integrating by parts in all variables, we get that
\begin{equation*}
\left|I^k_{p,q}(x,\xi) \right| \lesssim 2^{k(p+q)} \int \langle y \rangle^{-3}\langle z\rangle^{-3} \langle \eta\rangle^{-3} \langle \zeta\rangle^{-3} \ dy dz d\eta d\zeta \lesssim 2^{k(p+q)},
\end{equation*}
without any loss in $h^{-\delta}$ due to the fact that we are considering symbols $A(x,\xi) \in S_{\delta, \sigma}(1)$ with $\delta\in \{0, 1/2-\sigma, 1/2\}$, and the derivatives of $A(x +\sqrt{h}y, \xi + \sqrt{h}\eta)$ (resp. of $A(x +t\sqrt{h}z, \xi +t\sqrt{h}\zeta)$) with respect to $y$ and $\eta$ (resp. with respect to $z$ and $\zeta$). In a similar way one can also prove that $|\partial^\alpha_x\partial^\beta_\xi I^k_{p,q}(x,\xi)|\lesssim_{\alpha,\beta} h^{-\frac{1}{2}(|\alpha| + |\beta|)} 2^{k(p+q)}$, for any $\alpha,\beta\in\mathbb{N}^2$.
Theorem 7.11 of \cite{dimassi:spectral} implies then the statement in this case. 

For indices $k\in K$ such that $h \lesssim 2^k \le h^{1/2-\sigma}$ we observe that 
\begin{gather*}
\gamma\Big(\frac{x|\xi|}{h^{1/2-\sigma}} - h^\sigma\xi\Big)\psi(2^{-k}\sqrt{h}\xi)a\Big(\frac{x}{\sqrt{h}}\Big)b_p(\sqrt{h}\xi) \in 2^{kp}S_{\frac{1}{2}, \sigma}(1), \\
\widetilde{\psi}(2^{-k}\sqrt{h}\xi)a'\Big(\frac{x}{\sqrt{h}}\Big)b_q(\sqrt{h}\xi) \in 2^{kq}S_{\frac{1}{2},\sigma}(1).
\end{gather*}
Then $\widetilde{I}^k_{p,q}(x,\xi)=I^k_{p,q}\Big(\frac{x}{\sqrt{h}},\sqrt{h}\xi\Big) \in 2^{k(p+q)}S_{\frac{1}{2},0}(1)$ and theorem 7.11 of \cite{dimassi:spectral} along with equality \eqref{unitary_transf_Thetah} imply that $\oph(I^k_{p,q}):L^2\rightarrow L^2$ is bounded, uniformly in $h$.

The last part of the statement can be proved following an analogous scheme, after having previously made an integration in $dzd\eta$ (or in $dyd\zeta$ if dealing with $J^k_{p,0}$).
\endproof
\end{prop}

\begin{prop}[Continuity on $L^p$] \label{Prop: Continuity_Lp_wave}
Let $1\le p\le +\infty$, $\gamma\in C^\infty_0(\mathbb{R}^2)$ be radial, $\psi\in C^\infty_0(\mathbb{R}^2\setminus \{0\})$, $a(x)$ be a smooth function, bounded together with all its derivatives. Let also $\sigma>0$ be small, $k\in K$ with $K$ given by \eqref{set_frequencies_K} and $q\in\mathbb{Z}$.
Then $\oph\big(\gamma\big(\frac{x|\xi|-\xi}{h^{1/2-\sigma}}\big)\psi(2^{-k}\xi)a(x)b_q(\xi)\big):L^p\rightarrow L^p$ is a bounded operator with
\begin{equation*}
\left\|\oph\Big(\gamma\Big(\frac{x|\xi|-\xi}{h^{1/2-\sigma}}\Big)\psi(2^{-k}\xi)a(x)b_q(\xi)\Big)\right\|_{\mathcal{L}(L^p)}\lesssim 2^{kq} .
\end{equation*}
\proof
In order to prove the result of the statement we need to show that kernel $K^k(x,y)$ associated to $\oph\big(\gamma\big(\frac{x|\xi|-\xi}{h^{1/2-\sigma}}\big)\psi(2^{-k}\xi)a(x)b_q(\xi)\big)$, i.e.
\begin{equation}\label{kernel_Kk}
K^k(x,y):=\frac{1}{(2\pi h)^2}\int e^{\frac{i}{h}(x-y)\cdot\xi} \gamma\Big(\frac{\big(\frac{x+y}{2}\big)|\xi|-\xi}{h^{1/2-\sigma}}\Big)\psi(2^{-k}\xi)a\Big(\frac{x+y}{2}\Big)b_q(\xi) d\xi,
\end{equation}
is such that
\begin{equation*}
\sup_x \int |K^k(x,y)| dy \lesssim 2^{kq}, \quad \sup_y \int |K^k(x,y)|dx \lesssim 2^{kq}.
\end{equation*}
From the symmetry between variables $x,y$, it will be enough to show that one of the two above inequalities is satisfied. 
To do that we study $K^k$ separately in different spatial regions, distinguishing also between indices $k\in K$ such that $2^k\le h^{1/2-\sigma}$ and $2^k>h^{1/2-\sigma}$.
We hence introduce three smooth cut-off functions $\theta_s, \theta_b, \theta$, supported respectively for $|x|\le m\ll 1$, $|x|\ge M\gg 1$, $0<m'\le |x|\le M'<+\infty$, for some constants $m, m', M, M'$, and such that $\theta_s+\theta_b+\theta\equiv 1$.
Denoting concisely by $A^k(x,\xi)$ the multiplier in \eqref{kernel_Kk}, we split it as follows
\begin{equation*}
A^k(x,\xi) = A^k_s(x,\xi) + A^k_b(x,\xi)+A^k_1(x,\xi),
\end{equation*}
with $A^k_s(x,\xi):=A^k(x,\xi)\theta_s(x)$, $A^k_b(x,\xi):=A^k(x,\xi)\theta_b(x)$ and $A^k_1(x,\xi):=A^k(x,\xi)\theta(x)$.

\smallskip
\underline{\textbf{Case I:}}
Let us consider $k\in K$ such that $h\lesssim 2^k\le h^{1/2-\sigma}$.
According to the above decomposition we have that
\begin{equation*}
K^k(x,y)=K^k_s(x,y)+K^k_b(x,y)+K^k_1(x,y),
\end{equation*}
with clear meaning of kernels $K^k_s,K^k_b,K^k_1$.
Let us first prove that
\begin{equation} \label{L1norm_Ks_Kb}
\sup_x \int |K^k_s(x,y)| dy +  \sup_x \int |K^k_b(x,y)| dy \lesssim 2^{kq}.
\end{equation}
For $|x|\ll 1$ (resp. $|x|\gg 1$), $\big|\frac{x|\xi| - \xi}{h^{1/2-\sigma}}\big|\gtrsim h^{-1/2+\sigma}|\xi|$ (resp. $\big|\frac{x|\xi| - \xi}{h^{1/2-\sigma}}\big|\gtrsim h^{-1/2+\sigma}|\xi||x|\gtrsim h^{-1/2+\sigma}|\xi|$) so by lemma \ref{Lem : est on gamma for wave}
\begin{equation}\label{derivative_xi_gamma}
\left| \partial^\beta_\xi\Big[\gamma\Big(\frac{x|\xi| - \xi}{h^{1/2-\sigma}}\Big)\Big]\right|\lesssim \sum_{j=0}^{|\beta|}h^{-j(\frac{1}{2}-\sigma)}\left\langle \frac{x|\xi|-\xi}{h^{1/2-\sigma}}\right\rangle^{-j}|b_{j-|\beta|}(\xi)| \lesssim |\xi|^{-|\beta|}.
\end{equation}
Therefore
\begin{equation} \label{partial_beta_Aks}
\left|\partial^\beta_\xi A^k_s(x, 2^k\xi)\right|\lesssim \sum_{|\beta_1|\le |\beta|} 2^{k|\beta|}|2^k\xi|^{-|\beta_1|}2^{-k(|\beta|-|\beta_1|)+kq}\mathds{1}_{|\xi|\sim 1} \lesssim 2^{kq}\mathds{1}_{|\xi|\sim 1},
\end{equation}
so making a change of coordinates $\xi\mapsto 2^k\xi$ and some integration by parts we derive that
\begin{equation*}
|K^k_s(x,y)|\lesssim 2^{kq}(2^kh^{-1})^2\left\langle 2^kh^{-1} (x-y)\right\rangle^{-3},
\end{equation*}
for every $(x,y)\in\mathbb{R}^2\times \mathbb{R}^2$. The same argument applies to $K^k_b(x,y)$, hence taking the $L^1$ norm we obtain \eqref{L1norm_Ks_Kb}.

As concerns kernel $K^k_1(x,y)$, we deduce from lemma \ref{Lem : est on gamma for wave}, the fact that $\theta_1(x)$ is supported for $|x|\sim 1$, and that $2^k\lesssim h^{1/2-\sigma}$, the following inequality:
\begin{equation*}
\left| \partial^\beta_\xi \Big[A^k_1\Big(\frac{x+y}{2}, 2^k\xi\Big)\Big]\right|\lesssim 2^{k|\beta|}\Big[2^{k(q-|\beta|)}+\sum_{j=1}^{|\beta|}h^{-j(\frac{1}{2}-\sigma)}|b_{j-|\beta|+q}(2^k\xi)|\Big]\lesssim 2^{kq}.
\end{equation*}
Performing a change of coordinates $\xi\mapsto 2^k\xi$ and making some integration by parts one finds that
\begin{equation*}
|K^k_1(x,y)|\lesssim 2^{kq}(2^kh^{-1})^2 \left\langle 2^k h^{-1}(x-y) \right\rangle^{-3}, \qquad \forall (x,y),
\end{equation*}
and consequently that
\begin{equation*}
\sup_x \int |K^k(x,y)|dy\lesssim 2^{kq}.
\end{equation*}
Summing up with \eqref{L1norm_Ks_Kb}, this gives us that
\begin{equation*}
\oph(A^k(x,\xi))=\oph(A^k_s(x,\xi))+\oph(A^k_b(x,\xi))+\oph(A^k_1(x,\xi))
\end{equation*}
is a bounded operator on $L^p$, for every $1\le p\le +\infty$, with norm $O(2^{kq})$.

\underline{\textbf{Case II:}}
Let us now suppose that $k\in K$ is such that $h^{1/2-\sigma}< 2^k\le h^{-\sigma}$.
From \eqref{derivative_xi_gamma} we have that $\widetilde{A}^k_s(x,\xi) = A^k_s\big(\frac{x}{\sqrt{h}},\sqrt{h}\xi\big)$ satisfies
\begin{equation*}
\left|\partial^\beta_\xi \widetilde{A}^k_s(x,\xi)\right| \lesssim \sum_{|\beta_1|\le |\beta|} h^{\frac{|\beta|}{2}}|\sqrt{h}\xi|^{-|\beta_1|} 2^{-k(|\beta|-|\beta_1|)+kq} \mathds{1}_{|\xi|\sim 2^k h^{-1/2}},
\end{equation*}
for every $(x,\xi)\in\mathbb{R}^2\times\mathbb{R}^2$, and hence
\begin{equation*}
\left| \partial^\beta_\xi \widetilde{A}^k_s(x,2^kh^{-1/2}\xi)\right|\lesssim \sum_{|\beta_1|\le |\beta|}2^{k|\beta|}|2^k\xi|^{-|\beta_1|}2^{-k(|\beta|-|\beta_1|)+kq}\mathds{1}_{|\xi|\sim1}\lesssim 2^{kq}\mathds{1}_{|\xi|\sim1}.
\end{equation*}
By making a change of coordinates $\xi\mapsto 2^kh^{-1/2}\xi$, some integrations by parts and using the above inequality, one can show that kernel $\widetilde{K}^k_s(x,y)$ associated to $\oph(\widetilde{A}^k_s(x,\xi))$, i.e.
\begin{equation*}
\widetilde{K}^k_s(x,y)=\frac{1}{(2\pi h)^2}\int e^{\frac{i}{h}(x-y)\cdot\xi} \widetilde{A}^k_s\Big(\frac{x+y}{2},\xi\Big) d\xi,
\end{equation*}
is such that
\begin{equation*}
|\widetilde{K}^k_s(x,y)|\lesssim 2^{kq}(2^kh^{-\frac{3}{2}})^2\left\langle 2^kh^{-\frac{3}{2}}(x-y)\right\rangle^{-3} \qquad\forall (x,y),
\end{equation*}
which implies that $sup_x \int |\widetilde{K}^k_s(x,y)| dy\lesssim 2^{kq}$. The same argument and result hold for $\widetilde{K}^k_b(x,y)$ so $\oph(A^k_s)$ and $\oph(A^k_b)$ verify the statement.

The last thing to prove is that $\oph(A_1(x,\xi))\in\mathcal{L}(L^p)$ for every $1\le p\le +\infty$. Let $K^k_1(x,y)$ be its associated kernel, i.e.
\begin{equation}
K^k_1(x,y) =\frac{1}{(2\pi h)^2}\int e^{\frac{i}{h}(x-y)\cdot\xi} \gamma\Big(\frac{\big(\frac{x+y}{2}\big)|\xi|-\xi}{h^{1/2-\sigma}}\Big) \psi(2^{-k}\xi)a\Big(\frac{x+y}{2}\Big)b_q(\xi) d\xi,
\end{equation}
and assume, without loss of generality, that $\gamma(x)=\gamma(|x|^2)$. Set
\begin{equation*}
\frac{x+y}{2} = r[\cos\theta, \sin\theta],
\end{equation*}
with $m'\le r\le M'$ on the support of $\theta_1\big(\frac{x+y}{2}\big)$, and for fixed $r,\theta$ let
\begin{equation} \label{change_of_coordinates}
\xi= \rho[\cos\theta, \sin\theta] + r\Omega [-\sin\theta, \cos\theta].
\end{equation}
We immediately notice that $[\frac{\partial(\xi_1,\xi_2)}{\partial(\rho,\Omega)}]=r\sim 1$ and that $|\xi|^2=\rho^2+r^2\Omega^2$. Moreover,
\begin{equation*}
\Big|\Big(\frac{x+y}{2}\Big)|\xi|-\xi \Big|^2 = \left[r\sqrt{\rho^2 + r^2\Omega^2}-\rho\right]^2 + r^2\Omega^2.
\end{equation*} 
If the support of $\gamma$ is of size $0<\alpha\ll 1$ sufficiently small, from the above equality and the fact that $|\xi|\sim 2^k$ on the support of $\psi(2^{-k}\xi)$, with $h^{1/2-\sigma}< 2^k\lesssim h^{-\sigma}$, we deduce that
\begin{equation*}
r\Omega \le \sqrt{\alpha}h^{1/2-\sigma} \quad \text{and} \quad |\rho|\sim |\xi|\sim 2^k \quad \text{and} \quad \frac{r\Omega}{|\rho|}\le \sqrt{\alpha}.
\end{equation*}
Consequently 
\begin{equation*}
\alpha h^{1-2\sigma}\ge \left[r\sqrt{\rho^2 + r^2\Omega^2}-\rho\right]^2 \gtrsim \rho^2|r-1|^2.
\end{equation*}
The above left inequality implies that $\rho>0$, inferring so the right one.
Moreover
\begin{equation*}
\begin{split}
\alpha h^{1-2\sigma}\ge \left[r\sqrt{\rho^2 + r^2\Omega^2}-\rho\right]^2 +r^2\Omega^2 &= \rho^2 \left[(r-1) + r \left[\sqrt{1+\frac{r^2\Omega^2}{\rho^2}}-1\right]\right]^2 + r^2\Omega^2\\
& = \rho^2|r-1|^2 + r^2\Omega^2\left[1+ a(r,\Omega,\rho)\right],
\end{split}
\end{equation*} 
with $a(r,\Omega,\rho)$ bounded such that, for any $l,m,n\in\mathbb{N}$,
\begin{equation*}
|\partial^l_r \partial^m_\Omega \partial^n_\rho a(r,\Omega, \rho)|=O(\rho^{-(m+n)}).
\end{equation*}
If 
\begin{equation*}
\Gamma_h:=\gamma\Big(\frac{ \rho^2|r-1|^2}{h^{1-2\sigma}} + \frac{r^2\Omega^2}{h^{1-2\sigma}}\left[1+ a(r,\Omega,\rho)\right]\Big) \psi(2^{-k}\sqrt{\rho^2+r^2\Omega^2}) a(r,\theta)b_q(\rho),
\end{equation*}
from all the observations made above along with the fact that $h^{-1/2+\sigma}\lesssim \rho^{-1}$ we deduce that, for any $m,n\in\mathbb{N}$,
\begin{equation}\label{deriv_Gammah}
\left|\partial^m_\rho \Gamma_h\right| =O(2^{kq}\rho^{-m}) \quad \text{and} \quad \left|\partial^n_\Omega \Gamma_h\right|=O(2^{kq}\rho^{-n}).
\end{equation}
With the change of coordinates considered in \eqref{change_of_coordinates}, and setting $w:=x-y$, $e_\theta:=[\cos\theta, \sin\theta]$, kernel $K^k_1(x,y)$ transforms into
\begin{equation*}
\frac{1}{(2\pi h)^2}\int e^{\frac{i}{h}\rho w\cdot e_\theta + \frac{i}{h}r\Omega w\cdot e^\perp_\theta}\Gamma_h \, rd\rho d\Omega
\end{equation*}
and is restricted to $|\rho|\sim 2^k$, $|\Omega|\lesssim h^{1/2-\sigma}$, so by making some integrations by parts, using \eqref{deriv_Gammah}, and reminding that $|r-1|\ll 2^{-k}h^{1/2-\sigma}\ll 1$ on the support of $\Gamma_h$, we find that, for any $N\in\mathbb{N}$,
\begin{equation*}
|K^k_1(x,y)|\lesssim h^{-\frac{3}{2}-\sigma}2^k \left\langle \frac{2^k}{h}w\cdot e_\theta \right\rangle^{-N} \left\langle\frac{2^k}{h}w\cdot e^\perp_\theta \right\rangle^{-N} \mathds{1}_{||\frac{x+y}{2}|-1|\ll 1}.
\end{equation*}
Now, as $w=(x-y)$, $e_\theta = \frac{x+y}{|x+y|}$, and $|x+y|=2r\sim 1$ on the support of $\Gamma_h$, we have that $|w\cdot e_\theta|\sim ||x|^2-|y|^2|$, $|w\cdot e^\perp_\theta|\sim |(x+y)(x+y)^\perp|\sim 2|x\cdot y^\perp| = 2|x_1 y_2-x_2y_1|$, and consequently
\begin{equation*}
|K^k_1(x,y)|\lesssim h^{-\frac{3}{2}-\sigma}2^{k(1+q)} \left\langle \frac{2^k}{h}\big||x|^2-|y|^2\big| \right\rangle^{-N} \left\langle\frac{2^k}{h}(x_1y_2 - x_2y_1)\right\rangle^{-N} \mathds{1}_{||\frac{x+y}{2}|-1|\ll 1}.
\end{equation*}
Successively, taking the $L^1(dy)$ norm of $K^k_1(x,y)$ and using the above estimate we find that:

$\bullet$ if $|x|\ll |y|$ or $|x|\gg |x|$,
\begin{equation*}
\left\langle \frac{2^k}{h}\big||x|^2-|y|^2\big| \right\rangle^{-N} \mathds{1}_{||\frac{x+y}{2}|-1|\ll 1}\lesssim h^{N(\frac{1}{2}+\sigma)},
\end{equation*}
as follows from the fact that $h2^{-k}< h^{1/2+\sigma}$. We obtain that 
\begin{equation*}
\sup_x\int |K^k_1(x,y)|dy \lesssim h^{-\frac{3}{2}}2^{k(1+q)} h^{N(\frac{1}{2}+\sigma)} \lesssim 1
\end{equation*}
by taking $N\in\mathbb{N}$ sufficiently large (e.g. $N>3$) and $\sigma>0$ small.

$\bullet$ if $|x|\sim |y|$, we deduce that $|x|\ge c>0$ from the fact that $\big| \big|\frac{x+y}{2}\big|-1\big|\le \sqrt{\alpha}h^{1/2-\sigma}2^{-k}$ on the support of $\Gamma_h$. Without loss of generality we can assume that $x=\lambda e_1$ (this always being possible by making a rotation) and $|\lambda|\ge c>0$. 
If $w:=x+y$, 
\begin{equation*}
|x|^2-|y|^2 = w\cdot(x-y) = w\cdot (2x-w)=w\cdot (2\lambda e_1 -w) = 2\lambda w_1-w_1^2-w_2^2,
\end{equation*}
and then
\begin{equation*}
\frac{\big||x|^2-|y|^2\big|}{h} = -\frac{(w_1-\lambda)^2-\lambda^2}{h}+ \Big(\frac{w_2}{\sqrt{h}}\Big)^2
\end{equation*}
while
\begin{equation*}
x_1y_2-x_2y_1 = \lambda w_2.
\end{equation*}
This implies that
\begin{equation*}
|K^k_1(x,y)|\lesssim h^{-\frac{3}{2}-\sigma}2^{k(1+q)}\left\langle \frac{2^k}{h}\left((w_1-\lambda)^2-\lambda^2\right) \right\rangle^{-N} \left\langle  \frac{2^k}{h} w_2\right\rangle^{-N}.
\end{equation*}
Since $\int |K^k_1(x,y)| dy = \int |K^k_1(x,y)| dw$, from the above estimate (with a fixed $N\in\mathbb{N}$ sufficiently large) this integral is bounded by $2^{kq}$ when restricted to $|x|\sim |y|$. Indeed, when the integral is taken in a neighbourhood of $w_1=0$ or $w_1=2\lambda$, $(w_1-\lambda^2)-\lambda^2$ can be considered as the variable of integration, and by a change of coordinates along with the fact that $2^{-k}<h^{-1/2+\sigma}$ one deduces that
\begin{equation*}
\int_{U_0\cup U_{2\lambda}} |K^k_1(x,y)| dw \lesssim h^{-\frac{3}{2}-\sigma}2^{k(1+q)}h^22^{-2k}\lesssim 2^{kq},
\end{equation*}
where $U_0$ (resp. $U_{2\lambda}$) is a neighbourhood of $w_1=0$ (resp. of $w_1=2\lambda$). Outside of $U_0\cup U_{2\lambda}$, 
\begin{equation*}
\left\langle \frac{2^k}{h}\left((w_1-\lambda)^2-\lambda^2\right) \right\rangle^{-N}\lesssim (h2^{-k})^N \langle w_1\rangle^{-N}\lesssim h^{N(\frac{1}{2}+\sigma)}\langle w_1\rangle^{-N},
\end{equation*}
so
\begin{equation*}
\int_{(U_0\cup U_{2\lambda})^\complement} |K^k_1(x,y)| dw \lesssim h^{-\frac{3}{2}-\sigma}2^{k(1+q)}h2^{-k}h^{N(\frac{1}{2}+\sigma)}\lesssim 2^{kq}.
\end{equation*}

This finally proves that also $\oph(A^k_1(x,\xi))$ is a bounded operator on $L^p$ with norm $O(2^{kq})$.
\endproof
\end{prop}

Let us introduce \index{Omegah@$\Omega_h$, semi-classical Euclidean rotation} the Euclidean rotation in the semi-classical setting
\begin{equation}\label{def_Omega_h}
\Omega_h :=x_1 hD_2 - x_2 hD_1 = \oph(x_1\xi_2-x_2\xi_1).
\end{equation}

\begin{prop} \label{Prop : continuity of Op(gamma1):X to L2}
Under the same assumptions as in proposition \ref{Prop : continuity Op(gamma) L2 to L2}, with $\gamma$ replaced by $\gamma_1$, we have that for any $w\in L^2(\mathbb{R}^2)$ such that $\Omega_h w\in L^2_{loc}(\mathbb{R}^2)$  
\begin{equation} \label{est x-Linfty of Op(gamma1)}
\Big\| \oph\Big(\gamma_1\Big(\frac{x|\xi|-\xi}{h^{1/2-\sigma}}\Big)\psi(2^{-k}\xi)a(x)b_p(\xi)\Big)w\Big\|_{L^\infty}
 \lesssim 2^{kp} h^{-\frac{1}{2} -\sigma}\left(\|w\|_{L^2}+ \|\theta_0\Omega_hw\|_{L^2}\right),
\end{equation}
where $\theta_0$ is a smooth function supported in some annulus centred in the origin.
\proof
We prove the statement distinguishing between three spatial regions.
For that, we introduce three cut-off functions: $\theta_s(x)$ supported for $|x|\le m \ll 1$; $\theta_b(x)$ supported for $|x| \ge M'\gg 1$; $\theta(x)$ supported for $m'\le |x| \le M'$, for some $0<m'\ll 1, M\gg 1$,
such that $\theta_s + \theta_b + \theta \equiv 1$. We define respectively $A^k_s(x,\xi):= \gamma_1\Big(\frac{x|\xi|-\xi}{h^{1/2-\sigma}}\Big)\psi(2^{-k}\xi)a(x)b_p(\xi) \theta_s(x)$, $A^k_b(x,\xi):=\gamma_1\Big(\frac{x|\xi|-\xi}{h^{1/2-\sigma}}\Big)\psi(2^{-k}\xi)a(x)b_p(\xi)\theta_b(x)$, and $A^k(x,\xi):=\gamma_1\Big(\frac{x|\xi|-\xi}{h^{1/2-\sigma}}\Big)\psi(2^{-k}\xi)a(x)b_p(\xi)\theta(x)$, so that
$$\gamma_1\Big(\frac{x|\xi|-\xi}{h^{1/2-\sigma}}\Big)\psi(2^{-k}\xi)a(x)b_p(\xi) = A^k_s(x,\xi) + A^k_b(x,\xi) + A^k(x,\xi).$$
The fact that $\oph(A^k_s), \oph(A^k_b)\in\mathcal{L}(L^2;L^\infty)$ and their norm is a $O(2^{kp}h^{-1/2-\sigma})$ follows from lemmas \ref{Lemma on inequalities for Op(A)} and \ref{Lem : est on gamma for wave}. Indeed, when $|x|\ll 1$ (resp. $|x|\gg 1$) we have that $\big|\frac{x|\xi| - \xi}{h^{1/2-\sigma}}\big|\gtrsim h^{-1/2+\sigma}|\xi|$ (resp. $\big|\frac{x|\xi| - \xi}{h^{1/2-\sigma}}\big|\gtrsim h^{-1/2+\sigma}|\xi||x|\gtrsim h^{-1/2+\sigma}|\xi|$), so from lemma \ref{Lem : est on gamma for wave} we derive that
\begin{equation*}
\left| \partial^\alpha_x \partial^\beta_\xi\Big[\gamma_1\Big(\frac{x|\xi| - \xi}{h^{1/2-\sigma}}\Big)\Big]\right|\lesssim \sum_{j=0}^{|\beta|}h^{-(|\alpha|+j)(\frac{1}{2}-\sigma)}\left| \frac{x|\xi|-\xi}{h^{1/2-\sigma}}\right|^{-1-|\alpha|-j}|b_{|\alpha|+j-|\beta|}(\xi)| \lesssim h^{\frac{1}{2}-\sigma}|\xi|^{-1-|\beta|}.
\end{equation*}

Consequently, as $2^{-k}h\le 1$, we deduce that $|\partial^\alpha_x \partial^\beta_\xi \big[A^k_s(\frac{x+y}{2}, h\xi)\big]|\lesssim 2^{kp}h^{-1/2-\sigma}|\xi|^{-1}$ for any $\alpha,\beta\in\mathbb{N}^2$. Therefore
\begin{equation*} 
\left\|\partial_y^{\alpha}\partial_\xi^{\beta}\Big[A^k_s\Big(\frac{x+y}{2},h\xi\Big)\Big]\right\|_{L^2(d\xi)} \lesssim 2^{kp}h^{-\frac{1}{2}-\sigma} \left(\int_{|\xi|\sim 2^kh^{-1}}| \xi|^{-2}d\xi\right)^{\frac{1}{2}}\lesssim 2^{kp}h^{-\frac{1}{2}-\sigma}.
\end{equation*}
The same holds for $A^k_b(x,\xi)$ so, injecting these estimates in inequality \eqref{inequalities |Op(A)|}, we derive that $\|\oph(A^k_s(x,\xi))w\|_{L^\infty} + \|\oph(A^k_b(x,\xi))w\|_{L^\infty} \le C 2^{kp}h^{-\frac{1}{2}-\sigma}\|w\|_{L^2}$. 

A different analysis is needed for $\oph(A^k(x,\xi))w$, since it is no longer true that there exists a positive constant $C$ such that $|x|\xi|-\xi|\ge C |\xi|$ on the support of $A^k(x,\xi)$. 
In this case we exploit the fact that $A^k(x,\xi)$ is supported
in an annulus to perform a change of variables.
If $\theta_0\in C^\infty_0(\mathbb{R}^2\setminus\{0\})$ is a cut-off function equal to 1 on the support of $\theta$ we have that, for any $N\in\mathbb{N}$, $A^k(x,\xi) = \theta_0(x)\sharp A^k(x,\xi) + r_N^k(x,\xi)$ by means of proposition \ref{Prop: a sharp b}, where 
\begin{multline*}
r^k_N(x,\xi) = \left(\frac{h}{2i}\right)^N \frac{N}{(\pi h)^4}   \sum_{|\alpha|= N} \frac{(-1)^{|\alpha|}}{\alpha!}\int e^{\frac{2i}{h}(\eta\cdot z - y\cdot\zeta)} \int_0^1 \partial^\alpha_x \theta_0(x+tz)(1-t)^{N-1} dt\\
\times (\partial^\alpha_\xi A^k)(x, \xi+\eta) \,dy dz d\eta d\zeta.
\end{multline*}
If we take $N$ sufficiently large it turns out that the quantization of $r^k_N$ satisfies a better estimate than \eqref{est x-Linfty of Op(gamma1)}. Indeed, using lemma \ref{Lem : est on gamma for wave} and integrating in $dyd\zeta$, it can be rewritten as
\begin{multline} \label{rkN in prop 1.3}
r^k_N(x,\xi) = \sum_{j\le N} \frac{h^{N-j(\frac{1}{2}-\sigma)}}{(\pi h)^2}\int e^{\frac{2i}{h}\eta\cdot z} \int_0^1 \theta_0(x+tz)(1-t)^{N-1} dt \\
\times \gamma_{1+j}\Big(\frac{x|\xi+\eta|-(\xi+\eta)}{h^{1/2-\sigma}}\Big)\psi(2^{-k}(\xi+\eta))\theta_j(x)a(x)b_{p+j-N}(\xi + \eta)\, dz d\eta,
\end{multline}
for some new functions $\theta_0, \gamma_{1+j}, \psi,\theta_j, a, b_{p+j-N}$. As it is compactly supported in $x$, by lemma \ref{Lemma on inequalities for Op(A)} there is a new cut-off function (that we still call $\theta$) such that
\begin{equation*}
|\oph(r^k_N(x,\xi))w| \lesssim \|w\|_{L^2}\int \Big|\theta\Big(\frac{x+y}{2}\Big)\Big|\sum_{|\alpha'|\le 3}\Big\|\partial^{\alpha'}_y\Big[ r^k_N\Big(\frac{x+y}{2}, h\xi\Big)\Big]\Big\|_{L^2(d\xi)}dy.
\end{equation*}
One can check that the action of $\partial^{\alpha'}_y$ on $r^k_N(\frac{x+y}{2}, h\xi)$ makes appear factors $(h^{-1/2+\sigma}h|\xi+\eta|)^i$, for $i\le |\alpha'|$, without changing the underlining structure of $r^k_N$, and these are bounded by $(h^{-1/2+\sigma}2^k)^i$ on the support of $\psi(2^{-k}h(\xi+\eta))$. After a change of variables $\eta \mapsto h\eta$ in \eqref{rkN in prop 1.3}, we use that $e^{2i \eta\cdot z}=\left(\frac{1-2i\eta\cdot\partial_z}{1+4|\eta|^2}\right)^3\left(\frac{1-2i z\cdot \partial_\eta}{1+4|z|^2}\right)^3e^{2i \eta\cdot z}$, integrate by parts, apply Young's inequality for convolutions, and fix $N>7$, in order to deduce the following chain of inequalities:
\begin{multline*}
 \left\|\partial^{\alpha'}_y r^k_N\Big(\frac{x+y}{2}, h\xi\Big)\right\|_{L^2(d\xi)}^2\\
 \lesssim \sum_{i\le |\alpha'|, j\le N} h^{2N-2j(\frac{1}{2}-\sigma)}\big( h^{-\frac{1}{2}+\sigma}2^k\big)^{2i} 2^{2k(p+j-N)}\int d\xi \left|\int \langle z \rangle^{-3} \langle \eta\rangle^{-3} |\psi(2^{-k}h(\xi+\eta))|dzd\eta\right|^2 \\
  \lesssim  \sum_{i\le |\alpha'|, j\le N} h^{2N-2j(\frac{1}{2}-\sigma)}\big( h^{-\frac{1}{2}+\sigma}2^k\big)^{2i} 2^{2k(p+j-N)} \int |\psi(2^{-k}h\xi)|^2 d\xi \\
\lesssim \sum_{i\le |\alpha'|, j\le N}   h^{2N-2j(\frac{1}{2}-\sigma)}\big( h^{-\frac{1}{2}+\sigma}2^k\big)^{2i} 2^{2k(p+j-N)}\big( h^{-1}2^k\big)^2 \lesssim 2^{2kp},
\end{multline*}
and that $\|\oph(r^k_N)\|_{\mathcal{L}(L^2;L^\infty)}\lesssim 2^{kp}$.
We can then focus on the analysis of the $L^\infty$ norm of $\theta_0(x)\oph(A^k(x,\xi))w$. In polar coordinates $x=\rho e^{i\alpha}$ operator $\Omega_h$ reads as $D_\alpha$, so using the classical one-dimensional Sobolev injection with respect to variable $\alpha$, the one-dimensional semi-classical Sobolev injection with respect to variable $\rho$, and successively returning back to coordinates $x$, we deduce that
\begin{align*}
\big| \theta_0(x)\oph(A^k(x,\xi))w \big|  & \lesssim h^{-\frac{1}{2}} \Big[\|\oph(A^k)w\|_{L^2(dx)} + \|\oph(\xi)\oph(A^k)w\|_{L^2(dx)} \\
& \hspace{1cm}+ \|\Omega_h \theta_0 \oph(A^k)w\|_{L^2(dx)} + \|\oph(\xi)\Omega_h \theta_0 \oph(A^k)w\|_{L^2(dx)}\Big] \\
& \lesssim 2^{kp}h^{-\frac{1}{2}-\sigma}[\|w\|_{L^2} + \|\theta_0\Omega_h w\|_{L^2}]\,.
\end{align*}
The latter of above inequalities is derived observing that the commutator between $\Omega_h$ and $\oph(A^k)$ is a semi-classical pseudo-differential operator whose symbol is linear combination of terms of the form
\begin{equation*}
\gamma_1\Big(\frac{x|\xi|-\xi}{h^{1/2-\sigma}}\Big)\psi(2^{-k}\xi)a(x)\theta(x)b_p(\xi),
\end{equation*}
for some new $\gamma_1, \psi, a, \theta, b_p$, and from the fact that operators $\oph(A^k(x,\xi))$, $\oph(\xi A^k(x,\xi))$ are bounded on $L^2$ (see proposition \ref{Prop : continuity Op(gamma) L2 to L2}), with norm $O(2^{kp})$, $O(2^{k(p+1)})$ respectively, and that $2^k\le h^{-\sigma}$.
\endproof
\end{prop}

\begin{prop} \label{Prop : Linfty est of integral remainders}
Under the same hypothesis as proposition \ref{Prop: L2 est of integral remainders}, $\oph(I^k_{p,q}(x,\xi))$ and $\oph(J^k_{p,q}(x,\xi))$ are bounded operators from $L^\infty$ to $L^2$, with
\begin{equation} \label{norm_Linfty-L2_OP(Ikpq)}
\left\|\oph(I^k_{p,q}(x,\xi))\right\|_{\mathcal{L}(L^2;L^\infty)} + \left\|\oph(J^k_{p,q}(x,\xi))\right\|_{\mathcal{L}(L^2;L^\infty)} \lesssim \sum_{i\le 6} 2^{k(p+q)} (h^{-\frac{1}{2}+\sigma}2^k)^{i} (h^{-1}2^k).
\end{equation}
The same result holds if $q=0$ and $\widetilde{\psi}(2^{-k}\xi)b_q(\xi)\equiv 1$.
\proof
As in proposition \ref{Prop: L2 est of integral remainders}, we prove the statement only for $\oph(I^k_{p,q})$, leaving to the reader to check that the result is true also for $\oph(J^k_{p,q})$.

Let $w\in L^2$. After lemma \ref{Lemma on inequalities for Op(A)} we should prove that $\left\|\partial^\alpha_y \partial^\beta_\xi \Big[I^k_{p,q}(\frac{x+y}{2}, h\xi)\Big]\right\|_{L^2(d\xi)}$ is estimated by the right hand side of \eqref{norm_Linfty-L2_OP(Ikpq)}, for any $|\alpha|, |\beta|\le 3$. A change of variables $\eta \mapsto h\eta$, $\zeta \mapsto h\zeta$ allows us to write $I^k_{p,q}(\frac{x+y}{2}, h\xi)$ as
\begin{multline*}
\frac{1}{\pi^4} \int e^{2i (\eta\cdot z - y'\cdot\zeta)} \left[ \int_0^1  \Big(\gamma\big(h^{\frac{1}{2}+\sigma}(x|\xi|-\xi)\big)\psi(2^{-k}h\xi)a(x)b_p( h\xi)\Big)|_{(\frac{x+y}{2}+tz, \xi+t\zeta)} f(t)dt \right. \\ 
\left. \times \widetilde{\psi}(2^{-k}h(\xi + \eta))a'\Big(\frac{x+y}{2}+y'\Big)b_q( h(\xi+\eta))\right] dy' dz d\eta d\zeta.
\end{multline*}
We observe that, while on the one hand the action of $\partial^\alpha_y$ on the above integral makes appear a factor $(h^{-\frac{1}{2}+\sigma} |h(\xi+t\zeta)|)^i$, with $i\le |\alpha|$, on the other hand that of $\partial^\beta_\xi$ has basically no effect on the $L^2$ norm that we want to estimate as one can check using lemma \ref{Lem : est on gamma for wave} and the fact that $2^{-k}h\le 1$. With this in mind, we can reduce to the study of the $L^2(d\xi)$ norm of an integral function as\small
\begin{multline*}
\sum_{i \le 3} (h^{-\frac{1}{2}+\sigma}2^k)^i \int e^{2i (\eta\cdot z - y'\cdot\zeta)} \left[ \int_0^1  \Big(\gamma\Big(h^{\frac{1}{2}+\sigma}(x|\xi|-\xi)\Big)\psi(2^{-k}h\xi)a(x)b_p(h\xi)\Big)|_{(\frac{x+y}{2}+tz, \xi +t\zeta)}f(t) dt\right. \\
\times  \widetilde{\psi}(2^{-k}h(\xi + \eta))a'\Big(\frac{x+y}{2}+y'\Big)b_q(h(\xi + \eta))\Big] dy' dz d\eta d\zeta,
\end{multline*}\normalsize
for some new functions $\gamma, \psi, a, b_p, \widetilde{\psi}, a', b_q$, with the same properties as their previous homonyms. We use that
$$e^{2i (\eta z - y'\cdot\zeta)}=\left(\frac{1+2iy'\cdot\partial_\zeta}{1+4|y'|^2}\right)^3\left(\frac{1-2i\eta\cdot\partial_z}{1+4|\eta|^2}\right)^3\left(\frac{1-2i z\cdot \partial_\eta}{1+4|z|^2}\right)^3\left(\frac{1+2i\zeta\cdot\partial_{y'}}{1+4|\zeta|^2}\right)^3e^{2i (\eta\cdot z - y'\cdot\zeta)}$$ 
and make some integration by parts to obtain the integrability in $dy'dz d\eta d\zeta$, up to new factors $(h^{-\frac{1}{2}+\sigma} |h(\xi+t\zeta)|)^{j}$, with $j\le 3$, coming out from the derivation of the integrand with respect to $z$. Then, using that functions $h^jb_{p-j}(h(\xi+t\zeta))$ (resp. $h^j b_{q-j}(h(\xi+\eta))$), $j\le 3$, appearing from the derivation of $b_p(h(\xi+t\zeta))$ with respect to $\zeta$ (resp. the derivation of $b_q(h(\xi+\eta))$ with respect to $\eta$), are such that $|h^jb_{p-j}(h(\xi+t\zeta))|\le h^j2^{k(p-j)}\lesssim 2^{kp}$ on the support of $\psi(2^{-k}h(\xi+t\zeta))$ (resp. $|h^jb_{q-j}(h(\xi +\eta))|\le 2^{kq}$ on the support of $\widetilde{\psi}(2^{-k}h(\xi+\eta))$), and the fact that 
\[\Big\| \int \langle\eta\rangle^{-3} |\widetilde{\psi}(2^{-k}h(\xi + \eta))|d\eta\Big\|_{L^2(d\xi)}\le \|\widetilde{\psi}(2^{-k}h\cdot)\|_{L^2} \lesssim h^{-1}2^k,\] we obtain the result of the statement.

The last part of the statement can be proved following an analogous scheme, after having previously made an integration in $dzd\eta$ (or in $dyd\zeta$ if dealing with $J^k_{p,0}$).
\endproof
\end{prop}

\begin{lem} \label{Lem : remainder r^k_N} 
Let $\sigma>0$ be sufficiently small, $k\in K$ with $K$ given by \eqref{set_frequencies_K} and $p,q\in\mathbb{N}$.
Let also $\psi, \widetilde{\psi}\in C^\infty_0(\mathbb{R}^2\setminus\{0\})$, $a(x)$ be either a smooth compactly supported function or $a\equiv 1$, and $f\in C(\mathbb{R})$.
For a fixed integer $N> 2(p+q)+9$ we define
\small
\begin{multline} \label{integral rkN}
r_{N,p}^k(x,\xi) :=  \frac{h^N}{(\pi h)^4} \sum_{|\alpha|=N} \int e^{\frac{2i}{h}(\eta\cdot z - y\cdot \zeta)} \left[\int_0^1 \partial^{\alpha}_x\Big( \gamma_1\Big(\frac{x|\xi|-\xi}{h^{1/2-\sigma}}\Big)\psi(2^{-k}\xi)a(x)b_p(\xi)\Big)|_{(x+tz,\xi+t\zeta)} \right.\\
\times f(t)dt\Big] \partial^{\alpha}_\xi\big(b_q(\xi)\widetilde{\psi}(2^{-k}\xi)\big)|_{(x+y,\xi + \eta)}\, dy dz d\eta d\zeta,
\end{multline}\normalsize
and
\small
\begin{multline} \label{integral rtilde kN}
\widetilde{r}^k_{N,p}(x,\xi) := \frac{h^N}{(\pi h)^4} \sum_{|\alpha_1|+|\alpha_2|=N} \int e^{\frac{2i}{h}(\eta\cdot z - y\cdot \zeta)} \left[\int_0^1 \partial^{\alpha_1}_x\partial^{\alpha_2}_\xi\Big( \gamma_1\Big(\frac{x|\xi|-\xi}{h^{1/2-\sigma}}\Big)\psi(2^{-k}\xi)a(x)b_p(\xi)\Big)|_{(x+tz,\xi+t\zeta)}\right. \\
\times f(t)dt \Big] \partial^{\alpha_2}_x \partial^{\alpha_1}_\xi\big(x_n b_q(\xi)\widetilde{\psi}(2^{-k}\xi)\big)|_{(x+y,\xi + \eta)}\, dy dz d\eta d\zeta\,.
\end{multline}
\normalsize
Then
\begin{equation} \label{est L2 and Linfty rkN r'kN}
\|\oph(r^k_{N,p})\|_{\mathcal{L}(L^2)} +\|\oph(\widetilde{r}^k_{N,p})\|_{\mathcal{L}(L^2)} + \|\oph(r^k_{N,p})\|_{\mathcal{L}(L^2;L^\infty)}+ \|\oph(\widetilde{r}^k_{N,p})\|_{\mathcal{L}(L^2;L^\infty)} \lesssim h^{p+q}.
\end{equation}
\proof
We remind definition \eqref{integral Ik} of integral $I^k_{p,q}(x,\xi)$ for general $k\in K, p,q\in\mathbb{Z}$.
After an explicit development of the derivatives appearing in \eqref{integral rkN} we find that $r^k_{N,p}(x,\xi)$ may be written as
\begin{equation*}
\sum_{j\le N}h^{N-j(\frac{1}{2}-\sigma)}I^k_{p+j,q-N}(x,\xi)
\end{equation*}
where $\gamma$ is replaced with $\gamma_1$ and $a'\equiv 1$ in $I^k_{p+j,q-N}$. Propositions \ref{Prop: L2 est of integral remainders} and \ref{Prop : Linfty est of integral remainders}, combined with the fact that $h\le 2^k\le h^{-\sigma}$, imply respectively that
\begin{align*}
&\|\oph(r^k_{N,p})\|_{\mathcal{L}(L^2)} \lesssim \sum_{j\le N}h^{N-j(\frac{1}{2}-\sigma)}2^{k(p+j+q-N)}\\
&\lesssim \sum_{\substack{j\le N \\ p+j+q\le N}}  h^{N-j(\frac{1}{2}-\sigma) +p+j+q-N} + \sum_{\substack{j\le N \\ p+j+q> N}}  h^{N-j(\frac{1}{2}-\sigma) -\sigma(p+j+q-N)}  \lesssim h^{p+q}
\end{align*}
and
\begin{multline*}
 \|\oph(r^k_{N,p})\|_{\mathcal{L}(L^2;L^\infty)}\lesssim \sum_{i\le 6, j\le N}h^{N-j(\frac{1}{2}-\sigma)}2^{k(p+j+q-N)}(h^{-\frac{1}{2}+\sigma}2^k)^i(h^{-1}2^k)\\
 \lesssim \sum_{\substack{i\le 6, j\le N \\ p+i+j+q\le N-1}} h^{N-1- (i+j)(\frac{1}{2}-\sigma)+p+i+j+q-N+1} +  \sum_{\substack{i\le 6, j\le N \\ p+i+j+q> N-1}} h^{N-1 - (i+j)(\frac{1}{2}-\sigma)-\sigma(p+i+j+q-N+1)} \\
\lesssim h^{p+q},
\end{multline*}
as $N>2(p+q)+9$.

As regards \eqref{integral rtilde kN}, we first observe that index $\alpha_2$ is such that $|\alpha_2|\le 1$ since $x_nb_q(\xi)\widetilde{\psi}(2^{-k}\xi)$ is linear in $x_n$. An explicit development of derivatives in \eqref{integral rtilde kN}, combined with lemma \ref{Lem : est on gamma for wave}, shows that $\widetilde{r}^k_{N,p}(x,\xi)$ splits into two contributions: \small
\begin{align*}
J_0(x,\xi) = \frac{h^N}{(\pi h)^4}\sum_{i\le N} h^{-i(\frac{1}{2}-\sigma)}
\int e^{\frac{2i}{h}(\eta\cdot z -y\cdot\zeta)}\int_0^1\Big( \gamma_{1+i}\Big(\frac{x|\xi| - \xi}{h^{1/2-\sigma}}\Big)\psi(2^{-k}\xi)a(x)b_{p+i}(\xi)\Big)|_{(x+tz, \xi+t\zeta)} f(t)dt\\
\times (x_n+y_n)b_{q-N}(\xi+\eta)\widetilde{\psi}(2^{-k}(\xi+\eta))\, dydzd\eta d\zeta,
\end{align*} \normalsize
for some new functions $a, \psi, \widetilde{\psi}$ and clear meaning for $\gamma_i, b_{p+i}$, $b_{q-N}$, coming out when $|\alpha_2|=0$;
\begin{multline*}
J_1(x,\xi) = \frac{h^N}{(\pi h)^4}\sum_{i\le N-1, j\le 1} h^{-(i+j)(\frac{1}{2}-\sigma)} \int e^{\frac{2i}{h}(\eta\cdot z -y\cdot\zeta)}\\ 
\times \int_0^1\Big( \gamma_{1+i+j}\Big(\frac{x|\xi| - \xi}{h^{1/2-\sigma}}\Big)\psi(2^{-k}\xi)a(x)b_{p+i+j-1}(\xi)\Big)|_{(x+tz, \xi+t\zeta)} f(t)dt \\
\times b_{q-N+1}(\xi+\eta)\widetilde{\psi}(2^{-k}(\xi+\eta))\, dydzd\eta d\zeta,
\end{multline*}\normalsize
for some new other $a, \psi, \widetilde{\psi}$, corresponding instead to $|\alpha_2|=1$. 
One has that
\begin{equation*}
J_1(x,\xi) = \sum_{i\le N-1, j\le 1}h^{N - (i+j)(\frac{1}{2}-\sigma)}I^k_{p+i+j-1,q-N+1}(x,\xi),
\end{equation*}
with $\gamma$ replaced with $\gamma_1$ and $a'\equiv 1$, so propositions \ref{Prop: L2 est of integral remainders} and \ref{Prop : Linfty est of integral remainders}, along with the fact that $N>2(p+q)+9$, imply
\begin{equation*}
\| \oph(J_1(x,\xi))\|_{\mathcal{L}(L^2)} \lesssim \sum_{i\le N-1, j\le 1} h^{N-(i+j)(\frac{1}{2}-\sigma)}2^{k(p+i+j+q-N)} \lesssim h^{p+q},
\end{equation*}
\begin{equation*}
\| \oph(J_1(x,\xi))\|_{\mathcal{L}(L^2; L^\infty)} \lesssim \sum_{\substack{i\le N-1, j\le 1\\ l\le 6}} h^{N-(i+j)(\frac{1}{2}-\sigma)}2^{k(p+i+j+q-N)}(h^{-\frac{1}{2}+\sigma}2^k)^l (h^{-1}2^k) \lesssim h^{p+q}.
\end{equation*}
In order to derive the same estimates for $J_0(x,\xi)$ we split the sum $x_n+y_n$ and analyse separately the two out-coming integrals, that we denote $J_{0,x}(x,\xi), J_{0,y}(x,\xi)$. In the latter one, we use that $y_ne^{-\frac{2i}{h}y\cdot\zeta} = -\frac{h}{2i}\partial_{\zeta_n} e^{-\frac{2i}{h}y\cdot\zeta}$ and successively integrate by parts in $d\zeta_n$ obtaining, with the help of lemma \ref{Lem : est on gamma for wave}, that
\begin{multline} \label{integral J0y}
J_{0,y}(x,\xi) = \sum_{i\le N, j\le 1}h^{N+1-(i+j)(\frac{1}{2}-\sigma)}\int e^{\frac{2i}{h}(\eta\cdot z - y\cdot\zeta)} \\
\times \int_0^1 \Big( \gamma_{1+i+j}\Big(\frac{x|\xi| - \xi}{h^{1/2-\sigma}}\Big)\psi(2^{-k}\xi)a(x)b_{p+i+j-1}(\xi)\Big)|_{(x+tz, \xi +t\zeta)}f(t) dt \\
\times  b_{q-N}(\xi-\eta)\widetilde{\psi}(2^{-k}(\xi+\eta))\ dydzd\eta d\zeta
\end{multline}
for some new functions $a, \psi, \widetilde{\psi}, f$.
Again by propositions \ref{Prop: L2 est of integral remainders}, \ref{Prop : Linfty est of integral remainders} and the fact that $h\le 2^k\le h^{-\sigma}$, $N>2(p+q)+9$, we deduce that: \small
\begin{subequations}\label{est J0y}
\begin{equation}
\| \oph(J_{0,y}(x,\xi))\|_{\mathcal{L}(L^2)} \lesssim \sum_{i\le N, j\le 1}h^{N+1-(i+j)(\frac{1}{2}-\sigma)}2^{k(p+i+j+q-N-1)} \lesssim h^{p+q},
\end{equation} 
\begin{equation}
\| \oph(J_{0,y}(x,\xi))\|_{\mathcal{L}(L^2; L^\infty)} \lesssim \sum_{\substack{i\le N, j\le 1\\ l\le 6}}h^{N+1-(i+j)(\frac{1}{2}-\sigma)}2^{k(p+i+j+q-N-1)} (h^{-\frac{1}{2}-\sigma}2^k)^l (h^{-1}2^k) \lesssim h^{p+q}.
\end{equation}
\end{subequations}\normalsize
In $J_{0,x}(x,\xi)$ we first integrate in $dyd\zeta$ and then we split the occurring integral into two other contributions, called $J_{0, x+tz}(x,\xi), J_{0,tz}(x,\xi)$, by writing $x_n= (x_n+tz_n) - tz_n$. Similarly to what done above, we use that $z_n e^{\frac{2i}{h}\eta\cdot z} = \frac{h}{2i}\partial_{\eta_n} e^{\frac{2i}{h}\eta\cdot z}$ in $J_{0,tz}$, and successively integrate by parts in $d\eta_n$: as $2^{-k}h\le 1$, we obtain that $J_{0, tz}$ has the same form as \eqref{integral J0y} for some new $b_{q-N},\widetilde{\psi}$, and verifies \eqref{est J0y}.
Finally, using that $x_n+tz_n = h^{\frac{1}{2}-\sigma}\big(\frac{(x_n+tz_n)|\xi| - \xi_n}{h^{1/2-\sigma}}\big)|\xi|^{-1} + \xi_n |\xi|^{-1}$, we derive that
\begin{equation*}
\begin{split}
& J_{0, x+tz}(x,\xi) \\
& = \sum_{i\le N} h^{N-(i-1)(\frac{1}{2}-\sigma)}\int e^{\frac{2i}{h}\eta\cdot z}\int \Big(\gamma_{i}\Big(\frac{x|\xi| -\xi}{h^{1/2-\sigma}}\Big) \psi(2^{-k}\xi) a(x)b_{p+i-1}(\xi)\Big)|_{(x+tz, \xi)} f(t) dt \\
& \hspace{10cm} \times b_{q-N}(\xi+\eta) \widetilde{\psi}(2^{-k}(\xi +\eta)) dz d\eta, \\
& + \sum_{i\le N} h^{N-i(\frac{1}{2}-\sigma)}\int e^{\frac{2i}{h}\eta\cdot z}\int \Big(\gamma_{1+i}\big(\frac{x|\xi| -\xi}{h^{1/2-\sigma}}\Big) \psi(2^{-k}\xi) a(x)b_{p+i}(\xi)\Big)|_{(x+tz, \xi)}f(t) dt \\
& \hspace{10cm}\times b_{q-N}(\xi+\eta) \widetilde{\psi}(2^{-k}(\xi +\eta)) dz d\eta,
\end{split}
\end{equation*}
so by propositions \ref{Prop: L2 est of integral remainders} and \ref{Prop : Linfty est of integral remainders} 
\begin{equation*}
\|\oph(J_{0,x+tz}(x,\xi))\|_{\mathcal{L}(L^2)} \lesssim \sum_{i\le N} h^{N-i(\frac{1}{2}-\sigma)}2^{k(p+i+q -N)}\lesssim  h^{p+q},
\end{equation*}
\begin{equation*}
\|\oph(J_{0,x+tz}(x,\xi))\|_{\mathcal{L}(L^2; L^\infty)} \lesssim \sum_{i\le N, l\le 3} h^{N-i(\frac{1}{2}-\sigma)}2^{k(p+i+q -N)}(h^{-\frac{1}{2}+\sigma}2^k)^l(h^{-1}2^k)\lesssim  h^{p+q}.
\end{equation*}
That concludes the proof as $\widetilde{r}^k_{N,p} = J_{0, x+tz}+J_{0, tz} + J_{0,y} + J_1$.
\endproof
\end{lem}

We introduce the following operator: \index{Mj@$\mathcal{M}_j$, operator} 
\begin{equation} \label{def Mj}
\mathcal{M}_j:=\frac{1}{h}\oph(x_j|\xi| - \xi_j), \quad j=1,2
\end{equation}
and use the notation $\|\mathcal{M}^\gamma w\| = \|\mathcal{M}^{\gamma_1}_1 \mathcal{M}^{\gamma_2}_2 w\|$ for any $\gamma=(\gamma_1,\gamma_2)\in\mathbb{N}^2$. 
We have now all the ingredients to state and prove the following two results.

\begin{lem} \label{Lemma : symbolic product development}
Let $\sigma,k,p,\psi,a$ be as in lemma \ref{Lem : remainder r^k_N} and $\widetilde{a}(x)$ such that
\begin{gather*}
(a\equiv 1) \Rightarrow  (\widetilde{a}\equiv 1),\\
(a \text{ compactly supported }) \Rightarrow [(\widetilde{a} \equiv 1) \text{ or } (\widetilde{a} \text{ compactly supported and } \widetilde{a}a\equiv a )].
\end{gather*}
We have that
\begin{multline} \label{symbolic dev 1}
\oph\Big(\gamma_1\Big(\frac{x|\xi| - \xi}{h^{1/2-\sigma}}\Big)\psi(2^{-k}\xi)a(x)b_p(\xi)(x_n|\xi| - \xi_n)\Big)
\\
= \oph\Big(\gamma_1\Big(\frac{x|\xi| - \xi}{h^{1/2-\sigma}}\Big)\psi(2^{-k}\xi)a(x)b_p(\xi)\Big)\widetilde{a}(x) h\mathcal{M}_n +  \oph( r^k_p(x,\xi)),
\end{multline} 
where
\begin{subequations} \label{est Op(rkp)}
\begin{equation} \label{est L2 Op(rkp)}
\big\| \oph(r^k_p (x,\xi))w\big\|_{L^2} \lesssim h^{1-\beta}\|w\|_{L^2} ,
\end{equation}
\begin{equation} \label{est Linfty Op(rkp)}
\big\| \oph(r^k_p(x,\xi))w\big\|_{L^\infty} \lesssim h^{\frac{1}{2}-\beta}(\|w\|_{L^2}+\|\theta_0 \Omega_h w\|_{L^2}),
\end{equation}
\end{subequations}
for some $\theta\in C^\infty_0(\mathbb{R}^2\setminus \{0\})$ and a small $\beta>0$, $\beta\rightarrow 0$ as $\sigma\rightarrow 0$.
Moreover
\begin{subequations} \label{est: L2 Linfty with L}
\begin{equation} \label{est: L2 of gamma1 with L}
\Big\|\oph\Big(\gamma_1\Big(\frac{x|\xi|-\xi}{h^{1/2-\sigma}}\Big)\psi(2^{-k}\xi)a(x)b_p(\xi)(x_n|\xi| - \xi_n)\Big)w\Big\|_{L^2} 
\lesssim h^{1-\beta} \big(\|w\|_{L^2} + \|\mathcal{M}_nw\|_{L^2} \big),
\end{equation}
\begin{multline} \label{est: Linfty of gamma1 with L}
\Big\|\oph\Big(\gamma_1\Big(\frac{x|\xi|-\xi}{h^{1/2-\sigma}}\Big)\psi(2^{-k}\xi)a(x)b_p(\xi) (x_n|\xi| - \xi_n) \Big)w\Big\|_{L^\infty}
\\
\lesssim h^{\frac{1}{2}-\beta}\sum_{\mu=0}^1\Big(\|(\theta_0\Omega_h)^\mu w\|_{L^2}+ \|  (\theta_0\Omega_h)^\mu \mathcal{M}_nw\|_{L^2}\Big).
\end{multline}
\end{subequations}
\proof
The proof of the statement is basically made of tedious calculations and the application of propositions \ref{Prop : continuity Op(gamma) L2 to L2}, \ref{Prop : continuity of Op(gamma1):X to L2} along with lemma \ref{Lem : remainder r^k_N}. 

Let $\widetilde{\psi}\in C^\infty_0(\mathbb{R}^2\setminus\{0\})$ such that $\widetilde{\psi}\equiv 1$ on the support of $\psi$.
From formulas \eqref{a sharp b asymptotic formula}, \eqref{r_N 1} and the hypothesis of the statement we derive that for a fixed $N\in\mathbb{N}$, and up to negligible multiplicative constants, 
\begin{equation} \label{dev of gamma1 with his argument}
\begin{split}
&\left[\gamma_1\Big(\frac{x|\xi|-\xi}{h^{1/2-\sigma}}\Big) \psi(2^{-k}\xi) a(x)b_p(\xi) \right] \sharp \left[(x_n|\xi| - \xi_n) \widetilde{a}(x)\widetilde{\psi}(2^{-k}\xi)\right]\\
&= \gamma_1\Big(\frac{x|\xi|-\xi}{h^{1/2-\sigma}}\Big) \psi(2^{-k}\xi) a(x)b_p(\xi) (x_n|\xi| - \xi_n)  \\
&+ \, h \left\{ \gamma_1\Big(\frac{x|\xi|-\xi}{h^{1/2-\sigma}}\Big)\psi(2^{-k}\xi)a(x)b_p(\xi) , (x_n|\xi| - \xi_n)\right\} \\
&+  \sum_{\substack{2\le |\alpha|<N\\ |\alpha_1| + |\alpha_2| = |\alpha|}} h^{|\alpha|}\partial^{\alpha_1}_x\partial^{\alpha_2}_\xi\Big[ \gamma_1\Big(\frac{x|\xi|-\xi}{h^{1/2-\sigma}}\Big)\psi(2^{-k}\xi)a(x)b_p(\xi)\Big]\partial^{\alpha_2}_x \partial^{\alpha_1}_\xi \big[(x_n|\xi|-\xi_n)\big]
+   r_{N,p}^k(x,\xi),
\end{split}
\end{equation} 
with \small
\begin{multline}\label{rkNp}
r_{N,p}^k(x,\xi) = \frac{h^N}{(\pi h)^4} \sum_{|\alpha_1|+|\alpha_2|=N} \int e^{\frac{2i}{h}(\eta\cdot z - y\cdot \zeta)} \left[\int_0^1 \partial^{\alpha_1}_x\partial^{\alpha_2}_\xi\Big[ \gamma_1\Big(\frac{x|\xi|-\xi}{h^{1/2-\sigma}}\Big)\psi(2^{-k}\xi)a(x)b_p(\xi)\Big]|_{(x+tz,\xi+t\zeta)}\right.\\
\times (1-t)^{N-1}dt\Big] 
\partial^{\alpha_2}_x\partial^{\alpha_1}_\xi\big[(x_n|\xi|-\xi_n)\widetilde{a}(x)\widetilde{\psi}(2^{-k}\xi)\big]|_{(x+y,\xi + \eta)}\, dy dz d\eta d\zeta\,.
\end{multline}\normalsize
If $\widetilde{a}\equiv 1$ above $r^k_{N,p}$ can be decomposed into the sum of integrals of the form \eqref{integral rkN} and \eqref{integral rtilde kN} with $q=1$, so
\begin{equation} \label{norms_Op_rkNp}
\left\|\oph(r^k_{N,p})\right\|_{\mathcal{L}(L^2)}+ \left\|\oph(r^k_{N,p})\right\|_{\mathcal{L}(L^2;L^\infty)} \lesssim h^{1+p}
\end{equation}
if $N$ is taken sufficiently large (e.g. $N>2p+11$).
The same is true if functions $a,\widetilde{a}$ are compactly supported as follows by propositions \ref{Prop: L2 est of integral remainders} and \ref{Prop : Linfty est of integral remainders}, since from lemma \ref{Lem : est on gamma for wave} and definition \eqref{integral Ik} of $I^k_{p,q}$ for general $k\in K, p,q\in \mathbb{Z}$
\begin{equation*}
r^k_{N,p}(x,\xi) = \sum_{\substack{|\alpha_1|+|\alpha_2| = N \\ i\le |\alpha_1|, 1\le j\le |\alpha_2|}} h^{N-(i+j)(\frac{1}{2}-\sigma)}I^k_{p+i+j -|\alpha_2|, 1-|\alpha_1|}(x,\xi).
\end{equation*}

An explicit computation of the Poisson bracket in \eqref{dev of gamma1 with his argument} shows that it is equal to
\begin{multline} \label{first order term symb dev}
h (\partial \gamma_1)\Big(\frac{x|\xi|-\xi}{h^{1/2-\sigma}}\Big)\Big(\frac{x_1\xi_2 - x_2\xi_1}{h^{1/2-\sigma}}\Big)\psi(2^{-k}\xi)a(x)b_p(\xi) \\
 +{\sum}' h \gamma_1\Big(\frac{x|\xi|-\xi}{h^{1/2-\sigma}}\Big)\psi(2^{-k}\xi)a(x)b_p(\xi),
\end{multline}
where $\sum'$ is a concise notation to indicate a linear combination, and $\psi, a, b_p$ are some new functions with the same features of their homonyms.
After writing
\begin{equation}\label{x1 xi2 - x2 xi1}
(x_1\xi_2 - x_2\xi_1) = (x_1|\xi| - \xi_1)\xi_2|\xi|^{-1} - (x_2|\xi| - \xi_2)\xi_1|\xi|^{-1},
\end{equation}
we recognize that the quantization of \eqref{first order term symb dev} verifies estimates \eqref{est Op(rkp)} thanks to propositions \ref{Prop : continuity Op(gamma) L2 to L2}, \ref{Prop : continuity of Op(gamma1):X to L2} and the fact that $2^{kp}\le h^{-\sigma p}$.

Let us denote concisely by $t^k_\alpha$ the $|\alpha|$-order contributions in \eqref{dev of gamma1 with his argument}, for $2\le |\alpha|<N$.
As factor $x_n|\xi| -\xi_n$ is affine in $x_n$, the length of multi-index $\alpha_2$ is less or equal than 1 and, using lemma \ref{Lem : est on gamma for wave}, $t^k_\alpha$ appears to be the sum of two terms. The first one corresponds to $|\alpha_2| = 0$ and has the form
\begin{equation*}
{\sum_{\substack{i\le |\alpha|\\ \mu=0,1}}}'h^{|\alpha| - i(\frac{1}{2}-\sigma)}  \gamma_{1+i}\Big(\frac{x|\xi|-\xi}{h^{1/2-\sigma}}\Big)\psi(2^{-k}\xi)a(x)b_{p+i+1 - |\alpha| }(\xi)\, x^\mu_n , 
\end{equation*}
for some new functions $\psi, a$. Observe that $\mu=0$ if $a\equiv 1$ because the derivation of $\gamma_1\big(\frac{x|\xi|-\xi}{h^{1/2-\sigma}}\big)$ $|\alpha_1|$-times with respect to $x$ makes appear, inter alia, a factor $|\xi|^{|\alpha_1|}$ that allows us to rewrite $\partial^{\alpha_1}_\xi(x_n|\xi|-\xi_n)$ from $(x_n|\xi|-\xi_n)+b_0(\xi)$, for some new $b_0$, and $\partial^{\alpha_1} _z\gamma_1(z)z_n$ is of the form $\gamma_{|\alpha_1|}(z)$).
The second term, corresponding instead to $|\alpha_2|=1$, is given by
\begin{equation*} 
{\sum_{i\le |\alpha|-1, j\le 1}}'h^{|\alpha| - (i+j)(\frac{1}{2}-\sigma)}  \gamma_{1+i+j}\Big(\frac{x|\xi|-\xi}{h^{1/2-\sigma}}\Big)\psi(2^{-k}\xi)a(x)b_{p+i+j+1 - |\alpha| }(\xi),
\end{equation*}
for some new other functions $\psi, a$. From propositions \ref{Prop : continuity Op(gamma) L2 to L2}, \ref{Prop : continuity of Op(gamma1):X to L2} we then deduce that
\begin{subequations}
\begin{gather}
\|\oph(t^k_\alpha)w\|_{L^2} \lesssim (h^{\frac{|\alpha|}{2}-\beta} + h^{1+p})\|w\|_{L^2},\label{tk_alpha} \\
\|\oph(t^k_\alpha)w\|_{L^\infty} \lesssim (h^{\frac{|\alpha|-1}{2}-\beta} + h^{\frac{1}{2}+p})(\|w\|_{L^2}+\| \theta \Omega_h w\|_{L^2}),
\end{gather}
\end{subequations}
which concludes that
\begin{multline*}
\left[\gamma_1\Big(\frac{x|\xi|-\xi}{h^{1/2-\sigma}}\Big) \psi(2^{-k}\xi) a(x)b_p(\xi) \right] \sharp \left[(x_n|\xi| - \xi_n) \widetilde{a}(x)\widetilde{\psi}(2^{-k}\xi)\right] \\ 
=  \gamma_1\Big(\frac{x|\xi|-\xi}{h^{1/2-\sigma}}\Big) \psi(2^{-k}\xi) a(x)b_p(\xi) (x_n|\xi| - \xi_n) + r^k_p(x,\xi),
\end{multline*}
with $r^k_p$ satisfying \eqref{est Op(rkp)}.

Finally, by symbolic calculus we have that, up to some multiplicative constants,
\begin{align*}
\oph\big((x_n|\xi| - \xi_n)\widetilde{a}(x)\widetilde{\psi}(2^{-k}\xi)\big) &= \widetilde{a}(x) \oph\big((x_n|\xi| - \xi_n)\widetilde{\psi}(2^{-k}\xi)\big) + \oph(r^k(x,\xi)) \\
&=\oph(\widetilde{\psi}(2^{-k}\xi))\widetilde{a}(x) h\mathcal{M}_n + h\widetilde{a}(x) \oph((\partial\widetilde{\psi})(2^{-k}\xi)(2^{-k}|\xi|)) \\
&+ \oph(\widetilde{r}^k(x,\xi)) h\mathcal{M}_n
+ \oph(r^k(x,\xi)),
\end{align*}
where
\begin{align*}
r^k(x,\xi)&=\frac{h}{(2\pi)^2}\int e^{\frac{2i}{h}\eta\cdot z}\int \partial_x\widetilde{a}(x+tz) dt\ \partial_\xi\big[(x_n|\xi| -\xi_n)\widetilde{\psi}(2^{-k}\xi)\big]|_{(x,\xi+\eta)} dzd\eta, \\
\widetilde{r}^k(x,\xi)&=\frac{h2^{-k}}{(2\pi)^2}\int e^{\frac{2i}{h}\eta\cdot z}\int \partial_x\widetilde{a}(x+tz) dt\ (\partial_\xi\widetilde{\psi})(2^{-k}(\xi+\eta)) dzd\eta,
\end{align*}
are such that $\|\oph(r^k_1)\|_{\mathcal{L}(L^2)}=O(h)$, $\|\oph(\widetilde{r}^k_1)\|_{\mathcal{L}(L^2)}=O(1)$.
An explicit computation shows also that $\|[\Omega_h, \oph(r^k)]\|_{\mathcal{L}(L^2)}=O(h)$ and $\|[\Omega_h, \oph(\widetilde{r}^k)]\|_{\mathcal{L}(L^2)}=O(1)$.
Therefore, since $\widetilde{\psi}\equiv 1$ on the support of $\psi$, $\widetilde{a}\equiv 1$ on the support of $a$, one can use remark \ref{Remark:symbols_with_null_support_intersection} together with propositions \ref{Prop: L2 est of integral remainders}, \ref{Prop : Linfty est of integral remainders}, and also propositions \ref{Prop : continuity Op(gamma) L2 to L2}, \ref{Prop : continuity of Op(gamma1):X to L2}, to show that
\begin{multline*}
\oph\Big(\gamma_1\Big(\frac{x|\xi|-\xi}{h^{1/2-\sigma}}\Big) \psi(2^{-k}\xi) a(x)b_p(\xi)\Big)\oph\big((x_n|\xi| - \xi_n) \widetilde{a}(x)\widetilde{\psi}(2^{-k}\xi)\big)\\
 = \oph\Big(\gamma_1\Big(\frac{x|\xi|-\xi}{h^{1/2-\sigma}}\Big) \psi(2^{-k}\xi) a(x)b_p(\xi)\Big)\widetilde{a}(x)h\mathcal{M}_n +\oph(r^k_p(x,\xi)),
\end{multline*}
for a new $\oph(r^k_p(x,\xi))$ satisfying \eqref{est L2 Op(rkp)}. This concludes the proof of \eqref{symbolic dev 1} and of the entire statement by applying propositions \ref{Prop : continuity Op(gamma) L2 to L2}, \ref{Prop : continuity of Op(gamma1):X to L2} to the first operator in the above right hand side.
\endproof
\end{lem}

\begin{lem}\label{Lem: Gamma with double argument-wave}
Let $\sigma>0$ be small, $k\in K$ with $K$ given by \eqref{set_frequencies_K} and $p\in \mathbb{N}$. Let also $\gamma\in C^\infty_0(\mathbb{R}^2)$ be equal to 1 in a neighbourhood of the origin, $\psi\in C^\infty_0(\mathbb{R}^2\setminus\{0\})$, and $a\in C^\infty_0(\mathbb{R}^2)$. 
For any function $w\in L^2(\mathbb{R}^2)$ such that $\mathcal{M}w\in L^2(\mathbb{R}^2)$, any $m,n=1,2$, we have that
\begin{multline*}
\oph\Big(\gamma\Big(\frac{x|\xi|-\xi}{h^{1/2-\sigma}}\Big)\psi(2^{-k}\xi)a(x)b_p(\xi)(x_m|\xi|-\xi_m)(x_n|\xi|-\xi_n) \Big)w\\= \oph\Big(\gamma\Big(\frac{x|\xi|-\xi}{h^{1/2-\sigma}}\Big)\psi(2^{-k}\xi)a(x)b_p(\xi)(x_m|\xi|-\xi_m) \Big)[h\mathcal{M}_n w] +O_{L^2}\big(h^{2-\beta}(\|w\|_{L^2}+\|\mathcal{M}w\|_{L^2})\big),
\end{multline*} 
with $\beta>0$ small, $\beta\rightarrow 0$ as $\sigma\rightarrow0$.
\proof
Let $\widetilde{\gamma}(z):=\gamma(z)z_m$ and $\widetilde{\psi}\in C^\infty_0(\mathbb{R}^2\setminus\{0\})$ be identically equal to 1 on the support of $\psi$.
We saw in the proof of the previous lemma that the symbolic product
\begin{equation*}
\left[\widetilde{\gamma}\Big(\frac{x|\xi|-\xi}{h^{1/2-\sigma}}\Big)\psi(2^{-k}\xi)a(x)b_p(\xi)\right] \sharp [(x_n|\xi|-\xi_n)\widetilde{\psi}(2^{-k}\xi)]
\end{equation*}
develops as in \eqref{dev of gamma1 with his argument}, \eqref{rkNp}, with $\gamma_1$ replaced with $\widetilde{\gamma}$ and $\widetilde{a}\equiv 1$.
From \eqref{first order term symb dev}, the fact that
\begin{equation*}
\{x_m|\xi|-\xi_m, x_n|\xi|-\xi_n\}=
\begin{cases}
0 \quad &\text{if } m=n,\\
(-1)^{m+1} (x_1\xi_2-\xi_2x_1) &\text{if } m\ne n,
\end{cases}
\end{equation*}
and that $(x_1\xi_2-\xi_2x_1) = (x_1|\xi|-\xi_1)\xi_2|\xi|^{-1}-(x_2|\xi|-\xi_2)\xi_1|\xi|^{-1}$, we derive that the first order term of the mentioned symbolic development is a linear combination of products of the form
\begin{equation*}
h^{\frac{3}{2}}\gamma \Big(\frac{x|\xi|-\xi}{h^{1/2-\sigma}}\Big)\psi(2^{-k}\xi)a(x)b_p(\xi)(x_j|\xi|-\xi_j),
\end{equation*}
for some new functions $\gamma, \psi, a$, and its quantization acting on $w$ is a remainder as in the statement after estimate \eqref{est: L2 of gamma1 with L}.

The second order term is given, up to some negligible multiplicative constants, by
\begin{multline*}
h^{1+2\sigma}\sum_{|\alpha|=2}(\partial^\alpha\gamma)\Big(\frac{x|\xi|-\xi}{h^{1/2-\sigma}}\Big)\psi(2^{-k}\xi)a_1(x)b_{p+1}(\xi)(x_m|\xi|-\xi_m)\\
+ h^{\frac{3}{2}+\sigma}\sum_{|\alpha|=1}(\partial^\alpha\gamma)\Big(\frac{x|\xi|-\xi}{h^{1/2-\sigma}}\Big)\psi_2(2^{-k}\xi)a_2(x)b_{p+1}(\xi)\\
+ h^2{\sum}' \gamma\Big(\frac{x|\xi|-\xi}{h^{1/2-\sigma}}\Big)\psi_3(2^{-k}\xi)a_3(x)b_{p+1}(\xi),
\end{multline*}
for some new smooth compactly supported $\psi_2,\psi_3, a_1, a_2, a_3$, and as the derivatives of $\gamma$ vanish in a neighbourhood of the origin we can replace $(\partial^\alpha \gamma)(z)$ with $\sum_j \gamma_1^j(z)z_j$, $\gamma_j^1(z):=(\partial^\alpha\gamma)(z)z_j|z|^{-2}$, when $|\alpha|=1$.
The third order one is given by
\begin{multline*}
h^{\frac{3}{2}+3\sigma}\sum_{|\alpha|=3}(\partial^\alpha\gamma)\Big(\frac{x|\xi|-\xi}{h^{1/2-\sigma}}\Big)\psi(2^{-k}\xi)a_1(x)b_{p+1}(\xi)(x_m|\xi|-\xi_m)\\
 + h^2{\sum}' \gamma_1\Big(\frac{x|\xi|-\xi}{h^{1/2-\sigma}}\Big)\psi_1(2^{-k}\xi)a_2(x)b_{p+1}(\xi),
\end{multline*}
for some other $\psi_1, a_1, a_2$ and a new $\gamma_1\in C^\infty_0(\mathbb{R}^2)$.
Using estimate \eqref{est: L2 of gamma1 with L} for the summations in $\alpha$ and proposition \ref{Prop : continuity Op(gamma) L2 to L2} for the remaining terms in the above expressions we obtain that the quantizations of the second and third order term are also a $O_{L^2}\big(h^{2-\beta}(\|w\|_{L^2}+\|\mathcal{M}w\|_{L^2})\big)$ when acting on $w$, for a small $\beta>0$, $\beta\rightarrow 0$ as $\sigma\rightarrow 0$.

In all the other $|\alpha|$-order terms, with $4\le |\alpha|\le N-1$, and in integral remainder $r^k_{N,p}$, we look at $\gamma\big(\frac{x|\xi|-\xi}{h^{1/2-\sigma}}\big)\psi(2^{-k}\xi)a(x)b_p(\xi)(x_m|\xi|-\xi_m)$ as a symbol of the form 
\begin{equation*}
\gamma\Big(\frac{x|\xi|-\xi}{h^{1/2-\sigma}}\Big)\psi(2^{-k}\xi)a(x)b_{p+1}(\xi)
\end{equation*}
for a new $a_1\in C^\infty_0(\mathbb{R}^2)$. From \eqref{tk_alpha} and \eqref{norms_Op_rkNp} when $N>11$, we derive that the quantizations of these terms are also a $O_{L^2}\big(h^{2-\beta}(\|w\|_{L^2}+\|\mathcal{M}w\|_{L^2})\big)$ when acting on $w$.

We finally obtain that
\begin{multline*}
\oph\Big(\gamma\Big(\frac{x|\xi|-\xi}{h^{1/2-\sigma}}\Big)\psi(2^{-k}\xi)a(x)b_p(\xi)(x_m|\xi|-\xi_m)(x_n|\xi|-\xi_n) \Big)w\\= \oph\Big(\gamma\Big(\frac{x|\xi|-\xi}{h^{1/2-\sigma}}\Big)\psi(2^{-k}\xi)a(x)b_p(\xi)(x_m|\xi|-\xi_m) \Big)\oph\big((x_n|\xi|-\xi_n)\widetilde{\psi}(2^{-k}\xi)\big) \\
+O_{L^2}\big(h^{2-\beta}(\|w\|_{L^2}+\|\mathcal{M}w\|_{L^2})\big),
\end{multline*} 
and the conclusion of the proof comes then from the fact that, by symbolic calculus,
\begin{equation*} 
\oph\big((x_n|\xi| - \xi_n)\widetilde{\psi}_1(2^{-k}\xi)\big)= h \oph(\widetilde{\psi}_1(2^{-k}\xi))\mathcal{M}_n - \frac{h}{2i} \oph\big((\partial\widetilde{\psi}_1)(2^{-k}\xi)\cdot(2^{-k}\xi)\big),
\end{equation*}
and by remark \ref{Remark:symbols_with_null_support_intersection} as all derivatives of $\widetilde{\psi}$ vanish on the support of $\psi$.
\endproof
\end{lem}
The following lemma is introduced especially for the proof of lemma \ref{Lem: preliminary on Op(e)}. Even if quite similar to lemma \ref{Lemma : symbolic product development}, we are going to see that the particular structure of symbolic product in the left hand side of \eqref{symboli_dev_enhanced} allows for a remainder $r^k_p$ satisfying enhanced estimate \eqref{est Linfty Op(rkp) enhanced} rather than \eqref{est Linfty Op(rkp)}.

\begin{lem} \label{Lemma : on the enhanced symbolic product}
Let us take $\sigma>0$ sufficiently small, $k\in K$ and $p,q\in\mathbb{N}$.
Let also $\gamma\in C^\infty_0(\mathbb{R}^2)$ such that $\gamma\equiv 1$ in a neighbourhood of the origin, $\psi, \widetilde{\psi}\in C^\infty_0(\mathbb{R}^2\setminus\{0\})$ such that $\psi \equiv 1$ on the support of $\widetilde{\psi}$, $a(x)$ be a smooth compactly supported function.
Then
\begin{multline} \label{symboli_dev_enhanced}
\Big[ (x_n|\xi| - \xi_n)\widetilde{\psi}(2^{-k}\xi)a(x)b_p(\xi)\Big] \sharp \, \Big[\gamma\Big(\frac{x|\xi|-\xi}{h^{1/2-\sigma}}\Big)\psi(2^{-k}\xi)\Big]\\
= \gamma\Big(\frac{x|\xi|-\xi}{h^{1/2-\sigma}}\Big)\widetilde{\psi}(2^{-k}\xi)a(x)b_p(\xi)(x_n|\xi| - \xi_n) + r^k_p(x,\xi),
\end{multline}
where
\begin{subequations} \label{est Op(rkp) enhanced}
\begin{equation} \label{est L2 Op(rkp) enhanced}
\left\|\oph(r^k_p(x,\xi))w \right\|_{L^2}\lesssim h^{\frac{3}{2}-\beta} (\|w\|_{L^2} + \| \mathcal{M}w\|_{L^2}) +h^{1+p}\|w\|_{L^2},
\end{equation}
\begin{equation} \label{est Linfty Op(rkp) enhanced}
\left\|\oph(r^k_p(x,\xi))w \right\|_{L^\infty}
\lesssim h^{1-\beta}\sum_{\mu=0}^1\Big(\|(\theta_0\Omega_h)^\mu w\|_{L^2}+ \|  (\theta_0\Omega_h)^\mu \mathcal{M}w\|_{L^2}\Big),
\end{equation}
\end{subequations}
for some $\theta\in C^\infty_0(\mathbb{R}^2\setminus\{0\})$, and a small $\beta>0$, $\beta\rightarrow 0$ as $\sigma\rightarrow 0$.
\proof
Using proposition \ref{Prop: a sharp b}, for a fixed $N\in\mathbb{N}$ and up to multiplicative constants independent of $h, k,$ we have the following symbolic development:
\begin{equation} \label{symb dev 2}
\begin{split}
&\Big[ (x_n|\xi| - \xi_n) \widetilde{\psi}(2^{-k}\xi)a(x)b_p(\xi)\Big] \sharp \, \Big[\gamma\Big( \frac{x|\xi|-\xi}{h^{1/2-\sigma}}\Big)\psi(2^{-k}\xi)\Big] \\
& = \gamma\Big(\frac{x|\xi|-\xi}{h^{1/2-\sigma}}\Big)\widetilde{\psi}(2^{-k})a(x)b_p(\xi)(x_n|\xi| - \xi_n) \\
& + h \left\{(x_n|\xi| - \xi_n)\widetilde{\psi}(2^{-k}\xi)a(x)b_p(\xi),  \gamma\Big(\frac{x|\xi|-\xi}{h^{1/2-\sigma}}\Big)\right\} \\
& + \sum_{\substack{\alpha=(\alpha_1,\alpha_2)\\ 2\le |\alpha| <N}}h^{|\alpha|}\partial^{\alpha_1}_x \partial^{\alpha_2}_\xi\Big[(x_n|\xi| - \xi_n)\widetilde{\psi}(2^{-k}\xi)a(x)b_p(\xi) \Big]   \partial^{\alpha_2}_x \partial^{\alpha_1}_\xi\Big[\gamma\Big(\frac{x|\xi|-\xi}{h^{1/2-\sigma}}\Big)\Big]+ r^k_{N,p}(x,\xi),
\end{split}
\end{equation}
with \small
\begin{multline*}
r^k_{N,p}(x,\xi) = \frac{h^N}{(\pi h)^4}\sum_{|\alpha_1|+|\alpha_2| =N}\int e^{\frac{2i}{h}(\eta\cdot z - y\cdot \zeta)}\left[\int_0^1 \partial^{\alpha_1}_x \partial^{\alpha_2}_\xi\big[(x_n|\xi| - \xi_n)a(x)b_p(\xi)\widetilde{\psi}(2^{-k}\xi) \big]|_{(x+tz, \xi+t\zeta)}\right.\\
\times (1-t)^{N-1}dt \Big]  \partial^{\alpha_2}_x \partial^{\alpha_1}_\xi\Big[ \gamma\Big(\frac{x|\xi| - \xi}{h^{1/2-\sigma}}\Big)\psi(2^{-k}\xi)\Big]|_{(x+y,\xi+\eta)} \ dydz d\eta d\zeta.
\end{multline*}\normalsize

For sake of simplicity, we denote by $t^k_1$ (resp. $t^k_\alpha$, $|\alpha|=2,\dots, N-1$) the Poisson brackets (resp. the $|\alpha|$-th contribution) in \eqref{symb dev 2}.
An explicit computation of $t^k_1$, combined with the fact that $x_1\xi_2-x_2\xi_1 = (x_1|\xi|-\xi_1)\xi_2|\xi|^{-1} - (x_2|\xi|-\xi_2)\xi_1|\xi|^{-1}$, shows that it is linear combination of terms of the form
$$h(\partial\gamma)\Big(\frac{x|\xi|-\xi}{h^{1/2-\sigma}}\Big)\Big(\frac{x_j|\xi| - \xi_j}{h^{1/2-\sigma}}\Big)\widetilde{\psi}(2^{-k}\xi)a(x)b_p(\xi),$$ 
for $j\in\{1,2\}$ and some new functions $\widetilde{\psi}, a, b_p$, so by inequalities \eqref{est: L2 Linfty with L} we derive that
\begin{subequations} \label{est L2 Linfty Op(tk1)}
\begin{equation}
\left\|\oph(t^k_1)w \right\|_{L^2}\lesssim h^{\frac{3}{2}-\beta} \left(\|w\|_{L^2}+\|\mathcal{M}w\|_{L^2}\right),
\end{equation}
\begin{equation}
\left\|\oph(t^k_1)w \right\|_{L^\infty}\lesssim h^{1-\beta}\sum_{\mu=0}^1(\|(\theta_0\Omega_h)^\mu w\|_{L^2} + \|(\theta_0\Omega_h)^\mu \mathcal{M}w\|_{L^2}).
\end{equation}
\end{subequations}
The improvement of these estimates with respect to \eqref{est Op(rkp)} is attributable to the choice of $\psi$ identically equal to 1 on the support of $\widetilde{\psi}$.
All derivatives of $\psi$ vanish against $\widetilde{\psi}$, so in the development of $t^k_1$ we avoid terms like $\gamma\big(\frac{x|\xi|-\xi|}{h^{1/2-\sigma}}\big)\widetilde{\psi}(2^{-k}\xi)a(x)b_p(\xi)(\partial\psi)(2^{-k}\xi)(2^{-k}|\xi|)$, coming out from $\{x_n|\xi| -\xi_n, \psi(2^{-k}\xi)\}\gamma\big(\frac{x|\xi|-\xi|}{h^{1/2-\sigma}}\big)\widetilde{\psi}(2^{-k}\xi)a(x)b_p(\xi)$, that do not enjoy estimates like \eqref{est L2 Linfty Op(tk1)}.

Using formula \eqref{derivatives of gamma_n 2} and looking at $(x_n|\xi|-\xi_n)\widetilde{\psi}(2^{-k}\xi)a(x)b_p(\xi)$ as a linear combination of terms $\widetilde{\psi}(2^{-k}\xi)a(x)b_{p+1}(\xi)$, for some new $\widetilde{\psi},a,b_{p+1}$, we realize that, for any $2\le |\alpha|<N$, 
\begin{equation*}
t^k_\alpha = \sum_{\substack{|\alpha_1| + |\alpha_2| = |\alpha| \\ 1\le j \le |\alpha_1|}} h^{|\alpha| - (j+|\alpha_2|)(\frac{1}{2}-\sigma)}\gamma_{j+|\alpha_2|}\Big(\frac{x|\xi| - \xi}{h^{1/2-\sigma}}\Big)\widetilde{\psi}(2^{-k}\xi)a_j(x) b_{p+j+1-|\alpha_1|}(\xi),
\end{equation*}
for some new other $\widetilde{\psi}, a_j$, with $a_j$ compactly supported, and then that
\begin{equation*}
\|\oph(t^k_\alpha)w\|_{L^2}\lesssim \sum_{\substack{|\alpha_1| + |\alpha_2| = |\alpha| \\ 1\le j\le |\alpha_1|}} h^{|\alpha| - (j+|\alpha_2|)(\frac{1}{2}-\sigma)}2^{k(p+j+1 - |\alpha_1|)} \|w\|_{L^2}, 
\end{equation*}
\begin{multline*}
\|\oph(t^k_\alpha)w\|_{L^\infty}\\
\lesssim \sum_{\substack{|\alpha_1| + |\alpha_2| = |\alpha| \\ 1\le j\le |\alpha_1|}} h^{|\alpha| - (j+|\alpha_2| )(\frac{1}{2}-\sigma)}2^{k(p+j+1 - |\alpha_1|)}h^{-\frac{1}{2}-\sigma} \sum_{\mu=0}^1(\|(\theta_0\Omega_h)^\mu w\|_{L^2} + \|(\theta_0\Omega_h)^\mu \mathcal{M}w\|_{L^2}),
\end{multline*}
after propositions \ref{Prop : continuity Op(gamma) L2 to L2}, \ref{Prop : continuity of Op(gamma1):X to L2}.
For $|\alpha|\ge 3$, the above estimates imply $\|\oph(t^k_\alpha)\|_{\mathcal{L}(L^2)}\lesssim h^{\frac{3}{2}-\beta}$ and $\|\oph(t^k_\alpha)w\|_{L^\infty}\lesssim h^{1-\beta} \sum_{\mu=0}^1(\|(\theta_0\Omega_h)^\mu w\|_{L^2} + \|(\theta_0\Omega_h)^\mu \mathcal{M}w\|_{L^2})$. 
For $|\alpha| =2$, we exploit the fact that functions $\gamma_{j+|\alpha_2|}$ vanish in a neighbourhood of the origin, as they come from $\gamma$'s derivatives, and define $\gamma^i_{j+|\alpha_2| }(z):= \gamma_{j+|\alpha_2|}(z) z_i|z|^{-2}$, $i=1,2$, so that \small
\begin{equation*}
t^k_\alpha = \sum_{\substack{|\alpha_1| + |\alpha_2| = |\alpha| \\ 1\le j \le |\alpha_1|, i=1,2}} h^{|\alpha| - (j+|\alpha_2|)(\frac{1}{2}-\sigma)}\gamma^i_{j+|\alpha_2|}\Big(\frac{x|\xi| - \xi}{h^{1/2-\sigma}}\Big)\Big(\frac{x_i|\xi| - \xi_i}{h^{1/2-\sigma}}\Big)\widetilde{\psi}(2^{-k}\xi)a_j(x) b_{p+j+1-|\alpha_1|}(\xi),
\end{equation*}\normalsize
to which we can then apply lemma \ref{Lemma : symbolic product development}.
After inequalities \eqref{est: L2 Linfty with L}, $\oph(t^k_\alpha)$ with $|\alpha|=2$ also satisfies \eqref{est L2 Linfty Op(tk1)}.

Finally, reminding definition \eqref{integral Jk} of $J^k_{p,q}(x,\xi)$ for general $k\in K, p,q\in\mathbb{Z}$, and developing derivatives in $r^k_{N,p}$ using lemma \ref{Lem : est on gamma for wave}, we observe that 
\begin{equation*}
r^k_{N,p} = \sum_{\substack{|\alpha_1| + |\alpha_2| = N \\ 0\le j\le |\alpha_1|}} h^{N-(|\alpha_2| +j)(\frac{1}{2}-\sigma)}J^k_{p+1-|\alpha_2|, |\alpha_2|+j-|\alpha_1|}(x,\xi),
\end{equation*}
hence propositions \ref{Prop: L2 est of integral remainders} and \ref{Prop : Linfty est of integral remainders} give that
\begin{equation*}
\|\oph(r^k_{N,p})\|_{\mathcal{L}(L^2)}\lesssim \sum_{\substack{|\alpha_1| + |\alpha_2| = N \\ 0\le j\le |\alpha_1|}} h^{N-(|\alpha_2| +j)(\frac{1}{2}-\sigma)}2^{k(p+1 +j- |\alpha_1|)}\lesssim h^{1+p},
\end{equation*}
\begin{equation*}
\|\oph(r^k_{N,p})\|_{\mathcal{L}(L^2;L^\infty)}\lesssim \sum_{\substack{|\alpha_1| + |\alpha_2| = N\\ 0\le j\le |\alpha_1|, i\le 6}} h^{N-(|\alpha_2| +j)(\frac{1}{2}-\sigma)}2^{k(p+1+ j- |\alpha_1|)}(h^{-\frac{1}{2}+\sigma}2^k)^i(h^{-1}2^k)\lesssim h^{1+p},
\end{equation*}
if $N$ is chosen sufficiently large (e.g. $N>10+2p$). 
We should also highlight the fact that, at the difference of \eqref{est Linfty Op(rkp) enhanced}, \eqref{est L2 Op(rkp) enhanced} does not improve \eqref{est L2 Op(rkp)}: if we get a $h^{\frac{3}{2}-\beta}$ factor in front of the first term in the right hand side, the second term $h^{1+p}\|w\|_{L^2}$ is just a $O(h^{1-\beta})$ in the case $p=0$, coming from $|\alpha_1|=N$, $j=|\alpha_2|=0$, $p=0$ above.
\endproof
\end{lem}

\subsection{Operators for the Klein-Gordon solution: some estimates} \label{Subsection: Some Technical Estimates II}
This subsection is mostly devoted to the introduction of some symbols and operators, along with their properties, that we will often use in the paper when dealing with the Klein-Gordon component of the solution to starting system \eqref{wave KG system}. From now on we will use the notation $p(\xi):=\sqrt{1+|\xi|^2}$ (thus, $p'(\xi)$ denotes the gradient of $p(\xi)$, $p''(\xi)=(\partial^2_{ij}p(\xi))_{ij}$ the $2\times 2$ Hessian matrix of $p(\xi)$). \index{pxi@$p(\xi)$, function}

Proposition \ref{Prop : Continuity on H^s} is a general result about continuity on spaces $H^s_h(\mathbb{R}^2)$ of operators with symbols of order $r\in\mathbb{R}$ and generalises theorem 7.11 in \cite{dimassi:spectral}.
Proposition \ref{Prop : Continuity from $L^2$ to L^inf} is a result of continuity from $L^2$ to $H_h^{\rho, \infty}$ of a particular class of operators that will act on the Klein-Gordon component.
In the spirit of \cite{ifrim_tataru:global_bounds} for the Schr\"{o}dinger equation, it allows to pass from uniform norms to the $L^2$ norm losing only a power $h^{-\frac{1}{2}-\beta}$ for a small $\beta>0$ instead of a $h^{-1}$ as for the semi-classical Sobolev injection. 
Proposition \ref{Prop:Continuity Lp-Lp} is, instead, a result of uniform $L^p-L^p$ continuity of such operators, for every $1\le p\le +\infty$.
\begin{prop}[Continuity on $H^s_h$]\label{Prop : Continuity on H^s}
Let $s\in \mathbb{R}$.
Let $a \in S_{\delta,\sigma}(\langle\xi\rangle^r)$, $r\in\mathbb{R}$, $\delta \in [0, \frac{1}{2}]$, $\sigma \ge 0$.
Then $\oph(a)$ is uniformly bounded : $H^s_h(\mathbb{R}^2)\rightarrow H^{s-r}_h(\mathbb{R}^2)$ and there exists a positive constant $C$ independent of $h$ such that 
\begin{equation*}
\|\oph(a)\|_{\mathcal{L}(H^s_h;H^{s-r}_h)}\le C\, ,  \qquad \forall h\in ]0,1]\, .
\end{equation*}
\end{prop}
\begin{prop} [Continuity from $L^2$ to $H_h^{\rho ,\infty}$]\label{Prop : Continuity from $L^2$ to L^inf}
Let $\rho \in \mathbb{N}$.
Let $a \in S_{\delta,\sigma}(\langle \frac{x-p'(\xi)}{\sqrt{h}}\rangle^{-1})$, $\delta \in [0, \frac{1}{2}]$, $\sigma>0$.
Then $\oph(a)$ is bounded : $L^2(\mathbb{R}^2)\rightarrow H_h^{\rho, \infty}(\mathbb{R}^2)$ and there exists a positive constant $C$ independent of $h$ such that 
\begin{equation*}
\|\oph(a)\|_{\mathcal{L}(L^2;H_h^{\rho, \infty})}\le C h^{-\frac{1}{2}-\beta}\, ,  \qquad \forall h\in ]0,1]\, ,
\end{equation*}
where $\beta>0$ depends linearly on $\sigma$. 
\proof
We first remark that, after definition \ref{def of h-Sobolev spaces} $(i)$ of the $H_h^{\rho,\infty}$ norm,
$$\|\oph(a)w\|_{H^{\rho,\infty}_h} = \|\langle hD_x\rangle^\rho \oph(a)w\|_{L^\infty},$$
and that, by symbolic calculus of lemma \ref{Lem : a sharp b}, $\langle \xi\rangle^\rho \sharp a(x,\xi)$ belongs to
$S_{\delta,\sigma}(\langle\xi\rangle^\rho \big\langle \frac{x-p'(\xi)}{\sqrt{h}}\big\rangle^{-1})\subset h^{-\rho\sigma } S_{\delta,\sigma}( \big\langle \frac{x-p'(\xi)}{\sqrt{h}}\big\rangle^{-1})$. 
This means that estimating the $H_h^{\rho,\infty}$ norm of an operator whose symbol is rapidly decaying in $|h^{\sigma}\xi|$ corresponds actually to estimate the $L^{\infty}$ norm of an operator associated to another symbol (namely, $\tilde{a}(x,\xi)= \langle \xi \rangle^\rho \sharp a(x,\xi) $) which is still in the same class as $a$, up to a small loss $h^{-\rho\sigma}$.

From definition \ref{Def: Weyl and standard quantization} $(i)$ of $\oph(a)w$, and using a change of coordinates $y\mapsto \sqrt{h}y$, $\xi\mapsto \sqrt{h}\xi$, integration by part, Cauchy-Schwarz inequality, and Young's inequality for convolutions, we derive what follows:
\begin{equation} \label{form 3.19}
\begin{split}
& |\oph(a)w| = \\
& = \left| \frac{1}{(2\pi)^2}\int\int e^{i(\frac{x}{\sqrt{h}}-y)\cdot\xi}a\Big(\frac{x+\sqrt{h} y}{2},\sqrt{h}\xi\Big) w(\sqrt{h}y ) \, dyd\xi \right|\\
& = \left|\frac{1}{(2\pi)^4 h}\int\hat{w}\Big(\frac{\eta}{\sqrt{h}}\Big)d\eta \int\int e^{i(\frac{x}{\sqrt{h}}-y)\cdot\xi + i\eta\cdot y}a\Big(\frac{x+\sqrt{h} y}{2},\sqrt{h}\xi\Big)\, dyd\xi \right| \\
& = \left|\frac{1}{(2\pi)^4 h}\int \hat{w}\Big(\frac{\eta}{\sqrt{h}}\Big) \int\int\left(\frac{1-i\big(\frac{x}{\sqrt{h}}-y\big)\cdot\partial_{\xi}}{1+|\frac{x}{\sqrt{h}}-y|^2}\right)^{3}\left(\frac{1+i(\xi -\eta)\cdot\partial_y}{1+|\xi -\eta|^2}\right)^{3}e^{i(\frac{x}{\sqrt{h}}-y)\cdot\xi + i\eta\cdot y} \right. \\
& \left. \hspace{0.5 cm} \times \, a\Big(\frac{x+\sqrt{h} y}{2},\sqrt{h}\xi\Big)\, dyd\xi d\eta  \right| \\
& \lesssim \frac{1}{h} \int \left|\hat{w}\Big(\frac{\eta}{\sqrt{h}}\Big)\right| \int\int \Big\langle \frac{x}{\sqrt{h}}- y\Big\rangle^{-3} \langle \xi - \eta \rangle^{-3}\langle h^{\sigma}\sqrt{h}\xi \rangle^{-N} \Big\langle\frac{\frac{x+\sqrt{h} y}{2}- p'(\sqrt{h}\xi)}{\sqrt{h}}\Big\rangle^{-1} dy d\xi d\eta\\
& \lesssim \frac{1}{h}\left\|\hat{w}\Big(\frac{\cdot}{\sqrt{h}}\Big)\right\|_{L^2} \|\langle \eta \rangle^{-3}\|_{L^1(\eta)}\, \left\| \int \Big\langle \frac{x}{\sqrt{h}}-y\Big\rangle^{-3}\langle h^{\sigma}\sqrt{h}\xi \rangle^{-N} \Big\langle \frac{\frac{x+ \sqrt{h}y}{2}-p'(\sqrt{h}\xi)}{\sqrt{h}}\Big\rangle^{-1} dy \right\|_{L^2(d\xi)} \\
& \lesssim h^{-\frac{1}{2}}\|w\|_{L^2} \int \left\langle \frac{x}{\sqrt{h}}-y \right\rangle^{-3} \Big\|\langle h^{\sigma}\sqrt{h}\xi \rangle^{-N}\Big\langle \frac{\frac{x+ \sqrt{h}y}{2}-p'(\sqrt{h}\xi)}{\sqrt{h}}\Big\rangle^{-1} \Big\|_{L^2(\xi)} dy \, ,
\end{split}
\end{equation}
where $N>0$ will be properly chosen later.
We draw attention to two facts: in the third equality in \eqref{form 3.19} we use that 
$$\left(\frac{1-i(\frac{x}{\sqrt{h}}-y)\cdot\partial_{\xi}}{1+(\frac{x}{\sqrt{h}}-y)^2}\right)^{3}\left(\frac{1+i(\xi -\eta)\cdot\partial_y}{1+(\xi -\eta)^2}\right)^{3} \left[e^{i(\frac{x}{\sqrt{h}}-y)\cdot\xi + i\eta\cdot y}\right]=e^{i(\frac{x}{\sqrt{h}}-y)\cdot\xi + i\eta\cdot y}$$
so, integrating by part, derivatives $\partial_y, \partial_\xi$ fall on $\langle \frac{x}{\sqrt{h}}-y\rangle^{-1}$, $\langle \xi - \eta\rangle^{-1}$ (giving rise to more decreasing factors) and/or on $a\left(\frac{x+\sqrt{h}y}{2},\sqrt{h}\xi\right)$;
symbol $a$ belongs to $S_{\delta,\sigma}(1)$ with $\delta \le \frac{1}{2}$, but the loss of $h^{-\delta}$ is offset by the factor $\sqrt{h}$ coming from the derivation of $a(\frac{x+\sqrt{h}y}{2},\sqrt{h}\xi)$ with respect to $y$ and $\xi$.

In order to estimate $\big\|\langle h^{\sigma}\sqrt{h}\xi \rangle^{-N}\big\langle\frac{\frac{x+\sqrt{h} y}{2}- p'(\sqrt{h}\xi)}{\sqrt{h}}\big\rangle^{-1}\big\|_{L^2_\xi}$ we first introduce a smooth cut-off function $\chi(\frac{x+\sqrt{h}y}{2})$, with $\chi$ supported in some ball $B_C(0)$, to distinguish between the case when $\frac{x+\sqrt{h}y}{2}$ is bounded from the one where $|\frac{x+\sqrt{h}y}{2}|\rightarrow +\infty$.
In the latter situation, say for $|\frac{x + \sqrt{h}y}{2}|>2$, we have $\big\langle\frac{\frac{x+\sqrt{h} y}{2}-p'(\sqrt{h}\xi)}{\sqrt{h}}\big\rangle^{-1} \lesssim \sqrt{h}$ and 
\begin{equation*}
\Big|(1-\chi)\Big(\frac{x+\sqrt{h}y}{2}\Big) \Big|\Big\|\langle h^{\sigma}\sqrt{h}\xi \rangle^{-N}\Big\langle\frac{\frac{x+\sqrt{h} y}{2}-p'(\sqrt{h}\xi)}{\sqrt{h}}\Big\rangle^{-1}\Big\|_{L^2(d\xi)} \lesssim h^{-\sigma}.
\end{equation*}

On the other hand, when $\frac{x+\sqrt{h}y}{2}$ is bounded we consider a Littlewood-Paley decomposition and write
\begin{equation} \label{summation over k}
\begin{split}
\left\|\langle h^{\sigma}\sqrt{h}\xi \rangle^{-N}\Big\langle\frac{\frac{x+\sqrt{h} y}{2}-p'(\sqrt{h}\xi)}{\sqrt{h}}\Big\rangle^{-1}\right\|_{L^2(\xi)}^2 & = h^{-1}\sum_{k\ge 0} \int  \langle h^{\sigma}\xi \rangle^{-2N}\Big\langle\frac{\frac{x+\sqrt{h} y}{2}-p'(\xi)}{\sqrt{h}}\Big\rangle^{-2} \varphi_k(\xi) d\xi \\
&= h^{-1}\sum_{k\ge 0} I_k 
\end{split}
\end{equation}
where
\begin{equation*}
I_0= \int \langle h^{\sigma}\xi \rangle^{-2N} \Big\langle\frac{\frac{x+\sqrt{h} y}{2}-p'(\xi)}{\sqrt{h}}\Big\rangle^{-2} \varphi_0(\xi) d\xi 
\end{equation*}
and
\begin{equation} \label{I_k}
\begin{split}
I_k &= \int \langle h^{\sigma}\xi \rangle^{-2N} \Big\langle\frac{\frac{x+\sqrt{h} y}{2}-p'(\xi)}{\sqrt{h}}\Big\rangle^{-2} \varphi(2^{-k}\xi) d\xi \\
& = 2^{2k} \int \langle h^{\sigma}2^k \xi \rangle^{-2N}\Big\langle\frac{\frac{x+\sqrt{h} y}{2}-p'(2^k\xi)}{\sqrt{h}}\Big\rangle^{-2} \varphi(\xi) d\xi  \\
& \lesssim 2^{(-2N+2)k}h^{-2\sigma N} \int \Big\langle\frac{\frac{x+\sqrt{h} y}{2}-p'(2^k\xi)}{\sqrt{h}}\Big\rangle^{-2} \varphi(\xi) d\xi \, .
\end{split} \qquad k\ge 1
\end{equation}
For a fixed $k_0$ and any $k\le k_0$, $|\det(p''(2^k\xi))|\ge C> 0$ on the support of $\varphi$.
For $k\ge k_0$, function $\xi \rightarrow g_k(\xi)= 2^{3k}(\frac{x+\sqrt{h} y}{2})- 2^{3k}p'(2^k\xi)$ is such that $\text{det}(g_k'(\xi)) = \frac{2^{4k}}{(1+|2^k\xi|^2)^2}$ and $|\text{det}(g_k'(\xi))|\sim 1$ on the support of $\varphi$.
We may thus split the $d\xi$ integral in a finite number (independent of $k$) of integrals, computed on compact domains, on which $\xi \mapsto g_k(\xi)$ is a change of variables with Jacobian of size 1. 
We are then reduced to estimate $2^{(-2N+2)k} h^{-2\sigma N}\int_{|z|\le C} \langle \frac{z+g_k(\xi_0)}{2^{3k}\sqrt{h}}\rangle^{-2} dz$, where $C$ is a positive constant and $\xi_0$ is in $supp\varphi$. Since we assumed that $\frac{x+\sqrt{h} y}{2}$ is bounded, $|g_k(\xi_0)| = O(2^{3k})$ and we get
\begin{equation*}
\begin{split}
I_k &\lesssim 2^{(-2N+2)k} h^{-2\sigma N}\int_{|z|\lesssim 2^{3k}} \Big\langle\frac{z}{2^{3k}\sqrt{h}}\Big\rangle^{-2} dz \\
& \lesssim 2^{(-2N+8)k}h^{-2\sigma N}h \int_{|z|\lesssim h^{-1/2}} \langle z \rangle^{-2}  dz \\
& \lesssim 2^{(-2N+8)k}h^{-2\sigma N + 1} \log (h^{-1})\, .
\end{split}
\end{equation*}
Taking the sum of all $I_k$ for $k\ge 0$ we then deduce that
\begin{equation*}
\left\|\langle h^{\sigma}\sqrt{h}\xi\rangle^{-N}\Big\langle\frac{\frac{x+\sqrt{h} y}{2}-p'(\sqrt{h}\xi)}{\sqrt{h}}\Big\rangle^{-1}\right\|_{L^2(\xi)}\lesssim h^{-\sigma N -\delta}\Big(\sum_{k\ge 0}2^{(-2N +8)k}\Big)^\frac{1}{2} \lesssim h^{-\sigma N - \delta} \, ,  
\end{equation*}
for $\delta>0$ as small as we want, if we choose $N>0$ such that $-2N +8<0$ (e.g. $N=5$).
Finally
\begin{equation*}
\|\oph(a)\|_{\mathcal{L}(L^2; H_h^{\rho, \infty})} = O( h^{-\frac{1}{2}- \beta})\,,
\end{equation*}
where $\beta(\sigma) = (N + \rho)\sigma + \delta$.
\endproof
\end{prop}
The following lemma is as simple as useful and will be largely recalled from subsection \ref{Subsection : The Derivation of the ODE Equation} on. It is also useful to introduce now the following manifold\index{Lambdakg@$\Lkg$, manifold associated to the Klein-Gordon equation}
\begin{equation}\label{def_Lkg}
\Lkg:=\{(x,\xi)\in\mathbb{R}^2\times\mathbb{R}^2 : x-p'(\xi)=0\}
\end{equation}
which appears to be the graph of function $\xi=-d\phi(x)$, with $\phi(x)=\sqrt{1-|x|^2}$ (see picture \ref{picture: Lkg}).

\begin{lem} \label{Lem:family_thetah}
Let $\gamma,\chi\in C^\infty_0(\mathbb{R}^2)$ be equal to 1 in a neighbourhood of the origin and with sufficiently small support, and $\sigma>0$ be small.
There exists a family of smooth functions $\theta_h(x)$, real valued, equal to 1 for $|x|\le 1-ch^{2\sigma}$ and supported for $|x|\le 1-c_1h^{2\sigma}$, for some $0<c_1<c$, with $\|\partial^{\alpha}\theta_h\|_{L^{\infty}}=O(h^{-2|\alpha|\sigma})$ and $(h\partial_h)^k\theta_h$ bounded for every $k\in\mathbb{N}$, such that
\begin{equation} \label{cut-off-thetah}
\gamma\left(\frac{x-p'(\xi)}{\sqrt{h}}\right)\chi(h^\sigma\xi)= \theta_h(x) \gamma\left(\frac{x-p'(\xi)}{\sqrt{h}}\right)\chi(h^\sigma\xi).
\end{equation}
\proof
Straightforward after observing that function $\gamma\big(\frac{x-p'(\xi)}{\sqrt{h}}\big)\chi(h^\sigma\xi)$ is localized around manifold $\Lkg$, meaning that its support is included in $\{(x,\xi) | |\xi| \lesssim h^{-\sigma}, |x| \le 1 - ch^{2\sigma}\}$, for a small $c>0$.
\endproof
\end{lem}

\begin{prop}[Continuity from $L^p$ to $L^p$]\label{Prop:Continuity Lp-Lp}
Let $\gamma,\chi\in C^\infty_0(\mathbb{R}^2)$ be equal to 1 in a neighbourhood of the origin and with sufficiently small support, $\Sigma(\xi)=\langle\xi\rangle^\rho$ with $\rho\in\mathbb{N}$, and $\sigma>0$.
Then $\oph\big(\gamma\big(\frac{x-p'(\xi)}{\sqrt{h}}\big)\chi(h^\sigma\xi)\Sigma(\xi)\big):L^p\rightarrow L^p$ is bounded and its $\mathcal{L}(L^p)$ norm is estimated by $h^{-\sigma\rho-\beta}$, for a small $\beta>0$, $\beta\rightarrow 0$ as $\sigma\rightarrow 0$, for every $1\le p\le +\infty$.
\proof
From lemma \ref{Lem:family_thetah} and the fact that $\gamma\left(\frac{x-p'(\xi)}{\sqrt{h}}\right)\chi(h^\sigma\xi)$ is supported in a neighbourhood of $\Lambda_{kg}$ introduced above, we can find a new smooth cut-off function $\gamma_1$, suitably supported, so that\small
\begin{equation*}
\oph\left(\gamma\left(\frac{x-p'(\xi)}{\sqrt{h}}\right)\chi(h^\sigma\xi)\Sigma(\xi)\right)= \oph\left( \gamma\left(\frac{x-p'(\xi)}{\sqrt{h}}\right)\chi(h^\sigma\xi)\Sigma(\xi) \gamma_1\Big(\frac{\xi+d\phi(x)}{h^{1/2-\beta}}\Big)\theta_h(x)\right)
\end{equation*}\normalsize
where $\beta>0$ is a small constant, $\beta\rightarrow 0$ as $\sigma\rightarrow 0$, that takes into account the degeneracy of the equivalence between the two equations of $\Lambda_{kg}$ when approaching the boundary of $supp \theta_h$.
Denoting $\gamma\left(\frac{x-p'(\xi)}{\sqrt{h}}\right)\chi(h^\sigma\xi)\Sigma(\xi)$ concisely by $A(x,\xi)$ and looking at the kernel associated to the above operator
\begin{equation*}
\begin{split}
K(x,y)&:= \frac{1}{(2\pi h)^2}\int e^{\frac{i}{h}(x-y)\cdot\xi} A\left(\frac{x+y}{2},\xi \right) \gamma_1\Big(\frac{\xi+d\phi(\frac{x+y}{2})}{h^{1/2-\beta}}\Big)\theta_h\Big(\frac{x+y}{2}\Big) d\xi\\
&= \frac{e^{-\frac{i}{h}(x-y)\cdot d\phi(\frac{x+y}{2})}}{(2 \pi h)^2}\theta_h\Big(\frac{x+y}{2}\Big) \int e^{\frac{i}{h}(x-y)\cdot\xi} A\left(\frac{x+y}{2},\xi-d\phi\Big(\frac{x+y}{2}\Big) \right) \gamma_1\Big(\frac{\xi}{h^{1/2-\beta}}\Big) d\xi,
\end{split}
\end{equation*}
we observe that, since 
\begin{equation*}
\Big(\frac{x}{\sqrt{h}}\Big)^\alpha e^{\frac{i}{h}(x-y)\cdot\xi} = \Big(\frac{\sqrt{h}}{i}\Big)^{|\alpha|}\partial^\alpha_\xi e^{\frac{i}{h}(x-y)\cdot\xi}
\end{equation*}
and $h^{|\alpha|/2}\partial^\alpha_\xi A(\frac{x+y}{2},\xi)$ is bounded by $h^{-\sigma\rho}$ for any $\alpha\in\mathbb{N}^2$, by making some integration by parts
\begin{equation*}
\left|\Big(\frac{x}{\sqrt{h}}\Big)^\alpha K(x,y) \right| \lesssim h^{-2-\sigma\rho}\int_{|\xi|\lesssim h^{1/2-\beta}} d\xi \lesssim h^{-1-\sigma\rho-2\beta}, \quad \forall (x,y)\in\mathbb{R}^2\times\mathbb{R}^2.
\end{equation*}
This means in particular that
\begin{equation*}
|K(x,y)|\lesssim h^{-1-\sigma\rho-2\beta}\Big\langle \frac{x}{\sqrt{h}}\Big\rangle^{-3}, \quad |K(x,y)|\lesssim h^{-1-\sigma\rho-2\beta}\Big\langle \frac{y}{\sqrt{h}}\Big\rangle^{-3},\quad \forall (x,y)
\end{equation*}
implying that
\begin{equation*}
\sup_x \int |K(x,y)| dy \lesssim h^{-\sigma\rho-2\beta}, \quad \sup_y\int |K(x,y)|dx \lesssim h^{-\sigma\rho-2\beta}.
\end{equation*}
The operator associated to $K(x,y)$ is hence bounded on $L^p$ with norm $O(h^{-\sigma\rho-2\beta})$, for every $1\le p\le +\infty$. 
\endproof
\end{prop}

The following lemma shows that we have nice upper bounds for operators whose symbol is supported for large frequencies $|\xi|\ge h^{-\sigma}$, $\sigma >0$, when acting on functions $w$ that belong to $H^s_h$, for some large $s$. We state it in space dimension 2 but it is clear that it holds true in general space dimension $d\ge 1$. This result is useful when we want to reduce to symbols rapidly decaying in $|h^{\sigma}\xi|$, for example in the intention of using proposition \ref{Prop : Continuity from $L^2$ to L^inf} or when we want to pass from a symbol of a certain positive order to another one of order zero, up to small losses of order $O(h^{-\beta})$, $\beta>0$ depending linearly on $\sigma$. We can always split a symbol using that $1= \chi(h^{\sigma}\xi) + (1-\chi)(h^{\sigma}\xi)$, for a smooth $\chi$ equal to 1 close to the origin, and consider as remainders all contributions coming from the latter. 

\begin{lem} \label{Lem : new estimate 1-chi}
Let $s'\ge 0$ and $\chi \in C^{\infty}_0(\mathbb{R}^2)$, $\chi \equiv 1$ in a neighbourhood of zero.
Then
\begin{equation*}
\|\oph((1-\chi)(h^{\sigma}\xi))w\|_{H^{s'}_h}\le C h^{\sigma(s-s')}\|w\|_{H^s_h} \, , \qquad \qquad \forall s> s'\, .
\end{equation*}
\proof
The result is a simple consequence of the fact that $(1-\chi)(h^{\sigma}\xi)$ is supported for $|\xi|\gtrsim h^{-\sigma}$, because
\begin{equation*}
\begin{split}
\|\oph((1-\chi)(h^{\sigma}\xi))w \|_{H^{s'}_h}^2 & = \int (1+ |h\xi|^2)^{s'}|(1-\chi)(h^{\sigma}h\xi)|^2 |\hat{w}(\xi)|^2 d\xi \\
& = \int (1+ |h\xi|^2)^s (1+ |h\xi|^2)^{s'-s}|(1-\chi)(h^{\sigma}h\xi)|^2 |\hat{w}(\xi)|^2 d\xi \\
& \le C h^{2\sigma(s-s')}\|w\|^2_{H^s_h}\, ,
\end{split}
\end{equation*}
where the last inequality follows from an integration on $|h\xi|\gtrsim h^{-\sigma}$ and from the fact that $s'-s< 0$, $(1+|h\xi|^2)^{s'-s}\le C h^{-2\sigma(s'-s)}$.
\endproof
\end{lem}

We introduce the following operator: \index{Lj@$\mathcal{L}_j$, operator}
\begin{equation}\label{def Lj}
\mathcal{L}_j : = \frac{1}{h}\oph(x - p'_j(\xi)), \quad j=1,2,
\end{equation}
and use the notation $\|\mathcal{L}^\gamma w\| = \|\mathcal{L}^{\gamma_1}_1 \mathcal{L}^{\gamma_2}_2 w\|$ for any $\gamma=(\gamma_1,\gamma_2)\in\mathbb{N}^2$.

\begin{lem} \label{Lem:class of gamma c(x,xi)}
Let $\gamma\in C^\infty_0(\mathbb{R}^2)$ be equal to 1 in a neighbourhood of the origin, $c(x,\xi)\in S_{\delta,\sigma}(1)$ with $\delta\in [0,\frac{1}{2}[$ and $\sigma>0$. Then $\gamma(\frac{x-p'(\xi)}{\sqrt{h}})c(x,\xi)$ belongs to $S_{\frac{1}{2},\sigma}(1)\big(\langle\frac{x-p'(\xi)}{\sqrt{h}}\rangle^{-N}\big)$, for all $N\ge 0$. 
\proof
Straightforward.
\endproof
\end{lem}

\begin{lem} \label{Lem : composition gamma-1 and its argument}
Let $n\in\mathbb{N}$ and $\gamma_n(z)$ be a smooth function such that $|\partial^{\alpha}\gamma_n(z)| \lesssim \langle z\rangle^{-|\alpha|-n}$ for all $\alpha\in\mathbb{N}^2$.
Let also $c(x,\xi) \in S_{\delta,\sigma}(1)$, with $\delta \in [0,\frac{1}{2}[$, $\sigma>0$, be supported for $|\xi|\lesssim h^{-\sigma}$.
Up to some multiplicative constants independent of $h$, we have the following equality:
\begin{multline} \label{sharp gammatilde and its argument}
\left[c(x,\xi)\gamma_n\Big(\frac{x-p'(\xi)}{\sqrt{h}}\Big) \right]\sharp \big(x_j- p_j'(\xi)\big) =c(x,\xi)\gamma_n\Big(\frac{x-p'(\xi)}{\sqrt{h}}\Big) \big(x_j- p_j'(\xi)\big) \\
+ h \gamma_n\Big(\frac{x-p'(\xi)}{\sqrt{h}}\Big)\big[(\partial_{\xi_j} c) + (\partial_xc)\cdot (\partial_\xi p'_j)\big] + h \sum_{|\alpha|=2}(\partial^\alpha\gamma_n)\Big(\frac{x-p'(\xi)}{\sqrt{h}}\Big)c(x,\xi)(\partial^\alpha_\xi p'_j)(\xi) + r(x,\xi),
\end{multline} 
with $r \in h^{3/2-\delta} S_{\frac{1}{2},\sigma}\big(\langle\frac{x-p'(\xi)}{\sqrt{h}}\rangle^{-n}\big)$, and if $\chi \in C^\infty_0(\mathbb{R}^2)$ is such that $\chi(h^\sigma\xi)\equiv 1$ on the support of $c(x,\xi)$,
\begin{subequations} \label{est L2 Linfty for GammaTilde}
\begin{equation} \label{est L2 for GammaTilde}
\left\|\oph\Big(c(x,\xi)\gamma_n\Big(\frac{x-p'(\xi)}{\sqrt{h}}\Big)(x_j-p'_j(\xi))\Big)\widetilde{v} \right\|_{L^2} \lesssim \sum_{|\gamma|=0}^1 h^{1-\beta} \|\oph(\chi(h^\sigma\xi))\mathcal{L}^\gamma\widetilde{v}\|_{L^2} \, , 
\end{equation}
\begin{equation} \label{est Linf for GammaTilde}
\left\|\oph\Big(c(x,\xi)\gamma_n\Big(\frac{x-p'(\xi)}{\sqrt{h}}\Big)(x_j-p'_j(\xi))\Big)\widetilde{v} \right\|_{L^{\infty}}  \lesssim \sum_{|\gamma|=0}^1  h^{\frac{1}{2}\delta_n-\beta} \|\oph(\chi(h^\sigma\xi))\mathcal{L}^\gamma\widetilde{v}\|_{L^2} \, ,
\end{equation}
\end{subequations}
where $\delta_n=1$ if $n>0$, 0 otherwise, and $\beta>0$ is small, $\beta\rightarrow 0$ as $\delta,\sigma\rightarrow 0$.

Moreover, if $n\in\mathbb{N}^*$ and $\partial^\alpha\gamma_n$ vanishes in a neighbourhood of the origin whenever $|\alpha|\ge 1$, we also have that
\begin{subequations}  \label{est:L2 Linfty Op(gamma) LLv}
\begin{multline} \label{est:L2 Op(gamma)LLv}
\left\|\oph\Big(c(x,\xi)\gamma_n\Big(\frac{x-p'(\xi)}{\sqrt{h}}\Big)(x_i-p'_i(\xi))(x_j-p'_j(\xi))\Big)\widetilde{v} \right\|_{L^2} \lesssim\\
\sum_{0\le |\gamma| \le 2}  h^{2-\beta} \|\oph(\chi(h^\sigma\xi))\mathcal{L}^\gamma \widetilde{v}\|_{L^2},
\end{multline}
\begin{multline}
\left\|\oph\Big(c(x,\xi)\gamma_n\Big(\frac{x-p'(\xi)}{\sqrt{h}}\Big)(x_i-p'_i(\xi))(x_j-p'_j(\xi))\Big)\widetilde{v} \right\|_{L^{\infty}}  \lesssim\\
\sum_{0\le |\gamma|\le 2}  h^{\frac{3}{2}-\beta}\|\oph(\chi(h^\sigma\xi))\mathcal{L}^\gamma\widetilde{v}\|_{L^2}.
\end{multline}
\end{subequations}
\proof
As $c(x,\xi)\gamma_n\Big(\frac{x-p'(\xi)}{\sqrt{h}}\Big) \in S_{\frac{1}{2},\sigma}\big(\langle\frac{x-p'(\xi)}{\sqrt{h}}\rangle^{-n}\big)$ and $\partial^\alpha_{x,\xi}( x_j-p'_j(\xi))\in S_{0,0}(1)$ for any $|\alpha|\ge 1$, equality \eqref{sharp gammatilde and its argument} follows from the last part of lemma \ref{Lem : a sharp b} and symbolic development \eqref{a sharp b asymptotic formula} until order 2, after having observed that
\begin{equation} \label{poisson brackets gamma and its symbol}
\left\{c(x,\xi)\gamma_n\Big(\frac{x-p'(\xi)}{\sqrt{h}}\Big), x_j-p'_j(\xi)\right\} = \gamma_n\Big(\frac{x-p'(\xi)}{\sqrt{h}}\Big) \big[(\partial_{\xi_j}c) + (\partial_x c)\cdot (\partial_\xi p'_j)\big],
\end{equation}
and that, up to some multiplicative negligible, \small
\begin{multline*}
h^2\sum_{|\alpha|=2}\partial^\alpha_x\left[c(x,\xi)\gamma_n\Big(\frac{x-p'(\xi)}{\sqrt{h}}\Big)\right](\partial^\alpha_\xi p'_j)(\xi) = h \sum_{|\alpha|=2}(\partial^\alpha\gamma_n)\Big(\frac{x-p'(\xi)}{\sqrt{h}}\Big)c(x,\xi)(\partial^\alpha_\xi p'_j)(\xi) \\
+\underbrace{ h^{\frac{3}{2}}\sum_{\substack{|\alpha|=2\\ |\alpha_1|,|\alpha_2|= 1}} (\partial^{\alpha_1}\gamma_n)\Big(\frac{x-p'(\xi)}{\sqrt{h}}\Big)(\partial^{\alpha_2}_xc)(x,\xi)(\partial^\alpha_\xi p'_j)(\xi) + h^2 \sum_{|\alpha|=2}\gamma_n\Big(\frac{x-p'(\xi)}{\sqrt{h}}\Big)(\partial^\alpha_xc)(x,\xi)(\partial^\alpha_\xi p'_j)(\xi) }_{\in h^{\frac{3}{2}-\delta}S_{\frac{1}{2},\sigma}\big(\big\langle\frac{x-p'(\xi)}{\sqrt{h}}\big\rangle^{-n}\big)}.
\end{multline*}\normalsize

If $\chi$ is a cut-off function as in the statement, its derivatives vanish on the support of $c(x,\xi)$, and from remark \ref{Remark:symbols_with_null_support_intersection}
\begin{equation} \label{gamma_n-sharp chi}
c(x,\xi)\gamma_n\Big(\frac{x-p'(\xi)}{\sqrt{h}}\Big) = \left[c(x,\xi)\gamma_n\Big(\frac{x-p'(\xi)}{\sqrt{h}}\Big)\right]\sharp \chi(h^\sigma\xi) + r_\infty(x,\xi)
\end{equation}
with $r_\infty\in h^NS_{\frac{1}{2},\sigma}(\langle \frac{x-p'(\xi)}{\sqrt{h}}\rangle^{-n})$, $N\in\mathbb{N}$ as large as we want. 
Estimates \eqref{est L2 Linfty for GammaTilde} follow then as a straight consequence of \eqref{sharp gammatilde and its argument},
 definition \eqref{def Lj} of $\mathcal{L}_j$, proposition \ref{Prop : Continuity on H^s} and semi-classical Sobolev's injection \eqref{semi-classical Sobolev injection} (resp. proposition \ref{Prop : Continuity from $L^2$ to L^inf}) when $n=0$ (resp. $n>0$).

In order to prove the last part of the statement (estimates \eqref{est:L2 Linfty Op(gamma) LLv}) we use equality \eqref{sharp gammatilde and its argument} with $\gamma_n$ replaced by $\widetilde{\gamma}_{n-1}(z)=\gamma_n(z)z_i$, where $|\partial^\alpha \widetilde{\gamma}_{n-1}(z)|\lesssim \langle z\rangle^{-|\alpha|-(n-1)}$, which gives that
\begin{align*}
 c(x,\xi)\gamma_n\Big(\frac{x-p'(\xi)}{\sqrt{h}}\Big)(x_i-p'_i(\xi))(x_j-p'_j(\xi)) = \left[c(x,\xi)\gamma_n\Big(\frac{x-p'(\xi)}{\sqrt{h}}\Big)(x_i-p'_i(\xi))\right]\sharp (x_j-p'_j(\xi))\\
- h\gamma_n\Big(\frac{x-p'(\xi)}{\sqrt{h}}\Big)(x_i-p'_i(\xi))\left[(\partial_{\xi_j}c) + (\partial_xc)\cdot(\partial_\xi p'_j)\right] \\
 - h^\frac{3}{2}\sum_{|\alpha|=2}(\partial^\alpha \widetilde{\gamma}_{n-1})\Big(\frac{x-p'(\xi)}{\sqrt{h}}\Big)c(x,\xi)(\partial^\alpha_\xi p'_j)(\xi) - \sqrt{h}r(x,\xi),
\end{align*}
with $r\in h^{\frac{3}{2}-\delta}S_{\frac{1}{2},\sigma}\big(\langle\frac{x-p'(\xi)}{\sqrt{h}}\rangle^{-(n-1)}\big)$.
As $\partial^\alpha \widetilde{\gamma}_{n-1}$ vanishes in a neighbourhood of the origin for $|\alpha|=2$ by the hypothesis made on $\gamma_n$, we can rewrite it as $\sum_{l=1}^2 \widetilde{\gamma}^l_{n+2}(z) z_l$, where $\widetilde{\gamma}^l_{n+2}(z):=(\partial^\alpha \widetilde{\gamma}_{n-1})(z)z_l|z|^{-2}$ is such that $|\partial^\beta \widetilde{\gamma}^l_{n+2}(z)|\lesssim \langle z\rangle^{-|\beta|-(n+2)}$.
Then, using again equality \eqref{sharp gammatilde and its argument} for all products different from $r(x,\xi)$ in the above right hand side (with $c$ replaced with $h^\delta[(\partial_{\xi_j}c) - (\partial_xc)\cdot(\partial_\xi p'_j)]$ in the second addend, and $\gamma_n$ and $c$ replaced with
$\widetilde{\gamma}^l_{n+2}$ and $c(\partial^\alpha_\xi p'_j)$ respectively in the third one, $l=1,2$) we find that
\begin{multline*}
c(x,\xi)\gamma_n\Big(\frac{x-p'(\xi)}{\sqrt{h}}\Big)(x_i-p'_i(\xi))(x_j-p'_j(\xi)) = \\
\left[c(x,\xi)\gamma_n\Big(\frac{x-p'(\xi)}{\sqrt{h}}\Big)\right] \sharp (x_i-p'_i(\xi))\sharp  (x_j-p'_j(\xi)) + h r_1(x,\xi)\sharp (x_j-p'_j(\xi)) - \sqrt{h}r(x,\xi),
\end{multline*}
for a new $r_1\in h^{-\delta}S_{\frac{1}{2},\sigma}\big(\langle\frac{x-p'(\xi)}{\sqrt{h}}\rangle^{-n}\big)$.
Estimates \eqref{est:L2 Linfty Op(gamma) LLv} are then obtained using \eqref{gamma_n-sharp chi} and propositions \ref{Prop : Continuity on H^s}, \ref{Prop : Continuity from $L^2$ to L^inf}).
\endproof
\end{lem}
We will also need the following result, which is detailed in lemma 1.2.6 in \cite{delort:semiclassical} for the one-dimensional case.

\begin{lem} \label{Lem : on e and etilde}
Let $\gamma \in C^{\infty}_0(\mathbb{R}^2)$, and $\phi(x)=\sqrt{1-|x|^2}$.
If the support of $\gamma$ is sufficiently small, 
\begin{subequations}
\begin{gather}
(x_k-p'_k(\xi))\gamma\big(\langle\xi\rangle^2(x-p'(\xi))\big) = \sum_{l=1}^2 e^k_l(x,\xi)(\xi_l + d_l\phi(\xi)), \\
(\xi_k + d_k\phi(x)) \gamma\big(\langle\xi\rangle^2(x-p'(\xi))\big) = \sum_{l=1}^2 \widetilde{e}^k_l(x,\xi)(x_l-p'_l(\xi)),
\end{gather}
\end{subequations}
for any $k=1,2$, where functions $e^k_l(x,\xi), \widetilde{e}^k_l(x,\xi)$ are such that, for any $\alpha,\beta\in\mathbb{N}^2$,
\begin{subequations}
\begin{align}
|\partial^{\alpha}_x \partial_{\xi}^{\beta}e^k_l(x,\xi)| &\lesssim_{\alpha\beta} \langle\xi\rangle^{-3+2|\alpha|-|\beta|}\,,\\
|\partial^{\alpha}_x \partial_{\xi}^{\beta}\widetilde{e}^k_l(x,\xi)| &\lesssim_{\alpha\beta} \langle\xi\rangle^{3+2|\alpha|-|\beta|} \,,\label{etildekj}
\end{align}
\end{subequations}
for any $k,l=1,2$.
\end{lem}

\chapter{Energy Estimates} \label{Chap:Energy estimates}

The aim of this chapter is to write an energy inequality for $E_n(t;u_\pm, v_\pm)$ and $E^k_3(t;u_\pm, v_\pm)$ respectively, which allows us to propagate the a-priori energy estimates made in theorem \ref{Thm: bootstrap argument}, i.e. to pass from \eqref{est: bootstrap argument a-priori est} to \eqref{est: bootstrap enhanced Enn}, \eqref{est: boostrap enhanced E02}.
Such an inequality is in general derived by computing and estimating the derivative in time of the energy, i.e. of the $L^2$ norm to the square of $u^I_\pm, v^I_\pm$.
As this computation makes use of the system of equations satisfied by $(u^I_\pm, v^I_\pm)$ (see \eqref{system for uI+-, vI+-}), two main difficulties arise due to the quasi-linear nature of the starting problem and the very slow decay in time \eqref{est: bootstrap upm} of the wave solution.

On the one hand, among all quadratic terms appearing in the right hand side of \eqref{system for uI+-, vI+-} we find the quasi-linear ones $Q^\mathrm{w}_0(v_\pm, D_1 v^I_\pm)$ and $Q^{\mathrm{kg}}_0(v_\pm, D_1u^I_\pm)$, whose $L^2$ norm is bounded by $\|v_\pm(t,\cdot)\|_{H^{1,\infty}}(\|u^I_\pm(t,\cdot)\|_{H^1}+ \|v^I_\pm(t,\cdot)\|_{H^1})$, as usual for this kind of terms. 
This means that they are at the wrong energy level, in the sense that they cannot be controlled in $L^2$ by $E_n(t;u_\pm, v_\pm)$ or $E^k_3(t;u_\pm, v_\pm)$.
This causes a "loss of derivatives" in the energy inequality if we roughly estimate
\[\frac{1}{2}\partial_t \Big(\|u^I_\pm(t,\cdot)\|^2_{L^2}+ \|v^I_\pm(t,\cdot)\|^2_{L^2}\Big) = -\Im \Big[\langle Q^\mathrm{w}_0(v_\pm, D_1 v^I_\pm), u^I_\pm\rangle + \langle Q^{\mathrm{kg}}_0(v_\pm, D_1u^I_\pm), v^I_\pm\rangle + \dots\Big]\]
using the Cauchy-Schwarz inequality.
This issue is however only technical.
In fact, by writing system \eqref{system for uI+-, vI+-} in a vectorial fashion and para-linearising it in order to stress out the very troublesome terms (see subsection \ref{Subsection: Paralinearization}) we are able to \emph{symmetrize} it, i.e. to derive an equivalent system in which the quasi-linear contribution is represented by a self-adjoint operator of order 1 (see subsection \ref{Subs: Symmetrization}, proposition \ref{Prop: equation of WIs}).
As this operator is self-adjoint it essentially disappears in the energy inequality, replaced with an operator of order 0 whose action on $u^I_\pm, v^I_\pm$ is bounded in $L^2$ by $E_n(t;u_\pm, v_\pm)$ or $E^k_3(t;u_\pm, v_\pm)$, depending on the multi-index $I$ we are dealing with.

On the other hand, the $L^2$ norm of some semi-linear contributions to the right hand side of \eqref{system for uI+-, vI+-} decays very slowly in time. It is the case, for instance, of $Q^\mathrm{kg}_0(v^I_\pm, D_1u_\pm)$, whose $L^2$ norm is bounded by $\|u_\pm(t,\cdot)\|_{H^{2,\infty}}\|v^I_\pm(t,\cdot)\|_{L^2}$ and only has the slow decay \eqref{est: bootstrap upm} of the wave component $u_\pm$.
Since we want to prove that
\[\partial_t E_n(t;u_\pm, v_\pm) = O\big(\varepsilon t^{-1+\frac{\delta}{2}} E_n(t;u_\pm, v_\pm)^\frac{1}{2}\big), \quad \partial_t E^k_3(t;u_\pm, v_\pm) = O\big(\varepsilon t^{-1+\frac{\delta_k}{2}} E^k_3(t;u_\pm, v_\pm)^\frac{1}{2}\big)\]
we need to get rid of such terms by means of normal forms (see section \ref{Section : Normal Forms for system}).
Because of the quasi-linear nature of our problem, some of them will be eliminated by an adapted quasi-linear normal form argument (see subsection \ref{sub: a first normal form transformation}), while the remaining ones can be treated with an usual semi-linear one (see subsection \ref{sub: second normal form}).
At that point we will be able to prove proposition \ref{Prop: Propagation of the energy estimate} and to derive estimates \eqref{est: bootstrap enhanced Enn}, \eqref{est: boostrap enhanced E02}.

\section{Paralinearization and Symmetrization} \label{Section: Paralinearization and Symmetrization}

As anticipated above, the first step towards the derivation of the right energy inequality is to handle the quasi-linear terms appearing in the right hand side of \eqref{system for uI+-, vI+-} in order to avoid any loss of derivatives.
We realize that the very quasi-linear contribution to our system appears in equation \eqref{equation WI-1} through a para-differential operator whose symbol is a real \emph{non symmetric} matrix. As we need this operator to be self-adjoint (up to an operator of order 0), we \emph{symmetrize} equation \eqref{equation WI-1} by defining a new function $W^I_s$ in terms of $W^I$, that will be solution to a new equation in which the symbol of the quasi-linear contribution is a real symmetric matrix (see subsection \ref{Subs: Symmetrization}). 
Also, we set aside subsection \ref{Subs: Estimate of quadratic terms} to the estimate of the $L^2$ norms of the non-linear terms in the right hand side of \eqref{equation WI-1}.

\subsection{Paralinearization}  \label{Subsection: Paralinearization}

Let us remind definitions \eqref{def uIpm vIpm} and \eqref{def mathcal(I)}.
Since admissible vector fields considered in $\mathcal{Z}=\{\Omega, Z_j, \partial_j, j=1,2\}$ exactly commute with the linear part of system \eqref{wave KG system}, we deduce from remark \ref{Remark:Vector_field_on_null_structure} and \eqref{Gamma_I_nonlinearity} that, for any multi-index $I$, $(\Gamma^Iu, \Gamma^Iv)$ is solution to
\begin{equation*}
\begin{cases}
\left(\partial^2_t - \Delta_x\right)\Gamma^Iu = \displaystyle\sum_{\substack{(I_1,I_2)\in\mathcal{I}(I)\\ |I_1|+|I_2|=|I|}} Q_0(\Gamma^{I_1}v, \partial_1\Gamma^{I_2}v)+ \sum_{\substack{(I_1,I_2)\in\mathcal{I}(I)\\ |I_1|+|I_2|<|I|}} c_{I_1,I_2}Q_0(\Gamma^{I_1}v, \partial\Gamma^{I_2}v) ,\\
\left(\partial^2_t - \Delta_x+ 1\right)\Gamma^Iv = \displaystyle\sum_{\substack{(I_1,I_2)\in\mathcal{I}(I)\\ |I_1|+|I_2|=|I|}} Q_0(\Gamma^{I_1}v, \partial_1\Gamma^{I_2}u)+ \sum_{\substack{(I_1,I_2)\in\mathcal{I}(I)\\ |I_1|+|I_2|<|I|}} c_{I_1,I_2}Q_0(\Gamma^{I_1}v, \partial\Gamma^{I_2}u) ,
\end{cases}
\end{equation*}
with coefficients $c_{I_1,I_2}\in \{-1,0,1\}$ such that $c_{I_1,I_2}=1$ for $|I_1|+|I_2|=|I|$, in which case the derivative $\partial$ acting on $\Gamma^{I_2}v$ (resp. on $\Gamma^{I_2}u$) is equal to $\partial_1$, and $\partial$ representing one of the partial derivatives $\partial_a$, $a\in \{0,1,2\}$.
Let us remind that, if $\Gamma^I$ contains at least $k$ ($\le |I|$) space derivatives, above summations are taken over indices $I_1,I_2$ such that $k \le |I_1|+|I_2|\le |I|$.
Hence, introducing from \eqref{null form Q0}, \eqref{def u+- v+-}, \index{Qw0@$Q^{\mathrm{w}}_0(v_\pm, D_a v_\pm)$, null form} \index{Qkg0@$Q^{\mathrm{kg}}_0(v_\pm, D_a u_\pm)$, null form}
\begin{equation}  \label{Q0_pm}
\begin{split}
Q^{\mathrm{w}}_0(v_\pm, D_a v_\pm) & := \frac{i}{4}\left[(v_+ + v_{-})D_a(v_+ + v_{-}) - \frac{D_x}{\langle D_x \rangle}(v_+ - v_{-})\cdot\frac{D_x D_a}{\langle D_x \rangle}(v_+ - v_{-})\right] ,\\
 Q^{\mathrm{kg}}_0(v_\pm, D_a u_\pm) & := \frac{i}{4} \left[(v_+ + v_{-})D_a(u_+ + u_{-}) - \frac{D_x}{\langle D_x\rangle}(v_+ - v_{-})\cdot\frac{D_x D_a}{|D_x|}(u_+ - u_{-})\right],
\end{split}
\end{equation}
for any $a=0,1,2$, we deduce that $(u^I_+, v^I_+, u^I_{-}, v^I_{-})$ is solution to 
\begin{equation} \label{system for uI+-, vI+-}
\begin{cases}
& (D_t - |D_x|)u^I_+(t,x) = \displaystyle\sum_{\substack{(I_1,I_2)\in\mathcal{I}(I)\\ |I_1| + |I_2| = |I| }}Q_0^{\mathrm{w}}(v^{I_1}_\pm, D_1 v^{I_2}_\pm) + \sum_{\substack{(I_1,I_2)\in\mathcal{I}(I)\\|I_1| + |I_2| < |I|} } c_{I_1, I_2}Q_0^{\mathrm{w}}(v^{I_1}_\pm, D v^{I_2}_\pm)  \\
& (D_t - \langle D_x\rangle)v^I_+(t,x) =  \displaystyle\sum_{\substack{(I_1,I_2)\in\mathcal{I}(I)\\|I_1| + |I_2| = |I|} } Q_0^{\mathrm{kg}}(v^{I_1}_\pm, D_1 u^{I_2}_\pm) + \sum_{\substack{(I_1,I_2)\in\mathcal{I}(I)\\|I_1| + |I_2| < |I| }} c_{I_1, I_2}Q_0^{\mathrm{kg}}(v^{I_1}_\pm, D u^{I_2}_\pm)  \\
& (D_t + |D_x|)u^I_{-}(t,x) = \displaystyle\sum_{\substack{(I_1,I_2)\in\mathcal{I}(I)\\ |I_1| + |I_2| = |I| }}Q_0^{\mathrm{w}}(v^{I_1}_\pm, D_1 v^{I_2}_\pm) + \sum_{\substack{(I_1,I_2)\in\mathcal{I}(I)\\|I_1| + |I_2| < |I|} } c_{I_1, I_2}Q_0^{\mathrm{w}}(v^{I_1}_\pm, D v^{I_2}_\pm)    \\
& (D_t + \langle D_x \rangle)v^I_{-}(t,x) = \displaystyle\sum_{\substack{(I_1,I_2)\in\mathcal{I}(I)\\|I_1| + |I_2| = |I|} } Q_0^{\mathrm{kg}}(v^{I_1}_\pm, D_1 u^{I_2}_\pm) + \sum_{\substack{(I_1,I_2)\in\mathcal{I}(I)\\|I_1| + |I_2| < |I| }} c_{I_1, I_2}Q_0^{\mathrm{kg}}(v^{I_1}_\pm, D u^{I_2}_\pm)
\end{cases}
\end{equation}
The quasi-linear structure of the above system can be emphasized by using \eqref{dev in paraproduct} and decomposing $Q_0^{\mathrm{w}}(v_\pm, D_1 v^I_\pm)$, $Q_0^{\mathrm{kg}}(v_\pm, D_1 u^I_\pm)$ as follows: 
\begin{equation} \label{dec quasi-linear term}
 Q_0^{\mathrm{w}}(v_\pm, D_1 v^I_\pm) = (QL)_1 + (SL)_1,\quad
 Q_0^{\mathrm{kg}}(v_\pm, D_1 u^I_\pm) =(QL)_2 + (SL)_2 ,
\end{equation}
with
\begin{equation*}
\begin{split}
& (QL)_1 := \frac{i}{4}\left[Op^B\big((v_+ + v_{-})\eta_1\big)(v^I_+ + v^I_{-}) - Op^B\Big(\frac{D_x}{\langle D_x\rangle}(v_+ - v_{-})\cdot\frac{\eta \eta_1}{\langle\eta\rangle}\Big)(v^I_+- v^I_{-}) \right] ,\\
& (SL)_1 := \frac{i}{4}\left[Op^B\big(D_1(v^I_+ + v^I_{-})\big)(v_+ + v_{-}) - Op^B\Big(\frac{D_x D_1}{\langle D_x\rangle}(v^I_+ - v^I_{-})\cdot\frac{\eta}{\langle\eta\rangle}\Big)(v_+ - v_{-})  \right. \\
& \left. \hspace{2cm }+ Op^B_R\big((v_+ + v_{-})\eta_1\big)(v^I_+ + v^I_{-})\big) - Op^B_R\Big(\frac{D_x}{\langle D_x\rangle}(v_+ - v_{-})\cdot\frac{\eta \eta_1}{\langle \eta\rangle}\Big)(v^I_+ - v^I_{-})\right], \\
& (QL)_2 := \frac{i}{4}\left[Op^B\big((v_+ + v_{-})\eta_1\big)(u^I_+ + u^I_{-}) - Op^B\Big(\frac{D_x}{\langle D_x\rangle}(v_+ - v_{-})\cdot\frac{\eta \eta_1}{|\eta|}\Big)(u^I_+- u^I_{-}) \right] ,\\
& (SL)_2 := \frac{i}{4}\left[Op^B\big(D_1(u^I_+ + u^I_{-})\big)(v_+ + v_{-}) - Op^B\Big(\frac{D_x D_1}{|D_x|}(u^I_+ - u^I_{-})\cdot\frac{\eta}{\langle\eta\rangle}\Big)(v_+ - v_{-})  \right. \\
& \left. \hspace{2cm }+ Op^B_R\big((v_+ + v_{-})\eta_1\big)(u^I_+ + u^I_{-}) - Op^B_R\Big(\frac{D_x}{\langle D_x\rangle}(v_+ - v_{-})\cdot\frac{\eta \eta_1}{|\eta|}\Big)(u^I_+- u^I_{-}) \right]\,,
\end{split}
\end{equation*}
where the Bony quantization $Op^B$ (resp. $Op^B_R$) has been defined in \ref{Def: Paradiff_operator} (resp. in \eqref{operator OpB_R}).
We do a similar decomposition also for the semi-linear contribution $Q^{\mathrm{kg}}_0(v^I_\pm, D_1 u_\pm)$, for this term will thereafter be the object of the two normal forms mentioned at the beginning of this section: 
\begin{equation} \label{dec semi-linear term}
\begin{split}
Q_0^{\mathrm{kg}}(v^I_\pm, D_1 u_\pm) &=  \frac{i}{4}\left[Op^B\big((v^I_+ + v^I_{-})\eta_1\big)(u_+ + u_{-}) - Op^B\Big(\frac{D_x}{\langle D_x\rangle}(v^I_+ - v^I_{-})\cdot\frac{\eta \eta_1}{|\eta|}\Big)(u_+- u_{-}) \right] \\
& +\frac{i}{4}\left[Op^B\big(D_1(u_+ + u_{-})\big)(v^I_+ + v^I_{-}) - Op^B\Big(\frac{D_x D_1}{|D_x|}(u_+ - u_{-})\cdot\frac{\eta}{\langle\eta\rangle}\Big)(v^I_+ - v^I_{-})  \right] \\
& +\frac{i}{4}\left[ Op^B_R\big((v^I_+ + v^I_{-})\eta_1\big)(u_+ + u_{-}) - Op^B_R\Big(\frac{D_x}{\langle D_x\rangle}(v^I_+ - v^I_{-})\cdot\frac{\eta \eta_1}{|\eta|}\Big)(u_+- u_{-}) \right]\,. 
\end{split}
\end{equation}

In order to handle system \eqref{system for uI+-, vI+-} in the most efficient way we proceed to write it in a vectorial fashion. To this purpose, we introduce the following matrices: 
\begin{equation} \label{matrices A A'}
A(\eta)=
\begin{bmatrix}
|\eta| & 0 & 0 & 0 \\
0 & \langle\eta\rangle & 0 & 0\\
0 & 0 & -|\eta| & 0 \\
0 & 0 & 0 & - \langle\eta\rangle
\end{bmatrix}, \quad
A'(V;\eta) :=
\begin{bmatrix}
0 & a_k\eta_1 & 0 & b_k\eta_1 \\
a_0\eta_1 & 0 & b_0\eta_1 & 0 \\
0 & a_k\eta_1 & 0 & b_k\eta_1\\
a_0\eta_1 & 0 & b_0\eta_1 & 0
\end{bmatrix},
\end{equation}
\begin{equation}\label{matrix A''(VI)}
A''(V^I;\eta) :=
\begin{bmatrix}
0 & 0 & 0 & 0 \\
a^I_0\eta_1 & 0 & b^I_0 \eta_1 & 0 \\
0 & 0 & 0 & 0\\
a^I_0\eta_1 & 0 & b^I_0 \eta_1 & 0
\end{bmatrix},
\end{equation}
\begin{equation} \label{matrices C' C''}
C'(W^I ;\eta):=
\begin{bmatrix}
0 & c^I_0 & 0 & d^I_0 \\
0 & e_0^I & 0 & f_0^I \\
0 & c^I_0 & 0 & d^I_0 \\
0 & e_0^I & 0 & f_0^I
\end{bmatrix}, \quad
C''(U ;\eta):=
\begin{bmatrix}
0 & 0 & 0 & 0 \\
0 & e_0 & 0 & f_0 \\
0 & 0 & 0 & 0 \\
0 & e_0 & 0 & f_0
\end{bmatrix}
\end{equation}
where
\begin{equation} \label{def ak, bk, a0, b0}
\begin{cases}
& a_k=a_k(v_{\pm}; \eta) := \frac{i}{4}\big[(v_+ + v_{-}) - \frac{D_x}{\langle D_x\rangle}(v_+ - v_{-})\cdot \frac{\eta}{\langle\eta\rangle}\big] \\
& b_k=b_k(v_{\pm}; \eta) := \frac{i}{4}\big[(v_+ + v_{-}) + \frac{D_x}{\langle D_x\rangle}(v_+ - v_{-})\cdot \frac{\eta}{\langle\eta\rangle}\big] \\
& a_0=a_0(v_{\pm}; \eta) := \frac{i}{4}\big[(v_+ + v_{-}) - \frac{D_x}{\langle D_x\rangle}(v_+ - v_{-})\cdot \frac{\eta}{|\eta|}\big] \\
& b_0=b_0(v_{\pm}; \eta) := \frac{i}{4}\big[(v_+ + v_{-}) + \frac{D_x}{\langle D_x\rangle}(v_+ - v_{-})\cdot \frac{\eta}{|\eta|}\big]
\end{cases}
\end{equation}
\begin{equation}  \label{def c0 d0 e0 f0}
\begin{cases}
& c_0=c_0(v_{\pm};\eta) := \frac{i}{4}\big[D_1(v_+ + v_{-}) - \frac{D_x D_1}{\langle D_x\rangle}(v_+ - v_{-})\cdot\frac{\eta}{\langle\eta\rangle}\big] \\
& d_0=d_0(v_{\pm};\eta) := \frac{i}{4}\big[D_1(v_+ + v_{-}) + \frac{D_x D_1}{\langle D_x\rangle}(v_+ - v_{-})\cdot\frac{\eta}{\langle\eta\rangle}\big] \\
& e_0=e_0(u_{\pm};\eta) := \frac{i}{4}\big[D_1(u_+ + u_{-}) - \frac{D_x D_1}{|D_x|}(u_+ - u_{-})\cdot \frac{\eta}{\langle\eta \rangle}\big] \\
& f_0=f_0(u_{\pm}; \eta) := \frac{i}{4}\big[D_1(u_+ + u_{-}) + \frac{D_x D_1}{|D_x|}(u_+ - u_{-})\cdot\frac{\eta}{\langle\eta\rangle}\big]
\end{cases}
\end{equation}
\begin{equation} \label{def_aI0 bI0 cI0 dI0 eI0 fI0}
\begin{gathered}
a^I_0 = a_0(v^I_\pm;\eta), \quad b^I_0 = b_0(v^I_\pm; \eta), 
c^I_0 = c_0(v^I_\pm; \eta), \quad d^I_0 = d_0(v^I_\pm;\eta),
\\
e^I_0 = e_0(u^I_\pm;\eta), \quad f_0^I(u^I_\pm;\eta),
\end{gathered}
\end{equation}
vectors $W, U, V$:\index{W@$W$, wave-Klein-Gordon vector}\index{V@$V$, Klein-Gordon vector}\index{U@$U$, wave vector}
\begin{equation} \label{def W,V, U}
W:=
\begin{bmatrix}
u_+ \\
v_+ \\
u_{-} \\
v_{-}
\end{bmatrix}, \quad
V:=
\begin{bmatrix}
0 \\
v_+ \\
0 \\
v_ {-}
\end{bmatrix}, \quad
U:=
\begin{bmatrix}
u_+ \\
0 \\
u_{-} \\
0
\end{bmatrix},
\end{equation}
along with $W^I$\index{WI@$W^I$, wave-Klein-Gordon vector with admissible vector fields} (resp. $V^I, U^I$)\index{VI@$V^I$, Klein-Gordon vector with admissible vector fields}\index{UI@$U^I$, wave vector with admissible vector fields} defined from $W$ (resp. $V,U$) by replacing $u_\pm, v_\pm$ with $u^I_\pm, v^I_\pm$; and finally
\begin{equation} \label{matrix QI}
Q^I_0(V,W) = 
\begin{bmatrix}
\sum_{\substack{(I_1,I_2)\in\mathcal{I}(I) \\ |I_2|<|I|}} c_{I_1, I_2} Q^{\mathrm{w}}_0(v^{I_1}_\pm, D v^{I_2}_\pm) \\
\sum_{\substack{(I_1,I_2)\in\mathcal{I}(I) \\ |I_1|, |I_2|<|I|}} c_{I_1, I_2} Q^{\mathrm{kg}}_0(v^{I_1}_\pm, D u^{I_2}_\pm) \\
\sum_{\substack{(I_1,I_2)\in\mathcal{I}(I) \\ |I_2|<|I|}} c_{I_1, I_2} Q^{\mathrm{w}}_0(v^{I_1}_\pm, D v^{I_2}_\pm) \\
\sum_{\substack{(I_1,I_2)\in\mathcal{I}(I) \\ |I_1|, |I_2|<|I|}} c_{I_1, I_2} Q^{\mathrm{kg}}_0(v^{I_1}_\pm, D u^{I_2}_\pm) 
\end{bmatrix}
\end{equation}
The quantization $Op^B$ (resp. $Op^B_R$) of a matrix $A=(a_{ij})_{1\le i,j\le n}$ is meant as a matrix of operators $Op^B(A) = (Op^B(a_{ij}))_{1\le i,j\le n}$ (resp. $Op^B_R(A) = (Op^B_R(a_{ij}))_{1\le i,j\le n}$), and for a vector $Y=[y_1, \dots, y_n]$, 
\[
Op^B(A)Y^\dagger =
\begin{bmatrix}
\displaystyle\sum_{j=1}^n Op^B(a_{1j})y_j \\
\vdots \\
\displaystyle\sum_{j=1}^n Op^B(a_{nj})y_j
\end{bmatrix},
\]
$Y^\dagger$ being the transpose of $Y$. We also remind that
\[\|A\|_{L^2}=\Big(\sum_{i,j}|a_{ij}|^2\Big)^{\frac{1}{2}}, \quad \|A\|_{L^\infty}=\sup_{ij}|a_{ij}|.\]

With notations introduced above, system \eqref{system for uI+-, vI+-} writes in the following compact fashion which has the merit to well highlight, among all non-linear terms, the very quasi-linear contributions $(QL)_1, (QL)_2$, represented below by $Op^B(A'(V;\eta))W^I$:
\begin{equation} \label{equation WI}
\begin{split}
D_t W^I & = A(D) W^I + Op^B(A'(V;\eta))W^I + Op^B(C'(W^I;\eta))V + Op^B_R(A'(V;\eta))W^I \\
& + Op^B(A''(V^I;\eta))U + Op^B(C''(U;\eta))V^I + Op^B_R(A''(V^I;\eta))U + Q^I_0(V,W).
\end{split}
\end{equation}
The energies defined in \eqref{def_generalized_energy} take the form 
\begin{subequations} \label{energy Ekm(t,W)}
\begin{gather}
E_n(t;u_\pm, v_\pm) = \sum_{|\alpha|\le n} \|D^\alpha_x W(t,\cdot)\|_{L^2}, \quad \forall\, n\in\mathbb{N}, n\ge 3, \\
E^k_3(t;u_\pm, v_\pm) = \sum_{\substack{ |\alpha|+|I|\le 3\\  |I|\le 3-k}}\|D^\alpha_x W^I(t,\cdot)\|^2_{L^2}, \quad\forall \, 0\le k \le 2,  \label{energy_Ek2}
\end{gather}
\end{subequations}
and we can refer to them, respectively, as $E_n(t;W), E^k_3(t;W)$. \index{EnW@$E_n(t;W)$, generalized energy} \index{EkW@$E^k_3(t;W)$, generalized energy}
We also notice that, since
\begin{subequations} \label{commutators_Z}
\begin{equation} \label{commutator_Z_Dt-|D|}
 [\Gamma, D_t \pm |D_x|]=
\begin{cases} 
 0 \quad  &\text{if } \Gamma\in \{\Omega, \partial_j, j=1,2\} ,\\
  \mp \frac{D_m}{|D_x|}(D_t \pm |D_x|) \quad &\text{if } \Gamma =Z_m, m=1,2,
 \end{cases}
\end{equation}
\begin{equation}\label{commutator_Z_Dt-<D>}
 [\Gamma, D_t \pm \langle D_x\rangle]=
\begin{cases} 
 0 \quad  &\text{if } \Gamma\in \{\Omega, \partial_j, j=1,2\} ,\\
  \mp \frac{D_m}{\langle D_x\rangle}(D_t \pm \langle D_x\rangle) \quad &\text{if } \Gamma =Z_m, m=1,2,
 \end{cases}
\end{equation}
\end{subequations}
and operators $D_m|D_x|^{-1}, D_m\langle D_x\rangle^{-1}$ are continuous on $L^2$ for $m=1,2$,
there exists a constant $C>0$ such that
\begin{equation} \label{equivalence GammaIW-Ekn}
C^{-1}\sum_{I\in\mathcal{I}^k_3}\|\Gamma^I W(t,\cdot)\|^2_{L^2}\le E^k_3(t;W)\le C\sum_{I\in\mathcal{I}^k_3}\|\Gamma^I W(t,\cdot)\|^2_{L^2},
\end{equation}
where, for any integer $0\le k \le 2$, \index{Ik3@$\mathcal{I}^k_3$, set of multi-indices}
\begin{equation} \label{set_Ik3}
\mathcal{I}^k_3:=\left\{|I|\le 3: \Gamma^I = D^\alpha_x \Gamma^J \text{ with } |\alpha|+|J|=|I|,  |J|\le 3-k\right\}.
\end{equation}
For convenience, we also introduce the following set: \index{In@$\mathcal{I}_n$, set of multi-indices}
\begin{equation} \label{set_In}
\mathcal{I}_n:=\left\{|I|\le n : \Gamma^I= D^\alpha_x \text{ with } |\alpha|=|I|\right\}, \quad n\in\mathbb{N}, n\ge 3.
\end{equation}

Matrices $A(\eta), A'(V;\eta), A''(V^I;\eta)$ are of order 1 and $A'(V;\eta), A''(V^I;\eta)$ are singular at $\eta = 0$ (i.e. some of their elements are singular at $\eta = 0$), while $C'(W^I;\eta), C''(U;\eta)$ are of order 0. 
Since we will need to do some symbolic calculus on $A'(V;\eta)$, we need to isolate the mentioned singularity. We hence define
\begin{equation} \label{matrices A'1 A'-1}
A'_1(V;\eta) := 
\begin{bmatrix}
0 & a_0\eta_1 & 0 & b_0\eta_1 \\
a_0\eta_1 & 0 & b_0\eta_1 & 0 \\
0 & a_0\eta_1 & 0 & b_0\eta_1\\
a_0\eta_1 & 0 & b_0\eta_1 & 0
\end{bmatrix}, \quad 
A'_{-1}(V;\eta) : = 
\begin{bmatrix}
0 & (a_k-a_0)\eta_1 & 0 & (b_k-b_0)\eta_1 \\
0 & 0 & 0 & 0 \\
0 & 0 & 0 &0\\
0 & (a_k-a_0)\eta_1 & 0 & (b_k-b_0)\eta_1 
\end{bmatrix},
\end{equation}
$A'_1(V;\eta)$ being a matrix of order 1, $A'_{-1}(V;\eta)$ of order $-1$, both singular at $\eta = 0$, and write $A'_1(V;\eta) = A'_1(V;\eta)(1-\chi)(\eta) + A'_1(V;\eta)\chi(\eta)$, where $\chi\in C^\infty_0(\mathbb{R}^2)$ is equal to 1 in the unit ball. 
Equation \eqref{equation WI} can be the rewritten as follows
\begin{equation} \label{equation WI-1}
\begin{split}
D_t W^I & = A(D) W^I + Op^B(A'_1(V;\eta)(1-\chi)(\eta))W^I + Op^B(A'_1(V;\eta)\chi(\eta))W^I \\ &
+ Op^B(A'_{-1}(V;\eta))W^I + Op^B(C'(W^I;\eta))V + Op^B_R(A'(V;\eta))W^I 
+ Op^B(A''(V^I;\eta))U \\&
+ Op^B(C''(U;\eta))V^I + Op^B_R(A''(V^I;\eta))U + Q^I_0(V,W),
\end{split}
\end{equation}
and the symbol $A'_1(V;\eta)(1-\chi)(\eta)$ associated to the quasi-linear contribution is no longer singular at $\eta=0$.
We observe that this matrix is real since $i(v_+ + v_{-}) = 2\partial_t v$, $i\frac{D_x}{\langle D_x\rangle}(v_+ - v_{-}) = 2\partial_x v$ and $v$ is a real solution, but it is not symmetric and such a lack of symmetry could lead to a loss of derivatives when writing an energy inequality for $W^I$. The issue is however only technical, in the sense that $A_1(V;\eta)(1-\chi)(\eta)$ can be replaced with a real, symmetric matrix, as explained in subsection \ref{Subs: Symmetrization} (see proposition \ref{Prop: equation of WIs}). 
Before proving such result, we need to derive some $L^2$ estimates for the semi-linear terms in the right hand side of \eqref{equation WI-1}.

\subsection{Estimates of quadratic terms} \label{Subs: Estimate of quadratic terms}

In this subsection we recover some estimates for the $L^2$ norm of the non-linear terms in the right hand side of equation \eqref{equation WI-1}.

\begin{lem}  \label{Lemma: L2 estimate of semilinear terms}
Let $I$ be a fixed multi-index and $\chi\in C^\infty_0(\mathbb{R}^2)$ equal to 1 in a neighbourhood of the origin. The following estimates hold:
\begin{subequations}
\begin{equation} \label{L2 est on Op(A'1 + A'0)WI}
\left\|\Big[Op^B\big(A'_1(V;\eta)\chi(\eta)\big) + Op^B\big(A'_{-1}(V;\eta)\big)\Big]W^I(t,\cdot)\right\|_{L^2} \lesssim \|V(t,\cdot)\|_{H^{1,\infty}} \|W^I(t,\cdot)\|_{L^2};
\end{equation}
\begin{equation} \label{L2 est on Op(C(WI,eta))V}
\left\|Op^B(C'(W^I;\eta))V(t,\cdot) \right\|_{L^2}\lesssim \|V(t,\cdot)\|_{H^{6,\infty}} \|W^I(t,\cdot)\|_{L^2} ;
\end{equation}
\begin{equation} \label{L2 est on OpBR(A')WI}
\| Op^B_R(A'(V;\eta))W^I(t,\cdot)\|_{L^2} \lesssim \|V(t,\cdot)\|_{H^{7,\infty}} \|W^I(t,\cdot)\|_{L^2};
\end{equation} 
\begin{multline} \label{L2 est on OpBR(A'')W}
\| Op^B(A''(V^I;\eta))U(t,\cdot)\|_{L^2} + \| Op^B_R(A''(V^I;\eta))U(t,\cdot)\|_{L^2}  \\
\lesssim \big( \| R_1U (t,\cdot)\|_{H^{6,\infty}} + \|U(t,\cdot)\|_{H^{6,\infty}}\big)\|V^I(t,\cdot)\|_{L^2} ;
\end{multline}
\begin{equation} \label{L2 est on OpB(C")VI}
\| Op^B (C''(U;\eta))V^I(t,\cdot)\|_{L^2} \lesssim  \big( \| R_1U (t,\cdot)\|_{H^{2,\infty}} + \|U(t,\cdot)\|_{H^{2,\infty}}\big)\|W^I(t,\cdot)\|_{L^2} ;
\end{equation}
\end{subequations}
\proof 

$\bullet$
Inequality \eqref{L2 est on Op(A'1 + A'0)WI} follows applying proposition \ref{Prop : Paradiff action on Sobolev spaces-NEW} to $Op^B\left(A'_{-1}(V;\eta)(1-\chi)(\eta)\right)W^I$ whose symbol $A'_{-1}(V;\eta)(1-\chi)(\eta)$ is of order $-1$ and has $M^{-1}_0$ seminorm bounded from above by $\|V(t,\cdot)\|_{H^{1,\infty}}$, after definitions \eqref{def: seminorm Mmo}, 
\eqref{def ak, bk, a0, b0} and \eqref{matrices A'1 A'-1}.

$\bullet$ Since from definition \eqref{matrices C' C''} of matrix $C'(W^I;\eta)$
\begin{multline*}
\left\|Op^B(C'(W^I;\eta))V \right\|_{L^2}\lesssim \left\|Op^B(D_1(v^I_+ + v^I_{-}))v_\pm\right\|_{L^2} + \left\|Op^B\Big(\frac{D_xD_1}{\langle D_x\rangle}(v^I_+ -v^I_{-})\cdot\frac{\eta}{\langle\eta\rangle}\Big)v_\pm\right\|_{L^2} \\
+ \left\|Op^B\big(D_1(u^I_+ + u^I_{-})\big)v_\pm\right\|_{L^2} +  \left\|Op^B\Big(\frac{D_xD_1}{| D_x|}(u^I_+ -u^I_{-})\cdot\frac{\eta}{\langle\eta\rangle}\Big)v_\pm\right\|_{L^2},
\end{multline*}
we reduce to prove inequality \eqref{L2 est on Op(C(WI,eta))V} for $Op^B\big(\frac{D_xD_1}{\langle D_x\rangle}(v^I_+ -v^I_{-})\cdot\frac{\eta}{\langle\eta\rangle}\big)v_+$, the same argument being applicable to all other $L^2$ norms appearing in the above right hand side.
Using equality \eqref{fourier transform of paradiff op}, and considering a new admissible cut-off function $\chi_1$ identically equal to 1 on the support of $\chi$, we first derive that
\begin{align*}
& \reallywidehat{Op^B\Big(\frac{D_x D_1}{\langle D_x \rangle}(v^I_+ + v^I_{-})\cdot\frac{\eta}{\langle\eta\rangle}\Big)v_+}(\xi)  = \frac{1}{(2\pi)^2}\int \chi\Big(\frac{\xi-\eta}{\langle\eta\rangle}\Big) \reallywidehat{\frac{D_x D_1}{\langle D_x \rangle}(v^I_+ + v^I_{-})}(\xi-\eta) \cdot \reallywidehat{\frac{D_x}{\langle D_x \rangle} v_+}(\eta) d\eta \\
& = \frac{1}{(2\pi)^2}\int \chi\Big(\frac{\xi-\eta}{\langle\eta\rangle}\Big)\Big(\frac{\xi_1 - \eta_1}{\langle\eta\rangle}\Big) \reallywidehat{\frac{D_x}{\langle D_x\rangle}(v^I_+ + v^I_{-})}(\xi-\eta) \cdot\reallywidehat{D_x v_+}(\eta) d\eta \\
& = \frac{1}{(2\pi)^2}\int \chi_1\Big(\frac{\xi-\eta}{\langle\eta\rangle}\Big)\reallywidehat{ \left[\chi\Big(\frac{D_x}{\langle\eta\rangle}\Big)\frac{D_1}{\langle\eta\rangle}\frac{D_x}{\langle D_x\rangle}(v^I_+ + v^I_{-})\right]}(\xi-\eta)\cdot \reallywidehat{D_x v_+}(\eta) d\eta \\
&=
\reallywidehat{Op^B\Big(\chi\Big(\frac{D_x}{\langle\eta\rangle}\Big)\frac{D_1}{\langle\eta\rangle}\frac{D_x}{\langle D_x\rangle}(
v^I_+ + v^I_{-})\Big) D_x v_+}(\xi).
\end{align*}
Successively, by decomposition \eqref{dev in paraproduct} and the fact that $R(u,v)$ is symmetric in $(u,v)$, we have that
\begin{multline*}
Op^B\Big(\chi\Big(\frac{D_x}{\langle\eta\rangle}\Big)\frac{D_1}{\langle\eta\rangle}\frac{D_x}{\langle D_x\rangle}(
v^I_+ + v^I_{-})\Big) D_x v_+ = \chi\Big(\frac{D_x}{\langle\eta\rangle}\Big)\frac{D_1}{\langle\eta\rangle}\frac{D_x}{\langle D_x \rangle}(v^I_+ + v^I_{-}) \cdot D_x v_+\\
 - \left[Op^B(D_x v_+) + Op^B_R(D_x v_+)\right]\left[\chi\Big(\frac{D_x}{\langle\eta\rangle}\Big)\frac{D_1}{\langle\eta\rangle}\frac{D_x}{\langle D_x\rangle}(v^I_+ + v^I_{-})\right],
\end{multline*}
so propositions \ref{Prop : Paradiff action on Sobolev spaces-NEW}, \ref{Prop: Paradiff action with non smooth symbols and R(u,v)} $(ii)$, and the fact that $\chi\Big(\frac{D_x}{\langle\eta\rangle}\Big)\frac{D_1}{\langle\eta\rangle}\frac{D_x}{\langle D_x\rangle}$ is an operator uniformly bounded on $L^2$, imply that
\[\left\| Op^B\Big(\frac{D_x D_1}{\langle D_x \rangle}(v^I_+ + v^I_{-})\cdot\frac{\eta}{\langle\eta\rangle}\Big)v_+\right\|_{L^2} \lesssim \|V(t,\cdot)\|_{H^{6,\infty}} \|V^I(t,\cdot)\|_{L^2}.\]

$\bullet$ By definition \eqref{matrices A A'} of $A'(V;\eta)$,
\begin{multline*}
\left\| Op^B_R\big(A'(V;\eta)\big)W^I(t,\cdot) \right\|_{L^2} \lesssim \left\| Op^B_R(v_++v_{-})v^I_\pm \right\|_{L^2} + \left\| Op^B_R\Big(\frac{D_x}{\langle D_x\rangle}(v_+-v_{-})\cdot\frac{\eta}{\langle\eta\rangle}\Big)v^I_\pm \right\|_{L^2}\\
+ \left\|Op^B_R(v_++v_{-})u^I_\pm \right\|_{L^2} + \left\| Op^B_R\Big(\frac{D_x}{\langle D_x\rangle}(v_+-v_{-})\cdot\frac{\eta}{|\eta|}\Big)u^I_\pm\right\|_{L^2}.
\end{multline*}
Let us only show that inequality \eqref{L2 est on OpBR(A')WI} holds for $Op^B\big(\frac{D_x}{\langle D_x\rangle}(v_+-v_{-})\cdot\frac{\eta \eta_1}{|\eta|}\big)u^I_+$.
For a smooth cut-off function $\phi$ equal to 1 in the unit ball we write
\begin{multline*}
Op^B_R\left(\frac{D_x}{\langle D_x\rangle}(v_+ - v_{-})\cdot \frac{\eta\eta_1}{|\eta|}\right)u^I_+ = Op^B_R\Big(\frac{D_x}{\langle D_x\rangle}(v_+ - v_{-})\cdot \frac{\eta\eta_1}{|\eta|} \phi(\eta)\Big)u^I_+ \\
+ Op^B_R\Big(\frac{D_x}{\langle D_x\rangle}(v_+ - v_{-})\cdot \frac{\eta\eta_1}{|\eta|}(1-\phi)(\eta)\Big)u^I_+,
\end{multline*}
where by proposition \ref{Prop: Paradiff action with non smooth symbols and R(u,v)} $(i)$
\begin{multline*}
\left\| Op^B_R\Big(\frac{D_x}{\langle D_x\rangle}(v_+ - v_{-})\cdot \frac{\eta\eta_1}{|\eta|} \phi(\eta)\Big)u^I_+\right\|_{L^2}\lesssim \left\| \frac{D_x}{\langle D_x\rangle}(v_+ - v_{-})(t,\cdot)\right\|_{L^\infty} \|u^I_+(t,\cdot)\|_{L^2}\\ \lesssim \|V(t,\cdot)\|_{H^{1,\infty}} \|W^I(t,\cdot)\|_{L^2}.
\end{multline*}
On the other hand
\begin{equation*}
 Op^B_R\Big(\frac{D_x}{\langle D_x\rangle}(v_+ - v_{-})\cdot \frac{\eta\eta_1}{|\eta|}(1-\phi)(\eta)\Big)u^I_+ = \int e^{ix\cdot\xi} m(\xi,\eta) \left[\langle D_x\rangle^7(\hat{v}_+ - \hat{v}_{-})(\xi-\eta)\right] \hat{u}^I_+(\eta) d\xi d\eta,
\end{equation*} 
where
\begin{equation*}
m(\xi,\eta):= \frac{1}{(2\pi)^2} \left(1 - \chi\left(\frac{\xi-\eta}{\langle\eta\rangle}\right) - \chi\left(\frac{\eta}{\langle\xi - \eta\rangle}\right)\right)(1-\phi)(\eta) \frac{\xi-\eta}{\langle\xi -\eta\rangle^8}\cdot \frac{\eta\eta_1}{|\eta|}
\end{equation*}
and frequencies $\xi-\eta$ and $\eta$ are either bounded or equivalent on the support of $m(\xi,\eta)$. Therefore $m(\xi,\eta)$ satisfies the hypothesis of lemma \ref{Lem_appendix: Kernel with 1 function} $(i)$ $|\partial^\alpha_\xi \partial^\beta_\eta m(\xi,\eta)|\lesssim \langle\xi\rangle^{-3}\langle\eta\rangle^{-3}$ for any $\alpha,\beta\in\mathbb{N}^2$, and by inequality \eqref{est: corollary_app L2 norm}
\begin{equation*}
\left\| Op^B_R\Big(\frac{D_x}{\langle D_x\rangle}(v_+ - v_{-})\cdot \frac{\eta\eta_1}{|\eta|}(1-\phi)(\eta)\Big)u^I_+\right\|_{L^2}\lesssim \|V(t,\cdot)\|_{H^{7,\infty}}\|W^I(t,\cdot)\|_{L^2}.
\end{equation*}

$\bullet$ From definition \eqref{matrix A''(VI)} of $A''(V;\eta)$,
\begin{equation*}
\left\| Op^B\big(A''(V;\eta)\big)U(t,\cdot) \right\|_{L^2}\lesssim \left\|Op^B\big((v^I_+ + v^I_{-})\eta_1\big)u_\pm \right\|_{L^2} + \left\|Op^B\Big(\frac{D_x}{\langle D_x\rangle}(v^I_+ -v^I_{-})\cdot\frac{\eta\eta_1}{|\eta|}\Big)u_\pm\right\|_{L^2},
\end{equation*}
(the same inequality holds evidently when $Op^B$ is replaced by $Op^B_R$). As done for previous cases, we reduce to show \eqref{L2 est on OpBR(A'')W} for $Op^B\big(\frac{D_x}{\langle D_x\rangle}(v^I_+ - v^I_{-} )\cdot \frac{\eta\eta_1}{|\eta|}\big)u_+$ (resp. for $Op^B$ replaced with $Op^B_R$).
Using decomposition \eqref{dev in paraproduct} and the fact that $R(u,v)$ is symmetric in $(u,v)$ we find that
\begin{multline*}
Op^B\Big(\frac{D_x}{\langle D_x\rangle}(v^I_+ - v^I_{-} )\cdot \frac{\eta\eta_1}{|\eta|}\Big)u_+  = \frac{D_x}{\langle D_x\rangle}(v^I_+ - v^I_{-}) \cdot\frac{D_xD_1}{|D_x|}u_+\\
- Op^B\Big(\frac{D_x D_1}{|D_x|}u_+\cdot\frac{\eta}{\langle\eta\rangle}\Big)(v^I_+ - v^I_{-}) 
 - Op^B_R\Big(\frac{D_x D_1}{|D_x|}u_+\cdot\frac{\eta}{\langle\eta\rangle}\Big)(v^I_+ - v^I_{-}),
\end{multline*}
and
\begin{equation*}
Op^B_R\Big(\frac{D_x}{\langle D_x\rangle}(v^I_+ - v^I_{-} )\cdot \frac{\eta\eta_1}{|\eta|}\Big)u_+  = Op^B_R\Big(\frac{D_x D_1}{|D_x|}u_+\cdot\frac{\eta}{\langle\eta\rangle}\Big)(v^I_+ - v^I_{-}),
\end{equation*}
so a direct application of propositions \ref{Prop : Paradiff action on Sobolev spaces-NEW} and \ref{Prop: Paradiff action with non smooth symbols and R(u,v)} $(ii)$ gives that the $L^2$ norm of the above right hand sides is bounded by $\left\| \frac{D_x D_1}{|D_x|}u_+\right\|_{H^{4,\infty}} \|V^I(t,\cdot)\|_{L^2}$, and hence by $\|R_1U(t,\cdot)\|_{H^{6,\infty}}\|V^I(t,\cdot)\|_{L^2}$, which gives inequality \eqref{L2 est on OpBR(A'')W}.

$\bullet$ From definition \eqref{matrices C' C''} of matrix $C''(U;\eta)$,
\begin{multline*}
\|Op^B(C''(U;\eta))V^I\|_{L^2} \lesssim \\ \left\|Op^B(D_1(u_+ + u_{-}))(v^I_+ + v^I_{-})\right\|_{L^2} + \left\| Op^B\left(\frac{D_x D_1}{|D_x|}(u_+ - u_{-})\cdot\frac{\eta}{\langle\eta\rangle}\right)(v^I_+ - v^I_{-})\right\|_{L^2},
\end{multline*} so estimate \eqref{L2 est on OpB(C")VI} follows immediately from proposition \ref{Prop : Paradiff action on Sobolev spaces-NEW}.
\endproof
\end{lem}

Lemmas \ref{Lem: L2 est nonlinearities} and \ref{Lem:L2 est nonlinearity Dt} below are introduced with the aim of deriving an estimate of the $L^2$ norm of vector $Q^I_0(V,W)$ defined in \eqref{matrix QI} (see corollary \ref{Cor: L2 est QI0(V,W)}).
We remind that the summations defining $Q^I_0(V,W)$ come from the action of family $\Gamma^I$ of admissible vector fields on the quadratic non-linearities $Q_0(v, \partial_1v)$ and  $Q_0(v, \partial_1u)$ in \eqref{wave KG system} (or, in terms of $u_\pm, v_\pm$, on $Q^\mathrm{w}_0(v_\pm, D_1v_\pm)$ and $ Q^\mathrm{kg}_0(v_\pm, D_1u_\pm)$).
According to remark \ref{Remark:Vector_field_on_null_structure}, if $I\in \mathcal{I}_n$ and $\Gamma^I$ is a product of spatial derivatives only the action of $\Gamma^I$ on $Q^\mathrm{w}_0(v_\pm, D_1v_\pm)$ (resp. on $Q^\mathrm{kg}_0(v_\pm, D_1 u_\pm)$) "distributes" entirely on its factors, meaning that
\begin{equation*}
\Gamma^I Q^\mathrm{w}_0(v_\pm, D_1v_\pm) = \sum_{\substack{(I_1,I_2)\in \mathcal{I}(I)\\ |I_1|+|I_2|=|I|}}Q^\mathrm{w}_0(v^{I_1}_\pm, D_1 v^{I_2}_\pm),
\end{equation*}
(the same for $\Gamma^I Q^\mathrm{kg}_0(v_\pm, D_1u_\pm)$), and all coefficients $c_{I_1,I_2}$ in the right hand side of \eqref{system for uI+-, vI+-} are equal to 0.
On the contrary, if $I\in\mathcal{I}^k_3$ for $0\le k\le 2$ and $\Gamma^I$ contains some Klainerman vector fields $\Omega, Z_m, m=1,2$, the commutation between $\Gamma^I$ and the null structure gives rise to new quadratic contributions in which the derivative $D_1$ is eventually replaced with $D_2, D_t$. 
As already seen in \eqref{Gamma_I_nonlinearity}, in this case we have
\begin{equation}\label{GammaI_Q0}
\Gamma^I Q^\mathrm{w}_0(v_\pm, D_1v_\pm) = \sum_{\substack{(I_1,I_2)\in \mathcal{I}(I)\\ |I_1|+|I_2|=|I|}}Q^\mathrm{w}_0(v^{I_1}_\pm, D_1 v^{I_2}_\pm) + \sum_{\substack{(I_1,I_2)\in \mathcal{I}(I)\\ |I_1|+|I_2|<|I|}}c_{I_1,I_2} Q^\mathrm{w}_0(v^{I_1}_\pm, D v^{I_2}_\pm),
\end{equation}
with some of the coefficients $c_{I_1,I_2}$ being equal to $1$ or $-1$, and $D\in\{D_1,D_2,D_t\}$ depending on the addend we are considering (similarly for $\Gamma^I Q^\mathrm{kg}_0(v_\pm, D_1u_\pm)$).
For our scopes, there will be no difference between the case $D=D_1$ and $D=D_2$, the two associated quadratic contributions enjoying the same $L^2$ and $L^\infty$ estimates.
When $D=D_t$, we should make use of the equation satisfied by $v^{I_2}_\pm$ (resp. by $u^{I_2}_\pm$) in system \eqref{system for uI+-, vI+-} to replace $Q^\mathrm{w}_0(v^{I_1}_\pm, D_tv^{I_2}_\pm)$ (resp. $Q^\mathrm{kg}_0(v^{I_1}_\pm, D_tu^{I_2}_\pm)$) with
\begin{equation}\label{Q(Dt)}
\begin{gathered}
Q^\mathrm{w}_0(v^{I_1}_\pm, \langle D_x\rangle v^{I_2}_\pm) + Q^\mathrm{w}_0\left(v^{I_1}_\pm, \Gamma^{I_2}
Q^\mathrm{kg}_0(v_\pm, D_1u_\pm)\right), \\
\left(\text{resp. with } Q^\mathrm{kg}_0(v^{I_1}_\pm, |D_x| u^{I_2}_\pm) + Q^\mathrm{kg}_0\left(v^{I_1}_\pm, \Gamma^{I_2} Q^\mathrm{w}_0(v_\pm, D_1v_\pm)\right)\right),
\end{gathered}
\end{equation}
where the left hand side quadratic terms are given by
\begin{equation} \label{Qw0-Qk0- |Dx|}
\begin{gathered}
Q^\mathrm{w}_0(v^{I_1}_\pm, \langle D_x\rangle v^{I_2}_\pm) = (v^{I_1}_+ + v^{I_1}_{-})\langle D_x\rangle(v^{I_2}_+ - v^{I_2}_{-}) - \frac{D_x}{\langle D_x\rangle}(v^{I_1}_+ - v^{I_1}_{-}) \cdot D_x(v^{I_2}_+ + v^{I_2}_{-}),\\
\left(\text{resp. } Q^\mathrm{kg}_0(v^{I_1}_\pm, |D_x| u^{I_2}_\pm) = (v^{I_1}_+ + v^{I_1}_{-})| D_x| (u^{I_2}_+ - u^{I_2}_{-}) - \frac{D_x}{\langle D_x\rangle}(v^{I_1}_+ - v^{I_1}_{-}) \cdot D_x(u^{I_2}_+ + u^{I_2}_{-})\right),
\end{gathered}
\end{equation}
while the right hand side ones in \eqref{Q(Dt)} are cubic.
On the Fourier side, these new quadratic contributions write as
\begin{equation*}
\begin{gathered}
\sum_{j_1,j_2\in \{+,-\}}\int j_2\left(1-j_1j_2\frac{\xi-\eta}{\langle \xi-\eta\rangle}\cdot\frac{\eta}{\langle\eta\rangle}\right)\langle\eta\rangle \hat{v}^{I_1}_{j_1}(\xi-\eta)\hat{v}^{I_2}_{j_2}(\eta) d\xi d\eta, \\
\left(\text{resp. } \sum_{j_1,j_2\in \{+,-\}}\int j_2\left(1-j_1j_2\frac{\xi-\eta}{\langle \xi-\eta\rangle}\cdot\frac{\eta}{|\eta|}\right)|\eta| \hat{v}^{I_1}_{j_1}(\xi-\eta)\hat{u}^{I_2}_{j_2}(\eta) d\xi d\eta\right),
\end{gathered}
\end{equation*}
and have basically the same nature of the starting ones, as
\begin{equation*}
\begin{gathered}
\reallywidehat{Q^\mathrm{w}_0(v^{I_1}_\pm, D_1 v^{I_2}_\pm)}(\xi)= \sum_{j_1,j_2\in \{+,-\}}\int \left(1-j_1j_2\frac{\xi-\eta}{\langle \xi-\eta\rangle}\cdot\frac{\eta}{\langle\eta\rangle}\right)\eta_1 \hat{v}^{I_1}_{j_1}(\xi-\eta)\hat{v}^{I_2}_{j_2}(\eta) d\xi d\eta, \\
\left(\text{resp. } \reallywidehat{Q^\mathrm{kg}_0(v^{I_1}_\pm, D_1 u^{I_2}_\pm)}(\xi)=\sum_{j_1,j_2\in \{+,-\}}\int \left(1-j_1j_2\frac{\xi-\eta}{\langle \xi-\eta\rangle}\cdot\frac{\eta}{|\eta|}\right)\eta_1 \hat{v}^{I_1}_{j_1}(\xi-\eta)\hat{u}^{I_2}_{j_2}(\eta) d\xi d\eta\right).
\end{gathered}
\end{equation*}
For this reason, as long as we can neglect the cubic terms in \eqref{Q(Dt)}, we will not pay attention to the value of $D\in \{D_1,D_2,D_t\}$ in the second sum in the right hand side of \eqref{GammaI_Q0}.
Lemma \ref{Lem:L2 est nonlinearity Dt} is meant to show that the mentioned cubic terms are, indeed, remainders.

Before proving lemmas \ref{Lem: L2 est nonlinearities}, \ref{Lem:L2 est nonlinearity Dt}, we need to introduce a new set of indices. According to the order established in $\mathcal{Z}$ at the beginning of section \ref{sec: statement of the main results} (see \eqref{order_Z}), we define \index{K@$\mathcal{K}$, set of multi-indices}
\begin{equation} \label{set_K}
\mathcal{K}:=\{I=(i_1,i_2) : i_1,i_2 =1,2, 3\}
\end{equation}
as the set of indices $I$ such that $\Gamma^I$ is the product of two Klainerman vector fields only, together with \index{Vk@$\mathcal{V}^k$, set of multi-indices}
\begin{equation}\label{set_V}
\mathcal{V}^k:=\{I\in \mathcal{I}^k_3 : \exists (I_1,I_2)\in\mathcal{I}(I) \text{ with } I_1\in\mathcal{K} \},
\end{equation}
which is evidently empty when $k=2$.
We also warn the reader that, in inequality \eqref{est:L2_norm_Rk3(t,x)} with $k=2$, $E^3_3(t;W)$ stands for $E_3(t;W)$, this double notation allowing us to combine in one line all cases $k=0,1,2$.

\begin{lem} \label{Lem: L2 est nonlinearities}
$(i)$ Let $n\in\mathbb{N}, n\ge 3$ and $I\in \mathcal{I}_n$.
Then
\begin{equation} \label{est: L2 Qw0 (vI1 vI2)-only derivatives}
\sum_{\substack{(I_1,I_2)\in\mathcal{I}(I)\\ |I_2|<n}}\left\|Q^\mathrm{w}_0(v^{I_1}_\pm, D_x v^{I_2}_\pm)\right\|_{L^2}+\sum_{\substack{(I_1,I_2)\in\mathcal{I}(I)\\ |I_1|\le [\frac{n}{2}], |I_2|<n}}\left\|Q^\mathrm{kg}_0(v^{I_1}_\pm, D_xu^{I_2}_\pm)\right\|_{L^2} \lesssim \|V(t,\cdot)\|_{H^{[\frac{n}{2}]+2, \infty}}E_n(t;W)^\frac{1}{2}, 
\end{equation}
\begin{equation} \label{est: L2 Qkg0 (vI1 vI2)-only derivatives}
\sum_{\substack{(I_1,I_2)\in \mathcal{I}(I)\\ |I_1|> [\frac{n}{2}]}}\left\|Q^\mathrm{kg}_0(v^{I_1}_\pm, D_xu^{I_2}_\pm)\right\|_{L^2}
\lesssim \left(\|U(t,\cdot)\|_{H^{[\frac{n}{2}]+2, \infty}}+ \|\mathrm{R}_1 U(t,\cdot)\|_{H^{[\frac{n}{2}]+2,\infty}}\right)E_n(t;W)^\frac{1}{2}.
\end{equation}

$(ii)$ Let $0\le k \le 2$ and $I\in \mathcal{I}^k_3$.
There exists a constant $C>0$ such that, if we assume a-priori estimates \eqref{est: bootstrap upm}, \eqref{est: boostrap vpm} satisfied and $0<\varepsilon_0<(2A+B)^{-1}$ small, for any $\chi\in C^\infty_0(\mathbb{R}^2)$ equal to 1 in a neighbourhood of the origin and $\sigma>0$ small we have
\begin{subequations} \label{decomposition_lemma_Qwo-Qkg0}
\begin{gather}
\sum_{\substack{(I_1,I_2)\in\mathcal{I}(I)\\ |I_2|<3}} Q^\mathrm{w}_0(v^{I_1}_\pm, D_x v^{I_2}_\pm) = \mathfrak{R}^k_3(t,x),\label{decomposition_Qw0} \\
\sum_{\substack{(I_1,I_2)\in\mathcal{I}(I)\\|I_1|, |I_2|<3}} Q^\mathrm{kg}_0(v^{I_1}_\pm, D_x u^{I_2}_\pm) = \delta_{\mathcal{V}^k}\sum_{\substack{(I_1,I_2)\in\mathcal{I}(I)\\ I_1\in\mathcal{K},  |I_2|\le 1}} Q^\mathrm{kg}_0\left(v^{I_1}_\pm, \chi(t^{-\sigma}D_x)  D_x u^{I_2}_\pm\right) +\mathfrak{R}^k_3(t,x), \label{decomposition_Qkg0}
\end{gather}
\end{subequations}
where $\delta_{\mathcal{V}^k}=1$ if $I\in\mathcal{V}^k$, 0 otherwise, and
\begin{equation} \label{est:L2_norm_Rk3(t,x)}
\|\mathfrak{R}^k_3(t,\cdot)\|_{L^2}\le C(A+B)\varepsilon t^{-1}  E^k_3(t,W)^\frac{1}{2} + CB\varepsilon t^{-\frac{5}{4}},
\end{equation}
with $\beta>0$ small, $\beta\rightarrow 0$ as $\sigma\rightarrow 0$, for all $t\in[1,T]$.
The same result holds with $D_x v^{I_2}_\pm$ (resp. $D_x u^{I_2}_\pm$) replaced with $\langle D_x\rangle v^{I_2}_\pm$ (resp. $|D_x| u^{I_2}_\pm$).
\proof
$(i)$ The proof of follows straightly from \eqref{Q0_pm} with $a=1,2$, by bounding the $L^2$ norm of each product with the $L^\infty$ norm of the factor indexed in $J\in \{I_1,I_2\}$ such that $|J|\le \big[\frac{|I|}{2}\big]$, times the $L^2$ norm of the remaining one.

$(ii)$ Let $I\in\mathcal{I}^k_3$. One immediately sees that:
\begin{multline} \label{ineq:I1,0- 0,I2}
\sum_{(J,0)\in\mathcal{I}(I)} \|Q^\mathrm{w}_0(v^J_\pm, D_x v_\pm)\|_{L^2} + \sum_{\substack{(J,0)\in\mathcal{I}(I)\\ |J|<3 }}\Big(\|Q^\mathrm{w}_0(v_\pm, D_x v^J_\pm)\|_{L^2}+ \|Q^\mathrm{kg}_0(v_\pm, D_x u^J_\pm)\|_{L^2}\Big) \\
\lesssim \|V(t,\cdot)\|_{H^{2,\infty}}E^k_3(t;W)^\frac{1}{2};
\end{multline}
if $(I_1,I_2)\in\mathcal{I}(I)$ is such that $|I_2|<3$ and either $\Gamma^{I_1}$ or $\Gamma^{I_2}$ is a product of spatial derivatives only
\begin{equation} \label{est:quadratic-terms-derivative1}
\|Q^\mathrm{w}_0(v^{I_1}_\pm, D_x v^{I_2}_\pm)\|_{L^2}\lesssim \|V(t,\cdot)\|_{H^{4,\infty}}E^k_3(t;W)^\frac{1}{2};
\end{equation}
if $(I_1,I_2)\in\mathcal{I}(I)$ is such that $|I_2|<3$ and $\Gamma^{I_1}$ is a product of spatial derivatives only
\begin{equation} \label{est:quadratic-term-derivative2}
\|Q^\mathrm{kg}_0(v^{I_1}_\pm, D_x u^{I_2}_\pm)\|_{L^2}\lesssim \|V(t,\cdot)\|_{H^{3,\infty}}E^k_3(t;W)^\frac{1}{2}.
\end{equation}
Hence, the remaining quadratic contributions to be estimated are those corresponding to indices $(I_1,I_2)\in\mathcal{I}(I)$, with $|I_2|<3$, such that:
both $\Gamma^{I_1}$ and $\Gamma^{I_2}$ contain at least one Klainerman vector field, in the left hand side of \eqref{decomposition_Qw0}; $\Gamma^{I_1}$ contains one or two Klainerman vector fields, in the left hand side of \eqref{decomposition_Qkg0}.

The idea to estimate the $L^2$ norm of the $Q^\mathrm{w}_0(v^{I_1}_\pm, D v^{I_2}_\pm)$, for indices $I_1,I_2$ just mentioned above, is to decompose the Klein-Gordon component carrying exactly one Klainerman vector field in frequencies, by means of a truncation $\chi(t^{-\sigma}D_x)$ for some smooth cut-off function $\chi$ and $\sigma>0$ small.
Basically, the $L^\infty$ norm of the contribution truncated for large frequencies $|\xi|\gtrsim t^\sigma$ can be bounded by making appear a power of $t$ as negative as we want, while that of the remaining one, localized for $|\xi|\lesssim t^{\sigma}$, enjoys the sharp Klein-Gordon decay $t^{-1}$ as proved in lemma \ref{Lem_appendix: sharp_est_VJ} in appendix \ref{Appendix B}. The same argument can be applied to $Q^\mathrm{kg}_0(v^{I_1}_\pm, D_xu^{I_2}_\pm)$ with $I_1$ such that $\Gamma^{I_1}$ contains exactly one Klainerman vector field.
Then, by lemma \ref{Lem_app:products_Gamma} in appendix \ref{Appendix B} with $L=L^2$
we find that, for some $\chi\in C^\infty_0(\mathbb{R}^2)$, the following:
 if $\Gamma^{I_1}$ contains exactly one Klainerman vector field,
\begin{multline*}
\left\|Q^\mathrm{w}_0(v^{I_1}_\pm, D_x v^{I_2}_\pm)(t,\cdot)\right\|_{L^2}\lesssim \left\|\chi(t^{-\sigma}D_x)v^{I_1}_\pm(t,\cdot)\right\|_{H^{1,\infty}}\|v^{I_2}_\pm(t,\cdot)\|_{H^1}\\
+ t^{-N(s)}\left(\|v_\pm(t,\cdot)\|_{H^s}+\|D_t v_\pm(t,\cdot)\|_{H^s}\right)\Big(\sum_{|\mu|=0}^1\|x^\mu v^{I_2}_\pm(t,\cdot)\|_{H^1}+t\|v^{I_2}_\pm(t,\cdot)\|_{H^1}\Big)
\end{multline*}
and
\begin{multline*}
\left\|Q^\mathrm{kg}_0(v^{I_1}_\pm, D_xu^{I_2}_\pm)(t,\cdot)\right\|_{L^2}\lesssim \left\|\chi(t^{-\sigma}D_x)v^{I_1}_\pm(t,\cdot)\right\|_{H^{1,\infty}}\|u^{I_2}_\pm(t,\cdot)\|_{H^1}\\
+ t^{-N(s)}\left(\|v_\pm(t,\cdot)\|_{H^s}+\|D_t v_\pm(t,\cdot)\|_{H^s}\right)\Big(\sum_{|\mu|=0}^1\|x^\mu D_xu^{I_2}_\pm(t,\cdot)\|_{L^2}+t\|u^{I_2}_\pm(t,\cdot)\|_{H^1}\Big);
\end{multline*}
if $\Gamma^{I_2}$ contains exactly one Klainerman vector field,
\begin{multline*}
\left\|Q^\mathrm{w}_0(v^{I_1}_\pm, D_xv^{I_2}_\pm)(t,\cdot)\right\|_{L^2}\lesssim \left\|\chi(t^{-\sigma}D_x)v^{I_2}_\pm(t,\cdot)\right\|_{H^{2,\infty}}\|v^{I_1}_\pm(t,\cdot)\|_{L^2}\\
+ t^{-N(s)}\left(\|v_\pm(t,\cdot)\|_{H^s}+\|D_t v_\pm(t,\cdot)\|_{H^s}\right)\Big(\sum_{|\mu|=0}^1\|x^\mu v^{I_1}_\pm(t,\cdot)\|_{L^2}+t\|v^{I_1}_\pm(t,\cdot)\|_{L^2}\Big),
\end{multline*}
where, in all above inequalities, $N(s)\ge 3$ if $s>0$ is large enough.
From inequalities \eqref{Hs_norm_DtU}, \eqref{Hs norm DtV}, estimates \eqref{norm_L2_xj-GammaIv-}, lemma \ref{Lem_appendix: sharp_est_VJ} and the boostrap assumptions \eqref{est: bootstrap argument a-priori est}, together with the fact $\delta,\delta_j\ll 1$ are small, for $j=0,1,2$, we derive that there is a positive constant $C$ such that, for multi-indices $I_1,I_2$ considered in above inequalities,
\begin{equation*}
\left\|Q^\mathrm{w}_0(v^{I_1}_\pm, Dv^{I_2}_\pm)(t,\cdot)\right\|_{L^2} + \left\|Q^\mathrm{kg}_0(v^{I_1}_\pm, Du^{I_2}_\pm)(t,\cdot)\right\|_{L^2}\le CB\varepsilon t^{-1}E^k_3(t;W)^\frac{1}{2} + CB\varepsilon t^{-\frac{5}{4}}.
\end{equation*}

The remaining quadratic terms are $Q^\mathrm{kg}_0(v^{I_1}_\pm, D_x u^{I_2}_\pm)$ with $I_1\in \mathcal{K}$ (and hence $|I_2|\le 1$) if $\mathcal{V}^k$ is non empty. Applying lemma \ref{Lem_app:products_Gamma} with $L=L^2$, $w=u$ and the same $s$ as before, and making use of estimates \eqref{est: bootstrap argument a-priori est}, \eqref{norm_L2_xj-GammaIv-}, together with inequality \eqref{Hs_norm_DtU}, we see that
\begin{multline*}
\left\|Q^\mathrm{kg}_0(v^{I_1}_\pm, D_x u^{I_2}_\pm)(t,\cdot)\right\|_{L^2}\lesssim \left\|Q^\mathrm{kg}_0\left(v^{I_1}_\pm, \chi(t^{-\sigma}D_x) D_xu^{I_2}_\pm\right)(t,\cdot)\right\|_{L^2} \\
+ t^{-3}\Big(\sum_{|\mu|=0}^1 \|x^\mu v^{I_1}_\pm(t,\cdot)\|_{L^2}+t\|v^{I_1}_\pm(t,\cdot)\|_{L^2}\Big)\left(\|u_\pm(t,\cdot)\|_{H^s}+\|D_tu_\pm(t,\cdot)\|_{H^s}\right)\\
\lesssim \left\|Q^\mathrm{kg}_0\left(v^{I_1}_\pm, \chi(t^{-\sigma}D_x) D_xu^{I_2}_\pm\right)(t,\cdot)\right\|_{L^2} + CB\varepsilon t^{-\frac{5}{4}},
\end{multline*}
which hence concludes the proof of $(ii)$.
We should highlight the fact that the quadratic contribution in the above left hand side is treated differently from the previous ones, because we do not have a sharp decay $O(t^{-1})$ for $v^{I_1}_\pm$ when $I_1\in\mathcal{K}$ (neither when truncated for moderate frequencies), but only a control in $O(t^{-1+\beta'})$, for some small $\beta'>0$ (see lemma \ref{Lem_appendix:est vI I=2}). Moreover, the decay enjoyed by the uniform norm of $\chi(t^{-\sigma}D_x)D_x u^{I_2}_\pm$, appearing in the quadratic term in the above right hand side, is very weak (only $t^{-1/2+\beta'}$, see lemma \ref{Lem_appendix: est UJ}). Such terms, that contribute to the energy and decay slowly in time, will be successively eliminated by a normal form argument (see subsection \ref{sub: second normal form}).
\endproof
\end{lem}

\begin{lem} \label{Lem:L2 est nonlinearity Dt}
Let $0\le k \le 2$ and $I\in\mathcal{I}^k_3$. For any $\chi\in C^\infty_0(\mathbb{R}^2)$ equal to 1 in a neighbourhood of the origin and $\sigma>0$ small
\begin{subequations}
\begin{gather}
\sum_{\substack{(I_1,I_2)\in\mathcal{I}(I)\\ |I_1|+|I_2|\le 2}}Q^\mathrm{w}_0(v^{I_1}_\pm, D_tv^{I_2}_\pm) = \mathfrak{R}^k_3(t,x),\label{eq:Qw0(Dt)-statement} \\
\sum_{\substack{(I_1,I_2)\in\mathcal{I}(I)\\ |I_1|+|I_2|\le 2}}Q^\mathrm{kg}_0(v^{I_1}_\pm, D_tu^{I_2}_\pm) =\delta_{\mathcal{V}^k}\sum_{\substack{(J,0)\in\mathcal{I}(I)\\ J\in\mathcal{K}}}Q^\mathrm{kg}_0(v^J_\pm,  \chi(t^{-\sigma}D_x) |D_x| u_\pm)+ \mathfrak{R}^k_3(t,x), \label{eq:Qkg0(Dt)-statement}
\end{gather}
\end{subequations}
with $\delta_{\mathcal{V}^k}=1$ if $I\in\mathcal{V}^k$, 0 otherwise, and $\mathfrak{R}^k_3(t,x)$ satisfying \eqref{est:L2_norm_Rk3(t,x)}.
\proof
Using the equation satisfied by $v^{I_2}_\pm$ and $u^{I_2}_\pm$ respectively in system \eqref{system for uI+-, vI+-} with $I=I_2$ we see that
\begin{subequations}\label{eq:non-lin-Dt}
\begin{multline} \label{eq:Qw0(Dt)}
\sum_{\substack{(I_1,I_2)\in\mathcal{I}(I)\\ |I_1|+|I_2|\le 2}}Q^\mathrm{w}_0(v^{I_1}_\pm, D_tv^{I_2}_\pm) = \sum_{\substack{(I_1,I_2)\in\mathcal{I}(I)\\ |I_1|+|I_2|\le 2}}Q^\mathrm{w}_0(v^{I_1}_\pm, \langle D_x\rangle v^{I_2}_\pm)\\
 +\sum_{\substack{(I_1,I_2)\in\mathcal{I}(I)\\ |I_1|+|I_2|\le 2}} \sum_{(J_1,J_2)\in\mathcal{I}(I_2)}c_{J_1,J_2} Q^\mathrm{w}_0\left(v^{I_1}_\pm, Q^\mathrm{kg}_0(v^{J_1}_\pm, D u^{J_2}_\pm)\right),
\end{multline}
\begin{multline}\label{eq:Qkg0(Dt)}
\sum_{\substack{(I_1,I_2)\in\mathcal{I}(I)\\ |I_1|+|I_2|\le 2}}Q^\mathrm{kg}_0(v^{I_1}_\pm, D_tu^{I_2}_\pm) = \sum_{\substack{(I_1,I_2)\in\mathcal{I}(I)\\ |I_1|+|I_2|\le 2}}Q^\mathrm{kg}_0(v^{I_1}_\pm, | D_x| u^{I_2}_\pm)\\
 +\sum_{\substack{(I_1,I_2)\in\mathcal{I}(I)\\ |I_1|+|I_2|\le 2}} \sum_{(J_1,J_2)\in\mathcal{I}(I_2)}c_{J_1,J_2} Q^\mathrm{kg}_0\left(v^{I_1}_\pm, Q^\mathrm{w}_0(v^{J_1}_\pm, D v^{J_2}_\pm)\right),
\end{multline}
\end{subequations}
with coefficients $c_{J_1,J_2}\in \{-1,0,-1\}$ such that $c_{J_1,J_2}=1$ whenever $|J_1|+|J_2|=|I_2|$, and $Q^\mathrm{w}_0(v^{I_1}_\pm, \langle D_x\rangle v^{I_2}_\pm)$ (in which case $D=D_1$), $Q^\mathrm{kg}_0(v^{I_1}_\pm, |D_x| u^{I_2}_\pm)$ given explicitly by \eqref{Qw0-Qk0- |Dx|}.
After lemma \ref{Lem: L2 est nonlinearities} $(ii)$ we know that
\begin{equation*}
\sum_{\substack{(I_1,I_2)\in\mathcal{I}(I)\\ |I_1|+|I_2|\le 2}}\left[Q^\mathrm{w}_0(v^{I_1}_\pm, \langle D_x\rangle v^{I_2}_\pm) + Q^\mathrm{kg}_0(v^{I_1}_\pm, | D_x| u^{I_2}_\pm)\right]= \sum_{\substack{(J,0)\in \mathcal{I}(I)\\ J\in\mathcal{K}}}Q^\mathrm{kg}_0(v^J_\pm, | D_x| u_\pm) + \mathfrak{R}^k_3(t,x),
\end{equation*}
with $\mathfrak{R}^k_3$ verifying \eqref{est:L2_norm_Rk3(t,x)}.
The only thing to prove is that
the cubic terms in the right hand side of \eqref{eq:non-lin-Dt} are remainders $\mathfrak{R}^k_3$. 
We focus on those in the right hand side of \eqref{eq:Qw0(Dt)} as the same argument applies to the ones in \eqref{eq:Qkg0(Dt)}.

First, let us consider cubic terms corresponding to indices $I_1,I_2$ such that $|I_1|=2$ and $|I_2|=0$.
In this case we evidently have that $|J_1|=|J_2|=0$, and by \eqref{est Hsinfty NLkg-new} with $s=1$ and $\theta\ll 1$ small, together with a-priori estimate \eqref{est: bootstrap argument a-priori est},
\begin{equation*}
\left\|Q^\mathrm{w}_0\left(v^{I_1}_\pm , Q^\mathrm{kg}_0(v_\pm, D_1 u_\pm)\right)\right\|_{L^2}\lesssim \|v^{I_1}_\pm(t,\cdot) \|_{L^2} \|Q^\mathrm{kg}_0(v_\pm, D_1 u_\pm)\|_{H^{1,\infty}}\le CB\varepsilon t^{-\frac{3}{2}+\beta'},
\end{equation*}
for some $\beta'>0$ small as long as $\sigma, \delta_0$ are small.

Let us now consider indices $I_1,I_2$ such that $\Gamma^{I_1}\in\{\Omega, Z_m, m=1,2\}$. As we also require that $(I_1,I_2)\in\mathcal{I}(I)$ with $|I_2|\le 2$, we have in this case that $|I_2|\le 1$ and consequently, for each $(J_1,J_2)\in\mathcal{I}(I_2)$, either $|J_1|=0$ or $|J_2|=0$.
Using lemma \ref{Lem_app:products_Gamma} in appendix \ref{Appendix B} with $L=L^2$ and $w=v$, we derive that for any $\chi\in C^\infty_0(\mathbb{R}^2)$ as in the statement and $\sigma>0$ small
\begin{align*}
\sum_{(J_1,J_2)\in\mathcal{I}(I_2)}&\left\|Q^\mathrm{w}_0\left(v^{I_1}_\pm, Q^\mathrm{kg}_0(v^{J_1}_\pm, D u^{J_2}_\pm)\right)(t,\cdot)\right\|_{L^2}\\
& \lesssim \sum_{(J_1,J_2)\in\mathcal{I}(I_2)}\left\|\chi(t^{-\sigma}D_x)v^{I_1}_\pm(t,\cdot)\right\|_{L^\infty}\left\|Q^\mathrm{kg}_0(v^{J_1}_\pm, D u^{J_2}_\pm)(t,\cdot)\right\|_{L^2}\\
&+\sum_{(J_1,J_2)\in\mathcal{I}(I_2)}t^{-N(s)}\left(\|v_\pm(t,\cdot)\|_{H^s} + \|D_tv_\pm(t,\cdot)\|_{H^s}\right)\\
&\hspace{2cm}\times\Big(\sum_{|\mu|=0}^1\left\|x^\mu Q^\mathrm{kg}_0(v^{J_1}_\pm, D u^{J_2}_\pm)(t,\cdot)\right\|_{L^2}+ t \left\|Q^\mathrm{kg}_0(v^{J_1}_\pm, D u^{J_2}_\pm)(t,\cdot)\right\|_{L^2}\Big),
\end{align*}
with $N(s)\ge 3$ is $s>0$ is sufficiently large. Here
\begin{align*}
& \sum_{(J_1,J_2)\in\mathcal{I}(I_2)} \left\|x Q^\mathrm{kg}_0(v^{J_1}_\pm, D u^{J_2}_\pm)(t,\cdot)\right\|_{L^2}
\\
&\hspace{0.5cm}\lesssim \sum_{\substack{|\mu|=0,1 \\ |J|\le 1}}\left[\left\|x\Big(\frac{D_x}{\langle D_x\rangle}\Big)^\mu v_\pm(t,\cdot)\right\|_{L^\infty}\left(\|u^J_\pm(t,\cdot)\|_{H^1} + \|D_tu^J_\pm(t,\cdot)\|_{L^2}\right)\right. \\
&\hspace{1cm} \left. + \|xv^J_\pm(t,\cdot)\|_{L^2}\left(\|\mathrm{R}^\mu u_\pm(t,\cdot)\|_{H^{2,\infty}}+ \|D_t\mathrm{R}^\mu u_\pm(t,\cdot)\|_{H^{1,\infty}}\right)\right] \le C(A+B)B\varepsilon^2 t^{\frac{1}{2}+\frac{\delta_2}{2}}
\end{align*}
and 
\begin{multline}\label{est_J1J2_less 2}
\sum_{(J_1,J_2)\in\mathcal{I}(I_2)}\left\|Q^\mathrm{kg}_0(v^{J_1}_\pm, D u^{J_2}_\pm)(t,\cdot)\right\|_{L^2}\\
 \lesssim \sum_{|J|\le 1}\Big[ \|v^J_\pm(t,\cdot)\|_{L^2}\Big(\sum_{|\mu|=0}^1\|\mathrm{R}^\mu u_\pm(t,\cdot)\|_{H^{2,\infty}}+ \|D_t\mathrm{R}^\mu u_\pm(t,\cdot)\|_{H^{1,\infty}}\Big) \\+ \|v_\pm(t,\cdot)\|_{H^{1,\infty}}\big(\|u^J_\pm(t,\cdot)\|_{H^1}+ \|D_tu^J_\pm(t,\cdot)\|_{L^2}\big)\Big]\le C(A+B)B\varepsilon^2 t^{-\frac{1}{2}+\frac{\delta_2}{2}}
\end{multline}
by \eqref{Hs_norm_DtU}, \eqref{Hsinfty_norm_DtU}, \eqref{Hsinfty norm Dt R1U-new}, \eqref{DtUI} and estimates \eqref{est: bootstrap argument a-priori est}, \eqref{norms_H1_Linfty_xv-}, \eqref{norm_L2_xj-GammaIv-}, so together with lemma \ref{Lem_appendix: sharp_est_VJ} and \eqref{Hs norm DtV}, these inequalities give
\[\sum_{(J_1,J_2)\in\mathcal{I}(I_2)}\left\|Q^\mathrm{w}_0\left(v^{I_1}_\pm, Q^\mathrm{kg}_0(v^{J_1}_\pm, D u^{J_2}_\pm)\right)(t,\cdot)\right\|_{L^2} \le CB\varepsilon t^{-\frac{3}{2}+\beta'},\]
for some new $\beta'>0$ small, $\beta'\rightarrow 0$ as $\sigma,\delta_0\rightarrow 0$.

Finally, for indices $I_1,I_2$ such that $\Gamma^{I_1}\in \{D^\alpha_x , |\alpha|\le 1\}$
\begin{equation}\label{last_J1J2}
\sum_{(J_1,J_2)\in\mathcal{I}(I_2)}\left\|Q^\mathrm{w}_0\left(v^{I_1}_\pm, Q^\mathrm{kg}_0(v^{J_1}_\pm, D u^{J_2}_\pm)\right)\right\|_{L^2} \lesssim \sum_{(J_1,J_2)\in\mathcal{I}(I_2)}\|v_\pm(t,\cdot)\|_{H^{2,\infty}}\left\|Q^\mathrm{kg}_0(v^{J_1}_\pm, D u^{J_2}_\pm)\right\|_{L^2}.
\end{equation}
For $(J_1,J_2)\in\mathcal{I}(I_2)$ such that $|J_1|+|J_2|=|I_2|$ we have by lemma \ref{Lem: L2 est nonlinearities} $(ii)$ and a-priori estimates \eqref{est: bootstrap argument a-priori est} that
\begin{equation*}
\begin{split}
\| Q^\mathrm{kg}_0(v^{J_1}_\pm, D_1 u^{J_2}_\pm)\|_{L^2} &\lesssim \|\mathfrak{R}^k_3(t,\cdot)\|_{L^2} + \sum_{J\in\mathcal{K}}\|Q^\mathrm{kg}_0(v^J_\pm, D_1 \chi(t^{-\sigma}D_x)u_\pm)\|_{L^2}\\
&\lesssim \|\mathfrak{R}^k_3(t,\cdot)\|_{L^2} +t^\beta \sum_{|\mu|=0}^1\|\mathrm{R}_1^\mu u_\pm(t,\cdot)\|_{L^\infty} E^1_3(t;W)^\frac{1}{2}\\
&\le CB\varepsilon t^{-\frac{1}{2}+\beta +\frac{\delta_1}{2}},
\end{split}
\end{equation*}
with $\beta>0$ small, $\beta\rightarrow 0$ as $\sigma\rightarrow 0$, while for $(J_1,J_2)\in\mathcal{I}(I_2)$ such that $|J_1|+|J_2|<|I_2|$ (hence $<2$) an estimate such as \eqref{est_J1J2_less 2} holds. These estimates, together with \eqref{est: boostrap vpm}, imply that the right hand side of \eqref{last_J1J2} is bounded by $CAB\varepsilon^2 t^{-\frac{3}{2}+\beta'}$, for a new small $\beta'>0$, $\beta'\rightarrow 0$ as $\sigma, \delta_0\rightarrow 0$, and that concludes the proof of the statement.
\endproof
\end{lem}

\begin{cor} \label{Cor: L2 est QI0(V,W)}
Let $Q^I_0(V,W)$ be the vector defined in \eqref{matrix QI}.
There exists a constant $C>0$ such that, if we assume that a-priori estimates \eqref{est: bootstrap argument a-priori est} are satisfied in interval $[1,T]$, for some fixed $T>1$, with $\varepsilon_0<(2A+B)^{-1}$ small:

$(i)$ if $I\in\mathcal{I}_n$ with $n\ge 3$:
\begin{equation}\label{est:QI0-In}
\|Q^I_0(V,W)\|_{L^2}\le  C A\varepsilon t^{-\frac{1}{2}+\frac{\delta}{2}};
\end{equation}
$(ii)$ if $I\in\mathcal{I}^k_3$, with $0\le k\le 2$,
\begin{equation} \label{est:QI0-Ik3}
\|Q^I_0(V,W)\|_{L^2} \le C(A+B)\varepsilon t^{-\frac{1}{2}+\frac{\delta_k}{2}}.
\end{equation}
\proof
$(i)$ Inequality \eqref{est:QI0-In} is straightforward after definition \eqref{matrix QI} (all coefficients $c_{I_1,I_2}$ are equal to 0 when $I\in\mathcal{I}_n$), lemma \ref{Lem: L2 est nonlinearities} $(i)$, and a-priori estimates \eqref{est: bootstrap upm}, \eqref{est: boostrap vpm}.

$(ii)$ If $I\in\mathcal{I}^k_3$ for a fixed $0\le k\le 2$ we have by definition \eqref{matrix QI} and lemmas \ref{Lem: L2 est nonlinearities}, \ref{Lem:L2 est nonlinearity Dt} that
\begin{equation}\label{dec_1}
\sum_{\substack{(I_1,I_2)\in\mathcal{I}(I)\\ |I_2|<|I|}} Q^\mathrm{w}_0(v^{I_1}_\pm, D_x v^{I_2}_\pm) + \sum_{\substack{(I_1,I_2)\in\mathcal{I}(I)\\|I_1|+|I_2|\le 2, |I_2|<|I|}} Q^\mathrm{w}_0(v^{I_1}_\pm, D_t v^{I_2}_\pm)=\mathfrak{R}^k_3(t,x),
\end{equation}
with $\mathfrak{R}^k_3(t,x)$ satisfying \eqref{est:L2_norm_Rk3(t,x)}.
Moreover, for some smooth $\chi\in C^\infty_0(\mathbb{R}^2)$, equal to 1 in a neighbourhood of the origin and $\sigma>0$ small,
\begin{equation}\label{dec_Qkg0}
\begin{gathered}
\sum_{\substack{(I_1,I_2)\in\mathcal{I}(I)\\ |I_2|<|I|}} Q^\mathrm{kg}_0(v^{I_1}_\pm, D_x u^{I_2}_\pm)=\delta_{\mathcal{V}^k}\sum_{\substack{(I_1,I_2)\in \mathcal{I}(I)\\ I_1\in\mathcal{K}, |I_2|\le 1}}Q^\mathrm{kg}_0\left(v^{I_1}_\pm, \chi(t^{-\sigma}D_x)D_x u^{I_2}_\pm\right)+ \mathfrak{R}^k_3(t,x), \\
\sum_{\substack{(I_1,I_2)\in\mathcal{I}(I)\\|I_1|+|I_2|\le 2, |I_2|<|I|}} Q^\mathrm{kg}_0(v^{I_1}_\pm, D_t u^{I_2}_\pm)=\delta_{\mathcal{V}^k}\sum_{\substack{(J,0)\in \mathcal{I}(I)\\J\in\mathcal{K}}}Q^\mathrm{kg}_0\left(v^J_\pm, \chi(t^{-\sigma}D_x)|D_x| u_\pm\right)+ \mathfrak{R}^k_3(t,x), 
\end{gathered}
\end{equation}
with sets $\mathcal{K}, \mathcal{V}^k$ given, respectively, by \eqref{set_K}, \eqref{set_V}, $\delta_{\mathcal{V}^k}=1$ if $I\in\mathcal{V}^k$, 0 otherwise (remind that $\mathcal{V}^2$ is empty).
Observe that, if $k=0,1$, $I\in\mathcal{I}^k_3$ and $(I_1,I_2)\in\mathcal{I}(I)$ with $I_1\in\mathcal{K}$, two situations may occur: if $\Gamma^{I_2}\in\{D^\alpha_x, |\alpha|\le 1\}$ then product $\Gamma^{I_1}$ contains exactly the same number of Klainerman vector fields as in $\Gamma^I$ and $V^{I_1}$ would be at the same energy level as $V^I$ (i.e. its $L^2$ norm being controlled by $E^k_3(t;W)^{1/2}$). In this case, from a-priori estimates \eqref{est: bootstrap upm}
\begin{equation}\label{product vI1, uI2_1}
\begin{split}
\|v^{I_1}_\pm(t,\cdot)\|_{L^2}\left(\|\chi(t^{-\sigma}D_x)u^{I_2}_\pm(t,\cdot)\|_{H^{\rho,\infty}}+ \|\chi(t^{-\sigma}D_x)\mathrm{R} u^{I_2}_\pm(t,\cdot)\|_{H^{\rho,\infty}}\right) &\le A\varepsilon t^{-\frac{1}{2}}E^k_3(t;W)^\frac{1}{2}.
\end{split}
\end{equation}
If instead $I_2$ is such that $\Gamma^{I_2}\in \{\Omega, Z_m, m=1,2\}$ is a Klainerman vector field, we automatically have that $\Gamma^I$ is a product of three Klainerman vector fields and that $V^{I_1}$ is at an energy level strictly lower than $V^I$ (i.e. its $L^2$ norm is controlled by $E^1_3(t;W)^{1/2}$ whereas that of $V^I$ is bounded by $E^0_3(t;W)^{1/2}$). From lemma
\ref{Lem_appendix: est UJ} we deduce that
\begin{multline}\label{product_vI1-uI2_2}
\|v^{I_1}_\pm(t,\cdot)\|_{L^2}\left(\|\chi(t^{-\sigma}D_x)u^{I_2}_\pm(t,\cdot)\|_{H^{\rho,\infty}}+ \|\chi(t^{-\sigma}D_x)\mathrm{R}u^{I_2}_\pm(t,\cdot)\|_{H^{\rho,\infty}}\right)\\
 \le C(A+B)\varepsilon t^{-\frac{1}{2}+\beta+\frac{\delta_1}{2}}E^1_3(t;W)^\frac{1}{2},
\end{multline}
for a small $\beta>0$, $\beta\rightarrow 0$ as $\sigma\rightarrow 0$.
Summing up \eqref{dec_1} to \eqref{product_vI1-uI2_2} and using \eqref{est:L2_norm_Rk3(t,x)} we obtain that there is a positive constant $C$ such that
\begin{equation} \label{ineq_QI0_L2}
\|Q^I_0(V,W)\|_{L^2} \le \delta_k C(A+B)\varepsilon t^{-\frac{1}{2}}\Big[E^k_3(t;W)^\frac{1}{2} +\delta_0 t^{\beta+\frac{\delta_1}{2}}E^1_3(t;W)^\frac{1}{2}\Big]+ CB\varepsilon t^{-\frac{5}{4}},
\end{equation}
with $\delta_k=1$ for $k=0,1$, equal to 0 when $k=2$, and $\delta_0=1$ only when $k=0$, 0 otherwise.
Finally, taking $\sigma>0$ small so that $\beta+\delta_1/2\ll \delta_0/2$ and using a-priori estimates \eqref{est: bootstrap E02} we deduce estimate \eqref{est:QI0-Ik3} from \eqref{ineq_QI0_L2}.
\endproof
\end{cor}

\subsection{Symmetrization}\label{Subs: Symmetrization}

\begin{prop} \label{Prop: equation of WIs}
As long as $H^{1,\infty}$ norm of $V(t,\cdot)$ is sufficiently small, there exists a real matrix $P(V;\eta)$ of order 0 and a real symmetric matrix $\widetilde{A}_1(V;\eta)$ of order 1, vanishing at order 1 at $V=0$, such that 
\begin{equation}\label{def_WIs}
W^I_s:=Op^B\big(P(V;\eta)\big)W^I
\end{equation}
is solution to \index{WIs@$W^I_s$, wave-Klein-Gordon vector after symmetrization}
\begin{equation} \label{equation WIs}
\begin{split}
D_t W^I_s & = A(D) W^I_s + Op^B(\widetilde{A}_1(V;\eta))W^I_s + Op^B(A''(V^I;\eta))U \\
&+ Op^B(C''(U;\eta))V^I + Op^B_R(A''(V^I;\eta))U + Q^I_0(V,W) + \mathfrak{R}(U, V),
\end{split}
\end{equation}
where $\mathfrak{R}(U, V)$ satisfies, for any $\theta\in ]0,1[$,
\begin{equation} \label{L2 est of R(U,V)}
\begin{split}
\|\mathfrak{R}(U, V) (t,\cdot)\|_{L^2}  &\lesssim \Big[\|V(t,\cdot)\|_{H^{7,\infty}} 
+ \|V(t,\cdot)\|^{1-\theta}_{H^{1,\infty}} \|V(t,\cdot)\|^\theta_{H^{3}}\left(\| U(t,\cdot)\|_{H^{2,\infty}}+ \|\mathrm{R}_1U(t,\cdot)\|_{H^{2,\infty}}\right)\\
&+ \|V(t,\cdot)\|_{L^\infty}\left(\| U(t,\cdot)\|^{1-\theta}_{H^{2,\infty}}+ \|\mathrm{R}_1U(t,\cdot)\|^{1-\theta}_{H^{2,\infty}}\right)\|U(t,\cdot)\|^\theta_{H^{4}} \Big] \|W^I(t,\cdot)\|_{L^2} \\
& + \|V(t,\cdot)\|_{H^{1,\infty}} \Big(\|W(t,\cdot)\|_{H^{7,\infty}} + \|\mathrm{R}U(t,\cdot)\|_{H^{6,\infty}}\Big) \|W^I(t,\cdot)\|_{L^2} \\
& + \|V(t,\cdot)\|_{H^{1,\infty}}\|Q^I_0(V,W)\|_{L^2}.
\end{split}
\end{equation}
Moreover, for any $n,r\in\mathbb{N}$,
\begin{gather}
M^0_r\left(P(V;\eta)-I_4;n\right)\lesssim \|V(t,\cdot)\|_{H^{1+r,\infty}},  \label{seminorm P-I4} \\
M^1_r\big(\widetilde{A}_1(V;\eta);n\big)\lesssim \|V(t,\cdot)\|_{H^{1+r,\infty}}, \label{seminorm Atilde 1}
\end{gather} 
and as long as the $H^{2,\infty}$ norm of $V(t,\cdot)$ is small there is a constant $C>0$ such that 
\begin{equation}\label{equivalence WIs WI}
C^{-1}\|W^I(t,\cdot)\|_{L^2}\le \|W^I_s(t,\cdot)\|_{L^2}\le C \|W^I(t,\cdot)\|_{L^2}.
\end{equation}
\end{prop}
In order to prove proposition \ref{Prop: equation of WIs}, we first need to introduce the following lemma.

\begin{lem} \label{Lemma: existence of P(alpha, beta)}
Let $\alpha,\beta\in\mathbb{R}$, $L\in M_2(\mathbb{R})$ and $M_0, N(\alpha,\beta)\in M_4(\mathbb{R})$ given by
\begin{equation*}
L =
\begin{bmatrix}
0 & 1 \\
1 & 0
\end{bmatrix},\quad
M_0 = 
\begin{bmatrix}
I_2 & 0 \\
0 & -I_2
\end{bmatrix}, \quad 
N(\alpha,\beta)= 
\begin{bmatrix}
\alpha L & \beta L \\
\alpha L & \beta L 
\end{bmatrix}
=
\begin{bmatrix}
0 & \alpha & 0 & \beta \\
\alpha & 0 & \beta & 0 \\
0 & \alpha & 0 & \beta \\
\alpha & 0 & \beta & 0
\end{bmatrix} .
\end{equation*}
There exist a small $\delta>0$ and a smooth function defined on open ball $B_\delta(0)$ of radius $\delta$, $$(\alpha,\beta)\in B_\delta(0)\rightarrow P(\alpha,\beta)\in Sym_4(\mathbb{R}),$$
with values in the space of real, symmetric, $4\times 4$ matrices $Sym_4(\mathbb{R})$, such that $P(0,0)=I_4$, $P(\alpha,\beta)=I_4 + O(|\alpha|+|\beta|)$ and $P(\alpha,\beta)^{-1}\big(M_0 + N(\alpha, \beta)\big)P(\alpha,\beta)$ is symmetric for any $(\alpha, \beta)\in B_\delta(0)$.
Furthermore $P^{-1}(\alpha,\beta)= I_4 + O(|\alpha| + |\beta|)$.
\proof
Let $\mathcal{E}$ be the vector space of $2\times 2$ matrices $B(\alpha, \beta)=\alpha I_2 + \beta L$ and $\mathcal{F}$ be the set of $4\times 4$ matrices of the form 
\[
\begin{bmatrix}
F_{11} & F_{12} \\
F_{21} & F_{22}
\end{bmatrix}
\]
with $F_{ij}\in \mathcal{E}$.
We look for a matrix $P$ of the form
\begin{equation}\label{matrix P(B)}
P(B) = (I_2 - B^2)^{-\frac{1}{2}}
\begin{bmatrix}
I_2 & -B \\
-B & I_2
\end{bmatrix}
\end{equation}
with $B\in \mathcal{E}$ close to zero (so that in particular $(I_2-B^2)^{1/2}$ is well defined). 
We remark that matrix $P(B)^{-1}$ has the form
\[
P(B)^{-1}= (I_2 - B^2)^{-\frac{1}{2}}
\begin{bmatrix}
I_2 & B \\
B & I_2
\end{bmatrix}
\]
and that $P(0)=P^{-1}(0)=I_4$.
We consider $\Phi : \mathbb{R}^2\times \mathcal{E}\rightarrow \mathcal{F}$ defined by $\Phi(\alpha,\beta, B):= P(B)^{-1}\big[M_0 + N(\alpha, \beta)\big]P(B) = \big(\Phi_{ij}(\alpha, \beta, B)\big)_{1\le i,j\le 2}$, where $\Phi_{ij}\in \mathcal{E}$ as $\mathcal{E}$ is a commutative sub-algebra of $M_2(\mathbb{R})$. We also define $\Psi(\alpha, \beta, B):= \Phi_{12}(\alpha, \beta, B) - \Phi_{21}^\dagger(\alpha,\beta, B)$ with $\Phi_{21}^\dagger$ denoting the transpose of $\Phi_{21}$. We have that $\Psi(0,0, 0)=0$ and
\begin{equation*}
D_B\Phi(0,0,0)\cdot B = 
\begin{bmatrix}
0 & B \\
B & 0 
\end{bmatrix}M_0  - M_0
\begin{bmatrix}
0 & B \\
B & 0 
\end{bmatrix} = 2
\begin{bmatrix}
0 & -B \\
B & 0 
\end{bmatrix}
\end{equation*}
from which follows that $D_B\Psi(0,0,0)\cdot B = -4B$, i.e. $D_B\Psi(0,0,0) = -4I$. Therefore, there exist a small $\delta>0$ and a smooth function $(\alpha, \beta)\in B_\delta(0)\rightarrow B(\alpha,\beta)\in \mathcal{E}$ such that $B(0,0)=0$ (which implies $P(B(0,0))=I_4$), and $\Psi(\alpha, \beta, B(\alpha, \beta)) = 0$ $\forall (\alpha,\beta)\in B_\delta(0)$. This is equivalent to say that $\Phi(\alpha,\beta, B(\alpha,\beta))$ is symmetric and
moreover $P(B(\alpha,\beta)), P(B(\alpha,\beta))^{-1} = I_4 + O(|\alpha|+|\beta|)$. Defining $P(\alpha,\beta):=P(B(\alpha,\beta))$ concludes the proof of the statement.
\endproof
\end{lem}

\proof[Proof of proposition \ref{Prop: equation of WIs}]

With notations introduced in lemma \ref{Lemma: existence of P(alpha, beta)} and in \eqref{matrices A A'}, \eqref{matrices A'1 A'-1}, $A(\eta) = \langle \eta\rangle M_0 + S(\eta)$ and $A'_1(V;\eta)(1-\chi)(\eta) = \langle \eta\rangle N(\alpha,\beta)$, with
\[
S(\eta) = 
\begin{bmatrix}
|\eta| -\langle\eta\rangle & 0 & 0 & 0 \\
0 & 0 & 0 & 0 \\
0 & 0 & -(|\eta| - \langle\eta\rangle) & 0 \\
0 & 0 & 0 & 0
\end{bmatrix} \text{whose elements are $O(|\eta|^{-1}), |\eta|\rightarrow +\infty$},
\]
and $\alpha = a_0(v_{\pm};\eta)\frac{\eta_1}{\langle\eta\rangle}(1-\chi)(\eta)$, $\beta= b_0(v_{\pm};\eta)\frac{\eta_1}{\langle\eta\rangle}(1-\chi)(\eta)$, $a_0, b_0$ defined in \eqref{def ak, bk, a0, b0}. 
Since $\sup_\eta\big(|\alpha| + |\beta|\big) \lesssim \|V(t,\cdot)\|_{H^{1,\infty}}$, by lemma \ref{Lemma: existence of P(alpha, beta)} we have that, as long as $\|V(t,\cdot)\|_{H^{1,\infty}}$ is sufficiently small, there exists a real symmetric matrix $P=P(V;\eta)$ of the form \eqref{matrix P(B)} such that $P(V;\eta)^{-1}\big[M_0 + N(\alpha, \beta)\big]P(V;\eta)$ is real and symmetric. Moreover $P = I_4 + Q(V;\eta)$ and $P^{-1} = I_4 + Q'(V;\eta)$, where $Q(V;\eta)$, $Q'(V;\eta)$ are matrices depending smoothly on $\alpha,\beta$ (which are symbols of order 0), null at order 1 at $V=0$, verifying for any $n,r\in\mathbb{N}$
\begin{equation*}
M^0_r\left(Q(V;\eta);n\right)+ M^0_r\left(Q'(V;\eta);n\right)\lesssim \|V(t,\cdot)\|_{H^{1+r,\infty}}.
\end{equation*}
We define the following matrix of order 1
\begin{equation*}
\widetilde{A}_1(V;\eta) := P(V;\eta)^{-1}\big[\langle\eta\rangle\big(M_0 + N(\alpha, \beta)\big)\big]P(V;\eta) - \langle\eta\rangle M_0
\end{equation*}
and $W^I_s:= Op^B(P^{-1}(V;\eta))W^I$. 
From the fact that $\widetilde{A}_1(V;\eta)$ also writes as
\begin{equation*}
\langle\eta\rangle \left[Q'(V;\eta)M_0 + P^{-1}(V;\eta)M_0Q(V;\eta) + P^{-1}(V;\eta)N(\alpha,\beta)P(V;\eta)\right]
\end{equation*}
we see that it vanishes at order 1 at $V=0$ and is such that $M^1_r(\widetilde{A}_1(V;\eta);n)\lesssim \|V(t,\cdot)\|_{H^{1+r,\infty}}$. 
Moreover, from proposition \ref{Prop: paradifferential symbolic calculus} $(ii)$ with $r=1$ it follows that 
\begin{equation} \label{dec_operator_I}
I = Op^B(P(V;\eta))Op^B(P^{-1}(V;\eta)) + T_{-1}(V),
\end{equation}
where operator $T_{-1}(V)$ is of order less or equal than $-1$ whose $\mathcal{L}(L^2)$ norm is a $O(\|V(t,\cdot)\|_{H^{2,\infty}})$.
Therefore $W^I = Op^B(P(V;\eta))W^I_s + T_{-1}(V)W^I$ and from proposition \ref{Prop : Paradiff action on Sobolev spaces-NEW} the $L^2$ norms of $W^I, W^I_s$ are equivalent as long as the $H^{2,\infty}$ norm of $V$ is small.
Using equation \eqref{equation WI-1} we find that:
\begin{equation} \label{eq: preliminary equation WIs}
\begin{split}
D_t W^I_s &= Op^B(P^{-1}(V;\eta))Op^B\big(A(\eta) + A'_1(V;\eta)(1-\chi)(\eta)\big)W^I  \\
& + Op^B(P^{-1}(V;\eta))\Big[Op^B\big(A'_1(V;\eta)\chi(\eta)\big) + Op^B\big(A'_{-1}(V;\eta)\big)\Big]W^I \\
& + Op^B(P^{-1}(V;\eta))\Big[Op^B(C'(W^I;\eta))V + Op^B_R(A'(V;\eta))W^I \Big]\\
& +  Op^B(P^{-1}(V;\eta))\Big[  Op^B(A''(V^I;\eta))U 
+ Op^B(C''(U;\eta))V^I + Op^B_R(A''(V^I;\eta))U\Big] \\
& +  Op^B(P^{-1}(V;\eta))Q^I_0(V,W) + Op^B(D_tP^{-1}(V;\eta))W^I
\end{split}
\end{equation}
where
\small
\begin{equation} \label{equation WIs - calculations}
\begin{split}
& Op^B(P^{-1}(V;\eta))Op^B\big(A(\eta) + A'_1(V;\eta)(1-\chi)(\eta)\big)W^I \\
&=   Op^B(P^{-1}(V;\eta))Op^B\left(\langle\eta\rangle\big(M_0 + N(\alpha,\beta)\big)\right)W^I + Op^B(S(\eta))W^I + Op^B(Q'(V;\eta))Op^B(S(\eta))W^I \\
&= Op^B(P^{-1}(V;\eta))Op^B\left(\langle\eta\rangle\big(M_0 + N(\alpha,\beta)\big)\right)Op^B(P(V;\eta))W^I_s \\
& +   Op^B(P^{-1}(V;\eta))Op^B\left(\langle\eta\rangle\big(M_0 + N(\alpha,\beta)\big)\right)T_{-1}(V)W^I + Op^B(S(\eta))W^I_s \\
& +  Op^B(S(\eta))Op^B(Q(V;\eta))W^I_s+ Op^B(S(\eta))T_{-1}(V)W^I + Op^B(Q'(V;\eta))Op^B(S(\eta))W^I \\
= & Op^B(A(\eta) + \widetilde{A}_1(V;\eta))W^I_s + \widetilde{T}_0(V)W^I_s + \widetilde{T'}_0(V)W^I
\end{split}
\end{equation}
\normalsize
with $\widetilde{T}_0(V), \widetilde{T'}_0(V)$ operators of order 0 and $\mathcal{L}(L^2)$ norm $O(\|V(t,\cdot)\|_{H^{2,\infty}})$. Last equality follows indeed from the fact that, by proposition \ref{Prop: paradifferential symbolic calculus} $(ii)$ with $r=1$ and proposition \ref{Prop : Paradiff action on Sobolev spaces-NEW},
\begin{multline*}
 Op^B(P^{-1}(V;\eta))Op^B\big[\langle\eta\rangle\big(M_0 + N(\alpha,\beta)\big)\big]Op^B(P(V;\eta)) \\
=  Op^B\big(P(V;\eta)^{-1}\big[\langle\eta\rangle\big(M_0 + N(\alpha, \beta)\big)\big]P(V;\eta)\big) + \widetilde{T}_0(V)
\end{multline*}
and $Op^B(S(\eta))Op^B(Q(V;\eta))$, $Op^B(Q'(V;\eta))Op^B(S(\eta))$ are operator of order 0, too (the former of the form $\widetilde{T}_0(V)$, the latter of the form $T_0(V)$), while $Op^B(S(\eta))T_{-1}(V)$ is of order $-1$ (and can be included in $T_0(V)$).
The equivalence between the $L^2$ norms of $W^I_s$ and $W^I$ implies that $\widetilde{T}_0(V)W^I_s + \widetilde{T'}_0(V)W^I$ in \eqref{equation WIs - calculations} is a remainder $\mathfrak{R}(U, V)$. 

All operators appearing in the second and third line of \eqref{eq: preliminary equation WIs} are also remainders $\mathfrak{R}(U,V)$ because, from proposition \ref{Prop : Paradiff action on Sobolev spaces-NEW}, the fact that $M^0_0(P^{-1}(V;\eta);2)=O(1)$ and lemma \ref{Lemma: L2 estimate of semilinear terms}, their $L^2$ norm is bounded by $\|V(t,\cdot)\|_{H^{7,\infty}}\|W^I(t,\cdot)\|_{L^2}$.
Last term in \eqref{eq: preliminary equation WIs} also contributes to $\mathfrak{R}(U,V)$ for matrix $D_tP^{-1}(V;\eta)$ is of order 0, its $M^0_0(\cdot,2)$ seminorm is bounded by $\|D_tV(t,\cdot)\|_{H^{1,\infty}}$ and for any $\theta \in[0,1]$
\begin{equation*}
\begin{split}
\|D_tV(t,\cdot)\|_{H^{1,\infty}}&\lesssim \|V(t,\cdot)\|_{H^{2,\infty}} + \|V(t,\cdot)\|^{1-\theta}_{H^{1,\infty}} \|V(t,\cdot)\|^\theta_{H^{3}}\left(\| U(t,\cdot)\|_{H^{2,\infty}}+ \|\mathrm{R}_1U(t,\cdot)\|_{H^{2,\infty}}\right)\\
&+ \|V(t,\cdot)\|_{L^\infty}\left(\| U(t,\cdot)\|^{1-\theta}_{H^{2,\infty}}+ \|\mathrm{R}_1U(t,\cdot)\|^{1-\theta}_{H^{2,\infty}}\right)\|U(t,\cdot)\|^\theta_{H^{4}},
\end{split}
\end{equation*}
as follows from \eqref{est: Hsinfty Dt V} with $s=1$.
Finally, in remaining contributions in \eqref{eq: preliminary equation WIs} we replace $Op^B(P^{-1}(V;\eta))$ with $I+Op^B(Q'(V;\eta))$ and observe that the terms on which $Op^B(Q'(V;\eta))$ acts are remainders $\mathfrak{R}(U,V)$ after proposition \ref{Prop : Paradiff action on Sobolev spaces-NEW}, the fact that $M^0_0\left(Q'(V;\eta);2\right) = O(\|V(t,\cdot)\|_{H^{1,\infty}})$ and lemma \ref{Lemma: L2 estimate of semilinear terms}. 
Interchanging the notation of $P(V;\eta)$ and $P^{-1}(V;\eta)$, we obtain the result of the statement.
\endproof

\section{Normal forms and energy estimates} \label{Sec: Quasi-linear normal forms and energy estimates}

Before going further in writing an energy inequality for $W^I_s$ we should make few remarks.
As we previously anticipated, the $L^2$ norm of some of the semi-linear terms appearing in equation \eqref{equation WIs} have a very slow decay in time. 
On the one hand, it is the case of $Op^B(A''(V^I;\eta))U$, $Op^B(C''(U;\eta))V^I$ and $Op^B_R(A''(V;\eta))U$, whose $L^2$ norms are estimated in \eqref{L2 est on OpBR(A'')W}, \eqref{L2 est on OpB(C")VI} in terms of the uniform norms of $U, \mathrm{R}_1U$. On the other hand, also some of the contributions to $Q^I_0(V,W)$ are only a $O_{L^2}(t^{-1/2+\beta'})$,
for some small $\beta'>0$, after corollary \ref{Cor: L2 est QI0(V,W)}.
Nevertheless, we are going to see that $Op^B(A''(V^I;\eta))U$, $Op^B_R(A''(V;\eta))U$ and the mentioned contributions to $Q^I_0(V,W)$ can be easily eliminated by performing a semi-linear normal form argument in the energy inequality (see subsection \ref{sub: second normal form}). 
Such an argument is however not well adapted to handle $Op^B(C''(U;\eta))V^I$, for it leads to a loss of derivatives linked to the quasi-linear nature of the problem, i.e. to the fact that matrix $\widetilde{A}_1(V;\eta)$ in the right hand side of \eqref{equation WIs} is of order 1. This latter contribution should instead be eliminated through a suitable normal form applied directly on equation \eqref{equation WIs}, which is the object of the subsection \ref{sub: a first normal form transformation}.

\subsection{A first normal forms transformation and the energy inequality} \label{sub: a first normal form transformation}

First of all, we replace $Op^B(C''(U;\eta))V^I$ in equation \eqref{equation WIs} with $Op^B(C''(U;\eta))V^I_s$, having defined $V^I_s:= Op^B(P^{-1}(V;\eta))V^I$, and remind that from \eqref{seminorm P-I4} with $r=0$ and \eqref{dec_operator_I} the $L^2$ norm of $V^I$ and $V^I_s$ are equivalent as long as the $H^{2,\infty}$ norm of $V(t,\cdot)$ is small (assumption compatible with \eqref{est: boostrap vpm} if $\rho\ge 2$).
We will rather deal with 
\begin{equation} \label{equation WIs_1} 
\begin{split}
(D_t - A(D)) W^I_s  & =  Op^B(\widetilde{A}_1(V;\eta))W^I_s  
+ Op^B(A''(V^I;\eta))U + Op^B(C''(U;\eta))V^I_s\\ 
& + Op^B_R(A''(V^I;\eta))U + Q^I_0(V,W) + \mathfrak{R}(U, V),
\end{split}
\end{equation}
for a new $\mathfrak{R}(U,V)$ satisfying \eqref{L2 est of R(U,V)} and show how to get rid of $Op^B(C''(U;\eta))V^I_s$ in the above right hand side.
More precisely, we are going to prove the following result:

\begin{prop} \label{Prop: a first normal form}
Let $N\in\mathbb{N}^*$.
There exist three matrices $E^0_d(U;\eta), E^{-1}_d(U;\eta)$, $E_{nd}(U;\eta)$ linear in $(u_+,u_{-})$, with $E^0_d(U;\eta)$ real diagonal of order 0 and $E^{-1}_d(U;\eta), E_{nd}(U;\eta)$ of order -1, and, as long as $\|\mathrm{R}_1U(t,\cdot)\|_{H^{2,\infty}}$ is small, a real diagonal matrix $F^0_d(U;\eta)$ of order 0 such that, if
\begin{equation}\label{def_WtildeI}
\begin{gathered}
\widetilde{W}^I_s:= Op^B(I_4 + E(U;\eta))W^I_s ,\\
\text{with } E(U;\eta) := E^0_d(U;\eta)+ E^{-1}_d(U;\eta) + E_{nd}(U;\eta),
\end{gathered}
\end{equation}
then
\begin{equation} \label{equation Wtilde-Is}
\begin{split}
&(D_t -A(D))\widetilde{W}^I_s =  Op^B\left((I_4+ E^0_d(U;\eta))\widetilde{A}_1(V;\eta)(I_4 + F^0_d(U;\eta))\right)\widetilde{W}^I_s \\
& + Op^B(A''(V^I;\eta))U + Op^B_R(A''(V^I;\eta))U  + Q^I_0(V,W) + T_{-N}(U)W^I_s + \mathfrak{R}'(U, V).
\end{split}
\end{equation}
In the above right hand side $T_{-N}(U)=(\sigma_{ij}(U,D_x))_{ij}$ is a pseudo-differential operator of order less or equal than $-N$, with
\begin{equation} \label{norm of T-3 in propositon}
\|T_{-N}(U)\|_{\mathcal{L}(H^{s-N};H^s)}\lesssim  \|\mathrm{R}_1U(t,\cdot)\|_{H^{N+2,\infty}} + \|U(t,\cdot)\|_{H^{N+6,\infty}},
\end{equation}
for any $s\in\mathbb{R}$ and such that
\begin{subequations}  \label{sigmaij in proposition Normal Forms}
\begin{equation}
\mathcal{F}_{x\mapsto\xi}(\sigma_{ij}(U,\eta))(\xi) = 
\begin{cases}
\sigma^+_{ij}(\xi, \eta)\hat{u}_+(\xi) + \sigma^{-}_{ij}(\xi,\eta)\hat{u}_{-}(\xi), \qquad & i,j \in\{2,4\}, \\
 0, &\text{otherwise},
\end{cases}
\end{equation} 
with $\sigma^\pm_{ij}(\xi,\eta)$ supported for $|\xi|\le \varepsilon\langle\eta\rangle$ for a small $\varepsilon>0$, and for any $\alpha,\beta\in\mathbb{N}^2$
\begin{equation}
|\partial^\alpha_\xi \partial^\beta_\eta \sigma^\pm_{ij}(\xi,\eta)|\lesssim_{\alpha,\beta}
|\xi|^{N+1-|\alpha|}\langle\eta\rangle^{-N-|\beta|}, \quad i,j \in \{2,4\}.
\end{equation}
\end{subequations}
Also, $\mathfrak{R}'(U,V)$ is a remainder satisfying, for any $\theta\in ]0,1[$
\begin{equation} \label{L2 norm R'(U,V)}
\begin{split}
\|\mathfrak{R}' (U, V)& (t,\cdot)\|_{L^2}  \lesssim (1 + \|U(t,\cdot)\|_{H^{5,\infty}})\|\mathfrak{R}(U,V)\|_{L^2} \\
&+ \left(\|\mathrm{R}_1U(t,\cdot)\|_{H^{1,\infty}} + \|U(t,\cdot)\|_{H^{5,\infty}} \right) \Big[ \|Q^I_0(V,W)\|_{L^2} \\
& +\left(\|\mathrm{R} U(t,\cdot)\|_{H^{6,\infty}}+ \|U(t,\cdot)\|_{H^{6,\infty}}\right)\|W^I(t,\cdot)\|_{L^2} \Big]\\
&   + \|V(t,\cdot)\|^{2-\theta}_{H^{5,\infty}} \|V(t,\cdot)\|^\theta_{H^7} \|W^I(t,\cdot)\|_{L^2},
\end{split}
\end{equation} 
with $\mathfrak{R}(U,V)$ verifying \eqref{L2 est of R(U,V)}.

For any $n,r\in\mathbb{N}$, any $\chi\in C^\infty_0(\mathbb{R}^2)$ equal to 1 close to the origin and supported in open ball $B_{\varepsilon}(0)$, with $\varepsilon>0$ sufficiently small, we have that
\begin{subequations} \label{seminorms E}
\begin{equation} \label{seminorm E0d}
M^0_r\left(E^0_d\left(\chi\left(\frac{D_x}{\langle\eta\rangle}\right) U;\eta\right);n\right) \lesssim \|\mathrm{R}_1 U(t,\cdot)\|_{H^{1+r,\infty}},
\end{equation}
\begin{equation}\label{seminorm E-1d}
M^{-1}_r\left(E^{-1}_d\left(\chi\left(\frac{D_x}{\langle\eta\rangle}\right) U;\eta\right);n\right) \lesssim \|U(t,\cdot)\|_{H^{5+r,\infty}},
\end{equation}
\begin{equation}  \label{seminorm End}
M^{-1}_r\left(E_{nd}\left(\chi\left(\frac{D_x}{\langle\eta\rangle}\right) U;\eta\right);n\right) \lesssim \|U(t,\cdot)\|_{H^{5+r,\infty}};
\end{equation}
\end{subequations}
and
\begin{equation} \label{seminorm F0}
M^0_r\left(F^0_d\left(\chi\left(\frac{D_x}{\langle\eta\rangle}\right) U;\eta\right);n\right) \lesssim \|\mathrm{R}_1 U(t,\cdot)\|_{H^{1+r,\infty}}.
\end{equation}
Finally, as long as $\|\mathrm{R}_1U(t,\cdot)\|_{H^{2,\infty}}+\|U(t,\cdot)\|_{H^{5,\infty}}$ is small, there is a constant $C>0$ such that
\begin{equation} \label{equivalence WtildeIs WIs}
C^{-1}\|W^I_s(t,\cdot)\|_{L^2}\le \|\widetilde{W}^I_s(t,\cdot)\|_{L^2}\le C\|W^I_s(t,\cdot)\|_{L^2}.
\end{equation}
\end{prop}

\begin{remark}  
From propositions \ref{Prop: equation of WIs} and \ref{Prop: a first normal form} it follows that, as long as $\|\mathrm{R}_1U(t,\cdot)\|_{H^{2,\infty}}$, $\|U(t,\cdot)\|_{H^{5,\infty}}$ and $\|V(t,\cdot)\|_{H^{2,\infty}}$ are small, there is a constant $C>0$ such that
\begin{equation} \label{remark on equivalence between L2 norms}
C^{-1}\|W^I(t,\cdot)\|_{L^2}\le \|\widetilde{W}^I_s(t,\cdot)\|_{L^2}\le C\|W^I(t,\cdot)\|_{L^2}.
\end{equation}
This implies that, if \index{Etilden@$\widetilde{E}_n(t;W)$, first modified energy} \index{Etildek@$\widetilde{E}^k_3(t;W)$, first modified energy}
\begin{subequations} \label{modified_energies_Etilde}
\begin{gather}
\widetilde{E}_n(t;W):= \sum_{|\alpha|\le n} \left\|Op^B(I_4+E(U;\eta))Op^B(P(V;\eta))D^\alpha_x W(t,\cdot)\right\|_{L^2}, \quad \forall \, n\in\mathbb{N}, n \ge 3,\label{energy E_tilde_n}\\
\widetilde{E}^k_3(t;W):=\sum_{\substack{|\alpha|+|I|\le 3\\ |I|\le 3-k}} \left\|Op^B(I_4+E(U;\eta))Op^B(P(V;\eta))D^\alpha_x W^I(t,\cdot)\right\|_{L^2}, \forall\, 0\le k\le 2,\label{energy E_tilde_k2}
\end{gather}
\end{subequations}
there exists a constant $C_1>0$ such that
\begin{equation}\label{energy_equivalence Ekm Etilde_km}
\begin{gathered}
C_1^{-1}E_n(t;W)\le \widetilde{E}_n(t;W)\le C_1 E_n(t;W), \quad \forall \, n\ge 3, \\
C_1^{-1}E^k_3(t;W)\le \widetilde{E}^k_3(t;W)\le C_1 E^k_3(t;W), \quad \forall \, 0\le k\le 2.
\end{gathered}
\end{equation}
Thanks to the above equivalence, the propagation of some suitable estimates on $\widetilde{E}_n(t;W)$ and $\widetilde{E}^k_3(t;W)$ will provide us with \eqref{est: bootstrap enhanced Enn} and \eqref{est: boostrap enhanced E02} respectively, so we can rather focus on the derivation of an energy inequality for $\widetilde{E}_n(t;W), \widetilde{E}^k_3(t;W)$.
\end{remark}

In order to get rid of $Op^B(C''_d(U;\eta))V^I_s$ in \eqref{equation WIs_1}we introduce matrices
\begin{equation}
 C''_d(U;\eta) = 
\begin{bmatrix}
0 & 0 & 0 & 0 \\
0 & e_0 & 0 & 0 \\
0 & 0 & 0 & 0 \\
0 & 0 & 0 & f_0
\end{bmatrix}, \quad
C''_{nd}(U;\eta) = 
\begin{bmatrix}
0 & 0 & 0 & 0 \\
0 & 0 & 0 & f_0 \\
0 & 0 & 0 & 0 \\
0 & e_0 & 0 & 0
\end{bmatrix}
\end{equation}
so that
\[C''(U;\eta) = C''_d(U;\eta) + C''_{nd}(U;\eta),\]
and proceed to eliminate $Op^B(C''_d(U;\eta))V^I_s$ and $Op^B(C''_{nd}(U;\eta))V^I_s$ separately.

\begin{lem} \label{Lem: Normal Forms on C''d(U,eta)}
Let $N\in\mathbb{N}^*$. There exists a diagonal matrix $E_d(U;\eta)$ of order 0, linear in $(u_+, u_{-})$, such that
\begin{multline} \label{equation for Op(Ed)}
Op^B(C''_d(U;\eta))V^I_s + Op^B(D_t E_d(U;\eta))W^I_s - [A(D),Op^B(E_d(U;\eta))]W^I_s \\ = T_{-N}(U)W^I_s  + \mathfrak{R}'(V,V),
\end{multline}
where $\mathfrak{R}'(V,V)$ satisfies, for any $\theta\in ]0,1[$,
\begin{equation} \label{L2 est R'(V,V)}
\|\mathfrak{R}'(V,V)(t,\cdot)\|_{L^2}\lesssim \|V(t,\cdot)\|^{2-\theta}_{H^{5,\infty}}\|V(t,\cdot)\|^\theta_{H^7}\|V^I(t,\cdot)\|_{L^2},
\end{equation}
and $T_{-N}(U)$ is a pseudo-differential operator of order less or equal than $-N$ such that, for any $s\in\mathbb{R}$,
\begin{equation} \label{norm of T-3WI}
\|T_{-N}(U) \|_{\mathcal{L}(H^{s-N};H^s)}\lesssim \|\mathrm{R}_1U(t,\cdot)\|_{H^{N+2,\infty}} + \|U(t,\cdot)\|_{H^{N+6,\infty}},
\end{equation}
whose symbol $\sigma(U,\eta) = \left(\sigma_{ij}(U,\eta)\right)_{1\le i,j\le 4}$ is such that
\begin{subequations} \label{sigma ij of operator T-3}
\begin{equation}
\mathcal{F}_{x\mapsto \xi}(\sigma_{ij}(U,\eta))(\xi) = 
\begin{cases}
\sigma^+_{ii}(\xi, \eta)\hat{u}_+(\xi) + \sigma^{-}_{ii}(\xi,\eta)\hat{u}_{-}(\xi), \qquad & i=j \in\{2,4\}, \\
 0, &\text{otherwise},
\end{cases}
\end{equation} 
with $\sigma^\pm_{ii}(\xi,\eta)$ supported for $|\xi|\le \varepsilon\langle\eta\rangle$ for a small $\varepsilon>0$, and verifying, for any $\alpha,\beta\in\mathbb{N}^2$,
\begin{equation}
|\partial^\alpha_\xi \partial^\beta_\eta \sigma^\pm_{ii}(\xi,\eta)|\lesssim_{\alpha,\beta}
|\xi|^{N+1-|\alpha|}\langle\eta\rangle^{-N-|\beta|}, \quad \text{for } i=2,4.
\end{equation}
\end{subequations}
Moreover, if $\chi\in C^\infty_0(\mathbb{R}^2)$ is equal to 1 close to the origin and has a sufficiently small support, 
\begin{equation}\label{decomposition Ed}
E_d\left(\chi\left(\frac{D_x}{\langle\eta\rangle}\right)U;\eta\right) = E^0_d\left(\chi\left(\frac{D_x}{\langle\eta\rangle}\right)U;\eta\right) + E^{-1}_d\left(\chi\left(\frac{D_x}{\langle\eta\rangle}\right)U;\eta\right),
\end{equation}
the former matrix in the above right hand side being real of order 0 and satisfying \eqref{seminorm E0d}, the latter being of order $-1$ and verifying \eqref{seminorm E-1d}.
\proof
Because of the diagonal structure of $A(\eta)$ and $C''_d(U;\eta)$ we look for a matrix $E_d = (e_{ij})_{1\le i,j\le 4}$ satisfying \eqref{equation for Op(Ed)} such that $e_{ij}=0$ for all $i,j$ but $i=j\in\{2,4\}$, and we also require symbols $e_{22}, e_{44}$ to be of order 0 and linear in $(u_+, u_{-})$.
If we remind that matrix $A(\eta)$ in \eqref{matrices A A'} is of order 1 and make the ansatz that $E_d$ is of order 0, then by symbolic calculus of proposition \ref{Prop: paradifferential symbolic calculus} we have that
\begin{equation} \label{commutator A(D) Ed}
-[A(D), Op^B(E_d(U;\eta))] = - \sum_{ |\alpha|=1}^N\frac{1}{\alpha!} Op^B\left(\partial^\alpha_\eta A(\eta)D^\alpha_xE_d(U;\eta)\right) + T_{-N}(U)
\end{equation}
with $T_{-N}(U)$ pseudo-differential operator of order less or equal than $-N$ such that, for any $s\in\mathbb{R}$, \small
\begin{equation} \label{norm T-N}
\|T_{-N}(U)\|_{\mathcal{L}(H^{s-N}; H^s)} \lesssim M^1_{N+1}(A(\eta);N+3)M^0_0(E_d(U;\eta);2) + M^1_{0}(A(\eta);N+3)M^0_{N+1}(E_d(U;\eta);2)
\end{equation}\normalsize
and whose symbol $\sigma(U,\eta)=\left(\sigma_{ij}(U,\eta)\right)_{ij}$ is such that $\sigma_{ij}(U,\eta)=0$ for all $i,j$ but $i=j\in \{2,4\}$.
Therefore, for any fixed $\chi\in C^\infty_0(\mathbb{R}^2)$ equal to 1 in $\overline{B_{\varepsilon_1}(0)}$ and supported in $B_{\varepsilon_2}(0)$, for some $0<\varepsilon_1<\varepsilon_2\ll 1$, we look for $E_d(U;\eta)$ such that
\begin{equation*} 
\chi\left(\frac{D_x}{\langle\eta\rangle}\right)\left[C''_d(U;\eta) + D_tE_d(U;\eta) - \sum_{ |\alpha| =1}^N \frac{1}{\alpha!} \partial^\alpha_\eta A(\eta) D^\alpha_x E_d(U;\eta)\right] = 0.
\end{equation*}
Since $E_d(U;\eta)$ is required to be linear in $(u_+, u_{-})$, we should write it rather as $E_d(u_+, u_{-}; \eta)$ to then realize that, as $u_+$ (resp. $u_{-}$) is solution to the first (resp. to the third) equation in \eqref{system for uI+-, vI+-} with $|I|=0$, 
\begin{align*}
D_t E_d(u_+, u_{-};\eta) &= E_d(|D_x| u_+, - |D_x|u_{-};\eta) + E_d\big(Q^{\mathrm{w}}_0(v_\pm, D_1 v_\pm), Q^{\mathrm{w}}_0(v_\pm, D_1 v_\pm);\eta\big), \\
D^\alpha_x E_d(u_+, u_{-};\eta) & = E_d(D^\alpha_x u_+, D^\alpha_x u_{-};\eta), \quad \forall \alpha\in\mathbb{N}^2.
\end{align*}
If we temporarily neglecting contribution $E_d\big(Q_0^{\mathrm{w}}(v_\pm, D_1 v_\pm), Q_0^{\mathrm{w}}(v_\pm, D_1 v_\pm);\eta\big)$, we are lead to solve the following equation
\begin{equation*} 
\chi\left(\frac{D_x}{\langle\eta\rangle}\right)\left[C''_d(U;\eta) + E_d(|D_x| u_+, - |D_x|u_{-};\eta) - \sum_{ |\alpha| =1}^N \frac{1}{\alpha!} \partial^\alpha_\eta A(\eta) E_d(D^\alpha_x u_+, D^\alpha_x u_{-};\eta)\right] = 0,
\end{equation*}
which is equivalent to system
\small
\begin{equation*}
\begin{cases}
e_{22}\left(\chi\left(\frac{D_x}{\langle\eta\rangle}\right)\Big(|D_x|  - \displaystyle\sum_{|\alpha|=1}^N\frac{1}{\alpha!}\partial^\alpha_\eta(\langle\eta\rangle)D^\alpha_x\Big) u_+, -\chi\left(\frac{D_x}{\langle\eta\rangle}\right)\Big(|D_x| + \sum_{|\alpha|=1}^N\frac{1}{\alpha!}\partial^\alpha_\eta(\langle\eta\rangle) D^\alpha_x\Big)u_{-} ; \eta \right)
\\ \hspace{11cm} = -\chi\left(\frac{D_x}{\langle\eta\rangle}\right)e_0 \\
 e_{44}\left(\chi\left(\frac{D_x}{\langle\eta\rangle}\right)\Big(|D_x| + \displaystyle\sum_{|\alpha|=1}^N\frac{1}{\alpha!}\partial^\alpha_\eta(\langle\eta\rangle) D^\alpha_x\Big) u_+, -\chi\left(\frac{D_x}{\langle\eta\rangle}\right)\Big(|D_x| - \sum_{|\alpha|=1}^N\frac{1}{\alpha!}\partial^\alpha_\eta(\langle\eta\rangle) D^\alpha_x\Big)u_{-}; \eta \right)  \\
 \hspace{11cm} = -\chi\left(\frac{D_x}{\langle\eta\rangle}\right) f_0,
\end{cases}
\end{equation*}
\normalsize
with $e_0, f_0$ defined in \eqref{def c0 d0 e0 f0}. Then, if we look for $e_{ii}$ of the form
\begin{equation} \label{form of eii}
e_{ii}(u_+, u_{-};\eta) = \int e^{ix\cdot\xi} \alpha_{ii}(\xi, \eta)\hat{u}_+(\xi) d\xi + \int e^{ix\cdot\xi} \beta_{ii}(\xi, \eta)\hat{u}_{-}(\xi) d\xi,
\end{equation}
this system implies, inter alia, that 
\begin{multline*}
\int e^{ix\cdot\xi}\chi\left(\frac{\xi}{\langle\eta\rangle}\right)\Bigl( |\xi| - \sum_{|\alpha|=1}^N \frac{1}{\alpha!}\partial^\alpha_\eta(\langle\eta\rangle)\xi^\alpha\Bigr)\alpha_{22}(\xi,\eta)\hat{u}_+(\xi)d\xi =\\
 - \frac{i}{4}\int e^{ix\cdot\xi}\chi\left(\frac{\xi}{\langle\eta\rangle}\right)\Bigl(1-\frac{\eta}{\langle\eta\rangle}\cdot\frac{\xi}{|\xi|}\Bigr)\xi_1 \hat{u}_+(\xi) d\xi.
\end{multline*}
As
\begin{equation*}
\left(1 \mp \sum_{|\alpha|=1}^N \frac{1}{\alpha!}\partial^\alpha_\eta (\langle\eta\rangle) \frac{\xi^\alpha}{|\xi|}\right) = 1\mp \sum_{k=1}^N\frac{1}{k!}(\xi\cdot\nabla_\eta)^k(\langle\eta\rangle)
\end{equation*}
and
\begin{equation*}
(\partial_{\eta_1}\xi_1 + \partial_{\eta_2}\xi_2)^k \langle\eta\rangle = \frac{|\xi|^k}{\langle \eta\rangle^{k-1}}\left(1- \left(\frac{\eta}{\langle\eta\rangle}\cdot\frac{\xi}{|\xi|}\right)^2\right)b_k(\xi,\eta), \quad 2\le k\le N,
\end{equation*}
with $b_k(\xi,\eta)$ polynomial of degree $k-2$ in $\frac{\eta}{\langle\eta\rangle}\cdot\frac{\xi}{|\xi|}$, we derive that
\begin{equation}\label{formula_1}
\left(1 \mp \sum_{|\alpha|=1}^N \frac{1}{\alpha!}\partial^\alpha_\eta (\langle\eta\rangle) \frac{\xi^\alpha}{|\xi|}\right) = \left(1\mp \frac{\eta}{\langle\eta\rangle}\cdot\frac{\xi}{|\xi|}\right)\left(1 \mp b_\pm(\xi,\eta)\right)
\end{equation}
with 
\begin{equation} \label{derivatives bpm}
\begin{gathered}
b_\pm(\xi,\eta):= \sum_{k=2}^N \frac{1}{k!}|\xi|^{k-1}\langle\eta\rangle^{-(k-1)}\Big(1\pm \frac{\eta}{\langle\eta\rangle}\cdot\frac{\xi}{|\xi|}\Big)b_k(\xi,\eta),\\
|\partial^\mu_\xi \partial^\nu_\eta b_\pm(\xi,\eta)|\lesssim_{\mu, \nu} |\xi|^{1-|\mu|}\langle\eta\rangle^{-1-|\nu|}, \quad \forall \mu,\nu\in\mathbb{N}^2,
\end{gathered}
\end{equation}
and we can then choose $\alpha_{22}(\xi,\eta)$ in \eqref{form of eii} such that, when $|\xi|\le \varepsilon_2\langle\eta\rangle$,
\begin{equation}
\alpha_{22}(\xi,\eta) =  - \frac{i}{4} \left(1 - b_+(\xi,\eta)\right)^{-1}\frac{\xi_1}{|\xi|}.
\end{equation}
Similarly, we choose multipliers $\beta_{22}, \alpha_{44}, \beta_{44}$ such that, as long as $|\xi|\le \varepsilon_2\langle\eta\rangle$, 
\begin{align*}
&\beta_{22}(\xi, \eta)= \frac{i}{4} \left(1 + b_{-}(\xi,\eta)\right)^{-1}\frac{\xi_1}{|\xi|}, \\
&\alpha_{44}(\xi, \eta) = -\frac{i}{4} \left(1 + b_{-}(\xi,\eta)\right)^{-1}\frac{\xi_1}{|\xi|}, \quad
\beta_{44}(\xi, \eta) = \frac{i}{4} \left(1 - b_+(\xi,\eta)\right)^{-1}\frac{\xi_1}{|\xi|}.
\end{align*}
These multipliers are all well defined for $|\xi|\le \varepsilon_2\langle\eta\rangle$ as $b_\pm(\xi,\eta)=O(|\xi|\langle\eta\rangle^{-1})$.
Moreover, using that $(1\pm b_\mp(\xi,\eta))^{-1}=1\mp b_\mp(\xi,\eta) + O(|\xi|^2\langle\eta\rangle^{-2})$ as long as $|\xi|\le \varepsilon_2\langle\eta\rangle$, we have that
$$\alpha_{22}(\xi,\eta) = -\frac{i}{4} \frac{\xi_1}{|\xi|} + \alpha^{-1}_{22}(\xi, \eta), \quad \beta_{22}(\xi,\eta) = \frac{i}{4} \frac{\xi_1}{|\xi|} + \beta^{-1}_{22}(\xi, \eta),$$
$$\alpha_{44}(\xi,\eta) = -\frac{i}{4} \frac{\xi_1}{|\xi|} + \alpha^{-1}_{44}(\xi, \eta), \quad \beta_{44}(\xi,\eta) = \frac{i}{4} \frac{\xi_1}{|\xi|} + \beta^{-1}_{44}(\xi, \eta),$$
with $|\partial^\mu_\xi\partial^\nu_\eta \alpha^{-1}_{ii}| + |\partial^\mu_\xi\partial^\nu_\eta \beta^{-1}_{ii}| \lesssim_{\mu,\nu} |\xi|^{1-|\mu|}\langle\eta\rangle^{-1-|\nu|}$ for any $\mu,\nu\in\mathbb{N}^2$. 
Injecting the above $\alpha_{ii},\beta_{ii}$, $i\in\{2,4\}$, in \eqref{form of eii} we find that
\begin{align*}
e_{22}\left(\chi\left(\frac{D_x}{\langle\eta\rangle}\right) u_+, \chi\left(\frac{D_x}{\langle\eta\rangle}\right)u_{-};\eta\right) &= -\frac{i}{4}R_1(u_+- u_{-}) + e^{-1}_{22}\left(\chi\left(\frac{D_x}{\langle\eta\rangle}\right) u_+, \chi\left(\frac{D_x}{\langle\eta\rangle}\right)u_{-};\eta\right), \\
e_{44}\left(\chi\left(\frac{D_x}{\langle\eta\rangle}\right) u_+, \chi\left(\frac{D_x}{\langle\eta\rangle}\right)u_{-};\eta\right) &= -\frac{i}{4}R_1(u_+ - u_{-})+ e^{-1}_{44}\left(\chi\left(\frac{D_x}{\langle\eta\rangle}\right) u_+, \chi\left(\frac{D_x}{\langle\eta\rangle}\right)u_{-};\eta\right),
\end{align*}
where, for $i\in\{2,4\}$,
\begin{multline*}
e^{-1}_{ii}\left(\chi\left(\frac{D_x}{\langle\eta\rangle}\right) u_+, \chi\left(\frac{D_x}{\langle\eta\rangle}\right)u_{-};\eta\right) =\\
 \int e^{ix\cdot\xi} \chi\left(\frac{\xi}{\langle\eta\rangle}\right)\alpha^{-1}_{ii}(\xi,\eta) \hat{u}_+(\xi) d\xi +  \int e^{ix\cdot\xi} \chi\left(\frac{\xi}{\langle\eta\rangle}\right)\beta^{-1}_{ii}(\xi,\eta) \hat{u}_{-}(\xi) d\xi.
\end{multline*}
After lemma \ref{Lem_appendix: Kernel with 1 function} $(i)$ and above estimates for $\alpha^{-1}_{ii}, \beta^{-1}_{ii}$, kernels
$$K^i_+(x,\eta):= \int e^{ix\cdot\xi} \chi\left(\frac{\xi}{\langle\eta\rangle}\right)\alpha^{-1}_{ii}(\xi,\eta) \langle\xi\rangle^{-4}d\xi, \quad K^i_{-}(x,\eta):=\int e^{ix\cdot\xi} \chi\left(\frac{\xi}{\langle\eta\rangle}\right)\beta^{-1}_{ii}(\xi,\eta) \langle\xi\rangle^{-4}d\xi$$ 
are such that, for any $\beta\in\mathbb{N}^2$, $|\partial^\beta_\eta K^i_\pm(x,\eta)|\lesssim |x|^{-1}\langle x \rangle^{-2}\langle\eta\rangle^{-1-|\beta|}$ for every $(x,\eta)$. 
This implies that
\small
\begin{multline*}
\left|\partial^\beta_\eta e^{-1}_{ii}\left(\chi\left(\frac{D_x}{\langle\eta\rangle}\right) u_+, \chi\left(\frac{D_x}{\langle\eta\rangle}\right)u_{-};\eta\right)\right| \le \\
\left|\int \partial^\beta_\eta K^i_+(x-y,\eta)[\langle D_x\rangle^4 u_+](y)dy\right| + \left|\int \partial^\beta_\eta K^i_{-}(x-y,\eta) [\langle D_x\rangle^4 u_{-}](y)dy\right| 
\lesssim \|U(t,\cdot)\|_{H^{4,\infty}}\langle\eta\rangle^{-1-|\beta|}
\end{multline*} \normalsize
and $e^{-1}_{ii}$ is a symbol of order $-1$, for $i=2,4$.
Moreover, using definition \eqref{def: seminorm Mmr} and the fact that space $W^{r,\infty}$ injects in $H^{r+1,\infty}$, one can check that for any $r,n\in\mathbb{N}$, 
\begin{equation*}
M^{-1}_r\left(  e^{-1}_{ii}\left(\chi\left(\frac{D_x}{\langle\eta\rangle}\right) u_+, \chi\left(\frac{D_x}{\langle\eta\rangle}\right)u_{-};\eta\right);n \right) \lesssim \|U(t,\cdot)\|_{H^{5+r,\infty}}
\end{equation*}
and therefore that 
\begin{equation*}
M^0_r\left(  e_{ii}\left(\chi\left(\frac{D_x}{\langle\eta\rangle}\right) u_+, \chi\left(\frac{D_x}{\langle\eta\rangle}\right)u_{-};\eta\right);n \right)\lesssim \|\mathrm{R}_1U(t,\cdot)\|_{H^{1+r,\infty}}+\|U(t,\cdot)\|_{H^{5+r,\infty}}.
\end{equation*}
Defining
\begin{equation*}
E^0_d(U;\eta) =
\begin{bmatrix}
0 & 0 & 0 & 0 \\
0 & -\frac{i}{4}R_1(u_+-u_{-}) & 0 & 0 \\
0 & 0 & 0& 0 \\
0 & 0 & 0 & -\frac{i}{4}R_1(u_+-u_{-})
\end{bmatrix}, \quad
E^{-1}_d(U;\eta) = 
\begin{bmatrix}
0 & 0 & 0 & 0 \\
0 & e^{-1}_{22} & 0 & 0\\
0 & 0 & 0 & 0 \\
0 & 0 & 0& e^{-1}_{44}
\end{bmatrix},
\end{equation*}
decomposition \eqref{decomposition Ed} and estimate \eqref{seminorm E0d}, \eqref{seminorm E-1d} hold. Consequently, as
\begin{equation*}
E_d\big(Q_0^{\mathrm{w}}(v_\pm, D_1 v_\pm), Q_0^{\mathrm{w}}(v_\pm, D_1 v_\pm);\eta\big) = E^{-1}_d\big(Q_0^{\mathrm{w}}(v_\pm, D_1 v_\pm), Q_0^{\mathrm{w}}(v_\pm, D_1 v_\pm);\eta\big)
\end{equation*}
for any $n\in\mathbb{N}$ and $\theta\in ]0,1[$, we derive from \eqref{est Hsinfty for NLw-New} with $s=4$ that
\begin{multline*}
M^0_0\left(E_d\big(Q_0^{\mathrm{w}}(v_\pm, D_1 v_\pm), Q_0^{\mathrm{w}}(v_\pm, D_1 v_\pm);\eta\big);n\right) \lesssim \|Q_0^{\mathrm{w}}(v_\pm, D_1 v_\pm)\|_{H^{4,\infty}}\\
\lesssim \|V(t,\cdot)\|^{2-\theta}_{H^{5,\infty}}\|V(t,\cdot)\|^\theta_{H^7},
\end{multline*}
and hence that the quantization of $E_d\big(Q_0^{\mathrm{w}}(v_\pm, D_1 v_\pm), Q_0^{\mathrm{w}}(v_\pm, D_1 v_\pm);\eta\big)$ acting on $V^I_s$ verifies \eqref{L2 est R'(V,V)} after proposition \ref{Prop : Paradiff action on Sobolev spaces-NEW}. Also, \eqref{norm of T-3WI} is deduced from \eqref{norm T-N} while properties \eqref{sigma ij of operator T-3} are obtained using essentially \eqref{Fourier transform composition symbol}.
\endproof
\end{lem}

\begin{lem} \label{Lem: Normal Form on C''nd}
Let $N\in\mathbb{N}^*$.
There exists a purely imaginary matrix $E_{nd}(U;\eta)$, linear in $(u_+, u_{-})$ and of order $-1$, satisfying estimate \eqref{seminorm End}, such that
\begin{multline} \label{equation for Op(End)}
Op^B(C''_{nd}(U;\eta))V^I_s + Op^B(D_t E_{nd}(U;\eta))W^I_s - [A(D),Op^B(E_{nd}(U;\eta))]W^I_s \\ = T_{-N}(U)W^I_s  + \mathfrak{R}'(V,V),
\end{multline}
where $ \mathfrak{R}'(V,V)$ is a remainder satisfying \eqref{L2 est R'(V,V)} and
$T_{-N}(U)$ is a pseudo-differential operator of order less or equal than $-N$ such that, for any $s\in\mathbb{R}$,
\begin{equation}\label{norm_T-N_lemmaCnd}
\|T_{-N}(U)\|_{\mathcal{L}(H^{s-N};H^s)} \lesssim \|U(t,\cdot)\|_{H^{N+6,\infty}}.
\end{equation}
Moreover, its symbol $\sigma(U,\eta) = \left(\sigma_{ij}(U,\eta)\right)_{1\le i,j\le 4}$ is such that
\begin{subequations} \label{sigma ij of operator T-3nd}
\begin{equation}
\mathcal{F}_{x\mapsto\xi}(\sigma_{ij}(U,\eta))(\xi) = 
\begin{cases}
\sigma^+_{ij}(\xi, \eta)\hat{u}_+(\xi) + \sigma^{-}_{ij}(\xi,\eta)\hat{u}_{-}(\xi), \qquad & (i,j)\in \{(2,4), (4,2)\}, \\
 0, &\text{otherwise},
\end{cases}
\end{equation} 
with $\sigma^\pm_{ij}$ supported for $|\xi|\le \varepsilon\langle\eta\rangle$ for a small $\varepsilon>0$, and verifying, for any $\alpha,\beta\in\mathbb{N}^2$,
\begin{equation}
|\partial^\alpha_\xi \partial^\beta_\eta \sigma^\pm_{ij}(\xi,\eta)|\lesssim_{\alpha,\beta}
|\xi|^{N+2-|\alpha|}\langle\eta\rangle^{-N-1-|\beta|}, 
\end{equation}
\end{subequations}
for $(i,j)\in \{(2,4), (4,2)\}$. 
\proof
Because of the structure of $C''_{nd}(U;\eta)$, we seek for a matrix $E_{nd}(U;\eta)$ satisfying \eqref{equation for Op(End)}, of the form $E_{nd}(U;\eta) = (e_{ij})_{1\le i,j\le 4}$ with $e_{ij}=0$  for all $i,j$, except $(i,j)\in\{ (2,4), (4,2)\}$. If we make the ansatz that $E_{nd}(U;\eta)$ is linear in $(u_+,u_{-})$, of order $-1$, and remind that $A(\eta)$ in \eqref{matrices A A'} is of order 1, from symbolic calculus of proposition \ref{Prop: paradifferential symbolic calculus} we have that
\begin{align*}
- [A(D),Op^B(E_{nd}(U;\eta))] =& - Op^B(A(\eta)E_{nd}(U;\eta) - E_{nd}(U;\eta)A(\eta)) \\ 
&- \sum_{|\alpha|=1}^N \frac{1}{\alpha !}Op^B(\partial^\alpha_\eta A(\eta) \cdot D^\alpha_x E_{nd}(U;\eta)) 
+ T_{-N}(U),
\end{align*}
where $T_{-N}(U)$ is a pseudo-differential operator of order less or equal than $-N$, such that, for any $s\in\mathbb{R}$,\small
\begin{equation} \label{norm Tnd-N}
\|T_{-N}(U)\|_{\mathcal{L}(H^{s-N}; H^s)} \lesssim M^1_{N+1}(A(\eta);N+3)M^{-1}_0(E_{nd}(U;\eta);2) + M^1_{0}(A(\eta);N+3)M^{-1}_{N+1}(E_{nd}(U;\eta);2),
\end{equation}\normalsize
and whose symbol $\sigma(U,\eta)=(\sigma_{ij}(U,\eta))_{ij}$ is such that $\sigma_{ij}=0$ for all $i,j$ but $(i,j)\in\{(2,4), (4,2)\}$.
Hence, for any fixed $\chi\in \mathbb{R}^2$ equal to 1 in $\overline{B_{\varepsilon_1}(0)}$ and supported in $B_{\varepsilon_2}(0)$, for some $0<\varepsilon_1<\varepsilon_2\ll 1$, we look for $E_{nd}(U;\eta)$ such that
\begin{multline} \label{equation for End}
\chi\left(\frac{D_x}{\langle\eta\rangle}\right)\Big[C''_{nd}(U;\eta) + D_tE_{nd}(U;\eta) - A(\eta)E_{nd}(U;\eta) + E_{nd}(U;\eta)A(\eta) 
\\ - \sum_{|\alpha|=1}^N \frac{1}{\alpha!} \partial^\alpha_\eta A(\eta)\cdot D^\alpha_x E_{nd}(U;\eta)\Big] = 0.
\end{multline}
Furthermore, as $E_{nd}(U;\eta) = E_{nd}(u_+, u_{-};\eta)$ is linear in $(u_+, u_{-})$ and $u_+$ (resp. $u_{-}$) is solution to the first (resp. the third) equation in \eqref{system for uI+-, vI+-} with $|I|=0$, we have that
\begin{align*}
D_t E_{nd}(u_+, u_{-};\eta) &= E_{nd}(|D_x| u_+, - |D_x|u_{-};\eta) + E_{nd}\big(Q^{\mathrm{w}}_0(v_\pm, D_1 v_\pm), Q^{\mathrm{w}}_0(v_\pm, D_1 v_\pm);\eta\big), \\
D^\alpha_x E_{nd}(u_+, u_{-};\eta) & = E_{nd}(D^\alpha_x u_+, D^\alpha_x u_{-};\eta), \quad \forall \alpha\in\mathbb{N}^2
\end{align*}
while
\begin{equation*}
- A(\eta)E_{nd}(U;\eta) + E_{nd}(U;\eta)A(\eta)  = 
\begin{bmatrix}
0 & 0 & 0 & 0 \\
0 & 0 & 0 & -2\langle\eta\rangle e_{24} \\
0 & 0 & 0 & 0 \\
0 & 2\langle\eta\rangle e_{42} & 0 & 0
\end{bmatrix}.
\end{equation*}
Then we rather search for symbols $e_{24}$ and $e_{42}$ such that
\small
\begin{equation*}
\begin{cases}
& \chi\left(\frac{D_x}{\langle\eta\rangle}\right)  e_{2,4}\left( \Big(|D_x| - \displaystyle\sum_{|\alpha|=1}^N \frac{1}{\alpha!}\partial^\alpha(\langle\eta\rangle) D^\alpha_x - 2\langle\eta\rangle\Big) u_+, -  \Big(|D_x|+ \sum_{|\alpha|=1}^N\frac{1}{\alpha!}\partial^\alpha(\langle\eta\rangle)  D^\alpha_x + 2\langle\eta\rangle\Big)u_{-};\eta\right)\\
& \hspace{11cm} = -\chi\left(\frac{D_x}{\langle\eta\rangle}\right)f_0, \\
& \chi\left(\frac{D_x}{\langle\eta\rangle}\right)e_{4,2}\left( \Big(|D_x| + \displaystyle\sum_{|\alpha|=1}^N\frac{1}{\alpha!}\partial^\alpha(\langle\eta\rangle) D^\alpha_x + 2\langle\eta\rangle\Big) u_+, - \Big(|D_x| - \sum_{|\alpha|=1}^N\frac{1}{\alpha!}\partial^\alpha(\langle\eta\rangle)D^\alpha_x - 2\langle\eta\rangle\Big)u_{-};\eta\right)\\
& \hspace{11cm} = -\chi\left(\frac{D_x}{\langle\eta\rangle}\right)e_0,
\end{cases}
\end{equation*}\normalsize
with $e_0, f_0$ given by \eqref{def c0 d0 e0 f0}, neglecting contribution $E_{nd}\big(Q^{\mathrm{w}}_0(v_\pm, D_1 v_\pm), Q^{\mathrm{w}}_0(v_\pm, D_1 v_\pm);\eta\big)$ whose quantization acting on $W^I_s$ gives rise to a remainder $\mathfrak{R}'(V,V)$, as we will see at the end of the proof.
We look for $e_{ij}$ of the form
$$e_{ij}(u_+, u_{-};\eta) = \int e^{ix\cdot\xi} \alpha_{ij}(\xi, \eta)\hat{u}_+(\xi) d\xi + \int e^{ix\cdot\xi} \beta_{ij}(\xi, \eta)\hat{u}_{-}(\xi) d\xi,$$
for $(i,j)\in\{(2,4), (4,2)\}$, and reminding \eqref{formula_1}, \eqref{derivatives bpm} we choose the above multipliers such that, as long as $|\xi|\le \varepsilon_2\langle\eta\rangle$,
\begin{align*}
& \alpha_{24}(\xi,\eta)= -\frac{i}{4}\left(1+ \frac{\eta}{\langle\eta\rangle}\cdot\frac{\xi}{|\xi|}\right)\left(\left(1- \frac{\eta}{\langle\eta\rangle}\cdot\frac{\xi}{|\xi|}\right)(1-b_+(\xi,\eta)) - 2\frac{\langle\eta\rangle}{|\xi|}\right)^{-1} \frac{\xi_1}{|\xi|},
 \\
& \beta_{24}(\xi,\eta)=-\frac{i}{4} \left(1-\frac{\eta}{\langle\eta\rangle}\cdot\frac{\xi}{|\xi|}\right)\left(\left(1+ \frac{\eta}{\langle\eta\rangle}\cdot\frac{\xi}{|\xi|}\right)(1+b_{-}(\xi,\eta)) + 2\frac{\langle\eta\rangle}{|\xi|}\right)^{-1} \frac{\xi_1}{|\xi|} ,
\end{align*}
\begin{equation*}
\alpha_{42}(\xi,\eta) =\beta_{24}, \quad \beta_{42}(\xi,\eta) =\alpha_{24}(\xi,\eta).
\end{equation*} 
One can check that, on the support of $\chi\big(\frac{\xi}{\langle\eta\rangle}\big)$ and for any $\mu,\nu\in\mathbb{N}^2$, $|\partial^\mu_\xi \partial^\nu_\eta \alpha_{ij}|+|\partial^\mu_\xi \partial^\nu_\eta \beta_{ij}|\lesssim_{\mu,\nu}|\xi|^{1-|\mu|}\langle\eta\rangle^{-1-|\nu|}$,
and then that, if
\begin{equation*}
K^{ij}_+(x,\eta):=\int e^{ix\cdot\eta} \chi\Big(\frac{\xi}{\langle\eta\rangle}\Big)\alpha_{ij}(\xi,\eta)\langle\xi\rangle^{-4} d\xi , \quad K^{ij}_{-}(x,\eta):=\int e^{ix\cdot\eta} \chi\Big(\frac{\xi}{\langle\eta\rangle}\Big)\beta_{ij}(\xi,\eta)\langle\xi\rangle^{-4} d\xi,
\end{equation*}
for $(i,j)\in \{(2,4), (4,2)\}$, $|\partial^\beta_\eta K^{ij}_\pm(x,\eta)|\lesssim |x|^{-1}\langle x \rangle^{-2}\langle\eta\rangle^{-1-|\beta|}$, for any $\beta\in\mathbb{N}^2$ and $(x,\eta)\in\mathbb{R}^2\times\mathbb{R}^2$, as a consequence of lemma \ref{Lem_appendix: Kernel with 1 function}. Therefore
\small
\begin{multline*}
\left|\partial^\beta_\eta e_{ij}\left(\chi\left(\frac{D_x}{\langle\eta\rangle}\right) u_+, \chi\left(\frac{D_x}{\langle\eta\rangle}\right)u_{-};\eta\right)\right| \le \\
\left|\int \partial^\beta_\eta K^{ij}_+(x-y,\eta)[\langle D_x\rangle^4 u_+](y)dy\right| + \left|\int \partial^\beta_\eta K^{ij}_{-}(x-y,\eta) [\langle D_x\rangle^4 u_{-}](y)dy\right| 
\lesssim \|U(t,\cdot)\|_{H^{4,\infty}}\langle\eta\rangle^{-1-|\beta|},
\end{multline*} \normalsize
which implies that $e_{24}, e_{42}$ are symbols of order $-1$. Also, for $(i,j)\in\{(2,4), (4,2)\}$ and any $n,r\in\mathbb{N}$, one can prove that
\begin{equation*}
M^{-1}_r\left(e_{ij}\left(\chi\left(\frac{D_x}{\langle\eta\rangle}\right) u_+, \chi\left(\frac{D_x}{\langle\eta\rangle}\right)u_{-};\eta\right);n\right) \lesssim \|U(t,\cdot)\|_{H^{5+r,\infty}}
\end{equation*}
using definition \eqref{def: seminorm Mmr} and the fact that space $W^{r,\infty}$ injects in $H^{r+1}$ for any $r\in\mathbb{N}$.
Estimate \eqref{norm_T-N_lemmaCnd} follows from \eqref{norm Tnd-N} and symbol $\sigma(U;\eta)$ associated to $T_{-N}(U)$ satisfies \eqref{sigma ij of operator T-3nd}, as one can check using \eqref{Fourier transform composition symbol} and the estimates derived above for $\alpha_{ij},\beta_{ij}$.
Finally, from \eqref{est Hsinfty for NLw-New} with $s=4$ we deduce that, for any $\theta \in ]0,1[$,
$$M^{-1}_0\big(E_{nd}\big(Q^{\mathrm{w}}_0(v_\pm, D_1 v_\pm), Q^{\mathrm{w}}_0(v_\pm, D_1 v_\pm);\eta\big);n\big)\lesssim \|V(t,\cdot)\|^{2-\theta}_{H^{5,\infty}}\|V(t,\cdot)\|^\theta_{H^7},$$ 
and the quantization of $E_{nd}\big(Q^{\mathrm{w}}_0(v_\pm, D_1 v_\pm), Q^{\mathrm{w}}_0(v_\pm, D_1 v_\pm)$ acting on $W^I_s$ is a remainder verifying \eqref{L2 est R'(V,V)} by proposition \ref{Prop : Paradiff action on Sobolev spaces-NEW}. 
\endproof
\end{lem}

\proof[Proof of Proposition \ref{Prop: a first normal form}]
Lemmas \ref{Lem: Normal Forms on C''d(U,eta)} and \ref{Lem: Normal Form on C''nd} show that there exist two matrices $E_d(U;\eta)$ and $E_{nd}(U;\eta)$, linear in $(u_+, u_{-})$, satisfying equations \eqref{equation for Op(Ed)} and \eqref{equation for Op(End)} respectively. After definition \eqref{def_WtildeI} of $\widetilde{W}^I_s$ and equalities \eqref{equation WIs_1}, \eqref{equation for Op(Ed)} and \eqref{equation for Op(End)} we deduce that 
\begin{align*}
& (D_t - A(D))\widetilde{W}^I_s =  Op^B(\widetilde{A}_1(V;\eta))W^I_s + Op^B(A''(V^I;\eta))U  + Op^B_R(A''(V^I;\eta))U  \\
& + Q^I_0(V,W) + \mathfrak{R}(U,V) + Op^B(E(U;\eta))\Big[ Op^B(\widetilde{A}_1(V;\eta))W^I_s  + Op^B(A''(V^I;\eta))U  \\
& + Op^B(C''(U;\eta))V^I_s +  Op^B_R(A''(V^I;\eta))U   + Q^I_0(V,W) + \mathfrak{R}(U,V)\Big] + T_{-N}(U)W^I_s + \mathfrak{R}'(V,V)
\end{align*}
where $\mathfrak{R}(U,V)$ satisfies \eqref{L2 est of R(U,V)}, $\mathfrak{R}'(V,V)$ satisfies \eqref{L2 est R'(V,V)}, and $T_{-N}(U)$ is a pseudo-differential operator of order less or equal than $-N$ verifying \eqref{sigmaij in proposition Normal Forms}, \eqref{norm of T-3 in propositon}.
Contribution
\begin{multline*}
 Op^B(E(U;\eta))\Big[  Op^B(A''(V^I;\eta))U   + Op^B(C''(U;\eta))V^I_s + Op^B_R(A''(V^I;\eta))U\\ + Q^I_0(V,W) + \mathfrak{R}(U,V)\Big] 
\end{multline*}
is a remainder of the form $\mathfrak{R}'(U,V)$ satisfying estimate \eqref{L2 norm R'(U,V)} as a consequence of proposition \ref{Prop : Paradiff action on Sobolev spaces-NEW}, 
estimates \eqref{seminorms E} with $r=0$, lemma \ref{Lemma: L2 estimate of semilinear terms}, and the fact that the $L^2$ norms of $V^I_s$ and $V^I$ are equivalent as long as $\|V(t,\cdot)\|_{H^{2,\infty}}$ is small. 

According to the definition of $E(U;\eta)$ and decomposition \eqref{decomposition Ed}
\begin{align*}
Op^B(E(U;\eta))Op^B(\widetilde{A}_1(V;\eta)) &= Op^B(E^0_d(U;\eta))Op^B(\widetilde{A}_1(V;\eta)) \\&
+ Op^B\big(E^{-1}_d(U;\eta) +E_{nd}(U;\eta) \big)Op^B(\widetilde{A}_1(V;\eta)).
\end{align*}
Proposition \ref{Prop : Paradiff action on Sobolev spaces-NEW} and estimates \eqref{seminorm Atilde 1}, \eqref{seminorm E-1d}, \eqref{seminorm End} with $r=0$, imply that the latter addend in the above right hand side is a bounded operator on $L^2$ whose $\mathcal{L}(L^2)$ norm is estimated by $\|U(t,\cdot)\|_{H^{5,\infty}}\|V(t,\cdot)\|_{H^{1,\infty}}$. The former one writes instead as $Op^B(E^0_d(U;\eta)\widetilde{A}_1(V;\eta))+ T_0(U,V)$, for an operator $T_0(U,V)$ of order less or equal than 0 and $\mathcal{L}(L^2)$ norm controlled by $\|\mathrm{R}_1U(t,\cdot)\|_{H^{2,\infty}}\|V(t,\cdot)\|_{H^{2,\infty}}$, as follows from corollary \ref{Cor : paradiff symbolic calculus at order 1} and estimates \eqref{seminorm Atilde 1}, \eqref{seminorm E0d} with $r=1$. Hence, 
\begin{equation*}
Op^B(E(U;\eta))Op^B(\widetilde{A}_1(V;\eta))W^I_s = Op^B(E^0_d(U;\eta)\widetilde{A}_1(V;\eta))W^I_s + \mathfrak{R}'(U,V),
\end{equation*}
for a new $\mathfrak{R}'(U,V)$ satisfying \eqref{L2 norm R'(U,V)}.

After \eqref{seminorm E0d} matrix $I_4+ E^0_d\big(\chi(\frac{D_x}{\langle\eta\rangle})U;\eta\big)$ is invertible as long as $\|\mathrm{R}_1U(t,\cdot)\|_{H^{1,\infty}}$ is small and
$F^0_d(U;\eta):=\big[I_4+ E^0_d\big(\chi(\frac{D_x}{\langle\eta\rangle})U;\eta\big)\big]^{-1}-I_4$ is such that, for any $n,r\in\mathbb{N}$,
\begin{equation*}
M^0_r\left(F^0_d\Bigl(\chi\Bigl(\frac{D_x}{\langle\eta\rangle}\Bigr)U;\eta\Bigr);n\right)\lesssim \|\mathrm{R}_1U(t,\cdot)\|_{H^{1+r,\infty}}.
\end{equation*}
Moreover, $F^0_d(U;\eta)$ is a real diagonal matrix of order 0, and by corollary \ref{Cor : paradiff symbolic calculus at order 1} with $r=1$
\begin{equation*}
Op^B(I_4 + F^0_d(U;\eta)) Op^B(I_4+E^0_d(U;\eta))= Id + T_{-1}(U),
\end{equation*}
with $T_{-1}(U)$ of order less or equal than 0 and $\mathcal{L}(H^{s-1};H^s)$ norm bounded by $\|\mathrm{R}_1U(t,\cdot)\|_{H^{2,\infty}}$, for any $s\in\mathbb{R}$. This implies that
\begin{equation*}
\begin{gathered}
Op^B(I_4+F^0_d(U;\eta))\widetilde{W}^I_s = W^I_s + \widetilde{T}_{-1}(U)W^I_s,
\widetilde{T}_{-1}(U)=T_{-1}(U) + Op^B(E^{-1}_d(U;\eta) + E_{nd}(U;\eta))
\end{gathered}
\end{equation*}
with $\widetilde{T}_{-1}(U)$ of order less or equal than $-1$ and 
\begin{equation} \label{norm_Ttilde-1U}
\|\widetilde{T}_{-1}(U)\|_{\mathcal{L}(H^{s-1};H^s)}\lesssim \|\mathrm{R}_1U(t,\cdot)\|_{H^{2,\infty}}+ \|U(t,\cdot)\|_{H^{5,\infty}}
\end{equation}
for any $s\in\mathbb{R}$. Hence, as long as this quantity is small, there exists a positive constant $C$ such that \eqref{equivalence WtildeIs WIs} holds. Also,
\begin{align*}
Op^B(I_4+E^0_d(U;\eta))&Op^B(\widetilde{A}_1(V;\eta))W^I_s \\
&= Op^B(I_4+E^0_d(U;\eta))Op^B(\widetilde{A}_1(V;\eta)) Op^B(I_4+F^0_d(U;\eta))\widetilde{W}^I_s\\
 &- Op^B(I_4+E^0_d(U;\eta))Op^B(\widetilde{A}_1(V;\eta))\widetilde{T}_{-1}(U)W^I_s,
\end{align*}
where from proposition \ref{Prop : Paradiff action on Sobolev spaces-NEW}, \eqref{seminorm Atilde 1}, \eqref{seminorm E0d}, \eqref{norm_Ttilde-1U} and \eqref{equivalence WIs WI} the $L^2$ norm of the latter term in the above right hand side is estimated by
\begin{equation} \label{est_norm_Ttilde0}
\|V(t,\cdot)\|_{H^{1,\infty}}\left(\|\mathrm{R}_1U(t,\cdot)\|_{H^{2,\infty}}+\|U(t,\cdot)\|_{H^{5,\infty}}\right)\|W^I(t,\cdot)\|_{L^2}.
\end{equation}
On the other hand, by corollary \ref{Cor : paradiff symbolic calculus at order 1} with $r=1$ we get that
\begin{align*}
 Op^B(I_4+E^0_d(U;\eta))Op^B(\widetilde{A}_1(V;\eta)) &Op^B(I_4+F^0_d(U;\eta))\widetilde{W}^I_s \\
 &= Op^B\big((I_4+E^0_d(U;\eta))\widetilde{A}_1(V;\eta)(I_4+F^0_d(U;\eta))\big)\widetilde{W}^I_s\\
 & +  Op^B(I_4+E^0_d(U;\eta))T_0(U,V)\widetilde{W}^I_s + \widetilde{T}_0(U,V)\widetilde{W}^I_s,
\end{align*}
with $T_0(U,V), \widetilde{T}_0(U,V)$ operators of order less or equal than 0 and $\mathcal{L}(L^2)$ norm controlled, respectively, by $\|\mathrm{R}_1U(t,\cdot)\|_{H^{2,\infty}}\|V(t,\cdot)\|_{H^{2,\infty}}$ and \eqref{est_norm_Ttilde0}, so the last two terms in above right hand side are also remainders $\mathfrak{R}'(U,V)$ by proposition \ref{Prop : Paradiff action on Sobolev spaces-NEW} and estimate \eqref{seminorm E0d}.
That concludes the proof of the statement.
\endproof

\subsection{A second normal forms transformation.} \label{sub: second normal form}

In proposition \ref{Prop: a first normal form} in previous subsection we showed that one can get rid of the slow-decaying-in-time semi-linear contribution $Op^B\big(C''(U;\eta)\big)V^I_s$ in \eqref{equation WIs} by introducing a new function $\widetilde{W}^I_s$, defined in \eqref{def_WtildeI} in terms of $W^I_s$ and solution to equation \eqref{equation Wtilde-Is}.
That naturally led us to the introduction of new energies $\widetilde{E}_n(t;W)$, for $n\in\mathbb{N}, n\ge 3$, and $\widetilde{E}^k_3(t;W)$, for $k\in\mathbb{N}, 0\le k\le 2$, (see \eqref{energy E_tilde_n}) which are respectively equivalent to starting $E_n(t;W)$ and $E^k_3(t;W)$ whenever some uniform norms of $U,V$ are sufficiently small.
However, these new energies do not allow us yet to recover enhanced estimates \eqref{est: bootstrap enhanced Enn} and \eqref{est: boostrap enhanced E02} as it is not true that
\begin{equation} \label{wished_energy_estimate}
\left|\partial_t \widetilde{E}_n(t;W)\right| = O\left(\varepsilon t^{-1+\frac{\delta}{2}}E_n(t;W)^\frac{1}{2}\right), \quad
\left|\partial_t \widetilde{E}^k_3(t;W)\right| = O\left(\varepsilon t^{-1+\frac{\delta}{2}}E^k_3(t;W)^\frac{1}{2}\right).
\end{equation}
This is do to the fact that we still have to deal with semi-linear slow-decaying contributions $Op^B(A''(V^I;\eta))U$, $Op^B_R(A''(V^I;\eta))U$, $Q^I_0(V,W)$ to the right hand side of \eqref{equation Wtilde-Is}, together with the new $T_{-N}(U)W^I_s$ whose $L^2$ norm is also a $O(t^{-\frac{1}{2}}\|W^I(t,\cdot)\|_{L^2})$ after \eqref{norm of T-3 in propositon} and \eqref{est: bootstrap upm}.
The aim of the current subsection is hence to perform a new normal form argument to replace the mentioned terms with more decaying ones.
This is actually done at the energy level, meaning that we are going to add some suitable cubic perturbations to $\widetilde{E}_n(t;W)$ and $\widetilde{E}^k_3(t;W)$ so that the new energies so defined satisfy estimates as in \eqref{wished_energy_estimate}.

Let us first focus on the slow decaying terms that appear when computing
\[\partial_t \widetilde{E}_n(t;W) = \sum_{I\in\mathcal{I}_n}\left\langle \partial_t \widetilde{W}^I_s, \widetilde{W}^I_s\right\rangle\]
for any integer $n\ge 3$.
Using equation \eqref{equation Wtilde-Is} and rewriting $\widetilde{W}^I_s$ in terms of $W^I$ we find, on the one hand, the contribution
\begin{equation} \label{term_getridoff_Etilden}
-\sum_{I\in \mathcal{I}_n}\Im\Bigl[\langle Op^B(A''(V^I;\eta))U + Op^B_R(A''(V^I;\eta))U , W^I\rangle + \langle T_{-N}(U)W^I, W^I\rangle\Bigr],
\end{equation}
which is a $O(\varepsilon t^{-1/2}E_n(t;W))$ after Cauchy-Schwarz inequality, lemma \ref{Lemma: L2 estimate of semilinear terms} and a-priori estimates \eqref{est: bootstrap argument a-priori est}.
But we also have
\begin{equation} \label{contribution-to-get-ridoff}
-\sum_{I\in \mathcal{I}_n}\sum_{\substack{(I_1,I_2)\in \mathcal{I}(I) \\ [\frac{|I|}{2}]<|I_1|<|I|}}\Im \left[\langle Q^{\mathrm{kg}}_0(v^{I_1}_\pm, D_1u^{I_2}_\pm), v^I_+ + v^I_{-}\rangle \right],
\end{equation}
which enjoys the same decay as the previous one, as can be immediately seen using again Cauchy-Schwarz inequality along with \eqref{est: L2 Qkg0 (vI1 vI2)-only derivatives} and \eqref{est: bootstrap upm}.
From definition \eqref{matrix A''(VI)} of matrix $A''(V^I,\eta)$, Plancherel's formula, \eqref{fourier transform of paradiff op} and the fact that $\overline{v^I_+} = -v^I_{-}$
\begin{align*}
&\langle Op^B(A''(V^I;\eta))U, W^I \rangle  = \langle Op^B(a_0(v^I_\pm;\eta)\eta_1)u_+ + Op^B(b_0(v^I_\pm;\eta)\eta_1)u_{-}, v^I_+ + v^I_{-}\rangle  \\
& = 
-\frac{i}{4(2\pi)^2} \int \chi\left(\frac{\xi - \eta}{\langle\eta\rangle}\right) \left[(\reallywidehat{v_+^{I}+v_{-}^{I}})(\xi-\eta)(\reallywidehat{u_+ + u_{-}})(\eta) - \frac{\xi-\eta}{\langle\xi-\eta\rangle}\cdot\frac{\eta}{|\eta|}(\reallywidehat{v_+^{I}-v_{-}^{I}})(\xi-\eta) \right. \\
& \left. \hspace{8,5cm} \times (\reallywidehat{u_+ - u_{-}})(\eta)\right]\eta_1( \reallywidehat{v^I_{-} + v^I_{-}})(-\xi) d\xi d\eta,
\end{align*}
with $\chi$ denoting a smooth function equal to 1 in $\overline{B_{\varepsilon_1}(0)}$ and supported in $B_{\varepsilon_2}(0)$, for some $0<\varepsilon_1<\varepsilon_2\ll 1$. Hence
\begin{equation} \label{sum of CIk}
 - \Im\left[\langle Op^B(A''(V^I;\eta))U, W^I \rangle\right]=  \sum_{j_k\in\{+,-\}} C^I_{(j_1,j_2,j_3)}
\end{equation}
with
\begin{equation} \label{integral_CI}
C^I_{(j_1,j_2,j_3)}= \frac{1}{4(2\pi)^2}\displaystyle\int \chi\Big(\frac{\xi-\eta}{\langle\eta\rangle}\Big)\left(1-j_1j_2\frac{\xi-\eta}{\langle\xi-\eta\rangle}\cdot\frac{\eta}{|\eta|}\right)\eta_1 \hat{v}^I_{j_1}(\xi-\eta) \hat{u}_{j_2}(\eta) \hat{v}^I_{j_3}(-\xi) d\xi d\eta,
\end{equation}
for any $j_1,j_2,j_3\in\{+,-\}$.
Analogously, from equality \eqref{fourier transform of R(v,w)}
\begin{equation} \label{sum of CIRk}
 - \Im\left[\langle Op^B_R(A''(V^I;\eta))U, W^I \rangle\right]=\sum_{j_k\in\{+,-\}} C^{I,R}_{(j_1,j_2,j_3)}
\end{equation}
with
\begin{multline} \label{integral_CIR}
C^{I,R}_{(j_1,j_2,j_3)}= \frac{1}{4(2\pi)^2}\displaystyle\int \left[1-\chi\Big(\frac{\xi-\eta}{\langle\eta\rangle}\Big) - \chi\Big(\frac{\eta}{\langle\xi-\eta\rangle}\Big)\right]\left(1-j_1j_2\frac{\xi-\eta}{\langle\xi-\eta\rangle}\cdot\frac{\eta}{|\eta|}\right)\eta_1 \\
\times \hat{v}^I_{j_1}(\xi-\eta) \hat{u}_{j_2}(\eta) \hat{v}^I_{j_3}(-\xi) d\xi d\eta.
\end{multline}
After proposition \ref{Prop: a first normal form}, $T_{-N}(U)=(\sigma_{ij}(U,D_x))_{ij}$ with symbols $\sigma_{ij}(U,\eta)$ satisfying \eqref{sigmaij in proposition Normal Forms}. 
Introducing $\rho:\{+,-\}\rightarrow \{2,4\}$ such that $\rho(+)=2, \rho(-) =4$ and using the convention that $-j_k\in \{+,-\}\setminus \{j_k\}$, we have that
\begin{equation} \label{expression T-N in normal forms}
\begin{split}
\langle T_{-N}(U)W^I, W^I\rangle & = \sum_{i,j\in\{+,-\}} \langle \sigma_{\rho(i) \rho(j)}(U,D_x) v^I_j, v^I_i\rangle \\
& = - \frac{1}{(2\pi)^2}\sum_{j_k\in \{+,-\}}\int \sigma^{j_2}_{\rho(j_3), \rho(j_1)}(\eta, \xi-\eta) \hat{v}^I_{j_1}(\xi-\eta)\hat{u}_{j_2}(\eta)\hat{v}^I_{-j_3}(-\xi)d\xi d\eta ,
\end{split}
\end{equation}
where multipliers $\sigma^{j_2}_{\rho(j_3), \rho(j_1)}(\eta, \xi-\eta)$ are supported for $|\eta|\le\varepsilon |\xi-\eta|$ and such that, for any $\alpha,\beta\in\mathbb{N}^2$,
\begin{equation*}
\left|\partial^\alpha_\xi \partial^\beta_\eta \sigma^{j_2}_{\rho(j_3), \rho(j_1)}(\eta, \xi-\eta)\right|\lesssim_{\alpha,\beta} |\eta|^{N+1-|\beta|}\langle \xi-\eta\rangle^{-N-|\alpha|}, 
\end{equation*}
for any $(\xi,\eta)\in\mathbb{R}^2\times\mathbb{R}^2$, any $j_1,j_2,j_3\in \{+,-\}$.
Moreover, by \eqref{Q0_pm} we have that
\begin{gather} \label{sum CI1I2}
-\Im \left[\langle Q^{\mathrm{kg}}_0(v^{I_1}_\pm, D_1u^{I_2}_\pm), v^I_+ + v^I_{-}\rangle \right] = \sum_{j_k\in\{+,-\}} C^{I_1,I_2}_{(j_1,j_2,j_3)}
\end{gather}
with
\begin{equation} \label{integral_CI1I2}
C^{I_1,I_2}_{(j_1,j_2,j_3)}:= \frac{1}{4(2\pi)^2}\displaystyle\int \left(1-j_1j_2\frac{\xi-\eta}{\langle\xi-\eta\rangle}\cdot\frac{\eta}{|\eta|}\right)\eta_1\ \hat{v}^{I_1}_{j_1}(\xi-\eta) \hat{u}^{I_2}_{j_2}(\eta) \hat{v}^I_{j_3}(-\xi) d\xi d\eta.
\end{equation}
The above equalities lead us to introduce the following multipliers
\begin{equation} \label{def of B(i1,i2,i3)} 
B^k_{(j_1,j_2,j_3)}(\xi,\eta):= \frac{1}{j_1\langle\xi-\eta\rangle+j_2|\eta| + j_3\langle\xi\rangle}\left(1-j_1j_2\frac{\xi-\eta}{\langle\xi-\eta\rangle}\cdot\frac{\eta}{|\eta|}\right)\eta_k, \quad k=1,2
\end{equation}
and
\begin{equation}\label{def_multiplier_sigmatildeN}
\widetilde{\sigma}^N_{(j_1,j_2,j_3)}(\xi,\eta):= \frac{ \sigma^{j_2}_{\rho(j_3), \rho(j_1)}(\eta, \xi-\eta)}{j_1\langle\xi-\eta\rangle + j_2|\eta| -j_3\langle\xi\rangle},
\end{equation}
together with the following integrals
\begin{subequations} \label{integral Dk DR_k DT}
\begin{equation}
D^I_{(j_1,j_2,j_3)} := \frac{i}{4(2\pi)^2}\displaystyle\int  \chi\left(\frac{\xi-\eta}{\langle\eta\rangle}\right) B^1_{(j_1,j_2,j_3)}(\xi, \eta)\hat{v}^I_{j_1}(\xi-\eta) \hat{u}_{j_2}(\eta)\hat{v}^I_{j_3}(-\xi)\, d\xi d\eta,  \label{def of D_k}
\end{equation}
\begin{multline}
D^{I,R}_{(j_1, j_2, j_3)} :=
 \frac{i}{4(2\pi)^2}\displaystyle\int\left[1 - \chi\left(\frac{\xi-\eta}{\langle\eta\rangle}\right) - \chi\left(\frac{\eta}{\langle\xi-\eta\rangle}\right)\right] B^1_{(j_1,j_2,j_3)}(\xi, \eta )\\
 \times \hat{v}^I_{j_1}(\xi-\eta) \hat{u}_{j_2}(\eta)\hat{v}^I_{j_3}(-\xi)d\xi d\eta , \label{def of DR_k}
 \end{multline}
 \begin{equation}
D^{I,T_{-N}}_{(j_1, j_2, j_3)}  := \Re \left[\frac{1}{(2\pi)^2}\int\widetilde{\sigma}^N_{(j_1,j_2,j_3)}(\xi, \eta)\hat{v}^I_{j_1}(\xi-\eta)\hat{u}_{j_2}(\eta)\hat{v}^I_{-j_3}(-\xi)d\xi d\eta\right] \label{def of DIT}
\end{equation}
\end{subequations}
and
\begin{equation} \label{def of DI1I2} 
D^{I_1,I_2}_{(j_1,j_2,j_3)} := \frac{i}{4(2\pi)^2}\displaystyle\int B^1_{(j_1,j_2,j_3)}(\xi, \eta)\hat{v}^{I_1}_{j_1}(\xi-\eta) \hat{u}^{I_2}_{j_2}(\eta)\hat{v}^I_{j_3}(-\xi)\, d\xi d\eta
\end{equation}
for any triplet $(j_1,j_2,j_3)\in \{+,-\}^3$.
We warn the reader that definitions \eqref{integral Dk DR_k DT} and \eqref{def of DI1I2} are given here for any general multi-indices $I, I_1, I_2$.

\begin{defn}\label{def_Edag_n}
For every integer $n\ge 3$ we define the second modified energy $\widetilde{E}^\dagger_n(t;W)$ as
\begin{multline} \label{energy_dag_En}
\widetilde{E}^\dagger_n(t;W):= \widetilde{E}_n(t;W)+ \sum_{\substack{I\in\mathcal{I}_n \\ j_i\in\{+,-\}}} \left(D^I_{(j_1,j_2,j_3)} + D^{I,R}_{(j_1,j_2,j_3)} + D^{I,T_{-N}}_{(j_1,j_2,j_3)}\right)\\
 + \sum_{\substack{I\in\mathcal{I}_n \\ j_i\in\{+,-\}}}\sum_{\substack{(I_1,I_2)\in\mathcal{I}(I)\\ [\frac{|I|}{2}]<|I_1|<|I|}} D^{I_1,I_2}_{(j_1,j_2,j_3)}.
\end{multline}
\end{defn}

Let us now analyse the time derivative of $\widetilde{E}^k_3(t;W)$ for integers $0\le k\le 2$.
As in the previous case, from equation \eqref{equation Wtilde-Is} we see appear the same contribution as in \eqref{term_getridoff_Etilden}, but with the sum over $\mathcal{I}_n$ replaced with that on $\mathcal{I}^k_3$.
We also find 
\begin{equation} \label{term_getridoff_QI0}
-\sum_{I\in\mathcal{I}^k_3}\Im[\langle Q^I_0(V,W), W^I\rangle]
\end{equation}
which is a $O(\varepsilon t^{-(1+\delta_k)/2}E^k_3(t;W)^{1/2})$ from Cauchy-Schwarz inequality and estimate \eqref{est:QI0-Ik3}.
To be more precise, the slow decay in time of the above scalar product is due to some particular quadratic term appearing in $Q^I_0(V,W)$. 
In fact, according to definition \eqref{matrix QI} and to \eqref{decomposition_Qw0}, \eqref{est:L2_norm_Rk3(t,x)}, \eqref{eq:Qw0(Dt)-statement}, \eqref{est: bootstrap E02}, for any $I\in\mathcal{I}^k_3$
\begin{multline} \label{est_Qw0}
\sum_{\substack{(I_1,I_2)\in\mathcal{I}(I)\\ |I_2|<|I|}} \left|\langle Q^\mathrm{w}_0(v^{I_1}_\pm, D_x v^{I_2}_\pm), u^I_+ + u^I_{-}\rangle\right| + \sum_{\substack{(I_1,I_2)\in\mathcal{I}(I)\\ |I_1|+|I_2|\le 2, |I_2|<|I|}} \left|\langle Q^\mathrm{w}_0(v^{I_1}_\pm, D_t v^{I_2}_\pm), u^I_+ + u^I_{-}\rangle\right| \\
\lesssim \|\mathfrak{R}^k_3(t,\cdot)\|_{L^2}\|U^I(t,\cdot)\|_{L^2} \le C(A+B)\varepsilon t^{-1+\frac{\delta_k}{2}} E^k_3(t;W)^\frac{1}{2}.
\end{multline}
Also, after \eqref{decomposition_Qkg0} and \eqref{eq:Qkg0(Dt)-statement} we have that for all $I\notin \mathcal{V}^k$, with $\mathcal{V}^k$ defined in \eqref{set_V},
\begin{multline} \label{est_Qkg0-Ik3}
\sum_{\substack{(I_1,I_2)\in\mathcal{I}(I)\\ |I_1|,|I_2|<|I|}}\left|\left\langle Q^\mathrm{kg}_0(v^{I_1}_\pm, D_xu^{I_2}_\pm), v^I_+ +v^I_{-}\right\rangle\right| +\sum_{\substack{(I_1,I_2)\in\mathcal{I}(I)\\ |I_1|+|I_2|\le 2, |I_1|,|I_2|<|I|}}\left|\left\langle Q^\mathrm{kg}_0(v^{I_1}_\pm, D_t u^{I_2}_\pm), v^I_+ +v^I_{-}\right\rangle\right|  \\ 
\lesssim \|\mathfrak{R}^k_3(t,\cdot)\|_{L^2}\|V^I(t,\cdot)\|_{L^2} \le C(A+B)\varepsilon t^{-1+\frac{\delta_k}{2}} E^k_3(t;W)^\frac{1}{2}.
\end{multline}
Observe that the decay rate $O(t^{-1+\delta_k/2})$ in the right hand side of the two above inequalities is the slowest one that allows us to propagate a-priori estimate \eqref{est: bootstrap E02} and it gives us back exactly the slow growth in time $t^{\delta_k/2}$ enjoyed by $E^k_3(t;W)^{1/2}$, for $0\le k\le 2$.
On the other hand, for $I\in\mathcal{V}^k$ with $k=0,1,$ we have that, for some smooth cut-off function $\chi$ and some $\sigma>0$ small,
\begin{equation*}
\begin{gathered}
\sum_{\substack{(I_1,I_2)\in\mathcal{I}(I)\\ |I_1|,|I_2|<|I|}}c_{I_1,I_2} Q^\mathrm{kg}_0(v^{I_1}_\pm, D_xu^{I_2}_\pm) = \sum_{\substack{(I_1,I_2)\in\mathcal{I}(I)\\ I_1\in \mathcal{K}, |I_2|\le 1}}c_{I_1,I_2} Q^\mathrm{kg}_0\left( v^{I_1}_\pm,\chi(t^{-\sigma}D_x) D_x u^{I_2}_\pm\right) + \mathfrak{R}^k_3(t,x),\\
\sum_{\substack{(I_1,I_2)\in\mathcal{I}(I)\\|I_1|+|I_2|\le 2, |I_1|,|I_2|<|I|}}c_{I_1,I_2} Q^\mathrm{kg}_0(v^{I_1}_\pm, D_tu^{I_2}_\pm) =\sum_{\substack{(J,0)\in\mathcal{I}(I)\\ J\in\mathcal{K}}} c_{J,0} Q^\mathrm{kg}_0(v^J_\pm,  \chi(t^{-\sigma}D_x) |D_x| u_\pm)+ \mathfrak{R}^k_3(t,x).
\end{gathered}
\end{equation*}
The $L^2$ norms of the summations in the above right hand sides are bounded by
\begin{equation*}
\sum_{J|\le 1}\left( \|\chi(t^{-\sigma} D_x) u^J_\pm(t,\cdot)\|_{H^{2,\infty}} +  \|\chi(t^{-\sigma} D_x)\mathrm{R} u^J_\pm(t,\cdot)\|_{H^{2,\infty}}\right) E^k_3(t;W)^\frac{1}{2}
\end{equation*} 
and hence by $\varepsilon t^{-1/2}E^k_3(t;W)^{1/2}$ as follows by sharp estimate \eqref{Linfty_est_UJ} derived in appendix \ref{Appendix B}.
Therefore, the very contribution to \eqref{term_getridoff_QI0} that has to be eliminated from $\partial_t \widetilde{E}^k_3(t;W)$ appears only for $k=0,1$ and is
\begin{multline} \label{decomposition-Qkg0-complete}
-\sum_{I\in \mathcal{V}^k}\sum_{\substack{(I_1,I_2)\in\mathcal{I}(I)\\ I_1\in\mathcal{K}, |I_2|\le 1}}c_{I_1,I_2} \Im\left[\left\langle Q^\mathrm{kg}_0\left( v^{I_1}_\pm,\chi(t^{-\sigma}D_x) D_x u^{I_2}_\pm\right), v^I_+ + v^I_{-}\right\rangle\right] \\
 -\sum_{I\in \mathcal{V}^k}\sum_{\substack{(J,0)\in\mathcal{I}(I)\\ J\in\mathcal{K}}}c_{J,0} \Im\left[\left\langle Q^\mathrm{kg}_0\left( v^J_\pm,\chi(t^{-\sigma}D_x) |D_x| u_\pm\right), v^I_+ + v^I_{-}\right\rangle\right] .
\end{multline}
As
\begin{equation} \label{sum_FI1I2}
\begin{gathered}
-\Im \left[\langle Q^\mathrm{kg}_0(v^{I_1}_\pm,  \chi(t^{-\sigma}D_x)D_l u^{I_2}_\pm), v^I_+ + v^I_{-} \rangle\right] =\sum_{j_i\in \{+,-\}} F^{I_1,I_2,l}_{(j_1,j_2,j_3)}, \quad l=1,2\\
-\Im \left[\langle Q^\mathrm{kg}_0(v^{I_1}_\pm,  \chi(t^{-\sigma}D_x)|D_x| u^{I_2}_\pm), v^I_+ + v^I_{-} \rangle\right] =\sum_{j_i\in \{+,-\}} F^{I_1,I_2,3}_{(j_1,j_2,j_3)},
\end{gathered}
\end{equation}
with
\begin{equation} \label{integral_FI1I2} 
F^{I_1,I_2, l}_{(j_1,j_2,j_3)} = \frac{1}{4(2\pi)^2}\displaystyle\int \left(1-j_1j_2\frac{\xi-\eta}{\langle\xi-\eta\rangle}\cdot\frac{\eta}{|\eta|}\right)\eta_l\ \hat{v}^{I_1}_{j_1}(\xi-\eta) \reallywidehat{\chi(t^{-\sigma}D_x) u^{I_2}_{j_2}}(\eta) \hat{v}^I_{j_3}(-\xi) d\xi d\eta,
\end{equation}
for any $j_i\in\{+,-\}$, $l=1,2,3$, and $\eta_3:=j_2|\eta|$,
we introduce a new multiplier
\begin{equation}\label{multiplier_B3}
B^3_{(j_1,j_2,j_3)}(\xi,\eta):= \frac{j_2}{j_1\langle\xi-\eta\rangle+j_2|\eta| + j_3\langle\xi\rangle}\left(1-j_1j_2\frac{\xi-\eta}{\langle\xi-\eta\rangle}\cdot\frac{\eta}{|\eta|}\right)|\eta|
\end{equation}
together with integrals
\begin{equation}\label{integral_GI1I2}
G^{I_1,I_2,l}_{(j_1,j_2,j_3)} = \frac{i}{4(2\pi)^2}\displaystyle\int B^l_{(j_1,j_2,j_3)}(\xi,\eta)\ \hat{v}^{I_1}_{j_1}(\xi-\eta) \reallywidehat{\chi(t^{-\sigma}D_x) u^{I_2}_{j_2}}(\eta) \hat{v}^I_{j_3}(-\xi) d\xi d\eta
\end{equation}
for any $l=1,2,3$, $(j_1,j_2,j_3)\in\{+,-\}^3$, with multipliers $B^l_{(j_1,j_2,j_3)}$ given by \eqref{def of B(i1,i2,i3)} when $l=1,2$, and by \eqref{multiplier_B3} when $l=3$.
We warn the reader that in what follows we will sometimes refer to multipliers $B^l_{(j_1,j_2,j_3)}$ (resp. integrals $F^{I_1,I_2, l}_{(j_1,j_2,j_3)}$ and $G^{I_1,I_2,l}_{(j_1,j_2,j_3)}$) simply as $B_{(j_1,j_2,j_3)}$ (resp. $F^{I_1,I_2}_{(j_1,j_2,j_3)}$ and $G^{I_1,I_2}_{(j_1,j_2,j_3)}$) forgetting about superscript $l$.
This choice reveals to be convenient when we do not need to distinguish between $l=1,2,3$.
\begin{defn} \label{def_Edag_k3}
For every integer $0\le k\le 2$ we define the second modified energy $\widetilde{E}^{k,\dagger}_3(t;W)$ as
\begin{multline}\label{energy_dag_Ek2}
\widetilde{E}^{k,\dagger}_3(t;W):= \widetilde{E}^k_3(t;W)+ \sum_{\substack{I\in\mathcal{I}^k_3\\ j_i\in \{+,-\}}}\left(D^I_{(j_1,j_2,j_3)} + D^{I,R}_{(j_1,j_2,j_3)} + D^{I,T_{-N}}_{(j_1,j_2,j_3)}\right) \\
+\delta_{k<2}\sum_{\substack{I\in\mathcal{V}^k\\ j_i\in \{+,-\}}}\sum_{\substack{(I_1,I_2)\in\mathcal{I}(I)\\ I_1\in\mathcal{K}, |I_2|\le 1}}c_{I_1,I_2} G^{I_1,I_2}_{(j_1,j_2,j_3)},
\end{multline}
with $\delta_{k<2}=1$ if $k=0,1,$ 0 otherwise, and coefficients $c_{I_1,I_2}\in \{-1,0,1\}$.
\end{defn}

In view of the lemmas to follow it is useful  to remind that, after system \eqref{system for uI+-, vI+-}, for any multi-index $I$ vector $(\hat{u}^I_+, \hat{v}^I_+, \hat{u}^I_{-}, \hat{v}^I_{-})$ is solution to
\begin{equation} \label{system for hat(u)I+-, hat(v)I+-}
\begin{cases}
& (D_t - |\xi|)\hat{u}^I_+(t,\xi) = \sum_{|I_1| + |I_2| = |I| }\reallywidehat{Q_0^{\mathrm{w}}(v^{I_1}_\pm, D_1 v^{I_2}_\pm)} + \sum_{|I_1| + |I_2| < |I| } c_{I_1, I_2}\reallywidehat{Q_0^{\mathrm{w}}(v^{I_1}_\pm, D v^{I_2}_\pm)}  \\
& (D_t - \langle\xi\rangle)\hat{v}^I_+(t,\xi) =  \sum_{|I_1| + |I_2| = |I| } \reallywidehat{Q_0^{\mathrm{kg}}(v^{I_1}_\pm, D_1 u^{I_2}_\pm)} + \sum_{|I_1| + |I_2| < |I| } c_{I_1, I_2}\reallywidehat{Q_0^{\mathrm{kg}}(v^{I_1}_\pm, D u^{I_2}_\pm)}  \\
& (D_t + |\xi|)\hat{u}^I_{-}(t,x) = \sum_{|I_1| + |I_2| = |I| }\reallywidehat{Q_0^{\mathrm{w}}(v^{I_1}_\pm, D_1 v^{I_2}_\pm)} + \sum_{|I_1| + |I_2| <|I| } c_{I_1, I_2}\reallywidehat{Q_0^{\mathrm{w}}(v^{I_1}_\pm, D v^{I_2}_\pm)}  \\
& (D_t + \langle \xi \rangle)\hat{v}^I_{-}(t,x) = \sum_{|I_1| + |I_2| = |I| } \reallywidehat{Q_0^{\mathrm{kg}}(v^{I_1}_\pm, D_1 u^{I_2}_\pm)} + \sum_{|I_1| + |I_2| < |I| } c_{I_1, I_2}\reallywidehat{Q_0^{\mathrm{kg}}(v^{I_1}_\pm, D u^{I_2}_\pm)}
\end{cases}
\end{equation}
with coefficients $c_{I_1, I_2}\in \{-1,0,1\}$ and indices $I_1,I_2$ in above right hand side such that $(I_1,I_2)\in\mathcal{I}(I)$. 
In lemmas \ref{Lem:Quartic_terms_1} and \ref{Lem:Quartic_Terms_II} we will check that, with definitions \ref{def_Edag_n}, \ref{def_Edag_k3}, the slow decaying contributions highlighted in \eqref{term_getridoff_Etilden} are replaced in $\partial_t \widetilde{E}^\dagger_n(t;W)$, $\partial_t\widetilde{E}^{k,\dagger}_3(t;W)$ by some new quartic terms.
These latter ones are obtained from integrals \eqref{integral Dk DR_k DT} by replacing each factor $\hat{v}^I_{j_1}, \hat{u}_{j_2},\hat{v}^I_{j_3}$ at a time with the non-linearity appearing in the equation that factor satisfies in \eqref{system for hat(u)I+-, hat(v)I+-}.
Lemma \ref{Lem:Quartic_terms_III} (resp. lemma \ref{Lem:Analysis quartic terms IV}) shows that the same is for troublesome contributions \eqref{contribution-to-get-ridoff} in $\partial_t\widetilde{E}^\dagger_n(t;W)$ (resp. for \eqref{decomposition-Qkg0-complete} in $\partial_t\widetilde{E}^{k,\dagger}_3(t;W)$).
We are also going to see that, if $N\in\mathbb{N}^*$ is chosen sufficiently large (e.g. $N=18$), all these quartic terms suitably decay in time, and that modified energies $\widetilde{E}^\dagger_n(t;W), \widetilde{E}^{k,\dagger}_3(t;W)$ are equivalent, respectively, to $E_n(t;W), E^k_3(t,W)$.
We point out the fact that the normal form's step performed in previous section was necessary to avoid here some problematic quartic contributions coming from quasi-linear terms in \eqref{system for hat(u)I+-, hat(v)I+-} and that could lead to some loss of derivatives.
Before proving the mentioned lemmas, we need to introduce two preliminary results, that will be useful in the proof of lemmas \ref{Lem:Quartic_terms_1}, \ref{Lem:Quartic_terms_III}.
\begin{lem} \label{Lem:Est_integrals_quartic-terms}
For any $j_i\in\{+,-\}$, $i=1,2,3$, let $B^k_{(j_1,j_2,j_3)}(\xi,\eta)$ be the multiplier defined in \eqref{def of B(i1,i2,i3)} when $k=1,2$ and in \eqref{multiplier_B3} when $k=3$, and $\psi_1,\psi_2, \psi_3$ be three smooth cut-off functions such that $\psi_1(x)$ is supported for $|x|\le c$, $\psi_2(x)$ is supported for $c'\le |x|\le C'$, $\psi_3(x)$ is supported for $|x|\ge C$, for some $0<c,c'\ll 1$, $C,C'\gg 1$, and $\psi_1+\psi_2+\psi_3 \equiv 1$.
Let also $\delta_k$ be equal to 1 for $k=1,2$, 0 for $k=3$.

$(i)$
For any $j_1,\dots,j_5\in\{+,-\}$, $i=1,2$, and any $u_1,u_2,u_3,u_4$ such that $u_1\in H^{4,\infty}(\mathbb{R}^2)$, $u_2,u_4\in L^2(\mathbb{R}^2)$, $u_3\in H^{11,\infty}(\mathbb{R}^2)$ and $\delta_k\mathrm{R}_ku_3\in H^{7,\infty}(\mathbb{R}^2)$, \small
\begin{multline} \label{ineq:psi1-psi2}
\left| \int \psi_i\Big(\frac{\xi-\eta}{\langle\eta\rangle}\Big)B^k_{(j_1,j_2,j_3)}(\xi,\eta) \left(1-j_4j_5\frac{\xi-\eta-\zeta}{\langle\xi-\eta-\zeta\rangle}\cdot\frac{\zeta}{|\zeta|}\right)\zeta_1 \hat{u}_1(\xi-\eta-\zeta)\hat{u}_2(\zeta)\hat{u}_3(\eta) \hat{u}_4(-\xi) d\xi d\eta d\zeta \right|\\
\lesssim \|u_1\|_{H^{4,\infty}}\|u_2\|_{L^2}\left(\|u_3\|_{H^{11,\infty}}+ \delta_k\|\mathrm{R}_ku_3\|_{H^{7,\infty}}\right)\|u_4\|_{L^2};
\end{multline}\normalsize
$(ii)$ For any $j_1,\dots,j_5\in\{+,-\}$, and any $u_1,u_2,u_3,u_4$ such that $u_1\in H^{7,\infty}(\mathbb{R}^2)$, $u_2\in H^1(\mathbb{R}^2)$, $u_4\in L^2(\mathbb{R}^2)$, $u_3\in H^{4,\infty}(\mathbb{R}^2)$ and $\delta_k\mathrm{R}_k u_3\in L^\infty(\mathbb{R}^2)$,
\small
\begin{multline} \label{ineq:psi3}
\left| \int \psi_3\Big(\frac{\xi-\eta}{\langle\eta\rangle}\Big)B^k_{(j_1,j_2,j_3)}(\xi,\eta) \left(1-j_4j_5\frac{\xi-\eta-\zeta}{\langle\xi-\eta-\zeta\rangle}\cdot\frac{\zeta}{|\zeta|}\right)\zeta_1 \hat{u}_1(\xi-\eta-\zeta)\hat{u}_2(\zeta)\hat{u}_3(\eta) \hat{u}_4(-\xi) d\xi d\eta d\zeta \right|\\
\lesssim \|u_1\|_{H^{7,\infty}}\|u_2\|_{H^1}\left(\|u_3\|_{H^{4,\infty}}+ \delta_k\|\mathrm{R}_ku_3\|_{L^\infty}\right)\|u_4\|_{L^2}.
\end{multline}\normalsize
\proof
Let $k=1,2$. We are going to refer to $B^k_{(j_1,j_2,j_3)}$ (resp. $\eta_k$ and $\mathrm{R}_k$) simply as $B_{(j_1,j_2,j_3)}$ (resp. $\eta$ and $\mathrm{R}$) and rather use a superscript to define a decomposition of this multiplier (see \eqref{def_B0-B1})

Let us observe that, as
\begin{equation*}
B_{(j_1,j_2,j_3)}(\xi,\eta) =\frac{j_1\langle\xi-\eta\rangle+j_2|\eta| - j_3\langle\xi\rangle}{2j_1j_2\langle\xi-\eta\rangle |\eta|}\eta,
\end{equation*}
we can write
\begin{equation}\label{B=B0+B1}
B_{(j_1,j_2,j_3)}(\xi,\eta) = B^0_{(j_1,j_2,j_3)}(\xi,\eta)\frac{\eta}{|\eta|} + B^1_{(j_1,j_2,j_3)}(\xi,\eta)\langle\eta\rangle^4,
\end{equation}
where for any smooth cut-off function $\phi$, equal to 1 in a neighbourhood of the origin,
\begin{equation}\label{def_B0-B1}
\begin{gathered}
B^0_{(j_1,j_2,j_3)}(\xi,\eta) :=\frac{j_1\langle\xi-\eta\rangle+j_2|\eta| - j_3\langle\xi\rangle}{2j_1j_2\langle\xi-\eta\rangle }\phi(\eta),\\
B^1_{(j_1,j_2,j_3)}(\xi,\eta) :=\frac{j_1\langle\xi-\eta\rangle+j_2|\eta| - j_3\langle\xi\rangle}{2j_1j_2\langle\xi-\eta\rangle |\eta|}\eta\langle\eta\rangle^{-4}(1-\phi)(\eta).
\end{gathered}
\end{equation}
According to decomposition \eqref{B=B0+B1} we have that, for any $i=1,2,3$, \small
\begin{equation}
\begin{split}
&\int \psi_i\Big(\frac{\xi-\eta}{\langle\eta\rangle}\Big)B_{(j_1,j_2,j_3)}(\xi,\eta) \left(1-j_4j_5\frac{\xi-\eta-\zeta}{\langle\xi-\eta-\zeta\rangle}\cdot\frac{\zeta}{|\zeta|}\right)\zeta_1 \hat{u}_1(\xi-\eta-\zeta)\hat{u}_2(\zeta)\hat{u}_3(\eta) \hat{u}_4(-\xi) d\xi d\eta d\zeta \\
&= \int \psi_i\Big(\frac{\xi-\eta}{\langle\eta\rangle}\Big)B^0_{(j_1,j_2,j_3)}(\xi,\eta) \left(1-j_4j_5\frac{\xi-\eta-\zeta}{\langle\xi-\eta-\zeta\rangle}\cdot\frac{\zeta}{|\zeta|}\right)\zeta_1 \hat{u}_1(\xi-\eta-\zeta)\hat{u}_2(\zeta)\widehat{\mathrm{R} u}_3(\eta) \hat{u}_4(-\xi) d\xi d\eta d\zeta\\
& + \int \psi_i\Big(\frac{\xi-\eta}{\langle\eta\rangle}\Big)B^1_{(j_1,j_2,j_3)}(\xi,\eta) \left(1-j_4j_5\frac{\xi-\eta-\zeta}{\langle\xi-\eta-\zeta\rangle}\cdot\frac{\zeta}{|\zeta|}\right)\zeta_1 \hat{u}_1(\xi-\eta-\zeta)\hat{u}_2(\zeta)\widehat{\langle D_x\rangle^4 u}_3(\eta) \hat{u}_4(-\xi) \\
& \hspace{14cm} d\xi d\eta d\zeta\\
& =: I^0_i + I^1_i.
\end{split}
\end{equation}\normalsize
$(i)$ The first thing we observe concerning integral $I^k_i$ for $k=0,1$, $i=1,2$, is that $|\xi-\eta|, |\xi|\lesssim \langle\eta\rangle$ on the support of $\psi_i\big(\frac{\xi-\eta}{\langle\eta\rangle}\big)$ and that $|\zeta|\le \langle \xi-\eta-\zeta\rangle \langle\eta\rangle$.
Therefore, introducing the following multipliers
\begin{equation*}
B^{i,k}_{(j_1,\dots,j_5)}(\xi,\eta,\zeta):=\psi_i\Big(\frac{\xi-\eta}{\langle\eta\rangle}\Big) B^k_{(j_1,j_2,j_3)}(\xi,\eta)\left(1-j_4j_5\frac{\xi-\eta-\zeta}{\langle\xi-\eta-\zeta\rangle}\cdot\frac{\zeta}{|\zeta|}\right)\zeta_1 \langle\eta\rangle^{-7}\langle\xi-\eta-\zeta\rangle^{-4},
\end{equation*}
for any $j_1,\dots j_5\in\{+,-\}$, $k=0,1$, $i=1,2,$ a straight computation shows that, for any $\alpha,\beta,\gamma\in\mathbb{N}^2$,
\begin{equation} \label{est_multipliers_Bik}
\begin{gathered}
\left|\partial^\alpha_\xi \partial^\beta_\eta B^{i,k}_{(j_1,\dots,j_5)}(\xi,\eta,\zeta)\right| \lesssim \langle\zeta\rangle^{-3}|g_{\alpha,\beta}(\xi,\eta)|,\\
\left|\partial^\alpha_\xi \partial^\beta_\eta \partial^\gamma_\zeta B^{i,k}_{(j_1,\dots,j_5)}(\xi,\eta,\zeta)\right| \lesssim (|\zeta|\langle\zeta\rangle^{-1})^{1-|\gamma|}\langle\zeta\rangle^{-3}|g_{\alpha,\beta}(\xi,\eta)|, \ |\gamma|\ge 1,\\
\end{gathered}
\end{equation}
with
\begin{equation} \label{est_multipliers_g-alpha-beta}
\begin{gathered}
|g_{\alpha,0}(\xi,\eta)|\lesssim_{\alpha} \langle\eta\rangle^{-3}\langle\xi\rangle^{-3},\\
|g_{\alpha,\beta}(\xi,\eta)|\lesssim_{\alpha,\beta} (|\eta|\langle\eta\rangle^{-1})^{1-|\beta|}\langle\eta\rangle^{-3}\langle\xi\rangle^{-3}, \ |\beta|\ge 1.
\end{gathered}
\end{equation}
If
\begin{equation*}
K^{i,k}_{(j_1,\dots, j_5)}(x,y,z):=\int e^{ix\cdot\xi + iy\cdot\eta + iz\cdot\zeta} B^{i,k}_{(j_1,\dots,j_5)}(\xi,\eta,\zeta) d\xi d\eta d\zeta,
\end{equation*}
by lemma \ref{Lem_appendix: Kernel with 1 function} $(i)$ we first find that, for any $\alpha,\beta\in\mathbb{N}^2$,
\begin{equation*}
\left|\partial^\alpha_\xi\partial^\beta_\eta \int e^{iz\cdot\zeta}B^{i,k}_{(j_1,\dots,j_5)}(\xi,\eta,\zeta) d\zeta\right|\lesssim |z|^{-1}\langle z\rangle^{-2}|g_{\alpha,\beta}(\xi,\eta)|
\end{equation*}
and successively that
\begin{equation*}
\left|\partial^\alpha_\xi \int e^{iy\cdot\eta+ iz\cdot\zeta}B^{i,k}_{(j_1,\dots,j_5)}(\xi,\eta,\zeta)d\eta d\zeta\right|\lesssim |y|^{-1}\langle y\rangle^{-2} |z|^{-1}\langle z\rangle^{-2}\langle \xi\rangle^{-3},
\end{equation*}
for every $\xi\in\mathbb{R}^2$, $(y,z)\in\mathbb{R}^2\times\mathbb{R}^2$.
Corollary \ref{Cor_appendix: decay of integral operators} $(i)$ hence implies that 
\[|K^{i,k}_{(j_1,\dots,j_5)}(x,y,z)|\lesssim \langle x\rangle^{-3}  |y|^{-1}\langle y\rangle^{-2} |z|^{-1}\langle z\rangle^{-2}, \quad\forall (x,y,z)\in(\mathbb{R}^2)^3.\]
As for $i=1,2$
\begin{equation*}
\begin{split}
I^0_i & = \int B^{i,0}_{(j_1,\dots,j_5)}(\xi,\eta,\zeta) \widehat{\langle D_x\rangle^4 u_1}(\xi-\eta-\zeta) \hat{u}_2(\zeta) \reallywidehat{\langle D_x\rangle^7 \mathrm{R}u_3}(\eta) \hat{u}_4(-\xi)\, d\xi d\eta d\zeta, \\
& = \int K^{i,0}_{(j_1,\dots,j_5)}(t-x,x-z,x-y) [\langle D_x\rangle^4 u_1](x) u_2(y) [\langle D_x\rangle^7 \mathrm{R}u_3](z) u_4(t) dxdydzdt,
\end{split}
\end{equation*}
\begin{equation*}
\begin{split}
I^1_i &= \int B^{i,1}_{(j_1,\dots,j_5)}(\xi,\eta,\zeta) \widehat{\langle D_x\rangle^4 u_1}(\xi-\eta-\zeta) \hat{u}_2(\zeta) \widehat{\langle D_x\rangle^{11}u_3}(\eta) \hat{u}_4(-\xi)\, d\xi d\eta d\zeta \\
& = \int K^{i,1}_{(j_1,\dots,j_5)}(t-x,x-z,x-y) [\langle D_x\rangle^4 u_1](x) u_2(y) [\langle D_x\rangle^{11}u_3](z) u_4(t) dxdydzdt,
\end{split}
\end{equation*}
inequality \eqref{ineq:psi1-psi2} follows by the fact that, for any $\widetilde{u}_1, \dots \widetilde{u}_4\in L^2\cap L^\infty$, any $f,g,h\in L^1$, integrals such as
\begin{equation} \label{example_integral_u1-u4}
 \int f(t-x)g(x-z)h(x-y) |\widetilde{u}_1(x)|  |\widetilde{u}_2(y)|  |\widetilde{u}_3(z)|  |\widetilde{u}_4(t)| dxdydzdt 
\end{equation} 
can be bounded from above by the product of the $L^2$ norm of any two functions $\widetilde{u}_k$ times the $L^\infty$ norm of the remaining ones.

$(ii)$ For a cut-off function $\phi$ as the one introduced at the beginning of the proof we decompose integral $I^k_3$, $k=0,1,$ distinguishing between $|\zeta|\lesssim 1$ and $|\zeta|\gtrsim 1$.
On the one hand, for any $j_1,\dots,j_5, k=0,1,$ we consider 
\begin{equation*}
B^{3,k}_{(j_1,\dots,j_5)}(\xi,\eta,\zeta):=\psi_3\Big(\frac{\xi-\eta}{\langle\eta\rangle}\Big)\phi(\zeta)B^k_{(j_1,j_2,j_3)}(\xi,\eta)\left(1-j_4j_5\frac{\xi-\eta-\zeta}{\langle\xi-\eta-\zeta\rangle}\cdot\frac{\zeta}{|\zeta|}\right)\zeta_1 \langle\xi-\eta-\zeta\rangle^{-3}
\end{equation*}
and observe that, since $|\xi|\le \langle\xi-\eta-\zeta\rangle$ on the support of $\psi_3\big(\frac{\xi-\eta}{\langle\eta\rangle}\big)\phi(\zeta)$, the above multiplier satisfies estimates \eqref{est_multipliers_Bik}, \eqref{est_multipliers_g-alpha-beta}. From the same argument as before this implies that
\begin{multline} \label{ineq: J03}
\left|J_3^0:= \int B^{3,0}_{(j_1,\dots,j_5)}(\xi,\eta,\zeta) \widehat{\langle D_x\rangle^3u_1}(\xi-\eta-\zeta) \hat{u}_2(\zeta) \widehat{\mathrm{R}u_3}(\eta) \hat{u}_4(-\xi) d\xi d\eta d\zeta \right| \\
\lesssim \|u_1\|_{H^{3,\infty}}\|u_2\|_{L^2}\|\mathrm{R}u_3\|_{L^\infty}\|u_4\|_{L^2}
\end{multline}
together with
\begin{multline} \label{ineq: J13}
\left|J_3^1:= \int B^{3,1}_{(j_1,\dots,j_5)}(\xi,\eta,\zeta) \widehat{\langle D_x\rangle^3u_1}(\xi-\eta-\zeta) \hat{u}_2(\zeta) \widehat{\langle D_x\rangle^4 u_3}(\eta) \hat{u}_4(-\xi) d\xi d\eta d\zeta \right| \\
\lesssim \|u_1\|_{H^{3,\infty}}\|u_2\|_{L^2}\|u_3\|_{H^{4,\infty}}\|u_4\|_{L^2}.
\end{multline}
On the other hand, we make a further decomposition on the integral restricted to $|\zeta|\gtrsim 1$ by means of functions $\psi_i, i=1,2,3$, distinguishing between three regions: for $|\zeta|\le c \langle\xi-\eta\rangle$, for $c'\langle \xi-\eta\rangle\le |\zeta|\le C'\langle\xi-\eta\rangle$ and $|\zeta|>C\langle\xi-\eta\rangle$.
For any $j_1,\dots,j_5\in\{+,-\}$, $k=0,1$, we hence introduce the following multipliers
\begin{multline*}
\widetilde{B}^{3,i,k}_{(j_1,\dots,j_5)}(\xi,\eta,\zeta):= \psi_3\Big(\frac{\xi-\eta}{\langle\eta\rangle}\Big)(1-\phi)(\zeta)\psi_i\Big(\frac{\zeta}{\langle\xi-\eta\rangle}\Big)\\
\times B^k_{(j_1,j_2,j_3)}(\xi,\eta)\left(1-j_4j_5\frac{\xi-\eta-\zeta}{\langle\xi-\eta-\zeta\rangle}\cdot\frac{\zeta}{|\zeta|}\right)\zeta_1 \langle\xi-\eta-\zeta\rangle^{-7};
\end{multline*}
for $i=1m,3$, and
\begin{multline}\label{multiplier_Btilde,3,2,k}
\widetilde{B}^{3,2,k}_{(j_1,\dots,j_5)}(\xi,\eta,\zeta):= \psi_3\Big(\frac{\xi-\eta}{\langle\eta\rangle}\Big)(1-\phi)(\zeta)\psi_2\Big(\frac{\zeta}{\langle\xi-\eta\rangle}\Big)\\
\times B^k_{(j_1,j_2,j_3)}(\xi,\eta)\left(1-j_4j_5\frac{\xi-\eta-\zeta}{\langle\xi-\eta-\zeta\rangle}\cdot\frac{\zeta}{|\zeta|}\right)\zeta_1 \langle\zeta\rangle^{-1}.
\end{multline}
Since $|\xi|\sim |\xi-\eta|\sim|\xi-\eta-\zeta|$ on the support of $\psi_3\big(\frac{\xi-\eta}{\langle\eta\rangle}\big)(1-\phi)(\zeta)\psi_1\big(\frac{\zeta}{\langle\xi-\eta\rangle}\big)$ (resp.  $|\xi|\sim |\xi-\eta|\lesssim |\zeta|\sim |\xi-\eta-\zeta|$ on the support of $\psi_3\big(\frac{\xi-\eta}{\langle\eta\rangle}\big)(1-\phi)(\zeta)\psi_3\big(\frac{\zeta}{\langle\xi-\eta\rangle}\big)$), a straight computation shows that above multipliers verify \eqref{est_multipliers_Bik}, \eqref{est_multipliers_g-alpha-beta}, from which follows that
\begin{multline} \label{ineq: Jtilde-i0}
\left|\widetilde{J}^{i,0}_3:= \int \widetilde{B}^{3,i, 0}_{(j_1,\dots,j_5)}(\xi,\eta,\zeta) \widehat{\langle D_x\rangle^7u_1}(\xi-\eta-\zeta) \hat{u}_2(\zeta) \widehat{\mathrm{R}u_3}(\eta) \hat{u}_4(-\xi) d\xi d\eta d\zeta \right| \\
\lesssim \|u_1\|_{H^{7,\infty}}\|u_2\|_{L^2}\|\mathrm{R}u_3\|_{L^\infty}\|u_4\|_{L^2}
\end{multline}
along with
\begin{multline} \label{ineq: Jtilde-i1}
\left|\widetilde{J}^{i,1}_3:= \int \widetilde{B}^{3,i, 1}_{(j_1,\dots,j_5)}(\xi,\eta,\zeta) \widehat{\langle D_x\rangle^7u_1}(\xi-\eta-\zeta) \hat{u}_2(\zeta) \widehat{\langle D_x\rangle^4 u_3}(\eta) \hat{u}_4(-\xi) d\xi d\eta d\zeta \right| \\
\lesssim \|u_1\|_{H^{7,\infty}}\|u_2\|_{L^2}\|u_3\|_{H^{4,\infty}}\|u_4\|_{L^2},
\end{multline}
for $i=1,3$.
Finally, on the support of $\psi_3\big(\frac{\xi-\eta}{\langle\eta\rangle}\big)(1-\phi)(\zeta)\psi_2\big(\frac{\zeta}{\langle\xi-\eta\rangle}\big)$ we have that $|\xi|\sim |\xi-\eta|\sim |\zeta|$ and $|\xi-\eta-\zeta|\lesssim |\zeta|$. 
Replacing $\zeta$ with $\xi-\zeta$ by a change of coordinates
we find that, for any $\alpha,\beta,\gamma\in\mathbb{N}^2$,
\begin{equation} \label{est_Btilde-3,2k}
\begin{gathered}
\left|\partial^\alpha_\xi \partial^\gamma_\zeta \widetilde{B}^{3,2,k}_{(j_1,\dots,j_5)}(\xi,\eta,\xi-\zeta) \right|\lesssim_{\alpha,\gamma} \langle\eta\rangle^{-3}\langle\xi\rangle^{-|\alpha|}, \\
\left|\partial^\alpha_\xi \partial^\beta_\eta \partial^\gamma_\zeta \widetilde{B}^{3,2,k}_{(j_1,\dots,j_5)}(\xi,\eta,\xi-\zeta) \right|\lesssim (|\eta|\langle\eta\rangle^{-1})^{1-|\beta|}\langle\eta\rangle^{-3}\langle\xi\rangle^{-|\alpha|}, \ |\beta|\ge 1.
\end{gathered}
\end{equation}
If we introduce a Littlewood-Paley decomposition such that
\begin{equation*}
\widetilde{B}^{3,2,k}_{(j_1,\dots,j_5)}(\xi,\eta,\xi-\zeta) =\sum_{l\ge 1} \widetilde{B}^{3,2,k}_{(j_1,\dots,j_5)}(\xi,\eta,\xi-\zeta)\varphi(2^{-l}\xi),
\end{equation*}
one can check, using lemma \ref{Lem_appendix: Kernel with 1 function} $(i)$ to obtain the decay in $y$, making a change of coordinates $\xi\mapsto 2^l\xi$, some integration by parts, and using inequalities \eqref{est_Btilde-3,2k}, that
\begin{equation*}
K^{k,l}_{(j_1,\dots,j_5)}(x,y,z):=\int e^{ix\cdot\xi + iy\cdot\eta + iz\cdot\zeta}\widetilde{B}^{3,2,k}_{(j_1,\dots,j_5)}(\xi,\eta,\xi-\zeta) \varphi(2^{-l}\xi) d\xi d\eta d\zeta
\end{equation*}
is such that 
\begin{equation} \label{est_kernel_Kkl}
|K^{k,l}_{(j_1,\dots,j_5)}(x,y,z)|\lesssim 2^{2l}\langle 2^l x\rangle^{-3}|y|^{-1}\langle y\rangle^{-2}\langle z\rangle^{-3}, \quad \forall (x,y,z)\in (\mathbb{R}^2)^3.
\end{equation}
Moreover, since $|\xi|\sim|\xi-\zeta|$ on the support of $\widetilde{B}^{3,2,k}_{(j_1,\dots,j_5)}(\xi,\eta,\zeta)$ there are two other suitably supported cut-off functions $\varphi_1, \varphi_2$ such that $\varphi(2^{-l}\xi)=\varphi(2^{-l}\xi)\varphi_1(2^{-l}\xi)\varphi_2(2^{-l}(\xi-\zeta))$, for any $l\ge 1$. 
If $\Delta^l_j w:=\varphi_j(2^{-l}D_x)w$, we finally obtain that
\begin{equation*}
\begin{split}
 \widetilde{J}^{2,0}_3& :=\int \widetilde{B}^{3,2,0}_{(j_1,\dots,j_5)}(\xi,\eta,\zeta) \hat{u}_1(\xi-\eta-\zeta) \widehat{\langle D_x\rangle u}_2(\zeta) \widehat{\mathrm{R}u}_3(\eta)\hat{u}_4(-\xi) d\xi d\eta d\zeta \\
 & = \int \widetilde{B}^{3,2,0}_{(j_1,\dots,j_5)}(\xi,\eta,\xi-\zeta) \hat{u}_1(\zeta-\eta) \widehat{\langle D_x\rangle u}_2(\xi-\zeta) \widehat{\mathrm{R}u}_3(\eta)\hat{u}_4(-\xi) d\xi d\eta d\zeta\\
& =\sum_{l\ge 1}\int K^{0,l}_{(j_1,\dots,j_5)}(t-y,x-z,y-x) u_1(x) [\Delta^l_1\langle D_x\rangle u_2](y) [\mathrm{R}u_3](z) [\Delta^l_2u_4](t) dxdydzdt,
\end{split}
\end{equation*}
and by \eqref{est_kernel_Kkl} together with Cauchy-Schwarz inequality we derive that
\begin{equation} \label{ineq: Jtilde-20}
| \widetilde{J}^{2,0}_3|\lesssim \|u_1\|_{L^\infty}\|\mathrm{R}_1u_3\|_{L^\infty}\sum_{l\ge 1}\|\Delta^l_1 \langle D_x\rangle u_2\|_{L^2}\|\Delta^l_2 u_4\|_{L^2}\lesssim  \|u_1\|_{L^\infty}\|u_2\|_{H^1} \|\mathrm{R}_1u_3\|_{L^\infty} \|u_4\|_{L^2}.
\end{equation}
In a similar way we obtain that
\begin{equation*} 
\widetilde{J}^{2,1}_3 :=\int \widetilde{B}^{3,2,1}_{(j_1,\dots,j_5)}(\xi,\eta,\xi-\zeta) \hat{u}_1(\zeta-\eta) \widehat{\langle D_x\rangle u}_2(\xi-\zeta) \widehat{\langle D_x\rangle^4 u}_3(\eta)\hat{u}_4(-\xi) d\xi d\eta d\zeta
\end{equation*}
satisfies
\begin{equation}\label{ineq: Jtilde-21}
| \widetilde{J}^{2,1}_3|\lesssim \|u_1\|_{L^\infty}\|u_2\|_{H^1} \|u_3\|_{H^{4,\infty}} \|u_4\|_{L^2}.
\end{equation}
The result of statement $(ii)$ follows then from inequalities \eqref{ineq: J03}, \eqref{ineq: J13}, \eqref{ineq: Jtilde-i0}, \eqref{ineq: Jtilde-i1}, \eqref{ineq: Jtilde-20}, \eqref{ineq: Jtilde-21}, after having recognized that \small
\begin{multline*}
 \int \psi_3\Big(\frac{\xi-\eta}{\langle\eta\rangle}\Big)B_{(j_1,j_2,j_3)}(\xi,\eta) \left(1-j_4j_5\frac{\xi-\eta-\zeta}{\langle\xi-\eta-\zeta\rangle}\cdot\frac{\zeta}{|\zeta|}\right)\zeta_1 \hat{u}_1(\xi-\eta-\zeta)\hat{u}_2(\zeta)\hat{u}_3(\eta) \hat{u}_4(-\xi) d\xi d\eta d\zeta \\
 = \sum_{k=0}^1J^k_3 + \sum_{k=0}^1\sum_{i=1}^3 \widetilde{J}^{i,k}_3.
\end{multline*}\normalsize
In conclusion, the same proof of above applies to multiplier $B^3_{(j_1,j_2,j_3)}$ introduced in \eqref{multiplier_B3}, which can be decomposed as 
\begin{equation*}
j_2 B^0_{(j_1,j_2,j_3)}(\xi,\eta) + \widetilde{B}^1_{(j_1,j_2,j_3)}(\xi,\eta)\langle\eta\rangle^4
\end{equation*}
with the same $B^0_{(j_1,j_2,j_3)}$ as in \eqref{def_B0-B1} and 
\begin{equation*}
\widetilde{B}^1_{(j_1,j_2,j_3)}(\xi,\eta):=\frac{j_1\langle\xi-\eta\rangle+j_2|\eta| - j_3\langle\xi\rangle}{2j_1\langle\xi-\eta\rangle }\langle\eta\rangle^{-4}(1-\phi)(\eta).
\end{equation*}
The lack of factor $\eta_1|\eta|^{-1}$ against $B^0_{(j_1,j_2,j_3)}$, in comparison to decomposition \eqref{B=B0+B1}, is the reason why inequality \eqref{ineq:psi1-psi2} (resp. \eqref{ineq:psi3}) holds with $\|u_3\|_{H^{11,\infty}}+\|\mathrm{R}u_3\|_{H^{7,\infty}}$ (resp. $\|u_3\|_{H^{4,\infty}}+\|\mathrm{R}u_3\|_{L^\infty}$) replaced with $\|u_3\|_{H^{11,\infty}}$ (resp. with $\|u_3\|_{H^{4,\infty}}$).
\endproof
\end{lem}

\begin{lem} \label{Lem:Est_integrals_quartic-terms-2}
Under the same assumptions as in lemma \ref{Lem:Est_integrals_quartic-terms} we have that:

$(i)$
for any $j_1,\dots,j_5\in\{+,-\}$, $i=1,2$, and any $u_1,u_2,u_3,u_4$ such that $u_1\in H^{4,\infty}(\mathbb{R}^2)$, $u_2,u_4\in L^2(\mathbb{R}^2)$, $u_3\in H^{11,\infty}(\mathbb{R}^2)$ and $\delta_k\mathrm{R}_ku_3\in H^{7,\infty}(\mathbb{R}^2)$, \small
\begin{multline} \label{ineq:psi1-psi2-second lemma}
\left| \int \psi_i\Big(\frac{\xi-\eta}{\langle\eta\rangle}\Big)B^k_{(j_1,j_2,j_3)}(\xi,\eta) \left(1+j_4j_5\frac{\xi+\zeta}{\langle\xi+\zeta\rangle}\cdot\frac{\zeta}{|\zeta|}\right)\zeta_1 \hat{u}_1(-\xi-\zeta)\hat{u}_2(\zeta)\hat{u}_3(\eta) \hat{u}_4(\xi-\eta) d\xi d\eta d\zeta \right|\\
\lesssim \|u_1\|_{H^{4,\infty}}\|u_2\|_{L^2}\left(\|u_3\|_{H^{11,\infty}}+\delta_k \|\mathrm{R}_ku_3\|_{H^{7,\infty}}\right)\|u_4\|_{L^2};
\end{multline}\normalsize
$(ii)$ for any $j_1,\dots,j_5\in\{+,-\}$, and any $u_1,u_2,u_3,u_4$ such that $u_1\in H^{7,\infty}(\mathbb{R}^2)$, $u_2\in L^2(\mathbb{R}^2)$, $u_4\in H^1(\mathbb{R}^2)$, $u_3\in H^{4,\infty}(\mathbb{R}^2)$ and $\delta_k\mathrm{R}_ku_3\in L^\infty(\mathbb{R}^2)$,
\small
\begin{multline} \label{integral_psi3}
\left| \int \psi_3\Big(\frac{\xi-\eta}{\langle\eta\rangle}\Big)B^k_{(j_1,j_2,j_3)}(\xi,\eta) \left(1+j_4j_5\frac{\xi+\zeta}{\langle\xi+\zeta\rangle}\cdot\frac{\zeta}{|\zeta|}\right)\zeta_1 \hat{u}_1(-\xi-\zeta)\hat{u}_2(\zeta)\hat{u}_3(\eta) \hat{u}_4(\xi-\eta)  d\xi d\eta d\zeta \right|\\
\lesssim \|u_1\|_{H^{7,\infty}}\|u_2\|_{L^2}\left(\|u_3\|_{H^{4,\infty}}+ \delta_k \|\mathrm{R}_k u_3\|_{L^\infty}\right)\|u_4\|_{H^1}.
\end{multline}\normalsize
\proof
The proof of the statement is analogous to that of lemma \ref{Lem:Est_integrals_quartic-terms} after a change of coordinates $-\xi\mapsto \xi-\eta$. In \eqref{integral_psi3} we take the $H^1$ norm on $u_4$ instead of $u_2$, as done in \eqref{ineq:psi3}, by replacing multiplier $\widetilde{B}^{3,2,k}_{(j_1,j_2,j_3)}$ in \eqref{multiplier_Btilde,3,2,k} with
\begin{equation*}
 \psi_3\Big(\frac{\xi-\eta}{\langle\eta\rangle}\Big)(1-\phi)(\zeta)\psi_2\Big(\frac{\zeta}{\langle\xi-\eta\rangle}\Big) B^k_{(j_1,j_2,j_3)}(\xi,\eta)\left(1-j_4j_5\frac{\xi-\eta-\zeta}{\langle\xi-\eta-\zeta\rangle}\cdot\frac{\zeta}{|\zeta|}\right)\zeta_1 \langle\xi\rangle^{-1}.
\end{equation*}
\endproof
\end{lem}

\begin{lem} [Analysis of quartic terms. I] \label{Lem:Quartic_terms_1}
For any general multi-index $I$, any $j_k\in \{+,-\}$, $k=1,2,3$, let $C^I_{(j_1,j_2,j_3)}$, $C^{I,R}_{(j_1,j_2,j_3)}$ be the integrals defined in
\eqref{integral_CI}, \eqref{integral_CIR} respectively, and $D^I_{(j_1,j_2,j_3)}$, $D^{I,R}_{(j_1,j_2,j_3)}$ introduced in \eqref{def of D_k}, \eqref{def of DR_k}. Then
\begin{equation} \label{derivative DIs}
\partial_t \left[D^I_{(j_1,j_2,j_3)}+D^{I,R}_{(j_1,j_2,j_3)} \right] = -C^I_{(j_1,j_2,j_3)}-C^{I,R}_{(j_1,j_2,j_3)} + \mathfrak{D}^{I}_{\text{quart}},
\end{equation}
where $\mathfrak{D}^{I}_{\text{quart}}$ satisfies \small
\begin{equation} \label{est_DI1_quart}
\begin{split}
& \left|\mathfrak{D}^{I}_{\text{quart}}(t)\right| \\
&\lesssim \left[\|V(t,\cdot)\|^{\frac{7}{4}}_{H^{10,\infty}}\|V(t,\cdot)\|^{\frac{1}{4}}_{H^{12}}+ \|V(t,\cdot)\|_{H^{4,\infty}}\left(\|U(t,\cdot)\|_{H^{11,\infty}}+\|\mathrm{R}_1U(t,\cdot)\|_{H^{7,\infty}}\right) \right]\|W^I(t,\cdot)\|^2_{L^2}\\
&+ \sum_{\substack{(I_1,I_2)\in\mathcal{I}(I)\\ |I_2|<|I|}} \|Q^{\mathrm{kg}}_0(v^{I_1}_\pm, Du^{I_2}_\pm)(t,\cdot)\|_{L^2} \left(\|U(t,\cdot)\|_{H^{11,\infty}}+\|\mathrm{R}_1U(t,\cdot)\|_{H^{7,\infty}}\right) \|V^I(t,\cdot)\|_{L^2}.
\end{split}
\end{equation}\normalsize
\proof
Using definitions \eqref{integral_CI}, \eqref{def of D_k}, \eqref{def of B(i1,i2,i3)} with $k=1$, and system \eqref{system for hat(u)I+-, hat(v)I+-}, we find that
\begin{equation} \label{dt DI}
\begin{split}
 &-4(2\pi)^2\left[ \partial_t  D^I_{(j_1, j_2, j_3)} + C^I_{(j_1, j_2, j_3)}\right]  \\
& = \int \chi\left(\frac{\xi-\eta}{\langle\eta\rangle}\right)B^1_{(j_1,j_2,j_3)}(\xi,\eta) \left[\sum_{(I_1,I_2)\in\mathcal{I}(I)}c_{I_1,I_2}\reallywidehat{Q^{\mathrm{kg}}_0(v^{I_1}_\pm, D u^{I_2}_\pm)}\right](\xi-\eta) \hat{u}_{j_2}(\eta) \hat{v}^I_{j_3}(-\xi) d\xi d\eta \\
& + \int \chi\left(\frac{\xi-\eta}{\langle\eta\rangle}\right)B^1_{(j_1,j_2,j_3)}(\xi,\eta)\ \hat{v}^I_{j_1}(\xi-\eta) \reallywidehat{Q^{\mathrm{w}}_0(v_\pm, D_1v_\pm)}(\eta) \hat{v}^I_{j_3}(-\xi) d\xi d\eta \\
&+ \int\chi\left(\frac{\xi-\eta}{\langle\eta\rangle}\right) B^1_{(j_1,j_2,j_3)}(\xi,\eta) \hat{v}^I_{j_1}(\xi-\eta)\hat{u}_{j_2}(\eta) \left[\sum_{(I_1,I_2)\in\mathcal{I}(I)}c_{I_1,I_2}\reallywidehat{Q^{\mathrm{kg}}_0(v^{I_1}_\pm, Du^{I_2}_\pm)}\right](-\xi) d\xi d\eta \\
& =: S_1 + S_2 + S_3,
\end{split}
\end{equation}
where coefficients $c_{I_1,I_2}\in \{-1,0,1\}$ are such that $c_{I_1,I_2}=1$ when $|I_1| + |I_2| = |I|$ (in which case $D=D_1$) and $\chi\in C^\infty_0(\mathbb{R}^2)$ is equal to 1 close to the origin and has a sufficiently small support.
All integrals in the above right hand side are quartic terms for they involve the quadratic non-linearities of \eqref{system for hat(u)I+-, hat(v)I+-}.

The fact that $S_2$ is a remainder $\mathfrak{D}^{I}_{\text{quart}}$ satisfying \eqref{est_DI1_quart} follows by inequalities \eqref{estimate integral B uvw-new}, \eqref{est Hsinfty for NLw-New} with $s=7$, and the fact that
\begin{equation} \label{est_R1Qw0}
\|\mathrm{R}_1 Q^{\mathrm{w}}_0(v_\pm, D_1v_\pm)\|_{H^{7,\infty}} 
\lesssim \|V(t,\cdot)\|^{2-(2-\theta)\theta}_{H^{10,\infty}} \|V(t,\cdot)\|^{(2-\theta)\theta}_{H^{12}},
\end{equation}
for any $\theta\in ]0,1[$.
The above inequality is justified by the fact that, for any function $w\in W^{1,\infty}\cap H^1$, $\rho\in\mathbb{N}$ and any $\theta \in ]0,1[$, setting $p=\frac{2}{\theta}\in ]2,\infty[$,
\begin{multline} \label{injection_R1w}
\|\langle D_x\rangle^\rho\mathrm{R}_1 w\|_{L^\infty}\lesssim \|\langle D_x\rangle^\rho\mathrm{R}_1 w\|_{W^{1,p}} \lesssim \|\langle D_x\rangle^\rho w\|_{W^{1,p}} \lesssim \|\langle D_x\rangle^\rho w\|^{1-\theta}_{W^{1,\infty}}\|\langle D_x\rangle^\rho w\|^\theta_{H^1}\\
\lesssim \|\langle D_x\rangle^\rho w\|^{1-\theta}_{H^{2,\infty}}\|\langle D_x\rangle^\rho w\|^\theta_{H^1},
\end{multline}
as a consequence of Morrey's inequality, continuity of $\mathrm{R}_1:L^p\rightarrow L^p$ for $p<+\infty$, interpolation inequality, and the injection of $W^{1,\infty}$ into $H^{2,\infty}$. 
This implies that
\begin{equation} \label{ineq_for R1Qw0}
\|\mathrm{R}_1 Q^{\mathrm{w}}_0(v_\pm, D_1v_\pm)\|_{H^{\rho,\infty}}\lesssim \| Q^{\mathrm{w}}_0(v_\pm, D_1v_\pm)\|^{1-\theta}_{H^{\rho+2,\infty}} \| Q^{\mathrm{w}}_0(v_\pm, D_1v_\pm)\|^{\theta}_{H^{\rho+1}},
\end{equation}
for any $\rho\in\mathbb{N}$, and gives \eqref{est_R1Qw0} when $\rho=7$ after inequalities \eqref{est Hs NLw-New} with $s=8$, \eqref{est Hsinfty for NLw-New} with $s=9$. 
Therefore, for any $\theta\in ]0,1[$,
\begin{equation*}
|S_2|\lesssim \left(\|V(t,\cdot)\|^{2-\theta}_{H^{8,\infty}}\|V(t,\cdot)\|^\theta_{H^{10}} + \|V(t,\cdot)\|^{2-(2-\theta)\theta}_{H^{10,\infty}} \|V(t,\cdot)\|^{(2-\theta)\theta}_{H^{12}}\right) \|V^I(t,\cdot)\|^2_{L^2},
\end{equation*}
so choosing $\theta\ll 1$ small (e.g. $\theta\le 1/8$) and keeping in mind estimates \eqref{est: boostrap vpm}, \eqref{est: bootstrap Enn} we deduce that $S_2$ is controlled by the first term in the right hand side of \eqref{est_DI1_quart}.

Inequality \eqref{estimate integral B uvw-new} allows also to estimate all integrals in summations $S_1, S_3$ corresponding to indices $(I_1,I_2)\in\mathcal{I}(I)$ with $|I_2|<|I|$, and to bound them with the latter term in the right hand side of \eqref{est_DI1_quart}. 
This is not the case for integrals with $I_2=I$ involving quasi-linear term $Q^{\mathrm{kg}}_0(v_\pm, D_1u^I_\pm)$, because a straight application of that inequality would give a bound at the wrong energy level $n+1$, as $\|Q^{\mathrm{kg}}_0(v_\pm, D_1u^{I}_\pm)\|_{L^2}\lesssim \|V(t,\cdot)\|_{H^{1,\infty}}\|D_1U^I(t,\cdot)\|_{L^2}$.
Instead, since
\begin{equation} \label{Fourier transform Qkg0}
\reallywidehat{Q^{\mathrm{kg}}_0(v_\pm, D_1u^I_\pm)}(\xi) = \frac{i}{4}\sum_{j_4,j_5\in\{+,-\}}\int \Big(1 - j_4j_5 \frac{\xi-\zeta}{\langle\xi - \zeta\rangle}\cdot\frac{\zeta}{|\zeta|} \Big)\zeta_1\hat{v}_{j_4}(\xi - \zeta) \hat{u}^I_{j_5}(\zeta)d\zeta,
\end{equation}
we can rather write those integrals as the sum over $j_k\in \{+,-\}, k=1,\dots 4$, of the following:
\begin{subequations} \label{integrals B1 B2} 
\begin{equation} \label{integral B-1}
\int \chi\Big(\frac{\xi-\eta}{\langle\eta \rangle}\Big)B^1_{(j_1,j_2,j_3)}(\xi,\eta) \Big(1 - j_4j_5 \frac{\xi-\eta-\zeta}{\langle\xi-\eta - \zeta\rangle}\cdot\frac{\zeta}{|\zeta|} \Big)\zeta_1 \hat{v}_{j_4}(\xi-\eta - \zeta) \hat{u}^I_{j_5}(\zeta) \hat{u}_{j_2}(\eta)\hat{v}^I_{j_3}(-\xi)\ d\xi d\eta d\zeta,
\end{equation} 
\begin{equation}\label{integral B-2}
\int \chi\Big(\frac{\xi-\eta}{\langle\eta \rangle}\Big)B^1_{(j_1,j_2,j_3)}(\xi,\eta) \Big(1 + j_4j_5 \frac{\xi+\zeta}{\langle\xi+\zeta\rangle}\cdot\frac{\zeta}{|\zeta|} \Big)\zeta_1 \hat{v}^I_{j_1}(\xi-\eta) \hat{u}_{j_2}(\eta) \hat{v}_{j_4}(-\xi-\zeta) \hat{u}^I_{j_5}(\zeta)\ d\xi d\eta d\zeta,
\end{equation} 
\end{subequations}
and estimate them by using inequalities \eqref{ineq:psi1-psi2} and \eqref{ineq:psi1-psi2-second lemma} respectively. 
We hence obtain that
\begin{equation*}
\begin{split}
|S_1|+|S_3|&\lesssim  \|V(t,\cdot)\|_{H^{4,\infty}}\left(\|U(t,\cdot)\|_{H^{11,\infty}}+\|\mathrm{R}_1U(t,\cdot)\|_{H^{7,\infty}}\right) \|W^I(t,\cdot)\|^2_{L^2} \\
&+ \sum_{\substack{(I_1,I_2)\in\mathcal{I}(I)\\ |I_2|<|I|}} \|Q^{\mathrm{kg}}_0(v^{I_1}_\pm, Du^{I_2}_\pm)(t,\cdot)\|_{L^2} \left(\|U(t,\cdot)\|_{H^{11,\infty}}+\|\mathrm{R}_1U(t,\cdot)\|_{H^{7,\infty}}\right) \|V^I(t,\cdot)\|_{L^2},
\end{split}
\end{equation*}
and, since the same argument applies to $\partial_t D^{I,R}_{(j_1,j_2,j_3)}$, this also concludes the proof of the statement.
\endproof
\end{lem} 

\begin{lem}[Analysis of quartic terms. II] \label{Lem:Quartic_Terms_II}
For any general multi-index $I$, any $j_k\in \{+,-\}$, $k=1,2,3$, let $D^{I,T_{-N}}_{(j_1,j_2,j_3)}$ be defined as in \eqref{def of DIT}. Then
\begin{equation} \label{derivative DIT}
\partial_t D^{I,T_{-N}}_{(j_1,j_2,j_3)} = \Im\left[\langle T_{-N}(U)W^I, W^I\right] + \mathfrak{D}^{I,N}_{\text{quart}}
\end{equation}
and if $N\ge 18$ $\mathfrak{D}^{I,N}_{\text{quart}}$satisfies
\begin{equation} \label{est_DIN_quart}
\begin{split}
\left|\mathfrak{D}^{I,N}_{\text{quart}} \right| &\lesssim \|V(t,\cdot)\|_{H^{N+4,\infty}}^{\frac{7}{4}}\|V(t,\cdot)\|^{\frac{1}{4}}_{H^{N+6}}\|W^I(t,\cdot)\|^2_{L^2}\\
& +\sum_{\substack{(I_1,I_2)\in\mathcal{I}(I)\\ |I_2|<|I|}}\left\| Q^{\mathrm{kg}}_0(v^{I_1}_\pm, D u^{I_2}_\pm)\right\|_{L^2} \|U(t,\cdot)\|_{H^{N+3,\infty}}\|V^I(t,\cdot)\|_{L^2}.
\end{split}
\end{equation}
\proof
For any triplet $(j_1,j_2,j_3)$, we compute the time derivative of $D^{I,T_{-N}}$ by making use of system \eqref{system for hat(u)I+-, hat(v)I+-}.
Recalling \eqref{expression T-N in normal forms} and \eqref{def_multiplier_sigmatildeN}, we find that
\begin{equation}
\begin{split}
& \partial_t \left[\sum_{j_k\in\{+,-\}} D^{I,T_{-N}}_{(j_1,j_2,j_3)}\right] - \Im[\langle T_{-N}(U)W^I,W^I\rangle] = \\
& =\Re \left[ \frac{1}{(2\pi)^2} \int \widetilde{\sigma}^N_{(j_1,j_2,j_3)}(\xi,\eta) \left[\sum_{(I_1,I_2)\in \mathcal{I}(I)}c_{I_1,I_2}\reallywidehat{Q^{\mathrm{kg}}_0(v^{I_1}_\pm, D u^{I_2}_\pm)}(\xi-\eta)\right] \hat{u}_{j_2}(\eta) \hat{v}^I_{-j_3}(-\xi) d\xi d\eta \right. \\
&  \hspace{0,8cm}\left.+ \frac{1}{(2\pi)^2} \int \widetilde{\sigma}^N_{(j_1,j_2,j_3)}(\xi,\eta)  \hat{v}^I_{j_1}(\xi-\eta) \reallywidehat{Q^{\mathrm{w}}_0(v_\pm, D_1v_\pm)}(\eta) \hat{v}^I_{-j_3}(-\xi) d\xi d\eta \right. \\
&  \hspace{0,8cm}\left. +  \frac{1}{(2\pi)^2} \int \widetilde{\sigma}^N_{(j_1,j_2,j_3)}(\xi,\eta)  \hat{v}^I_{j_1}(\xi-\eta) \hat{u}_{j_2}(\eta) \left[\sum_{(I_1,I_2)\in \mathcal{I}(I)}c_{I_1,I_2}\reallywidehat{Q^{\mathrm{kg}}_0(v^{I_1}_\pm, D u^{I_2}_\pm)}(-\xi)\right] d\xi d\eta \right] \\
& =: S^{T_{-N}}_1 + S^{T_{-N}}_2 + S^{T_{-N}}_3,
\end{split}
\end{equation}
with coefficients $c_{I_1,I_2}\in\{-1,0,1\}$ such that $c_{I_1,I_2}=1$ whenever $|I_1|+|I_2|=|I|$ (in which case $D=D_1$).
After lemma \ref{Lem_appendix: integral sigma_tilde_N} and inequality \eqref{est Hsinfty for NLw-New} with $s=N+3$ we deduce that, if $N\ge 15$, for any $\theta\in ]0,1[$
\begin{equation*}
|S^{T_{-N}}_2|\lesssim  \|V(t,\cdot)\|_{H^{N+4,\infty}}^{2-\theta}\|V(t,\cdot)\|^\theta_{H^{N+6}}\|V^I(t,\cdot)\|^2_{L^2}.
\end{equation*}
Choosing $\theta\ll 1$ small (e.g. $\theta\le 1/8$) we then obtain that $S^{T_{-N}}_2$ is a remainder $\mathfrak{D}^{I,N}_{\text{quart}}$ satisfying \eqref{est_DIN_quart}.
Also, the same lemma implies that each contribution in $S^{T_{-N}}_1, S^{T_{-N}}_3$ corresponding to $(I_1,I_2)\in \mathcal{I}(I)$ with $|I_2|<|I|$ is bounded by
\begin{equation*}
\left\| Q^{\mathrm{kg}}_0(v^{I_1}_\pm, D u^{I_2}_\pm)\right\|_{L^2} \|U(t,\cdot)\|_{H^{N+3,\infty}}\|V^I(t,\cdot)\|_{L^2}.
\end{equation*}
Reminding instead \eqref{Fourier transform Qkg0}, we find that the remaining contribution to $S^{T_{-N}}_1$, corresponding to $I_2=I$, is equal to the sum over $j_1,\dots,j_5\in \{+,-\}$ of the (imaginary part) of the following integrals:
\begin{equation} \label{integral sigma_tilde_N-1}
\int \widetilde{\sigma}^N_{(j_1,j_2,j_3)}(\xi,\eta)\Big(1 - j_4j_5 \frac{\xi-\eta-\zeta}{\langle\xi-\eta - \zeta\rangle}\cdot\frac{\zeta}{|\zeta|} \Big)\zeta_1 \hat{v}_{j_4}(\xi-\eta - \zeta) \hat{u}^I_{j_5}(\zeta) \hat{u}_{j_2}(\eta)\hat{v}^I_{j_3}(-\xi)\ d\xi d\eta d\zeta.
\end{equation} 
Analogously, the contribution corresponding to $I_2=I$ in $S^{T_{-N}}_3$ is the sum over $j_k\in\{+,-\}, k=1,\dots,5$ of
\begin{equation} \label{integral sigma_tilde_N-2}
\int  \widetilde{\sigma}^N_{(j_1,j_2,j_3)}(\xi,\eta) \Big(1 + j_4j_5 \frac{\xi+\zeta}{\langle\xi+\zeta\rangle}\cdot\frac{\zeta}{|\zeta|} \Big)\zeta_1 \hat{v}^I_{j_1}(\xi-\eta) \hat{u}_{j_2}(\eta) \hat{v}_{j_4}(-\xi-\zeta) \hat{u}^I_{j_5}(\zeta)\ d\xi d\eta d\zeta.
\end{equation} 
Since $\widetilde{\sigma}^N_{(j_1,j_2,j_3)}(\xi,\eta)$ satisfies \eqref{derivatives sigma tilde N} and is supported for $|\eta|\le \varepsilon |\xi-\eta|$, for a small $0<\varepsilon\ll 1$, we rewrite above integrals, respectively, as
\begin{multline} \label{new_integral sigma_tilde_N-1}
\int \widetilde{\sigma}^N_{(j_1,j_2,j_3)}(\xi,\eta) \langle\eta\rangle^{-N-3}\Big(1 - j_4j_5 \frac{\xi-\eta-\zeta}{\langle\xi-\eta - \zeta\rangle}\cdot\frac{\zeta}{|\zeta|} \Big)\zeta_1 \langle\xi-\eta-\zeta\rangle^{-4}\\
\times  \reallywidehat{\langle D_x\rangle^4 v}_{j_4}(\xi-\eta - \zeta) \hat{ u}^I_{j_5}(\zeta) \reallywidehat{\langle D_x\rangle^{N+3} u}_{j_2}(\eta)\hat{v}^I_{j_3}(-\xi)\ d\xi d\eta d\zeta,
\end{multline}
and
\begin{multline} \label{new_integral sigma_tilde_N-2}
\int  \widetilde{\sigma}^N_{(j_1,j_2,j_3)}(\xi,\eta)\langle\eta\rangle^{-N-7} \Big(1 + j_4j_5 \frac{\xi+\zeta}{\langle\xi+\zeta\rangle}\cdot\frac{\zeta}{|\zeta|} \Big)\zeta_1\langle\xi+\zeta\rangle^{-4}\\
\times \hat{ v}^I_{j_1}(\xi-\eta) \reallywidehat{\langle D_x\rangle^{N+7} u}_{j_2}(\eta) \reallywidehat{\langle D_x\rangle^4v}_{j_4}(-\xi-\zeta) \hat{u}^I_{j_5}(\zeta)\ d\xi d\eta d\zeta.
\end{multline}
With such a choice, the new multipliers, that we denote concisely by $\widetilde{\sigma}^{N,k}_{(j_1,\dots, j_5)}(\xi,\eta,\zeta)$, $k=0,1$, verify, for any $\alpha,\beta,\gamma\in\mathbb{N}^2$,
\begin{gather*}
\left|\partial^\alpha_\xi\partial^\beta_\eta \widetilde{\sigma}^{N,k}_{(j_1,\dots, j_5)}(\xi,\eta,\zeta)\right| \lesssim \langle\zeta\rangle^{-3} |g^N_{\alpha,\beta}(\xi)|, \\
\left|\partial^\alpha_\xi\partial^\beta_\eta \partial^\gamma_\zeta \widetilde{\sigma}^{N,k}_{(j_1,\dots, j_5)}(\xi,\eta,\zeta)\right| \lesssim (|\zeta|\langle \zeta\rangle^{-1})^{1-|\gamma|}\langle\zeta\rangle^{-3} |g^N_{\alpha,\beta}(\xi)|, \quad |\gamma|\ge 1,
\end{gather*}
with $g^N_{\alpha,\beta}(\xi,\eta)$ supported for $|\eta|\le \varepsilon |\xi-\eta|$ and such that
\begin{gather*}
|g^N_{\alpha,\beta}(\xi,\eta)|\lesssim \langle \xi-\eta\rangle^{6-N+|\alpha| +2|\beta|} |\eta|^{N-|\beta|}\langle \eta\rangle^{-N-3}, \quad \forall (\xi,\eta)\in\mathbb{R}^2\times\mathbb{R}^2.
\end{gather*}
If $N\in\mathbb{N}^*$ is sufficiently large (e.g. $N\ge 18$), the above estimate implies that,
for any $\alpha,\beta\in\mathbb{N}^2$ of length less or equal than 3, 
\[|g^N_{\alpha,\beta}(\xi,\eta)|\lesssim \langle\eta\rangle^{-3}\langle\xi\rangle^{-3},\]
so by lemma \ref{Lem_appendix: Kernel with 1 function} $(i)$ together with corollary \ref{Cor_appendix: decay of integral operators} $(i)$ we obtain that, for any $k=0,1$,
\begin{equation*}
K^{N,k}_{(j_1,\dots,j_5)}(x,y,z):=\int e^{ix\cdot\xi + i y\cdot \eta + iz\cdot\zeta} \widetilde{\sigma}^{N,k}_{(j_1,\dots, j_5)}(\xi,\eta,\zeta) d\xi d\eta d\zeta
\end{equation*} 
is such that 
\begin{equation} \label{est_KNk}
|K^{N,k}_{(j_1,\dots,j_5)}(x,y,z)|\lesssim \langle x\rangle^{-3}|y|^{-1}\langle y \rangle^{-2}|z|^{-1}\langle z \rangle^{-2}, \quad \forall (x,y,z)\in(\mathbb{R}^2)^3.
\end{equation}
By \eqref{new_integral sigma_tilde_N-1}, \eqref{new_integral sigma_tilde_N-2}, integrals \eqref{integral sigma_tilde_N-1}, \eqref{integral sigma_tilde_N-2} are respectively equal to
\begin{equation*}
\int K^{N,0}_{(j_1,\dots,j_5)}(t-x, x-z, x-y) [\langle D_x\rangle^4v_{j_4}](x) u^I_{j_5}(y) [\langle D_x\rangle^{N+3} u_{j_2}](z) v^I_{j_3}(t) dxdydzdt 
\end{equation*}
and
\begin{equation*}
\int K^{N,1}_{(j_1,\dots,j_5)}(z-x, x-y, z-t) v^I_{j_1}(x) [\langle D_x\rangle^{N+7} u_{j_2}](y) [\langle D_x\rangle^4v_{j_4}](z) u^I_{j_5}(t) dxdydzdt.
\end{equation*}
Using \eqref{est_KNk} and the fact that integrals such as \eqref{example_integral_u1-u4} can be bounded from above by the product of the $L^2$ norm of any two functions $\widetilde{u}_k$ times the $L^\infty$ norm of the remaining ones, they are estimated by 
\begin{equation*}
\|V(t,\cdot)\|_{H^{4,\infty}}\|U(t,\cdot)\|_{H^{N+7,\infty}}\|W^I(t,\cdot)\|^2_{L^2},
\end{equation*}
which concludes the proof of the statement.
\endproof
\end{lem}

\begin{lem}[Analysis of quartic terms. III]\label{Lem:Quartic_terms_III}
Let $n\in\mathbb{N}, n\ge 3$, $I\in\mathcal{I}_n$ and $(I_1,I_2)\in\mathcal{I}(I)$ be such that $[\frac{|I|}{2}]<|I_1|<|I|$.
Let also $C^{I_1,I_2}_{(j_1,j_2,j_3)}$, $D^{I_1,I_2}_{(j_1,j_2,j_3)}$ be the integrals defined, respectively, in \eqref{integral_CI1I2}, \eqref{def of DI1I2}, for any $j_k\in \{+,-\}, k=1,2,3$.
Then
\begin{equation} \label{derivative DI1I2}
\partial_t D^{I_1,I_2}_{(j_1,j_2,j_3)} =- C^{I_1,I_2}_{(j_1,j_2,j_3)}+ \mathfrak{D}^{I_1,I_2}_{\text{quart}},
\end{equation}
where $\mathfrak{D}^{I_1,I_2}_{\text{quart}}$ satisfies
\begin{multline} \label{est_DI1I2_quart}
\left|\mathfrak{D}^{I_1,I_2}_{\text{quart}}(t) \right|\lesssim \left[\left(\|W(t,\cdot)\|_{H^{[\frac{n}{2}]+12,\infty}}+ \|\mathrm{R}_1U(t,\cdot)\|_{H^{[\frac{n}{2}]+8,\infty}}\right)^2 \right. \\ \left. +  \|V(t,\cdot)\|^{\frac{7}{4}}_{H^{[\frac{n}{2}]+11,\infty}} \|V(t,\cdot)\|^{\frac{1}{4}}_{H^{[\frac{n}{2}]+12}}\right] E_n(t;W).
\end{multline}
\proof
We compute the time derivative of $D^{I_1,I_2}_{(j_1,j_2,j_3)}$ by making use of system \eqref{system for hat(u)I+-, hat(v)I+-}. 
We remind that, after remark \ref{Remark:Vector_field_on_null_structure} and definition \eqref{def mathcal(I)}, if $\Gamma^I$ is a product of spatial derivatives then all couples of indices $(I_1,I_2)$ belonging to $\mathcal{I}(I)$ are such that $|I_1|+|I_2|=|I|$ and $\Gamma^{I_1},\Gamma^{I_2}$ are also products of spatial derivatives. Therefore, all coefficients $c_{I_1,I_2}$ appearing in the right hand side of \eqref{system for hat(u)I+-, hat(v)I+-} are equal to 0.
By definitions \eqref{def of B(i1,i2,i3)} with $k=1$, \eqref{integral_CI1I2}, \eqref{def of DI1I2}, we find that
\begin{equation}
\begin{split}
-4(2\pi)^2&\Big[ \partial_t  D^{I_1,I_2}_{(j_1, j_2, j_3)}+ C^{I_1,I_2}_{(j_1, j_2, j_3)} \Big]= \\
&   \int B^1_{(j_1,j_2,j_3)}(\xi,\eta) \left[\sum_{(J_1,J_2)\in\mathcal{I}(I_1)}\reallywidehat{Q^{\mathrm{kg}}_0(v^{J_1}_\pm, D_1 u^{J_2}_\pm)}(\xi-\eta)\right]\hat{u}^{I_2}_{j_2}(\eta) \hat{v}^I_{j_3}(-\xi) d\xi d\eta \\
+ &  \int B^1_{(j_1,j_2,j_3)}(\xi,\eta)\ \hat{v}^{I_1}_{j_1}(\xi-\eta) \left[\sum_{(J_1,J_2)\in\mathcal{I}(I_2)}\reallywidehat{Q^{\mathrm{w}}_0(v^{J_1}_\pm, D_1v^{J_2}_\pm)}(\eta)\right] \hat{v}^I_{j_3}(-\xi) d\xi d\eta \\
+ &  \int B^1_{(j_1,j_2,j_3)}(\xi,\eta) \hat{v}^{I_1}_{j_1}(\xi-\eta)\hat{u}^{I_2}_{j_2}(\eta) \left[\sum_{(J_1,J_2)\in\mathcal{I}(I)}\reallywidehat{Q^{\mathrm{kg}}_0(v^{J_1}_\pm, D_1u^{J_2}_\pm)}\right](-\xi) d\xi d\eta \\
 =: & S^{I_1,I_2}_1 + S^{I_1,I_2}_2 + S^{I_1,I_2}_3.
\end{split}
\end{equation}
Since $|J_1|+|J_2|=|I_1|<|I|\le n$ in $S^{I_1,I_2}_1$, we can estimate all its contributions using inequality \eqref{estimate integral B uvw-new}. Using lemma \ref{Lem: L2 est nonlinearities} $(i)$, the fact that $|I_2|\le [\frac{n}{2}]$ by the hypothesis and, hence, that 
\begin{equation*}
\|u^{I_2}_\pm(t,\cdot)\|_{H^{7,\infty}}+\|\mathrm{R}_1u^{I_2}_\pm(t,\cdot)\|_{H^{7,\infty}}\lesssim \|U(t,\cdot)\|_{H^{[\frac{n}{2}]+8,\infty}},
\end{equation*}
we deduce that
\begin{equation*}
\left|S^{I_1,I_2}_1\right|\lesssim \left(\|W(t,\cdot)\|_{H^{[\frac{n}{2}]+2}} + \|\mathrm{R}_1U(t,\cdot)\|_{H^{[\frac{n}{2}]+2,\infty}}\right)\|U(t,\cdot)\|_{H^{[\frac{n}{2}]+8,\infty}}E_n(t;W),
\end{equation*}
and above estimate holds also for all integrals in $S^{I_1,I_2}_3$ corresponding to $|J_2|<|I|$.
The same inequality \eqref{estimate integral B uvw-new}, combined with \eqref{ineq_for R1Qw0} applied to $Q^\mathrm{w}_0(v^{J_1}_\pm, D_1v^{J_2}_\pm)$ and with corollary \ref{Cor_appendixA:Hs-Hsinfty norm of bilinear expressions} in appendix \ref{Appendix A}, gives that, for any $\theta\in ]0,1[$,
\begin{equation*}
\begin{split}
&|S^{I_1,I_2}_2| \\
&\lesssim \sum_{|J_1|+|J_2|= |I_2|}\left[\left\|Q^\mathrm{w}_0(v^{J_1}_\pm, D_1v^{J_2}_\pm) \right\|_{H^{7, \infty}} + \left\|Q^\mathrm{w}_0(v^{J_1}_\pm, D_1v^{J_2}_\pm) \right\|^{1-\theta}_{H^{9, \infty}}\left\|Q^\mathrm{w}_0(v^{J_1}_\pm, D_1v^{J_2}_\pm) \right\|^\theta_{H^8} \right] E_n(t;W)\\
&\lesssim \|V(t,\cdot)\|^{2-(2-\theta)\theta}_{H^{[\frac{n}{2}]+11,\infty}} \|V(t,\cdot)\|^{(2-\theta)\theta}_{H^{[\frac{n}{2}]+12}} E_n(t;W).
\end{split}
\end{equation*}
Finally, the last remaining integral in $S^{I_1,I_2}_3$, corresponding to indices $J_1=0,J_2=I$, can be written using \eqref{Fourier transform Qkg0} as
\begin{equation*}
\sum_{j_4, j_4\in \{+,-\}}\int B^1_{(j_1,j_2,j_3)}(\xi,\eta)\left(1+j_4j_5\frac{\xi+\zeta}{\langle\xi+\zeta\rangle}\cdot\frac{\zeta}{|\zeta|}\right)\zeta_1 \hat{v}^{I_1}_{j_1}(\xi-\eta)\hat{u}^{I_2}_{j_2}(\eta)\hat{v}_{j_4}(-\xi-\zeta)\hat{u}^I_{j_5}(\zeta) 
d\xi d\eta d\zeta,
\end{equation*} 
and is estimated, after lemma \ref{Lem:Est_integrals_quartic-terms-2} and the fact that $|I_1|<|I|$, by
\begin{equation*}
\|V(t,\cdot)\|_{H^{7,\infty}}\left(\|U(t,\cdot)\|_{H^{[\frac{n}{2}]+12,\infty}}+ \|\mathrm{R}_1U(t,\cdot)\|_{H^{[\frac{n}{2}]+8,\infty}}\right)E_n(t;W).
\end{equation*}
This gives that
\begin{equation*}
\left| S^{I_1,I_2}_3\right|\lesssim \left(\|W(t,\cdot)\|_{H^{[\frac{n}{2}]+12,\infty}}+ \|\mathrm{R}_1U(t,\cdot)\|_{H^{[\frac{n}{2}]+8,\infty}}\right)^2 E_n(t;W)
\end{equation*}
and concludes the proof of the statement.
\endproof
\end{lem}

\begin{lem}[Analysis of quartic terms. IV] \label{Lem:Analysis quartic terms IV}
Let $k=0,1$, $\mathcal{K},\mathcal{V}_k$ be the sets introduced in \eqref{set_K}, \eqref{set_V} respectively, $I\in \mathcal{V}^k$ and $(I_1,I_2)\in\mathcal{I}(I)$ be such that $I_1\in\mathcal{K}$, $|I_2|\le 1$. Let also $F^{I_1,I_2,l}_{(j_1,j_2,j_3)}$, $G^{I_1,I_2,l}_{(j_1,j_2,j_3)}$ be the integrals defined in \eqref{integral_FI1I2}, \eqref{integral_GI1I2}, for any $l=1,2,3$, $j_i\in\{+,-\}, i=1,2,3$.
For any $l=1,2,3$, any triplet $(j_1,j_2,j_3)$, we have that
\begin{equation} \label{derivative_GI1I2}
\partial_t G^{I_1,I_2,l}_{(j_1,j_2,j_3)} = -F^{I_1,I_2,l}_{(j_1,j_2,j_3)} + \mathfrak{G}^{I_1,I_2}_{\text{quart}},
\end{equation}
and there is a constant $C>0$ such that, if a-priori estimates \eqref{est: bootstrap argument a-priori est} are satisfied in interval $[1,T]$ for a fixed $T>1$, with $\varepsilon_0<(2A+B)^{-1}$ small, 
\begin{equation} \label{est_GI1I2quart}
|\mathfrak{G}^{I_1,I_2}_{\text{quart}}(t)|\le C(A+B)^2\varepsilon^2 t^{-1+\frac{\delta_k}{2}}\left[E^k_3(t;W)^\frac{1}{2}+\delta_{\mathcal{V}^0}t^{\beta+\frac{\delta_1}{2}}E^1_3(t;W)^\frac{1}{2}+ t^{-\frac{1}{4}-\frac{\delta_k}{2}}\right],
\end{equation}
for every $t\in [1,T]$, with $\delta_{\mathcal{V}^0}=1$ if $I\in\mathcal{V}^0$, 0 otherwise, and $\beta>0$ as small as we want.
\proof
First of all, it is useful to remind that from \eqref{product vI1, uI2_1}, \eqref{product_vI1-uI2_2} and a-priori estimate \eqref{est: bootstrap E02}, for any $k=0,1$, $I\in\mathcal{I}^k_3$, $(I_1,I_2)\in\mathcal{I}(I)$ such that $I_1\in\mathcal{K}$, $|I_2|\le 1$, and $\sigma>0$ sufficiently small
\begin{equation} \label{est_VI1UI2}
\|V^{I_1}(t,\cdot)\|_{L^2}\left(\|\chi(t^{-\sigma}D_x)U^{I_2}(t,\cdot)\|_{H^{\rho,\infty}}+ \|\chi(t^{-\sigma}D_x)\mathrm{R} U^{I_2}(t,\cdot)\|_{H^{\rho,\infty}}\right) \le C(A+B)B\varepsilon^2 t^{-\frac{1}{2}+\frac{\delta_k}{2}},
\end{equation}
for every $t\in [1,T]$.

For any fixed $(j_1,j_2,j_3)$, any $l=1,2,3$, we compute $\partial_t G^{I_1,I_2,l}_{(j_1,j_2,j_3)}$ recurring to system \eqref{system for hat(u)I+-, hat(v)I+-} along with its compact form
\begin{equation*}
\begin{cases}
(D_t \mp \langle D_x\rangle)v^I_\pm = \Gamma^I Q^\mathrm{w}_0(v_\pm, D_1 v_\pm),\\
(D_t \mp | D_x|)u^I_\pm = \Gamma^I Q^\mathrm{kg}_0(v_\pm, D_1 u_\pm),
\end{cases}
\end{equation*}
and using that $[D_t, \chi(t^{-\sigma}D_x)]=t^{-1}\chi_1(t^{-\sigma}D_x)$ with $\chi_1(\xi):=i\sigma (\partial\chi)(\xi)\cdot\xi$. 
We find that\small
\begin{equation*}
\begin{split}
&-4(2\pi)^2 \left[\partial_t G^{I_1,I_2,l}_{(j_1,j_2,j_3)} + F^{I_1,I_2,l}_{(j_1,j_2,j_3)}\right] \\
&= \int B^l_{(j_1,j_2,j_3)}(\xi,\eta) \Big[\reallywidehat{\Gamma^{I_1}Q^\mathrm{kg}_0(v_\pm, D_1u_\pm)}(\xi-\eta)\Big] \reallywidehat{\chi(t^{-\sigma}D_x)u^{I_2}_{j_2}}(\eta)\hat{v}^I_{j_3}(-\xi) d\xi d\eta \\
&+ \int B^l_{(j_1,j_2,j_3)}(\xi,\eta)  \hat{v}^{I_1}_{j_1}(\xi-\eta) \Big[\reallywidehat{\chi(t^{-\sigma}D_x)\Gamma^{I_2}Q^\mathrm{w}_0(v_\pm, D_1 v_\pm)}(\eta) + t^{-1}\reallywidehat{\chi_1(t^{-\sigma}D_x)u^{I_2}_{j_2}}(\eta)\Big]\hat{v}^I_{j_3}(-\xi) d\xi d\eta \\
& + \int B^l_{(j_1,j_2,j_3)}(\xi,\eta)  \hat{v}^{I_1}_{j_1}(\xi-\eta)  \reallywidehat{\chi(t^{-\sigma}D_x)u^{I_2}_{j_2}}(\eta)\Big[\sum_{(J_1,J_2)\in\mathcal{I}(I)}c_{J_1,J_2}\reallywidehat{Q^\mathrm{kg}_0(v^{J_1}_\pm, Du^{J_2}_\pm)}(-\xi)\Big] d\xi d\eta\\
&=: S^{I_1,I_2,l}_1 + S^{I_1,I_2,l}_2 + S^{I_1,I_2,l}_3,
\end{split}
\end{equation*}\normalsize
with $B^l_{(j_1,j_2,j_3)}$ given by \eqref{def of B(i1,i2,i3)} when $l=1,2$ or \eqref{multiplier_B3} when $l=3$. 

Applying \eqref{estimate integral B uvw-new} to $S^{I_1,I_2,l}_2$, using \eqref{injection_R1w} with $w=\Gamma^{I_2}Q^\mathrm{w}_0(v_\pm, D_1v_\pm)$ and $\rho=7$, together with the fact that operators $\chi(t^{-\sigma}D_x)$, $\chi_1(t^{-\sigma}D_x)$ are bounded from $L^\infty$ to $H^{\rho,\infty}$ for any $\rho\ge 0$ with norm $O(t^{\sigma\rho})$, and from $L^2$ to $H^s$ for any $s\ge0$ with norm $O(t^{\sigma s})$, we deduce that, for any $\theta\in ]0,1[$,
\begin{multline} \label{preliminary_SI1I2_2}
|S^{I_1,I_2,l}_2|\lesssim t^\beta \|V^{I_1}(t,\cdot)\|_{L^2}\|V^I(t,\cdot)\|_{L^2}\\
\times \left[\|\Gamma^{I_2}Q^\mathrm{w}_0(v_\pm, D_1v_\pm)\|_{L^\infty}+ \delta_l \|\Gamma^{I_2}Q^\mathrm{w}_0(v_\pm, D_1v_\pm)\|_{L^\infty}^{1-\theta}\|\Gamma^{I_2}Q^\mathrm{w}_0(v_\pm, D_1v_\pm)\|^\theta_{L^2} \right.\\
\left.+ t^{-1} \left(\|\chi_1(t^{-\sigma}D_x)u^{I_2}_\pm(t,\cdot)\|_{L^\infty}+\|\chi_1(t^{-\sigma}D_x)\mathrm{R} u^{I_2}_\pm(t,\cdot)\|_{L^\infty}\right)\right],
\end{multline}
for some $\beta>0$ small, $\beta\rightarrow 0$ as $\sigma\rightarrow 0$, and with $\delta_l=1$ if $l=1,2$, 0 otherwise.
When $|I_2|=0$ the above right hand side can be estimated using \eqref{est L2 NLw}, \eqref{est Linfty NLw} and a-priori estimates \eqref{est: bootstrap argument a-priori est}.
When $|I_2| = 1$ we derive from \eqref{Gamma_nonlinearity} that
\begin{equation*}
\Gamma^{I_2}Q^\mathrm{w}_0(v_\pm, D_1v_\pm) = Q^\mathrm{w}_0(v^{I_2}_\pm, D_1v_\pm) + Q^\mathrm{w}_0(v_\pm, D_1v^{I_2}_\pm) + G^w_1(v_\pm, D v_\pm)
\end{equation*}
with $G^w_1(v_\pm, Dv_\pm) = G_1(v, \partial v)$ given by \eqref{def_G1}. 
Using lemma \ref{Lem_app:products_Gamma} in appendix \ref{Appendix B} with $L=L^\infty$ to estimate the $L^\infty$ norm of the first two quadratic terms in the above right hand side, we find that, for some new $\chi\in C^\infty_0(\mathbb{R}^2)$ and $\sigma>0$ small, there is a constant $C>0$ such that
\begin{equation*}
\begin{split}
\|\Gamma^{I_2}& Q^\mathrm{w}_0(v_\pm, D_1v_\pm)\|_{L^\infty}\lesssim \left\|\chi(t^{-\sigma}D_x)v^{I_2}_\pm(t,\cdot)\right\|_{H^{2,\infty}}\|v_\pm(t,\cdot)\|_{H^{2,\infty}} \\
&+ t^{-N(s)}\left(\|v_\pm(t,\cdot)\|_{H^s}+ \|D_tv_\pm(t,\cdot)\|_{H^s}\right)\Big(\sum_{|\mu|=0}^1\|x^\mu v_\pm(t,\cdot)\|_{H^1}+ t\|v_\pm(t,\cdot)\|_{H^1}\Big) \\
& + \|v_\pm(t,\cdot)\|_{H^{1,\infty}}\left(\|v_\pm(t,\cdot)\|_{H^{2,\infty}}+\|D_tv_\pm(t,\cdot)\|_{H^{1,\infty}}\right)\\
&\le CAB\varepsilon^2 t^{-2},
\end{split}
\end{equation*}
where last inequality is obtained by picking $s>0$ sufficiently large so that $N(s)\ge 4$ and using \eqref{Hs norm DtV}, \eqref{est: Hsinfty Dt V}, \eqref{norm_xv-}, lemma \ref{Lem_appendix: sharp_est_VJ}, together with a-priori estimates.
Also, by \eqref{Hs norm DtV} with $s=0$ and a-priori estimates
\begin{equation*}
\|\Gamma^{I_2}Q^\mathrm{w}_0(v_\pm, D_1v_\pm)\|_{L^2} \lesssim \|V(t,\cdot)\|_{H^{2,\infty}}\left(\|V^{I_2}(t,\cdot)\|_{H^1}+\|D_tV(t,\cdot)\|_{L^2}\right)\le CAB\varepsilon^2t^{-1+\frac{\delta_2}{2}}.
\end{equation*}
Therefore, using lemma \ref{Lem_appendix: est UJ} and taking $\theta,\sigma>0$ sufficiently small we deduce from \eqref{preliminary_SI1I2_2} and the above estimates that, for any $l=1,2,3$ and a new $C>0$,
\begin{equation} \label{est_SI1I2_2}
|S^{I_1,I_2,l}_2|\le CAB\varepsilon^2 t^{-\frac{5}{4}}E^k_3(t;W)^\frac{1}{2}.
\end{equation}
We make use of inequality \eqref{estimate integral B uvw-new} to estimate $S^{I_1,I_2,l}_1$, too.
From \eqref{Gamma_I_nonlinearity} we have that
\begin{equation*}
\Gamma^{I_1}Q^\mathrm{kg}_0(v_\pm, D_1u_\pm) = Q^\mathrm{kg}_0(v^{I_1}_\pm, D_1u_\pm)+ \sum_{\substack{(J_1,J_2)\in\mathcal{I}(I_1)\\ |J_1|<|I_1|}}c_{J_1,J_2} Q^\mathrm{kg}_0(v^{J_1}_\pm, Du^{J_2}_\pm)
\end{equation*}
with $c_{J_1,J_2}\in \{-1,0,1\}$,
and then from \eqref{decomposition_Qkg0}, \eqref{eq:Qkg0(Dt)-statement} and the fact that $I_1\in\mathcal{K}$,
\begin{equation*}
\Gamma^{I_1}Q^\mathrm{kg}_0(v_\pm, D_1u_\pm)= Q^\mathrm{kg}_0(v^{I_1}_\pm, \chi(t^{-\sigma}D_x)D_1u_\pm) +\mathfrak{R}^k_3(t,x),
\end{equation*}
with $\mathfrak{R}^k_3$ satisfying \eqref{est:L2_norm_Rk3(t,x)} and
\begin{equation*}
\| Q^\mathrm{kg}_0(v^{I_1}_\pm, \chi(t^{-\sigma}D_x)D_1 u_\pm)\|_{L^2}\le \left(\|U(t,\cdot)\|_{H^{2,\infty}}+\|\mathrm{R}U(t,\cdot)\|_{H^{2,\infty}}\right)\|V^{I_1}(t,\cdot)\|_{L^2}.
\end{equation*}
So from \eqref{est:L2_norm_Rk3(t,x)}, \eqref{est_VI1UI2}, lemma \ref{Lem_appendix: est UJ} and priori estimates \eqref{est: bootstrap argument a-priori est}
\begin{equation} \label{est_SI1I1_1}
\begin{split}
|S^{I_1,I_2,l}_1|&\lesssim \left[\left(\|U(t,\cdot)\|_{H^{2,\infty}}+\|\mathrm{R}U(t,\cdot)\|_{H^{2,\infty}}\right)\|V^{I_1}(t,\cdot)\|_{L^2} + \|\mathfrak{R}^k_3(t,\cdot)\|_{L^2}\right]\\
&\times \left(\|\chi(t^{-\sigma}D_x)U^{I_2}(t,\cdot)\|_{H^{7,\infty}}+ \|\chi(t^{-\sigma}D_x)\mathrm{R} U^{I_2}(t,\cdot)\|_{H^{7,\infty}} \right)\|V^I(t,\cdot)\|_{L^2}\\
& \le CAB\varepsilon^2 t^{-1+\frac{\delta_k}{2}}E^k_{3}(t;W)^\frac{1}{2}.
\end{split}
\end{equation}
Let us now consider all the addends in $S^{I_1,I_2,l}_3$ with $|J_2|<|I|$, which by inequality \eqref{estimate integral B uvw-new} are bounded by
\[\|V^{I_1}(t,\cdot)\|_{L^2}\Big(\sum_{|\mu|=0}^1 \|\chi(t^{-\sigma}D_x)\mathrm{R}^\mu U^{I_2}(t,\cdot)\|_{H^{7,\infty}}\Big)\sum_{\substack{(J_1,J_2)\in\mathcal{I}(I)\\ |J_2|<|I|}}\Big\|c_{J_1,J_2}Q^\mathrm{kg}_0(v^{J_1}_\pm, Du^{J_2}_\pm)\Big\|_{L^2}.\]
As the latter above factor is bounded by the $L^2$ norm of $Q^I_0(V,W)$ (see definition \eqref{matrix QI}), inequalities \eqref{ineq_QI0_L2} and \eqref{est_VI1UI2} imply that those integrals are remainders $\mathfrak{G}^{I_1,I_2}_{\text{quart}}$ satisfying \eqref{est_GI1I2quart}.
Finally, the last contribution to $S^{I_1,I_2,l}_3$, corresponding to $|J_1|=0, J_2=I$, for which $D=D_1$, can be rewritten using \eqref{Fourier transform Qkg0} as the sum over $j_4,j_5\in \{+,-\}$ of
\begin{equation*}
\int B^1_{(j_1,j_2,j_3)}(\xi,\eta) \left(1+j_4j_5\frac{\xi+\zeta}{\langle\xi+\zeta\rangle}\cdot\frac{\zeta}{|\zeta|}\right)\zeta_1 \hat{v}_{4}(-\xi-\zeta) \hat{u}^I_{j_5}(\zeta) \reallywidehat{\chi(t^{-\sigma}D_x)u^{I_2}_{j_2}}(\eta) \hat{v}^{I_1}_{j_1}(\xi-\eta) d\xi d\eta.
\end{equation*}
By means of lemma \ref{Lem:Est_integrals_quartic-terms-2} it is bounded by
\begin{equation*}
\|V(t,\cdot)\|_{H^{7,\infty}}\Big(\sum_{|\mu|=0}^1\|\chi(t^{-\sigma}D_x)D_1\mathrm{R}^\mu U^{J_2}(t,\cdot)\|_{H^{11,\infty}} \Big)\|V^{I_1}(t,\cdot)\|_{H^1}\|U^I(t,\cdot)\|_{L^2}
\end{equation*}
for every $t\in [1,T]$, and hence by $CA(A+B)\varepsilon^2 t^{-\frac{3}{2}+\beta'} E^k_3(t;W)$, with $\beta'>0$ small as long as $\sigma,\delta_0$ are small, as follows by a-priori estimate \eqref{est: boostrap vpm} and lemma \ref{Lem_appendix: est UJ}.
\endproof
\end{lem}

\subsection{Propagation of the energy estimate}

\begin{prop}[Propagation of the energy estimate] \label{Prop: Propagation of the energy estimate}
Let us fix $K_2>0$.
There exist two integers $n\gg \rho\gg 1$ sufficiently large, two constants $A,B>1$ sufficiently large, $\varepsilon_0\in ]0,(2A+B)^{-1}[$ sufficiently small, and some $0<\delta \ll \delta_2\ll \delta_1\ll \delta_0\ll 1$ small such that, for any $0<\varepsilon<\varepsilon_0$, if $(u,v)$ is solution to \eqref{wave KG system}-\eqref{initial data} in some interval $[1,T]$ for a fixed $T>1$, and $u_\pm, v_\pm$ defined in \eqref{def u+- v+-} satisfy a-priori estimates \eqref{est: bootstrap argument a-priori est} for every $t\in [1,T]$, then they also verify \eqref{est: bootstrap enhanced Enn}, \eqref{est: boostrap enhanced E02} on the same interval $[1,T]$.
\proof
For any integer $k,n\in\mathbb{N}$, with $n\ge 3$ and $0\le k\le 2$, let $\widetilde{E}_n(t;W)$, $\widetilde{E}^k_3(t;W)$ be the first modified energies introduced in \eqref{modified_energies_Etilde} and $\widetilde{E}^\dagger_n(t;W)$, $\widetilde{E}^{k,\dagger}_3(t;W)$ be the second modified energies, introduced in \eqref{energy_dag_En} and \eqref{energy_dag_Ek2} respectively.
Let also $D^I_{(j_1,j_2,j_3)}, D^{I,R}_{(j_1,j_2,j_3)}$, $D^{I,T_{-N}}_{(j_1,j_2,j_3)}$ be the integrals defined in \eqref{integral Dk DR_k DT}, $D^{I_1,I_2}_{(j_1,j_2,j_3)}$ in \eqref{def of DI1I2}, and $G^{I_1,I_2,l}_{(j_1,j_2,j_3)}$ in \eqref{integral_GI1I2}. Fix $N=18$.

The first thing we observe is that, as long as estimates \eqref{est: bootstrap upm}, \eqref{est: boostrap vpm} are satisfied and $\rho\in\mathbb{N}$ is sufficiently large (e.g. $\rho\ge \max\{[\frac{n}{2}]+8,21\}$), there is a constant $C>0$ such that for every $t\in [1,T]$
\begin{subequations}
\begin{gather}
C^{-1}E_n(t;W)\le \widetilde{E}^{\dagger}_n(t;W)\le C E_n(t;W),\label{equivalence_En-Edaggern}\\
C^{-1}E^k_3(t;W)\le \widetilde{E}^{k,\dagger}_3(t;W)\le C E^k_3(t;W).\label{equivalence_Ek2-Edaggerk2}
\end{gather}
\end{subequations}
Above equivalences follow from \eqref{energy_equivalence Ekm Etilde_km}, a-priori estimates \eqref{est: bootstrap upm}, \eqref{est: boostrap vpm}, the fact that for a general multi-index $I$ ($I\in\mathcal{I}_n$ or $I\in\mathcal{I}^k_3$ for $0\le k\le 2$)
\begin{equation*}
\sum_{j_i\in\{+,-\}} \left| D^I_{(j_1,j_2,j_3)}\right| + \left| D^{I,R}_{(j_1,j_2,j_3)}\right| \lesssim \left(\|U(t,\cdot)\|_{H^{7,\infty}} + \|\mathrm{R}_1U(t,\cdot)\|_{H^{7,\infty}}\right) \|V^I(t,\cdot)\|^2_{L^2}
\end{equation*}
by inequality \eqref{estimate integral B uvw-new}, 
\begin{equation*}
\sum_{j_k\in\{+,-\}} \left| D^{I,T_{-18}}_{(j_1,j_2,j_3)}\right| \lesssim \|U(t,\cdot)\|_{H^{21,\infty}} \|W^I(t,\cdot)\|^2_{L^2}
\end{equation*}
by inequality \eqref{estimate_integral_sigmatildeN}, and:

$\bullet$ as concerns especially \eqref{equivalence_En-Edaggern}, from the fact that for any $I\in\mathcal{I}_n$, any $(I_1,I_2)\in\mathcal{I}(I)$ with $[\frac{|I|}{2}]<|I_1|<|I|$, by \eqref{estimate integral B uvw-new}
\begin{equation*}
\begin{split}
\sum_{j_i\in\{+,-\}} \left| D^{I_1,I_2}_{(j_1,j_2,j_3)}\right|&
\lesssim \left(\|U^{I_2}(t,\cdot)\|_{H^{7,\infty}}+\|\mathrm{R}_1U^{I_2}(t,\cdot)\|_{H^{7,\infty}}\right)\|V^{I_1}(t,\cdot)\|_{L^2}\|V^I(t,\cdot)\|_{L^2}\\
&\lesssim \left(\|U(t,\cdot)\|_{H^{[\frac{n}{2}]+8,\infty}} + \|\mathrm{R}_1U(t,\cdot)\|_{H^{[\frac{n}{2}]+8,\infty}}\right)E_n(t;W);\\
\end{split}
\end{equation*}

$\bullet$ as concerns especially \eqref{equivalence_Ek2-Edaggerk2}, the fact that for any $I\in\mathcal{V}^k$ (see definition \eqref{set_V}), any $(I_1,I_2)\in\mathcal{I}(I)$ with $I_1\in\mathcal{K}$ (see \eqref{set_K}) and $|I_2|\le 1$, and any $l=1,2,3$, by \eqref{estimate integral B uvw-new} and \eqref{est_VI1UI2}
\begin{equation}
\begin{split}
\sum_{j_i\in\{+,-\}}\left| G^{I_1,I_2,l}_{(j_1,j_2,j_3)}\right| & \lesssim \sum_{|\mu|=0}^1\|\chi(t^{-\sigma}D_x)\mathrm{R}^\mu U^{I_2}(t,\cdot)\|_{H^{7,\infty}} \|V^{I_1}(t,\cdot)\|_{L^2}\|V^I(t,\cdot)\|_{L^2}\\
&\le C(A+B)B\varepsilon^2 t^{-\frac{1}{2}+\frac{\delta_k}{2}}E^k_3(t;W)^\frac{1}{2}.
\end{split}
\end{equation}

Let us now consider a general multi-index $I$. From equation \eqref{equation Wtilde-Is} we deduce the following equality:
\begin{equation} \label{partial t || WIs ||}
\begin{split}
& \frac{1}{2}\partial_t \|\widetilde{W}^I_s(t,\cdot)\|^2_{L^2}= -\Im\left[\langle D_t \widetilde{W}^I_s,\widetilde{W}^I_s\rangle\right]\\
& = -\Im \left[\langle A(D)\widetilde{W}^I_s, \widetilde{W}^I_s\rangle + \Bigl\langle Op^B\Bigl((I_4+ E^0_d(U;\eta))\widetilde{A}_1(V;\eta)(I_4 + F^0_d(U;\eta))\Bigr)\widetilde{W}^I_s , \widetilde{W}^I_s \Bigr\rangle \right.\\
& \hspace{0.7cm}  \left. +\langle Op^B(A''(V^I;\eta))U + Op^B_R(A''(V^I;\eta))U , \widetilde{W}^I_s\rangle + \langle Q^I_0(V,W) , \widetilde{W}^I_s\rangle \right. \\
& \hspace{0.7cm} \left.  + \langle T_{-18}(U)W^I_s, \widetilde{W}^I_s\rangle + \langle\mathfrak{R}'(U, V),\widetilde{W}^I_s\rangle\right]
\end{split}
\end{equation}
and immediately notice that $\Im[\langle A(D)\widetilde{W}^I_s, \widetilde{W}^I_s\rangle] = 0$ because of the fact that $A(\eta)$, introduced in \eqref{matrices A A'}, is real diagonal matrix and its quantization is a self-adjoint operator. 

Matrix $(I_4+ E^0_d(U;\eta))\widetilde{A}_1(V;\eta)(I_4 + F^0_d(U;\eta))$ is real, symmetric, of order 1, with semi-norm
\begin{equation*}
M^1_1\Bigl(\big(I_4+ E^0_d(U;\eta))\widetilde{A}_1(V;\eta)(I_4 + F^0_d(U;\eta)\big), 3\Bigr) \lesssim
 (1+ \|\mathrm{R}_1U(t,\cdot)\|_{H^{2,\infty}})^2\|V(t,\cdot)\|_{H^{2,\infty}}
\end{equation*}
as follows by estimate \eqref{seminorm E0d} on $E^0_d$, \eqref{seminorm F0} of $F^0_d$, and \eqref{seminorm Atilde 1} on $\widetilde{A}_1(V;\eta)$.
Corollary \ref{Cor : paradiff ajoint at order 1} and a-priori estimates \eqref{est: bootstrap upm}, \eqref{est: boostrap vpm} imply then that the second term in the right hand side of \eqref{partial t || WIs ||} reduces to $\langle T_0(U,V) \widetilde{W}^I_s, \widetilde{W}^I_s\rangle$, with $T_0(U,V)$ operator of order less or equal than 0 such that 
\begin{equation*}
\|T_0(U,V)\|_{\mathcal{L}(L^2)}\lesssim M^1_1\Bigl(\big(I_4+ E^0_d(U;\eta))\widetilde{A}_1(V;\eta)(I_4 + F^0_d(U;\eta)\big), 3\Bigr) \le CA\varepsilon t^{-1},
\end{equation*}
so after Cauchy-Schwarz inequality and equivalence \eqref{remark on equivalence between L2 norms} it is a remainder $R(t)$ satisfying, for every $t\in [1,T]$
\begin{equation} \label{remainder_R(t)}
\left| R(t) \right|\le CA\varepsilon t^{-1}\|W^I(t,\cdot)\|^2_{L^2}.
\end{equation}
Observe that, by the definition of $\widetilde{W}^I_s$ in \eqref{def_WtildeI} and of $W^I_s$ in \eqref{def_WIs}, we have that
\begin{equation} \label{est_L2_Witildes-WI}
\begin{split}
\left\|(\widetilde{W}^I_s - W^I)(t,\cdot) \right\|_{L^2}&\le
\|Op^B(P(V;\eta)-I_4)W^I\|_{L^2} + \|Op^B(E(U;\eta))W^I_s\|_{L^2}\\
 &\lesssim \left(\|V(t,\cdot)\|_{H^{1,\infty}}+\|U(t,\cdot)\|_{H^{5,\infty}}+\|\mathrm{R}_1U(t,\cdot)\|_{H^{1,\infty}}\right)\|W^I(t,\cdot)\|_{L^2},
\end{split}
\end{equation}
the latter inequality following from proposition \ref{Prop : Paradiff action on Sobolev spaces-NEW}, estimate \eqref{seminorm P-I4}, the fact that $E(U;\eta)$ verifies, after \eqref{seminorms E} and for any admissible cut-off function $\chi$,
\begin{equation*}
M^0_0\left(E\Big(\chi\Big(\frac{D_x}{\langle \eta\rangle}\Big) U;\eta\Big);n \right)\lesssim \|U(t,\cdot)\|_{H^{5,\infty}}+\|\mathrm{R}_1U(t,\cdot)\|_{H^{1,\infty}},
\end{equation*}
and equivalence \eqref{equivalence WIs WI}.
Therefore, third and fifth brackets in the right hand side of \eqref{partial t || WIs ||} can be replaced with
\begin{equation*}
\langle Op^B(A''(V^I;\eta))U + Op^B_R(A''(V^I;\eta))U , W^I\rangle +  \langle T_{-18}(U)W^I, W^I\rangle
\end{equation*}
up to some new remainders $R(t)$, satisfying \eqref{remainder_R(t)} after Cauchy-Schwarz inequality, estimates \eqref{L2 est on OpBR(A'')W}, \eqref{norm of T-3 in propositon}, \eqref{est_L2_Witildes-WI} and \eqref{est: bootstrap upm}, \eqref{est: boostrap vpm}. 

Summing up, equality \eqref{partial t || WIs ||} reduces to:
\begin{equation}\label{partial WIs-2}
\begin{split}
\frac{1}{2}\partial_t \|\widetilde{W}^I_s(t,\cdot)\|_{L^2} = -&\Im\Big[ \langle Op^B(A''(V^I;\eta))U + Op^B_R(A''(V^I;\eta))U , W^I\rangle\\
&+ \langle Q^I_0(V,W), \widetilde{W}^I_s\rangle + \langle T_{-18}(U)W^I_s, \widetilde{W}^I_s\rangle + \langle\mathfrak{R}'(U,V), \widetilde{W}^I_s\rangle \Big] +R(t).
\end{split}
\end{equation}
In order to analyse the behaviour of the second and fourth brakets in above right hand side we need, at this point, to distinguishing between indices $I\in \mathcal{I}_n$ and $I\in\mathcal{I}^k_3$. 

\underline{\textbf{Propagation of a-priori estimate \eqref{est: bootstrap Enn}}}: Let us suppose that $I\in\mathcal{I}_n$.
Using \eqref{est_L2_Witildes-WI} and estimate \eqref{est:QI0-In} we find that
\begin{equation}\label{QI0,WtildeIS}
\langle Q^I_0(V,W), \widetilde{W}^I_s\rangle = \langle Q^I_0(V,W), W^I\rangle +R_n(t)
\end{equation}
where, for a new constant $C>0$ and every $t\in [1,T]$,
\begin{equation}  \label{remainder_Rn(t)}
|R_n(t)|\le CA\varepsilon t^{-1+\frac{\delta}{2}}E_n(t;W)^\frac{1}{2}.
\end{equation}
Reminding definition \eqref{matrix QI} of $Q^I_0(V,W)$ and the fact that coefficients $c_{I_1,I_2}$ are all equal to 0 when $I\in\mathcal{I}_n$, we notice that some of the contributions to the scalar product in the right hand side of \eqref{QI0,WtildeIS} are also remainders $R_n(t)$.
These are precisely the following ones:
\begin{equation*}
\sum_{(I_1,I_2)\in\mathcal{I}(I)}\langle Q^\mathrm{w}_0(v^{I_1}_\pm, D_1v^{I_2}_\pm), u^I_+ + u^I_{-}\rangle + \sum_{\substack{(I_1,I_2)\in\mathcal{I}(I)\\ |I_1|\le [\frac{|I|}{2}]}}\langle Q^\mathrm{kg}_0(v^{I_1}_\pm, D_1u^{I_2}_\pm), v^I_+ + v^I_{-}\rangle
\end{equation*}
in consequence of Cauchy-Schwarz inequality and estimates \eqref{est: L2 Qw0 (vI1 vI2)-only derivatives}, \eqref{est: boostrap vpm}, \eqref{est: bootstrap Enn}.
Moreover, $\langle \mathfrak{R}'(U,V), \widetilde{W}^I\rangle$ in the right hand side of \eqref{partial WIs-2} is also a remainder $R_n(t)$ because of Cauchy-Schwarz, \eqref{est_L2_Witildes-WI}, a-priori estimates \eqref{est: bootstrap upm}, \eqref{est: boostrap vpm}, and the fact that
\begin{equation*}
\|\mathfrak{R}'(U,V)\|_{L^2}\le CA\varepsilon t^{-1+\frac{\delta}{2}},
\end{equation*} 
which follows choosing $\theta\ll 1$ in \eqref{L2 norm R'(U,V)}, using \eqref{est:QI0-In} and \eqref{est: bootstrap upm}-\eqref{est: bootstrap Enn}.

Since remainder $R(t)$ in \eqref{partial WIs-2} (verifying \eqref{remainder_R(t)}) can be enclosed in $R_n(t)$ after \eqref{est: bootstrap Enn}, we obtain that equality \eqref{partial WIs-2} can be further reduced to
\begin{multline*}
\frac{1}{2}\partial_t \|\widetilde{W}^I_s(t,\cdot)\|^2_{L^2} = -\Im\Big[\langle Op^B(A''(V^I;\eta))U + Op^B_R(A''(V^I;\eta))U , W^I\rangle\\
+ \sum_{\substack{(I_1,I_2)\in\mathcal{I}(I)\\ [\frac{|I|}{2}]<|I_1|<|I|}} \langle  Q^\mathrm{kg}_0(v^{I_1}_\pm, D_1u^{I_2}_\pm), v^I_+ + v^I_{-}\rangle +  \langle T_{-18}(U)W^I, W^I\rangle \Big] + R_n(t).
\end{multline*}
From definition \eqref{energy_dag_En}, equalities \eqref{sum of CIk}, \eqref{sum of CIRk}, \eqref{expression T-N in normal forms} with $N=18$, \eqref{sum CI1I2}, together with \eqref{derivative DIs}, \eqref{derivative DIT} with $N=18$, \eqref{derivative DI1I2}, we deduce that 
\begin{multline*}
\frac{1}{2}\left| \partial_t \widetilde{E}^{\dagger}_n(t;W) \right| \lesssim |R_n(t)| + \sum_{I\in\mathcal{I}_n}\left( \left|\mathfrak{D}^{I}_{\text{quart}}(t)\right| +  \left|\mathfrak{D}^{I,18}_{\text{quart}}(t)\right| \right)
+ \sum_{I\in\mathcal{I}_n}\sum_{\substack{(I_1,I_2)\in\mathcal{I}(I)\\ [\frac{|I|}{2}]<|I_1|<|I|}} \left| \mathfrak{D}^{I_1, I_2}_{\text{quart}}(t)\right|,
\end{multline*}
where quartic terms $\mathfrak{D}^{I}_{\text{quart}}, \mathfrak{D}^{I,18}_{\text{quart}}, \mathfrak{D}^{I_1, I_2}_{\text{quart}}$ satisfy, respectively, \eqref{est_DI1_quart}, \eqref{est_DIN_quart} with $N=18$, \eqref{est_DI1I2_quart}.
These latter ones can also be considered as remainders $R_n(t)$ thanks to lemma \ref{Lem: L2 est nonlinearities} $(i)$ and a-priori estimates \eqref{est: bootstrap argument a-priori est}, which implies that, for some new $C>0$ and every $t\in [1,T]$,
\begin{equation*}
\left|\partial_t \widetilde{E}^{\dag}_n(t;W)\right| \le CA\varepsilon t^{-1+\frac{\delta}{2}} E_n(t;W)^\frac{1}{2}.
\end{equation*}
Then
\begin{equation*}
\widetilde{E}^{\dagger}_n(t;W)^\frac{1}{2}\le \widetilde{E}^{\dagger}_n(1;W)^\frac{1}{2} + \int_1^t C
A\varepsilon \tau^{-1+\frac{\delta}{2}}  d\tau ,
\end{equation*}
so after equivalence \eqref{equivalence_En-Edaggern} and a-priori estimate \eqref{est: bootstrap Enn}
\begin{equation*}
\begin{split}
E_n(t;W)^\frac{1}{2}& \le C E_n(1;W)^\frac{1}{2} + \int_1^t C A\varepsilon\tau^{-1+\frac{\delta}{2}} d\tau\\
& \le C E_n(1;W)^\frac{1}{2} + \frac{2C A\varepsilon}{\delta}t^\frac{\delta}{2},
\end{split}
\end{equation*}
again for a new $C>0$. 
Taking $B>1$ sufficiently large so that $E_n(1;W)^\frac{1}{2}\le \frac{B\varepsilon}{2C K_2}$ and $\frac{2CA}{\delta}<\frac{B}{2K_2}$ we finally obtain \eqref{est: bootstrap enhanced Enn}.

\smallskip
\underline{\textbf{Propagation of a-priori estimate \eqref{est: bootstrap E02}}}: Let us now suppose that $I\in\mathcal{I}^k_3$ for $0\le k\le 2$.
After \eqref{est:QI0-Ik3} and \eqref{est_L2_Witildes-WI} we have that
\begin{equation*}
\langle Q^I_0(V,W), \widetilde{W}^I_s\rangle = \langle Q^I_0(V,W), W^I\rangle +R^k_3(t)
\end{equation*}
with
\begin{equation}
|R^k_3(t)|\le CA(A+B)\varepsilon^2 t^{-1+\frac{\delta_k}{2}}E^k_3(t;W)^\frac{1}{2},
\end{equation}
and moreover 
\begin{equation}
\begin{split}
-\Im\left[\langle Q^I_0(V,W), W^I\rangle \right] = & -\delta_{\mathcal{V}^k}\sum_{\substack{(I_1,I_2)\in\mathcal{I}(I)\\ I_1\in\mathcal{K}, |I_2|\le 1}}c_{I_1,I_2} \Im\left[\left\langle Q^\mathrm{kg}_0\left( v^{I_1}_\pm,\chi(t^{-\sigma}D_x) D_x u^{I_2}_\pm\right), v^I_+ + v^I_{-}\right\rangle\right]\\
&  -\delta_{\mathcal{V}^k} \sum_{\substack{(J,0)\in\mathcal{I}(I)\\ J\in\mathcal{K}}}c_{J,0} \Im\left[\left\langle Q^\mathrm{kg}_0\left( v^J_\pm,\chi(t^{-\sigma}D_x) |D_x| u_\pm\right), v^I_+ + v^I_{-}\right\rangle\right]
 +R^k_3(t),
 \end{split}
\end{equation}
with $\delta_{\mathcal{V}^k}=1$ if $I\in\mathcal{V}^k$, 0 otherwise, as already seen in \eqref{decomposition-Qkg0-complete}.
Also, $\langle \mathfrak{R}'(U,V), \widetilde{W}^I_s\rangle$ in the right hand side of \eqref{partial WIs-2} and $R(t)$ are remainders $R^k_3(t)$ in consequence of the same argument used in the previous case, but with estimate \eqref{est:QI0-In} replaced with \eqref{est:QI0-Ik3}.
Therefore, we can further reduce \eqref{partial WIs-2} to the following equality:
\begin{equation*}
\begin{split}
\frac{1}{2}\partial_t \|\widetilde{W}^I_s(t,\cdot)\|^2_{L^2} =& -\Im\left[\langle Op^B(A''(V^I;\eta))U + Op^B_R(A''(V^I;\eta))U , \widetilde{W}^I\rangle + \langle T_{-18}(U)W^I, W^I\rangle\right]\\
& -\delta_{\mathcal{V}^k}\sum_{\substack{(I_1,I_2)\in\mathcal{I}(I)\\ I_1\in\mathcal{K}, |I_2|\le 1}}c_{I_1,I_2} \Im\left[\left\langle Q^\mathrm{kg}_0\left( v^{I_1}_\pm,\chi(t^{-\sigma}D_x) D u^{I_2}_\pm\right), v^I_+ + v^I_{-}\right\rangle\right] \\
&  -\delta_{\mathcal{V}^k} \sum_{\substack{(J,0)\in\mathcal{I}(I)\\ J\in\mathcal{K}}}c_{J,0} \Im\left[\left\langle Q^\mathrm{kg}_0\left( v^J_\pm,\chi(t^{-\sigma}D_x) |D_x| u_\pm\right), v^I_+ + v^I_{-}\right\rangle\right] +R^k_3(t),
\end{split}
\end{equation*}
and deduce from definition \eqref{energy_dag_Ek2}, equalities \eqref{sum of CIk}, \eqref{sum of CIRk}, \eqref{expression T-N in normal forms} with $N=18$, \eqref{sum_FI1I2}, together with \eqref{derivative DIs}, \eqref{derivative DIT} with $N=18$, and \eqref{derivative_GI1I2}, that 
\begin{equation*}
\left| \partial_t \widetilde{E}^{k,\dagger}_3(t;W) \right| \lesssim |R^k_3(t)| + \sum_{I\in\mathcal{I}^k_3}\left( \left|\mathfrak{D}^{I}_{\text{quart}}(t)\right| +  \left|\mathfrak{D}^{I,18}_{\text{quart}}(t)\right| \right) + \delta_{k<2}\sum_{\substack{I\in\mathcal{V}^k\\ j_i\in \{+,-\}}}\sum_{\substack{(I_1,I_2)\in\mathcal{I}(I)\\ I_1\in\mathcal{K}, |I_2|\le 1}} \left|\mathfrak{G}^{I_1,I_2}_{(j_1,j_2,j_3)}\right|
\end{equation*}
with $\delta_{k<2}=1$ for $k<2$, 0 otherwise. 
On the one hand, quartic terms $\mathfrak{D}^{I}_{\text{quart}}, \mathfrak{D}^{I,18}_{\text{quart}}$ satisfy, respectively, \eqref{est_DI1_quart} and \eqref{est_DIN_quart} with $N=18$, and are remainders $R^k_3(t)$ after \eqref{est:QI0-Ik3} and a-priori estimates. On the other hand, $\mathfrak{G}^{I_1,I_2}_{(j_1,j_2,j_3)}$ verifies estimate \eqref{est_GI1I2quart}.
Consequently, there is a constant $C>0$ such that
\begin{multline*}
\widetilde{E}^{k,\dagger}_3(t;W)\le \widetilde{E}^{k,\dagger}_3(1;W) + C(A+B)^2\varepsilon^2 \int_1^t  \tau^{-1+\frac{\delta_k}{2}}E^k_3(t\tau;W)^\frac{1}{2} d\tau \\
+ \delta_{k<2}\, C(A+B)^2\varepsilon^2 \left[ \delta_{k=0}\int_1^t \tau^{-1+\frac{\delta_0}{2}+\beta+\frac{\delta_1}{2}}E^1_3(\tau;W)^\frac{1}{2}d\tau+\int_1^t \tau^{-\frac{5}{4}}d\tau \right]
\end{multline*}
with $\delta_{k=0}=1$ if $k=0$, 0 otherwise, $\beta>0$ as small as we want, and after equivalence \eqref{equivalence_Ek2-Edaggerk2}
\begin{multline*}
E^k_3(t;W)\le C E^k_3(1;W) + C(A+B)^2\varepsilon^2 \int_1^t  \tau^{-1+\frac{\delta_k}{2}}E^k_3(\tau;W)^\frac{1}{2} d\tau \\
+ \delta_{k<2}\, C(A+B)^2\varepsilon^2 \left[ \delta_{k=0}\int_1^t \tau^{-1+\frac{\delta_0}{2}+\beta+\frac{\delta_1}{2}}E^1_3(\tau;W)^\frac{1}{2}d\tau  +\int_1^t \tau^{-\frac{5}{4}}d\tau \right] ,
\end{multline*}
for a new $C>0$. Injecting \eqref{est: bootstrap E02} in the above inequality and integrating in $d\tau$, we obtain that
\begin{equation*}
E^k_3(t;W)\le CE^k_3(1;W) +C(A+B)^2B\varepsilon^3\left[ \frac{1}{\delta_k} t^{\delta_k}+\delta_{k=0}\  \frac{1}{\frac{\delta_0}{2}+\beta+\delta_1}t^{\frac{\delta_0}{2}+\beta+\delta_1}\right],
\end{equation*} 
and taking $\beta$ sufficiently small so that $\beta+ \delta_1\le\delta_0/2$, $B>1$ sufficiently large so that $E^k_3(1;W)\le \frac{B^2\varepsilon^2}{2CK_2^2}$ and $B\ge A$, and $\varepsilon_0>0$ sufficiently small so that 
\[\varepsilon_0 \le \frac{1}{8BCK^2_2}\Big[\frac{1}{\delta_k}+\delta_{k=0}\frac{1}{\frac{\delta_0}{2}+\beta+\delta_1}\Big]^{-1},\]
we finally derive enhanced estimate \eqref{est: boostrap enhanced E02} and the conclusion of the proof.
\endproof
\end{prop}

\chapter{Uniform Estimates}

\section{Semilinear normal forms}\label{Section : Normal Forms for system}

In proposition \ref{Prop: Propagation of the energy estimate} of the previous chapter we proved the propagation of the a-priori the energy estimates, i.e. that there exist some constants $A,B>1$ large and $\varepsilon_0>0$ small, such that \eqref{est: bootstrap argument a-priori est} implies \eqref{est: bootstrap enhanced Enn}, \eqref{est: boostrap enhanced E02}.
To conclude the proof of theorem \ref{Thm: bootstrap argument} it only remains to show that \eqref{est: bootstrap argument a-priori est} also implies \eqref{est:bootstrap enhanced upm}, \eqref{est:bootstrap enhanced vpm}.
In particular, as $u_+ = -\overline{u_{-}}$ and $v_+ =-\overline{u_{-}}$, it will be enough to prove this result for $(u_{-},v_{-})$, which is solution to
\begin{equation} \label{wave-KG for u- v-}
\begin{cases}
& \left(D_t + |D_x|\right) u_{-} = Q_0^\mathrm{w}(v_\pm, D_1 v_\pm), \\
& \left(D_t + \langle D_x\rangle\right) v_{-} = Q_0^{\mathrm{kg}}(v_\pm, D_1 u_\pm),
\end{cases}
\end{equation}
with $Q_0^\mathrm{w}(v_\pm, D_1 v_\pm), Q_0^{\mathrm{kg}}(v_\pm, D_1 u_\pm)$ given by \eqref{Q0_pm}.

As for the simpler case of the one-dimensional Klein-Gordon equation (see \cite{stingo:1D_KG}), the main idea is to reformulate system \eqref{wave-KG for u- v-} in terms of two new functions $\widetilde{u}, \widetilde{v}$, defined from $u_{-}, v_{-}$ and living in a new framework (the \textit{semi-classical framework}), and to deduce a new simpler system, made of a transport equation and an ODE.
Through this new system we will be able to recover the required enhanced estimates \eqref{est:bootstrap enhanced upm}, \eqref{est:bootstrap enhanced vpm}.

Before introducing the semi-classical framework in which we will work for the rest of the paper, we need to replace almost all quadratic non-linearities in \eqref{wave-KG for u- v-} with cubic ones by a normal forms. This is the object of the following two subsections.
We highlight the fact that we do not make use directly of the normal forms obtained in the proof of the energy inequality, because our goals and constraints are henceforth different. In fact, we want to obtain a $L^\infty$ estimate for essentially $\rho$ derivatives of our solution, having a control on its $H^s$ norm for $s\gg \rho$.
Therefore, we are allowed to lose some derivatives in the normal form reduction, which means that we do not care any more about the quasi-linear nature of our problem.

We warn the reader that, for seek of compactness, we will often use the notation $\Nlw$ (resp. $\Nlkg$) when referring to $Q^\mathrm{w}_0(v_\pm, D_1 v_\pm)$ (resp. to $Q^\mathrm{kg}_0(v_\pm, D_1 u_\pm)$).

\subsection{Normal forms for the Klein-Gordon equation}

The aim of this subsection is to introduce a new unknown $\vnf$, defined in terms of $v_{-}$, in such a way it is solution to a cubic half Klein-Gordon equation instead of the quadratic one satisfied by $v_{-}$ in \eqref{wave-KG for u- v-}.
This normal form is motivated by the fact that the $L^2$ norm of $Q^\mathrm{kg}_0(v_\pm, D_1u_\pm)$ decays too slowly in time (only $t^{-1+\delta/2}$), as follows from \eqref{est L2 NLkg} and a-priori estimates \eqref{est: boostrap vpm}, \eqref{est: bootstrap Enn}, and this decay is not enough in view of proposition \ref{Prop:propagation_unif_est_V} (the required one being strictly faster than $t^{-3/2}$).

It is fundamental to observe that, after \eqref{est: bootstrap argument a-priori est} and inequality \eqref{est_Hsinfty_vnf-v} below with $\theta\ll 1$ small enough (e.g. $\theta<(2+\delta)^{-1}$), $\vnf$ and $v_{-}$ are comparable, in the sense that there is a positive constant $C$ such that
\begin{equation} \label{equivalence v- vnf}
\left| \|v_{-}(t,\cdot)\|_{H^{\rho,\infty}} -  \|\vnf(t,\cdot)\|_{H^{\rho,\infty}}\right| \le C\varepsilon^2 t^{-1}.
\end{equation}
Then bootstrap assumption \eqref{est: boostrap vpm} implies that the new unknown $\vnf$ disperses in time at the same rate $t^{-1}$ as $v_{-}$, and the propagation of a suitable estimate of the $H^{\rho,\infty}$ norm of $\vnf$ will provide us with enhanced \eqref{est:bootstrap enhanced vpm}.

\begin{prop} \label{Prop: normal forms on KG}
Assume that $(u,v)$ is solution to \eqref{wave KG system} in $[1,T]$ for a fixed $T>1$, consider $(u_+, v_+, u_{-}, v_{-})$ defined in \eqref{def u+- v+-} and solution to \eqref{system for uI+-, vI+-} with $|I|=0$, and remind definition \eqref{def W,V, U} of vectors $U,V$. Let \index{vnF@$v^{NF}$, Klein-Gordon component after a normal form}
\begin{equation} \label{def vNF}
v^{NF}:= v_{-} - \frac{i}{4(2\pi)^2}\sum_{j_1,j_2\in\{+,-\}} \int e^{ix\cdot\xi} B^1_{(j_1,j_2,+)}(\xi,\eta) \hat{v}_{j_1}(\xi-\eta) \hat{u}_{j_2}(\eta) d\xi d\eta,
\end{equation}
with $B^1_{(j_1,j_2,+)}(\xi,\eta)$ given by \eqref{def of B(i1,i2,i3)} with $k=1$ and $j_3=+$.
Then for every $t\in [1,T]$ $v^{NF}$ is solution to \index{rNFkg@$r^{NF}_{kg}$, cubic term in the Klein-Gordon equation after a normal form}
\begin{equation} \label{KG equation vNF}
\left(D_t + \langle D_x\rangle\right) v^{NF}(t,x) = r^{NF}_{kg}(t,x),
\end{equation}
where
\begin{multline} \label{def rNF-kg}
r^{NF}_{kg}(t,x) 
= -\frac{i}{4(2\pi)^2}\sum_{j_1,j_2\in\{+,-\}} \int e^{ix\cdot\xi}  B^1_{(j_1,j_2,+)}(\xi,\eta) \\
\times \left[\widehat{\textit{NL}_{kg}}(\xi - \eta) \hat{u}_{j_2}(\eta)+  \hat{v}_{j_1}(\xi-\eta)\widehat{\textit{NL}_{w}}(\eta) \right]d\xi d\eta
\end{multline} 
satisfies 
\begin{subequations} \label{est L2 Linfty (cut-off) rNF-kg-new}
\begin{equation} \label{est_rnfkg_L2}
\begin{split}
\|r^{NF}_{kg}(t,\cdot)\|_{L^2}& \lesssim \sum_{\mu=0}^1 \|V(t,\cdot)\|_{H^{1,\infty}}\|\mathrm{R}^\mu_1 U(t,\cdot)\|_{L^\infty}\|U(t,\cdot)\|_{H^1} + \|V(t,\cdot)\|^2_{H^{2,\infty}}\|V(t,\cdot)\|_{H^2},
\end{split}
\end{equation}\normalsize
\begin{equation} \label{est_rnfkg_Linfty}
\begin{split}
& \|\chi(t^{-\sigma}D_x)r^{NF}_{kg}(t,\cdot)\|_{L^\infty} \lesssim  \|V(t,\cdot)\|_{H^{1,\infty}}\Big(\sum_{\mu=0}^1\|\mathrm{R}^\mu_1 U(t,\cdot)\|_{H^{2,\infty}}\Big)^2 + t^\sigma\|V(t,\cdot)\|^3_{H^{2,\infty}},
 \end{split}
\end{equation}
\end{subequations} 
for any $\chi\in C^\infty_0(\mathbb{R}^2)$, $\sigma>0$.
Moreover, for every $s,\rho\ge 0$, any $\theta \in ]0,1[$,
\begin{subequations}\label{est_vNF-v-}
\begin{equation} \label{est_Hs_vnf-v}
\left\| (\vnf- v_{-})(t,\cdot)\right\|_{H^s}\lesssim \sum_{\mu=0}^1\|V(t,\cdot)\|_{H^s} \|\mathrm{R}_1^\mu U(t,\cdot)\|_{L^\infty} + \|V(t,\cdot)\|_{L^\infty}\|U(t,\cdot)\|_{H^{s+1}},
\end{equation}
\begin{equation}\label{est_Hsinfty_vnf-v}
\begin{split}
\left\| (\vnf- v_{-})(t,\cdot)\right\|_{H^{s, \infty}}& \lesssim \sum_{\mu=0}^1\|V(t,\cdot)\|^{1-\theta}_{H^{s,\infty}}\|V(t,\cdot)\|^\theta_{H^{s+2}}\|\mathrm{R}^\mu_1U(t,\cdot)\|_{L^\infty}\\
& +\sum_{\mu=0}^1 \|V(t,\cdot)\|_{L^\infty}\|\mathrm{R}^\mu_1 U(t,\cdot)\|^{1-\theta}_{H^{s+1,\infty}}\|U(t,\cdot)\|^\theta_{H^{s+3}},
\end{split}
\end{equation}
\begin{equation} \label{est_Omega_vnf-v}
\begin{split}
\left\| \Omega(\vnf- v_{-})(t,\cdot)\right\|_{L^2}&\lesssim \sum_{\mu,\nu=0}^1 \left[\|\Omega^\mu V(t,\cdot)\|_{L^2}\|\mathrm{R}^\nu_1U(t,\cdot)\|_{L^\infty} + \|V(t,\cdot)\|_{L^\infty}\|\Omega^\mu U(t,\cdot)\|_{H^1}\right]\\
& + \|\Omega V(t,\cdot)\|_{H^2}\|U(t,\cdot)\|_{H^1} + \|V(t,\cdot)\|_{L^2}\|\Omega U(t,\cdot)\|_{H^2},
\end{split}
\end{equation}
\end{subequations}
and
\begin{subequations}
\begin{equation}\label{est_chi_vnf-v}
\left\|\chi(t^{-\sigma}D_x)(\vnf - v_{-})(t,\cdot)\right\|_{L^2}\lesssim t^\sigma \|V(t,\cdot)\|_{H^{1,\infty}}\|U(t,\cdot)\|_{L^2},
\end{equation}
\begin{multline}\label{est_chi_Omega_vnf-v}
\left\|\chi(t^{-\sigma}D_x)\Omega(\vnf - v_{-})(t,\cdot)\right\|_{L^2}\\
\lesssim t^\sigma \left[\sum_{\mu=0}^1 \|\Omega V(t,\cdot)\|_{L^2}\|\mathrm{R}^\mu_1 U(t,\cdot)\|_{L^\infty} + \|V(t,\cdot)\|_{H^{1,\infty}}\|\Omega^\mu U(t,\cdot)\|_{L^2}\right].
\end{multline}
\end{subequations}
\proof
From definition \eqref{def vNF} of $v^{NF}$, system \eqref{system for uI+-, vI+-} with $|I|=0$, and the fact that
\begin{equation}
Q^{\mathrm{kg}}_0(v_\pm, D_1u_\pm) = \frac{i}{4(2\pi)^2}\sum_{j_1,j_2\in \{+,-\}}\int e^{ix\cdot\xi}\Big(1-j_1j_2 \frac{\xi-\eta}{\langle \xi-\eta\rangle}\cdot\frac{\eta}{|\eta|}\Big)\eta_1 \hat{v}_{j_1}(\xi-\eta) \hat{u}_{j_2}(\eta) d\xi d\eta,
\end{equation}
it immediately follows that $v^{NF}$ is solution to \eqref{KG equation vNF} with $r^{NF}_{kg}$ given by \eqref{def rNF-kg}.
We observe that, after formula \eqref{explicit integral B}, we have the following explicit expressions:
\begin{equation}\label{explicit_vNF-v_chapter5}
\begin{split}
v^{NF}-v_{-} =-\frac{i}{8}&\Big[(v_+ +v_{-})\mathrm{R}_1(u_+ - u_{-}) - \frac{D_1}{\langle D_x\rangle}(v_+ - v_{-}) (u_+ +u_{-}) \\
&+ D_1\big[[\langle D_x\rangle^{-1}(v_+-v_{-})](u_+ + u_{-})\big] 
- \langle D_x\rangle \big[[\langle D_x\rangle^{-1}(v_+ -v_{-})]\mathrm{R}_1 (u_+ - u_{-})\big]\Big]
\end{split}
\end{equation}
and
\begin{equation} \label{explicit_rnfkg_chapter5}
\rnfkg = -\frac{i}{4}\left[\Nlkg\, \mathrm{R}_1(u_+-u_{-}) - \frac{D_x}{\langle D_x\rangle}(v_+ - v_{-})\, \Nlw + D_1\big[\langle D_x\rangle^{-1}(v_+-v_{-})\, \Nlw\big] \right].
\end{equation}
Inequalities \eqref{est_Hs_vnf-v}, \eqref{est_Hsinfty_vnf-v} are straightforward from \eqref{explicit_vNF-v_chapter5} and corollary \ref{Cor_appendixA:Hs-Hsinfty norm of bilinear expressions} in appendix \ref{Appendix A}.
Inequality \eqref{est_Omega_vnf-v} is also obtained from corollary \ref{Cor_appendixA:Hs-Hsinfty norm of bilinear expressions}, after having applied $\Omega$ to \eqref{explicit_vNF-v_chapter5} and used the Leibniz rule, and from bounding the $L^\infty$ norm of $\Omega u_\pm, \Omega v_\pm$ with their $H^2$ norm by means of the classical Sobolev injection.
Inequalities \eqref{est_chi_vnf-v}, \eqref{est_chi_Omega_vnf-v} are also straightforward if one observes that operator $\chi(t^{-\sigma}D_x)$, with $\chi\in C^\infty_0(\mathbb{R}^2)$ and $\sigma>0$, is $L^2-H^1$ continuous with norm $O(t^\sigma)$.

As concerns $\rnfkg$, from \eqref{explicit_rnfkg_chapter5} and corollary \ref{Cor_appendixA:Hs-Hsinfty norm of bilinear expressions} we find that
\begin{equation*}
\begin{split}
\|r^{NF}_{kg}(t,\cdot)\|_{L^2} &\lesssim \sum_{\mu=0}^1\|\textit{NL}_{kg}(t,\cdot)\|_{L^2}\|\mathrm{R}^\mu_1U(t,\cdot)\|_{L^\infty}+ \|V(t,\cdot)\|_{L^2}\|\Nlw(t,\cdot)\|_{L^\infty}\\
& + \|V(t,\cdot)\|_{L^\infty}\|\textit{NL}_w(t,\cdot)\|_{H^1}
\end{split}
\end{equation*}
and
\begin{equation*}
\|\chi(t^{-\sigma}D_x) r^{NF}_{kg}(t,\cdot)\|_{L^\infty}\lesssim \sum_{\mu=0}^1\|\textit{NL}_{kg}(t,\cdot)\|_{L^\infty}\|\mathrm{R}^\mu_1U(t,\cdot)\|_{L^\infty} + t^\sigma \|V(t,\cdot)\|_{H^{1,\infty}}\|\Nlw(t,\cdot)\|_{L^\infty}.
\end{equation*}
Inequalities \eqref{est_rnfkg_L2} and \eqref{est_rnfkg_Linfty} follow then by \eqref{est Hs NLw-New} with $s=1$, \eqref{est Linfty NLw}, \eqref{est L2 NLkg} and \eqref{est Linfty NLkg}.
\endproof
\end{prop}

\subsection{Normal forms for the wave equation} \label{Subsection: Section : Normal Forms for the Wave Equation}

We now focus on the wave equation satisfied by $u_{-}$:
\begin{equation*} 
(D_t + |D_x| )u_{-}(t,x) = Q^{\mathrm{w}}_0(v_\pm, D_1 v_\pm),
\end{equation*}
and perform a normal form argument in order to replace (a part of) the quadratic non-linearity in the above right hand side with a cubic non-local one. The fundamental reason for that is to be found in lemma \ref{Lemma : estimate of e(x,xi)}, where we have to prove that the $L^2$ norm of some operator, acting on the non-linearity of the equation satisfied by $u_{-}$, decays like $t^{-1/2+\beta}$, for a small $\beta>0$. 
That decay is obtained if the $L^2$ norm of the mentioned non-linearity is a $O(t^{-3/2+\beta'})$, for some new small $\beta'>0$, which is not the case for $Q^{\mathrm{w}}_0(v_\pm, D_1 v_\pm)$, as follows from \eqref{est L2 NLw}, \eqref{est: boostrap vpm}, \eqref{est: bootstrap Enn}.
This normal form can be actually performed only on contributions depending on $(v_+, v_+)$ and $(v_{-}, v_{-})$ but not on $(v_+,v_{-})$, which are resonant.
Nevertheless, the structure of these latter contributions allows us to recover the right mentioned time decay for their $L^2$ norm (see lemmas \ref{Lem: hD|V|2} and \ref{Lem: (xi+dphi)Op(gamma)}).

Thanks to inequalities \eqref{Hrho-infty norm uNF- u-}, \eqref{Hrho-infty norm R(uNF-u-)} and a-priori estimates \eqref{est: bootstrap argument a-priori est}, the new unknown $\unf$ we define in \eqref{def uNF} below is equivalent to the former $u_{-}$, meaning that there exists a positive constant $C$ such that
\begin{equation}\label{equivalence u- unf}
\sum_{\kappa=0}^1 \left| \|\mathrm{R}^\kappa_1 u_{-}(t,\cdot)\|_{H^{\rho+1,\infty}} - \|\mathrm{R}^\kappa_1\unf(t,\cdot)\|_{H^{\rho,\infty}} \right| \le C\varepsilon^2 t^{-1+\frac{\delta}{2}}.
\end{equation}
After \eqref{est: bootstrap upm} this means that $\unf$ and $\mathrm{R}_1\unf$ decay in the $H^{\rho+1,\infty}$ norm at the same rate $t^{-1/2}$ as $u_{-},\mathrm{R}_1u_{-}$, and the propagation of a suitable estimate of this norm will provide us with enhanced \eqref{est:bootstrap enhanced upm}. 

Let us rewrite $Q^{\mathrm{w}}_0(v_\pm, D_1 v_\pm)$ as follows
\begin{equation} \label{decomposition Q0w}
\begin{split}
Q_0^{\mathrm{w}}(v_\pm, D_1 v_\pm) &= -\frac{1}{2}\Im\left[v_+\, D_1v_{-} +\frac{D_x}{\langle D_x\rangle}v_+ \cdot \frac{D_xD_1}{\langle D_x\rangle}v_{-}\right]\\
& + \frac{i}{4(2\pi)^2} \sum_{j\in\{+,-\}}\int e^{ix\cdot\xi}\left(1-\frac{\xi-\eta}{\langle\xi-\eta\rangle}\cdot\frac{\eta}{\langle\eta\rangle}\right)\eta_1 \hat{v}_j(\xi-\eta)\hat{v}_j(\eta) d\xi d\eta,
\end{split}
\end{equation}
and introduce, for any $j\in\{+,-\}$,
\begin{equation} \label{def Dj1j2(xi,eta)}
D_j(\xi,\eta):= \frac{\left(1-\frac{\xi-\eta}{\langle\xi-\eta\rangle}\cdot\frac{\eta}{\langle\eta\rangle}\right)\eta_1}{j\langle\xi-\eta\rangle + j\langle\eta\rangle + |\xi|}.
\end{equation}

\begin{prop} \label{Prop: NF on wave}
Assume that $(u,v)$ is solution to \eqref{wave KG system} in $[1,T]$ for a fixed $T>1$, consider $(u_+, v_+, u_{-}, v_{-})$ defined in \eqref{def u+- v+-} and solution to \eqref{system for uI+-, vI+-} with $|I|=0$, remind definition \eqref{def W,V, U} of vectors $U,V$ and \eqref{def vNF} of $\vnf$.
Let \index{uNF@$u^{NF}$, wave component after a normal form}
\begin{equation} \label{def uNF}
u^{NF} := u_{-} - \frac{i}{4(2\pi)^2}\sum_{j\in\{+,-\}} \int e^{i x\cdot \xi}D_j(\xi,\eta)\hat{v}_j(\xi - \eta)\hat{v}_j(\eta)d\xi d\eta,
\end{equation}
with multiplier $D_j$ defined in \eqref{def Dj1j2(xi,eta)}.
For every $t\in [1, T]$ $u^{NF}$ is solution to 
\begin{equation} \label{wave equation uNF}
(D_t + |D_x|)u^{NF}(t,x) = q_w(t,x)+c_w(t,x)+ r^{NF}_w(t,x),
\end{equation}
where quadratic term $q_w$ is given by\index{qw@$q_w$, quadratic term in the wave equation after a normal form}
\begin{equation} \label{def_qw}
q_w(t,x) = \frac{1}{2}\Im\left[\overline{\vnf}\,  D_1\vnf - \overline{\frac{D_x}{\langle D_x\rangle}\vnf}\cdot\frac{D_xD_1}{\langle D_x\rangle}\vnf\right],
\end{equation}
while cubic terms $c_w, r^{NF}_w$ are equal, respectively, to\index{cw@$c_w$, cubic term in the wave equation after a normal form}
\begin{equation}\label{def_cw}
\begin{split}
c_w(t,x)=\frac{1}{2}\Im& \left[\overline{(v_{-}-\vnf)}\, D_1v_{-} + \overline{\vnf}\, D_1(v_{-}-\vnf) \right.\\
&\left. - \overline{\frac{D_x}{\langle D_x\rangle}(v_{-}-\vnf)}\cdot\frac{D_xD_1}{\langle D_x\rangle}v_{-} - \overline{\frac{D_x}{\langle D_x\rangle}\vnf}\cdot \frac{D_xD_1}{\langle D_x\rangle}(v_{-} - \vnf)\right],
\end{split}
\end{equation}
and\index{rNFw@$ r^{NF}_w$, cubic term in the wave equation after a normal form}
\begin{equation} \label{def rNF}
 r^{NF}_w(t,x) =  -\frac{i}{4(2\pi)^2}\sum_{j\in\{+,-\}}\int e^{ix\cdot \xi} D_j(\xi, \eta) \left[\widehat{\textit{NL}_{kg}}(\xi - \eta) \hat{v}_j(\eta)+  \hat{v}_j(\xi-\eta)\widehat{\textit{NL}_{kg}}(\eta) \right]d\xi d\eta.
\end{equation}
For any $s,\rho\ge 0$, any $t\in [1,T]$,
\begin{subequations} \label{norms uNF - u-}
\begin{equation} \label{Hs norm uNF- u-}
\|u^{NF}(t,\cdot) - u_{-}(t,\cdot)\|_{H^s} \lesssim \|V(t,\cdot)\|_{L^\infty}\|V(t,\cdot)\|_{H^{s+15}},
\end{equation}
\begin{equation} \label{Hrho-infty norm uNF- u-}
\|u^{NF}(t,\cdot) - u_{-}(t,\cdot)\|_{H^{\rho+1,\infty}} \lesssim \|V(t,\cdot)\|_{L^\infty}\|V(t,\cdot)\|_{H^{\rho+18}} ,
\end{equation}
\begin{equation} \label{Hrho-infty norm R(uNF-u-)}
\|\mathrm{R}_j u^{NF}(t,\cdot) -\mathrm{R}_j u_{-}(t,\cdot)\|_{H^{\rho+1,\infty}} \lesssim \|V(t,\cdot)\|_{L^\infty}\|V(t,\cdot)\|_{H^{\rho+8}}, \quad j=1,2.
\end{equation}
\end{subequations}
Moreover, for any cut-off function $\chi\in C^\infty_0(\mathbb{R}^2)$ and $\sigma>0$ there exists some $\chi_1\in C^\infty_0(\mathbb{R}^2)$ and $s>0$ such that
\begin{subequations} \label{est_cw}
\begin{equation}\label{est_cw_L2}
\begin{split}
\left\|\chi(t^{-\sigma}D_x)c_w(t,\cdot)\right\|_{L^2}&\lesssim  t^\sigma \left\|\chi_1(t^{-\sigma}D_x)(\vnf-v_{-})(t,\cdot)\right\|_{L^2}\left(\|V(t,\cdot)\|_{H^{2,\infty}}+\|\vnf(t,\cdot)\|_{H^{1,\infty}}\right) \\
& + t^{-N(s)}\left\| (\vnf - v_{-})(t,\cdot)\right\|_{H^1}\left(\|V(t,\cdot)\|_{H^s}+\|\vnf(t,\cdot)\|_{H^s}\right),
\end{split}
\end{equation}
\begin{equation} \label{est_cw_Linfty}
\begin{split}
\left\|\chi(t^{-\sigma}D_x)c_w(t,\cdot)\right\|_{L^\infty}
& \lesssim  t^\sigma \left\|\chi_1(t^{-\sigma}D_x)\left(\vnf-v_{-}\right)(t,\cdot)\right\|_{L^\infty}\left(\|V(t,\cdot)\|_{H^{2,\infty}}+\|\vnf(t,\cdot)\|_{H^{1,\infty}}\right) \\
& + t^{-N(s)}\left\| (\vnf - v_{-})(t,\cdot)\right\|_{H^1}\left(\|V(t,\cdot)\|_{H^s}+\|\vnf(t,\cdot)\|_{H^s}\right)
\end{split}
\end{equation}
\begin{equation}\label{est_Omega_cw}
\begin{split}
\left\|\chi(t^{-\sigma}D_x)\Omega c_w(t,\cdot)\right\|_{L^2}&\lesssim t^\sigma \left\|\chi_1(t^{-\sigma}D_x)\Omega (\vnf - v_{-})(t,\cdot)\right\|_{L^2} \left(\|V(t,\cdot)\|_{H^{2,\infty}}+\|\vnf(t,\cdot)\|_{H^{1,\infty}}\right)\\
& + t^{-N(s)}\left\|\Omega (\vnf - v_{-})(t,\cdot)\right\|_{L^2} \left(\|V(t,\cdot)\|_{H^s}+\|\vnf(t,\cdot)\|_{H^s}\right)\\
& + t^\sigma \left\|(\vnf - v_{-})(t,\cdot)\right\|_{H^{1,\infty}}\times \sum_{\mu=0}^1\left(\|\Omega^\mu V(t,\cdot)\|_{H^1}+\|\Omega^\mu \vnf(t,\cdot)\|_{L^2}\right)
\end{split}
\end{equation}
\end{subequations}
with $N(s)>0$ as large as we want as long as $s>0$ is large, and
\begin{subequations} \label{est rNF}
\begin{equation} \label{est L2 rNF}
\|\chi(t^{-\sigma} D_x) r^{NF}_w(t,\cdot)\|_{L^2} \lesssim  \|V(t,\cdot)\|^2_{H^{13,\infty}}\|U(t,\cdot)\|_{H^1},
\end{equation}
\begin{equation} \label{est Linfty rNF}
\|\chi(t^{-\sigma} D_x) r^{NF}_w(t,\cdot)\|_{L^\infty}\lesssim  \|V(t,\cdot)\|^2_{H^{13,\infty}}\left(\|U(t,\cdot)\|_{H^{2,\infty}}+\|\mathrm{R}_1U(t,\cdot)\|_{H^{2,\infty}}\right),
\end{equation} 
and for any $\theta\in ]0,1[$,
\begin{equation} \label{est phi(D) Omega rNF-new}
\begin{split}
 \|\chi(t^{-\sigma} D_x)& \Omega r^{NF}_w(t,\cdot)\|_{L^2} \lesssim t^\beta\Big[ \|V(t,\cdot)\|^{1-\theta}_{H^{15,\infty}} \|V(t,\cdot)\|^\theta_{H^{17}}\left(\| U(t,\cdot)\|_{H^{1,\infty}}+ \|\mathrm{R}_1U(t,\cdot)\|_{H^{1,\infty}}\right)\\
&+ \|V(t,\cdot)\|_{L^\infty}\left(\| U(t,\cdot)\|^{1-\theta}_{H^{16,\infty}}+ \|\mathrm{R}_1U(t,\cdot)\|^{1-\theta}_{H^{16,\infty}}\right)\|U(t,\cdot)\|^\theta_{H^{18}}\Big]\|\Omega V(t,\cdot)\|_{L^2}\\
&+ t^\beta \Big[\|V(t,\cdot)\|_{H^{1,\infty}}\left(\|U(t,\cdot)\|_{H^1} + \|\Omega U(t,\cdot)\|_{H^1}\right)\\
& + \left(\|U(t,\cdot)\|_{H^{2,\infty}}+\| \mathrm{R}_1U(t,\cdot)\|_{H^{2,\infty}}\right)\left(\|V(t,\cdot)\|_{L^2}+\|\Omega V(t,\cdot)\|_{L^2}\right)\Big] \|V(t,\cdot)\|_{H^{17,\infty}}.
 \end{split}
\end{equation}
\end{subequations}\normalsize
\proof
By definition \eqref{def uNF} of $u^{NF}$, system \eqref{system for uI+-, vI+-} with $|I|=0$, \eqref{decomposition Q0w} and \eqref{def Dj1j2(xi,eta)}, it follows that $u^{NF}$ is solution to 
\begin{equation*}
(D_t + |D_x|)u^{NF}(t,x) =-\frac{1}{2}\Im\left[v_+\, D_1v_{-} +\frac{D_x}{\langle D_x\rangle}v_+ \cdot \frac{D_xD_1}{\langle D_x\rangle}v_{-}\right]+ r^{NF}_w(t,x),
\end{equation*}
with $r^{NF}_w$ given by \eqref{def rNF}. Reminding that $v_+=-\overline{v_{-}}$ and replacing each occurrence of $v_{-}$ in the quadratic contribution to the above right hand side, we find that $\unf$ is solution to \eqref{wave equation uNF}.

The first part of lemma \ref{Lem_Appendix: est on Dj1j2} and the fact that any $H^{\rho+1,\infty}$ injects into $H^{\rho+3}$ by Sobolev inequality immediately imply estimates \eqref{norms uNF - u-} and
\begin{equation*}
\|\chi(t^{-\sigma} D_x) r^{NF}_w(t,\cdot)\|_{L^2}\lesssim \|\textit{NL}_{kg}(t,\cdot)\|_{L^2}\|V(t,\cdot)\|_{H^{13,\infty}},
\end{equation*}
\begin{equation*}
\|\chi(t^{-\sigma} D_x) r^{NF}_w(t,\cdot)\|_{L^\infty}\lesssim  \|\textit{NL}_{kg}(t,\cdot)\|_{L^\infty}\|V(t,\cdot)\|_{H^{13,\infty}},
\end{equation*}
for any $s,\rho\ge 0$. Moreover, from \eqref{est: L2 Omega integral D with cut-off} we derive that
\begin{multline*}
\|\chi(t^{-\sigma} D_x) \Omega r^{NF}_w(t,\cdot)\|_{L^2} \lesssim t^\beta \left(\|\textit{NL}_{kg}(t,\cdot)\|_{L^2}+ \|\Omega \textit{NL}_{kg}(t,\cdot)\|_{L^2}\right)\|V(t,\cdot)\|_{H^{17,\infty}} \\
+ t^\beta \|\textit{NL}_{kg}(t,\cdot)\|_{H^{15,\infty}}\|\Omega V(t,\cdot)\|_{L^2},
\end{multline*}
so estimates \eqref{est rNF} are obtained using \eqref{est L2 NLkg}, \eqref{est Hsinfty NLkg-new} with $s=15$, and \eqref{est L2 Omega NLkg}.

Finally, inequality \eqref{est_cw_L2} (resp. \eqref{est_cw_Linfty}) is obtained using lemma \ref{Lem_appendix:L_estimate of products} in appendix \ref{Appendix B} with $L=L^2$ (resp. $L=L^\infty$), $w=v_{-}-\vnf$, and the fact that $\chi_1(t^{-\sigma}D_x)$ is continuous from $L^2$ to $H^1$ (resp. from $L^\infty$ to $H^{1,\infty}$) with norm $O(t^\sigma)$.
Inequality \eqref{est_Omega_cw} is deduced applying $\Omega$ to \eqref{def_cw} and using the Leibniz rule. The $L^2$ norm of products in which $\Omega$ is acting on $v_{-} -\vnf$ is estimated by means of lemma \ref{Lem_appendix:L_estimate of products} with $L=L^2$, $w=v_{-}-\vnf$, whereas the $L^2$ norm of the remaining products is simply estimated by taking the $L^\infty$ norm on $v_{-} - \vnf$ times the $L^2$ norm of the remaining factor.
\endproof
\end{prop}

\section{From PDEs to ODEs} \label{Sec: development of the PDE system}

In the previous section we showed that, if $(u_{-},v_{-})$ is solution to system \eqref{wave-KG for u- v-} in some interval $[1,T]$, for a fixed $T>1$, one can define two new functions, $v^{NF}$ as in \eqref{def vNF} and $u^{NF}$ as in \eqref{def uNF}, respectively comparable to $v_{-}$ and $u_{-}$ in the sense of \eqref{equivalence v- vnf} and \eqref{equivalence u- unf}, such that $(u^{NF}, v^{NF})$ is solution to a new wave-Klein-Gordon system:
\begin{equation} \label{wave_KG NF system}
\begin{cases}
\left(D_t + |D_x|\right)u^{NF}(t,x) =q_w(t,x)+c_w(t,x)+ r^{NF}_w(t,x),\\
\left(D_t + \langle D_x\rangle\right)v^{NF}(t,x) = r^{NF}_{kg}(t,x),
\end{cases}
\end{equation}
for every $(t,x)\in [1,T]\times\mathbb{R}^2$,
where quadratic inhomogeneous term $q_w$ is given by \eqref{def_qw} and cubic ones $c_w$, $r^{NF}_w$ and $r^{NF}_{kg}$ respectively by \eqref{def_cw}, \eqref{def rNF} and \eqref{def rNF-kg}.

As anticipated before, our aim is to deduce from \eqref{wave_KG NF system} a system made of a transport equation and an ODE, from which it will be possible to deduce suitable estimates on $(u^{NF},v^{NF})$ (and consequently on $(u_{-},v_{-})$). Thanks to \eqref{equivalence v- vnf} and \eqref{equivalence u- unf} these estimates will allow us to close the bootstrap argument and prove theorem \ref{Thm: bootstrap argument}.

In subsection \ref{Subsection : The Derivation of the ODE Equation} we focus on the deduction of the mentioned ODE starting from the Klein-Gordon equation satisfied by $v^{NF}$, while in subsection
\ref{Subsection : The Derivation of the Transport Equation} we show how to derive a transport equation from the wave equation satisfied by $u^{NF}$. The framework in which this plan takes place is the \textit{semi-classical framework}, introduced below. \index{h@$h$, semi-classical parameter}\index{utilde@$\widetilde{u}$, wave component in semi-classical setting}\index{vtilde@$\widetilde{v}$, Klein-Gordon component in semi-classical setting}

Let us introduce the \textit{semi-classical parameter} $h:=t^{-1}$ together with the following two new functions:
\begin{equation} \label{def utilde vtilde}
\widetilde{u}(t,x) := t u^{NF}(t,tx), \qquad \widetilde{v}(t,x) := t v^{NF}(t,tx),
\end{equation}
and observe that, from definition \eqref{def utilde vtilde} and inequalities \eqref{equivalence u- unf}, \eqref{equivalence v- vnf}, a-priori estimates \eqref{est: bootstrap upm}, \eqref{est: boostrap vpm} are equivalent respectively to
\begin{subequations}
\begin{gather}
\|\widetilde{u}(t,\cdot)\|_{H^{\rho+1,\infty}_h} + \left\|\oph(\xi|\xi|^{-1})\widetilde{u}(t,\cdot)\right\|_{H^{\rho+1,\infty}_h}\le C\varepsilon h^{-\frac{1}{2}}, \label{est:a-priori_ut}\\
\|\vt(t,\cdot)\|_{H^{\rho,\infty}_h}\le C\varepsilon,\label{est:a-priori_vt}
\end{gather}
\end{subequations}
for some positive constant $C$.
A suitable propagation of the above estimates will therefore provide us with \eqref{est:bootstrap enhanced upm} and \eqref{est:bootstrap enhanced vpm}.

A straight computation shows that $(\widetilde{u}, \widetilde{v})$ satisfies the following coupled system of semi-classical pseudo-differential equations:
\begin{equation} \label{wave-KG system in semi-classical coordinates}
\begin{cases}
& \big[D_t - \oph(x\cdot\xi - |\xi|)\big]\widetilde{u}(t,x) = h^{-1}\left[q_w(t,tx)+c_w(t,tx)+r^{NF}_w(t, tx) \right]\\
& \big[D_t - \oph(x\cdot\xi - \langle\xi\rangle)\big]\widetilde{v}(t,x) = h^{-1}r^{NF}_{kg}(t,tx),
\end{cases}
\end{equation}
where $\oph$ denotes the semi-classical Weyl quantization introduced in \ref{Def: Weyl and standard quantization} $(i)$. Moreover, if $\mathcal{M}_j$ (resp. $\mathcal{L}_j$), $j=1,2$, is the operator introduced in \eqref{def Mj} (resp. \eqref{def Lj}), $\mathcal{M}_j\widetilde{u}$ (resp. $\mathcal{L}_j\widetilde{v}$) can be expressed in term of $Z_ju^{NF}$ (resp. $Z_j v^{NF}$).  We have the following general result:

\begin{lem} \label{Lem: relation between Z and M/L}
$(i)$ Let $w(t,x)$ be a solution to the inhomogeneous half wave equation
\begin{equation} \label{half-wave-w}
\left[D_t + |D_x| \right] w(t,x) = f(t,x),
\end{equation}
and $\widetilde{w}(t,x)=tw(t,tx)$. For any $j=1,2$,
\begin{equation} \label{relation_Zjw_Mjwidetilde(w)}
Z_jw(t,y) =  ih \left[-\mathcal{M}_j\widetilde{w}(t,x) + \frac{1}{2i}\oph\left(\frac{\xi_j}{|\xi|}\right)\widetilde{w}(t,x) \right]|_{x=\frac{y}{t}} + iy_jf(t,y);
\end{equation}
$(ii)$ If $w(t,x)$ is solution to the inhomogeneous half Klein-Gordon equation
\begin{equation}\label{half KG}
\left[D_t + \langle D_x\rangle \right] w(t,x) = f(t,x),
\end{equation}
then
\begin{equation} \label{relation_Zjw_Ljwidetilde(w)}
Z_j  w(t,y) = ih\left[-\oph(\langle\xi\rangle)\mathcal{L}_j\widetilde{w}(t,x)+\frac{1}{i}\oph\Big(\frac{\xi_j}{\langle\xi\rangle}\Big)\widetilde{w}(t,x)\right]|_{x=\frac{y}{t}}+ iy_jf(t,y).
\end{equation}
\proof
$(i)$ If $w$ is solution to half wave equation \eqref{half-wave-w} then $\widetilde{w}(t,x)$ satisfies
\begin{equation*}
\big[D_t - \oph(x\cdot\xi - |\xi|)\big]\widetilde{w}(t,x) =  h^{-1}f(t, tx),
\end{equation*}
so, for any $i=1,2$,
\begin{equation*} 
\begin{split}
&Z_j  w(t,y) =\\
& i h^{-1}\left[x_j D_t + \oph(\xi_j - x_j x\cdot\xi) + \frac{3h}{2i}x_j\right]\left(\frac{1}{t}\widetilde{w}(t,x)\right)\Big|_{x=\frac{y}{t}} \\
& = i \left[ x_j D_t + \oph(\xi_j - x_jx\cdot\xi) + \frac{h}{2i}x_j\right]\widetilde{w}(t,x)\Big|_{x=\frac{y}{t}} 
\\
& = i \left[x_j \oph(x\cdot\xi - |\xi|)\widetilde{w}(t,x) + \oph(\xi_j - x_j x\cdot \xi)\widetilde{w}(t,x) + \frac{h}{2i}x_j\widetilde{u}(t,x) + h^{-1}x_j f(t,tx)\right]\big|_{x=\frac{y}{t}}   \\
& = ih \left[-\mathcal{M}_j\widetilde{w}(t,x) + \frac{1}{2i}\oph\left(\frac{\xi_j}{|\xi|}\right)\widetilde{w}(t,x)\right]|_{x=\frac{y}{t}}+ iy_j f(t,y).
\end{split}
\end{equation*}
We should specify that last equality is obtained by a trivial version of symbolic calculus \eqref{a sharp b asymptotic formula}, that applies also to symbols $b(\xi)$ singular at $\xi=0$.
Indeed, if symbol $a=a(x,\xi)$ is linear in $x$, and $b(\xi)$ is lipschitz, the development $a\sharp b$ is actually finite:
\begin{equation*}
a\sharp b(x,\xi) = a(x,\xi)b(\xi) - \frac{h}{2i}\partial_x a(x,\xi)\cdot\partial_\xi b(\xi).
\end{equation*}
$(ii)$ The result is analogous to the previous one, after observing that $\widetilde{w}$ satisfies
\begin{equation*}
\big[D_t - \oph(x\cdot\xi - \langle \xi\rangle)\big]\widetilde{w}(t,x) =  h^{-1}f(t, tx).
\end{equation*}
\endproof
\end{lem}

As a straight consequence of lemma \ref{Lem: relation between Z and M/L} and system \eqref{wave-KG system in semi-classical coordinates} we have that
\begin{subequations}
\begin{equation} \label{relation between Zju and Mj utilde-new}
Z_j u^{NF}(t,y) =  ih \left[-\mathcal{M}_j\widetilde{u}(t,x) + \frac{1}{2i}\oph\left(\frac{\xi_j}{|\xi|}\right)\widetilde{u}(t,x)\right]|_{x=\frac{y}{t}} + iy_j\left[q_w+c_w+r^{NF}_w\right](t,y),
\end{equation}
\begin{equation} \label{relation between Zjv and Lj vtilde}
Z_jv^{NF}(t,y) =  ih \left[ - \oph(\langle\xi\rangle)\mathcal{L}_j\widetilde{v}(t,x) + \frac{1}{i}\oph\Big(\frac{\xi_j}{\langle\xi\rangle}\Big)\widetilde{v}(t,x) \right]\Big|_{x=\frac{y}{t}} + iy_jr^{NF}_{kg}(t,y).
\end{equation}
\end{subequations}
In view of lemma \ref{Lemma : estimate of e(x,xi)}, it is also useful to write down the analogous relation between $(Z_m Z_nu)_{-}$ and $\mathcal{M}[t(Z_nu)_{-}(t,tx)]$.
As $(Z_nu)_{-}$ is solution to
\begin{equation*} 
\big(D_t + |D_x|\big)(Z_nu)_{-} = Z_n\Nlw(t,x),
\end{equation*}
from equality \eqref{relation_Zjw_Mjwidetilde(w)} with $w=(Z_nu)_{-}$ and the commutation between $Z_m$ and $D_t-|D_x|$ (see \eqref{commutator_Z_Dt-|D|}) we find that
\begin{multline} \label{relation ZmZnu Mutilde Zn-new}
(Z_mZ_nu)_{-}(t,y) = ih\Big[-\mathcal{M}_m\widetilde{u}^J(t,x) + \frac{1}{2i}\oph\Big(\frac{\xi_m}{|\xi|}\Big)\widetilde{u}^J(t,x)\Big]\big|_{x=\frac{y}{t}} + iy_mZ_n \Nlw(t,y) \\
-\frac{D_m}{|D_y|}(Z_nu)_{-}(t,y),
\end{multline}
where $J$ is the index such that $\Gamma^J=Z_n$ and $\widetilde{u}^J(t,x):= t(Z_nu)_{-}(t,tx)$.
Also, observe that from \eqref{Gamma_nonlinearity}, \eqref{def_G1}, \eqref{def u+- v+-} and \eqref{def uIpm vIpm}
\begin{equation*}
Z_n\Nlw = Q^{\mathrm{w}}_0\big((Z_nv)_\pm, D_1v_\pm\big) + Q^{\mathrm{w}}_0\big(v_\pm, D_1(Z_nv)_\pm\big) - \delta_n^1 Q^{\mathrm{w}}_0(v_\pm, D_tv_\pm)
\end{equation*}
with $\delta_n^1=1$ for $n=1$, and that from inequality \eqref{Hs norm DtV} with $s=0$,
\begin{multline} \label{L2 est NLwZn}
\|Z_n \Nlw(t,\cdot)\|_{L^2}\lesssim \|Z_nV(t,\cdot)\|_{H^1}\|V(t,\cdot)\|_{H^{2,\infty}}+ \big[\|V(t,\cdot)\|_{H^1} \\
+  \|V(t,\cdot)\|_{L^2}\left(\|U(t,\cdot)\|_{H^{1,\infty}}+ \|\mathrm{R}_1U(t,\cdot)\|_{H^{1,\infty}}\right) + \|V(t,\cdot)\|_{L^\infty}\|U(t,\cdot)\|_{H^1}\big]\|V(t,\cdot)\|_{H^{1,\infty}}.
\end{multline}

Moreover, from the definition of $\Mcal_j$ and $\Lcal_j$ we see that
\begin{gather*}
h\mathcal{M}_j\widetilde{w}(t,x) = \left[y_j|D_y| - tD_j + \frac{1}{2i}\frac{D_j}{|D_y|}\right]w(t,y)|_{y=tx},\\
h\oph(\langle \xi\rangle)\mathcal{L}_j \widetilde{w}(t,x) = \left[y_j\langle D_y\rangle - tD_j -i \frac{D_j}{\langle D_y\rangle}\right]w(t,y)|_{y=tx},
\end{gather*}
so lemma \ref{Lem: relation between Z and M/L} implies that, if $w$ is solution to half wave equation \eqref{half-wave-w} (resp. to half Klein-Gordon \eqref{half KG}),
\begin{subequations} \label{relation_w_Zjw}
\begin{gather}
\left[y_j|D_y| - tD_j + \frac{1}{2i}\frac{D_j}{|D_y|}\right]w(t,y)= iZ_jw(t,y)+ \frac{1}{2i}\frac{D_j}{|D_y|}w(t,y)+ y_jf(t,y),\label{relation_w_wave_Zjw-new} \\
\left(\text{resp. }\left[\langle D_y\rangle y_j - tD_j \right]w(t,y) =iZ_jw(t,y)-i \frac{D_j}{\langle D_y\rangle}w(t,y)+y_jf(t,y)\right). \label{relation_w_KG_Zjw}
\end{gather}
\end{subequations}

\subsection{Derivation of the ODE and propagation of the uniform estimate on the Klein-Gordon component} \label{Subsection : The Derivation of the ODE Equation} 

Let us firstly deal with the semi-classical Klein-Gordon equation satisfied by $\widetilde{v}$: 
\begin{equation} \label{semi-classical KG equation}
\big[D_t - \oph(x\cdot\xi - p(\xi))\big]\widetilde{v}(t,x) = h^{-1}r^{NF}_{kg}(t,tx),
\end{equation}
where $p(\xi) = \langle\xi\rangle$ and $r^{NF}_{kg}$ is given by \eqref{def rNF-kg} and satisfies \eqref{est L2 Linfty (cut-off) rNF-kg-new}.
We remind that $p'(\xi)$ denotes the gradient of $p(\xi)$ while $p''(\xi)$ is its $2\times 2$ Hessian matrix, and that $\Lcal_j$ is the operator introduced in \eqref{def Lj} for $j=1,2$.
We also remind definition \eqref{def_Lkg} of manifold $\Lkg$, represented in dimension 1 by picture \ref{picture: Lkg} below, and decompose $\vt$ into the sum of two contributions: one localized in a neighbourhood of $\Lkg$ of size $\sqrt{h}$ (in the spirit of \cite{ifrim_tataru:global_bounds}), the other localized out of this neighbourhood. 
\begin{figure}[h]
\begin{center}
\begin{tikzpicture}[scale=1.9]

\draw[->] (0,-1.7) -- (0,1.7);
\draw[->] (-2,0) -- (2,0);
\node[below] at (1.9,0) {\small $x$};
\node[left] at (0,1.6) {\small $\xi$};

\draw[dashed] (-1,-1.7) -- (-1,1.7);
\draw[dashed] (1,-1.7) -- (1,1.7);
\node[below left] at (-1,0) {$-1$};
\node[below right] at (1,0) {$1$};

\draw[blue,thick] [domain=0:0.87] plot(\x, {\x/((1-\x^2)^(1/2))});
\draw[blue,thick] [domain=0:0.87] plot(-\x, {-\x/((1-\x^2)^(1/2))});
\node[right] at (-0.8,-1.6) {$\Lkg$};
\end{tikzpicture}
\end{center}
\caption{Lagrangian for the Klein-Gordon equation}
\label{picture: Lkg}
\end{figure}
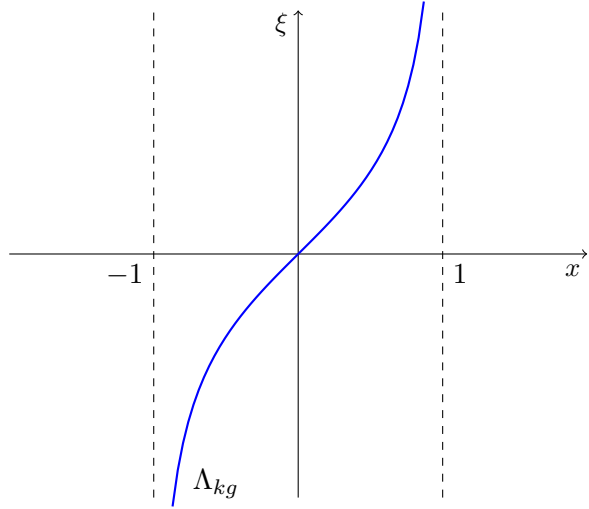
The contribution localized away from $\Lkg$ appears to be a $O(h^{1/2-0})$ if we assume a moderate growth for the $L^2$ norm of $\mathcal{L}^\mu\vt$, with $0\le |\mu|\le 2$, and has hence a better decay in time than the one expected for $\vt$ (remind $h=t^{-1}$).
Thus the main contribution to $\vt$ is the one localized around $\Lkg$. We are going to show that this latter one is solution to an ODE (see proposition \ref{Prop: ODE for vtilde Sigma_Lambda}) and that its $H^{\rho,\infty}_h$ norm is uniformly bounded in time, which will finally enable us to propagate \eqref{est:a-priori_vt} and obtain \eqref{est: boostrap vpm} (see proposition \ref{Prop:propagation_unif_est_V}).

For any fixed $\rho\in\mathbb{Z}$ let $\Sigma(\xi):=\langle\xi\rangle^{\rho}$, and for some $\gamma, \chi\in C^\infty_0(\mathbb{R}^2)$ equal to 1 close to the origin, $\sigma>0$ small (e.g. $\sigma<\frac{1}{4}$) let \index{Gammakg@$\Gamma^{kg}$, operator}
\begin{equation} \label{def Gamma_kg}
\Gamma^{kg} := \oph\left(\gamma\left(\frac{x - p'(\xi)}{\sqrt{h}}\right)\chi(h^\sigma\xi)\right).
\end{equation}
We also introduce the following notations:
\begin{equation}
\widetilde{v}^\Sigma:= \oph(\Sigma(\xi))\widetilde{v},
\end{equation}
together with \index{vtildeL@$\widetilde{v}^\Sigma_{\Lambda_{kg}}$, Klein-Gordon component in semi-classical setting, localised around $\Lambda_{kg}$}\index{vtildeLc@$\widetilde{v}^\Sigma_{\Lambda^c_{kg}}$, Klein-Gordon component in semi-classical setting, localised away from $\Lambda_{kg}$}
\begin{subequations}\label{def_both_vtilde_Sigma_Lambda}
\begin{equation} \label{def vtilde_Sigma_Lambda}
\widetilde{v}^\Sigma_{\Lambda_{kg}} := \Gamma^{kg}\widetilde{v}^\Sigma,
\end{equation}
\begin{equation} \label{def vtilde_Lambda,c}
\widetilde{v}^\Sigma_{\Lambda^c_{kg}}:= \oph\left(1 -\gamma\left(\frac{x-p'(\xi)}{\sqrt{h}}\right)\chi(h^\sigma\xi)\right)\widetilde{v}^\Sigma,
\end{equation}
\end{subequations}
so that $\widetilde{v}^\Sigma = \widetilde{v}^\Sigma_{\Lambda_{kg}} + \widetilde{v}^\Sigma_{\Lambda^c_{kg}}$, and remind that $\|\mathcal{L}^\gamma w\| = \|\mathcal{L}^{\gamma_1}_1\mathcal{L}^{\gamma_2}_2w\|$, for any $\gamma=(\gamma_1,\gamma_2)\in\mathbb{N}^2$.

\begin{lem} \label{Lem:Op(gammatilde)v}
Let $\widetilde{\gamma}\in C^\infty(\mathbb{R}^2)$ vanish in a neighbourhood of the origin and be such that $|\partial^\alpha_z\widetilde{\gamma}(z)|\lesssim \langle z\rangle^{-|\alpha|}$. Let $c(x,\xi)\in S_{\delta,\sigma}(1)$ with $\delta\in [0,\frac{1}{2}]$, $\sigma>0$, be supported for $|\xi|\lesssim h^{-\sigma}$. For any $\chi\in C^\infty_0(\mathbb{R}^2)$ such that $\chi(h^\sigma\xi)\equiv 1$ on the support of $c(x,\xi)$, 
\begin{subequations} \label{est_1L}
\begin{gather}
\left\| \oph\Big(\widetilde{\gamma}\Big(\frac{x-p'(\xi)}{\sqrt{h}}\Big)c(x,\xi)\Big)w\right\|_{L^2} \lesssim \sum_{|\mu|=0}^1 h^{\frac{1}{2}-\beta}\left\|\oph(\chi(h^\sigma\xi))\mathcal{L}^\mu w \right\|_{L^2}, \label{est_1L-L2}\\
\left\| \oph\Big(\widetilde{\gamma}\Big(\frac{x-p'(\xi)}{\sqrt{h}}\Big)c(x,\xi)\Big)w\right\|_{L^\infty} \lesssim \sum_{|\mu|=0}^1 h^{-\beta}\left\|\oph(\chi(h^\sigma\xi))\mathcal{L}^\mu w \right\|_{L^2}, \label{est_1L-Linfty}
\end{gather}
\end{subequations}
and
\begin{subequations}\label{est_2L}
\begin{gather}
\left\| \oph\Big(\widetilde{\gamma}\Big(\frac{x-p'(\xi)}{\sqrt{h}}\Big)c(x,\xi)\Big)w\right\|_{L^2} \lesssim \sum_{|\mu|=0}^2 h^{1-\beta}\left\|\oph(\chi(h^\sigma\xi))\mathcal{L}^\mu w \right\|_{L^2}, \label{est_2L-L2}\\
\left\| \oph\Big(\widetilde{\gamma}\Big(\frac{x-p'(\xi)}{\sqrt{h}}\Big)c(x,\xi)\Big)w\right\|_{L^\infty} \lesssim \sum_{|\mu|=0}^2 h^{\frac{1}{2}-\beta}\left\|\oph(\chi(h^\sigma\xi))\mathcal{L}^\mu w \right\|_{L^2}, \label{est_2L-Linfty}
\end{gather}
\end{subequations}
for a small $\beta>0$, $\beta\rightarrow 0$ as $\sigma\rightarrow 0$.
\proof
The proof of \eqref{est_1L} (resp. of \eqref{est_2L}) follows straightly by inequalities \eqref{est L2 Linfty for GammaTilde} (resp. \eqref{est:L2 Linfty Op(gamma) LLv}), after observing that, as $\widetilde{\gamma}$ vanishes in a neighbourhood of the origin,
\begin{equation*}
\widetilde{\gamma}\Big(\frac{x-p'(\xi)}{\sqrt{h}}\Big)c(x,\xi) = \sum_{j=1}^2\widetilde{\gamma}_1^j\Big(\frac{x-p'(\xi)}{\sqrt{h}}\Big)\Big(\frac{x_j-p'_j(\xi)}{\sqrt{h}}\Big)c(x,\xi),
\end{equation*}
where $\widetilde{\gamma}_1^j(z):=\widetilde{\gamma}(z)z_j|z|^{-2}$ is such that $|\partial^\alpha_z\widetilde{\gamma}_1^j(z)|\lesssim \langle z\rangle^{-1-|\alpha|}$ (resp.
\begin{equation*}
\widetilde{\gamma}\Big(\frac{x-p'(\xi)}{\sqrt{h}}\Big)c(x,\xi) = \sum_{j=1}^2\widetilde{\gamma}_2\Big(\frac{x-p'(\xi)}{\sqrt{h}}\Big)\Big(\frac{x-p'(\xi)}{\sqrt{h}}\Big)^2c(x,\xi),
\end{equation*}
where $\widetilde{\gamma}_2(z):=\widetilde{\gamma}(z)|z|^{-2}$ is such that $|\partial^\alpha_z\widetilde{\gamma}(z)|\lesssim\langle z\rangle^{-2-|\alpha|}$).
\endproof
\end{lem}

\begin{cor} \label{Prop : est on vtilde_Lambda,c}
There exists $s>0$ sufficiently large such that
\begin{subequations} \label{est widetildev_Lambda,c}
\begin{equation} \label{est H1h widetildev_Lambda,c}
\left\|\widetilde{v}^\Sigma_{\Lambda^c_{kg}}(t,\cdot)\right\|_{L^2} \lesssim h^{1-\beta}\left(\|\widetilde{v}(t,\cdot)\|_{H^s_h} +\sum_{1\le |\mu|\le 2}\|\oph(\chi(h^\sigma\xi))\mathcal{L}^\mu\widetilde{v}(t,\cdot)\|_{L^2}\right),
\end{equation}
\begin{equation} \label{est Linfty of widetildev_Lambda,c}
\left\|\widetilde{v}^\Sigma_{\Lambda^c_{kg}}(t,\cdot)\right\|_{L^\infty} \lesssim h^{\frac{1}{2}-\beta}\left(\|\widetilde{v}(t,\cdot)\|_{H^s_h} + \sum_{1\le |\mu|\le 2}\|\oph(\chi(h^\sigma\xi))\mathcal{L}^\mu\widetilde{v}(t,\cdot)\|_{L^2}\right).
\end{equation}
\end{subequations}
for a small $\beta>0$, $\beta\rightarrow 0$ as $\sigma\rightarrow 0$.
\proof
Since symbol $1- \gamma\big(\frac{x-p'(\xi)}{\sqrt{h}}\big)\chi(h^\sigma\xi)$ is supported for $|\frac{x-p'(\xi)}{\sqrt{h}}|\ge d_1>0$ or $|h^\sigma \xi|\ge d_2>0$, for some small $d_1,d_2>0$, we may consider a smooth cut-off function $\widetilde{\chi}$ equal to 1 close to the origin and such that $\widetilde{\chi}\chi \equiv \widetilde{\chi}$, so that $1- \gamma\big(\frac{x-p'(\xi)}{\sqrt{h}}\big)\chi(h^\sigma\xi)$ writes as
\begin{equation*}
\left[1-\gamma\left(\frac{x-p'(\xi)}{\sqrt{h}}\right)\right]\widetilde{\chi}(h^\sigma\xi) + \left[1- \gamma\left(\frac{x-p'(\xi)}{\sqrt{h}}\right)\chi(h^\sigma\xi)\right](1-\widetilde{\chi})(h^\sigma\xi),
\end{equation*}
the first symbol being supported in $\{(x,\xi): |\frac{x-p'(\xi)}{\sqrt{h}}|\ge d_1, |\xi|\lesssim h^{-\sigma}\}$, the second one for large frequencies $|\xi|\gtrsim h^{-\sigma}$. 

Using lemma \ref{Lem : a sharp b} and the fact that $\gamma\Big(\frac{x-p'(\xi)}{\sqrt{h}}\Big)\chi(h^\sigma\xi)\in S_{\frac{1}{2},\sigma}\big(\big\langle \frac{x-p'(\xi)}{\sqrt{h}} \big\rangle^{-M}\big)$, for any $M\in\mathbb{N}$, we have that, for a fixed $N\in\mathbb{N}^*$,
\begin{align*}
\left[1-\gamma\Big(\frac{x-p'(\xi)}{\sqrt{h}}\Big)\chi(h^\sigma\xi)\right]\big(1-\widetilde{\chi}(h^\sigma\xi)\big)&= \big(1-\widetilde{\chi}(h^\sigma\xi)\big)\sharp \left[1-\gamma\Big(\frac{x-p'(\xi)}{\sqrt{h}}\Big)\chi(h^\sigma\xi)\right]\\
& + \sum_{1\le j<N}\widetilde{\chi}_j(h^\sigma\xi)\sharp a_j(x,\xi) + r_N(x,\xi),
\end{align*}
where function $\widetilde{\chi}_j(h^\sigma\xi)$ is still supported for large frequencies $|\xi|\gtrsim h^{-\sigma}$, for every $1\le j<N$, up to negligible multiplicative constants,
$$a_j(x,\xi) = h^{j(\frac{1}{2}+\sigma)}\sum_{|\alpha|=j}(\partial^\alpha\gamma)\Big(\frac{x-p'(\xi)}{\sqrt{h}}\Big)\chi(h^\sigma\xi) \in h^{j(\frac{1}{2}+\sigma)}S_{\frac{1}{2},\sigma}\big(\big\langle \frac{x-p'(\xi)}{\sqrt{h}} \big\rangle^{-M}\big),$$
and $r_N\in h^{N(\frac{1}{2}+\sigma)}S_{\frac{1}{2},\sigma}\big(\big\langle \frac{x-p'(\xi)}{\sqrt{h}} \big\rangle^{-M}\big)$. 
Lemma \ref{Lem : new estimate 1-chi}, proposition \ref{Prop : Continuity on H^s}, and the semi-classical Sobolev injection imply that
\begin{gather*}
\left\|\oph\Big( \Big[1- \gamma\Big(\frac{x-p'(\xi)}{\sqrt{h}}\Big)\chi(h^\sigma\xi)\Big](1-\widetilde{\chi})(h^\sigma\xi)\Big)\widetilde{v}^\Sigma(t,\cdot)\right\|_{L^2} \lesssim h^{N(s)}\|\widetilde{v}(t,\cdot)\|_{H^s_h}, \\
\left\|\oph\Big( \Big[1- \gamma\Big(\frac{x-p'(\xi)}{\sqrt{h}}\Big)\chi(h^\sigma\xi)\Big](1-\widetilde{\chi})(h^\sigma\xi)\Big)\widetilde{v}^\Sigma (t,\cdot)\right\|_{L^\infty} \lesssim h^{N'(s)} \|\widetilde{v}(t,\cdot)\|_{H^s_h},
\end{gather*}
where $N(s), N'(s)\ge 1$ if $s>2$ is sufficiently large. 

On the other hand, as function $(1-\gamma)\big(\frac{x-p'(\xi)}{\sqrt{h}}\big)$ vanishes in a neighbourhood of the origin and is such that $|\partial^\alpha_z (1-\gamma)(z)|\lesssim \langle z\rangle^{-|\alpha|}$, by inequalities \eqref{est_2L} and the fact that, using symbolic calculus to commute $\mathcal{L}$ with $\Sigma(\xi)$,
\begin{equation} \label{Lgamma vSigma}
 \|\oph(\chi(h^\sigma\xi))\mathcal{L}^\mu \widetilde{v}^\Sigma(t,\cdot)\|_{L^2}\lesssim h^{-\nu}\sum_{|\mu_1|\le |\mu|}  \|\oph(\chi(h^\sigma\xi))\mathcal{L}^{\mu_1} \widetilde{v}(t,\cdot)\|_{L^2}
\end{equation}
with $\nu = \rho\sigma$ if $\rho\ge 0$, 0 otherwise, we have that
\begin{gather*}
\left\|\oph\Big((1-\gamma)\Big(\frac{x-p'(\xi)}{\sqrt{h}}\Big)\chi(h^\sigma\xi)\Big)\widetilde{v}^\Sigma(t,\cdot)\right\|_{L^2}\lesssim \sum_{|\mu|\le 2}h^{1-\beta} \left\|\oph(\chi(h^\sigma\xi))\mathcal{L}^\mu \widetilde{v}(t,\cdot)\right\|_{L^2},\\
\left\|\oph\Big((1-\gamma)\Big(\frac{x-p'(\xi)}{\sqrt{h}}\Big)\chi(h^\sigma\xi)\Big)\widetilde{v}^\Sigma(t,\cdot)\right\|_{L^\infty}\lesssim \sum_{|\mu|\le 2}h^{\frac{1}{2}-\beta} \left\|\oph(\chi(h^\sigma\xi))\mathcal{L}^\mu \widetilde{v}(t,\cdot)\right\|_{L^2},
\end{gather*}
for a small $\beta>0$, $\beta\rightarrow 0$ as $\sigma\rightarrow 0$.
\endproof
\end{cor}

In the following lemma we show how to develop the symbol $a(x,\xi)$ associated to an operator acting on $\Gamma^{kg}w$, for some suitable function $w$, at $\xi = -d\phi(x)$, where $\phi(x)=\sqrt{1-|x|^2}$.

\begin{lem} \label{Lem:dev of a symbol at xi = -dvarphi(x)}
Let $a(x,\xi)$ be a real symbol in $S_{\delta,0}(\langle\xi\rangle^q)$, $q \in\mathbb{R}$, for some $\delta>0$ small, $\Sigma(\xi)=\langle\xi\rangle^\rho$ for some fixed $\rho\in\Z$, $\Gamma^{kg}$ the operator introduced in \eqref{def Gamma_kg} and $w=w(t,x)$ such that $\Lcal^\mu w(t,\cdot)\in L^2(\R^2)$ for any $|\mu|\le 2$. Let us also introduce $w^\Sigma_{\Lkg}:=\Gamma^{kg}\oph(\Sigma)w$.
There exists a family $(\theta_h(x))_h$ of $C^{\infty}_0$ functions real valued, equal to 1 on the closed ball $\overline{B_{1-ch^{2\sigma}}(0)}$ and supported in $\overline{B_{1-c_1h^{2\sigma}}(0)}$, for some small $0<c_1<c,\sigma>0$, with $\|\partial^\alpha_x\theta_h\|_{L^\infty}=O(h^{-2|\alpha|\sigma})$ and $(h\partial_h)^k\theta_h$ bounded for every $k$, such that
\begin{equation} \label{Op(a)v development}
\oph(a)w^\Sigma_{\Lkg}= \theta_h(x)a(x,-d\phi(x))w^\Sigma_{\Lkg}+ R_1(w) \,,
\end{equation}
where $R_1(w)$ satisfies
\begin{subequations} \label{est L2 Linfty R1(widetildev)}
\begin{equation} \label{est L2 of R1(widetildev)}
\|R_1(w)(t,\cdot)\|_{L^2}\lesssim h^{1-\beta}\left(\|w(t,\cdot)\|_{H^s_h} + \sum_{|\gamma|=1}\|\oph(\chi(h^\sigma\xi))\mathcal{L}^\gamma w(t,\cdot)\|_{L^2}\right),
\end{equation}
\begin{equation} \label{est Linfty of R1(widetildev)}
\|R_1(w)(t,\cdot)\|_{L^\infty}\lesssim h^{\frac{1}{2}-\beta}\left(\|w(t,\cdot)\|_{H^s_h} + \sum_{|\gamma|=1}\|\oph(\chi(h^\sigma\xi)) \mathcal{L}^\gamma w(t,\cdot)\|_{L^2}\right),
\end{equation}
\end{subequations}
with $\beta=\beta(\sigma,\delta)>0$, $\beta\rightarrow 0$ as $\sigma, \delta \rightarrow 0$.
Moreover, if $\partial_\xi a(x,\xi)$ vanishes at $\xi=-d\phi(x)$, the above estimates can be improved and $R_1(w)$ is rather a remainder $R_2(w)$ such that
\begin{subequations} \label{est L2 Linfty of R2(widetildev)}
\begin{equation} \label{est L2 of R2(widetildev)}
\|R_2(w)(t,\cdot)\|_{L^2}\lesssim h^{2-\beta}\left(\|w(t,\cdot)\|_{H^s_h} +\sum_{1\le |\gamma|\le 2 } \|\oph(\chi(h^\sigma\xi))\mathcal{L}^\gamma w(t,\cdot)\|_{L^2}\right),
\end{equation}
\begin{equation}  \label{est Linfty of R2(widetildev)}
\|R_2(w)(t,\cdot)\|_{L^\infty}\lesssim h^{\frac{3}{2}-\beta}\left(\|w(t,\cdot)\|_{H^s_h} + \sum_{1\le |\gamma|\le 2}\|\oph(\chi(h^\sigma\xi))\mathcal{L}^\gamma w(t,\cdot)\|_{L^2}\right).
\end{equation}
\end{subequations}
\proof
After lemma \ref{Lem:family_thetah} we know that there exists a family of functions $\theta_h(x)$ as in the statement such that equality \eqref{cut-off-thetah} holds.
We highlight the fact that any derivative of $\theta_h$ vanishes on the support of $\gamma\big(\frac{x-p'(\xi)}{\sqrt{h}}\big)\chi(h^\sigma\xi)$ and its derivatives.
After remark \ref{Remark:symbols_with_null_support_intersection}, this implies that 
\[w^\Sigma_{\Lkg}= \theta_h(x)w^\Sigma_{\Lkg}+ r_\infty, \quad r_\infty\in h^NS_{\frac{1}{2},\sigma}(\langle x \rangle^{-\infty})\]
and hence that 
\[\oph(a)w^\Sigma_{\Lkg}= \oph(a)\theta_h(x)w^\Sigma_{\Lkg}+ \oph(  r^a_\infty)w^\Sigma_{\Lkg},\]
with $r^a_\infty=a\sharp r_\infty \in h^{N-\gamma}S_{\frac{1}{2},\sigma}(\langle x\rangle^{-\infty})$ and $\gamma=q\sigma$ if $q\ge 0$, 0 otherwise.
From proposition \ref{Prop : Continuity on H^s} and the semi-classical Sobolev injection it follows at once that $\oph( r^a_\infty)w^\Sigma_{\Lkg}$ satisfies enhanced estimates \eqref{est L2 Linfty of R2(widetildev)} if $N$ is taken sufficiently large.
Up to negligible multiplicative constants, a further application of symbolic calculus gives also that 
\begin{multline*}
\oph(a(x,\xi))\theta_h(x)w^\Sigma_{\Lkg}= \oph(a(x,\xi)\theta_h(x))w^\Sigma_{\Lkg}+\sum_{|\alpha|=1}^{N-1}h^{|\alpha|} \oph\big(\partial^\alpha_\xi a(x,\xi)\partial^\alpha_x\theta_h(x)\big)w^\Sigma_{\Lkg} \\
+ \oph(r_N(x,\xi))w^\Sigma_{\Lkg},
\end{multline*}
where $r_N\in h^{N-\beta} S_{\delta',0}(\langle\xi\rangle^{q-N}\langle x \rangle^{-\infty})$ for a new small $\beta=\beta(\delta,\sigma)$ and $\delta'=\max\{\delta,\sigma\}$.
From the same argument as above $\oph(r_N)w^\Sigma_{\Lkg}$ verifies enhanced estimates \eqref{est L2 Linfty of R2(widetildev)} if $N$ is suitably chosen. 
Also, since the support of $\partial^\alpha_\xi a(x,\xi)\cdot\partial^\alpha_x\theta_h(x)$ has empty intersection with that of $\gamma\Big(\frac{x-p'(\xi)}{\sqrt{h}}\Big)\chi(h^\sigma\xi)$ for any $|\alpha|\ge 1$, 
all the $|\alpha|$-order terms in the above equality are remainders $R_1(w)$.

Now, as symbol $a(x,\xi)\theta_h(x)$ is supported for $|x|\le 1-c_1h^{2\sigma}<1$, we are allowed to develop it at $\xi = -d\phi(x)$:
\begin{align*} \label{dev of a at xo = -dvarphi}
a(x,\xi)\theta_h(x) &= a(x,-d\phi(x))\theta_h(x) + \sum_{|\alpha| =1}\int_0^1 (\partial^\alpha_\xi a)(x, t\xi + (1-t)d\phi(x))dt\, \theta_h(x) (\xi + d\phi(x))^\alpha \\
& = a(x,-d\phi(x))\theta_h(x) + \sum_{j=1}^2 b_j(x,\xi) (x_j - p'_j(\xi)), \numberthis
\end{align*}
with
\begin{equation} \label{def bj in development}
b_j(x,\xi) = \sum_{|\alpha| =1}\int_0^1 (\partial^\alpha_\xi a)(x, t\xi + (1-t)d\phi(x)) dt\, \theta_h(x)\frac{(\xi + d\phi(x))^\alpha (x_j -p'_j(\xi))}{|x-p'(\xi)|^2}, \quad j=1,2.
\end{equation}
If $\chi_1\in C^\infty_0(\mathbb{R}^2)$ is a new cut-off function equal to 1 close to the origin, we can reduce ourselves to the analysis of symbol $b_j(x,\xi)(x_j - p'_j(\xi))\chi_1(h^\sigma \xi)$.
In fact, as $b_j(x,\xi)(x_j - p'_j(\xi))(1-\chi_1)(h^\sigma\xi)$ is supported for large frequencies, one can prove that its operator acting on $w^\Sigma_{\Lkg}$ is a $O_{L^2\cap L^\infty}(h^N\|w(t,\cdot)\|_{H^s_h})$ with $N>0$ large as long as $s>0$ is large, by using the semi-classical Sobolev injection, symbolic calculus of proposition \ref{Prop: a sharp b}, lemma \ref{Lem : new estimate 1-chi} and proposition \ref{Prop : Continuity on H^s}.
Furthermore, if we consider a smooth cut-off function $\widetilde{\gamma}\in C^\infty_0(\mathbb{R}^2)$, equal to 1 close to the origin and such that $\widetilde{\gamma}\big(\langle\xi\rangle^2 (x-p'(\xi))\big)\equiv 1$ on the support of $\gamma\big(\frac{x-p'(\xi)}{\sqrt{h}}\big)\chi(h^\sigma\xi)$ (which is possible if $\sigma<1/4$), we have that
\begin{align*}
b_j(x,\xi)(x_j - p'_j(\xi))\chi_1(h^\sigma \xi)  &= b_j(x,\xi)(x_j - p'_j(\xi))\chi_1(h^\sigma \xi) \widetilde{\gamma}\big(\langle\xi\rangle^2 (x-p'(\xi))\big) \\
&+ b_j(x,\xi)(x_j - p'_j(\xi))\chi_1(h^\sigma \xi)(1-\widetilde{\gamma})\big(\langle\xi\rangle^2 (x-p'(\xi))\big).
\end{align*}
Since $b_j(x,\xi)(x_j - p'_j(\xi))\chi_1(h^\sigma \xi)(1-\widetilde{\gamma})\big(\langle\xi\rangle^2 (x-p'(\xi))\big)\in h^{-\beta} S_{\delta,\sigma}(1)$, for some new small $\beta, \delta>0$, and its support has empty intersection with that of $\gamma\big(\frac{x-p'(\xi)}{\sqrt{h}}\big)$ (which instead belongs to class $S_{\frac{1}{2},0}(\langle \frac{x-p'(\xi)}{\sqrt{h}}\rangle^{-M})$, for $M\in\mathbb{N}$ as large as we want), its quantization acting on $w^\Sigma_{\Lkg}$ is also an enhanced remainder $R_2(w)$.

The very contribution that only enjoys estimates \eqref{est L2 Linfty R1(widetildev)} is 
$\oph\big(c(x,\xi)(x_j - p'_j(\xi))\big)w^\Sigma_{\Lkg}$, with $c(x,\xi):=b_j(x,\xi)\chi_1(h^\sigma \xi) \widetilde{\gamma}\big(\langle\xi\rangle^2 (x-p'(\xi))\big)\in h^{-\beta}S_{2\sigma,\sigma}(1)$ and $\beta$ depending linearly on $\sigma$.
In fact, if we assume that the support of $\chi_1$ is sufficiently small so that $\chi_1\chi \equiv \chi_1$ and all derivatives of $\chi$ vanish on that support, by using symbolic development \eqref{a sharp b asymptotic formula} until a sufficiently large order $N$ and observing that 
\begin{multline*}
\left\{ c(x,\xi)(x_j-p'_j(\xi)), \gamma\Big(\frac{x-p'(\xi)}{\sqrt{h}}\Big) \right\}= \left\{c(x,\xi), \gamma\Big(\frac{x-p'(\xi)}{\sqrt{h}}\Big) \right\} (x_j-p'_j(\xi)) \\
=\left[(\partial_\xi c)\cdot (\partial \gamma)\Big(\frac{x-p'(\xi)}{\sqrt{h}}\Big)+(\partial_x c)\cdot (\partial \gamma)\Big(\frac{x-p'(\xi)}{\sqrt{h}}\Big) p''(\xi)\right]\Big(\frac{x_j-p'_j(\xi)}{\sqrt{h}}\Big)
\end{multline*}
does not lose any power $h^{-1/2}$, we derive that, up to negligible constants,
\begin{multline*}
\Big[c(x,\xi)(x_j - p'_j(\xi))\Big] \sharp \Big[\gamma\Big(\frac{x-p'(\xi)}{\sqrt{h}}\Big)\chi(h^\sigma\xi)\Big] = \gamma\Big(\frac{x-p'(\xi)}{\sqrt{h}}\Big)\chi(h^\sigma\xi)c(x,\xi)(x_j - p'_j(\xi))\\
 + {\sum}'h\widetilde{\gamma}\Big(\frac{x-p'(\xi)}{\sqrt{h}}\Big)\widetilde{c}(x,\xi)+ r_N(x,\xi).
\end{multline*}
In the above equality ${\sum}'$ is a concise notation to indicate a linear combination, $\widetilde{\gamma}\in C^\infty_0(\mathbb{R}^2\setminus\{0\})$, $\widetilde{c}\in h^{-\beta}S_{\delta,\sigma}(1)$ for some new small $\beta,\delta>0$, and $r_N\in h^{(N+1)/2-\beta}S_{\frac{1}{2},\sigma}\big(\big\langle\frac{x-p'(\xi)}{\sqrt{h}}\big\rangle^{-(M-1)}\big)$ as $c(x,\xi)(x_j-p'_j(\xi))\in h^{1/2-\beta}S_{2\sigma,\sigma}\big(\langle \frac{x-p'(\xi)}{\sqrt{h}}\rangle\big)$.
From inequalities \eqref{est L2 Linfty for GammaTilde} and \eqref{Lgamma vSigma} we deduce that $\oph\big(\gamma\big(\frac{x-p'(\xi)}{\sqrt{h}}\big)\chi(h^\sigma\xi)c(x,\xi)(x_j - p'_j(\xi))\big)\oph(\Sigma)w$ is a remainder $R_1(w)$ satisfying \eqref{est L2 Linfty R1(widetildev)}. 
The quantization of all the addends in ${\sum}'$ acting on $\oph(\Sigma)w$
is estimated by using that $\widetilde{\gamma}(z)$ vanishes in a neighbourhood of the origin and can be rewritten as $\sum_{j=1,2}\widetilde{\gamma}_2(z)z^2_j$, with $\widetilde{\gamma}_2(z):=\widetilde{\gamma}(z)|z|^{-2}$ such that $|\partial^\alpha_z \widetilde{\gamma}_2(z)|\lesssim \langle z\rangle^{-2-|\alpha|}$. Inequalities \eqref{est:L2 Linfty Op(gamma) LLv} and the successive commutation of $\mathcal{L}^\gamma$ with $\Sigma$, for $|\gamma|=1,2$, give then that $h\oph\big(\widetilde{\gamma}\big(\frac{x-p'(\xi)}{\sqrt{h}}\big)\widetilde{c}(x,\xi)\big)\oph(\Sigma)w$ is a remainder $R_2(w)$. Finally, as 
\[r_N(x,\xi)\sharp \Sigma(\xi)\in h^{\frac{N}{2}-\beta-\mu}S_{\frac{1}{2},\sigma}(\langle\frac{x-p'(\xi)}{\sqrt{h}}\rangle^{-(M-1)})\]
with $\mu=\sigma\rho$ if $\rho\ge 0$, 0 otherwise, $\oph(r_N)\oph(\Sigma)w$ is also a remainder $R_2(w)$ just from \ref{Prop : Continuity on H^s}, \ref{Prop : Continuity from $L^2$ to L^inf}, fixing $N\in\mathbb{N}$ sufficiently large (e.g. $N=3$).

If symbol $a(x,\xi)$ is such that $\partial_\xi a|_{\xi = -d\phi}=0$, instead of equality \eqref{dev of a at xo = -dvarphi} with $b_j$ given by \eqref{def bj in development}, we have
\begin{equation*}
a(x,\xi)\theta_h(x) = a(x, -d\phi(x)) \theta_h(x)+ \sum_{j=1,2} b(x,\xi)(x_j - p'_j(\xi))^2,
\end{equation*}
with
\begin{equation*}
b(x,\xi) = \sum_{|\alpha|=2}\frac{2}{\alpha!}\int_0^1 (\partial^\alpha_\xi a)(t\xi - (1-t)d\phi(x))(1-t) dt\, \theta_h(x)\frac{(\xi +d\phi(x))^\alpha}{|x-p'(\xi)|^2}.
\end{equation*}
The same argument as before can be applied to $\oph\big(b(x,\xi)\theta_h(x)(x_j - p'_j(\xi))^2\big)w^\Sigma_{\Lkg}$ to show that it reduces to
\begin{equation*}
\oph\Big(b(x,\xi)\theta_h(x)(x_j - p'_j(\xi))^2\chi_1(h^\sigma \xi) \widetilde{\gamma}\big(\langle\xi\rangle^2 (x-p'(\xi))\big)\Big)w^\Sigma_{\Lkg} + R_2(w),
\end{equation*}
with $R_2(w)$ satisfying \eqref{est L2 Linfty of R2(widetildev)}. 
If 
\begin{equation*}
B(x,\xi) := b(x,\xi)\theta_h(x)\chi_1(h^\sigma \xi) \widetilde{\gamma}\big(\langle\xi\rangle^2 (x-p'(\xi))\big)
\end{equation*}
then $B(x,\xi)(x_j-p'_j(\xi))^2 \in h^{-\beta}S_{\delta',\sigma}(1)$ by lemma \ref{Lem : on e and etilde}, for some new small $\beta, \delta'>0$ depending on $\sigma,\delta$. 
Using lemma \ref{Lem : a sharp b}, symbolic development \eqref{a sharp b asymptotic formula} until order 4, and assuming that the support of $\chi_1$ is sufficiently small so that $\chi\chi_1\equiv \chi$,
we derive that
\begin{multline*}
\Big[B(x,\xi)(x_j - p'_j(\xi))^2\Big] \sharp \Big[\gamma\Big(\frac{x-p'(\xi)}{\sqrt{h}}\Big)\chi(h^\sigma\xi)\Big] = B(x,\xi) \gamma\Big(\frac{x-p'(\xi)}{\sqrt{h}}\Big)(x_j-p'_j(\xi))^2 \\
+ \frac{h}{i} \sum_{i=1}^2(\partial_i\gamma)\Big(\frac{x-p'(\xi)}{\sqrt{h}}\Big)\Big(\frac{x_j-p'_j(\xi)}{\sqrt{h}}\Big)\left[(\partial_{\xi_i}B) +\sum_j (\partial_{x_j}B) p''_{ij}(\xi)\right](x_j-p'_j(\xi))\\
+ {\sum_{2\le |\alpha|\le 3}}' h^{\frac{|\alpha|}{2}-2\delta'-\beta} \gamma_\alpha\Big(\frac{x-p'(\xi)}{\sqrt{h}}\Big)B_\alpha(x,\xi)+r_4(x,\xi),
\end{multline*}
where $\gamma_\alpha\in C^\infty_0(\mathbb{R}^2\setminus \{0\})$, $B_\alpha(x,\xi)\in S_{\delta',\sigma}(1)$, and $r_4(x,\xi)\in h^{2-4\delta'-\beta}S_{\frac{1}{2},\sigma}\big(\langle\frac{x-p'(\xi)}{\sqrt{h
}}\rangle^{-M}\big)$.
As $r_4(x,\xi)\sharp \Sigma(\xi) \in h^{2-\beta'}S_{\frac{1}{2},\sigma}\big(\langle\frac{x-p'(\xi)}{\sqrt{h
}}\rangle^{-M}\big)$, for $\beta'=2-4\delta'-\beta-\rho\sigma$ if $\rho\ge 0$, $\beta'=2-4\delta'-\beta$ otherwise, it immediately follows from propositions \ref{Prop : Continuity on H^s} and \ref{Prop : Continuity from $L^2$ to L^inf} that $\oph(r_4)\widetilde{v}^\Sigma$ is a remainder $R_2(w)$.
After inequalities \eqref{est:L2 Linfty Op(gamma) LLv} with $\gamma_n=\gamma$ and $c=B$ (resp. inequalities \eqref{est L2 Linfty for GammaTilde} with $\gamma_n(z) = \partial_i\gamma(z) z_j$ and $c = h^{\delta'}[(\partial_{\xi_i}B) + (\partial_x B)\cdot (\partial_\xi p'_1 + \partial_\xi p'_2)]\in S_{\delta',\sigma}(1)$, for $i,j=1,2$), and \eqref{Lgamma vSigma}, we deduce that the quantization of the first (resp. the second) contribution in above symbolic development is a remainder $R_2(w)$, when acting on $\oph(\Sigma)w$.
Finally, as $\gamma_\alpha$ vanishes in a neighbourhood of the origin, 
we write
\begin{gather*}
\gamma_\alpha\Big(\frac{x-p'(\xi)}{\sqrt{h}}\Big) = \sum_{k=1}^2 h^{-1}\underbrace{ \gamma_\alpha\Big(\frac{x-p'(\xi)}{\sqrt{h}}\Big)\Big|\frac{x-p'(\xi)}{\sqrt{h}}\Big|^{-2}}_{\widetilde{\gamma}_\alpha\big(\frac{x-p'(\xi)}{\sqrt{h}}\big)}\times (x_k-p'_k(\xi))^2, \quad |\alpha|=2,\\
\gamma_\alpha\Big(\frac{x-p'(\xi)}{\sqrt{h}}\Big) = \sum_{k=1}^2 h^{-\frac{1}{2}} \underbrace{\gamma_\alpha\Big(\frac{x-p'(\xi)}{\sqrt{h}}\Big)\Big(\frac{x_k-p'_k(\xi)}{\sqrt{h}}\Big)\Big|\frac{x-p'(\xi)}{\sqrt{h}}\Big|^{-2}}_{\widetilde{\gamma}^k_\alpha\big(\frac{x-p'(\xi)}{\sqrt{h}}\big)} \times (x_k-p'_k(\xi)), \quad |\alpha| =3
\end{gather*}
and obtain that the quantization of $\alpha$-th order term with $|\alpha|=2$ (resp. $|\alpha|=3$) is a remainder $R_2(w)$ when acting on $\oph(\Sigma)w$, after inequalities \eqref{est:L2 Linfty Op(gamma) LLv} (resp. \eqref{est L2 Linfty for GammaTilde}) with $\gamma_n = \widetilde{\gamma}_\alpha$ (resp. $\gamma_n=\widetilde{\gamma}^k_\alpha$, $k=1,2$) and $c=B_\alpha$.
\endproof
\end{lem}

The following two results allow us to finally derive the ODE satisfied by $\widetilde{v}^\Sigma_{\Lambda_{kg}}$.

\begin{lem} \label{Lem: Commutator Gamma-kg}
We have that
\begin{equation} \label{commutator Gamma kg with linear part}
\big[D_t - \oph(x\cdot\xi - p(\xi)), \Gamma^{kg}\big] = \oph(b),
\end{equation}
where
\begin{equation} \label{symbol of commutator}
\begin{split}
b(x,\xi) &= -\frac{h}{2i} (\partial\gamma)\Big(\frac{x-p'(\xi)}{\sqrt{h}}\Big)\cdot\Big(\frac{x-p'(\xi)}{\sqrt{h}}\Big) \chi(h^\sigma\xi) - \frac{\sigma h}{i}\gamma\Big(\frac{x-p'(\xi)}{\sqrt{h}}\Big) (\partial\chi)(h^\sigma\xi)\cdot(h^\sigma\xi) \\
&+\frac{i}{24}h^\frac{3}{2}\sum_{|\alpha|=3} (\partial^\alpha\gamma)\Big(\frac{x-p'(\xi)}{\sqrt{h}}\Big)(\partial^\alpha_\xi p'(\xi))\chi(h^\sigma\xi)+ r(x,\xi)
\end{split}
\end{equation}
and $r \in h^{5/2} S_{\frac{1}{2},\sigma}(\langle \frac{x-p'(\xi)}{\sqrt{h}}\rangle^{-N})$ for any $N\ge 0.$
Therefore, function $\widetilde{v}^\Sigma_{\Lambda_{kg}}$ is solution to
\begin{equation} \label{equation for vtilde Sigma Lambda}
\big[D_t - \oph(x\cdot\xi - p(\xi))\big]\widetilde{v}^\Sigma_{\Lambda_{kg}} = \Gamma^{kg}\oph(\Sigma(\xi))\big[h^{-1}r^{NF}_{kg}(t,tx)\big] + R_2(\widetilde{v})
\end{equation}
with $R_2(\widetilde{v})$ satisfying estimates \eqref{est L2 Linfty of R2(widetildev)}.
\proof
Recalling the definition \eqref{def Gamma_kg} of $\Gamma^{kg}$, one can prove by a straight computation that
\small
\begin{multline*}
\big[D_t, \Gamma^{kg}\big] = \frac{h}{i} \oph\Big((\partial\gamma)\Big(\frac{x-p'(\xi)}{\sqrt{h}}\Big)\cdot\frac{p''(\xi)\xi}{\sqrt{h}}\chi(h^\sigma\xi)\Big) \\
+ \frac{h}{2i}\oph\Big((\partial\gamma)\Big(\frac{x-p'(\xi)}{\sqrt{h}}\Big)\cdot\Big(\frac{x-p'(\xi)}{\sqrt{h}}\Big) \chi(h^\sigma\xi)\Big)  - \frac{(1+\sigma)h}{i}\oph\Big(\gamma\Big(\frac{x-p'(\xi)}{\sqrt{h}}\Big)(\partial\chi)(h^\sigma\xi)\cdot(h^\sigma\xi)\Big).
\end{multline*}\normalsize
Since the development of a commutator's symbol only contains odd-order terms, lemma \ref{Lem : a sharp b} gives that the symbol associated to $\big[\Gamma^{kg}, \oph(x\cdot\xi - p(\xi))\big]$ writes as
\begin{equation*}
\frac{h}{i}\Big\{\gamma\Big(\frac{x-p'(\xi)}{\sqrt{h}}\Big)\chi(h^\sigma\xi), x\cdot\xi - p(\xi)\Big\}+\frac{i}{24}h^\frac{3}{2}\sum_{|\alpha|=3} (\partial^\alpha\gamma)\Big(\frac{x-p'(\xi)}{\sqrt{h}}\Big)\chi(h^\sigma\xi)(\partial^\alpha_\xi p(\xi)) + r_5(x,\xi)
\end{equation*}
with $r_5\in h^{5/2}S_{\frac{1}{2},\sigma}(\langle \frac{x-p'(\xi)}{\sqrt{h}}\rangle^{-N})$ for any $N\ge 0$. Developing the above Poisson bracket one finds that
\begin{multline*}
\big[\Gamma^{kg}, \oph(x\cdot\xi - p(\xi))\big]= -\frac{h}{i}\oph\Big((\partial\gamma)\Big(\frac{x-p'(\xi)}{\sqrt{h}}\Big)\cdot\frac{p''(\xi)\xi}{\sqrt{h}}\chi(h^\sigma\xi)\Big) \\
- \frac{h}{i}\oph\Big((\partial\gamma)\Big(\frac{x-p'(\xi)}{\sqrt{h}}\Big)\cdot\Big(\frac{x-p'(\xi)}{\sqrt{h}}\Big) \chi(h^\sigma\xi)\Big) + \frac{h}{i}\oph\Big(\gamma\Big(\frac{x-p'(\xi)}{\sqrt{h}}\Big)(\partial\chi)(h^\sigma\xi)\cdot(h^\sigma\xi)\Big) \\+\frac{i}{24}h^\frac{3}{2}\sum_{|\alpha|=3} \oph\Big((\partial^\alpha\gamma)\Big(\frac{x-p'(\xi)}{\sqrt{h}}\Big)(\partial^\alpha_\xi p'(\xi))\chi(h^\sigma\xi)\Big) 
+ \oph(r_5(x,\xi)),
\end{multline*}
which summed to the previous commutator gives \eqref{symbol of commutator}.

The last part of the statement follows applying to equation \eqref{semi-classical KG equation} operators $\oph(\Sigma(\xi))$ (which commutes exactly with the linear part of the equation, evident in non semi-classical coordinates) and $\Gamma^{kg}$.
Since
\begin{multline*}
h\oph\Big((\partial\gamma)\Big(\frac{x-p'(\xi)}{\sqrt{h}}\Big)\cdot\Big(\frac{x-p'(\xi)}{\sqrt{h}}\Big) \chi(h^\sigma\xi)\Big)\widetilde{v}^\Sigma \\= \sum_{k=1}^2\oph\Big(\gamma^k\Big(\frac{x-p'(\xi)}{\sqrt{h}}\Big)\cdot (x-p'(\xi))(x_k-p'_k(\xi)) \Big)\widetilde{v}^\Sigma
\end{multline*}
with $\gamma^k(z):=(\partial\gamma)(z)z_k|z|^{-2}$, and 
\begin{equation*}
h^\frac{3}{2}\oph\Big((\partial^\alpha\gamma)\Big(\frac{x-p'(\xi)}{\sqrt{h}}\Big)(\partial^\alpha_\xi p'(\xi))\Big)  = h \oph\Big(\gamma^k_\alpha\Big(\frac{\xi-p'(\xi)}{\sqrt{h}}\Big)(\partial^\alpha_\xi p'(\xi))(x_k-p'_k(\xi))\Big)\widetilde{v}^\Sigma
\end{equation*}
with $\gamma^k_\alpha(z):=(\partial^\alpha\gamma)(z)z_k|z|^{-2}$, we obtain from inequalities \eqref{est:L2 Linfty Op(gamma) LLv} (resp. \eqref{est L2 Linfty for GammaTilde}) and \eqref{Lgamma vSigma} that $h\oph\big((\partial\gamma)\big(\frac{x-p'(\xi)}{\sqrt{h}}\big)\cdot\big(\frac{x-p'(\xi)}{\sqrt{h}}\big) \chi(h^\sigma\xi)\big)\widetilde{v}^\Sigma$ (resp. $h^{3/2}\oph\big((\partial^\alpha\gamma)\big(\frac{x-p'(\xi)}{\sqrt{h}}\big)(\partial^\alpha_\xi p'(\xi))\big)$, $|\alpha|=3$) is a remainder $R_2(\widetilde{v})$. The same holds true for
$\oph\big(\gamma\big(\frac{x-p'(\xi)}{\sqrt{h}}\big) (\partial\chi)(h^\sigma\xi)\cdot(h^\sigma\xi) \big)\widetilde{v}^\Sigma$, as follows combining symbolic calculus and lemma \ref{Lem : new estimate 1-chi}, because its symbol is supported for large frequencies $|\xi|\gtrsim h^{-\sigma}$.
From propositions \ref{Prop : Continuity on H^s} and \ref{Prop : Continuity from $L^2$ to L^inf} it immediately follows that $\oph(r_5)\widetilde{v}^\Sigma$ satisfies \eqref{est L2 of R2(widetildev)} and \eqref{est Linfty of R2(widetildev)}. 
\endproof
\end{lem}

\begin{prop}[Deduction of the ODE] \label{Prop: ODE for vtilde Sigma_Lambda}
There exists a family $(\theta_h(x))_h$ of $C^{\infty}_0$ functions, real valued, equal to 1 on the closed ball $\overline{B_{1-ch^{2\sigma}}(0)}$ and supported in $\overline{B_{1-c_1h^{2\sigma}}(0)}$, for some small $0<c_1<c$, $\sigma>0$, with $\|\partial^\alpha_x\theta_h\|_{L^\infty}=O(h^{-2|\alpha|\sigma})$ and $(h\partial_h)^k\theta_h$ bounded for every $k$, such that
\begin{equation} \label{dev linear part widetildev Sigma Lambda}
\oph(x\cdot\xi - p(\xi))\widetilde{v}^\Sigma_{\Lambda_{kg}} = -\phi(x)\theta_h(x)\widetilde{v}^\Sigma_{\Lambda_{kg}} + R_2(\widetilde{v}),
\end{equation}
where $\phi(x)=\sqrt{1-|x|^2}$ and $R_2(\widetilde{v})$ satisfies estimates \eqref{est L2 Linfty of R2(widetildev)}. Therefore,
$\widetilde{v}^\Sigma_{\Lambda_{kg}}$ is solution of the following non-homogeneous ODE:
\begin{equation} \label{ODE for widetildev Sigma Lambda}
D_t \widetilde{v}^\Sigma_{\Lambda_{kg}} = -\phi(x) \theta_h(x) \widetilde{v}^\Sigma_{\Lambda_{kg}} +  \Gamma^{kg}\oph(\Sigma(\xi))\big[h^{-1}r^{NF}_{kg}(t,tx)\big] + R_2(\widetilde{v}),
\end{equation}
with $r^{NF}_{kg}$ given by \eqref{def rNF-kg}.
\proof
The proof of the statement follows directly from lemma \ref{Lem:dev of a symbol at xi = -dvarphi(x)} if we observe that $\partial_\xi (x\cdot\xi - p(\xi)) = 0$ at $\xi = -d\phi(x)$ and $x\cdot (-d\phi(x))- p(-d\phi(x))=-\phi(x)$. Therefore, \eqref{dev linear part widetildev Sigma Lambda} holds and, injecting it in \eqref{equation for vtilde Sigma Lambda}, we obtain \eqref{ODE for widetildev Sigma Lambda}.
\endproof
\end{prop}

\begin{prop}[Propagation of the uniform estimate on $V$]\label{Prop:propagation_unif_est_V}
Let us fix $K_1>0$.
There exist two integers $n\gg \rho\gg 1$ sufficiently large, two constants $A,B>1$ sufficiently large, $\varepsilon_0\in ]0,(2A+B)^{-1}[$ sufficiently small, and $0\ll \delta\ll \delta_2\ll \delta_1\ll \delta_0\ll 1$ small, such that, for any $0<\varepsilon<\varepsilon_0$, if $(u,v)$ is solution to \eqref{wave KG system}-\eqref{initial data} in some interval $[1,T]$ for a fixed $T>1$, and $u_\pm, v_\pm$ defined in \eqref{def u+- v+-} satisfy a-priori estimates \eqref{est: bootstrap argument a-priori est} for every $t\in [1,T]$, then it also verify \eqref{est:bootstrap enhanced vpm} in the same interval $[1,T]$.
\proof
We warn the reader that, throughout the proof, we will denote by $C$, $\beta$ (resp. $\beta'$) two positive constants such that $\beta\rightarrow 0$ as $\sigma\rightarrow0$ (resp. $\beta'\rightarrow 0$ as $\delta_0,\sigma\rightarrow 0$). These constants may change line after line. We also remind that $h=1/t$.

In proposition \ref{Prop: normal forms on KG} we introduced function $v^{NF}$, defined from $v_{-}$ through \eqref{def vNF}, and proved that its $H^{\rho,\infty}$ norm differs from that of $v_{-}$ by a quantity satisfying \eqref{est_Hsinfty_vnf-v}.
Hence, from a-priori estimates \eqref{est: bootstrap upm}, \eqref{est: boostrap vpm}, \eqref{est: bootstrap Enn} and for $\theta\in ]0,1[$ sufficiently small (e.g. $\theta<1/4$) 
\begin{equation}\label{Hrho_infty_norm_v- in function of vNF}
\|v_{-}(t,\cdot)\|_{H^{\rho,\infty}}\le \|v^{NF}(t,\cdot)\|_{H^{\rho,\infty}}+ CA^{2-\theta}B^\theta\varepsilon^2t^{-\frac{5}{4}}.
\end{equation}
We successively introduced $\widetilde{v}$ in \eqref{def utilde vtilde} and decomposed it into the sum of functions $\widetilde{v}^\Sigma_{\Lambda_{kg}}$ and $\widetilde{v}^\Sigma_{\Lambda^c_{kg}}$ (see \eqref{def_both_vtilde_Sigma_Lambda}).
We will show in lemma \ref{Lem: from energy to norms in sc coordinates-KG} of appendix \ref{Appendix B} that, for any $s\le n$,
\begin{equation}\label{est_vtilde_energy}
\|\widetilde{v}(t,\cdot)\|_{H^s_h}+ \sum_{|\gamma|=1}^2\left\|\oph(\chi(h^\sigma\xi))\mathcal{L}^\gamma\widetilde{v}(t,\cdot) \right\|_{L^2}\le CB\varepsilon h^{-\beta'} 
\end{equation}
for all $t\in [1,T]$, so inequality \eqref{est Linfty of widetildev_Lambda,c} gives that
\begin{equation} \label{est_vSigma_Lambdakg}
\|\widetilde{v}^\Sigma_{\Lambda^c_{kg}}(t,\cdot)\|_{L^\infty}\le CB\varepsilon h^{\frac{1}{2}-\beta'}.
\end{equation}

As concerns $\vt^\Sigma_\Lkg$, we proved in proposition \ref{Prop: ODE for vtilde Sigma_Lambda} that it is solution to ODE \eqref{ODE for widetildev Sigma Lambda}, with $r^{NF}_{kg}$ given by \eqref{def rNF-kg} and satisfying \eqref{est L2 Linfty (cut-off) rNF-kg-new}, and $R_2(\widetilde{v})$ verifying \eqref{est L2 Linfty of R2(widetildev)}. From \eqref{est_vtilde_energy}, we then have that
\begin{equation*}
\|R_2(\widetilde{v})(t,\cdot)\|_{L^\infty}\le CB\varepsilon t^{-\frac{3}{2}+\beta'}.
\end{equation*}
We also have that
\begin{equation}\label{est_Gammakg_rNFkg}
\left\| \Gamma^{kg}\oph(\Sigma(\xi))[tr^{NF}_{kg}(t,tx)]\right\|_{L^\infty(dx)}\le C(A+B)AB\varepsilon^3 t^{-\frac{3}{2}+\beta'}.
\end{equation}
In fact, by symbolic calculus of lemma \ref{Lem : a sharp b} we derive that, for a fixed $N\in\mathbb{N}$ (e.g. $N>\rho$) and up to negligible multiplicative constants,
\begin{equation*}
\Gamma^{kg}\oph(\Sigma(\xi)) = \sum_{|\alpha|=0}^{N-1}h^\frac{|\alpha|}{2} \oph\Big((\partial^\alpha\gamma)\Big(\frac{x-p'(\xi)}{\sqrt{h}}\Big)\chi(h^\sigma\xi)(\partial^\alpha\Sigma)(\xi)\Big) + \oph(r_N(x,\xi)),
\end{equation*}
where $r_N\in h^\frac{N}{2}S_{\frac{1}{2},\sigma}(\langle\frac{x-p'(\xi)}{\sqrt{h}}\rangle^{-1})$. 
Choosing $N$ sufficiently large, we deduce from proposition \ref{Prop : Continuity from $L^2$ to L^inf}, the fact that $\|tw(t,t\cdot)\|_{L^2}=\|w(t,\cdot)\|_{L^2}$, inequality \eqref{est_rnfkg_L2} and a-priori estimates, that for every $t\in [1,T]$
\begin{equation*}
\Big\|\oph(r_N(x,\xi))[tr^{NF}_{kg}(t,tx)]\Big\|_{L^\infty(dx)}\le CA^2B\varepsilon^3 t^{-2}.
\end{equation*}
Using, instead, proposition \ref{Prop:Continuity Lp-Lp} with $p=+\infty$, inequality \eqref{est_Linfty_chi_rnfkg} in appendix \ref{Appendix B}, and that $h=t^{-1}$, we deduce that
\begin{multline*}
\sum_{|\alpha|=0}^{N-1}h^{\frac{|\alpha|}{2}}\left\|\oph\Big((\partial^\alpha\gamma)\Big(\frac{x-p'(\xi)}{\sqrt{h}}\Big)\chi(h^\sigma\xi)(\partial^\alpha\Sigma)(\xi)\Big) \oph(\chi_1(h^\sigma\xi))[tr^{NF}_{kg}(t,tx)]\right\|_{L^\infty} \\
\lesssim t^{1+\beta}\left\|\chi(t^{-\sigma}D_x) r^{NF}_{kg}(t,\cdot)\right\|_{L^\infty}\le C(A+B)AB\varepsilon^3 t^{-\frac{3}{2}+\beta'}.
\end{multline*}

Summing up, $\Gamma^{kg}\oph(\Sigma(\xi))[t^{-1}r^{NF}_{kg}(t,tx)]+ R_2(\widetilde{v}) = F_{kg}(t,x)$ with
\begin{equation*}
\|F_{kg}(t,\cdot)\|_{L^\infty}\le  [C(A+B)AB\varepsilon^3 +CB\varepsilon] t^{-\frac{3}{2}+\beta'},
\end{equation*}
Using equation \eqref{ODE for widetildev Sigma Lambda} we deduce that
\begin{equation}
\frac{1}{2}\partial_t |\widetilde{v}^\Sigma_{\Lambda_{kg}}(t,x)|^2 = \Im\left(\widetilde{v}^\Sigma_{\Lambda_{kg}} \overline{D_t \widetilde{v}^\Sigma_{\Lambda_{kg}}}\right)\le |\widetilde{v}^\Sigma_{\Lambda_{kg}}(t,x)| |F_{kg}(t,x)|
\end{equation}
and hence that
\begin{equation*}
\begin{split}
\|\widetilde{v}^\Sigma_{\Lambda_{kg}}(t,\cdot)\|_{L^\infty}&\le \|\widetilde{v}^\Sigma_{\Lambda_{kg}}(1,\cdot)\|_{L^\infty} + \int_1^t \|F_{kg}(\tau, \cdot)\|_{L^\infty} d\tau \\
& \le \|\widetilde{v}^\Sigma_{\Lambda_{kg}}(1,\cdot)\|_{L^\infty} + C(A+B)AB\varepsilon^3 +CB\varepsilon.
\end{split}
\end{equation*}
As $\|\widetilde{v}^\Sigma_{\Lambda_{kg}}(1,\cdot)\|_{L^\infty}\lesssim \|\widetilde{v}(1,\cdot)\|_{L^2}\le CB\varepsilon$ by proposition \ref{Prop : Continuity from $L^2$ to L^inf} and a-priori estimate \eqref{est: bootstrap Enn}, the above inequality together with\eqref{est_vSigma_Lambdakg} and definition \eqref{def utilde vtilde} of $\widetilde{v}$, gives that
\begin{equation*}
\|v^{NF}(t,\cdot)\|_{L^\infty}\le (C(A+B)AB\varepsilon^3 +CB\varepsilon)t^{-1},
\end{equation*}
which injected in \eqref{Hrho_infty_norm_v- in function of vNF} leads finally to \eqref{est:bootstrap enhanced vpm} if we take $A>1$ sufficiently large such that $CB<\frac{A}{3K_1}$, and $\varepsilon_0>0$ sufficiently small to verify $C(A+B)B\varepsilon^2_0 + CA^{1-\theta}B^\theta\varepsilon_0\le \frac{1}{3K_1}$.
\endproof
\end{prop}

\subsection{The derivation of the transport equation} \label{Subsection : The Derivation of the Transport Equation}

We now focus on the semi-classical wave equation satisfied by $\widetilde{u}$:
\begin{equation} \label{wave for utilde}
\big[D_t - \oph(x\cdot\xi - |\xi|)\big]\widetilde{u}(t,x) = h^{-1}\left[q_w(t,tx)+c_w(t,tx) + r^{NF}_w(t,tx)\right],
\end{equation}
with $q_w, c_w, r^{NF}_w$ given by \eqref{def_qw}, \eqref{def_cw}, \eqref{def rNF} respectively,
and on the derivation of the mentioned transport equation.
As we will make use several times of proposition \ref{Prop : continuity of Op(gamma1):X to L2} and inequalities \eqref{est: L2 Linfty with L}, we remind the reader about definition \eqref{def_Omega_h} of $\Omega_h$ and \eqref{def Mj} of $\Mcal_j$.
Also, $\theta_0(x)$ denotes a smooth radial cut-off function (often coming with operator $\Omega_h$) while $\chi\in C^\infty_0(\mathbb{R}^2)$ is equal to 1 in a neighbourhood of the origin and suitably supported. 

In order to recover a sharp estimate for $\widetilde{u}$ such as \eqref{est:a-priori_ut}, we study the behaviour of this function separately in different regions of the phase space $(x,\xi)\in\mathbb{R}^2\times\mathbb{R}^2$.
We start by fixing $\rho \in\mathbb{Z}$, and by introducing 
\begin{equation}\label{Sigma_j}
\Sigma_j(\xi):=
\begin{cases}
\langle\xi\rangle^\rho, \quad &\text{for } j=0,\\
\langle\xi\rangle^\rho \xi_j|\xi|^{-1}, \quad &\text{for } j=1,2.
\end{cases}
\end{equation}
Taking a smooth cut-off function $\chi_0$ equal to 1 in a neighbourhood of the origin, a Littlewood-Paley decomposition, and a small $\sigma>0$, we write the following for any $j\in \{0,1,2\}$: 
\begin{multline} \label{decomposition frequencies utilde}
\oph(\Sigma_j(\xi))\widetilde{u} =\oph(\Sigma_j(\xi)\chi_0(h^{-1}\xi))\widetilde{u} + \sum_k \oph\big(\Sigma_j(\xi)(1-\chi_0)(h^{-1}\xi)\varphi(2^{-k}\xi)\chi_0(h^\sigma\xi)\big)\widetilde{u}\\
+ \oph(\Sigma_j(\xi)(1-\chi_0)(h^\sigma\xi))\widetilde{u},
\end{multline}\normalsize
observing that the sum over $k$ is actually finite and restricted to set of indices $K:=\{k\in \mathbb{Z}: h\lesssim 2^k\lesssim h^{-\sigma}\}$.
From the classical Sobolev injection and the continuity on $L^2$ of the Riesz operator
\begin{equation} \label{utilde small frequencies}
\|\oph(\Sigma_j(\xi)\chi_0(h^{-1}\xi))\widetilde{u}(t,\cdot)\|_{L^\infty} = \|\Sigma_j(hD)\chi_0(D)\widetilde{u}(t,\cdot)\|_{L^\infty}\lesssim  \|\widetilde{u}(t,\cdot)\|_{L^2},
\end{equation}
while from the semi-classical Sobolev injection along with lemma \ref{Lem : new estimate 1-chi}
\begin{equation}\label{utilde large frequencies}
\|\oph(\Sigma_j(\xi)(1-\chi_0)(h^\sigma\xi))\|_{L^\infty}\lesssim h^N \|\widetilde{u}(t,\cdot)\|_{H^s_h},
\end{equation}
where $N=N(s)\ge 0$ if $s>0$ is sufficiently large.
The remaining terms in the right hand side of \eqref{decomposition frequencies utilde}, localised for frequencies $|\xi|\sim 2^k$, need a sharper analysis because a direct application of semi-classical Sobolev injection only gives that
\[\left\|\oph\big(\Sigma_j(\xi)(1-\chi_0)(h^{-1}\xi)\varphi(2^{-k}\xi)\chi_0(h^\sigma\xi)\big)\widetilde{u} \right\|_{L^\infty} \le 2^kh^{-1-\mu}\|\widetilde{u}\|_{L^2},\]
with $\mu=\sigma\rho$ if $\rho\ge 0$, 0 otherwise, and factor $2^k h^{-1-\mu}$ may grow too much when $h\rightarrow 0$. 

\begin{figure}[h]
\begin{center}
\begin{tikzpicture}[scale=1.9]

\draw[->] (0,-1.7) -- (0,1.7);
\draw[->] (-2,0) -- (2,0);
\node[below] at (1.9,0) {\small $x$};
\node[left] at (0,1.6) {\small $\xi$};

\draw[blue] (-1,-1.7) -- (-1,1.7);
\draw[blue] (1,-1.7) -- (1,1.7);
\node[below left] at (-1,0) {$-1$};
\node[below right] at (1,0) {$1$};
\node[right] at (1,1) {$\Lw$};

\draw[dashed] [domain=0:0.87] plot(\x, {\x/((1-\x^2)^(1/2))});
\draw[dashed] [domain=0:0.87] plot(-\x, {-\x/((1-\x^2)^(1/2))});
\end{tikzpicture}
\end{center}
\caption{Lagrangian for the wave equation}
\label{picture: Lw}
\end{figure}
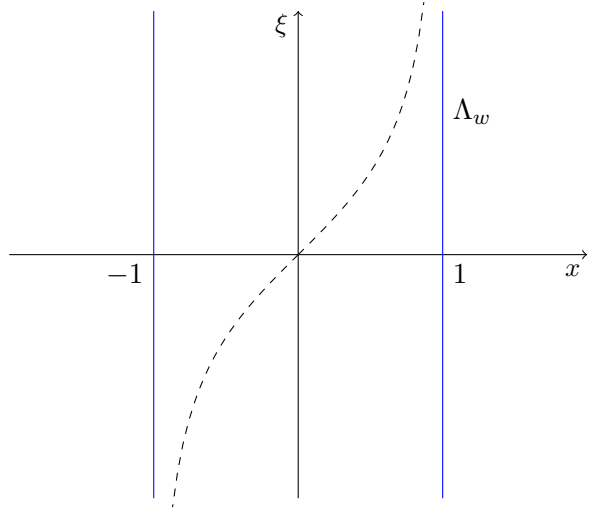

For any fixed $k\in K$, $\rho\in\mathbb{Z}$ and $j\in \{0,1,2\}$, let us introduce \index{uSigmak@$\widetilde{u}^{\Sigma,k}$, wave component in semi-classical setting, localised for frequencies $\sim 2^k$}
\begin{equation} \label{def utilde-Sigma,k}
\widetilde{u}^{\Sigma_j,k}(t,x) := \oph\big(\Sigma_j(\xi)(1-\chi_0)(h^{-1}\xi)\varphi(2^{-k}\xi)\chi_0(h^\sigma\xi)\big)\widetilde{u}(t,x)
\end{equation}
and observe that, from the commutation of the above operator with the linear part of equation \eqref{wave for utilde}, we get that $\widetilde{u}^{\Sigma_j,k}$ is solution to
\begin{equation} \label{wave equation u^k}
\begin{split}
& [D_t - \oph(x\cdot\xi - |\xi|)]\widetilde{u}^{\Sigma_j,k}(t,x)\\
 &=  h^{-1}\oph\big(\Sigma_j(\xi)(1-\chi_0)(h^{-1}\xi)\varphi(2^{-k}\xi)\chi_0(h^\sigma\xi)\big)\left[q_w(t,tx)+c_w(t,tx)+ r^{NF}_w(t,tx)\right]\\
&- i h \,\oph\big(\Sigma_j(\xi)(\partial\chi_0)(h^{-1}\xi)\cdot(h^{-1}\xi)\varphi(2^{-k}\xi)\big)\widetilde{u}-i\sigma h \, \oph\big(\Sigma_j(\xi)\varphi(2^{-k}\xi)(\partial\chi_0)(h^\sigma\xi))\cdot(h^\sigma\xi)\big)\widetilde{u}.
\end{split}
\end{equation}
We introduce the following manifold (see picture \ref{picture: Lw})\index{Lambdaw@$\Lw$, manifold associated to the wave equation}
\begin{equation} \label{def_Lw}
\Lw :=\left\{(x,\xi) : x-\frac{\xi}{|\xi|} = 0\right\},
\end{equation}
together with operator
\begin{equation} \label{def of Gamma_wk}
\Gamma^{w,k} :=
\oph\Big(\gamma\Big(\frac{x|\xi| - \xi}{h^{1/2-\sigma}}\Big)\psi(2^{-k}\xi)\Big),
\end{equation}
for some $\gamma\in C^\infty_0(\mathbb{R}^2)$ equal to 1 close to the origin and $\psi\in C^\infty_0(\mathbb{R}^2\setminus\{0\})$ equal to 1 on $supp\varphi$, whose symbol is localized in a neighbourhood of $\Lw \cap \{|\xi|\sim 2^k\}$ of size $h^{1/2-\sigma}$. 
We also define
\begin{subequations}\label{def_ukLambda-ukLambdac}
\begin{equation} \label{def of uk Lambda}
\widetilde{u}^{\Sigma_j,k}_{\Lambda_w} := \Gamma^{w,k} \widetilde{u}^{\Sigma_j,k}\,,
\end{equation}
\begin{equation} \label{def of uk Lambda complementary}
\widetilde{u}^{\Sigma_j,k}_{\Lambda^c_w} := \oph\Big(\big(1- \gamma\big)\Big(\frac{x|\xi| - \xi}{h^{1/2-\sigma}}\Big)\psi(2^{-k}\xi)\Big)\widetilde{u}^{\Sigma_j,k},
\end{equation}
\end{subequations}
so that $\widetilde{u}^{\Sigma_j,k} = \widetilde{u}^{\Sigma_j,k}_{\Lambda_w} + \widetilde{u}^{\Sigma_j,k}_{\Lambda^c_w}$.
We are going to prove that, if we suitably control the $L^2$ norm of $(\theta_0\Omega_h)^\mu\mathcal{M}^\nu\widetilde{u}^{\Sigma_j, k}$, for any $\mu, |\nu|\le 1$, then $\ut^{\Sigma_j,k}_{\Lw^c}$ is a $O_{L^\infty}(h^{-0})$ (see proposition \ref{Prop : est on widetildeu Lambda,c}).
As $h=t^{-1}$, this means that $\ut^{\Sigma_j,k}_{\Lambda^c_w}$ grows in time at a rate slower than the one expected for $\ut^{\Sigma_j,k}$ (that is $t^{1/2}$ after \eqref{est:a-priori_ut}).
Analogously to the Klein-Gordon case discussed in the previous subsection, the main contribution to $\ut^{\Sigma_j,k}$ is hence the one localized around $\Lw$ and represented by $\ut^{\Sigma_j,k}_\Lw$.
We will show that this function is solution to a transport equation (see proposition \ref{Prop: transport equation for uLambda}), from which we will be able to derive a suitable estimate of its uniform norm and to finally propagate \eqref{equivalence u- unf} (see proposition \ref{Prop: Propagation uniform estimate U,RU}).

\begin{prop} \label{Prop : est on widetildeu Lambda,c}
There exists a constant $C>0$ such that, for any $h\in ]0,1]$, $k\in K$, \index{utildeLambdaw@$\widetilde{u}^{\Sigma_j,k}_{\Lambda_w}$, wave component in semi-classical setting, localised around $\Lambda_w$}
\begin{subequations} \label{est L2 Linfty utilde Lambda,c}
\begin{equation}
\|\widetilde{u}^{\Sigma_j,k}_{\Lambda^c_w}(t,\cdot)\|_{L^2} \le C h^{\frac{1}{2}-\beta} \big(\|\widetilde{u}^{\Sigma_j, k}(t,\cdot)\|_{L^2} + \|\mathcal{M}\widetilde{u}^{\Sigma_j, k}(t,\cdot)\|_{L^2} \big)\,,
\end{equation} \index{utildeLambdaw_c@$\widetilde{u}^{\Sigma_j,k}_{\Lambda^c_w}$, wave component in semi-classical setting, localised away from $\Lambda_w$}
\begin{equation} \label{Linfy_norm_utilde_Lambdac}
\|\widetilde{u}^{\Sigma_j,k}_{\Lambda^c_w}(t,\cdot)\|_{L^\infty} \le C h^{-\beta}\sum_{\mu=0}^1\big(\|(\theta_0\Omega_h)^\mu\widetilde{u}^{\Sigma_j, k}(t,\cdot)\|_{L^2}+ \|(\theta_0\Omega_h)^\mu\mathcal{M}\widetilde{u}^{\Sigma_j, k}(t,\cdot)\|_{L^2}\big)\,,
\end{equation}
\end{subequations}
for a small $\beta>0$, $\beta\rightarrow 0$ as $\sigma\rightarrow 0$.
\proof
The proof is straightforward if one writes 
$$\widetilde{u}^{\Sigma_j,k}_{\Lambda^c_w} = \sum_{j=1}^2\oph\Big(\gamma^j_1\Big(\frac{x|\xi| - \xi}{h^{1/2-\sigma}}\Big)\Big(\frac{x_j|\xi| - \xi_j}{h^{1/2-\sigma}}\Big)\psi(2^{-k}\xi)\Big)\widetilde{u}^{\Sigma_j,k},$$
where
$\gamma^j_1(z):= \frac{(1-\gamma)(z)z_j}{|z|^2}$ is such that $|\partial^\alpha_z\gamma^j_1(z)|\lesssim \langle z \rangle^{-(|\alpha|+1)}$, and uses inequalities \eqref{est: L2 Linfty with L} with $a(x)=b_p(\xi)\equiv 1$.
\endproof
\end{prop}

\begin{lem} \label{Lem: commutator Gamma-wk with linear part}
Let $\widetilde{\varphi}\in C^\infty_0(\mathbb{R}^2\setminus\{0\})$ be such that $\widetilde{\varphi}\equiv 1$ on $supp\varphi$ and have a sufficiently small support so that $\psi\widetilde{\varphi}\equiv \psi$.
Then for any $k\in K$
\begin{equation} \label{commutator_Gammawk_equation}
\left[\Gamma^{w,k}, D_t - \oph\big((x\cdot \xi - |\xi|)\widetilde{\varphi}(2^{-k}\xi)\big)\right] \oph(\varphi(2^{-k}\xi)) = \oph(b(x,\xi)),
\end{equation}
where, for any $w\in L^2$ such that $\theta_0\Omega_h w, (\theta_0\Omega_h)^\mu\mathcal{M}w\in L^2(\mathbb{R}^2)$, for $\mu=0,1$,
\begin{subequations}
\begin{equation} \label{L2 norm Op(b)}
\| \oph(b(x,\xi))w\|_{L^2} \lesssim h^{1-\beta}\left(\|w\|_{L^2}+\|\mathcal{M}w\|_{L^2}\right),
\end{equation}
\begin{equation}\label{Linfty norm Op(b)}
\| \oph(b(x,\xi))w\|_{L^\infty} \lesssim h^{1-\beta}\sum_{\mu=0}^1\big(\|(\theta_0\Omega_h)^\mu w\|_{L^2}+\|(\theta_0\Omega_h)^\mu\mathcal{M}w\|_{L^2}\big),
\end{equation}
\end{subequations}
with $\beta>0$ small, $\beta\rightarrow 0$ as $\sigma\rightarrow 0$.
\proof
We warn the reader that most of the terms arising from the development of the commutator in the left hand side of \eqref{commutator_Gammawk_equation} satisfy a better $L^2$ estimate than \eqref{L2 norm Op(b)}, namely 
\begin{equation} \label{better L2 fwk}
\| \cdot\|_{L^2}\lesssim h^{\frac{3}{2}-\beta}\big(\|w\|_{L^2} + \|\mathcal{M}w\|_{L^2} \big).
\end{equation}
The only contribution whose $L^2$ norm is only a $O(h\|w\|_{L^2})$ is the integral remainder called $\widetilde{r}^k_N$, appearing in symbolic development \eqref{formula: symb dev}.

Since $\partial_t = -h^2\partial_h$, an easy computation shows that
\begin{equation} \label{commutator Dt Gamma^wk}
\begin{split}
[\Gamma^{w,k}, D_t] = & \Big(\frac{1}{2}+\sigma\Big)\frac{h}{i}\oph\Big((\partial\gamma)\Big(\frac{x|\xi|-\xi}{h^{1/2-\sigma}}\Big)\cdot\Big(\frac{x|\xi|-\xi}{h^{1/2-\sigma}}\Big)\psi(2^{-k}\xi)\Big) \\
& + \frac{h}{i} \oph\Big(\gamma\Big(\frac{x|\xi|-\xi}{h^{1/2-\sigma}}\Big)(\partial\psi)(2^{-k}\xi)\cdot(2^{-k}\xi)\Big).
\end{split}
\end{equation}
The first term in the above right hand side satisfies \eqref{better L2 fwk} and \eqref{Linfty norm Op(b)} after inequalities \eqref{est: L2 Linfty with L}.
The same estimates hold also for the latter one when it acts on $\oph(\varphi(2^{-k}\xi))w$, for the derivatives of $\psi$ vanish on the support of $\widetilde{\varphi}$ (and then of $\varphi$) as a consequence of our assumptions. In fact, if we introduce a smooth function $\widetilde{\psi}\in C^\infty_0(\mathbb{R}^2\setminus\{0\})$, equal to 1 on the support of $\partial\psi$ and such that $supp \widetilde{\psi}\cap supp\varphi = \emptyset$,
and use symbolic calculus we find that, for any fixed $N\in\mathbb{N}$,
\begin{multline*}
\oph\left(\gamma\Big(\frac{x|\xi|-\xi}{h^{1/2-\sigma}}\Big)(\partial\psi)(2^{-k}\xi)\cdot(2^{-k}\xi)\right)\oph(\varphi(2^{-k}\xi))\\
 = \oph\left(\gamma\Big(\frac{x|\xi|-\xi}{h^{1/2-\sigma}}\Big)(\partial\psi)(2^{-k}\xi)\cdot(2^{-k}\xi)\right) \oph\big(\widetilde{\psi}(2^{-k}\xi)\varphi(2^{-k}\xi)\big)- \oph(r^k_N),
\end{multline*}
where the first term in the above right hand side is 0,
and integral remainder $r^k_N$ is given by \small
\begin{multline*}
r^k_N = \Big(\frac{h}{2i}\Big)^N\sum_{|\alpha| = N}\frac{N (-1)^{|\alpha|}}{\alpha! (\pi h )^4}\int e^{\frac{2i}{h}(\eta \cdot z - y\cdot\zeta)}\int_0^1 \partial^\alpha_x \Big[\gamma\Big(\frac{x|\xi|-\xi}{h^{1/2-\sigma}}\Big)(\partial\psi)(2^{-k}\xi)\cdot (2^{-k}\xi)\Big]|_{(x+tz, \xi + t\zeta)} dt \\
\times \partial^\alpha_\xi \big(\widetilde{\psi}(2^{-k}\xi)\big)|_{(\xi + \eta)}\, dydz d\eta d\zeta.
\end{multline*}\normalsize
Developing explicitly the above derivatives and reminding definition \eqref{integral Ik} of integrals $I^k_{p,q}$, for general $k\in K$, $p,q\in\mathbb{Z}$, one recognizes that, up to some multiplicative constants, $r^k_N$ has the form
\begin{equation*}
h^{N-N(\frac{1}{2}-\sigma)}2^{-kN}I^k_{N,0}(x,\xi),
\end{equation*}
with $a,a',b_q\equiv 1$, $p=N$ and $\psi(2^{-k}\xi)$ replaced with $(\partial\psi)(2^{-k}\xi)\cdot(2^{-k}\xi)$.
Propositions \ref{Prop: L2 est of integral remainders} and \ref{Prop : Linfty est of integral remainders} imply then that 
\[\|\oph(r^k_N)\|_{\mathcal{L}(L^2)} + \|\oph(r^k_N)\|_{\mathcal{L}(L^2;L^\infty)} \lesssim h\]
if $N\in\mathbb{N}$ is chosen sufficiently large (e.g. $N>9$), which implies that the $\mathcal{L}(L^2)$ and $\mathcal{L}(L^2;L^\infty)$ norms of the latter operator in the right hand side of \eqref{commutator Dt Gamma^wk} is bounded by $h^2$.

As regards $[\Gamma^{w,k}, \oph((x\cdot\xi - |\xi|)\widetilde{\varphi}(2^{-k}\xi))]$, we first remind that the symbolic development of a commutator's symbol only contains odd order terms. Consequently, for a new fixed $N\in\mathbb{N}$ and up to multiplicative constants independent of $h, k$, the symbol of the considered commutator writes as
\begin{multline} \label{formula: symb dev}
h\Big\{\gamma\Big(\frac{x|\xi| - \xi}{h^{1/2-\sigma}}\Big), (x\cdot\xi - |\xi|)\widetilde{\varphi}(2^{-k}\xi)\Big\} \\
+ \sum_{\substack{3\le |\alpha|< N \\ |\alpha|=|\alpha_1|+|\alpha_2|}} h^{|\alpha|} \partial^{\alpha_1}_x \partial^{\alpha_2}_{\xi}\Big[\gamma\Big(\frac{x|\xi| - \xi}{h^{1/2-\sigma}}\Big)\Big]\partial^{\alpha_2}_x \partial^{\alpha_1}_{\xi}\big[(x\cdot\xi - |\xi|)\widetilde{\varphi}(2^{-k}\xi)\big] +  \widetilde{r}^k_N(x,\xi),
\end{multline}
with \small
\begin{multline*}
\widetilde{r}^k_N(x, \xi) = \Big(\frac{h}{2i}\Big)^N\sum_{|\alpha_1|+|\alpha_2| = N}\frac{N (-1)^{|\alpha_1|}}{\alpha! (\pi h )^4}\int e^{\frac{2i}{h}(\eta\cdot z - y\cdot\zeta)}\int_0^1 \partial^{\alpha_1}_x \partial^{\alpha_2}_{\xi}\Big[\gamma\Big(\frac{x|\xi| - \xi}{h^{1/2-\sigma}}\Big) \psi(2^{-k}\xi)\Big]\big|_{(x+tz,\xi+t\zeta)} dt \\ 
\times \partial^{\alpha_2}_x \partial^{\alpha_1}_{\xi}\big[(x\cdot\xi - |\xi|)\widetilde{\varphi}(2^{-k}\xi)\big]|_{(x+y,\xi +\eta)} \ dydzd\eta d\zeta \,.
\end{multline*}\normalsize
Since $\big\{\gamma\Big(\frac{x|\xi| - \xi}{h^{1/2-\sigma}}\Big), x\cdot \xi - |\xi|\big\} =0$ the Poisson braket in the above sum reduces to
\begin{equation*}
h\sum_{j,l}(\partial_j\gamma)\Big(\frac{x|\xi| - \xi}{h^{1/2-\sigma}}\Big)(\partial_j\widetilde{\varphi})(2^{-k}\xi)\big(\frac{x_l|\xi|-\xi_l}{h^{1/2-\sigma}}\big)(2^{-k}\xi_l)
\end{equation*}
and its quantization acting on $\oph(\varphi(2^{-k}\xi))w$ satisfies \eqref{better L2 fwk}, \eqref{Linfty norm Op(b)} because $\partial\widetilde{\varphi}$ vanishes on the support of $\varphi$.

An explicit calculation of terms of order $3\le |\alpha|<N$, with the help of lemma \ref{Lem : est on gamma for wave} and the observation that $|\alpha_2|\le 1$ because $(x\cdot\xi-|\xi|)\widetilde{\varphi}(2^{-k}\xi)$ is affine in $x$, shows that they are linear combination of products
\[h^{|\alpha| - |\alpha|(\frac{1}{2}-\sigma)}\gamma_{|\alpha|}\Big(\frac{x|\xi| - \xi}{h^{1/2-\sigma}}\Big)\widetilde{\varphi}(2^{-k}\xi)x^\nu b_1(\xi)\]
and 
\[h^{|\alpha| - (|\alpha|-1)(\frac{1}{2}-\sigma)}\widetilde{\gamma}\Big(\frac{x|\xi| - \xi}{h^{1/2-\sigma}}\Big)\widetilde{\varphi}(2^{-k}\xi)b_0(\xi)\]
for two new cut-off functions $\widetilde{\gamma},\widetilde{\varphi}$, $|\partial^\beta b_0(\xi)|\lesssim_\beta |\xi|^{-|\beta|}$, and $\nu\in\mathbb{N}^2$ of length at most 1. 
Furthermore, for $j=1,2$,
\begin{multline*}
h^{|\alpha| - |\alpha|(\frac{1}{2}-\sigma)}\gamma_{|\alpha|}\Big(\frac{x|\xi| - \xi}{h^{1/2-\sigma}}\Big)\widetilde{\varphi}(2^{-k}\xi)x_jb_1(\xi) = h^{|\alpha| - (|\alpha|-1)(\frac{1}{2}-\sigma)}\widetilde{\gamma}^j_{|\alpha|}\Big(\frac{x|\xi| - \xi}{h^{1/2-\sigma}}\Big)\widetilde{\varphi}(2^{-k}\xi)b_0(\xi)\\ + h^{|\alpha| - |\alpha|(\frac{1}{2}-\sigma)}\gamma_{|\alpha|}\Big(\frac{x|\xi| - \xi}{h^{1/2-\sigma}}\Big)\widetilde{\varphi}(2^{-k}\xi)\xi_jb_0(\xi),
\end{multline*}
with $\widetilde{\gamma}^j_{|\alpha|}(z):=\gamma_{|\alpha|}(z)z_j$. From propositions \ref{Prop : continuity Op(gamma) L2 to L2}, \ref{Prop : continuity of Op(gamma1):X to L2}, the fact that $|\alpha|\ge 3$ and $2^k\le h^{-\sigma}$ we deduce that the quantization of these $|\alpha|$-order terms acting on $\oph(\varphi(2^{-k}\xi))w$ satisfies \eqref{better L2 fwk}, \eqref{Linfty norm Op(b)}.

Finally, we notice that integral remainder $\widetilde{r}^k_N$ can be actually seen as the sum of two contributions, one of the form \eqref{integral rkN}, the other like \eqref{integral rtilde kN}, with $a\equiv 1$ and $p=1$.
Lemma \ref{Lem : remainder r^k_N} implies then that the $\mathcal{L}(L^2)$ and $\mathcal{L}(L^2;L^\infty)$ norms of $\oph(\widetilde{r}^k_N)$ are bounded by $h$ as foretold, which concludes the proof of the statement. 
\endproof
\end{lem}

\begin{lem}
Function $\widetilde{u}^{\Sigma_j,k}_{\Lambda_w}$ is solution to the following equation:
\begin{equation} \label{wave equation uk-Lambda}
\begin{split}
&\big[D_t - \oph\big((x\cdot \xi - |\xi|)\widetilde{\varphi}(2^{-k}\xi)\big)\big] \widetilde{u}^{\Sigma_j,k}_{\Lambda_w}(t,x) = f^w_k(t,x) \\
& +h^{-1} \Gamma^{w,k} \oph\big(\Sigma_j(\xi)(1-\chi_0)(h^{-1}\xi)\varphi(2^{-k}\xi)\chi_0(h^\sigma\xi)\big)\left[q_w(t,tx)+c_w(t,tx)+r^{NF}_w(t,tx)\right]\\
&- i h \,\Gamma^{w,k}\oph\big(\Sigma_j(\xi)(\partial\chi_0)(h^{-1}\xi)\cdot(h^{-1}\xi)\varphi(2^{-k}\xi)\big)\widetilde{u}\\
&-i\sigma h \, \Gamma^{w,k}\oph\big(\Sigma_j(\xi)\varphi(2^{-k}\xi)(\partial\chi_0)(h^\sigma\xi))\cdot(h^\sigma\xi)\big)\widetilde{u},
\end{split}
\end{equation}
where $\widetilde{\varphi}\in C^\infty_0(\mathbb{R}^2\setminus\{0\})$ is equal to 1 on $supp\varphi$, and there exist two constants $C, C'>0$ such that, for any $h\in]0,1], k\in K$,
\begin{subequations} \label{est L2 Linfty fk}
\begin{equation} \label{est L2 of fk}
\|f^w_k(t,\cdot)\|_{L^2} \le C h^{1-\beta} \big(\|\widetilde{u}^{\Sigma_j,k}(t,\cdot)\|_{L^2} + \|\mathcal{M}\widetilde{u}^{\Sigma_j,k}(t,\cdot)\|_{L^2} \big)\,,
\end{equation}
\begin{equation} \label{est Linfty of fk}
\|f^w_k(t,\cdot)\|_{L^\infty} \le C' h^{1-\beta}\sum_{\mu=0}^1\big(\|(\theta_0\Omega_h)^\mu\widetilde{u}^{\Sigma_j,k}(t,\cdot)\|_{L^2} + \|(\theta_0\Omega_h)^\mu\mathcal{M}\widetilde{u}^{\Sigma_j,k}(t,\cdot)\|_{L^2} \big)\,,
\end{equation}
\end{subequations}
with $\beta>0$ small, $\beta\rightarrow 0$ as $\sigma\rightarrow 0$.
\proof
If we consider a cut-off function $\widetilde{\varphi}\in C^\infty_0(\mathbb{R}^2\setminus\{0\})$ such that $\widetilde{\varphi}\equiv 1$ on the support of $\varphi$ ($\varphi$ being the truncation on $\widetilde{u}^{\Sigma_j,k}$'s frequencies), we have the exact equality
\[\oph(x\cdot\xi - |\xi|)\widetilde{u}^{\Sigma_j,k} = \oph((x\cdot\xi - |\xi|)\widetilde{\varphi}(2^{-k}\xi))\widetilde{u}^{\Sigma_j,k}.\] 
Moreover, if we assume that its support is sufficiently small so that $\psi\widetilde{\varphi}\equiv \widetilde{\varphi}$, and apply operator $\Gamma^{w,k}$ to equation \eqref{wave equation u^k}, lemma \ref{Lem: commutator Gamma-wk with linear part} gives us the result of the statement.
\endproof
\end{lem}

The transport equation we talked about at the beginning of this section will be deduced from equation \eqref{wave equation uk-Lambda} by suitably developing symbol $(x\cdot \xi - |\xi|)\widetilde{\varphi}(2^{-k}\xi)$. To do that, we first need to restrict the support of that symbol to bounded values of $x$ through the introduction of a new cut-off function $\theta(x)$. 
We remind that $\Sigma'$ is a concise notation that we use to indicate a linear combination of a finite number of terms of the same form.

\begin{lem} \label{Lem: PDE equation for utilde-k}
Let $0<D_1<D_2$ and $\theta=\theta(x)$ be a smooth function equal to 1 for $|x|\le D_1$ and supported for $|x|\le D_2$. Then,
\begin{multline} \label{linear op with theta}
\oph\big((x\cdot\xi - |\xi|)\widetilde{\varphi}(2^{-k}\xi)\big) = \oph\big(\theta(x)(x\cdot\xi - |\xi|)\widetilde{\varphi}(2^{-k}\xi)\big) + (1-\theta)(x)\oph((x\cdot\xi - |\xi|)\widetilde{\varphi}(2^{-k}\xi))\\ +
{\sum}' \widetilde{\theta}(x)\oph(\widetilde{\varphi}_1(2^{-k}\xi)) + \oph(r(x,\xi)),
 \end{multline}
where $\widetilde{\theta}$ is a smooth function supported for $D_1<|x|<D_2$, $\widetilde{\varphi}_1\in C^\infty_0(\mathbb{R}^2\setminus\{0\})$ and 
\[\|\oph(r)\|_{\mathcal{L}(L^2)} + \|\oph(r)\|_{\mathcal{L}(L^2;L^\infty)}\lesssim h.\]
Therefore, $\widetilde{u}^{\Sigma_j,k}_{\Lambda_{kg}}$ verifies 
\begin{equation} \label{PDO equation for utilde-k}
\begin{split}
&\Big[D_t - \oph\big(\theta(x)(x\cdot\xi - |\xi|)\widetilde{\varphi}(2^{-k}\xi)\big)\Big]\widetilde{u}^{\Sigma_j,k}_{\Lambda_w}(t,x) = f^w_k(t,x)   \\
& + (1-\theta)(x)\oph((x\cdot\xi - |\xi|)\widetilde{\varphi}(2^{-k}\xi))\widetilde{u}^{\Sigma_j,k}_{\Lambda_w} +
 {\sum}'\widetilde{\theta}(x)\oph(\widetilde{\varphi}_1(2^{-k}\xi))\widetilde{u}^{\Sigma_j,k}_{\Lambda_w}\\
&+ h^{-1} \Gamma^{w,k} \oph\big(\Sigma(\xi)(1-\chi_0)(h^{-1}\xi)\varphi(2^{-k}\xi)\chi_0(h^\sigma\xi)\big)\left[q_w(t,tx)+c_w(t,tx)+r^{NF}_w(t,tx)\right]\\
&- i h \,\Gamma^{w,k}\oph\big(\Sigma(\xi)(\partial\chi_0)(h^{-1}\xi)\cdot(h^{-1}\xi)\varphi(2^{-k}\xi)\big)\widetilde{u}\\
&-i\sigma h \, \Gamma^{w,k}\oph\big(\Sigma(\xi)\varphi(2^{-k}\xi)(\partial\chi_0)(h^\sigma\xi))\cdot(h^\sigma\xi)\big)\widetilde{u},
\end{split}
\end{equation}
where $f^w_k$ satisfies estimates \eqref{est L2 Linfty fk}.
\proof
Let $\theta(x)$ be the cut-off function of the statement. 
By proposition \ref{Prop: a sharp b} we have that\small
\begin{equation} \label{dev with 1-theta}
\begin{split}
& (1-\theta)(x) (x\cdot\xi - |\xi|)\widetilde{\varphi}(2^{-k}\xi)= (1-\theta)(x) \sharp \, \big[(x\cdot\xi - |\xi|)\widetilde{\varphi}(2^{-k}\xi) \big]\\ & - \frac{h}{2i} \partial\theta(x)\cdot \Big(x - \frac{\xi}{|\xi|}\Big) \widetilde{\varphi}(2^{-k}\xi) -\frac{2^{-k}h}{2i} (x\cdot\xi - |\xi|)\partial\theta(x)\cdot(\partial \widetilde{\varphi})(2^{-k}\xi) + r_2^k(x,\xi) \\
 &= \, (1-\theta)(x) \sharp \, \big[(x\cdot\xi - |\xi|) \widetilde{\varphi}(2^{-k}\xi)\big] - \frac{h}{2i}\Big[\partial\theta(x)\cdot x \Big] \sharp \,\widetilde{\varphi}(2^{-k}\xi) + \frac{h}{2i}\sum_{l=1}^2\partial_l\theta(x)\sharp \Big[\frac{\xi_l}{|\xi|}\widetilde{\varphi}(2^{-k}\xi) \Big]\\
& -\frac{h}{2i}\sum_{j,l=1}^2\Big[\partial_j\theta(x)x_l \Big]\sharp \,\Big[(\partial_j\widetilde{\varphi})(2^{-k}\xi)(2^{-k}\xi_l)\Big] + \frac{h}{2i}\sum_{l=1}^2\partial_l\theta(x)\sharp \Big[(2^{-k}|\xi|)(\partial_l\widetilde{\varphi})(2^{-k}\xi)\Big] + r^k_2(x,\xi)+ \widetilde{r}^k_2(x,\xi),
\end{split}
\end{equation}\normalsize
where $\partial\theta$ is supported for $D_1<|x|<D_2$, and $r^k_2(t,x)$ (resp. $\widetilde{r}^k_2(t,x)$) is a linear combination of integrals of the form
\begin{equation*}
\frac{h^22^{-k}}{(\pi h)^2}\int e^{\frac{2i}{h}\eta\cdot z}\int_0^1 \theta(x+tz) (1-t)^2 dt \ x^\nu \widetilde{\varphi}(2^{-k}(\xi+\eta)) dzd\eta,
\end{equation*}
with $|\nu|=0,1$ (resp. $|\nu|=0$), for some new $\theta, \widetilde{\varphi} \in C^\infty_0(\mathbb{R}^2\setminus\{0\})$. 
By writing $x$ as $(x+tz)-tz$, using that $ze^{\frac{2i}{h}\eta\cdot z}= \big(\frac{h}{2i}\big)\partial_\xi e^{\frac{2i}{h}\eta\cdot z}$, and making an integration by parts, one can express$r^k_2(t,x)$ as the sum over $|\nu|=0,1$ of integrals such as
\begin{equation*}
\frac{h^22^{-k}(h2^{-k})^\nu}{(\pi h)^2}\int e^{\frac{2i}{h}\eta\cdot z}\int_0^1 \theta(x+tz) f(t) dt \ \widetilde{\varphi}(2^{-k}(\xi+\eta)) dzd\eta,
\end{equation*}
for some new smooth $\theta, f, \widetilde{\varphi}$, and show that for any $\alpha,\beta\in\mathbb{N}^2$
\[\big|\partial^\alpha_x\partial^\beta_\xi \big[(r^k_2+\widetilde{r}^k_2)(x,h\xi)\big]\big|\lesssim_{\alpha,\beta} h^22^{-k}\lesssim_{\alpha,\beta}h.\] 
Thus $(r^k_2+ \widetilde{r}^k_2)(x,h\xi) \in hS_{0}(1)$, which means, by classical results on pseudo-differential operators (see for instance \cite{hormander:the_analysis_III}), that
\[\oph((r^k_2 + \widetilde{r}^k_2)(x,\xi)) = Op^w((r^k_2+\widetilde{r}^k_2)(x, h\xi))\in \mathcal{L}(L^2)\] with norm $O(h)$.
Furthermore, one can also show that $\|\oph(r^k_2+\widetilde{r}^k_2)\|_{\mathcal{L}(L^2;L^\infty)}\lesssim h$ using lemma \ref{Lemma on inequalities for Op(A)} and the fact that, by making some integrations by parts, for any multi-indices $\alpha, \beta\in\mathbb{N}^2$ and a new $\widetilde{\varphi}\in C^\infty_0(\mathbb{R}^2\setminus\{0\})$
\begin{align*}
\left\|\partial^\alpha_y\partial^\beta_\xi \left[(r^k_2+ \widetilde{r}^k_2)\Big(\frac{x+y}{2},h\xi\Big)\right]\right\|_{L^2(d\xi)} \lesssim h^22^{-k} \left\|\int \langle \eta\rangle^{-3} |\widetilde{\varphi}(2^{-k}h(\xi+\eta))| d\eta \right\|_{L^2(d\xi)}\lesssim h.
\end{align*}
These considerations, along with the continuity of $\Gamma^{w,k}$ on $L^2$, uniformly in $h$ and $k$ (see proposition \ref{Prop : continuity Op(gamma) L2 to L2}), imply that $\oph(r^k_2 + \widetilde{r}_2^k)\widetilde{u}^{\Sigma_j,k}_{\Lambda_w}$ is a remainder $f^w_k$.
\endproof
\end{lem}

\begin{lem} \label{Lemma : development of linear part}
We have that
\[|\xi| - x\cdot\xi = \frac{1}{2}(1-|x|^2)x\cdot\xi + e(x,\xi)\]
with
\begin{equation} \label{def of e(x,xi)}
e(x,\xi) = \frac{1}{2}|\xi|\Big|x - \frac{\xi}{|\xi|}\Big|^2 + \frac{1}{2}\Big(\Big(x-\frac{\xi}{|\xi|}\Big)\cdot\xi\Big)\Big(x- \frac{\xi}{|\xi|}\Big)\cdot\Big(x + \frac{\xi}{|\xi|}\Big).
\end{equation}
\proof
\begin{equation*}
\begin{split}
|\xi| -x\xi & = \frac{1}{2}|\xi| \left|x-\frac{\xi}{|\xi|}\right|^2 + \frac{1}{2}|\xi| (1-|x|^2) \\
& =  \frac{1}{2}|\xi| \left|x-\frac{\xi}{|\xi|}\right|^2 + \frac{1}{2}(|\xi| - x\cdot\xi) (1-|x|^2) + \frac{1}{2}(1-|x|^2)x\cdot\xi \\
& = \underbrace{\frac{1}{2}|\xi| \left|x-\frac{\xi}{|\xi|}\right|^2 + \frac{1}{2}\left(\left(\frac{\xi}{|\xi|}-x\right)\cdot\xi\right)\left(\frac{\xi}{|\xi|} -x\right)\cdot\left(\frac{\xi}{|\xi|} + x\right)}_{e(x,\xi)} + \frac{1}{2}(1-|x|^2)x\cdot\xi\,.
\end{split}
\end{equation*}
\end{lem}

\begin{lem} \label{Lem: preliminary on Op(e)}
Let $\gamma, \theta\in C^\infty_0(\mathbb{R}^2)$ and $\widetilde{\varphi}\in C^\infty_0(\mathbb{R}^2\setminus\{0\})$ be such that $\widetilde{\varphi}\equiv 1$ on the support of $\varphi$ and have a sufficiently small support so that $\psi \widetilde{\varphi}\equiv \widetilde{\varphi}$. Let also
\begin{equation} \label{def B(x,xi)}
B(x,\xi):=\gamma\Big(\frac{x|\xi| - \xi}{h^{1/2-\sigma}}\Big)\widetilde{\varphi}(2^{-k}\xi)\theta(x)\Big(x_m - \frac{\xi_m}{|\xi|}\Big), \qquad m\in \{1,2\}.
\end{equation}
For any function $w\in L^2(\mathbb{R}^2)$ such that $\mathcal{M}w\in L^2(\mathbb{R}^2)$, any $m,n\in \{1,2\}$,
\begin{subequations}
\begin{equation} \label{est L2 Op(e)}
\left\|\oph\Big(\theta(x)\widetilde{\varphi}(2^{-k}\xi)\Big(x_m - \frac{\xi_m}{|\xi|}\Big)(x_n|\xi|-\xi_n)\Big)\Gamma^{w,k}w\right\|_{L^2}\lesssim 
h^{1-\beta} \big(\|w\|_{L^2} + \|\mathcal{M}w\|_{L^2}\big),
\end{equation}
\begin{multline} \label{est Linfty Op(e)}
\left\|\oph\Big(\theta(x)\widetilde{\varphi}(2^{-k}\xi)\Big(x_m - \frac{\xi_m}{|\xi|}\Big)(x_n|\xi|-\xi_n)\Big)\Gamma^{w,k}w\right\|_{L^\infty}\lesssim\\
h^{1-\beta} \big(\|w\|_{L^2} + \|\mathcal{M}w\|_{L^2}\big)
+ h^{-\beta}\|\oph\big(B(x,\xi)\xi\big)\mathcal{M}w\|_{L^2},
\end{multline}
\end{subequations}
with $\beta>0$ small, $\beta\rightarrow 0$ as $\sigma\rightarrow 0.$
\proof
After lemma \ref{Lemma : on the enhanced symbolic product} with $p=0$ we have that \small
\begin{equation*}
\oph\left(\theta(x)\widetilde{\varphi}(2^{-k}\xi)\Big(x_m - \frac{\xi_m}{|\xi|}\Big)(x_n|\xi|-\xi_n)\right)\Gamma^{w,k}w= \oph\left(B(x,\xi)(x_n|\xi| - \xi_n)\right)w+ \oph(r^k_0(x,\xi))w,
\end{equation*}\normalsize
and the $L^2$ (resp. $L^\infty$) norm of the latter term in the above right hand side is bounded by the right hand side of \eqref{est L2 Op(e)} (resp. of \eqref{est Linfty Op(e)}) after inequality \eqref{est L2 Op(rkp) enhanced} (resp. \eqref{est Linfty Op(rkp) enhanced}).
Moreover, the $L^2$ norm of $\oph\big(B(x,\xi)(x_n|\xi| - \xi_n)\big)w$ is also bounded by the right hand side of \eqref{est L2 Op(e)} as straightly follows from emma \ref{Lemma : symbolic product development}.
It only remains to prove that the $L^\infty$ norm of this term is bounded by the right hand side of \eqref{est Linfty Op(e)}.

We first consider a new cut-off function $\widetilde{\varphi}_1 \in C^\infty_0(\mathbb{R}^2\setminus\{0\})$, equal to 1 on $supp\widetilde{\varphi}$ so that its derivatives vanish against $\varphi$, and use symbolic calculus to write
\begin{equation*} 
\oph\big(B(x,\xi)(x_n|\xi| - \xi_n)\big) = \oph(\widetilde{\varphi}_1(2^{-k}\xi))\oph\big(B(x,\xi)(x_n|\xi| - \xi_n)\big) + \oph(r^k_{N,1}(x,\xi)),
\end{equation*}
where $r^k_{N,1}(x,\xi)$ is obtained using \eqref{r_N 2}.
Up to interchange the role of variables $y$ and $z$ (resp. $\eta$ and $\zeta$) and to consider $e^{\frac{2i}{h}(y\cdot\zeta-\eta\cdot\zeta)}$ instead of $e^{\frac{2i}{h}(\eta\cdot z-y\cdot\zeta)}$ (which does not affect estimate \eqref{est L2 and Linfty rkN r'kN}), $r^k_{N,1}$ is analogous to integral \eqref{integral rtilde kN} with $p=1$. 
Therefore, if $N\in\mathbb{N}$ is chosen sufficiently large (e.g. $N> 11$), lemma \ref{Lem : remainder r^k_N} implies that $\|\oph(r^k_{N,1})\|_{\mathcal{L}(L^2;L^\infty)}=O(h).$

Since $\widetilde{\varphi}_1$ localises frequencies $\xi$ in an annulus, the classical Sobolev injection gives that
\begin{multline*}
\left\|\oph(\widetilde{\varphi}_1(2^{-k}\xi))\oph\big(B(x,\xi)(x_n|\xi| - \xi_n)\big)w\right\|_{L^\infty} \\
\lesssim \log h \left\|\oph\big(B(x,\xi)(x_n|\xi| - \xi_n)\big)w\right\|_{L^2} + \left\|D_x \oph\big(B(x,\xi)(x_n|\xi| - \xi_n)\big)w\right\|_{L^2}.
\end{multline*}
As previously said, the former norm in the above right hand side satisfies inequality \eqref{est L2 Op(e)}.
As concerns the latter one, we remark that thanks to the specific structure of symbol $B(x,\xi)$ its first derivative with respect to $x$ does not lose any factor $h^{-1/2+\sigma}$, as
\begin{equation} \label{partial x B(x,xi)}
\partial_x \left[\gamma\Big(\frac{x|\xi|-\xi}{h^{1/2-\sigma}}\Big)\right]\widetilde{\varphi}(2^{-k}\xi)\theta(x)\Big(x_m - \frac{\xi_m}{|\xi|}\Big) = (\partial\gamma)\Big(\frac{x|\xi| -\xi}{h^{1/2-\sigma}}\Big)\widetilde{\varphi}(2^{-k}\xi)\theta(x)\Big(\frac{x_m|\xi| - \xi_m}{h^{1/2-\sigma}}\Big).
\end{equation}
Consequently, by using symbolic calculus we derive that
\begin{multline*}
D_x \oph\big(B(x,\xi)(x_n|\xi| - \xi_n)\big)w = h^{-1}\oph\big(B(x,\xi)(x_n|\xi| - \xi_n)\xi\big)w \\
+ {\sum}'  \oph\Big(\gamma\Big(\frac{x|\xi| - \xi}{h^{1/2-\sigma}}\Big)\widetilde{\varphi}(2^{-k}\xi)a(x)b_0(\xi)(x_j|\xi| - \xi_j) \Big)w,
\end{multline*} 
where ${\sum}'$ is a concise notation to indicate linear combinations, $j\in\{m,n\}$ and $\gamma, \widetilde{\varphi}, a$ are some new smooth functions with $a(x)$ compactly supported.
Again by lemma \ref{Lemma : symbolic product development} the $L^2$ norms of latter contributions in the above right hand side are bounded by $h^{1-\beta}(\|w\|_{L^2}+\|\mathcal{M}w\|_{L^2})$.

Finally, we observe that symbol $B(x,\xi)\xi$ can be seen as
\begin{equation} \label{B(x,xi)xi 1}
\gamma\Big(\frac{x|\xi| - \xi}{h^{1/2-\sigma}}\Big)(x_m|\xi| - \xi_m)\widetilde{\varphi}(2^{-k}\xi)\theta(x)b_0(x),
\end{equation}
which implies, after lemma \ref{Lem: Gamma with double argument-wave}, that
\begin{equation*}
h^{-1}\oph\big(B(x,\xi)(x_n|\xi|-\xi_n)\xi\big)w
 =\oph\big(B(x,\xi)\xi\big)\mathcal{M}_nw+ O_{L^2}(h^{1-\beta}(\|w\|_{L^2} + \|\mathcal{M}w\|_{L^2})).
\end{equation*}
\endproof
\end{lem}

\begin{lem} \label{Lemma : estimate of e(x,xi)}
Let $e(x,\xi)$ be the symbol defined in \eqref{def of e(x,xi)}, $\theta\in C^\infty_0(\mathbb{R}^2)$, and $\widetilde{\varphi}\in C^\infty_0(\mathbb{R}^2\setminus\{0\})$ with sufficiently small support so that $\psi \widetilde{\varphi}\equiv \widetilde{\varphi}$. 
If a-priori estimates \eqref{est: bootstrap argument a-priori est} are satisfied for every $t\in [1,T]$, for some fixed $T>1$, there exists a constant $C>0$ such that\small
\begin{equation} \label{est L2 Linfty e(x,xi)}
\left\|\oph\Big(\theta(x)\widetilde{\varphi}(2^{-k}\xi)e(x,\xi)\Big)\widetilde{u}^{\Sigma_j,k}_{\Lambda_w}(t,\cdot)\right\|_{L^2} + \left\|\oph\Big(\theta(x)\widetilde{\varphi}(2^{-k}\xi)e(x,\xi)\Big)\widetilde{u}^{\Sigma_j,k}_{\Lambda_w}(t,\cdot)\right\|_{L^\infty}\le CB\varepsilon h^{1-\beta'}
\end{equation}\normalsize
for every $t\in[1,T]$, with $\beta>0$ small, $\beta\rightarrow 0$ as $\sigma\rightarrow 0$.
\proof
We warn the reader that, throughout this proof, $C, \beta$ and $\beta'$ will denote three positive constants that may change line after line, with $\beta\rightarrow 0$ as $\sigma\rightarrow 0$ (resp. $\beta'\rightarrow 0$ as $\sigma, \delta_1\rightarrow 0$). 

Since symbol $e(x,\xi)$ writes as
\begin{equation*} 
e(x,\xi) = \frac{1}{2}\sum_{m=1}^2 \Big(x_m - \frac{\xi_m}{|\xi|}\Big)(x_m|\xi|-\xi_m)   + \frac{1}{2}\sum_{m,n=1}^2 \Big(x_m - \frac{\xi_m}{|\xi|}\Big)(x_n|\xi|-\xi_n)\Big(\frac{\xi_m}{|\xi|}\frac{\xi_n}{|\xi|} +  x_n \frac{\xi_m}{|\xi|} \Big)\,,
\end{equation*}
it follows that the $L^2$ norm of $\oph\big(\theta(x)\widetilde{\varphi}(2^{-k}\xi)e(x,\xi)\big)\widetilde{u}^{\Sigma_j,k}_{\Lambda_w}$ satisfies inequality \eqref{est L2 Linfty e(x,xi)} after lemmas \ref{Lem: preliminary on Op(e)} and \ref{Lem: from energy to norms in sc coordinates-WAVE} in appendix \ref{Appendix B}. Moreover, from lemma \ref{Lem: preliminary on Op(e)}
\begin{multline*}
\left\|\oph\Big(\theta(x)\widetilde{\varphi}(2^{-k}\xi)e(x,\xi)\Big)\widetilde{u}^{\Sigma_j,k}_{\Lambda_w}\right\|_{L^\infty} \lesssim h^{1-\beta}\left(\|\widetilde{u}^{\Sigma_j,k}(t,\cdot)\|_{L^2} + \|\mathcal{M}\widetilde{u}^{\Sigma_j,k}(t,\cdot)\|_{L^2}\right)\\
 + h^{-\beta}\|\oph\big(B(x,\xi)\xi\big)\mathcal{M}\widetilde{u}^{\Sigma_j,k}(t,\cdot)\|_{L^2},
\end{multline*}
with $B(x,\xi)$ defined in \eqref{def B(x,xi)}.
The aim of the proof is then to show that the $L^2$ norm of $\oph\big(B(x,\xi)\xi\big)\mathcal{M}\widetilde{u}^{\Sigma_j,k}$ is estimated by the right hand side of \eqref{est L2 Linfty e(x,xi)}.

First of all, we remind that $B(x,\xi)\xi$ can be seen as a symbol of the form \eqref{B(x,xi)xi 1}.
From proposition \ref{Prop : continuity Op(gamma) L2 to L2} we hence have that
\begin{subequations}
\begin{equation}\label{norm_operator B(x,xi)xi-1}
\|\oph\big(B(x,\xi)\xi\big)\|_{\mathcal{L}(L^2)}=O( h^{\frac{1}{2}-\beta}),
\end{equation}
while from inequality \eqref{est: L2 of gamma1 with L}
\begin{equation}\label{norm_operator B(x,xi)xi-2}
\|\oph\big(B(x,\xi)\xi\big)w\|_{L^2} \lesssim h^{1-\beta}(\|w\|_{L^2}+\|\mathcal{M}w\|_{L^2}).
\end{equation}
\end{subequations}
We also recall definition \eqref{def utilde-Sigma,k} of $\widetilde{u}^{\Sigma_j,k}$, use the concise notation $\phi^j_k(\xi)$ for its symbol $\Sigma_j(\xi)(1-\chi_0)(h^{-1}\xi)\varphi(2^{-k}\xi)\chi_0(h^\sigma\xi)$, and observe that
\begin{equation} \label{comm_Mcal_phik}
\begin{gathered}
\Big[\Mcal_n, \oph(\phi^j_k(\xi))\Big] = -\frac{1}{2i}\oph\big(|\xi|\partial_n\phi^j_k(\xi)\big),\\
\Big\| \Big[\Mcal_n, \oph(\phi^j_k(\xi))\Big] \Big\|_{\Lcal(L^2)} = O(h^{-\sigma}),
\end{gathered}
\end{equation} 
after propositions \ref{Prop: a sharp b} and \ref{Prop : continuity Op(gamma) L2 to L2}.

Using \eqref{comm_Mcal_phik} and recalling relation \eqref{relation between Zju and Mj utilde-new}, we find that for any $n=1,2$,
\begin{equation*}
\begin{split}
& \|\oph\big(B(x,\xi)\xi\big)\mathcal{M}_n\widetilde{u}^{\Sigma_j,k}(t,\cdot)\|_{L^2}\lesssim \|\oph\big(B(x,\xi)\xi\big)\oph(\phi^j_k(\xi))[t (Z_n u^{NF})(t,tx)]\|_{L^2(dx)}\\
&+ \left\|\oph\big(B(x,\xi)\xi\big)\oph\big(\xi_n|\xi|^{-1}\phi^j_k(\xi))\widetilde{u}(t,\cdot) \right\|_{L^2} + \left\|\oph\big(B(x,\xi)\xi\big)\oph\big(|\xi|\partial_n\phi^j_k(\xi)\big)\widetilde{u}(t,\cdot) \right\|_{L^2}\\
&+ \left\|\oph\big(B(x,\xi)\xi\big)\oph(\phi^j_k(\xi))\left[t(tx_n)\left[q_w(t,tx)+c_w(t,tx)+\rnfw(t,tx)\right]\right]\right\|_{L^2(dx)},
\end{split}
\end{equation*}
with $u^{NF}$ defined in \eqref{def uNF}, $q_w$, $c_w$ and $\rnfw$ given by \eqref{def_qw}, \eqref{def_cw} and \eqref{def rNF} respectively.
Evidently, after \eqref{norm_operator B(x,xi)xi-2} and a further commutation of $\mathcal{M}$ with $\oph\big(\xi_n|\xi|^{-1}\phi^j_k(\xi)\big)$ and $\oph\big(|\xi|\partial_n\phi^j_k(\xi)\big)$ respectively, the second and third $L^2$ norm in the above right hand side are estimated by 
\begin{equation*}
h^{1-\beta}(\|\widetilde{u}(t,\cdot)\|_{L^2}+ \|\oph(\chi(h^\sigma\xi))\mathcal{M}\widetilde{u}(t,\cdot)\|_{L^2}),
\end{equation*}
for some $\chi\in C^\infty_0(\mathbb{R}^2)$.
They are hence bounded by $CB\varepsilon h^{1-\beta'}$ by lemma \ref{Lem: from energy to norms in sc coordinates-WAVE}.

\smallskip
\underline{\textbf{$\bullet$ Estimate of $\| \oph\big(B(x,\xi)\xi)\oph(\phi^j_k(\xi))[t(Z_nu^{NF})(t,tx)]\|_{L^2}$:}}
This $L^2$ norm is basically estimated in terms of the $L^2$ norm of $(Z^\mu u)_{-}$, for $|\mu|\le 2$. In fact, after definition \eqref{def uNF} and equality \eqref{commutator_Z_Dt-|D|}
\begin{multline}\label{dev_Znu} 
(Z_nu^{NF})(t,tx) = (Z_nu)_{-}(t,tx)+\Big(\frac{D_n}{|D_x|}u_{-}\Big)(t,tx)\\ -\frac{i}{4(2\pi)^2}\sum_{l\in\{+,-\}} \Big[Z_n\int e^{iy\cdot\xi} D_l(\xi,\eta) \hat{v}_l(\xi-\eta)\hat{v}_l(\eta) d\xi d\eta\Big]\big|_{y=tx},
\end{multline}
with $D_l$ given by \eqref{def Dj1j2(xi,eta)}.
On the one hand, taking a new smooth cut-off function $\theta_1$ equal to 1 on the support of $\theta$, using \eqref{symbolic dev 1} with $\widetilde{a}=\theta_1$, together with \eqref{est L2 Op(rkp)}, proposition \ref{Prop : continuity Op(gamma) L2 to L2}, and \eqref{comm_Mcal_phik}, we deduce that
\begin{multline*}
\|\oph\big(B(x,\xi)\xi)\oph(\phi^j_k(\xi))[t(Z_nu)_{-}(t,tx)]\|_{L^2(dx)}\\
\lesssim \sum_{m=1}^2 h\|\theta_1(x)\oph(\phi^j_k(\xi))\mathcal{M}_m[t(Z_nu)_{-}(t,tx)]\|_{L^2(dx)}+ h^{1-\beta}\|(Z_nu)_{-}(t,\cdot)\|_{L^2}.
\end{multline*}
After relation \eqref{relation ZmZnu Mutilde Zn-new},
\begin{multline*}
\|\theta_1(x)\oph(\phi^j_k(\xi))\mathcal{M}_m[t(Z_nu)_{-}(t,tx)]\|_{L^2} \lesssim \|(Z_mZ_nu)_{-}(t,\cdot)\|_{L^2}+\|(Z_nu)_{-}(t,\cdot)\|_{L^2} \\
+ \Big\|\theta_1\Big(\frac{x}{t}\Big)\phi^j_k(D_x)\left[x_m Z_n\Nlw\right](t,\cdot) \Big\|_{L^2}.
\end{multline*}
Moreover,
\begin{equation*}
\theta_1\Big(\frac{x}{t}\Big)\phi^j_k(D_x)x_m =t \theta_{1,m}\Big(\frac{x}{t}\Big)\phi^j_k(D_x) + \theta_1\Big(\frac{x}{t}\Big)[\phi^j_k(D_x),x_m],
\end{equation*}
where $\theta_{1,m}(z)=\theta_1(z)z_m$, and $[\phi^j_k(D_x),x_m]$ is a bounded operator on $L^2$ with norm $O(t)$, as one can check computing its associated symbol and using that $2^{-k}\lesssim h^{-1}=t$. 
Therefore, using also inequality \eqref{L2 est NLwZn} together with a-priori estimates \eqref{est: bootstrap argument a-priori est} we deduce that
\begin{equation}\label{est_A}
\begin{split}
&\left\| \oph\big(B(x,\xi)\xi)\oph(\phi^j_k(\xi))\Big[t(Z_nu)_{-}(t,tx)\Big]\right\|_{L^2(dx)} \\
&\lesssim \sum_{|\mu|=1}^2h \|(Z^\mu u)_{-}(t,\cdot)\|_{L^2} + \|Z_nV(t,\cdot)\|_{H^1}\|V(t,\cdot)\|_{H^{2,\infty}}+ \big[\|V(t,\cdot)\|_{H^1} \\
&+  \|V(t,\cdot)\|_{L^2}\left(\|U(t,\cdot)\|_{H^{1,\infty}}+ \|\mathrm{R}_1U(t,\cdot)\|_{H^{1,\infty}}\right) + \|V(t,\cdot)\|_{L^\infty}\|U(t,\cdot)\|_{H^1}\big]\|V(t,\cdot)\|_{H^{1,\infty}}\\
& \le CB\varepsilon h^{1-\frac{\delta_1}{2}}.
\end{split}
\end{equation}
On the other hand, it is a straight consequence of \eqref{norm_operator B(x,xi)xi-2}, \eqref{comm_Mcal_phik} and lemma \ref{Lem: from energy to norms in sc coordinates-WAVE} that
\begin{multline}\label{est_B}
\left\| \oph\big(B(x,\xi)\xi)\oph(\phi^j_k(\xi)) [t (D_n|D_x|^{-1}u)_{-}(t,tx)]\right\|_{L^2}\\
\lesssim h^{1-\beta}(\|\widetilde{u}(t,\cdot)\|_{L^2}+\|\oph(\chi(h^\sigma\xi))\mathcal{M}\widetilde{u}(t,\cdot)\|_{L^2})\le CB\varepsilon h^{1-\frac{\delta_2}{2}}.
\end{multline}
Finally, by symbolic calculus and \eqref{partial x B(x,xi)} we have that
\begin{equation}\label{Op(B)hD}
\oph(B(x,\xi)\xi)=\oph(B(x,\xi))(hD_x) + \frac{h}{2i}
\oph\big(\partial_xB(x,\xi)\big),
\end{equation}
where $\partial_xB$ is of the form
\begin{equation}\label{derivative_B(x,xi)xi}
\gamma\Big(\frac{x|\xi|-\xi}{h^{1/2-\sigma}}\Big)\widetilde{\varphi}(2^{-k}\xi)\theta(x)b_0(\xi)
\end{equation}
for some new $\gamma, \theta\in C^\infty_0(\mathbb{R}^2)$. Consequently, by proposition \ref{Prop : continuity Op(gamma) L2 to L2}
\begin{multline}\label{Op(Bxi)integral_Zn}
\Big\|\oph(B(x,\xi)\xi)\oph(\phi^j_k(\xi))\Big[t Z_n\int e^{i y\cdot\xi} D_l(\xi,\eta) \hat{v}_l(\xi-\eta)\hat{v}_l(\eta) d\xi d\eta\Big]\big|_{y=tx}\Big\|_{L^2(dx)}  \\
\lesssim\Big\|\chi(t^{-\sigma}D_x) D_xZ_n\int e^{ix\cdot\xi} D_l(\xi,\eta) \hat{v}_l(\xi-\eta)\hat{v}_l(\eta) d\xi d\eta\Big\|_{L^2(dx)}\\+  
h\Big\|\chi(t^{-\sigma}D_x)  Z_n\int e^{ix\cdot\xi} D_l(\xi,\eta) \hat{v}_l(\xi-\eta)\hat{v}_l(\eta) d\xi d\eta\Big\|_{L^2(dx)}
\end{multline} 
and the above right hand side is bounded by 
\begin{multline*}
h^{-\beta} \left(\|V(t,\cdot)\|_{H^1}+ \|V(t,\cdot)\|_{L^2}(\|U(t,\cdot)\|_{H^{1,\infty}}+\|\mathrm{R}_1U(t,\cdot)\|_{H^{1,\infty}})+ \|V(t,\cdot)\|_{L^\infty}\|U(t,\cdot)\|_{H^1}\right)\\
\times \left(\|V(t,\cdot)\|_{H^{14,\infty}}+h\|V(t,\cdot)\|_{H^{13}}\right) + h^{-\beta}\|Z_nV(t,\cdot)\|_{L^2}\|V(t,\cdot)\|_{H^{17,\infty}}
\end{multline*}
after inequalities \eqref{est:L2 Z integral D with cut-off}, \eqref{est: L2 DjZ integral D with cut-off} and \eqref{Hs norm DtV} with $s=0$.
From a-priori estimates \eqref{est: bootstrap argument a-priori est} we then deduce that the left hand side of \eqref{Op(Bxi)integral_Zn} is bounded by $CB\varepsilon h^{1-\beta'}$, which implies, together with equality \eqref{dev_Znu} and estimates\eqref{est_A}, \eqref{est_B}, that the $L^2$ norm of contribution $ \oph\big(B(x,\xi)\xi)\oph(\phi^j_k(\xi))[t(Z_nu^{NF})(t,tx)]$ is estimated with the right hand side of \eqref{est L2 Linfty e(x,xi)}.

\underline{\textbf{$\bullet$ Estimate of} $\|\oph(B(x,\xi)\xi) \left[t(tx_n)q_w(t,tx)\right]\|_{L^2(dx)}$}:
After definition \eqref{def_qw} of $q_w(t,x)$ and \eqref{def utilde vtilde} of $\widetilde{v}$, we first notice that
\begin{equation}\label{def_qtilde_w}
tq_w(t,tx) = \frac{h}{2} \Im \left[\overline{\widetilde{v}}\, \oph(\xi_1)\vt - \overline{\oph\Big(\frac{\xi_1}{\langle \xi\rangle}\Big)\vt}\cdot \oph\Big(\frac{\xi\xi_1}{\langle \xi\rangle}\Big)\vt\right](t,x) =:\widetilde{q}_w(t,x),
\end{equation}
where \begin{equation}\label{norm_L2_qtilde_w}
\|\widetilde{q}_w(t,\cdot)\|_{L^2}\lesssim h\|\vt(t,\cdot)\|_{H^{1,\infty}}\|\vt(t,\cdot)\|_{H^1}.
\end{equation}
Then \small
\begin{equation*}
\|\oph(B(x,\xi)\xi)\oph(\phi^j_k(\xi)) \left[t(tx_n)q_w(t,tx)\right]\|_{L^2(dx)} = h^{-1}\|\oph(B(x,\xi)\xi)\oph(\phi^j_k(\xi))\left[x_n\widetilde{q}_w(t,x)\right]\|_{L^2(dx)}.
\end{equation*} \normalsize
Since $B(x,\xi)$ is compactly supported in $x$ and
\[\Big\|\Big[\oph\big(B(x,\xi)\xi\big)\oph(\phi^j_k(\xi)), x_n\Big]\Big\|_{\Lcal(L^2)} = O(h^{\frac{1}{2}-\beta}),\]
as follows from symbolic calculus, \eqref{norm_operator B(x,xi)xi-1}, equality \eqref{derivatives of gamma_n 2} and proposition \ref{Prop : continuity Op(gamma) L2 to L2}, we can morally reduce ourselves to the study of the $L^2$ norm of 
\begin{equation*}
h^{-1} \oph(B(x,\xi)\xi)\oph(\phi^j_k(\xi))\widetilde{q}_w(t,x)
\end{equation*}
up to a $O_{L^2}(h^{-1/2-\beta}\|\widetilde{q}_w\|_{L^2})$.
Using \eqref{Op(B)hD}, \eqref{derivative_B(x,xi)xi}, together with proposition \ref{Prop : continuity Op(gamma) L2 to L2}, we deduce that
\begin{equation*}
h^{-1}\left\|\oph(B(x,\xi)\xi)\oph(\phi^j_k(\xi))\widetilde{q}_w(t,\cdot)\right\|_{L^2}\lesssim h^{-1}\| \oph(\phi^j_k(\xi)) (hD_x)\widetilde{q}_w(t,\cdot)\|_{L^2}+ \|\widetilde{q}_w(t,\cdot)\|_{L^2},
\end{equation*}
so from lemma \ref{Lem: hD|V|2} below, estimates \eqref{norm_L2_qtilde_w}, \eqref{est:a-priori_vt}, and lemmas \ref{Lem: from energy to norms in sc coordinates-KG}, \ref{Lem_appendix: estimate L2vtilde} in appendix \ref{Appendix B}, we conclude that
\begin{multline}\label{norm_(hD)_qtilde}
h^{-1}\| \oph(\phi^j_k(\xi)) (hD_x)\widetilde{q}_w(t,\cdot)\|_{L^2}\\
\lesssim h^{1-\beta}\Big(\|\vt(t,\cdot)\|_{H^s}+\sum_{|\mu|=1}^2\|\oph(\chi(h^\sigma\xi))\mathcal{L}^\mu \vt(t,\cdot)\|_{L^2}\Big) \|\vt(t,\cdot)\|_{H^{1,\infty}}\le CB\varepsilon h^{1-\beta'};
\end{multline}

\underline{\textbf{$\bullet$ Estimate of} $\|\oph(B(x,\xi)\xi) \oph(\phi^j_k(\xi))(t(tx_n)c_w(t,tx)\|_{L^2(dx)}$}:
As for the previous estimate, we can reduce to the study of the $L^2$ norm of 
\begin{equation*}
\oph(B(x,\xi)\xi)\oph(\phi^j_k(\xi))[t^2 c_w(t,tx)],
\end{equation*}
up to a $O_{L^2}\big(h^{-1/2-\beta}\|\oph(\chi(h^\sigma\xi))[tc_w(t,tx)]\|_{L^2(dx)}\big)$ for some $\chi\in C^\infty_0(\mathbb{R}^2)$. 
So using \eqref{norm_operator B(x,xi)xi-1}, the fact that $\|tw(t,t\cdot)\|_{L^2}=\|w(t,\cdot)\|_{L^2}$, and \eqref{est_cw_L2} with $s>0$ sufficiently large so that $N(s)>2$, we obtain that for a new $\chi_1\in C^\infty_0(\mathbb{R}^2)$
\begin{equation}
\begin{split}
& \|\oph(B(x,\xi)\xi)\oph(\phi^j_k(\xi))[t^2c_w(t,t\cdot)]\|_{L^2}\lesssim h^{-\frac{1}{2}-\beta}\left\|\chi(t^{-\sigma}D_x)c_w(t,\cdot)\right\|_{L^2}\\
&\lesssim h^{-\frac{1}{2}-\beta}\left\|\chi_1(t^{-\sigma}D_x)(\vnf-v_{-})(t,\cdot)\right\|_{L^2}\left(\|V(t,\cdot)\|_{H^{2,\infty}}+\|\vnf(t,\cdot)\|_{H^{1,\infty}}\right) \\
& + h^{\frac{3}{2}}\left\| (\vnf - v_{-})(t,\cdot)\right\|_{H^1}\left(\|V(t,\cdot)\|_{H^s}+\|\vnf(t,\cdot)\|_{H^s}\right).
\end{split}
\end{equation}
Then inequalities \eqref{est_Hs_vnf-v} with $s=1$ and \eqref{est_chi_vnf-v}, together with a-priori estimates, give that
\[\|\oph(B(x,\xi)\xi)\oph(\phi^j_k(\xi))[t^2c_w(t,t\cdot)]\|_{L^2}\le CB\varepsilon h^{1-\beta'}.\]

\underline{\textbf{$\bullet$ Estimate of} $\|\oph(B(x,\xi)\xi) \oph(\phi^j_k(\xi))(t(tx_n)r^{NF}_w(t,tx)\|_{L^2(dx)}$}:
Analogously, from \eqref{est L2 rNF} and \eqref{est: bootstrap argument a-priori est} we obtain that
\begin{multline}\label{Op(B)rNFw}
\|\oph(B(x,\xi)\xi)\oph(\phi^j_k(\xi))[t^2 r^{NF}_w(t,t\cdot)]\|_{L^2}\lesssim h^{-\frac{1}{2}-\beta}\|\chi(t^{-\sigma}D_x) r^{NF}_w(t,\cdot)\|_{L^2}\\
\lesssim h^{-\frac{1}{2}-\beta}\|V(t,\cdot)\|^2_{H^{13,\infty}}\|U(t,\cdot)\|_{H^1}\le CB\varepsilon h^{\frac{3}{2}-\beta'}.
\end{multline}
\endproof
\end{lem}

\begin{lem}\label{Lem: hD|V|2}
Let $\varphi\in C^\infty_0(\mathbb{R}^2\setminus \{0\})$, $k\in K$ and $a_j(\xi)$ be two smooth real symbols of order $j=0,1$. Then
\begin{multline}\label{(hD)|V|2}
\left\|\oph(\varphi(2^{-k}\xi))(hD_x)\left[ \overline{a_0(hD_x)\vt} \, a_1(hD_x)\vt\right](t,\cdot)\right\|_{L^2}\\
\lesssim h^{1-\beta}\Big(\|\vt(t,\cdot)\|_{H^s_h}+\sum_{|\mu|=1}^2\|\oph(\chi(h^\sigma\xi))\mathcal{L}^\mu \vt(t,\cdot)\|_{L^2}\Big)\|\vt(t,\cdot)\|_{H^{1,\infty}_h}.
\end{multline}
\proof
Let us split both $\vt$ in the left hand side of \eqref{(hD)|V|2} into the sum $\vtl+\vtlc$, with $\vtl, \vtlc$ introduced in \eqref{def_both_vtilde_Sigma_Lambda} with $\Sigma_j\equiv 1$. Remind that $\vtlc$ satisfies inequality \eqref{est H1h widetildev_Lambda,c} and that
\begin{equation*}
\|a_0(hD_x)\vt(t,\cdot)\|_{L^\infty}+ \|a_0(hD_x)\vt_{\Lambda_{kg}}(t,\cdot)\|_{L^\infty} \lesssim h^{-\beta}\|\vt(t,\cdot)\|_{H^{1,\infty}_h},
\end{equation*}
for a small $\beta>0$, $\beta\rightarrow 0$ as $\sigma\rightarrow 0$, as follows
from lemma \ref{Prop:Continuity Lp-Lp} with $p=+\infty$ and the uniform continuity of $a_0(hD_x)$ from $H^{1,\infty}$ to $L^\infty$.
Therefore, using the continuity on $L^2$ of $\oph(\varphi(2^{-k}\xi))(hD_x)$ with norm $O(2^k)$ and the fact that $2^k\lesssim h^{-\sigma}$
we deduce that, for any $w_1,w_2\in \{\vt, \vtl,\vtlc\}$ with at least one $w_j$ equal to $\vtlc$,
\begin{multline*}
\left\|\oph(\varphi(2^{-k}\xi))(hD_x)\left[\overline{a_0(hD_x)w_1} a_1(hD_x)w_2\right]\right\|_{L^2}\\
\lesssim h^{1-\beta}\Big(\|\vt(t,\cdot)\|_{H^s_h}+\sum_{|\mu|=1}^2\|\oph(\chi(h^\sigma\xi))\mathcal{L}^\mu \vt(t,\cdot)\|_{L^2}\Big)\|\vt(t,\cdot)\|_{H^{1,\infty}_h}.
\end{multline*}
We are thus reduced to proving inequality \eqref{(hD)|V|2} for
\begin{equation*}
\left\|\oph(\varphi(2^{-k}\xi))(hD_x)\left[ \overline{a_0(hD_x)\vt_{\Lambda_{kg}}} \, a_1(hD_x)\vt_{\Lambda_{kg}}\right](t,\cdot)\right\|_{L^2}.
\end{equation*}
Furthermore, by means of lemma \ref{Lem:dev of a symbol at xi = -dvarphi(x)} we can replace the action of $a_j(hD_x)$ in the above $L^2$ norm, for $j=0,1$, with the multiplication operator by a real function, up to new remainders bounded in $L^2$ by the right hand side of \eqref{(hD)|V|2}. 
In fact,
\begin{equation*}
a_j(hD_x)\vt_{\Lambda_{kg}}=\theta_h(x)a_j(-d\phi(x))\vt_{\Lambda_{kg}}+R_1(\vt), \quad j=0,1,
\end{equation*}
where $\theta_h$ is a smooth cut-off function as in the statement of lemma \ref{Lem:dev of a symbol at xi = -dvarphi(x)} and $R_1(\vt)$ satisfies \eqref{est L2 of R1(widetildev)}.
Now
\begin{equation*}
hD_x |\vt_{\Lambda_{kg}}|^2 = \left[\oph(\xi+d\phi(x)\theta_h(x))\vt_{\Lambda_{kg}}\right]\overline{\vt_{\Lambda_{kg}}} - \vt_{\Lambda_{kg}}\left[\overline{\oph(\xi+d\phi(x)\theta_h(x))\vt_{\Lambda_{kg}}}\right],
\end{equation*}
and from lemma \ref{Lem: (xi+dphi)Op(gamma)} below
\begin{equation*}
\|\oph(\xi+d\phi(x)\theta_h(x))\vt_{\Lambda_{kg}}(t,\cdot)\|_{L^2}\lesssim h^{1-\beta}\sum_{|\mu|=0}^1\|\oph(\chi(h^\sigma\xi))\mathcal{L}^{\mu} \vt(t,\cdot)\|_{L^2}.
\end{equation*}
This implies, after having applied the Leibniz rule and proposition \ref{Prop:Continuity Lp-Lp}, that
\begin{multline*}
\left\|hD_x \left[a_0(-d\phi(x))a_1(-d\phi(x))\theta^2_h(x)|\vt_{\Lambda_{kg}}|^2(t,\cdot)\right]\right\|_{L^2} \\ \lesssim  h^{1-\beta}\Big(\|\vt(t,\cdot)\|_{H^s_h}+\sum_{|\mu|=1}^2\|\oph(\chi(h^\sigma\xi))\mathcal{L}^\mu \vt(t,\cdot)\|_{L^2}\Big)\|\vt(t,\cdot)\|_{L^\infty}
\end{multline*}
and the conclusion of the statement.
\endproof
\end{lem}

\begin{lem} \label{Lem: (xi+dphi)Op(gamma)}
Let $\gamma, \chi\in C^\infty_0(\mathbb{R}^2)$ be equal to 1 in a neighbourhood of the origin, $\sigma>0$ small, $(\theta_h(x))_h$ be a family of $C^{\infty}_0(B_1(0))$ functions, equal to 1 on the support of $\gamma\big(\frac{x-p'(\xi)}{\sqrt{h}}\big)\chi(h^\sigma\xi)$, with $\|\partial^\alpha_x\theta_h\|_{L^\infty}=O(h^{-2|\alpha|\sigma})$ and $(h\partial_h)^k\theta_h$ bounded for every $k$. Let also $\phi(x)=\sqrt{1-|x|^2}$.
Then for every $j=1,2$
\begin{multline*}
\left\|\oph(\xi_j+ d_j\phi(x)\theta_h(x))\oph\Big(\gamma\Big(\frac{x-p'(\xi)}{\sqrt{h}}\Big)\chi(h^\sigma\xi)\Big)\widetilde{v}(t,\cdot) \right\|_{L^2}\\
\lesssim h^{1-\beta}\sum_{|\mu|=0}^2\|\oph(\chi(h^\sigma\xi))\mathcal{L}^\mu \widetilde{v}(t,\cdot)\|_{L^2},
\end{multline*}
with $\beta>0$ small, $\beta\rightarrow 0$ as $\sigma\rightarrow 0$.
\proof
By symbolic calculus of lemma \ref{Lem : a sharp b} and the fact that $\theta_h\equiv 1$ on the support of $\gamma\big(\frac{x-p'(\xi)}{\sqrt{h}}\big)\chi(h^\sigma\xi)$, we have that, for any $j=1,2$,
\begin{equation}\label{(xi-dphi)V}
\begin{split}
\oph(\xi_j+ &d_j\phi(x)\theta_h(x))\oph\Big(\gamma\Big(\frac{x-p'(\xi)}{\sqrt{h}}\Big)\chi(h^\sigma\xi)\Big)\widetilde{v} = \oph\Big(\gamma\Big(\frac{x-p'(\xi)}{\sqrt{h}}\Big)\chi(h^\sigma\xi)(\xi_j+d_j\phi(x))\Big)\widetilde{v}  \\
&+\frac{\sqrt{h}}{2i}\oph\Big((\partial_j\gamma)\Big(\frac{x-p'(\xi)}{\sqrt{h}}\Big)\chi(h^\sigma\xi)\Big)\widetilde{v} \\
&-\frac{\sqrt{h}}{2i}\sum_{k,l=1}^2\oph\Big((\partial_l\gamma)\Big(\frac{x-p'(\xi)}{\sqrt{h}}\Big)p''_{k,l}(\xi)\partial_k( d_j\phi(x)\theta_h(x)) \chi(h^\sigma\xi)\Big)\widetilde{v}\\
& +\frac{h^{1+\sigma}}{2i}\sum_{k=1}^2\oph\Big(\gamma\Big(\frac{x-p'(\xi)}{\sqrt{h}}\Big)\partial_k( d_j\phi(x)\theta_h(x)) (\partial_k\chi)(h^\sigma\xi)\Big)\widetilde{v} + \oph(r_2(x,\xi))\widetilde{v},
\end{split}
\end{equation}
with $r_2\in h^{1-4\sigma}S_{\frac{1}{2},\sigma}(\langle \frac{x-p'(\xi)}{\sqrt{h}}\rangle^{-1})$.
On the one hand, as
\begin{equation*}
\oph\Big(\gamma\Big(\frac{x-p'(\xi)}{\sqrt{h}}\Big)\chi(h^\sigma\xi)(\xi_j-d_j\phi(x))\Big)\widetilde{v} = \sum_{k=1}^2\oph\Big(\gamma\Big(\frac{x-p'(\xi)}{\sqrt{h}}\Big)\chi(h^\sigma\xi)\widetilde{e}^j_k(x,\xi)(x_k-p'_k(\xi))\Big)\widetilde{v},
\end{equation*}
with $\widetilde{e}^j_k$ satisfying \eqref{etildekj} on the support of $\gamma\big(\frac{x-p'(\xi)}{\sqrt{h}}\big)\chi(h^\sigma\xi)$, the $L^2$ norm of the first term in the right hand side of \eqref{(xi-dphi)V} can be estimated using \eqref{est:L2 Op(gamma)LLv}.

On the other hand, as $\partial\gamma$ vanishes in a neighbourhood of the origin, the $L^2$ norm of the second and third term in the right hand side of \eqref{(xi-dphi)V} can be estimated using \eqref{est_1L-L2}. 

The two remaining contributions to the right hand side of \eqref{(xi-dphi)V}, that already carry the right power of $h$, can be estimated with $h^{1-\beta}\|\vt(t,\cdot)\|_{L^2}$ simply by proposition \ref{Prop : Continuity on H^s}.
\endproof
\end{lem}

We can finally state the following result:

\begin{prop}[Deduction of the transport equation] \label{Prop: transport equation for uLambda}
For any fixed $T>1$, $D>0$, let $\mathcal{C}^T_D:=\{(t,x):1\le t \le T,  |x|\le D\}$ be the truncated cylinder, and assume that estimates \eqref{est: bootstrap argument a-priori est} are satisfied in time interval $[1,T]$.
Then function 
\begin{equation} \label{def_u_Sigma_Lw}
\widetilde{u}^{\Sigma_j}_{\Lambda_w}(t,x):=\sum_k\widetilde{u}^{\Sigma_j,k}_{\Lambda_w}(t,x)
\end{equation}
is solution to the following transport equation:
\begin{equation} \label{eq: transport equation for uSigmaLambda}
\left[D_t +\frac{1}{2}(1-|x|^2)x\cdot (hD_x) + \frac{h}{2i}(1-2|x|^2)\right]\widetilde{u}^{\Sigma}_{\Lambda_w}(t,x) = F_w(t,x) , \quad \forall (t,x)\in \mathcal{C}^T_D,
\end{equation}
and there exists some constant $C>0$ such that
\begin{equation} \label{Linfty_norm_Fw}
\|F_w(t,\cdot)\|_{L^\infty} \le CB\varepsilon h^{1-\beta'}
\end{equation}
for some $\beta'>0$ small, $\beta'\rightarrow 0$ as $\sigma,\delta_1\rightarrow 0$.
\proof
By the assumption in the statement, all that we are going to say is to be meant in time interval $[1,T]$.
We remind the reader that, by the definition of $\widetilde{u}^{\Sigma_j,k}_{\Lambda_w}$ in \eqref{def of uk Lambda} and of $\widetilde{u}^{\Sigma_j,k}$ in \eqref{def utilde-Sigma,k}, the sum defining $\widetilde{u}^{\Sigma_j}_{\Lambda_w}$ is finite and restricted to indices $k\in K:=\{k\in\mathbb{Z}: h\lesssim 2^k\lesssim h^{-\sigma}\}$.
Also, we warn the reader that, throughout the proof, $C$ and $\beta$ will denote two positive constants that may change line after line, with $\beta\rightarrow 0$ as $\sigma\rightarrow 0$.

In lemma \ref{Lem: PDE equation for utilde-k} we proved that function $\widetilde{u}^{\Sigma_j,k}_{\Lambda_w}$ is solution to \eqref{PDO equation for utilde-k} with $f^w_k$ verifying \eqref{est L2 Linfty fk}.
Hence, by lemma \ref{Lem: from energy to norms in sc coordinates-WAVE} we derive that $f^w_k$ is a remainder of the form $F_w$ satisfying \eqref{Linfty_norm_Fw}.

For seek of compactness, we denote symbol $\Sigma_j(\xi)(1-\chi_0)(h^{-1}\xi)\varphi(2^{-k}\xi)\chi_0(h^\sigma\xi)$ in the right hand side of \eqref{PDO equation for utilde-k} by $\phi^j_k(\xi)$. 
On the one hand, reminding \eqref{def_qtilde_w} and using the $L^\infty-L^\infty$ continuity of operator $\Gamma^{w,k}$ (see proposition \ref{Prop: Continuity_Lp_wave}), together with the classical Sobolev injection, the fact that 
\begin{equation} \label{norm_O(mu)}
\left\|\oph(\phi^j_k(\xi))\right\|_{\Lcal(L^2)}=O(h^{-\mu}),
\end{equation}
with $\mu=\sigma\rho$ if $\rho\ge 0$, 0 otherwise, estimates \eqref{norm_L2_qtilde_w},\eqref{norm_(hD)_qtilde} and \eqref{est:a-priori_vt}, we find that
\begin{equation} \label{nl_term_1}
\begin{split}
&\left\|\Gamma^{w,k}\oph(\Sigma_j(\xi)(1-\chi_0)(h^{-1}\xi)\varphi(2^{-k}\xi)\chi_0(h^\sigma\xi))\left[t q_w(t,tx)\right]\right\|_{L^\infty}\\
&\lesssim h^{-\beta}\|\widetilde{q}_w(t,\cdot)\|_{L^2}+ h^{-1-\beta}\|\oph(\varphi(2^{-k}\xi))(hD_x)\widetilde{q}_w(t,\cdot)\|_{L^2}\le CB\varepsilon h^{1-\beta'}.
\end{split}
\end{equation}
On the other hand, using proposition \ref{Prop : continuity of Op(gamma1):X to L2}, estimate \eqref{norm_O(mu)}, the fact that the commutator between $\oph(\phi^j_k(\xi))$ and $\Omega_h$ is also continuous on $L^2$ with norm $O(h^{-\mu})$, equality $\|tw(t,t\cdot)\|_{L^2}=\|w(t,\cdot)\|_{L^2}$, and \eqref{est_cw_L2}, \eqref{est_Omega_cw} (in which we choose $s>0$ large enough to have, say, $N(s)\ge 2$), we deduce that there is a $\chi\in C^\infty_0(\mathbb{R}^2)$ such that
\begin{equation} \label{nl_term_1.5}
\begin{split}
&\|\Gamma^{w,k}\oph\big(\phi^j_k(\xi)\big)(h^{-1}c_w(t,tx))\|_{L^\infty(dx)}\\
& \lesssim t^{\frac{1}{2}+\beta}\left\|\chi(t^{-\sigma}D_x)(\vnf-v_{-})(t,\cdot)\right\|_{L^2}\left(\|V(t,\cdot)\|_{H^{2,\infty}}+\|\vnf(t,\cdot)\|_{H^{1,\infty}}\right) \\
& + t^{-\frac{3}{2}+\beta}\left\| (\vnf - v_{-})(t,\cdot)\right\|_{H^1}\left(\|V(t,\cdot)\|_{H^s}+\|\vnf(t,\cdot)\|_{H^s}\right) \\
& + t^{\frac{1}{2}+\beta}\left\|\chi(t^{-\sigma}D_x)\Omega (\vnf - v_{-})(t,\cdot)\right\|_{L^2} \left(\|V(t,\cdot)\|_{H^{2,\infty}}+\|\vnf(t,\cdot)\|_{H^{1,\infty}}\right)\\
& + t^{-\frac{3}{2}+\beta}\left\|\Omega (\vnf - v_{-})(t,\cdot)\right\|_{L^2} \left(\|V(t,\cdot)\|_{H^s}+\|\vnf(t,\cdot)\|_{H^s}\right)\\
& +  t^{\frac{1}{2}+\beta} \left\|(\vnf - v_{-})(t,\cdot)\right\|_{H^{1,\infty}}\times \sum_{\mu=0}^1\left(\|\Omega^\mu V(t,\cdot)\|_{H^1}+\|\Omega ^\mu\vnf(t,\cdot)\|_{L^2}\right).
\end{split}
\end{equation}
Also, from \eqref{est L2 rNF}, \eqref{est phi(D) Omega rNF-new} we get that for every $\theta\in ]0,1[$
\begin{equation}\label{nl_term_2}
\begin{split}
&\|\Gamma^{w,k}\oph\big(\phi^j_k(\xi)\big)(h^{-1}r^{NF}_w(t,tx))\|_{L^\infty} \lesssim  t^{\frac{1}{2}+\beta}\|V(t,\cdot)\|^2_{H^{13,\infty}}\|U(t,\cdot)\|_{H^1} \\
&+t^{\frac{1}{2}+\beta}\Big[ \|V(t,\cdot)\|^{1-\theta}_{H^{15,\infty}} \|V(t,\cdot)\|^\theta_{H^{17}}\left(\| U(t,\cdot)\|_{H^{1,\infty}}+ \|\mathrm{R}_1U(t,\cdot)\|_{H^{1,\infty}}\right)\\
&\hspace{20pt}+ \|V(t,\cdot)\|_{L^\infty}\left(\| U(t,\cdot)\|^{1-\theta}_{H^{16,\infty}}+ \|\mathrm{R}_1U(t,\cdot)\|^{1-\theta}_{H^{16,\infty}}\right)\|U(t,\cdot)\|^\theta_{H^{18}}\Big]\|\Omega V(t,\cdot)\|_{L^2}\\
&+ t^{\frac{1}{2}+\beta} \Big[\|V(t,\cdot)\|_{H^{1,\infty}}\left(\|U(t,\cdot)\|_{H^1} + \|\Omega U(t,\cdot)\|_{H^1}\right)\\
&\hspace{20pt} + \left(\|U(t,\cdot)\|_{H^{2,\infty}}+\| \mathrm{R}_1U(t,\cdot)\|_{H^{2,\infty}}\right)\|\Omega V(t,\cdot)\|_{L^2}\Big] \|V(t,\cdot)\|_{H^{17,\infty}}.
\end{split}
\end{equation}
Therefore, using \eqref{est_vNF-v-} with $s=1$, \eqref{est_chi_vnf-v}, \eqref{est: bootstrap argument a-priori est}, and choosing $\theta\ll 1$ sufficiently small, we derive that $h^{-1}\Gamma^{w,k}\oph\big(\phi^j_k(\xi)\big)(c_w(t,tx)+r^{NF}_w(t,tx))$ is a remainder $F_w(t,x)$ satisfying \eqref{Linfty_norm_Fw}.
Since function $(\partial\chi_0)(h^{-1}\xi)$ is localized for frequencies of size $h$, its product with $\psi(2^{-k}\xi)$ is non-zero only for values of $k\in\mathbb{Z}$ such that $2^k\sim h$.
In that case, by commutating $\Gamma^{w,k}$ with $\oph\big((\partial\chi_0)(h^{-1}\xi)\cdot(h^{-1}\xi)\psi(2^{-k}\xi)\big)$ and using the classical Sobolev injection, together with proposition \ref{Prop : continuity Op(gamma) L2 to L2}, we find that
\begin{equation}\label{nl_term_3}
\left\|i h \,\Gamma^{w,k}\oph\big((\partial\chi_0)(h^{-1}\xi)\cdot(h^{-1}\xi)\psi(2^{-k}\xi)\big)\widetilde{u}(t,\cdot) \right\|_{L^\infty}\lesssim h \|\widetilde{u}(t,\cdot)\|_{L^2}.
\end{equation}
Since $(\partial\chi_0)(h^\sigma\xi)$ is, instead, localized for frequencies larger than $h^{-\sigma}$, by applying the semi-classical Sobolev injection and lemma \ref{Lem : new estimate 1-chi} we find that
\begin{equation}\label{nl_term_4}
\left\| i\sigma h \, \Gamma^{w,k}\oph\big(\psi(2^{-k}\xi)(\partial\chi_0)(h^\sigma\xi))\cdot(h^\sigma\xi)\big)\widetilde{u}(t,\cdot)\right\|_{L^\infty}\lesssim h^N \|\widetilde{u}(t,\cdot)\|_{H^s_h},
\end{equation}
with $N=N(s)>1$ as long as $s>0$ is sufficiently large.
By lemma \ref{Lem: from energy to norms in sc coordinates-WAVE} we obtain that also the fifth and sixth addend in the right hand side of \eqref{PDO equation for utilde-k} are remainders $F_w(t,x)$.

Finally, after lemma \ref{Lemma : development of linear part} 
\begin{multline*}
-\oph\big(\theta(x)(x\cdot\xi - |\xi|)\widetilde{\varphi}(2^{-k}\xi)\big)\widetilde{u}^{\Sigma_j,k}_{\Lambda_w} = \frac{1}{2}\oph\big(\theta(x)(1-|x|^2)x\cdot\xi \widetilde{\varphi}(2^{-k}\xi)\big)\widetilde{u}^{\Sigma_j,k}_{\Lambda_w} \\
+ \oph(\theta(x)e(x,\xi)\widetilde{\varphi}(2^{-k}\xi))\widetilde{u}^{\Sigma_j,k}_{\Lambda_w}
\end{multline*}
with $e(x,\xi)$ given by \eqref{def of e(x,xi)}, and latter term in the above right hand side satisfies \eqref{est L2 Linfty e(x,xi)}. Using symbolic calculus of proposition \ref{Prop: a sharp b} until order $N\in\mathbb{N}$ we find that \small
\begin{multline*}
\frac{1}{2}\oph\big(\theta(x)(1-|x|^2)x\cdot\xi \widetilde{\varphi}(2^{-k}\xi)\big)\widetilde{u}^{\Sigma_j,k}_{\Lambda_w} 
= \theta(x)\Big[\frac{1}{2}(1-|x|^2)x\cdot(hD_x) + \frac{h}{2i}(1 -2|x|^2)\Big]\oph(\widetilde{\varphi}(2^{-k}\xi))\widetilde{u}^{\Sigma_j,k}_{\Lambda_w} \\
+ \frac{h}{4i}(\partial\theta)(x)\cdot x(1-|x|^2)\oph(\widetilde{\varphi}(2^{-k}\xi))\widetilde{u}^{\Sigma_j,k}_{\Lambda_w} 
+ {\sum}'h\theta_1(x) \oph(\widetilde{\varphi}_1(2^{-k}\xi))\widetilde{u}^{\Sigma_j,k}_{\Lambda_w} + \oph(r(_Nx,\xi))\widetilde{u}^{\Sigma_j,k}_{\Lambda_w},
\end{multline*}\normalsize
with ${\sum}'$ being a concise notation to indicate a linear combination, $\partial\theta (x)$ supported for $|x|>D_1$, $\theta_1\in C^\infty_0(\mathbb{R}^2)$, $\widetilde{\varphi}_1\in C^\infty_0(\mathbb{R}^2\setminus\{0\})$ coming out from the derivatives of $\widetilde{\varphi}$, and $r_N(x,\xi)$ integral remainder of the form
\begin{equation*}
\frac{h^N}{(\pi h)^2}\int e^{\frac{2i}{h}\eta\cdot z}\int_0^1 \theta_N(x+tz)(1-t)^{N-1}dt\ \widetilde{\varphi}_N(2^{-k}(\xi+\eta)) dzd\eta,
\end{equation*}
for some other $\theta_N\in C^\infty_0(\mathbb{R}^2)$, $ \widetilde{\varphi}_N\in C^\infty_0(\mathbb{R}^2\setminus\{0\})$, verifying that 
\begin{equation}\label{norm_Op(r)}
\|\oph(r_N(x,\xi))\|_{\mathcal{L}(L^2;L^\infty)}=O(h)
\end{equation}
if $N$ is taken sufficiently large. Therefore, from proposition \ref{Prop : continuity Op(gamma) L2 to L2}, \eqref{norm_O(mu)} and \eqref{est:utilde-Hs}
\begin{equation*}
\left\| \oph(r(x,\xi))\widetilde{u}^{\Sigma_j,k}_{\Lambda_w}(t,\cdot)\right\|_{L^\infty}\lesssim h^{1-\beta}\|\widetilde{u}(t,\cdot)\|_{L^2}\le CB\varepsilon h^{1-\beta'}.
\end{equation*}
Moreover, since $\widetilde{\varphi}\equiv 1$ on the support of $\varphi$ (which defines $\widetilde{u}^{\Sigma_j,k}$), by commutating $\oph(\widetilde{\varphi}(2^{-k}\xi))$ with $\Gamma^{w,k}$ and using remark \ref{Remark:symbols_with_null_support_intersection} we find that, for any $N\in\mathbb{N}$ as large as we want,
\[\oph(\widetilde{\varphi}(2^{-k}\xi))\widetilde{u}^{\Sigma_j,k}_{\Lambda_w} = \widetilde{u}^{\Sigma_j,k}_{\Lambda_w} + O_{L^\infty}(h^N\|\widetilde{u}\|_{L^2}).\]
Also, since $\widetilde{\varphi}_1$ is obtained from the derivatives of $\widetilde{\varphi}$ and vanishes on the support of $\varphi$, 
\[\theta_1(x)\oph(\widetilde{\varphi}_1(2^{-k}\xi))\widetilde{u}^{\Sigma_j,k}_{\Lambda_w} =O_{L^\infty}(h^N\|\widetilde{u}\|_{L^2}). \]
Therefore, again from \eqref{est:utilde-Hs} we deduce that
\begin{multline*}
-\oph\big(\theta(x)(x\cdot\xi - |\xi|)\widetilde{\varphi}(2^{-k}\xi)\big)\widetilde{u}^{\Sigma_j,k}_{\Lambda_w}=
\theta(x)\Big[\frac{1}{2}(1-|x|^2)x\cdot(hD_x) + \frac{h}{2i}(1 -2|x|^2)\Big]\widetilde{u}^{\Sigma_j,k}_{\Lambda_w} 
\\
+ \frac{h}{4i}(\partial\theta)(x)\cdot x(1-|x|^2)\widetilde{u}^{\Sigma_j,k}_{\Lambda_w} + \oph\big(\theta(x)\widetilde{\varphi}(2^{-k}\xi)e(x,\xi)\big)\widetilde{u}^{\Sigma_j,k}_{\Lambda_w}
 +  O_{L^\infty}(h^{1-\beta'}),
\end{multline*}
which implies, summed up with estimates from \eqref{nl_term_1} to \eqref{nl_term_4}, that $\widetilde{u}^{\Sigma_j,k}_{\Lambda_w}$ is solution to 
\begin{multline*}
\left[D_t + \theta(x)\frac{1}{2}(1-|x|^2)x\cdot (hD_x) + \theta(x)\frac{h}{2i}(1-2|x|^2)\right]\widetilde{u}^{\Sigma_j,k}_{\Lambda_w}(t,x) = F^k_w(t,x)  \\
 + \Big[(1-\theta)(x)\oph((x\cdot\xi - |\xi|)\widetilde{\varphi}(2^{-k}\xi))+
 \widetilde{\theta}(x)\oph(\widetilde{\varphi}_1(2^{-k}\xi))- \frac{h}{4i}(\partial\theta)(x)\cdot x(1-|x|^2)\Big]\widetilde{u}^{\Sigma_j,k}_{\Lambda_w}(t,x)\,,
\end{multline*}
where $F^k_w(t,x)$ satisfies \eqref{Linfty_norm_Fw}.
Choosing $D_1=D$, we obtain that $\widetilde{u}^{\Sigma_j}_{\Lambda_{w}}$ is solution to \eqref{eq: transport equation for uSigmaLambda} in cylinder $\mathcal{C}^T_D$, with $F_w(t,x):=\sum_k F^k_w(t,x)$ (this sum being finite and restricted to indices $k\in K$) satisfying the same $L^\infty$ estimate as $F^k_w$, up to an additional factor $h^{-\sigma}$.
\endproof
\end{prop}

\section{Analysis of the transport equation and end of the proof}

In previous section (see proposition \ref{Prop:propagation_unif_est_V}) we firstly showed how to propagate a-priori uniform estimate \eqref{est: boostrap vpm} on the Klein-Gordon component $v_{-}$, in the sense of deducing \eqref{est:bootstrap enhanced vpm} from estimates \eqref{est: bootstrap argument a-priori est}.
We then passed to the study of the wave equation and proved that, if $(u_{-},v_{-})$ is solution to \eqref{wave-KG for u- v-} in some interval $[1,T]$, function $\widetilde{u}^{\Sigma_j}_{\Lambda_w}$ defined in \eqref{def_u_Sigma_Lw} is solution to transport equation \eqref{eq: transport equation for uSigmaLambda} in truncated cylinder $\mathcal{C}^T_D:=\{(t,x):1\le t\le T, |x|\le D\}$, for any $D>0$.
The aim of this section is to study such a transport equation in order to deduce some information on the uniform norm of its solutions. This will allow us to finally propagate a-priori estimate \eqref{est: bootstrap upm} on the wave component $u_{-}$ and to close the bootstrap argument.
A short proof of main theorem \ref{Thm: Main theorem} is given at the end of this section.

\subsection{The inhomogeneous transport equation}
The aim of this subsection is to study the behaviour of a solution $w$ to the following transport equation
\begin{equation}
\left[D_t + \frac{1}{2}(1-|x|^2)x\cdot (hD_x) - \frac{i}{2t}(1-2|x|^2)\right]w = f\,, 
\end{equation}
in a cylinder $\mathcal{C} = \{(t,x) : t\ge 1, |x|\le D\}$ for a large constant $D\gg 1$, where the inhomogeneous term $f$ is a $O_{L^\infty}(\varepsilon t^{-1+\beta})$, for some $\varepsilon>0$ small and $0\le \beta<1/2$.
We distinguish in $\mathcal{C}$ two subregions:
\begin{equation*}
I_1 := \Big\{(t,x) : t\ge 1, |x|< \Big(\frac{t}{t-1}\Big)^{\frac{1}{2}}, |x|\le D\Big\}\,, \qquad I_2: = \Big\{(t,x) : t> 1, \Big(\frac{t}{t-1}\Big)^{\frac{1}{2}}\le |x|\le D\Big\},
\end{equation*}
and denote by $I_{1,t}, I_{2,t}$ their sections at a fixed time $t\ge 1$,
\begin{equation*}
I_{1,t} := \Big\{x : |x|< \Big(\frac{t}{t-1}\Big)^{\frac{1}{2}}, |x|\le D\Big\}\,, \qquad I_{2,t}: = \Big\{x : \Big(\frac{t}{t-1}\Big)^{\frac{1}{2}}\le |x|\le D\Big\}.
\end{equation*}


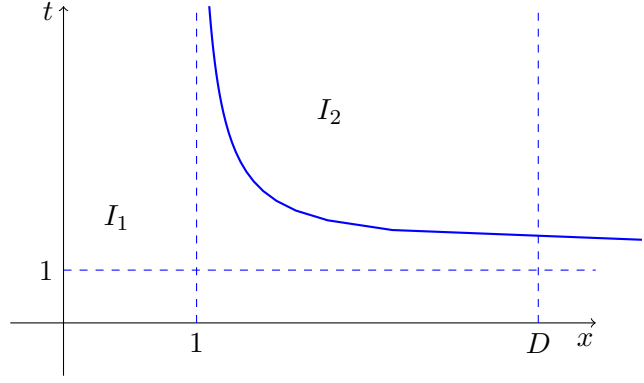
\begin{figure}
\begin{center}
\begin{tikzpicture}[scale=0.7]
\draw[->] (0,-1) -- (0,6);
\draw[->] (-1,0) -- (10,0);
\draw[blue,thick][scale=1.5] [domain=1.05:4] plot({((2.5)*\x/(\x-1))^(1/2)},\x);
\draw[blue,dashed] [domain=0:10]plot(\x, 1);
\draw[blue,dashed] [domain=0:6]plot(8.92, \x);
\draw[blue,dashed] [domain=0:6]plot(2.5, \x);
\node[below] at (2.5,0) {$1$};

\node[ left] at (0,1) {1};
\node[below] at (8.92,0) {$D$};
\node at (5,4) {$I_2$};
\node at (1, 2) {$I_1$};

\node[below] at (9.8,0) {$x$};
\node[left] at (0,5.9) {$t$};
\end{tikzpicture}
\end{center}
\caption{Regions $I_1$ and $I_2$ in space dimension 1} 
\end{figure}

The result we prove is the following.

\begin{prop} \label{Prop : estimate of solution of Tr in a cylinder}
Let $\varepsilon>0$ be small and $w$ be the solution to the following Cauchy problem
\begin{equation}  \label{Tr}
\begin{cases}
& \left[D_t + \frac{1}{2}(1-|x|^2)x\cdot (hD_x) - \frac{i}{2t}(1-2|x|^2)\right]w = f\,, \\
& w(1,x) = \varepsilon w_0(x)\,,
\end{cases}
\end{equation}
with $f =O_{L^\infty}(\varepsilon t^{-1+\beta})$, for some fixed $0\le \beta<1/2$. 
Let us suppose that $|w_0(x)|\lesssim \langle x \rangle^{-2}$ and that $|w(t,x)|\lesssim\varepsilon t^{\beta'}$ for some $\beta'>0$ whenever $|x|>D\gg 1$. Therefore,
\begin{equation} \label{est of w in the cylinder}
|w(t,x)|\lesssim \varepsilon \|w_0\|_{L^\infty} t^{\beta''} (1+|x|)^{-\frac{1}{2}}(t^{-1} + |1-|x||)^{-\frac{1}{2}+\beta''}\,,
\end{equation}
for every $(t,x)\in \mathcal{C}_D=\{(t,x) | t\ge 1, |x|\le D\}$,
with $\beta'' = \max\{\beta, \beta'\}$.
\end{prop}
We observe that, if $W(t,x)=t^{-1}w(t,t^{-1}x)$, the above inequality implies that
\begin{equation*}
|W(t,x)|\lesssim \varepsilon\|w_0\|_{L^\infty}(t+|x|)^{-\frac{1}{2}}(1+|t-|x||)^{-\frac{1}{2}+\beta''},
\end{equation*}
showing that the uniform norm of $W(t,\cdot)$ decays in time at a rate $t^{-1/2}$, enhanced to $t^{-1+\beta''}$ out of the light cone $t=|x|$.

In order to prove the result of proposition \ref{Prop : estimate of solution of Tr in a cylinder} we fix $T\ge 1$, $x\in B_D(0)$, and look for the characteristic curve of \eqref{Tr} with initial point $(T,x)$, i.e. map $t\mapsto X(t;T,x)$ solution of
\begin{equation} \label{eq of characteristics in I1}
\begin{cases}
\frac{d}{dt}X(t;T,x) = \frac{1}{2t}\big(1 - |X(t;T,x)|^2\big) X(t;T,x) \\
X(T;T,x) = x
\end{cases}\qquad t\ge T.
\end{equation}

\begin{lem} \label{Lem: characteristic equ in I1}
Solution $X(t;T,x)$ to \eqref{eq of characteristics in I1} writes explicitly as
\begin{equation} \label{expression of X(t;1,x) in I1}
X(t;T,x) = \frac{\sqrt{t}x}{(T - (T-t)|x|^2)^{\frac{1}{2}}}
\end{equation}
and it is well defined for all $t>T(1-|x|^{-2})$.
Moreover, for any fixed $t>T$, map $x\in\mathbb{R}^2 \mapsto X(t;T,x) \in \Big\{|x|< \big(\frac{t}{t-T}\big)^{\frac{1}{2}}\Big\}$ is a diffeomorphism of inverse $Y(t,y) = \frac{\sqrt{T} y}{(t + (T-t)|y|^2)^\frac{1}{2}}$. 
\proof
Multiplying equation \eqref{eq of characteristics in I1} by $2X(t;T,x)$ we deduce that $|X(t;T,x)|^2$ satisfies the equation
\begin{equation*}
\frac{d}{dt}|X(t;T,x)|^2 = \frac{1}{t}\big(1-|X(t;T,x)|^2\big)|X(t;T,x)|^2\,,
\end{equation*}
from which follows that $1- |X(t;T,x)|^2 = \frac{T(1-|x|^2)}{T- (T-t)|x|^2}$. Injecting this result in \eqref{eq of characteristics in I1} and integrating in time, we obtain expression \eqref{expression of X(t;1,x) in I1} and observe that the obtained map is well defined for all $t>T(1-|x|^{-2})$.

In order to prove the second part of the statement, we fix $t>T$, $y \in \Big\{|x|\le\big(\frac{t}{t-T}\big)^{\frac{1}{2}}\Big\}$ and look for $Y(t,y)$ such that $X(t;T, Y(t,y)) = y$. In other words,
$$y = \frac{\sqrt{t}Y(t,y)}{(T - (T-t)|Y(t,y)|^2)^{\frac{1}{2}}},$$ which implies that $Y(t,y) = \frac{\sqrt{T} y}{(t + (T-t)|y|^2)^\frac{1}{2}}$.
This map is well defined as long as $|y|<\big(\frac{t}{t-T}\big)^\frac{1}{2}$.
\endproof
\end{lem}

Along the characteristic curve $X(t;T,x)$ function $w$ satisfies
\begin{equation*}
\begin{split}
\frac{d}{dt}w\big(t, X(t;T,x)\big) & = -\frac{1}{2t}\big(1 - 2|X(t;T,x)|^2\big)\ w\big(t, X(t;T,x)\big) + i f\big(t, X(t;T,x)\big) \\
& = - \frac{1}{2t}\frac{T- T|x|^2 - t|x|^2}{T-(T-t)|x|^2}\ w\big(t, X(t;T,x)\big) + if\big(t, X(t;T,x)\big)
\end{split}
\end{equation*}
and hence
\begin{multline} \label{dt exp w(t, X(t,T,x))}
\frac{d}{dt}\left[\left(\exp \int_T^t \frac{1}{2\tau}\frac{T- T|x|^2 - \tau |x|^2}{T-(T-\tau)|x|^2}d\tau\right)w\big(t, X(t;T,x)\big)\right] \\ =i \left(\exp \int_T^t \frac{1}{2\tau}\frac{T- T|x|^2 - \tau |x|^2}{T-(T-\tau)|x|^2} d\tau\right)f\big(t, X(t;T,x)\big)\,.
\end{multline}

\begin{lem}
\begin{equation} \label{calculation exp of integral}
\exp\int_T^t \frac{1}{2\tau}\frac{T- T|x|^2 - \tau |x|^2}{T-(T-\tau)|x|^2}\, d\tau = \Big(\frac{t}{T}\Big)^\frac{1}{2}\Big(\frac{T - T|x|^2 + t|x|^2}{T}\Big)^{-1}.
\end{equation}
\proof
The result follows writing
$$\frac{1}{2\tau}\frac{T- T|x|^2 - \tau |x|^2}{T-(T-\tau)|x|^2} = \frac{1}{2\tau} - \frac{|x|^2}{T -T|x|^2 + \tau |x|^2}\,,$$
taking the integral over $\tau \in [T,t]$ and then passing to its exponential.
\endproof
\end{lem}

Let us first study the behaviour of $w$, solution to \eqref{Tr}, in region $I_1$. 
We fix $T=1$ and, integrating equality \eqref{dt exp w(t, X(t,T,x))} over $[1,t]$, we find that
\begin{multline} \label{expression of w along X}
\left( \exp \int_1^t \frac{1}{2\tau}\frac{1- |x|^2 - \tau |x|^2}{1-(1-\tau)|x|^2} d\tau\right) w\big(t, X(t;1,x)\big) \\ = w(1,x) + i \int_1^t \left(\exp \int_1^s \frac{1}{2\tau}\frac{1- |x|^2 - s |x|^2}{1-(1-s)|x|^2}ds \right)f\big(s, X(s;1 ,x)\big)ds.
\end{multline}
Using \eqref{calculation exp of integral} and the fact that $f = O_{L^\infty}(\varepsilon t^{-1+\beta})$, we then obtain that
\begin{multline} \label{inequality for w along charact}
\big|w(t, X(t;1,x))\big| \le t^{-\frac{1}{2}}(1-|x|^2 + t|x|^2)|w(1,x)| \\
+ C\varepsilon t^{-\frac{1}{2}}(1-|x|^2 + t|x|^2) \int_1^t \frac{ds}{(1-|x|^2+s|x|^2)s^{\frac{1}{2}-\beta}}\,,
\end{multline}
for some positive constant $C$.

\begin{lem} 
For any fixed $0\le \beta<1/2$
\begin{equation} \label{integral inequality}
 \int_1^t \frac{ds}{(1-|x|^2+s|x|^2)s^{\frac{1}{2}-\beta}} \lesssim \frac{t^{\frac{1}{2}+\beta}}{(1 + \sqrt{t}|x|)^{1+2\beta}}(1+|x|)^{-1+2\beta + \beta'}\,,
\end{equation}
for all $t\ge 1$ and $\beta'>0$ as small as we want.
\proof
For $\sqrt{t} |x|\le 1$, we have that
\begin{equation*}
\int_1^t \frac{ds}{(1-|x|^2+s|x|^2)s^{\frac{1}{2}-\beta}} \lesssim t^{\frac{1}{2}+\beta} \lesssim \frac{t^{\frac{1}{2}+\beta}}{(1 + \sqrt{t}|x|)^{1+2\beta}}(1+|x|)^{-1+2\beta + \beta'},
\end{equation*}
for any $\beta'\ge 0$. Suppose then that $\sqrt{t}|x| > 1$.
For $t\le 2$
\begin{equation*}
\begin{gathered}
\displaystyle\int_1^t \frac{ds}{(1-|x|^2+s|x|^2)s^{\frac{1}{2}-\beta}} \lesssim (1+|x|)^{-2}\log(1+ |x|^2) \\
\text{and}\quad |x|^{-2}\log(1+ |x|^2) \frac{(1+\sqrt{t}|x|)^{1+2\beta}}{t^{\frac{1}{2}+\beta}}\lesssim (1+|x|)^{-1+2\beta}\log(1+|x|^2),
\end{gathered}
\end{equation*}
which immediately implies inequality \eqref{integral inequality}.
For $t\ge 2$ 
\[ \int_1^t \frac{ds}{(1-|x|^2+s|x|^2)s^{\frac{1}{2}-\beta}} =  \int_1^2 \frac{ds}{(1-|x|^2+s|x|^2)s^{\frac{1}{2}-\beta}} +  \int_2^t \frac{ds}{(1-|x|^2+s|x|^2)s^{\frac{1}{2}-\beta}}\, ,\]
where the first integral is bounded from the right hand side of \eqref{integral inequality}.
The second one is less or equal than $\int_1^{t-1}\frac{ds}{(1+s|x|^2)s^{\frac{1}{2}-\beta}}$, so
for $|x|\ge 1$ it follows that
\begin{equation*}
\int_1^{t-1}\frac{ds}{(1+s|x|^2)s^{\frac{1}{2}-\beta}} \le |x|^{-2}\int_1^{t-1}\frac{ds}{s^{\frac{3}{2}-\beta}}\lesssim (1+|x|)^{-2}.
\end{equation*}
Since $\frac{(1+ \sqrt{t}|x|)^{1+2\beta}}{t^{\frac{1}{2}+\beta}}\le (1+|x|)^{1+2\beta}$, from the above inequality we deduce the right bound of the statement.
For $|x|< 1$, a change of variables gives that
\small
\begin{equation*}
\int_1^{t-1}\frac{ds}{(1+s|x|^2)s^{\frac{1}{2}-\beta}} = |x|^{-1-2\beta}\int_{|x|^2}^{(t-1)|x|^2}\frac{ds}{(1+s)s^{\frac{1}{2}-\beta}}\lesssim |x|^{-1-2\beta}\frac{(t|x|^2)^{\frac{1}{2}+\beta}}{(1+t|x|^2)^{\frac{1}{2}+\beta}}\le \frac{t^{\frac{1}{2}+\beta}}{(1+t|x|^2)^{\frac{1}{2}+\beta}}.
\end{equation*} \normalsize
\endproof
\end{lem}

If initial condition $w_0(x)$ is sufficiently decaying in space, e.g. $|w_0(x)|\lesssim \langle x\rangle^{-2}$, we deduce from inequalities \eqref{inequality for w along charact} and \eqref{integral inequality} the following bound for $w$ along the characteristic curve $X(t;1,x)$:
\begin{equation} \label{est for w along characteristic in I1}
\big|w(t, X(t;1,x))\big| \lesssim\varepsilon \|w_0\|_{L^\infty} t^\beta (1+\sqrt{t}|x|)^{1-2\beta}(1+|x|)^{-1+2\beta +\beta'}\,,
\end{equation}
for any $\beta'>0$ as small as we want.

\begin{figure}
\begin{center}
\begin{tikzpicture}[scale=0.7]
\draw[->] (0,-1) -- (0,7);
\draw[->] (-1,0) -- (12,0);

\node[below] at (11.8,0) {$x$};
\node[left] at (0,6.9) {$t$};

\draw[dashed] [domain=0:12]plot(\x, 1);
\node[ left] at (0,1) {1};

\draw[blue, thick] [domain=1.5:7]plot({((2.5/2)*\x^(1/2))/((9/8)+(1/4)*\x)^(1/2)},\x);
\node[blue,thick] at (1.25, 1.5) {$\bullet$};
\node[left] at (1.25, 1.5) {\tiny{$(T_1,x_1)$}};
\draw[dashed][domain=0:7] plot(2.5, \x);

\draw[blue, thick] [domain=1.58:7]plot({((2.5)*5*\x^(1/2))/((3/2)-(75/2)+(25)*\x)^(1/2)},\x);
\node[blue,thick] at (8.35, 1.58) {$\bullet$};
\node[below,right] at (8.35, 1.58) {\tiny{$(T_2,x_2)$}};

\node[below] at (2.5,0) {$1$};

\draw[red,thick][scale=1.8] [domain=1.064:4] plot({((2.5)*\x/(\x-1))^(1/2)},\x);
\node[above] at (10,2) {\tiny{$x=(\frac{t}{t-1})^\frac{1}{2}$}};

\end{tikzpicture}
\end{center}
 \caption{Characteristic curves of initial point $(T_i,x_i)\in I_1$, $i=1,2$, in space dimension 1} 
\end{figure}
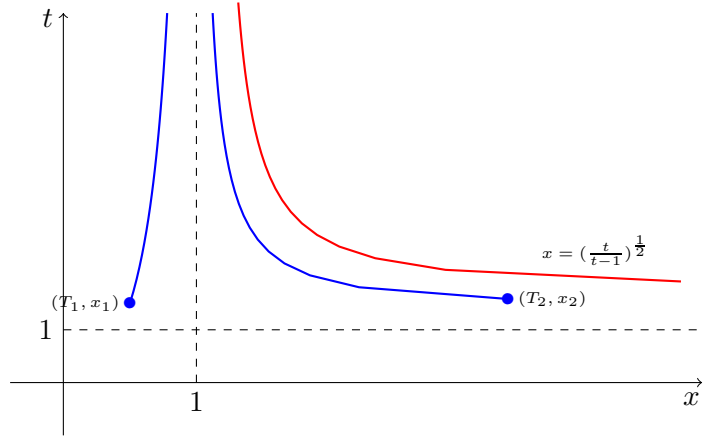

\begin{prop} \label{Prop: estimate for transport sol n I1}
Let $w$ be the solution to transport equation \eqref{Tr}, with $\|f(t,\cdot)\|_{L^\infty}\lesssim \varepsilon t^{-1+\beta}$ for some fixed $0\le \beta<1/2$, and initial condition $|w_0(x)|\lesssim \langle x \rangle^{-2}$, $\forall x\in\mathbb{R}^2$.
Then 
\begin{equation} \label{est of w in I1}
|w(t,x)|\lesssim\varepsilon t^\beta \big[t^{-1} + |1-|x||\big]^{-\frac{1}{2}+\beta}
\end{equation}
for every $(t,x)\in I_1=\{(t,x) : t\ge 1, |x|< \big(\frac{t}{t-1}\big)^{\frac{1}{2}}, |x|\le D\}$.
\proof
In lemma \ref{Lem: characteristic equ in I1} we proved that, for any fixed $t>T=1$, map $x\in\mathbb{R}^2\mapsto X(t;1,x)\in \big\{x : |x|<(\frac{t}{t-1})^\frac{1}{2}\big\}$ is a diffeomorphism with inverse $Y(t,y) = y(t + (1-t)|y|^2)^{-1/2}$.
From inequality \eqref{est for w along characteristic in I1} we hence deduce that, for any $y$ such that $|y|<\big(\frac{t}{t-1}\big)^\frac{1}{2}$,
\begin{equation*}
|w(t,y)|\lesssim\varepsilon t^\beta \big(1 + \sqrt{t}|Y(t,y)|\big)^{1-2\beta}\big(1 + |Y(t,y)|\big)^{-1+2\beta+\beta'}.
\end{equation*}
In particular, as $t(1-|y|^2)+|y|^2)\sim t|1-|y|^2| + |y|^2$ when $|y|<\big(\frac{t}{t-1}\big)^\frac{1}{2}$ and $t\ge t_0>1$, and $t|1-|y|^2| + |y|^2 \sim t|1-|y|| + |y|$ when $|y|\le D$, we find for those values of $(t,y)$ that
\begin{equation*}
|w(t,y)| \lesssim\varepsilon t^\beta \left(1 + \frac{\sqrt{t}|y|}{(t|1-|y|| + |y|)^{\frac{1}{2}}}\right)^{1-2\beta} \lesssim\varepsilon t^\beta\big[t^{-1} + |1-|y||\big]^{-\frac{1}{2}+\beta}\,,
\end{equation*}
simply using that $(1+|Y(t,y)|)^{-1+2\beta+\beta'}\le 1$. 
Moreover, for $t\rightarrow 1$ and $|y|\le D$,
\begin{equation*}
|w(t,y)|\lesssim \varepsilon \lesssim \varepsilon t^\beta \big[t^{-1} + |1-|y||\big]^{-\frac{1}{2}+\beta}.
\end{equation*}
\endproof
\end{prop}

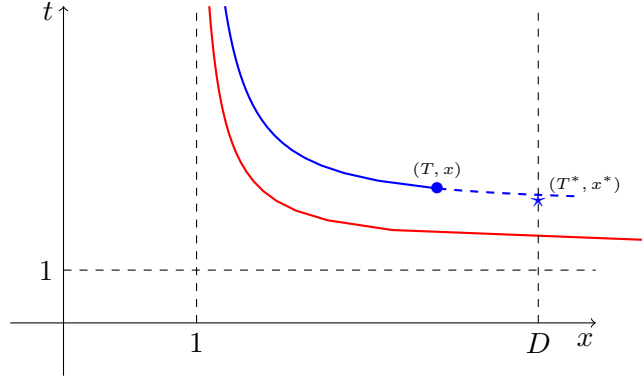
\begin{figure}
\begin{center}
\begin{tikzpicture}[scale=0.7]
\draw[->] (0,-1) -- (0,6);
\draw[->] (-1,0) -- (10,0);

\draw[red,thick][scale=1.5] [domain=1.05:4] plot({((2.5)*\x/(\x-1))^(1/2)},\x);
\draw[blue,thick][scale=1.5][domain=1.7:4]plot({1.60*(2*\x^(1/2))/(-6+4*\x)^(1/2)},\x);
\draw[blue, thick, dashed][scale=1.5][domain=1.6:1.7]plot({1.60*(2*\x^(1/2))/(-6+4*\x)^(1/2)},\x);

\draw[dashed] [domain=0:10]plot(\x, 1);
\draw[dashed] [domain=0:6]plot(8.92, \x);
\draw[dashed] [domain=0:6]plot(2.5, \x);
\node[below] at (2.5,0) {$1$};

\node[blue] at (7.02,2.55) {$\bullet$};
\node[above] at (7.02,2.55) {\tiny{$(T,x)$}};

\node[above right] at (8.92, 2.27) {\tiny{$(T^*, x^*)$}};
\node[blue] at (8.92, 2.32) {$\star$};
\node[ left] at (0,1) {1};
\node[below] at (8.92,0) {$D$};

\node[below] at (9.8,0) {$x$};
\node[left] at (0,5.9) {$t$};
\end{tikzpicture}
\end{center}
\caption{Characteristic curve of initial point $(T,x)\in I_2$} 
\end{figure}

\begin{prop} \label{Prop : est for transport solution in I2}
Let $\varepsilon>0$ be small and $w$ be the solution to transport equation \eqref{Tr}, with $\|f(t,\cdot)\|_{L^\infty}\lesssim \varepsilon t^{-1+\beta}$ for some fixed $0\le \beta<1/2$, and suppose that $|w(t,x)|\lesssim\varepsilon t^{\beta'}$ for some $\beta'>0$ whenever $|x|\ge D$. Then
\begin{equation*}
|w(t,x)| \lesssim\varepsilon t^{\beta''} (|x|^2 - 1)^{\beta''-\frac{1}{2}}\,,
\end{equation*}
for every $(t,x)\in I_2 = \{(t,x) : t>1, \big(\frac{t}{t-1}\big)^\frac{1}{2}\le |x| \le D\}$, where $\beta'' = \max\{\beta, \beta'\}$.
\proof
For a fixed $(T,x)\in I_2$ we look at $X(t; T,x)$, solution to \eqref{eq of characteristics in I1} and given by the explicit expression \eqref{expression of X(t;1,x) in I1}. We observe that there exists a time $T^*$, $1< T^*<T$, such that $X(t;T,x)$ hits the boundary $|y| =D$ at $t =T^*$. In other words, $t = T^*$ is the first time when $X(t;T,x)$ enters in the region $\{(t,x) : t\ge 1, |x|\le D\}$, to never leave it again for function $t\mapsto |X(t;T,x)|$ is strictly decreasing. A simple computation shows that
\begin{equation} \label{T* in terms of T}
T^* = \frac{D^2}{D^2 - 1}(1-|x|^{-2})T <T\,.
\end{equation}
Integrating expression \eqref{dt exp w(t, X(t,T,x))} over $[T^*, T] $ and using \eqref{calculation exp of integral}, we find that
\begin{multline} \label{w(T,x)}
w(T,x) = \Big(\frac{T^*}{T}\Big)^\frac{1}{2}\Big(\frac{T - T(1 -|x|^{-2})}{T^* - T(1-|x|^{-2})}\Big) w(T^*, X(T^*;T,x)) \\
+ i \int_{T^*}^T \Big(\frac{t}{T}\Big)^\frac{1}{2}\Big(\frac{T - T(1 -|x|^{-2})}{t - T(1-|x|^{-2})}\Big) f\big(t, X(t;T,x)\big) dt\,.
\end{multline}
From \eqref{T* in terms of T}
\[T^* - T(1-|x|^{-2}) = \frac{1}{D^2-1}(1-|x|^{-2})T \quad\text{and}\quad \frac{T^*}{T} = \frac{D^2}{D^2-1}(1-|x|^{-2})\]
so since $|w(t,x)|\lesssim \varepsilon t^{\beta'}$ whenever $|x|\ge D$, for some $\beta'>0$ by the hypothesis, we find that the first term in right hand side of \eqref{w(T,x)} is bounded by $C\varepsilon (|x|^2-1)^{-\frac{1}{2}}(T^*)^{\beta'}$, for some constant $C>0$.
Setting $c= \frac{1}{D^2-1}$, by the hypothesis on $f$ we derive that
\begin{equation*}
\begin{split}
 \Big|\int_{T^*}^T \Big(\frac{t}{T}\Big)^\frac{1}{2}\Big(\frac{T - T(1 -|x|^{-2})}{t - T(1-|x|^{-2})}\Big) f\big(t, X(t;T,x)\big) &  dt\Big|  \lesssim\varepsilon T^\frac{1}{2}\int_{T^*}^T \big(t - T (1-|x|^{-2})\big)^{-1} t^{-\frac{1}{2}+\beta} dt \\
& =\varepsilon T^\frac{1}{2}\int_{T^*}^T \big(t-T^* + c (1-|x|^{-2})T\big)^{-1} t^{-\frac{1}{2}+\beta} dt \\
& \le\varepsilon T^{\frac{1}{2}}\int_0^{T-T^*}\frac{dt}{\big(t + c (1-|x|^{-2})T\big) t^{\frac{1}{2}-\beta}} \\
& \lesssim\varepsilon T^\frac{1}{2}\big((1-|x|^{-2})T\big)^{\beta-\frac{1}{2}} =\varepsilon T^\beta(1-|x|^{-2})^{\beta - \frac{1}{2}}\,.
\end{split}
\end{equation*}
\endproof
\end{prop}

\subsection{Propagation of the uniform estimate on the wave component} \label{Subs: propagation of the unif est wave}

\begin{prop}[Propagation of the a-priori estimate on $U, RU$] \label{Prop: Propagation uniform estimate U,RU}
Let us fix $K_1>0$.
There exist two integers $n\gg\rho\gg 1$ sufficiently large,, two constants $A,B>1$ sufficiently large, some small $0<\delta\ll \delta_2\ll \delta_1\ll \delta_0$, and $\varepsilon_0\in ]0,1[$ sufficiently small, such that, for any $0<\varepsilon<\varepsilon_0$, if $(u,v)$ is solution to \eqref{wave KG system}-\eqref{initial data} in some interval $[1,T]$, for a fixed $T>1$, and $u_\pm, v_\pm$ defined in \eqref{def u+- v+-} satisfy a-priori estimates \eqref{est: bootstrap argument a-priori est}, for every $t\in [1,T]$, then it also verify \eqref{est:bootstrap enhanced upm} in the same interval $[1,T]$.
\proof
We warn the reader that, throughout the proof, $C,\beta, \beta'$ will denote some positive constants that may change line after line, such that $\beta\rightarrow 0$ as $\sigma\rightarrow0$ (resp. $\beta'\rightarrow 0$ as $\delta_1,\sigma\rightarrow 0$). We also remind that $h=1/t$.

In proposition \ref{Prop: NF on wave} we introduced function $u^{NF}$, defined from $u_{-}$ through \eqref{def uNF}, and showed that its $H^{\rho+1,\infty}$ norm (resp. the $H^{\rho+1,\infty}$ norm of $\mathrm{R}u^{NF}$) differs from that of $u_{-}$ (resp. of $\mathrm{R}u_{-}$) by a quantity satisfying \eqref{Hrho-infty norm uNF- u-} (resp. \eqref{Hrho-infty norm R(uNF-u-)}).
If $n$ is sufficiently large with respect to $\rho$ (at least $n\ge \rho+18$), a-priori estimates \eqref{est: boostrap vpm}, \eqref{est: bootstrap Enn} give that, for every $t\in [1,T]$,
\begin{multline} \label{Hrhoinfty_norm_u- in terms of uNF}
\|u_{-}(t,\cdot)\|_{H^{\rho+1,\infty}} + \sum_{j=1}^2\|\mathrm{R}_j u_{-}(t,\cdot)\|_{H^{\rho+1,\infty}}\\
\le \|u^{NF}(t,\cdot)\|_{H^{\rho+1,\infty}}+\sum_{j=1}^2 \|\mathrm{R}_j u^{NF}(t,\cdot)\|_{H^{\rho,\infty}} + 2AB\varepsilon^2 t^{-1+\frac{\delta}{2}}.
\end{multline}
We successively considered $\widetilde{u}(t,x):=t\widetilde{u}^{NF}(t,tx)$ and decomposed it as in \eqref{decomposition frequencies utilde}, with $\Sigma_j$ given by \eqref{Sigma_j}, showing that it satisfies \eqref{utilde small frequencies} (resp. \eqref{utilde large frequencies}) when restricted to small frequencies $|\xi|\lesssim t^{-1}$ (resp. large frequencies $|\xi|\gtrsim t^\sigma)$. 
We then focused on $\widetilde{u}^{\Sigma_j,k}$ defined in \eqref{def utilde-Sigma,k}, which is localized for frequencies supported in an annulus of size $2^k$ with $k\in K=\{k\in\mathbb{Z}: h\lesssim 2^k\lesssim h^{-\sigma}\}$, and further split it into the sum of functions $\widetilde{u}^{\Sigma_j,k}_{\Lambda_w}$, $\widetilde{u}^{\Sigma_j,k}_{\Lambda^c_w}$ (see \eqref{def_ukLambda-ukLambdac}). 
On the one hand, from inequality \eqref{Linfy_norm_utilde_Lambdac} and lemma \ref{Lem: from energy to norms in sc coordinates-WAVE} we have that there is a positive constant $C$ such that, for every $t\in [1,T]$,
\begin{equation*}
\|\widetilde{u}^{\Sigma_j,k}_{\Lambda^c_w}(t,\cdot)\|_{L^\infty}\le C\varepsilon t^{\beta'}.
\end{equation*}
On the other hand, we proved in proposition \ref{Prop: transport equation for uLambda} that, for any $D>0$ and any $(t,x)$ in truncated cylinder $\mathcal{C}^T_D=\{(t,x): 1\le t\le T, |x|\le D\}$, $\widetilde{u}^{\Sigma_j}_{\Lw}(t,x)$ defined in \eqref{def_u_Sigma_Lw} is solution to inhomogeneous transport equation \eqref{eq: transport equation for uSigmaLambda}, with inhomogeneous term $F_w(t,x)$ satisfying \eqref{Linfty_norm_Fw}.
Observe that, by definition \eqref{def Mj} of $\mathcal{M}$, symbolic calculus, and proposition \ref{Prop : Continuity on H^s}, we have that
\begin{equation*}
\|\widetilde{u}^{\Sigma_j}_{\Lw}(1,\cdot)\|_{L^2} + \|x \widetilde{u}^{\Sigma_j}_{\Lw}(1,\cdot)\|_{L^2} \lesssim \|\widetilde{u}(1,\cdot)\|_{L^2} + \|\oph(\chi(h^\sigma \xi))\mathcal{M}\widetilde{u}(1,\cdot)\|_{L^2} \le C\varepsilon,
\end{equation*}
which means that $\varepsilon^{-1}\langle x\rangle\widetilde{u}^{\Sigma_j}_{\Lw}(1,x)\in L^2$.
That hence implies that $|\widetilde{u}^{\Sigma_j}_{\Lw}(1,x)|\lesssim \varepsilon \langle x\rangle^{-2}$ for every $x\in\mathbb{R}^2$ (if not, we would have $\|\langle\cdot\rangle^{-1}\|_{L^2}\le \varepsilon^{-1} \|\langle \cdot\rangle \widetilde{u}^{\Sigma_j}_{\Lw}(1,\cdot)\|_{L^2}$).
Moreover, if $D\gg 1$ is sufficiently large, from lemma \ref{Lem: est utilde Sigma,k large x} below and \ref{Lem: from energy to norms in sc coordinates-WAVE} in appendix \ref{Appendix B} we deduce that
\begin{multline} \label{est_utilde_xBig}
|\mathds{1}_{|x|\ge D} \widetilde{u}^{\Sigma_j}_{\Lw}(t,x)|\le C\frac{\log{|x|}}{|x|} h^{-\beta}\big(\|\oph(\chi(h^\sigma\xi))\widetilde{u}(t,\cdot)\|_{L^2} + \|\oph(\chi(h^\sigma\xi))\mathcal{M}\widetilde{u}(t,\cdot)\|_{L^2}\big)\\
\le C\varepsilon \frac{\log |x|}{|x|} t^{\beta'}.
\end{multline}
Therefore, from proposition \ref{Prop : estimate of solution of Tr in a cylinder} we obtain that
\begin{equation*}
|\widetilde{u}^{\Sigma_j}_{\Lw}(t,x)|\lesssim C\varepsilon t^{\beta'} (1+|x|)^{-\frac{1}{2}}\big(t^{-1}+ |1-|x||\big)^{-\frac{1}{2}+\beta'}, \quad \forall (t,x)\in  \mathcal{C}^T_D.
\end{equation*}
Summing up, denoting by $\mathds{1}_{\mathcal{C}^T_D}$ the characteristic function of cylinder $\mathcal{C}^T_D$,
\begin{equation*}
|\widetilde{u}^\Sigma(t,x)|\le C\varepsilon \mathds{1}_{\mathcal{C}^T_D} t^{\beta'} (1+|x|)^{-\frac{1}{2}}\big(t^{-1}+ |1-|x||\big)^{-\frac{1}{2}+\beta'} + C\varepsilon t^{\beta'}, \quad \forall (t,x)\in [1,T]\times\mathbb{R}^2.
\end{equation*}
Returning back to function $u^{NF}$ via \eqref{def utilde vtilde}, this means that, for every $(t,x)\in [1,T]\times \mathbb{R}^2$,
\begin{multline} \label{Hrhoinfty_bound_uNF}
\left|\langle D_x\rangle^\rho u^{NF}(t,x)\right| + \sum_{j=1}^2 \left|\langle D_x\rangle^\rho R_j u^{NF}(t,x)\right|\\
\le C\varepsilon \mathds{1}_{\{|x|\le Dt \}} (t+|x|)^{-\frac{1}{2}}(1+ |t-|x||)^{-\frac{1}{2}+\beta'}+ C\varepsilon t^{-1+\beta'}.
\end{multline}
Finally, reminding definition \ref{def Sobolev spaces-NEW} $(iii)$ of space $H^{\rho,\infty}$, injecting the above inequality in \eqref{Hrhoinfty_norm_u- in terms of uNF}, and choosing $A>1$ sufficiently large such that $C<\frac{A}{3K_1}$, $\varepsilon_0>0$ sufficiently small so that $CB\varepsilon_0<(3K_1)^{-1}$, we deduce enhanced estimate \eqref{est:bootstrap enhanced upm}.
\endproof
\end{prop}

\begin{remark}
Beside the propagation of estimate \eqref{est: bootstrap upm}, by combining inequalities \eqref{Hrhoinfty_norm_u- in terms of uNF}, \eqref{Hrhoinfty_bound_uNF}, and \eqref{def u+- v+-}, we also deduce the following inequality
\begin{equation*}
| \partial_t u (t,x)| + |\nabla_x u(t,x)|\le C\varepsilon \mathds{1}_{\{|x|\le Dt \}} (t+|x|)^{-\frac{1}{2}}(1+ |t-|x||)^{-\frac{1}{2}+\beta'}+ C\varepsilon t^{-1+\beta'},
\end{equation*}
with $\beta'>0$ small as long as $\sigma, \delta_1$ are small,
which almost corresponds to the optimal decay in time and space enjoyed by the linear wave in space dimension two.
\end{remark}

\begin{lem}\label{Lem: est utilde Sigma,k large x}
Let $\chi\in C^\infty_0(\mathbb{R}^2)$ be equal to 1 in a neighbourhood of the origin and $\sigma>0$ be small.
Let also $\varphi\in C^\infty_0(\mathbb{R}^2\setminus\{0\})$. There exists a constant $C>0$ such that for every $h\in ]0,1[, R\gg 1$, and any function $w(t,x)$ with $w(t,\cdot), \oph(\chi(h^\sigma\xi))\mathcal{M}w(t,\cdot)\in L^2(\mathbb{R}^2)$,
\begin{equation} \label{ineq:utilde xbig}
\Big\|\varphi\Big(\frac{\cdot}{R}\Big)\oph(\chi(h^\sigma\xi))w(t,\cdot)\Big\|_{L^\infty}\le C R^{-1}(\log{R}+|\log h|) \sum_{|\gamma|=0}^1\|\oph(\chi(h^\sigma\xi))\mathcal{M}^\gamma w(t,\cdot)\|_{L^2}\,.
\end{equation}
\proof
Let us fix $R\gg 1$ and, for seek of compactness, denote $\oph(\chi(h^\sigma\xi))w$ by $w^\chi$.
For a new smooth cut-off function $\chi_1$ equal to 1 on the support of $\chi$, we have that
\begin{equation*}
\varphi\Big(\frac{x}{R}\Big) \oph(\chi(h^\sigma\xi)) w= \oph(\chi_1(h^\sigma\xi)) \Big[\varphi\Big(\frac{x}{R}\Big) w^\chi\Big] + \left[\varphi\Big(\frac{x}{R}\Big), \oph(\chi_1(h^\sigma\xi))\right]w^\chi,
\end{equation*}
where the symbol associated to above commutator is given by
\begin{equation*}
r_R(x,\xi) = -\frac{h^{1+\sigma} R^{-1}}{i(\pi h)^2}\int e^{\frac{2i}{h}\eta\cdot z}\left[\int_0^1 (\partial \varphi)\Big(\frac{x+tz}{R}\Big) dt\right] (\partial \chi_1)(h^\sigma(\xi+\eta)) dz d\eta,
\end{equation*}
as follows from \eqref{r_N 1} and integration in $dy,d\zeta$.
Since $(\partial \chi_1)(h^\sigma\xi)$ is supported for frequencies $|\xi|\le h^{-\sigma}$, and $R^{-1}, h^{1+\sigma}\le 1$, by making a change of coordinates $\eta/h\mapsto \eta$ and using that $e^{2i \eta\cdot z} = \big(\frac{1-2i\eta\cdot\partial_z}{1+4|\eta|^2}\big)\big(\frac{1-2i z\cdot\partial_\eta}{1+4|z|^2}\big)e^{2i \eta\cdot z}$, together with some integration by parts, one can check that 
\[\left\|\partial^\alpha_y \partial^\beta_\xi \big[r_R(\frac{x+y}{2}, h\xi)\big]\right\|_{L^2(d\xi)}\lesssim R^{-1}\]
for any $\alpha,\beta\in\mathbb{N}^2$, and hence obtain from lemma \ref{Lemma on inequalities for Op(A)} that
\begin{equation*}
\|\oph(r^k_R(x,\xi))w^\chi(t,\cdot)\|_{L^\infty}\lesssim R^{-1}\|w^\chi(t,\cdot)\|_{L^2}.
\end{equation*}
Successively, taking a Littlewood-Paley decomposition such that
\begin{equation*}
\chi_1(h^\sigma\xi) \equiv \left[\phi\Big(\frac{R}{h}\xi\Big) + \sum_{hR^{-1}\le 2^j\le h^{-\sigma}}(1-\phi)\Big(\frac{R}{h}\xi\Big)\psi(2^{-j}\xi)\right] \chi_1(h^\sigma\xi),
\end{equation*}
with $\phi\in C^\infty_0(\mathbb{R}^2)$, equal to 1 close to the origin and $\psi\in C^\infty_0(\mathbb{R}^2\setminus\{0\})$, we derive that
\begin{multline}\label{Linfty_wchi_split1}
\left\|\oph(\chi_1(h^\sigma\xi))\Big[ \varphi\Big(\frac{x}{R}\Big) w^\chi\Big](t,\cdot) \right\|_{L^\infty} \lesssim \left\| \oph\Big(\phi\Big(\frac{R}{h}\xi\Big)\chi_1(h^\sigma\xi)\Big)\Big[ \varphi\Big(\frac{x}{R}\Big) w^\chi\Big](t,\cdot) \right\|_{L^\infty}\\
 + \sum_{hR^{-1}\le 2^j\le h^{-\sigma}} \left\|\oph\Big((1-\phi)\Big(\frac{R}{h}\xi\Big)\psi(2^{-j}\xi)\chi_1(h^\sigma\xi)\Big)\Big[ \varphi\Big(\frac{x}{R}\Big) w^\chi \Big](t,\cdot)\right\|_{L^\infty},
\end{multline}
and immediately notice that
\begin{multline}\label{Linfty_wchi_phi(Rh-1xi)}
\left\| \oph\Big(\phi\Big(\frac{R}{h}\xi\Big)\chi_1(h^\sigma\xi)\Big)\Big[ \varphi\Big(\frac{x}{R}\Big) w^\chi\Big](t,\cdot) \right\|_{L^\infty}  \\
=\left\|\phi(RD_x) \oph(\chi_1(h^\sigma\xi))\Big[ \varphi\Big(\frac{x}{R}\Big) w^\chi\Big](t,\cdot) \right\|_{L^\infty} \lesssim R^{-1}\|w^\chi(t,\cdot)\|_{L^2},
\end{multline}
just by the classical Sobolev injection and the uniform continuity of $\oph(\chi_1(h^\sigma\xi))\varphi\Big(\frac{x}{R}\Big)$ on $L^2$.
Introducing operators $\Theta_R, \Theta^{-1}_R$, where $\Theta_R u(x) := u(Rx)$, $\Theta^{-1}_Ru (x) := u\big(\frac{x}{R}\big)$, we have the following equality
\begin{multline} \label{operators Theta R}
\oph\Big((1-\phi)\Big(\frac{R}{h}\xi\Big)\psi(2^{-j}\xi)\chi_1(h^\sigma\xi)\Big)\Big[\varphi\Big(\frac{x}{R}\Big)w^\chi\Big]\\
 = \Big[\Theta^{-1}_R Op^w_{h_{Rj}}\Big( (1-\phi)\Big(\frac{\xi}{h_{Rj}}\Big)\psi(\xi)\chi_1(h^\sigma 2^j\xi)\Big)\varphi(x) \Theta_R\Big]w^\chi
\end{multline}
with $h_{Rj}:=\frac{h}{R2^j}\le 1$, and by $h_{Rj}$-symbolic calculus (that is proposition \ref{Prop: a sharp b} with $h$ replaced by $h_{Rj}$),
\begin{multline*}  
Op^w_{h_{Rj}}\Big( (1-\phi)\Big(\frac{\xi}{h_{Rj}}\Big)\psi(\xi)\chi_1(h^\sigma 2^j\xi)\Big)\varphi(x) =\\ Op^w_{h_{Rj}}\Big( (1-\phi)\Big(\frac{\xi}{h_{Rj}}\Big)\psi(\xi)\chi_1(h^\sigma 2^j\xi)\varphi(x)\Big) + Op^w_{h_{Rj}}(r(x,\xi))
\end{multline*}
with 
\begin{equation*}
r(x,\xi) = \frac{h_{Rj}}{2i (\pi h_{Rj})^2}\int e^{-\frac{2i}{h_{Rj}}y\cdot\zeta} \left[\int_0^1\partial_\xi\Big[(1-\phi)\Big(\frac{\xi}{h_{Rj}}\Big)\psi(\xi)\chi_1(h^\sigma 2^j\xi)\Big]\big|_{(\xi + t\zeta)} dt\right] (\partial\varphi)(x+y) dyd\zeta.
\end{equation*}
Similarly as before, one can prove that 
\[\left\|\partial^\alpha_x \partial^\beta_\xi \big[r(\frac{x+y}{2}, h\xi)\big]\right\|_{L^2(d\xi)}\lesssim 1\]
for any $\alpha,\beta\in\mathbb{N}^2$, observing that no negative power of $h_{Rj}$ appears in the right hand side of this inequality for the product of $\psi(\xi)$ with any derivative of $(1-\phi)(\frac{\xi}{h_{Rj}})$ is supported for $h_{Rj} \sim |\xi| \sim 1$. Hence lemma \ref{Lemma on inequalities for Op(A)} gives that operator $Op^w_{h_{Rj}}(r(x,\xi))$ is uniformly bounded from $L^2$ to $L^\infty$ and 
\begin{equation*} 
\big\|Op^w_{h_{Rj}}(r(x,\xi))\Theta_R w^\chi(t,\cdot)\big\|_{L^\infty}\lesssim \|\Theta_Rw^\chi(t,\cdot)\|_{L^2}\lesssim R^{-1}\|w^\chi(t,\cdot)\|_{L^2}\,.
\end{equation*}
Since symbol $(1-\phi)\big(\frac{\xi}{h_{Rj}}\big)\psi(\xi)\chi_1(h^\sigma 2^j\xi)\varphi(x)$ is supported for $|x|\sim|\xi|\sim 1$, 
\begin{multline*} 
(1-\phi)\Big(\frac{\xi}{h_{Rj}}\Big)\psi(\xi)\chi_1(h^\sigma 2^j\xi)\varphi(x) \\
= \sum_{l=1}^2\underbrace{\frac{(1-\phi)\big(\frac{\xi}{h_{Rj}}\big)\psi(\xi)\chi_1(h^\sigma 2^j\xi)\varphi(x)(Rx_l|2^j\xi| - 2^j\xi_l)}{|Rx|2^j\xi|-2^j\xi|^2}}_{a_l(x,\xi)}\big(Rx_l|2^j\xi| - 2^j\xi_l\big),
\end{multline*}
with $a_l(x,\xi) \in R^{-1}2^{-j}S_{0,0}(1)$ as long as $R\gg 1$, and by $h_{Rj}$-symbolic calculus
\begin{equation*} 
(1-\phi)\Big(\frac{\xi}{h_{Rj}}\Big)\psi(\xi)\chi_1(h^\sigma 2^j\xi)\varphi(x) = \sum_{l=1}^2 a_l(x,\xi)\sharp \Big[(Rx_l|2^j\xi| - 2^j\xi_l)\widetilde{\psi}(\xi)\Big] + r_{Rj}(x,\xi),
\end{equation*}
with $\widetilde{\psi}\in C^\infty_0(\mathbb{R}^2\setminus\{0\})$ such that $\widetilde{\psi}\psi \equiv \psi$, and $r_{Rj}\in h_{Rj}S_{0,0}(1)$.
From semi-classical Sobolev injection 
$$\|Op^w_{h_{Rj}}(r_{Rj}(x,\xi))\Theta_Rw^\chi(t,\cdot)\|_{L^\infty}\lesssim \|\Theta_Rw^\chi(t,\cdot)\|_{L^2}\le R^{-1}\|w^\chi(t,\cdot)\|_{L^2}$$
while
\begin{equation}\label{last formula Op(al)}
\begin{split}
& Op^w_{h_{Rj}}(a_l(x,\xi))Op_{h_{Rj}}^w\big((Rx_l|2^j\xi| - 2^j\xi)\widetilde{\psi}(\xi)\big)\Theta_Rw^\chi \\
&= Op_{h_{Rj}}^w(a_l(x,\xi))\Theta_R \Big[\oph\big((x_l|\xi| - \xi)\widetilde{\psi}(2^{-j}\xi)\big)w^\chi\Big] \\
&= Op_{h_{Rj}}^w(a_l(x,\xi))\Theta_R \Big[ \oph(\widetilde{\psi}(2^{-j}\xi))\oph(x_l|\xi| - \xi)w^\chi -\frac{h}{2i} \oph((2^{-j}\xi)\cdot(\partial\widetilde{\psi})(2^{-j}\xi))w^\chi\Big].
\end{split}
\end{equation}
The last thing to do to conclude the proof of the statement is to study continuity of operator $Op_{h_{Rj}}^w(a_l(x,\xi))$.

\begin{lem} \label{Lemma: continuity of Op(al)}
We have that $Op^w_{h_{Rj}}(a_l(x,\xi)) :L^2 \rightarrow L^\infty$ is bounded with norm
\begin{equation*}
\left\|Op^w_{h_{Rj}}(a_l(x,\xi))\right\|_{\mathcal{L}(L^2;L^\infty)} \lesssim h^{-1}.
\end{equation*}
\proof
The result comes straightly from lemma \ref{Lemma on inequalities for Op(A)}. Indeed, since symbol $a_l(x,\xi)$ is compactly supported in $x$ there is a smooth cut-off function $\varphi_1 \in C^\infty_0(\mathbb{R}^2\setminus\{0\})$, with $\varphi_1 \varphi \equiv \varphi$, such that
\begin{equation*} 
\left|Op^w_{h_{Rj}}(a_l(x,\xi))w\right| 
\lesssim\|w\|_{L^2(dx)}\int \Big|\varphi_1\Big(\frac{x+y}{2}\Big)\Big| \sum_{|\alpha|\le 3}\left\|\partial^\alpha_y \left[a_l\Big(\frac{x+y}{2},h_{R_j}\xi\Big)\right]\right\|_{L^2(d\xi)}dy,
\end{equation*}
and for $|\alpha|\le 3$
\begin{align*}
&\left\|\partial^\alpha_y \left[a_l\Big(\frac{x+y}{2},h_{R_j}\xi\Big)\right]\right\|_{L^2(d\xi)} \\
&\hspace{20pt}\lesssim\frac{R}{h}\left\|\partial^\alpha_y \left[\frac{(1-\phi)(\xi)\psi(h_{Rj}\xi)\chi_1( h_{Rj} h^\sigma2^j\xi)\varphi_1(\frac{x+y}{2})}{|R(\frac{x+y}{2})|\xi| -\xi|^2}\Big(R\Big(\frac{x_l+y_l}{2}\Big)|\xi| -\xi_l\Big)\right]\right\|_{L^2(d\xi)}\\ 
&\hspace{20pt} \lesssim \frac{|\widetilde{\varphi}(\frac{x+y}{2})|}{h} \left(\int \frac{|\psi(h_{Rj}\xi)|^2}{|\xi|^2} d\xi\right)^{\frac{1}{2}}  \lesssim \frac{|\widetilde{\varphi}(\frac{x+y}{2})|}{h},
\end{align*}
where $\widetilde{\varphi}\in C^\infty_0(\mathbb{R}^2\setminus\{0\})$.
\endproof
\end{lem}
Finally, summing up all formulas from \eqref{operators Theta R} to \eqref{last formula Op(al)} and using lemma \ref{Lemma: continuity of Op(al)}, we obtain that 
\begin{equation*}
\Big\|\oph\Big((1-\phi)\Big(\frac{R}{h}\xi\Big)\psi(2^{-j}\xi)\chi_1(h^\sigma\xi)\Big)\Big[\varphi\Big(\frac{x}{R}\Big)w^\chi(t,\cdot)\Big]\Big\|_{L^\infty}\lesssim R^{-1}(\|w^\chi(t,\cdot)\|_{L^2} + \|\mathcal{M}w^\chi(t,\cdot)\|_{L^2}),
\end{equation*}
for any index $j\in \mathbb{Z}$ such that $hR^{-1}\le 2^j\le h^{-\sigma}$.
Injecting \eqref{Linfty_wchi_phi(Rh-1xi)} and the above inequality in \eqref{Linfty_wchi_split1}, and using that $[\mathcal{M}, \oph(\chi(h^\sigma\xi))]=i \oph((\partial \chi)(h^\sigma\xi)(h^\sigma |\xi|))$ is uniformly continuous on $L^2$, we deduce \eqref{ineq:utilde xbig} (the loss in $\log{R}+|\log h|$ arising from the fact that we are considering a sum over indices $j$, with $\log h-\log{R}\lesssim j \lesssim \log (h^{-1})$).
\endproof  
\end{lem}

\subsection{Proof of the main theorems}

\proof[Proof of theorem \ref{Thm: bootstrap argument}]
Straightforward after propositions \ref{Prop: Propagation of the energy estimate}, \ref{Prop:propagation_unif_est_V}, \ref{Prop: Propagation uniform estimate U,RU}.
\endproof

\proof[Proof of theorem \ref{Thm: Main theorem}]
Let us prove that, for small enough data satisfying \eqref{condition_initial_data}, Cauchy problem \eqref{wave KG system}-\eqref{initial data} has a unique global solution.
This result follows by a local existence argument, after having proved that there exist two integers $n\gg \rho\gg 1$, two constants $A', B'>1$ sufficiently large, $\varepsilon_0>0$ sufficiently small, and $0<\delta\ll \delta_2\ll \delta_1\ll \delta_0$ small, such that, for any $0<\varepsilon<\varepsilon_0$, if $(u,v)$ is solution to 
\eqref{wave KG system}-\eqref{initial data} in $[1,T]\times\mathbb{R}^2$, for some $T>1$, with $\partial_{t,x}u\in C^0([1,T]; H^n(\mathbb{R}^2))$, $v\in C^0([1,T]; H^{n+1}(\mathbb{R}^2))\cap C^1([1,T];H^n(\mathbb{R}^2))$, and satisfies
\begin{subequations} \label{a-priori_estimate_u,v}
\begin{gather}
\|\partial_tu(t,\cdot)\|_{H^{\rho+1,\infty}} + \|\nabla_x u(t,\cdot)\|_{H^{\rho+1,\infty}} + \||D_x|u(t,\cdot)\|_{H^{\rho+1,\infty}}+ \sum_{j=1}^2\|\mathrm{R}_j\partial_tu(t,\cdot)\|_{H^{\rho+1,\infty}}\le A'\varepsilon t^{-\frac{1}{2}},\\
\|\partial_tv(t,\cdot)\|_{H^{\rho,\infty}} + \|v(t,\cdot)\|_{H^{\rho+1,\infty}}\le A'\varepsilon t^{-1}, \\
\|\partial_tu(t,\cdot)\|_{H^n} + \|\nabla_x u(t,\cdot)\|_{H^n}+ \|\partial_tv(t,\cdot)\|_{H^n}+\|\nabla_xv(t,\cdot)\|_{H^n}+\|v(t,\cdot)\|_{H^n}\le B'\varepsilon t^\frac{\delta}{2}, 
\end{gather}
\begin{multline}
\sum_{|I|=k} \left[\|\partial_t \Gamma^I u(t,\cdot)\|_{L^2} + \|\nabla_x \Gamma^I u(t,\cdot)\|_{L^2}+ \|\partial_t \Gamma^I v(t,\cdot)\|_{L^2}+\|\nabla_x \Gamma^Iv(t,\cdot)\|_{L^2}\right.\\
\left.+\|\Gamma^I v(t,\cdot)\|_{L^2}\right]\le B'\varepsilon t^\frac{\delta_{3-k}}{2}, \quad 1\le k \le 3,
\end{multline}
\end{subequations}\normalsize
for every $t\in [1,T]$, then in the same interval it satisfies
\begin{subequations}\label{a-priori_enhance_estimate_u,v}
\begin{gather}
\|\partial_tu(t,\cdot)\|_{H^{\rho+1,\infty}} + \|\nabla_x u(t,\cdot)\|_{H^{\rho+1,\infty}} + \||D_x|u(t,\cdot)\|_{H^{\rho+1,\infty}}+ \sum_{j=1}^2\|\mathrm{R}_j\partial_tu(t,\cdot)\|_{H^{\rho+1,\infty}}\le  \frac{A'}{2}\varepsilon t^{-\frac{1}{2}},\\
\|\partial_tv(t,\cdot)\|_{H^{\rho,\infty}} + \|v(t,\cdot)\|_{H^{\rho+1,\infty}}\le \frac{A'}{2}\varepsilon t^{-1}, \\
\|\partial_tu(t,\cdot)\|_{H^n} + \|\nabla_x u(t,\cdot)\|_{H^n}+ \|\partial_tv(t,\cdot)\|_{H^n}+\|\nabla_xv(t,\cdot)\|_{H^n}+\|v(t,\cdot)\|_{H^n}\le  \frac{B'}{2}\varepsilon t^\frac{\delta}{2}, \\
\end{gather}
\begin{multline} 
\sum_{|I|=k} \left[\|\partial_t \Gamma^I u(t,\cdot)\|_{L^2} + \|\nabla_x \Gamma^I u(t,\cdot)\|_{L^2}+ \|\partial_t \Gamma^I v(t,\cdot)\|_{L^2}+\|\nabla_x \Gamma^Iv(t,\cdot)\|_{L^2}\right.\\
\left. +\|\Gamma^I v(t,\cdot)\|_{L^2}\right]\le  \frac{B'}{2}\varepsilon t^\frac{\delta_{3-k}}{2}, \quad 1\le k\le 3.
\end{multline}
\end{subequations}\normalsize
We remind that, if $I=(i_1,\dots, i_n)$ is a multi-index of length $|I|=n$, with $i_j\in \{1,\dots,5\}$, $\Gamma^I = \Gamma_{i_1}\cdots \Gamma_{i_n}$ is a product of vector fields in family $\mathcal{Z}=\{\Omega, Z_j, \partial_j |j=1,2\}$.

We can immediately observe that the above bounds are verified at time $t=1$ after \eqref{condition_initial_data} and Sobolev injection.
By definition \eqref{def u+- v+-} we also notice that
\begin{subequations}\label{dependence_unif_norm_u,v_upm,vpm}
\begin{multline} 
\|u_\pm (t,\cdot)\|_{H^{\rho+1,\infty}}+ \sum_{j=1}^2\|\mathrm{R}_ju_\pm (t,\cdot)\|_{H^{\rho+1,\infty}}\le 2 \|\partial_tu(t,\cdot)\|_{H^{\rho+1,\infty}}+ 2\||D_x|u(t,\cdot)\|_{H^{\rho+1,\infty}}\\
 + 2 \sum_{j=1}^2\left( \|\partial_j u(t,\cdot)\|_{H^{\rho+1,\infty}}+\|\mathrm{R}_j\partial_tu(t,\cdot)\|_{H^{\rho+1,\infty}}\right),
\end{multline}
\begin{equation}
\|v_\pm (t,\cdot)\|_{H^{\rho,\infty}}\le 2\|\partial_tv(t,\cdot)\|_{H^{\rho,\infty}}+ 2\|v(t,\cdot)\|_{H^{\rho+1,\infty}},
\end{equation}
\end{subequations}
and, conversely,
\begin{subequations} \label{dependence u,v in terms of upm vpm}
\begin{multline}
\|\partial_tu(t,\cdot)\|_{H^{\rho+1,\infty}} + \||D_x|u(t,\cdot)\|_{H^{\rho+1,\infty}}+ \sum_{j=1}^2\left(\|\partial_ju(t,\cdot)\|_{H^{\rho+1,\infty}}+ \|\mathrm{R}_j\partial_tu(t,\cdot)\|_{H^{\rho+1,\infty}}\right)  \\
\le \|u_+(t,\cdot)\|_{H^{\rho+1,\infty}}+ \|u_{-}(t,\cdot)\|_{H^{\rho+1,\infty}} + \sum_{j=1}^2 \left(\|\mathrm{R}_j u_+(t,\cdot)\|_{H^{\rho+1,\infty}}+ \|\mathrm{R}_j u_{-}(t,\cdot)\|_{H^{\rho+1,\infty}}\right),
\end{multline}
\begin{equation}
\|\partial_tv(t,\cdot)\|_{H^{\rho,\infty}}+ \|v(t,\cdot)\|_{H^{\rho+1,\infty}}\le \|v_+(t,\cdot)\|_{H^{\rho,\infty}}+\|v_{-}(t,\cdot)\|_{H^{\rho,\infty}}.
\end{equation}
\end{subequations}
Moreover, reminding definition \eqref{def_generalized_energy} of generalized energies $E_n(t;u_\pm, v_\pm)$, $E^k_3(t;u_\pm, v_\pm)$, for $n\ge 3$ and $0\le k \le 2$, and of set $\mathcal{I}^k_3$ in \eqref{set_Ik3}, there is a constant $C>0$
such that
\begin{subequations} \label{dependence_energy_u,v_energy upm vpm}
\begin{multline}
C^{-1}E_n(t;u_\pm, v_\pm) \le \left[\|\partial_tu(t,\cdot)\|^2_{H^n} + \|\nabla_x u(t,\cdot)\|^2_{H^n}\right.\\
\left.+ \|\partial_tv(t,\cdot)\|^2_{H^n}+\|\nabla_xv(t,\cdot)\|^2_{H^n}+\|v(t,\cdot)\|^2_{H^n}\right]  \le C E_n(t;u_\pm, v_\pm),
\end{multline}
and for any $0\le k\le 2$,
\begin{multline}
C^{-1}E^k_3(t;u_\pm, v_\pm) \le \sum_{I\in\mathcal{I}^k_3} \left[\|\partial_t \Gamma^I u(t,\cdot)\|^2_{L^2} + \|\nabla_x \Gamma^I u(t,\cdot)\|^2_{L^2}\right. \\
\left. + \|\partial_t \Gamma^I v(t,\cdot)\|^2_{L^2}+\|\nabla_x \Gamma^Iv(t,\cdot)\|^2_{L^2}+\|\Gamma^I v(t,\cdot)\|^2_{L^2}\right] \le C E^k_3(t;u_\pm, v_\pm).
\end{multline}
\end{subequations}
Therefore, after \eqref{dependence_unif_norm_u,v_upm,vpm}, \eqref{dependence_energy_u,v_energy upm vpm}, and \eqref{a-priori_estimate_u,v}, we deduce that estimates \eqref{est: bootstrap argument a-priori est} are satisfied with $A=2A'$, $B=C_1B'$, for some new $C_1>0$, so choosing for instance $K_1=4$ and $K_2$ sufficiently large, theorem \ref{Thm: bootstrap argument} and inequalities \eqref{dependence u,v in terms of upm vpm}, \eqref{dependence_energy_u,v_energy upm vpm} imply \eqref{a-priori_enhance_estimate_u,v}.
\endproof

\appendix

\chapter[Appendix A]{} \label{Appendix A}

\numberwithin{equation}{chapter}
\numberwithin{thm}{chapter}

The aim of this appendix is to prove the continuity of some trilinear integral operators (see lemmas \ref{Lem_appendix: est integrals Bj u v} and \ref{Lem_appendix: integral sigma_tilde_N}) that arise in subsection \ref{sub: second normal form} when performing a normal form argument at the energy level, and of some bilinear integral operators (see lemma \ref{Lem_Appendix: est on Dj1j2}) that instead appear in subsection \ref{Subsection: Section : Normal Forms for the Wave Equation} when we perform a normal form the wave equation (see proposition \ref{Prop: NF on wave}).
All the other results of this chapter are stated and proved in view of the above mentioned lemmas.

\begin{lem} \label{Lem_appendix: Kernel with 1 function}
Let $\check{a}(x)$ denote the inverse transform of a function $a(\xi)$.

$(i)$ If $a:\mathbb{R}^2\rightarrow \mathbb{C}$ is such that, for any $\alpha\in \mathbb{N}^2$ with $1\le |\alpha|\le 4$,
\[| a(\xi)|\lesssim \langle\xi\rangle^{-3}  \quad \text{and}\quad
|\partial^\alpha a(\xi)| \lesssim_{\alpha} (|\xi|\langle\xi\rangle^{-1})^{1-|\alpha|}\langle \xi\rangle^{-3}\quad \forall \xi\in\mathbb{R}^2
\]
then
\[
| \check{a}(x)|\lesssim |x|^{-1}\langle x \rangle^{-2}, \quad\forall x\in\mathbb{R}^2.
\]

$(ii)$ If $a$ is such that, for any $\alpha\in\mathbb{N}^2$ with $|\alpha|\le 3$,
\[|\partial^\alpha a(\xi)|\lesssim (|\xi|\langle\xi\rangle^{-1})^{-|\alpha|}\langle\xi\rangle^{-3},\quad \forall\xi\in\mathbb{R}^2 \]
then 
\[|\check{a}(x)|\lesssim\langle x \rangle^{-2}, \quad \forall x\in\mathbb{R}^2;\]

$(iii)$ Let $N\in\mathbb{N}$. If for any $\alpha\in\mathbb{N}^2$ with $|\alpha|\le N$ there exists $f_\alpha\in L^1(\mathbb{R}^2)$ such that $|\partial^\alpha a(\xi)|\lesssim_\alpha |f_\alpha(\xi)|$ then
\[|\check{a}(x)|\lesssim \langle x\rangle^{-N}, \quad \forall x\in\mathbb{R}^2.\]
\proof 
\textit{(i)} We consider a cut-off function $\phi\in C^\infty_0(\mathbb{R}^2)$ equal to 1 in the unit ball and write
\begin{equation}\label{splitting_check_a}
\begin{gathered}
\check{a}(x) = K_0(x) + K_1(x) \\
\text{with}\quad K_0(x) :=\frac{1}{(2\pi)^2} \int e^{ix\cdot\xi} a(\xi)\phi(\xi) d\xi, \quad K_1(x) :=\frac{1}{(2\pi)^2} \int e^{ix\cdot\xi} a(\xi)(1-\phi)(\xi) d\xi.
\end{gathered}
\end{equation}
On the one hand, since $|\partial^\alpha a(\xi)|\lesssim_{\alpha} \langle\xi\rangle^{-3}$ on the support of 
$(1-\phi)(\xi)$ for any $|\alpha|\le 4$, we immediately deduce by integration by parts that $| K_1(x)|\lesssim \langle x \rangle^{-4}$ for any $x\in \mathbb{R}^2$.
On the other hand, again an integration by parts gives that
\begin{equation*}
x  K_0(x) = \int e^{ix\cdot\xi} a_1(\xi) d\xi
\end{equation*}
with $a_1(\xi)$ supported for $|\xi|\lesssim 1$ and such that $|\partial^\alpha a_1(\xi)|\lesssim_\alpha |\xi|^{-|\alpha|}$ for any $\xi\in\mathbb{R}^2$, any $|\alpha|\le 3$. 
This implies that $|x  K_0(x)|\lesssim 1$ for any $x\in\mathbb{R}^2$.
Moreover, $|x^\alpha x\ K_0(x)|\lesssim_\alpha 1$ for any $|\alpha|\le 3$. 
This is obvious in the unit ball. Out of the unit ball we consider a Littlewood-Paley decomposition in frequencies so that
$$\phi(\xi) = \phi(\xi) \left[\varphi_0(2^{-L_0}\xi) + \sum_{k=L_0+1}^0\varphi(2^{-k}\xi)\right],$$
with $supp\varphi_0\subset B_1(0)$, $\varphi\in C^\infty_0(\mathbb{R}^2\setminus\{0\})$ and $L_0<0$ such that $2^{L_0}\sim |x|^{-1}$, and write
\begin{equation*}
\begin{gathered}
xK_0(x) = K_0^0(x) + \sum_{k=L_0+1}^0 K_0^k(x) \\
\text{with}\quad K^0_0(x) := \int e^{ix\cdot\xi} a_1(\xi) \varphi_0(2^{-L_0}\xi) d\xi, \quad K^k_0(x) := \int e^{ix\cdot\xi} a_1(\xi) \varphi_k(2^{-k}\xi) d\xi.
\end{gathered}
\end{equation*}
Performing a change of coordinates and making some integrations by parts we deduce that 
\[ |K^0_0(x)|\lesssim 2^{2L_0} \quad \text{and} \quad |K^k_0(x)|\lesssim 2^{2k}\langle 2^k x\rangle^{-3}, \quad L_0+1\le k\le 0 \]
for any $x\in\mathbb{R}^2$, which finally implies $|x  K_0(x)|\lesssim 2^{2L_0}\sim |x|^{-2} $. 

\textit{(ii)} The result follows splitting $\check{a}$ as in \eqref{splitting_check_a} and applying to $K_0(x)$ the same argument previously used for $xK_0(x)$.

$(iii)$ The result follows straightly from integration by parts and the fact that $f_\alpha\in L^1(\mathbb{R}^2)$ for any $|\alpha|\le N$.
\endproof
\end{lem}

\begin{cor} \label{Cor_appendix: decay of integral operators}
Let $d\in\mathbb{N}^*$, $N\in\mathbb{N}$ and $g_\beta \in L^1(\mathbb{R}^d)$ for every $|\beta|\le N$.

$(i)$ If $a(\xi,\eta):\mathbb{R}^2\times\mathbb{R}^d\rightarrow \mathbb{C}$ is such that, for any $\beta\in\mathbb{N}^d$ with $|\beta|\le N$,
\begin{equation} \label{ineq_a_1}
\begin{gathered}
|\partial^\beta_\eta a(\xi,\eta)|\lesssim_\beta \langle\xi\rangle^{-3}|g_\beta(\eta)|, \\
|\partial^\alpha_\xi \partial^\beta_\eta a(\xi,\eta)|\lesssim_{\alpha,\beta} (|\xi|\langle\xi\rangle^{-1})^{1-|\alpha|} \langle\xi\rangle^{-3}|g_\beta(\eta)|, \quad 1\le |\alpha|\le 4.
\end{gathered}
\end{equation}
for any $(\xi,\eta)\in\mathbb{R}^2\times\mathbb{R}^d$, then 
\begin{equation}\label{eq_cor_1}
\left|\int e^{ix\cdot\xi + iy\cdot\eta} a(\xi,\eta) d\xi d\eta\right| \lesssim |x|^{-1}\langle x \rangle^{-2}\langle y \rangle^{-N}, \quad \forall (x,y)\in\mathbb{R}^2\times \mathbb{R}^d.
\end{equation}
Moreover, if $d=2$ and $N=3$, for any $u,v\in L^2(\mathbb{R}^2)\cap L^\infty(\mathbb{R}^2)$
\begin{subequations}\label{est:coroll_app_L_norm}
\begin{equation}\label{est: corollary_app L2 norm}
\left\|\int e^{ix\cdot\xi}a(\xi,\eta)\hat{u}(\xi-\eta)\hat{v}(\eta) d\xi d\eta \right\|_{L^2(dx)} \lesssim \|u\|_{L^2}\|v\|_{L^\infty} \, (\text{or } \lesssim \|u\|_{L^\infty}\|v\|_{L^2})
\end{equation}
and
\begin{equation}
\left\|\int e^{ix\cdot\xi}a(\xi,\eta)\hat{u}(\xi-\eta)\hat{v}(\eta) d\xi d\eta \right\|_{L^\infty(dx)} \lesssim \|u\|_{L^\infty}\|v\|_{L^\infty}.
\end{equation}
\end{subequations}

$(ii)$ If $a(\xi,\eta)$ is such that, for any $\alpha\in\mathbb{N}^2$ with $|\alpha|\le 3$, $\beta\in\mathbb{N}^d$ with $|\beta|\le N$,
\begin{equation} \label{ineq_a_2}
|\partial^\alpha_\xi \partial^\beta_\eta a(\xi,\eta)|\lesssim_{\alpha,\beta} (|\xi|\langle\xi\rangle^{-1})^{-|\alpha|} \langle\xi\rangle^{-3}|g_\beta(\eta)|, 
\end{equation}
for any $(\xi,\eta)\in\mathbb{R}^2\times\mathbb{R}^d$, then
\begin{equation}\label{eq_cor_2}
\left|\int e^{ix\cdot\xi + iy\cdot\eta} a(\xi,\eta) d\xi d\eta\right| \lesssim\langle x \rangle^{-2}\langle y \rangle^{-N}, \quad \forall (x,y)\in\mathbb{R}^2\times \mathbb{R}^d.
\end{equation}
Moreover, if $d=2,N=3$, for any $u,v\in L^2(\mathbb{R}^2)$
\begin{subequations} \label{ineq_corA2}
\begin{equation}
\left\|\int e^{ix\cdot\xi}a(\xi,\eta)\hat{u}(\xi-\eta)\hat{v}(\eta) d\xi d\eta \right\|_{L^2(dx)} \lesssim \|u\|_{L^2}\|v\|_{L^2}
\end{equation}
while if $u\in L^2(\mathbb{R}^2), v\in L^\infty(\mathbb{R}^2)$,
\begin{equation}
\left\|\int e^{ix\cdot\xi}a(\xi,\eta)\hat{u}(\xi-\eta)\hat{v}(\eta) d\xi d\eta \right\|_{L^\infty(dx)} \lesssim \|u\|_{L^2}\|v\|_{L^\infty}.
\end{equation}
\end{subequations}
\proof
Let 
\begin{equation*}
K(x,\eta):=\int e^{ix\cdot\xi}a(\xi,\eta)d\xi \quad\text{and}\quad \widetilde{K}(x,y) := \int e^{ix\cdot\xi }K(x,\eta)  d\eta.
\end{equation*}
By the hypothesis on $a(\xi,\eta)$ and lemma \ref{Lem_appendix: Kernel with 1 function} $(i)$ (resp. $(ii)$) we derive that, for any $\beta\in\mathbb{N}^d$ with $|\beta|\le N$,
\[|\partial^\beta_\eta K(x,\eta)|\lesssim |x|^{-1}\langle x\rangle^{-2}|g_\beta(\eta)| \quad \left(\text{resp. } |\partial^\beta_\eta K(x,\eta)|\lesssim \langle x\rangle^{-2}|g_\beta(\eta)|\right) \quad \forall (x,\eta)\in\mathbb{R}^2\times\mathbb{R}^d.\]
Hence \eqref{eq_cor_1} (resp. \eqref{eq_cor_2}) follows applying lemma \ref{Lem_appendix: Kernel with 1 function} $(iii)$ to $\widetilde{K}(x,y)$.

$(i)$ If $d=2, N=3$, inequality \eqref{est: corollary_app L2 norm} from the fact that
\begin{equation*}
\int e^{ix\cdot\xi}a(\xi, \eta)\hat{u}(\xi-\eta) \hat{v}(\eta) d\eta  = \int \widetilde{K}(x-y, y-z) u(y)v(z) dydx,
\end{equation*}
and by \eqref{eq_cor_1}, for $L=L^2$ or $L=L^\infty$,
\begin{equation} \label{ineq: norm L(dx) kernel}
\begin{split}
\left\| \int \widetilde{K}(x-y,y-z) \widetilde{u}(y) \widetilde{v}(z) dydz\right\|_{L(dx)}& \lesssim \left\| \int |x-y|^{-1}\langle x-y\rangle^{-2}\langle y-z\rangle^{-3} \widetilde{u}(y) \widetilde{v}(z) dydz\right\|_{L(dx)} \\
& \lesssim \int |y|^{-1}\langle y \rangle^{-2}\langle z \rangle^{-3} \|\widetilde{u}(\cdot -y) \widetilde{v}(\cdot - y -z)\|_{L(dx)} dydz\\
& \lesssim \|\widetilde{u}\|_{L^\infty}\|\widetilde{v}\|_{L} \ (\text{or } \lesssim \|\widetilde{u}\|_{L}\|\widetilde{v}\|_{L^\infty} ).
\end{split}
\end{equation}

$(ii)$ By inequality \eqref{eq_cor_2}
\begin{equation*}
\begin{split}
&\left\| \int \widetilde{K}(x-y,y-z)u(y) v(z) dydz\right\|_{L^2(dx)}\lesssim \left\| \int\langle x-y\rangle^{-2}\langle y-z\rangle^{-3}| u(y) ||v(z)| dydz\right\|_{L^2(dx)}\\
& \lesssim \int  \langle y- z\rangle^{-3}|u(y)| |v(z)| dydz \lesssim \int |v(z)| \left(\int \langle y-z\rangle^{-3}dy\right)^\frac{1}{2} \left(\int \langle y-z\rangle^{-3}|u(y)|^2 dy\right)^\frac{1}{2} dz\\
& \lesssim \|v\|_{L^2} \left(\int \langle y-z\rangle^{-3}|u(y)|^2 dydz\right)^\frac{1}{2}\lesssim \|u\|_{L^2}\|v\|_{L^2}
 \end{split}
\end{equation*}
and
\begin{multline*}
\left\| \int \widetilde{K}(x-y,y-z)u(y) v(z) dydz\right\|_{L^\infty(dx)}\lesssim \left\| \int\langle x-y\rangle^{-2}\langle y-z\rangle^{-3}| u(y) ||v(z)| dydz\right\|_{L^\infty(dx)}\\
 \lesssim \|v\|_{L^\infty} \left\|\int \langle x-y\rangle^{-2}|u(y)| dy \right\|_{L^\infty(dx)}
 \lesssim \|u\|_{L^2}\|v\|_{L^\infty}.
\end{multline*}
\endproof
\end{cor}

\begin{lem}[Sobolev norm of a product] \label{Lem:Sobolev norm of products}
Let $s\in\mathbb{N}^*$. For any $u,v\in H^{s}(\mathbb{R}^2)\cap L^\infty(\mathbb{R}^2)$,
\begin{equation} \label{Hs_norm_product}
\|uv\|_{H^s}\lesssim \|u\|_{H^{s}}\|v\|_{L^\infty}+ \|u\|_{L^\infty}\|v\|_{H^{s}};
\end{equation}
for any $u,v\in H^{s,\infty}(\mathbb{R}^2)\cap H^{s+2}(\mathbb{R}^2)$, any $\theta\in ]0,1[$,
\begin{equation} \label{Hsinfty_norm_product}
\|uv\|_{H^{s,\infty}}\lesssim \|u\|^{1-\theta}_{H^{s,\infty}}\|u\|^\theta_{H^{s+2}}\|v\|_{L^\infty}+\|u\|_{L^\infty}\|v\|^{1-\theta}_{H^{s,\infty}}\|v\|^\theta_{H^{s+2}}.
\end{equation}
\proof
Inequality \eqref{Hs_norm_product} is a classical result (see, for instance, \cite{Alinhac-Gerard}).

In order to deduce \eqref{Hsinfty_norm_product} we decompose product $uv$ as follows:
\begin{equation}\label{dec_para-products}
uv= T_uv+ T_vu +R(u,v),
\end{equation}
where $T_uv$ is the para-product of $u$ times $v$ defined by
\begin{equation*}
T_u v:= S_{-3}u S_0 v +\sum_{k\ge 1}S_{k-3}u \Delta_kv,
\end{equation*}
with $S_k = \chi(2^{-k}D_x)$, $\chi\in C^\infty_0(\mathbb{R}^2)$ such that $\chi(\xi)=1$ for $|\xi|\le 1/2$, $\chi(\xi)=0$ for $|\xi|\ge 1$, $\Delta_0= S_0$ and $\Delta_k = S_k- S_{k-1}$ for $k\ge 1$, and $R(u,v)=\sum_{k}\Delta_k u\widetilde{\Delta}_k v,$
with $\widetilde{\Delta}_k=\Delta_{k-1}+\Delta_k +\Delta_{k+1}$.
Since
\begin{equation*}
T_uv =\sum_{j\ge 0}\Delta_j (T_uv) =\sum_{\substack{j,k\\ |j-k|\le N_0}} \Delta_j [S_{k-3}u \Delta_k v]
\end{equation*}
for a certain $N_0\in\mathbb{N}$, by definition \ref{def Sobolev spaces-NEW} $(iii)$ of the $H^{s,\infty}$ norm and the fact that $\|\Delta_k v\|_{L^\infty}\lesssim 2^k\|\Delta_kv\|_{L^2}$ we deduce that, for any fixed $\theta\in ]0,1[$,
\begin{equation}\label{Hsinfty_norm_Tuv}
\begin{split}
&\|T_uv\|_{H^{s,\infty}} = \|\langle D_x\rangle^s T_uv\|_{L^\infty} \le \sum_{\substack{j,k\\ |j-k|\le N_0}}2^{js}\|\Delta_j [S_{k-3}u \Delta_k v]\|_{L^\infty}\\
& \le \sum_{\substack{j,k\\ |j-k|\le N_0}}2^{js} \|S_{k-3}u\|_{L^\infty} \|\Delta_kv\|_{L^\infty} \le \sum_{\substack{j,k\\ |j-k|\le N_0}}2^{js} \|u\|_{L^\infty} (2^{-ks}\|\Delta_k \langle D_x\rangle^s v\|_{L^\infty})^{1-\theta}(2^k\|\Delta_k v\|_{L^2})^\theta \\
&\lesssim \sum_{\substack{j,k\\ |j-k|\le N_0}}2^{(j-k)s} \|u\|_{L^\infty}\|\Delta_k \langle D_x\rangle^s v\|_{L^\infty}^{1-\theta}\left(2^{-k}\|\Delta_k \langle D_x\rangle^{s+2}v\|_{L^2}\right)^\theta\\
&\lesssim \|u\|_{L^\infty}\|v\|^{1-\theta}_{H^{s,\infty}}\|v\|^\theta_{H^{s+2}}.
\end{split}
\end{equation}
Similarly, 
\begin{equation*}
\|T_vu\|_{H^{s,\infty}} + \|R(u,v)\|_{H^{s,\infty}}\lesssim \|u\|^{1-\theta}_{H^{s,\infty}}\|u\|^\theta_{H^{s+2}}\|v\|_{L^\infty}.
\end{equation*}
\endproof
\end{lem}

\begin{cor}\label{Cor_appendixA:Hs-Hsinfty norm of bilinear expressions}
Let $s\in \mathbb{N}^*$, $a_1(\xi)\in S^{m_1}_0(\mathbb{R}^2)$, $a_2(\xi)\in S^{m_2}_0(\mathbb{R}^2)$, for some $m_1,m_2\ge 0$.
For any $u\in H^{s+m_1}(\mathbb{R}^2)\cap H^{m_1,\infty}(\mathbb{R}^2)$, $v\in H^{s+m_2}(\mathbb{R}^2)\cap H^{m_2,\infty}(\mathbb{R}^2)$,
\begin{equation} \label{Hs_norm_bilinear-expressions}
\left\|[a_1(D_x)u]\, [a_2(D_x)v] \right\|_{H^s}\lesssim \|u\|_{H^{s+m_1}}\|v\|_{H^{m_2,\infty}}+ \|u\|_{H^{m_1,\infty}}\|v\|_{H^{s+m_2}};
\end{equation}
for any $u\in H^{s+m_1,\infty}(\mathbb{R}^2)\cap H^{s+m_1+2}(\mathbb{R}^2)$, $v\in H^{s+m_2,\infty}(\mathbb{R}^2)\cap H^{s+m_2+2}(\mathbb{R}^2)$, any $\theta\in ]0,1[$,
\begin{multline} \label{Hsinfty_norm_bilinear-expressions}
\left\|[a_1(D_x)u]\, [a_2(D_x)v]  \right\|_{H^{s,\infty}}\\
\lesssim \|u\|^{1-\theta}_{H^{s+m_1,\infty}} \|u\|^\theta_{H^{s+m_1+2}}\|v\|_{H^{m_2,\infty}}+ \|u\|_{H^{m_1,\infty}}\|v\|^{1-\theta}_{H^{s+m_2,\infty}}\|v\|^\theta_{H^{s+m_2+2}}.
\end{multline}
\proof
The result of the statement follows writing $[a_1(D_x)u]\, [a_2(D_x)v]$ in terms of para-products as in \eqref{dec_para-products}, and using that $T_{a_1(D)u}(a_2(D)v)$, $T_{a_2(D)v}(a_1(D)u)$ and remainder $R\big(a_1(D)u, a_2(D)v\big)$ can be written from $\widetilde{u}=\langle D_x\rangle^{m_1}u$, $\widetilde{v}=\langle D_x\rangle^{m_2}v$, as done below for the former of these terms:
\begin{equation*}
\begin{split}
T_{a_1(D)u}(a_2(D)v) &= [S_{-3}a_1(D)\langle D_x\rangle^{-m_1}\widetilde{u}] [S_0a_2(D)\langle D_x\rangle^{-m_2}\widetilde{v}] \\
& + \sum_{k}[S_{k-3}a_1(D)\langle D_x\rangle^{-m_1}\widetilde{u}] [\Delta_k a_2(D)\langle D_x\rangle^{-m_2}\widetilde{v}].
\end{split}
\end{equation*}
Since $a_j(\xi)\langle\xi\rangle^{-m_j}\in S^0_0(\mathbb{R}^2)$, $j=1,2$, operators $S_k a_j(D)\langle D_x\rangle^{-m_j}$, $\Delta_k a_j(D)\langle D_x\rangle^{-m_j}$ have the same spectrum (i.e. the support of the Fourier transform) up to a negligible constant of $S_k$ and $\Delta_k$ respectively.
\endproof
\end{cor}

In the following lemma we prove a result of continuity for a trilinear integral operator defined from multiplier $B^k_{(j_1,j_2,j_3)}(\xi,\eta)$ given by \eqref{def of B(i1,i2,i3)}
(resp. by \eqref{multiplier_B3}) for $k=1,2$ (resp. $k=3$), any $j_1,j_2,j_3\in \{+,-\}$.
It is useful to observe that, since
\[B^k_{(j_1,j_2,j_3)}(\xi,\eta) = \frac{j_1\langle\xi-\eta\rangle + j_2|\eta| - j_3\langle\xi\rangle}{2j_1 j_2\langle\xi-\eta\rangle|\eta|}\eta_k, \quad k=1,2\] 
from \eqref{def of B(i1,i2,i3)}, while 
\[B^3_{(j_1,j_2,j_3)}(\xi,\eta) = j_1\, \frac{j_1\langle\xi-\eta\rangle + j_2|\eta| - j_3\langle\xi\rangle}{2\langle\xi-\eta\rangle} \]
\eqref{multiplier_B3}, we have that
\begin{multline}\label{explicit integral B}
\frac{1}{(2\pi)^2}\int e^{ix\cdot\xi}B^k_{(j_1,j_2,j_3)}(\xi,\eta)\hat{u}(\xi-\eta)\hat{v}(\eta) d\xi d\eta = \frac{j_2}{2} (u\mathrm{R}_kv)(x) - \frac{j_1}{2} \Big[\Big(\frac{D_1}{\langle D_x\rangle}u\Big)v\Big](x) \\
+ \frac{j_1}{2}D_1 \left[(\langle D_x\rangle^{-1}u)v\right](x) - \frac{j_3}{2j_1j_2}\langle D_x\rangle [(\langle D_x\rangle^{-1}u)(\mathrm{R}_kv)](x)
\end{multline}
for $k=1,2$, while for $k=3$
\begin{multline} \label{explicit_integral_B3}
\frac{1}{(2\pi)^2}\int e^{ix\cdot\xi}B^3_{(j_1,j_2,j_3)}(\xi,\eta)\hat{u}(\xi-\eta)\hat{v}(\eta) d\xi d\eta = \frac{1}{2}(uv)(x) + \frac{j_1j_2}{2}[ (\langle D_x\rangle^{-1}u) |D_x|v](x) \\
- \frac{j_1j_3}{2}\langle D_x\rangle\left[(\langle D_x\rangle^{-1} u) v\right](x).
\end{multline}

\begin{lem} \label{Lem_appendix: est integrals Bj u v}
Let $B^k_{(j_1,j_2,j_3)}(\xi,\eta)$ be given by \eqref{def of B(i1,i2,i3)} when $k=1,2$, and by \eqref{multiplier_B3} when $k=3$, for any $j_1,j_2,j_3\in \{+,-\}$. Let also $\delta_k=1$ if $k\in \{1,2\}$, $\delta_k=0$ if $k=3$.
For any $u,w\in L^2(\mathbb{R}^2), v\in H^{2,\infty}(\mathbb{R}^2)$ such that $\delta_k\mathrm{R}_kv\in H^{2,\infty}(\mathbb{R}^2)$,
\begin{equation} \label{estimate integral B uvw-new}
\left|\int B^k_{(j_1,j_2,j_3)}(\xi,\eta)\hat{u}(\xi-\eta)\hat{v}(\eta)\hat{w}(-\xi) d\xi d\eta\right| \lesssim \|u\|_{L^2}\left(\|v\|_{H^{7,\infty}}+\delta_k \|\mathrm{R}_kv\|_{H^{7,\infty}}\right) \|w\|_{L^2}.
\end{equation}
\proof
First of all we observe that for $k\in\{1,2\}$\small
\begin{subequations}
\begin{multline} \label{int_Bk_app}
\int B^k_{(j_1,j_2,j_3)}(\xi,\eta)\hat{u}(\xi-\eta)\hat{v}(\eta)\hat{w}(-\xi) d\xi d\eta = \frac{j_2}{2} \int \reallywidehat{u(\mathrm{R}_kv)}(\xi) \hat{w}(-\xi)d\xi - \frac{j_1}{2} \int \reallywidehat{\big[\big(\frac{D_x}{\langle D_x\rangle}u\big)v\big]}(\xi) \hat{w}(-\xi) d\xi  \\
+ \frac{j_1}{2}\int \frac{\xi_1}{\langle\xi-\eta\rangle} \hat{u}(\xi-\eta)\hat{v}(\eta)\hat{w}(-\xi) d\xi d\eta -\frac{j_3}{2j_1j_2}\int \frac{\langle\xi\rangle}{\langle\xi-\eta\rangle}  \hat{u}(\xi-\eta)\widehat{\mathrm{R}_k v}(\eta)\hat{w}(-\xi) d\xi d\eta,
\end{multline}\normalsize
while for $k=3$,
\begin{multline}\label{int_B3_app}
\int B^3_{(j_1,j_2,j_3)}(\xi,\eta)\hat{u}(\xi-\eta)\hat{v}(\eta)\hat{w}(-\xi) d\xi d\eta = \frac{1}{2}\int \reallywidehat{uv}(\xi)\hat{w}(-\xi)d\xi \\+ \frac{j_1j_2}{2} \int \reallywidehat{\left[(\langle D_x\rangle^{-1}u)|D_x|v\right]}(\xi)\hat{w}(-\xi) d\xi
-\frac{j_1j_3}{2} \int \frac{\langle\xi\rangle}{\langle\xi-\eta\rangle} \hat{u}(\xi-\eta)\hat{v}(\eta)\hat{w}(-\xi)d\xi.
\end{multline}
\end{subequations}
H\"older's inequality shows immediately that the first two addends in both above right hand sides are bounded by the right hand side of \eqref{estimate integral B uvw-new}.
Then the result of the statement follows by proving that inequality \eqref{estimate integral B uvw-new} is satisfied by integrals such as
\begin{equation*}
\int a(\xi,\eta) \hat{u}_1(\xi-\eta) \hat{u}_2(\eta)\hat{u}_3(-\xi) d\xi d\eta
\end{equation*}
with $a(\xi,\eta) = \xi_1\langle\xi-\eta\rangle^{-1}$ or $a(\xi,\eta) = \langle\xi\rangle \langle\xi-\eta\rangle^{-1}$, and some general functions $u_1,u_3\in L^2(\mathbb{R}^2), u_2\in L^\infty(\mathbb{R}^2)$. 
By taking a Littlewood-Paley decomposition we can split the above integral as 
\begin{equation} \label{decomposition u1u2u3}
\sum_{k,l\ge 0}\int a(\xi,\eta)\varphi_k(\xi)\varphi_l(\eta)  \hat{u}_1(\xi-\eta) \hat{u}_2(\eta)\hat{u}_3(-\xi) d\xi d\eta, 
\end{equation}
with $\varphi_0\in C^\infty_0(\mathbb{R}^2)$, $\varphi \in C^\infty_0(\mathbb{R}^2\setminus\{0\})$ and $\varphi_k(\zeta) = \varphi(2^{-k}\zeta)$ for any $k\in\mathbb{N}^*$.
Since frequencies $\xi,\eta$ are bounded on the support of $\varphi_0(\xi)\varphi_0(\eta)$,
kernel 
\begin{equation*}
K_0(x,y):=\int e^{ix\cdot\xi+iy\cdot\eta}a(\xi,\eta)\varphi_0(\xi)\varphi_0(\eta) d\xi d\eta
\end{equation*}
is such that $|K_0(x,y)|\lesssim \langle x \rangle^{-3}\langle y\rangle^{-3}$ for any $(x,y)$, after the first part of corollary \ref{Cor_appendix: decay of integral operators} $(i)$.
Therefore
\begin{equation*}
\begin{split}
& \left| \int a(\xi,\eta)\varphi_0(\xi)\varphi_0(\eta)  \hat{u}_1(\xi-\eta) \hat{u}_2(\eta)\hat{u}_3(-\xi) d\xi d\eta\right| \\
&= \left| \int K_0(z-x, x-y) u_1(x) u_2(y) u_3(z) dxdydz\right| \\
&\lesssim  \int \langle z-x\rangle^{-3}\langle x-y\rangle^{-3} |u_1(x)| |u_2(y)| |u_3(z)| dx dy dz \\
&\lesssim  \|u_2\|_{L^\infty}\int \langle x\rangle^{-3}|u_1(z-x)||u_3(z)| dxdz \lesssim \|u_1\|_{L^2}\|u_2\|_{L^\infty}\|u_3\|_{L^2},
\end{split}
\end{equation*}
where last inequality obtained by H\"older inequality.

For positive indices $l,k$ such that $l>k+N_0\ge 0$ (resp. $|l-k|\le N_0$), for a suitably large integer $N_0>1$, we have that $|\xi|<|\eta|\sim|\xi-\eta|$ (resp. $|\xi|\sim |\eta|$) on the support of $\varphi_k(\xi)\varphi_l(\eta)$.
If we define $a_{l>k+N_0}(\xi,\eta):=a(\xi,\eta)\langle\eta\rangle^{-1}$ and $a_{|l-k|\le N_0}(\xi,\eta):=a(\xi,\eta)\langle\eta\rangle^{-7}$, it is a computation to check that, for any $\alpha,\beta\in\mathbb{N}^2$ with $|\alpha|, |\beta|\le 3$,
\[ |\partial^\alpha_\xi \partial^\beta_\eta [a_{l>k+N_0}(2^k\xi, 2^l\eta)]| + |\partial^\alpha_\xi \partial^\beta_\eta a_{|l-k|\le N_0}(\xi,\eta)|\lesssim 2^{-l}. \]
Hence, their associated kernels $K_{l>k+N_0}(x,y)$ and $K_{|l-k|\le N_0}(x,y)$ are such that 
\[|K_{l>k+N_0}(x,y)| + |K_{|l-k|\le N_0}(x,y)|\lesssim 2^{2k}2^l\langle 2^kx\rangle^{-3}\langle 2^l y\rangle^{-3}, \quad \forall (x,y)\in\mathbb{R}^2\times\mathbb{R}^2 \]
as follows after a change of coordinates and some integrations by parts, and for any $l>k+N_0$ 
\begin{equation}\label{eq_app_1}
\begin{split}
& \left| \int a(\xi,\eta)\varphi_k(\xi)\varphi_l(\eta)  \hat{u}_1(\xi-\eta) \hat{u}_2(\eta)\hat{u}_3(-\xi) d\xi d\eta\right| \\
=& \left| \int K_{l>k+N_0}(z-x, x-y) u_1(x) [\langle D_x\rangle u_2](y) u_3(z) dxdydz\right| \\
\lesssim& 2^{2k}2^l \left|\int \langle 2^k(z-x)\rangle^{-3}\langle 2^l(x-y)\rangle^{-3} |u_1(x)||\langle D_x\rangle u_2(y)| |u_3(z)| dxdydz \right| \\
\lesssim & 2^{-\frac{k}{2}}2^{-\frac{l}{2}}\|u_1\|_{L^2}\|u_2\|_{H^{1,\infty}}\|u_3\|_{L^2},
\end{split}
\end{equation}
while for $|l-k|\le N_0$
\begin{equation}\label{eq_app_2}
\begin{split}
& \left| \int a(\xi,\eta)\varphi_l(\xi)\varphi_l(\eta)  \hat{u}_1(\xi-\eta) \hat{u}_2(\eta)\hat{u}_3(-\xi) d\xi d\eta\right| \\
=& \left| \int K_{|l-k|\le N_0}(z-x, x-y) u_1(x) [\langle D_x\rangle^7 u_2](y) u_3(z) dxdydz\right| \\
\lesssim &2^{3l}\left|\int \langle 2^l(z-x)\rangle^{-3}\langle 2^l(x-y)\rangle^{-3} |u_1(x)||\langle D_x\rangle^7 u_2(y)| |u_3(z)| dxdydz \right| \\
\lesssim & 2^{-l}\|u_1\|_{L^2}\|u_2\|_{H^{7,\infty}}\|u_3\|_{L^2}.
\end{split}
\end{equation}
Finally, when positive indices $l,k$ are such that $k>l-N_0$ we observe that frequencies $\xi$ and $\xi-\eta$ are equivalent and of size $2^k$ on the support of $\varphi_k(\xi)\varphi_l(\eta)$. 
If we take $a_{k>l-N_0}(\xi,\eta)$ equal to $a_{l>k+N_0}(\xi,\eta)$, denote by $K_{k>l-N_0}(x,y)$ its associated kernel (which is hence equal to $K_{l>k+N_0}(x,y)$), and introduce two new smooth cut-off function $\varphi^1, \varphi^2\in C^\infty_0(\mathbb{R}^2)$ equal to 1 on the support of $\varphi$, together with operators $\Delta^1_k:=\varphi^1(2^{-k}D_x), \Delta^2_k:=\varphi^2(2^{-k}D_x)$,
we deduce that
\begin{equation*}
\begin{split}
& \left| \int a(\xi,\eta)\varphi_k(\xi)\varphi_l(\eta)  \hat{u}_1(\xi-\eta) \hat{u}_2(\eta)\hat{u}_3(-\xi) d\xi d\eta\right| \\
=& \left| \int K_{k>l-N_0}(z-x, x-y) [\Delta^1_k u_1](x) [\langle D_x\rangle u_2](y) [\Delta^2_k u_3](z) dxdydz\right| \\
\lesssim& 2^{2k}2^l \left|\int \langle 2^k(z-x)\rangle^{-3}\langle 2^l(x-y)\rangle^{-3} |[\Delta^1_k u_1](x)||\langle D_x\rangle u_2(y)| |[\Delta^2_k u_3](z)| dxdydz \right| \\
\lesssim &2^{-l}\|\Delta^1_k u_1\|_{L^2}\|u_2\|_{H^{1,\infty}}\|\Delta^2_k u_3\|_{L^2}.
\end{split}
\end{equation*}
Combining decomposition \eqref{decomposition u1u2u3} together with \eqref{eq_app_1}, \eqref{eq_app_2} and Cauchy-Schwarz inequality we finally obtain that
\begin{equation*}
\left| \int a(\xi,\eta) \hat{u}_1(\xi-\eta) \hat{u}_2(\eta)\hat{u}_3(-\xi) d\xi d\eta \right| \lesssim \|u_1\|_{L^2}\|u_2\|_{H^{7,\infty}}\|u_3\|_{L^2},
\end{equation*}
\endproof
\end{lem}

\begin{lem} \label{Lem_appendix: integral sigma_tilde_N}
Let $\varepsilon>0$ be small, $N\in\mathbb{N}^*$, and $\sigma^N(\xi,\eta):\mathbb{R}^2\times\mathbb{R}^2\rightarrow \mathbb{C}$ be supported for $|\xi|\le \varepsilon \langle\eta\rangle$ and such that, for any $\alpha,\beta\in\mathbb{N}^2$,
\[|\partial^\alpha_\xi \partial^\beta_\eta \sigma^N(\xi,\eta)|\lesssim |\xi|^{N+1-|\alpha|}\langle\eta\rangle^{-N-|\beta|}, \quad \forall (\xi,\eta)\in\mathbb{R}^2\times\mathbb{R}^2.\]
For any $(j_1,j_2,j_3)\in \{+,-\}^3$ let also
\begin{equation}\label{def_app_sigmatN}
\widetilde{\sigma}^N_{(j_1,j_2,j_3)}(\xi,\eta):=\frac{\sigma^N(\eta,\xi-\eta)}{j_1\langle\xi-\eta\rangle + j_2|\eta| -j_3\langle\xi\rangle}.
\end{equation}
Then for any $\alpha,\beta\in \mathbb{N}^2$
\begin{equation} \label{derivatives sigma tilde N}
\left|\partial^\alpha_\xi\partial^\beta_\eta \widetilde{\sigma}^N_{(j_1,j_2,j_3)}(\xi,\eta)\right|\lesssim_{\alpha,\beta} \langle \xi-\eta\rangle^{2-N+|\alpha|+2|\beta|}|\eta|^{N-|\beta|}, \quad \forall (\xi,\eta)\in\mathbb{R}^2\times\mathbb{R}^2
\end{equation}
and if $N\ge 15$, for any $u,w\in L^2(\mathbb{R}^2)$, $v\in H^{N+3,\infty}(\mathbb{R}^2)$,
\begin{equation} \label{estimate_integral_sigmatildeN}
\left|\int \widetilde{\sigma}^N_{(j_1,j_2,j_3)}(\xi,\eta)\hat{u}(\xi-\eta)\hat{v}(\eta)\hat{w}(-\xi) d\xi d\eta\right|\lesssim \|u\|_{L^2}\|v\|_{H^{N+3,\infty}}\|w\|_{L^2}.
\end{equation}
\proof
From definition \eqref{def_app_sigmatN} function $\widetilde{\sigma}^N_{(j_1,j_2,j_3)}$ can be written as follows
\begin{equation*}
\widetilde{\sigma}^N_{(j_1,j_2,j_3)}(\xi,\eta)=\frac{j_1\langle\xi-\eta\rangle +j_2|\eta| + j_3\langle\xi\rangle}{2j_1j_2 \langle\xi-\eta\rangle |\eta| - 2(\xi-\eta)\cdot\eta}\sigma^N_{(j_1,j_2,j_3)}(\eta,\xi-\eta).
\end{equation*}
We observe that
\[[j_1j_2  \langle\xi-\eta\rangle|\eta| - (\xi-\eta)\cdot\eta]^{-1}\lesssim \langle\xi-\eta\rangle|\eta|^{-1}, \quad \forall (\xi,\eta)\in\mathbb{R}^2\times\mathbb{R}^2\]
and that for any multi-indices $\alpha,\beta\in\mathbb{N}^2$ of positive length
\begin{multline*}
\left|\partial^\alpha_\xi \left[(j_1j_2 \langle\xi-\eta\rangle|\eta| - (\xi-\eta)\cdot\eta)^{-1} \right] \right|\\
\lesssim \sum_{1\le |\alpha_1|\le |\alpha|} \left| j_1j_2 \langle\xi-\eta\rangle|\eta| - (\xi-\eta)\cdot\eta\right|^{-1-|\alpha_1|}|\eta|^{|\alpha_1|}\langle\xi-\eta\rangle^{-(|\alpha| - |\alpha_1|)}, 
\end{multline*}
\begin{multline*}
\left|\partial^\beta_\eta \left[(j_1j_2 \langle\xi-\eta\rangle|\eta| - (\xi-\eta)\cdot\eta)^{-1} \right]\right| \\
\lesssim \sum_{0\le |\beta_1|<|\beta|}\left| j_1j_2 \langle\xi-\eta\rangle|\eta| - (\xi-\eta)\cdot\eta\right|^{-1-(|\beta|-|\beta_1|)}\sum_{\substack{i+j=|\beta| -2|\beta_1|\\ i,j\le |\beta|-|\beta_1|}}\langle\xi-\eta\rangle^i |\eta|^j.
\end{multline*}
From above inequalities we hence deduce that on the support of $\sigma^N_{(j_1,j_2,j_3)}(\eta,\xi-\eta)$ (i.e. for $|\eta|\le \varepsilon |\xi-\eta|$), for any $\alpha,\beta\in\mathbb{N}^2$, 
\begin{equation*}
\left| \partial^\alpha_\xi \partial^\beta_\eta \left[(j_1j_2 \langle\xi-\eta\rangle|\eta| - (\xi-\eta)\cdot\eta)^{-1} \right]\right|\lesssim_{\alpha,\beta} \langle \xi-\eta\rangle^{1+|\alpha| +2|\beta|}|\eta|^{-1-|\beta|},
\end{equation*}
and therefore that 
\small
\begin{equation*}
\left|\partial^\alpha_\xi \partial^\beta_\eta \left[\frac{j_1\langle\xi-\eta\rangle +j_2|\eta| + j_3\langle\xi\rangle}{2j_1j_2 \langle\xi-\eta\rangle |\eta| - 2(\xi-\eta)\cdot\eta}\right] \right| \lesssim_{\alpha,\beta} \langle\xi-\eta\rangle^{2+|\alpha|+2|\beta|}|\eta|^{-1-|\beta|} + \langle\xi\rangle\langle\xi-\eta\rangle^{1+|\alpha|+2|\beta|}|\eta|^{-1-|\beta|}.
\end{equation*} 
\normalsize
The above estimates, summed up with the fact 
\[|\partial^\alpha_\xi\partial^\beta_\eta [\sigma^N_{(j_1,j_2,j_3)}(\eta,\xi-\eta)]|\lesssim_{\alpha,\beta} \langle\xi-\eta\rangle^{-N-|\alpha|}|\eta|^{N+1-|\beta|},\]
gives the first part of the statement.

Let us now suppose that $N\ge 15$ and take $\chi\in C^\infty_0(\mathbb{R}^2)$ equal to 1 in a neighbourhood of the origin.
We have that 
\begin{multline*}
\int \widetilde{\sigma}^N_{(j_1,j_2,j_3)}(\xi,\eta)\hat{u}(\xi-\eta)\hat{v}(\eta)\hat{w}(-\xi) d\xi d\eta = \int K^N_0(z-x,x-y) u(x) v(y)w(z) dxdydz \\
+ \int K_1^N(z-x,x-y) u(x) [\langle D_x\rangle^{N+3}v](y) w(z) dxdydz,
\end{multline*}
with
\begin{equation*}
\begin{gathered}
K^N_k(x,y):= \int e^{ix\cdot\xi + iy\cdot\eta}\widetilde{\sigma}^{N,k}_{(j_1,j_2,j_3)}(\xi,\eta) d\xi d\eta,\\
\widetilde{\sigma}^{N,0}_{(j_1,j_2,j_3)}(\xi,\eta)=\widetilde{\sigma}^N_{(j_1,j_2,j_3)}(\xi,\eta)\chi(\eta) \quad\text{and} \quad \widetilde{\sigma}^{N,1}_{(j_1,j_2,j_3)}(\xi,\eta)=\widetilde{\sigma}^N_{(j_1,j_2,j_3)}(\xi,\eta)\langle \eta\rangle^{-N-3}(1-\chi)(\eta).
\end{gathered}
\end{equation*}
Then inequality \eqref{derivatives sigma tilde N} is obtained using the fact that, for any $\widetilde{u},\widetilde{w}\in L^2, \widetilde{v}\in L^\infty$,
\begin{equation*}
\begin{split}
\int \langle z-x\rangle^{-3}\langle x-y\rangle^{-3} |\widetilde{u}(x)| |\widetilde{v}(y)| |\widetilde{w}(z)| dxdydz &\lesssim \|v\|_{L^\infty}\int \langle z \rangle^{-3} |\widetilde{u}(x)| |\widetilde{w}(z-x)| dx dz \\
& \lesssim \|u\|_{L^2}\|v\|_{L^\infty}\|w\|_{L^2}.
\end{split}
\end{equation*}
\endproof
\end{lem}

In the following lemma we derive some results on the Sobolev continuity of the bilinear integral operator
$$(u,v)\mapsto \int e^{ix\cdot\xi}D_{(j_1,j_2)}(\xi,\eta)\hat{u}(\xi-\eta)\hat{v}(\eta) d\xi d\eta ,$$
with $D_{(j_1,j_2)}$ defined in \eqref{def Dj1j2(xi,eta)}.
We warn the reader that we are not going to take advantage of factor $(1-\frac{\xi-\eta}{\langle\xi-\eta\rangle}\cdot\frac{\eta}{\langle\eta\rangle})$ in $D_{(j_1,j_2)}(\xi,\eta)$ when deriving the estimates mentioned below, since the Sobolev continuity of the above integral operator does not depend on the null structure $Q_0(v,\partial_1v)$ we chose for the Klein-Gordon self-interaction in the wave equation in system \eqref{wave KG system}.

\begin{lem} \label{Lem_appendix: L2 Linfty inequalities integral D Dtilde}
Let $\rho\in \mathbb{N}$ and $D(\xi,\eta)$ a function satisfying, for any multi-indices $\alpha,\beta\in\mathbb{N}^2$, the following:

$(i)$ if $|\xi|\lesssim 1$, 
\begin{gather*}
|\partial^\beta_\eta D(\xi,\eta)|\lesssim_\beta \langle\eta\rangle^{\rho+|\beta|}, \\
|\partial^\alpha_\xi \partial^\beta_\eta D(\xi,\eta)|\lesssim_{\alpha,\beta} \langle\eta\rangle^{\rho+|\alpha|+|\beta|} + \sum_{|\alpha_1|+|\alpha_2|=|\alpha|}|\xi|^{1-|\alpha_1|}\langle\eta\rangle^{\rho+|\alpha_2|+|\beta|}, \quad |\alpha|\ge 1;
\end{gather*}
$(ii)$ for $|\xi| \gtrsim 1$, $|\eta|\lesssim \langle\xi - \eta\rangle$,
\begin{equation*}
|\partial^\alpha_\xi \partial^\beta_\eta D(\xi,\eta)|\lesssim_{\alpha,\beta} \langle\xi-\eta\rangle^{\rho+|\alpha|+|\beta|};
\end{equation*}
$(iii)$ for $|\xi|\gtrsim 1$, $|\eta| \gtrsim \langle\xi-\eta\rangle$:
\begin{equation*}
|\partial^\alpha_\xi \partial^\beta_\eta D(\xi,\eta)|\lesssim_{\alpha,\beta} \langle\eta\rangle^{\rho+|\alpha|+|\beta|}.
\end{equation*}
Then for any $s\ge0$, any $u,v\in H^{s+\rho+13}(\mathbb{R}^2)\cap L^\infty(\mathbb{R}^2)$ (resp. $u,v\in H^{s+\rho+13,\infty}(\mathbb{R}^2)\cap L^2(\mathbb{R}^2)$)
\begin{subequations} \label{est: L2 Linfty integral D(xi,eta)}
\begin{equation} \label{L2 est integral D(xi,eta)}
\begin{split}
\left\|\int e^{ix\cdot\xi} D(\xi,\eta) \hat{u}(\xi-\eta)\hat{v}(\eta) d\xi d\eta \right\|_{H^s(dx)} & \lesssim \|u\|_{H^{s+\rho+13}}\|v\|_{L^\infty} + \|u\|_{L^\infty}\|v\|_{H^{s+\rho+13}}\\
(\text{or } &\lesssim \|u\|_{H^{s+\rho+13,\infty}}\|v\|_{L^2} + \|u\|_{L^2}\|v\|_{H^{s+\rho+13,\infty}}),
\end{split}
\end{equation}
and for any $u,v\in H^{s+\rho+13,\infty}(\mathbb{R}^2)$
\begin{equation} \label{Linfty est integral D(xi,eta)}
\left\|\int e^{ix\cdot\xi} D(\xi,\eta) \hat{u}(\xi-\eta)\hat{v}(\eta) d\xi d\eta \right\|_{H^{s,\infty}(dx)} \lesssim \|u\|_{H^{s+\rho+13,\infty}}\|v\|_{L^\infty} + \|u\|_{L^\infty}\|v\|_{H^{s+\rho+13,\infty}} .
\end{equation}
\end{subequations}
Furthermore, if $\phi \in C^\infty_0(\mathbb{R}^2)$, $t\ge 1$, $\sigma>0$ small, there exists $\delta>0$ depending linearly on $\sigma$, such that
\begin{subequations} \label{est:L2 Linfty integral D with cut-off}
\begin{equation} \label{est:L2 integral D with cut-off}
\begin{split}
\left\|\phi(t^{-\sigma} D_x)\int e^{ix\cdot\xi} D(\xi,\eta) \hat{u}(\xi-\eta)\hat{v}(\eta) d\xi d\eta \right\|_{H^s(dx)} &\lesssim t^\delta \|u\|_{H^{\rho+13}}\|v\|_{L^\infty}  \\
(\text{or } &\lesssim t^\delta \|u\|_{H^{\rho+13,\infty}}\|v\|_{L^2})\\
(\text{or } &\lesssim t^\delta \|u\|_{L^\infty}\|v\|_{H^{\rho+13}}),\\
(\text{or } &\lesssim t^\delta \|u\|_{L^2}\|v\|_{H^{\rho+13,\infty}}),
\end{split}
\end{equation}
\begin{equation} \label{est: Linfty integral D with cut-off}
\begin{split}
\left\|\phi(t^{-\sigma} D_x)\int e^{ix\cdot\xi} D(\xi,\eta) \hat{u}(\xi-\eta)\hat{v}(\eta) d\xi d\eta \right\|_{H^{s,\infty}} &\lesssim t^\delta \|u\|_{H^{\rho+13,\infty}}\|v\|_{L^\infty} \\
(\text{or } &\lesssim t^\delta \|u\|_{L^\infty}\|v\|_{H^{\rho+13,\infty}}).
\end{split}
\end{equation}
\end{subequations}
Finally, if for any $\alpha,\beta\in\mathbb{N}^2$ $D(\xi,\eta)$ satisfies $(ii)$, $(iii)$ when $|\xi|\gtrsim 1$, together with:

$(\widetilde{i})$ if $|\xi|\lesssim 1$
\begin{equation*}
|\partial^\alpha_\xi\partial^\beta_\eta \widetilde{D}(\xi,\eta)|\lesssim_{\alpha,\beta} \langle\eta\rangle^{\rho+|\alpha|+|\beta|} + \sum_{|\alpha_1|+|\alpha_2|=|\alpha|}|\xi|^{-|\alpha_1|+1}\langle\eta\rangle^{\rho+|\alpha_2|+|\beta|},
\end{equation*}
then, for any $u,v\in H^{s+\rho+13}(\mathbb{R}^2)\cap L^\infty(\mathbb{R}^2)$, 
\begin{subequations} \label{est: L2 Linfty integral D-tilde(xi,eta)}\small
\begin{multline} \label{L2 est integral D-tilde(xi,eta)}
\left\|\int e^{ix\cdot\xi} D(\xi,\eta) \hat{u}(\xi-\eta)\hat{v}(\eta) d\xi d\eta \right\|_{H^s(dx)} \lesssim \|u\|_{H^{\rho+10}}\|v\|_{L^2}+ \|u\|_{H^{s+\rho+13}}\|v\|_{L^\infty} + \|u\|_{L^\infty}\|v\|_{H^{s+\rho+13}}\\
(\text{or } \lesssim\|u\|_{L^2}\|v\|_{H^{\rho+10}}+ \|u\|_{H^{s+\rho+13,\infty}}\|v\|_{L^2} + \|u\|_{L^2}\|v\|_{H^{s+\rho+13,\infty}}),
\end{multline}\normalsize
and for any $u,v\in H^{s+\rho+13,\infty}(\mathbb{R}^2)$, with $u\in H^{\rho+10}(\mathbb{R}^2)$ (or $u\in L^2(\mathbb{R}^2)$),
\begin{multline} \label{Linfty est integral D-tilde(xi,eta)}
\left\|\int e^{ix\cdot\xi} D(\xi,\eta) \hat{u}(\xi-\eta)\hat{v}(\eta) d\xi d\eta \right\|_{H^{s,\infty}(dx)} \lesssim\\ \|u\|_{H^{\rho+10}}\|v\|_{L^\infty}+ \|u\|_{H^{s+\rho+13,\infty}}\|v\|_{L^\infty} + \|u\|_{L^\infty}\|v\|_{H^{s+\rho+13,\infty}}\\
(\text{or } \lesssim \|u\|_{L^2}\|v\|_{H^{\rho+10,\infty}}+ \|u\|_{H^{s+\rho+13,\infty}}\|v\|_{L^\infty} + \|u\|_{L^\infty}\|v\|_{H^{s+\rho+13,\infty}}).
\end{multline}
\end{subequations}
\proof
Let $L(\mathbb{R}^2)$ denote either the $L^2(\mathbb{R}^2)$ space or the $L^\infty(\mathbb{R}^2)$ one.
After definition \ref{def Sobolev spaces-NEW} $(i)$ of space $H^s$ (resp. $(iii)$ of $H^{s,\infty}$), we should prove that the $L^2$ norm (resp. the $L^\infty$) norm of
\begin{equation} \label{integral Ds}
\int e^{ix\cdot\xi} D^s(\xi,\eta)\hat{u}(\xi-\eta)\hat{v}(\eta) d\xi d\eta,
\end{equation} 
with $D^s(\xi,\eta):=D(\xi,\eta)\langle\xi\rangle^s$, is bounded by the right hand side of \eqref{L2 est integral D(xi,eta)} and \eqref{est:L2 integral D with cut-off} (resp. \eqref{Linfty est integral D(xi,eta)} and \eqref{est: Linfty integral D with cut-off}).
Let us first take $\chi\in C^\infty_0(\mathbb{R}^2)$ equal to 1 in a neighbourhood of the origin and split the above integral, distinguishing between bounded and unbounded frequencies $\xi$, as
\begin{equation} \label{split integral Ds}
\int e^{ix\cdot\xi} D^s(\xi,\eta)\chi(\xi)\hat{u}(\xi-\eta)\hat{v}(\eta) d\xi d\eta +\int e^{ix\cdot\xi} D^s(\xi,\eta)(1-\chi)(\xi)\hat{u}(\xi-\eta)\hat{v}(\eta) d\xi d\eta.
\end{equation}
On the support of $\chi(\xi)$ frequencies $\xi-\eta, \eta$ are either bounded or equivalent, thus if
\begin{equation*}
a^s_0(\xi,\eta):=
\begin{cases}
& D^s(\xi,\eta)\chi(\xi) \langle \xi-\eta\rangle^{-\rho-10}\\
&\qquad \text{or}\\
& D^s(\xi,\eta)\chi(\xi) \langle \eta\rangle^{-\rho-10}
\end{cases}
\end{equation*} 
$a^s_0(\xi,\eta)$ satisfies \eqref{ineq_a_1} with $g_\beta(\eta)=\langle\eta\rangle^{-3}$ for any $|\beta|\le 3$, after hypothesis $(i)$ on $D(\xi,\eta)$.
Then by \eqref{est:coroll_app_L_norm} and depending on the choice of $a^s_0(\xi,\eta)$, we have that\small
\begin{subequations}\label{est:L integral D chi(xi)}
\begin{multline}\label{est:L integral D with derivatives on u}
\left\|\int e^{ix\cdot\xi} D^s(\xi,\eta)\chi(\xi)\hat{u}(\xi-\eta)\hat{v}(\eta) d\xi d\eta \right\|_{L(dx)} = \left\|\int  e^{ix\cdot\xi}  a^s_0(\xi,\eta)\reallywidehat{\langle D_x\rangle^{\rho+10}u}(\xi-\eta)\hat{v}(\eta)d\xi d\eta  \right\|_{L(dx)}\\
\lesssim \|\langle D_x\rangle^{\rho+10}u\|_{L}\|v\|_{L^\infty} (\text{or } \|\langle D_x\rangle^{\rho+10}u\|_{L^\infty}\|v\|_{L}),
\end{multline}\normalsize
or\small
\begin{multline}\label{est:L integral D with derivatives on v}
\left\|\int e^{ix\cdot\xi} D^s(\xi,\eta)\chi(\xi)\hat{u}(\xi-\eta)\hat{v}(\eta) d\xi d\eta \right\|_{L(dx)} = \left\| \int e^{ix\cdot\xi}  a^s_0(\xi,\eta) \hat{u}(\xi-\eta)\reallywidehat{\langle D_x\rangle^{\rho+10}v}(\eta)d\xi d\eta  \right\|_{L(dx)}\\
\lesssim \|u\|_{L^\infty}\|\langle D_x\rangle^{\rho+10} v\|_{L} (\text{or } \|u\|_{L}\|\langle D_x\rangle^{\rho+10} v\|_{L^\infty}).
\end{multline}
\end{subequations}\normalsize
Successively, we consider a Littlewood-Paley decomposition in order to write
\begin{multline}\label{LP decomposition integral D}
\int e^{ix\cdot\xi}D^s(\xi,\eta)(1-\chi)(\xi) \hat{u}(\xi-\eta)\hat{v}(\eta)d\xi d\eta \\
= \sum_{k\ge 1,l\ge 0} \int e^{ix\cdot\xi}D^s(\xi,\eta)(1-\chi)(\xi)\varphi_k(\xi)\varphi_l(\eta) \hat{u}(\xi-\eta)\hat{v}(\eta)d\xi d\eta,
\end{multline}
where $\varphi_0 \in C^\infty_0(\mathbb{R}^2)$, $\varphi_k(\zeta)=\varphi(2^{-k}\zeta)$ with $\varphi\in C^\infty_0(\mathbb{R}^2\setminus\{0\})$ for any $k\in\mathbb{N}^*$.
When positive indices $l,k$ are such that $k>l+N_0$ for a certain large $N_0\in\mathbb{N}^*$, we have that $|\eta|< |\xi-\eta|$ and $|\xi-\eta|\sim |\xi|\sim 2^k$ on the support of $\varphi_k(\xi)\varphi_l(\eta)$. If 
\[a^s_{k>l+N_0}(\xi,\eta):=D^s(\xi,\eta) \varphi_k(\xi)\varphi_l(\eta)\langle\xi-\eta\rangle^{-s-\rho-13},\]
by hypothesis $(ii)$ we deduce that, for any $\alpha,\beta\in\mathbb{N}^2$ of length less or equal than 3, 
\[|\partial^\alpha_\xi\partial^\beta_\eta[ a^s_{k>l+N_0}(2^k\xi,2^l\eta)]|\lesssim 2^{-k}, \quad\forall (\xi,\eta)\in\mathbb{R}^2\times\mathbb{R}^2\]
and its associated kernel
\begin{equation*}
K^s_{k>l+N_0}(x,y):=\int e^{ix\cdot\xi + iy\cdot\eta} a^s_{k>l+N_0}(\xi,\eta)d\xi d\eta 
\end{equation*}
verifies that 
\[|K^s_{k>l+N_0}(x,y)|\lesssim 2^k 2^{2l}\langle 2^kx\rangle^{-3}\langle 2^ly\rangle^{-3}, \quad \forall (x,y)\in\mathbb{R}^2\times\mathbb{R}^2\]
as one can check doing some integration by parts.
Therefore
\begin{equation} \label{est:L akl k>l}
\begin{split}
&\left\| \int e^{ix\cdot\xi} D^s(\xi,\eta)\varphi_k(\xi)\varphi_l(\eta)\hat{u}(\xi-\eta)\hat{v}(\eta)d\xi d\eta \right\|_{L(dx)}\\
& =  \left\|\int K^s_{k>l+N_0}(x-y,y-z) [\langle D_x\rangle^{s+\rho+13}u](y) v(z) dydz\right\|_{L(dx)} \\
&\lesssim \  2^k 2^{2l} \left\| \int \langle 2^k(x-y)\rangle^{-3}\langle 2^l (y-z)\rangle^{-3} |\langle D_x\rangle^{s+\rho +13}u(y)| |v(z)| dydz \right\|_{L(dx)} \\
&\lesssim  \ 2^k 2^{2l}\int \langle 2^ky\rangle^{-3} \langle 2^l z\rangle^{-3} \| [\langle D_x\rangle^{s+\rho+13} u](\cdot -y)v(\cdot -y-z)\|_{L(dx)} dydz  \\
& \lesssim  2^{-\frac{k}{2}}2^{-\frac{l}{2}} \|\langle D_x\rangle^{s+\rho+13}u\|_{L}\|v\|_{L^\infty} \ (\text{ or }  2^{-\frac{k}{2}}2^{-\frac{l}{2}}\|\langle D_x\rangle^{s+\rho+13} u\|_{L^\infty}\|v\|_L).
\end{split}
\end{equation}
For indices $l,k$ such that $1\le k\le l+N_0$ we have that $|\xi-\eta|\lesssim |\eta|$ on the support of $\varphi_k(\xi)\varphi_l(\eta)$. 
If
\[a^s_{k\le l+N_0}(\xi,\eta):=D^s(\xi,\eta)\varphi_k(\xi)\varphi_l(\eta)\langle\eta\rangle^{-s-\rho-13}\]
by hypothesis $(iii)$ for any multi-indices $\alpha,\beta$ of length less or equal than 3, 
\[|\partial^\alpha_\xi\partial^\beta_\eta [a^s_{k\le l+N_0}(2^k\xi,2^l\eta)]|\lesssim_{\alpha,\beta}2^{-l},\]
and its associated kernel $K^s_{k\le l+N_0}(x,y)$ is such that
\[K^s_{k\le l+N_0}(x,y)|\lesssim 2^{2k}2^l \langle 2^kx\rangle^{-3}\langle 2^ly\rangle^{-3}, \quad \forall (x,y)\in\mathbb{R}^2\times\mathbb{R}^2.\]
Consequently
\begin{equation} \label{est:L akl k<l}
\begin{split}
& \left\| \int e^{ix\cdot\xi}D^s(\xi,\eta)\varphi_k(\xi)\varphi_l(\eta)\hat{u}(\xi-\eta)\hat{v}(\eta)d\xi d\eta \right\|_{L(dx)}\\
& \lesssim  \ 2^{-\frac{k}{2}}2^{-\frac{l}{2}} \|u\|_{L^\infty}\|\langle D_x\rangle^{s+\rho+13}v\|_{L} \ (\text{or }  2^{-\frac{k}{2}}2^{-\frac{l}{2}} \|u\|_{L}\|\langle D_x\rangle^{s+\rho+13}v\|_{L^\infty}),
\end{split}
\end{equation}
and inequality \eqref{L2 est integral D(xi,eta)} (resp. \eqref{Linfty est integral D(xi,eta)}) is hence obtained by combining inequalities \eqref{est:L integral D chi(xi)}, \eqref{est:L akl k>l}, \eqref{est:L akl k<l} with $L=L^2$ (resp. $L=L^\infty$), and taking the sum over $k\ge 1, l\ge 0$.

In order to derive inequalities \eqref{est:L2 Linfty integral D with cut-off}, we first observe that we can reduce to study the $L^2$ and $L^\infty$ norm of \eqref{integral Ds} with $s=0$ and $D(\xi,\eta)$ multiplied by $\phi(t^{-\sigma}\xi)$, up to a factor $t^{s\sigma}$.
Here we use again decompositions \eqref{split integral Ds}, \eqref{LP decomposition integral D}, and only need to modify some of the multipliers defined above, depending on if we want derivatives falling entirely on $u$ or rather on $v$.
In fact, in order to prove the first two inequalities in \eqref{est:L2 integral D with cut-off} and the first one in \eqref{est: Linfty integral D with cut-off} we introduce
\begin{equation*}
\begin{split}
&a^\phi_{l\le k+N_0}(\xi,\eta):=D(\xi,\eta)\chi(t^{-\sigma}\xi)\varphi_k(\xi)\varphi_l(\eta)\\
 &a^\phi_{l>k+N_0}(\xi,\eta):=D(\xi,\eta)\chi(t^{-\sigma}\xi)\varphi_k(\xi)\varphi_l(\eta)\langle\xi-\eta\rangle^{-\rho-13}
\end{split}
\end{equation*} 
and deduce from hypothesis $(ii)-(iii)$ on $D(\xi,\eta)$ and the fact that $|\xi|\lesssim t^\sigma$ on the support of $\phi(t^{-\sigma}\xi)$ that, for any $\alpha,\beta\in \mathbb{N}^2$ of length less or equal than 3,
\[|\partial^\alpha_\xi \partial^\beta_\eta [a^\phi_{l\le k+N_0}(2^k\xi, 2^l\eta)]|\lesssim t^\delta 2^{-k} \quad \text{and} \quad |\partial^\alpha_\xi \partial^\beta_\eta [a^\phi_{ l>k+N_0}(2^k\xi, 2^l\eta)]|\lesssim 2^{-l}\]
with $\delta>0$, $\delta\rightarrow 0$ as $\sigma\rightarrow 0$.
On the one hand, kernel $K^\phi_{l\le k+N_0}(x,y)$ associated to $a^\phi_{l\le k+N_0}(\xi,\eta)$ verifies 
\[|K^\phi_{l\le k+N_0}(x,y)|\lesssim t^\delta 2^{k}2^{2l}\langle2^k x\rangle^{-3}\langle 2^l y\rangle^{-3}, \quad \forall(x, y)\in\mathbb{R}^2\times\mathbb{R}^2\]
and then for any $l,k$ such that $l\le k+N_0$
\begin{equation} \label{integral D with cut-off k bigger l}
\begin{split}
& \left\|\int e^{ix\cdot\xi} D(\xi,\eta)\phi(t^\sigma\xi)\varphi_k(\xi)\varphi_l(\eta)\hat{u}(\xi-\eta)\hat{v}(\eta) d\xi d\eta \right\|_{L(dx)} \\
= & \left\| \int K^\phi_{l\le k+N_0}(x-y,y-z) u(y) v(z) dydz \right\|_{L(dx)} \lesssim t^\delta 2^{-\frac{k}{2}}2^{-\frac{l}{2}} \|u\|_{L}\|v\|_{L^\infty}\\
& \hspace{7cm} (\text{or } \lesssim  t^\delta 2^{-\frac{k}{2}}2^{-\frac{l}{2}}\|u\|_{L^\infty}\|v\|_{L}).
\end{split}
\end{equation}
On the other hand, kernel $K^\phi_{l>k+N_0}(x,y)$
associated to $a^\phi_{l>k+N_0}(\xi,\eta)$ satisfies
\[|K^\phi_{l>k+N_0}(x,y)|\lesssim 2^{2k}2^l\langle2^k x\rangle^{-3}\langle 2^l y\rangle^{-3}, \quad \forall(x,y)\in\mathbb{R}^2\times\mathbb{R}^2\] 
so for indices $l,k$ such that $l>k+N_0$
\begin{equation*}
\begin{split}
& \left\|\int e^{ix\cdot\xi} D(\xi,\eta)\phi(t^\sigma\xi)\varphi_l(\xi)\varphi_l(\eta)\hat{u}(\xi-\eta)\hat{v}(\eta) d\xi d\eta \right\|_{L(dx)} \\
&=  \left\| \int K^\phi_{l>k+N_0}(x-y,y-z) [\langle D_x\rangle^{\rho+13} u](y) v(z) dydz \right\|_{L(dx)} 
\lesssim 2^{-\frac{k}{2}}2^{-\frac{l}{2}} \|\langle D_x\rangle^{\rho+13} u\|_{L}\|v\|_{L^\infty}\\
&\hspace{8.7cm}  (\text{or } \lesssim  2^{-\frac{k}{2}}2^{-\frac{l}{2}}\|\langle D_x\rangle^{\rho+13} u\|_{L^\infty}\|v\|_{L}) \Big).
\end{split}
\end{equation*}
Combining these two inequalities with \eqref{est:L integral D with derivatives on u} and taking the sum over $k\ge 1, l\ge 0$ we obtain the wished estimates.

Last two inequalities in \eqref{est:L2 integral D with cut-off} and last one in \eqref{est: Linfty integral D with cut-off} are instead obtained combining \eqref{est:L integral D with derivatives on v} with \eqref{est:L akl k<l} (that evidently holds for $D^s(\xi,\eta)$ replaced with $D(\xi,\eta)\phi(t^\sigma\xi)$) and \eqref{integral D with cut-off k bigger l}.

Finally, last part of the statement follows from the same argument of above, with the only difference that, after hypothesis $(\widetilde{i})$, multiplier $\widetilde{a}^s_0(\xi,\eta):=\widetilde{D}(\xi,\eta)\chi(\xi)\langle \eta\rangle^{-\rho-10}$ satisfies
\eqref{ineq_a_2} with $|g_\beta(\eta)|\lesssim \langle \eta\rangle^{-3}$ for any $|\beta|\le 3$, then by \eqref{ineq_corA2} we have that
\begin{multline*}
\left\|\int e^{ix\cdot\xi} \widetilde{D}^s(\xi,\eta)\chi(\xi)\hat{u}(\xi-\eta)\hat{v}(\eta) d\xi d\eta \right\|_{L(dx)} \\
= \left\|\int  e^{ix\cdot\xi}  \widetilde{a}^s_0(\xi,\eta)\reallywidehat{\langle D_x\rangle^{\rho+10}u}(\xi-\eta)\hat{v}(\eta)d\xi d\eta  \right\|_{L(dx)}
\lesssim \|\langle D_x\rangle^{\rho+10}u\|_{L^2}\|v\|_{L},
\end{multline*}\normalsize
or\small
\begin{multline*}
\left\|\int e^{ix\cdot\xi} \widetilde{D}^s(\xi,\eta)\chi(\xi)\hat{u}(\xi-\eta)\hat{v}(\eta) d\xi d\eta \right\|_{L(dx)} \\ =\left\| \int e^{ix\cdot\xi}  \widetilde{a}^s_0(\xi,\eta) \hat{u}(\xi-\eta)\reallywidehat{\langle D_x\rangle^{\rho+10}v}(\eta)d\xi d\eta  \right\|_{L(dx)}
\lesssim \|u\|_{L^2}\|\langle D_x\rangle^{\rho+10} v\|_{L} .
\end{multline*}
\endproof
\end{lem}

\begin{lem} \label{Lem_Appendix: est on Dj1j2}
Let $j\in \{+,-\}$, $\phi\in C^\infty_0(\mathbb{R}^2)$, $t\ge 1, \sigma>0$, and $D_j(\xi,\eta)$ be the multiplier introduced in \eqref{def Dj1j2(xi,eta)}.
For any $s\ge 0$, $i=1,2$, $D_j(\xi,\eta)$ and $\frac{\xi_i}{|\xi|}D_j(\xi,\eta)$ satisfy inequalities \eqref{est: L2 Linfty integral D(xi,eta)}, \eqref{est:L2 Linfty integral D with cut-off} with $\rho=2$, and\small
\begin{subequations} \label{est: L2 Linfty integral partial D}
\begin{multline}
\left\|\int e^{ix\cdot\xi} \partial_{\xi}D_j(\xi,\eta) \hat{u}(\xi-\eta)\hat{v}(\eta) d\xi d\eta \right\|_{H^s(dx)}  \lesssim \|u\|_{H^{13}}\|v\|_{L^2}+ \|u\|_{H^{s+16}}\|v\|_{L^\infty} + \|u\|_{L^\infty}\|v\|_{H^{s+16}}\\
(\text{resp. }\lesssim \|u\|_{H^{13}}\|v\|_{L^2}+ \|u\|_{H^{s+16,\infty}}\|v\|_{L^2} + \|u\|_{L^2}\|v\|_{H^{s+16,\infty}}),
\end{multline}
\begin{multline} 
\left\|\int e^{ix\cdot\xi} \partial_{\xi}D_j(\xi,\eta) \hat{u}(\xi-\eta)\hat{v}(\eta) d\xi d\eta \right\|_{H^{s,\infty}(dx)} \\
\lesssim \|u\|_{H^{13}}\|v\|_{L^\infty}+ \|u\|_{H^{s+16,\infty}}\|v\|_{L^\infty} + \|u\|_{L^\infty}\|v\|_{H^{s+16,\infty}},
\end{multline}
\end{subequations}\normalsize
together with
\begin{subequations} \label{est:app_partialD_cutoff}
\begin{multline} \label{est:L2 integral partialD with cut-off}
\left\|\phi(t^{-\sigma}D_x)\int e^{ix\cdot\xi} \partial_{\xi}D_j(\xi,\eta) \hat{u}(\xi-\eta)\hat{v}(\eta) d\xi d\eta \right\|_{H^s(dx)}  \lesssim t^\delta \|u\|_{H^{13}}\left(\|v\|_{L^2}+ \|v\|_{L^\infty}\right)\\
(\text{or }\lesssim t^\delta \|u\|_{L^2}\left(\|v\|_{H^{10}}+ \|v\|_{H^{13,\infty}}\right)),
\end{multline}
\begin{multline} 
\left\|\phi(t^{-\sigma}D_x) \int e^{ix\cdot\xi} \partial_{\xi}D_j(\xi,\eta) \hat{u}(\xi-\eta)\hat{v}(\eta) d\xi d\eta \right\|_{H^{s,\infty}(dx)} 
\lesssim t^\delta\left( \|u\|_{H^{13}}+ \|u\|_{H^{16,\infty}}\right)\|v\|_{L^\infty}\\
(\text{or } \lesssim t^\delta \left(\|u\|_{L^2}+\|u\|_{L^\infty}\right) \|v\|_{H^{16,\infty}}).
\end{multline}
\end{subequations}
Moreover, if $\Omega= x_1 \partial_2 - x_2\partial_1$ and $Z_n=x_n\partial_t +t\partial_n$, $n=1,2,$\small
\begin{subequations} \label{est: L2 integral Omega and Zn cut-off}
\begin{multline} \label{est: L2 Omega integral D with cut-off}
\left\|\phi(t^{-\sigma}D_x)\Omega \int e^{ix\cdot\xi}D_j(\xi,\eta) \hat{u}(\xi-\eta)\hat{v}(\eta)d\xi d\eta\right\|_{L^2(dx)}\\
\lesssim t^\delta \left[ \left(\|u\|_{L^2}+\|\Omega u\|_{L^2}\right)\|v\|_{H^{17,\infty}}+ \|u\|_{H^{15,\infty}}\|\Omega v\|_{L^2}\right],
\end{multline}
\begin{multline} \label{est:L2 Z integral D with cut-off}
\left\|\phi(t^{-\sigma}D_x)Z_n \int e^{ix\cdot\xi}D_j(\xi,\eta) \hat{u}(\xi-\eta)\hat{v}(\eta)d\xi d\eta\right\|_{L^2(dx)}\\
\lesssim t^\delta \left[\|\partial_tu\|_{L^2}\|v\|_{H^{13}}+\|u\|_{H^{13}}\|\partial_tv\|_{L^2}+ \|Z_nu\|_{L^2}\|v\|_{H^{15,\infty}}+\|u\|_{H^{15,\infty}}\| Z_nv\|_{L^2} \right],
\end{multline}
\begin{multline}\label{est: L2 DjZ integral D with cut-off}
\left\|\phi(t^{-\sigma}D_x)D_j Z_n \int e^{ix\cdot\xi}D_j(\xi,\eta) \hat{u}(\xi-\eta)\hat{v}(\eta)d\xi d\eta\right\|_{L^2(dx)}\\
\lesssim t^\delta \left[\|\partial_tu\|_{L^2}\|v\|_{H^{14,\infty}}+\|u\|_{H^{14,\infty}}\|\partial_tv\|_{L^2}+ \|Z_nu\|_{L^2}\|v\|_{H^{17,\infty}}+\|u\|_{H^{17,\infty}}\| Z_nv\|_{L^2} \right],
\end{multline}
\end{subequations}\normalsize
with $\delta>0$, $\delta\rightarrow 0$ as $\sigma\rightarrow 0$.
\proof
The statement follows essentially from the observation that, for $j\in \{+,-\}$, functions $D_j(\xi,\eta)$ and $[(\xi_i\partial_{\xi_j})^{k_1}(\eta_i\partial_{\eta_j})^{k_2}D_j](\xi,\eta)$ satisfy hypothesis $(i)-(iii)$ of lemma \ref{Lem_appendix: L2 Linfty inequalities integral D Dtilde} with $\rho=2$ and $\rho=2+2(k_1+k_2)$ respectively, while $\partial_\xi D_j(\xi,\eta)$ satisfies $(\widetilde{i}),(ii),(iii)$ with $\rho=3$.
In fact, we first remark that, for every $\xi, \eta$, denominator $1+\langle\xi-\eta\rangle\langle\eta\rangle - (\xi-\eta)\cdot\eta$ is bounded from below by a positive constant;
secondly, the derivation of that denominator gives rise to losses in $\langle\xi-\eta\rangle, \langle\eta\rangle$, as
\begin{align*}
\partial_{\xi_k}(1+\langle\xi-\eta\rangle\langle\eta\rangle - (\xi-\eta)\cdot\eta) &= \frac{\xi_k-\eta_k}{\langle\xi-\eta\rangle}\langle\eta\rangle + \eta_k, \\
\partial_{\eta_k}(1+\langle\xi-\eta\rangle\langle\eta\rangle - (\xi-\eta)\cdot\eta) & =  \frac{\xi_k-\eta_k}{\langle\xi-\eta\rangle}\langle\eta\rangle + \langle\xi-\eta\rangle\frac{\eta_k}{\langle\eta\rangle} + \eta_k - (\xi_k-\eta_k).
\end{align*}
For $|\xi|\lesssim 1$ we have that $\langle\xi-\eta\rangle \lesssim \langle\eta\rangle$, so for any $\alpha,\beta\in\mathbb{N}^2$
\begin{equation*}
\left|\partial^\alpha_\xi \partial^\beta_\eta \left[\frac{j\langle\xi-\eta\rangle + j\langle\eta\rangle }{1+\langle\xi-\eta\rangle\langle\eta\rangle - (\xi-\eta)\cdot\eta} \eta_1 \right]\right| \lesssim_{\alpha,\beta} \langle\eta\rangle^{2+|\alpha|+|\beta|},
\end{equation*}
while
\begin{align*}
\left|\partial^\beta_\eta \left[\frac{|\xi|}{1+\langle\xi-\eta\rangle\langle\eta\rangle - (\xi-\eta)\cdot\eta} \eta_1 \right]\right| &\lesssim_\beta \langle\eta\rangle^{1+ |\beta|}, \\
\left|\partial^\alpha_\xi\partial^\beta_\eta \left[\frac{|\xi|}{1+\langle\xi-\eta\rangle\langle\eta\rangle - (\xi-\eta)\cdot\eta}  \eta_1\right]\right| &\lesssim_{\alpha,\beta} \sum_{|\alpha_1|+|\alpha_2|=|\alpha|} |\xi|^{1-|\alpha_1|}\langle\eta\rangle^{1+ |\alpha_2|+|\beta|}, \quad |\alpha|\ge 1.
\end{align*}
For $|\xi|\gtrsim1$ and $|\eta|\lesssim \langle \xi-\eta\rangle$ (resp. $|\eta|\gtrsim \langle\xi -\eta\rangle$) we have that $|\xi|\lesssim |\xi-\eta|$ (resp.$|\xi|\lesssim |\eta|$), so each time a derivative hits the denominator of $D_j(\xi,\eta)$ we lose a factor $\langle\xi-\eta\rangle$ (resp. $\langle\eta\rangle$).
Hence lemma \ref{Lem_appendix: L2 Linfty inequalities integral D Dtilde} immediately implies inequalities \eqref{est: L2 Linfty integral D(xi,eta)}, \eqref{est:L2 Linfty integral D with cut-off} with $D= D_j$ and $\rho=2$, together with \eqref{est: L2 Linfty integral partial D}, \eqref{est:app_partialD_cutoff}, while inequalities \eqref{est: L2 integral Omega and Zn cut-off} follow from \eqref{est:L2 Linfty integral D with cut-off} and the fact that, after some integration by parts,
\begin{equation*} 
\begin{split}
\Omega & \int  e^{ix\cdot\xi}D_j(\xi,\eta)\hat{u}(\xi-\eta)\hat{v}(\eta) d\xi d\eta  \\
&=\sum_{k_1+k_2+k_3+k_4=1}\int  e^{ix\cdot\xi} [(\xi_1 \partial_{\xi_2} - \xi_2 \partial_{\xi_1})^{k_1}(\eta_1 \partial_{\eta_2} - \eta_2 \partial_{\eta_1})^{k_2}D_j](\xi,\eta)\widehat{\Omega^{k_3}u}(\xi-\eta)\reallywidehat{\Omega^{k_4}v}(\eta) d\xi d\eta ,\\
Z_n & \int  e^{ix\cdot\xi}D_j(\xi,\eta)\hat{u}(\xi-\eta)\hat{v}(\eta) d\xi d\eta \\
&=  \int  e^{ix\cdot\xi}[\partial_{\xi_n}D_j](\xi,\eta)D_t\big[\hat{u}(\xi-\eta)\hat{v}(\eta)\big] d\xi d\eta + \int  e^{ix\cdot\xi}[\partial_{\eta_n}D_j](\xi,\eta)\hat{u}(\xi-\eta)\widehat{D_t v}(\eta) d\xi d\eta  \\
& +  \int  e^{ix\cdot\xi}D_j(\xi,\eta)\widehat{Z_n u}(\xi-\eta)\hat{v}(\eta) d\xi d\eta +  \int  e^{ix\cdot\xi}D_j(\xi,\eta)\hat{u}(\xi-\eta)\widehat{Z_nv}(\eta) d\xi d\eta,
\end{split}
\end{equation*}
and, if $\delta_{jn}$ denotes the Kronecker delta,
\begin{equation*}
\begin{split}
D_j & Z_n  \int  e^{ix\cdot\xi}D_j(\xi,\eta)\hat{u}(\xi-\eta)\hat{v}(\eta) d\xi d\eta \\
& = \delta_{jn} \int e^{ix\cdot\xi}D_j(\xi,\eta) D_t\left[\hat{u}(\xi-\eta)\hat{v}(\eta)\right] d\xi d\eta \\
& + \int  e^{ix\cdot\xi}\partial_{\xi_n} [\xi_j D_j](\xi,\eta)D_t\big[\hat{u}(\xi-\eta)\hat{v}(\eta)\big] d\xi d\eta + \int  e^{ix\cdot\xi}\partial_{\eta_n}[\xi_j D_j](\xi,\eta)\hat{u}(\xi-\eta)\widehat{D_t v}(\eta) d\xi d\eta  \\
& +  \int  e^{ix\cdot\xi} \xi_j D_j(\xi,\eta)\widehat{Z_n u}(\xi-\eta)\hat{v}(\eta) d\xi d\eta +  \int  e^{ix\cdot\xi}\xi_j D_j(\xi,\eta)\hat{u}(\xi-\eta)\widehat{Z_nv}(\eta) d\xi d\eta.
\end{split}
\end{equation*}.
\endproof
\end{lem}

\chapter[Appendix B]{} \label{Appendix B}

\numberwithin{equation}{section}
\numberwithin{thm}{section}
\renewcommand{\thesection}{\thechapter.\arabic{section}}

The aim of this chapter is to show how, from the bootstrap assumptions \eqref{est: bootstrap argument a-priori est}, it is possible to derive a moderate growth in time for the $L^2$ norm of $\mathcal{L}^\mu\widetilde{v}$, with $0\le |\mu|\le 2$, and of $\Omega^\mu_h\mathcal{M}^\nu\widetilde{u}^{\Sigma,k}$, with $\mu, |\nu|=0,1$.
These estimates are fundamentally used in propositions \ref{Prop:propagation_unif_est_V} and \ref{Prop: Propagation uniform estimate U,RU}.
Moreover, we also prove in lemma \ref{Lem_appendix: sharp_est_VJ} a sharp decay estimate for the uniform norm of the Klein-Gordon solution when one Klainerman vector field is acting on it (and when considered for frequencies less or equal than $t^\sigma$, with $\sigma>0$ small).
We are hence going to assume for the rest of this chapter that a-priori estimates \eqref{est: bootstrap argument a-priori est} are satisfied in interval $[1,T]$, for some fixed $T>1$, and that $\varepsilon_0<(2A+B)^{-1}$.
We remind that $\Gamma$ generally denotes one of the admissible vector fields belonging to $\mathcal{Z}$ (see \eqref{order_Z}) and that, for a multi-index $I=(i_1,\dots,i_n)$ with $i_j\in \{1,\dots,5\}$ for $j=1,\dots, n$, $\Gamma^I = \Gamma_{i_1}\cdots \Gamma_{i_n}$. Also, we warn the reader that any norm $X$ ($X=L^\infty, H^s, H^s_h...$) of $w=w(t,x)$ is here considered with respect to spatial variable $x$. We will often write $\|\cdot\|_{X}$ in place of $\|\cdot\|_{X(dx)}$. 

\section{Some preliminary lemmas}\label{Sub: App_B1}

In the current section we list, on the one hand, some inequalities concerning the $H^s$ and $H^{s,\infty}$ norm of the quadratic non-linearities $Q^\mathrm{w}_0(v_\pm, D_1v_\pm), Q^\mathrm{kg}_0(v_\pm, D_1u_\pm)$ (see lemmas \ref{Lem: estimates NLkg NLw}, \ref{Lem: est DtU DtV}), as they are very frequently recalled in the second part of the paper. 
On the other hand, we introduce some preliminary small results that will be useful in sections \ref{sec_appB: first range of estimates} and \ref{sec_appB: second range of estimates}.

For seek of compactness, we denote $Q^\mathrm{w}_0(v_\pm, D_1v_\pm)$ and $Q^\mathrm{kg}_0(v_\pm, D_1 u_\pm)$ by $\textit{NL}_w$ and $\textit{NL}_{kg}$ respectively, i.e. \index{Nlw@$\Nlw$, quadratic non-linearity in the wave equation satisfied by $u_\pm$}\index{Nlkg@$\Nlkg$, quadratic non-linearity in the Klein-Gordon equation satisfied by $v_\pm$}
\begin{subequations}
\begin{align}
\Nlw & := \frac{i}{4}\left[(v_+ + v_{-})D_1(v_+ + v_{-}) - \frac{D_x}{\langle D_x \rangle}(v_+ - v_{-})\cdot\frac{D_x D_1}{\langle D_x \rangle}(v_+ - v_{-})\right] ,\\
\Nlkg & := \frac{i}{4} \left[(v_+ + v_{-})D_1(u_+ + u_{-}) - \frac{D_x}{\langle D_x\rangle}(v_+ - v_{-})\cdot\frac{D_x D_1}{|D_x|}(u_+ - u_{-})\right]. \label{def_app_Nlkg}
\end{align}
\end{subequations}
We recall the result of lemma \ref{Lem : new estimate 1-chi}, that can be also stated in the classical setting and says that, for any real positive $s>s'$ and $w\in H^s(\mathbb{R}^2)$,
\begin{equation} \label{ineq:1-chi}
\|(1-\chi)(t^{-\sigma}D_x)w\|_{H^{s'}}\le C t^{-\sigma(s-s')}\|w\|_{H^s}, \quad \forall s>s'.
\end{equation}
It is also useful to remind, in view of upcoming lemmas, that the $L^2$ norm of $(\Gamma^Iu)_\pm$ and $(\Gamma^Iv)_\pm$ is estimated with:
\begin{equation*}
\begin{split}
E_n(t;W)^\frac{1}{2},  &\ \text{whenever }|I|\le n \text{ and }\Gamma^I \text{ is a product of spatial derivatives;}\\
E^k_3(t;W)^\frac{1}{2}, &\ \text{whenever }|I|\le 3 \text{ and at most }3-k \text{ vector fields in }\Gamma^I \text{, with } 0\le k\le 2 \\
&\text{ belong to }\{\Omega, Z_m, m=1,2\}.
\end{split}
\end{equation*}
As assumed in \eqref{est: bootstrap Enn}, \eqref{est: bootstrap E02}, such energies have a moderate growth in time and a hierarchy is established among them in the sense that
\[0<\delta\ll \delta_2\ll \delta_1\ll \delta_0\ll 1.\]
We warn the reader that this hierarchy is often implicitly used throughout this chapter.

\begin{lem}  \label{Lem: estimates NLkg NLw}
For any $s\ge 0$, any $\theta\in ]0,1[$, $\textit{NL}_w$ satisfies the following inequalities:
\begin{subequations} \label{est Hs Hsinfty Nlw-new}
\begin{gather}
\|\textit{NL}_w(t,\cdot)\|_{L^2} \lesssim \|V(t,\cdot)\|_{H^{1,\infty}}\|V(t,\cdot)\|_{H^1}, \label{est L2 NLw} \\
\|\textit{NL}_w(t,\cdot)\|_{L^\infty} \lesssim \|V(t,\cdot)\|^2_{H^{2,\infty}},\label{est Linfty NLw} \\
\|\textit{NL}_w(t,\cdot)\|_{H^s} \lesssim \|V(t,\cdot)\|_{H^{s+1}}\|V(t,\cdot)\|_{H^{1,\infty}},  \label{est Hs NLw-New}\\
\|\textit{NL}_w(t,\cdot)\|_{H^{s,\infty}} \lesssim \|V(t,\cdot)\|_{H^{s+1,\infty}}^{2-\theta}\|V(t,\cdot)\|^\theta_{H^{s+3}}, \label{est Hsinfty for NLw-New}\\
\|\Omega \textit{NL}_w(t,\cdot)\|_{L^2}\lesssim \|V(t,\cdot)\|_{H^{2,\infty}}\left(\|V(t,\cdot)\|_{L^2}+ \|\Omega V(t,\cdot)\|_{H^1}\right), \label{est_Omega_NLw}
\end{gather}
\end{subequations}
while for $\textit{NL}_{kg}$ we have that:

\begin{subequations}\label{est on NLkg-new}
\begin{gather}
\|\textit{NL}_{kg}(t,\cdot)\|_{L^2} \lesssim \|V(t,\cdot)\|_{H^{1,\infty}}\|U(t,\cdot)\|_{H^1},  \label{est L2 NLkg} \\
\|\textit{NL}_{kg}(t,\cdot)\|_{L^\infty} \lesssim \|V(t,\cdot)\|_{H^{1,\infty}}\left(\|U(t,\cdot)\|_{H^{2,\infty}}+ \|\mathrm{R}_1U(t,\cdot)\|_{H^{2,\infty}}\right), \label{est Linfty NLkg}\\
\|\textit{NL}_{kg}(t,\cdot)\|_{H^s}  \lesssim \|V(t,\cdot)\|_{H^{s}}\left(\|U(t,\cdot)\|_{H^{1,\infty}}+ \|\mathrm{R}_1U(t,\cdot)\|_{H^{1,\infty}}\right) + \|V(t,\cdot)\|_{L^\infty}\|U(t,\cdot)\|_{H^{s+1}}, \label{est Hs NLkg-New}\\
\end{gather}
\begin{equation} \label{est Hsinfty NLkg-new}
\begin{split}
\|\textit{NL}_{kg}(t,\cdot)\|_{H^{s,\infty}}&\lesssim \|V(t,\cdot)\|^{1-\theta}_{H^{s,\infty}} \|V(t,\cdot)\|^\theta_{H^{s+2}}\left(\| U(t,\cdot)\|_{H^{1,\infty}}+ \|\mathrm{R}_1U(t,\cdot)\|_{H^{1,\infty}}\right)\\
&+ \|V(t,\cdot)\|_{L^\infty}\left(\| U(t,\cdot)\|^{1-\theta}_{H^{s+1,\infty}}+ \|\mathrm{R}_1U(t,\cdot)\|^{1-\theta}_{H^{s+1,\infty}}\right)\|U(t,\cdot)\|^\theta_{H^{s+3}},
\end{split}
\end{equation}
\begin{equation} \label{est L2 Omega NLkg}
\begin{split}
\|\Omega\textit{NL}_{kg}(t,\cdot)\|_{L^2}&\lesssim \left(\|V(t,\cdot)\|_{L^2}+\|\Omega V(t,\cdot)\|_{L^2}\right)\left(\|U(t,\cdot)\|_{H^{2,\infty}}+ \|\mathrm{R}_1U(t,\cdot)\|_{H^{2,\infty}}\right) \\
&+ \|V(t,\cdot)\|_{H^{1,\infty}}\|\Omega U(t,\cdot)\|_{H^1}.
\end{split}
\end{equation}
\end{subequations}
\proof
Inequalities \eqref{est L2 NLw}, \eqref{est Linfty NLw}, \eqref{est L2 NLkg}, and \eqref{est Linfty NLkg} are straightforward. The same is for \eqref{est_Omega_NLw} and \eqref{est L2 Omega NLkg} after commutation of $\Omega$ with the operators appearing in \eqref{Q0_pm}.
All other inequalities in the statement are rather derived using corollary \ref{Cor_appendixA:Hs-Hsinfty norm of bilinear expressions}.
\endproof
\end{lem}

\begin{lem} \label{Lem: est DtU DtV}
For any $s\ge 0$, any $\theta\in ]0,1[$,
\begin{subequations} 
\begin{gather}
\|D_tU(t,\cdot)\|_{H^s}\lesssim \|U(t,\cdot)\|_{H^{s+1}} + \|V(t,\cdot)\|_{H^{s+1}}\|V(t,\cdot)\|_{H^{1,\infty}},  \label{Hs_norm_DtU}\\
\|D_tU(t,\cdot)\|_{H^{s,\infty}}\lesssim \|U(t,\cdot)\|_{H^{s+2,\infty}} + \|V(t,\cdot)\|^{2-\theta}_{H^{s+1,\infty}}\|V(t,\cdot)\|^\theta_{H^{s+3}}, \label{Hsinfty_norm_DtU} \\
\|D_t \mathrm{R}_1U(t,\cdot)\|_{H^{s,\infty}}\lesssim \|\mathrm{R}_1 U(t,\cdot)\|_{H^{s+1,\infty}} + \|V(t,\cdot)\|_{H^{s+3}}\|V(t,\cdot)\|_{H^{1,\infty}},  \label{Hsinfty norm Dt R1U-new}\\
\|D_t\Omega U(t,\cdot)\|_{L^2}\le \|\Omega U(t,\cdot)\|_{H^1}+ \|V(t,\cdot)\|_{H^{2,\infty}}\left(\|V(t,\cdot)\|_{L^2}+ \|\Omega V(t,\cdot)\|_{H^1}\right), \label{L2_norm_DtOmegaU}
\end{gather}
\end{subequations}
while
\begin{subequations}
\begin{equation} \label{Hs norm DtV}
\begin{split}
\|D_tV(t,\cdot)\|_{H^s} &\lesssim \|V(t,\cdot)\|_{H^{s+1}} + \|V(t,\cdot)\|_{H^{s}}\left(\|U(t,\cdot)\|_{H^{1,\infty}}+ \|\mathrm{R}_1U(t,\cdot)\|_{H^{1,\infty}}\right)\\
& + \|V(t,\cdot)\|_{L^\infty}\|U(t,\cdot)\|_{H^{s+1}},
\end{split}
\end{equation}\small
\begin{equation}\label{est: Hsinfty Dt V}
\begin{split}
\|D_tV(t,\cdot)\|_{H^{s,\infty}}& \lesssim \|V(t,\cdot)\|_{H^{s+1,\infty}} + \|V(t,\cdot)\|^{1-\theta}_{H^{s,\infty}} \|V(t,\cdot)\|^\theta_{H^{s+1}}\left(\| U(t,\cdot)\|_{H^{1,\infty}}+ \|\mathrm{R}_1U(t,\cdot)\|_{H^{1,\infty}}\right)\\
&+ \|V(t,\cdot)\|_{L^\infty}\left(\| U(t,\cdot)\|^{1-\theta}_{H^{s+1,\infty}}+ \|\mathrm{R}_1U(t,\cdot)\|^{1-\theta}_{H^{s+1,\infty}}\right)\|U(t,\cdot)\|^\theta_{H^{s+3}},
\end{split}
\end{equation}\normalsize
\begin{equation} \label{L2_norm_DtOmegaV}
\begin{split}
\| D_t\Omega V(t,\cdot)\|_{L^2}& \le \|\Omega V(t,\cdot)\|_{H^1}+  \left(\|V(t,\cdot)\|_{L^2}+\|\Omega V(t,\cdot)\|_{L^2}\right)\left(\|U(t,\cdot)\|_{H^{2,\infty}}+ \|\mathrm{R}_1U(t,\cdot)\|_{H^{2,\infty}}\right) \\
&+ \|V(t,\cdot)\|_{H^{1,\infty}}\|\Omega U(t,\cdot)\|_{H^1}.
\end{split}
\end{equation}
\end{subequations}
\proof
Straight consequence of the previous lemma and the fact that $(u_+, v_+, u_{-},v_{-})$ is solution to system \eqref{wave-KG for u- v-}. Observe that inequality \eqref{Hsinfty norm Dt R1U-new} is derived using that
\begin{equation*}
\|\mathrm{R}_1 \textit{NL}_w(t,\cdot)\|_{H^{s,\infty}}\lesssim \|\textit{NL}_w(t,\cdot)\|_{H^{s+2}}
\end{equation*}
after classical Sobolev injection and continuity of $\mathrm{R}_1: H^s\rightarrow H^s$, for any $s\ge 0$.
\endproof
\end{lem}

\begin{lem}
Let $|I|=1$ be such that $\Gamma^I\in\{\Omega, Z_m, m=1,2\}$.
Then
\begin{multline} \label{DtUI}
\| D_t U^I(t,\cdot)\|_{L^2}\lesssim \|U^I(t,\cdot)\|_{H^1} + \|V(t,\cdot)\|_{H^{2,\infty}}\Big[\|V^I(t,\cdot)\|_{H^1} \\
+ \|V(t,\cdot)\|_{H^1}\Big(1+\sum_{\mu=0}^1\|\mathrm{R}_1^\mu U(t,\cdot))\|_{H^{1,\infty}}\Big) +\|V(t,\cdot)\|_{L^\infty} \|U(t,\cdot)\|_{H^1}\Big],
\end{multline}
\begin{multline}\label{DtVI}
\| D_t V^I(t,\cdot)\|_{L^2}\lesssim \|V^I(t,\cdot)\|_{H^1} + \sum_{\mu=0}^1\|\mathrm{R}_1^\mu U(t,\cdot)\|_{H^{2,\infty}}\|V^I(t,\cdot)\|_{L^2} \\
+ \|V(t,\cdot)\|_{H^{1,\infty}}\left(\|U^I(t,\cdot)\|_{H^1}+\|U(t,\cdot)\|_{H^1}+\|V(t,\cdot)\|_{H^{1,\infty}}\|V(t,\cdot)\|_{H^1}\right).
\end{multline}
\proof
The result of the statement follows using the equation satisfied, respectively, by $u^I_\pm$ and $v^I_\pm$, together with \eqref{Hs_norm_DtU}, \eqref{Hs norm DtV} with $s=0$.
In fact, by \eqref{Gamma_nonlinearity} with $|I|=1$,
\begin{gather*}
D_t u^I_\pm =\pm |D_x| u^I_\pm + Q^\mathrm{w}_0(v^I_\pm, D_1 v_\pm) + Q^\mathrm{w}_0(v_\pm, D_1 v^I_\pm) + G^\mathrm{w}_1(v_\pm,D v_\pm), \\
D_t v^I_\pm =\pm \langle D_x\rangle v^I_\pm + Q^\mathrm{kg}_0(v^I_\pm, D_1 u_\pm) + Q^\mathrm{kg}_0(v_\pm, D_1 u^I_\pm) + G^\mathrm{kg}_1(v_\pm,D u_\pm),
\end{gather*}
with $G^\mathrm{w}_1(v_\pm,\partial v_\pm)=G_1(v, \partial v)$, $G^\mathrm{kg}_1(v_\pm,D u_\pm)= G_1(v, \partial u)$ and $G_1$ given by \eqref{def_G1}.
Hence one can estimate the $L^2$ norm of the first two quadratic terms in above equalities with the $L^2$ norm of factors indexed in $I$ times the $L^\infty$ norm of the remaining one, while the $L^2$ norm of the latter quadratic terms can be instead bounded by taking the $L^2$ norm of one of the two factors times the $L^\infty$ norm of the remaining one, indifferently. We choose here to consider the $L^2$ norm of factors $Du_\pm, Dv_\pm$, and use \eqref{Hs_norm_DtU}, \eqref{Hs norm DtV} if the derivative $D$ is a time derivative.
\endproof
\end{lem}

It is useful to remind that, if $w(t,x)$ is solution to inhomogeneous half wave equation \eqref{half-wave-w} from \eqref{relation_w_wave_Zjw-new} we have that for any $j,k\in\{1,2\}$ and $|\mu|\le 1$
\begin{subequations}\label{equalities_app_xw}
\begin{equation} \label{equality_xDkw}
\begin{split}
 x_j D_k\Big(\frac{D_x}{|D_x|}\Big)^\mu w=& \frac{D_k}{|D_x|}\Big(\frac{D_x}{|D_x|}\Big)^\mu\left[x_j|D_x| -tD_j +\frac{1}{2i}\frac{D_j}{|D_x|}\right]w+t\frac{D_jD_k}{|D_x|}\Big(\frac{D_x}{|D_x|}\Big)^\mu w\\
&  -\frac{1}{2i}\frac{D_jD_k}{|D_x|^2}\Big(\frac{D_x}{|D_x|}\Big)^\mu w + i Op\Big(\partial_j\Big(\frac{\xi_k}{|\xi|}\Big(\frac{\xi}{|\xi|}\Big)^\mu\Big)|\xi|\Big)w\\
=&\,  i\frac{D_k}{|D_x|}\Big(\frac{D_x}{|D_x|}\Big)^\mu Z_jw + \frac{D_k}{|D_x|}\Big(\frac{D_x}{|D_x|}\Big)^\mu [x_jf(t,x)]\\
&+  t\frac{D_jD_k}{|D_x|}\Big(\frac{D_x}{|D_x|}\Big)^\mu w+ iOp\Big(\partial_j\Big(\frac{\xi_k}{|\xi|}\Big(\frac{\xi}{|\xi|}\Big)^\mu\Big)|\xi|\Big)w.
\end{split}
\end{equation} 
Analogously, if $w(t,x)$ is solution to inhomogeneous half Klein-Gordon \eqref{half KG}, from \eqref{relation_w_KG_Zjw} we have that
\begin{equation} \label{xjw_Zjw}
\begin{split}
x_j \Big(\frac{D_x}{\langle D_x\rangle}\Big)^\mu w &= \frac{1}{\langle D_x\rangle}\Big(\frac{D_x}{\langle D_x\rangle}\Big)^\mu\left[\langle D_x\rangle x_j -tD_j\right] w+ t \frac{D_j}{\langle D_x\rangle}\Big(\frac{D_x}{\langle D_x\rangle}\Big)^\mu w + i\oph\Big(\partial_j\Big(\frac{\xi}{\langle\xi\rangle}\Big)^\mu\Big)w \\
&= i\frac{1}{\langle D_x\rangle}\Big(\frac{D_x}{\langle D_x\rangle}\Big)^\mu Z_jw -i \frac{D_j}{\langle D_x\rangle^2}\Big(\frac{D_x}{\langle D_x\rangle}\Big)^\mu w + \frac{1}{\langle D_x\rangle}\Big(\frac{D_x}{\langle D_x\rangle}\Big)^\mu[x_j f(t,x)]\\
& +  t \frac{D_j}{\langle D_x\rangle}\Big(\frac{D_x}{\langle D_x\rangle}\Big)^\mu w+ i\oph\Big(\partial_j\Big(\frac{\xi}{\langle\xi\rangle}\Big)^\mu\Big)w.
\end{split}
\end{equation} 
\end{subequations}
We also remind the reader about equivalence \eqref{equivalence GammaIW-Ekn}, so we won't particularly care if we are dealing with $\Gamma^I u_\pm, \Gamma^I v_\pm$ instead of $(\Gamma^I u)_\pm, (\Gamma^I v)_\pm$, when we bound the $L^2$ norm of those terms with the energy defined in \eqref{def_generalized_energy}.

\begin{lem}\label{Lem_appendix: Technical_estimates}
There exists a positive constant $C>0$ such that, for every $j=1,2$, $t\in [1,T]$,
\begin{subequations} \label{norms_H1_Linfty_xv-}
\begin{gather}
\sum_{|\mu|=0}^1 \left\| x_j \Big(\frac{D_x}{\langle D_x\rangle}\Big)^\mu v_\pm(t,\cdot)\right\|_{H^1}\le CB\varepsilon t^{1+{\frac{\delta}{2}}}, \label{norm_xv-} \\
\sum_{|\mu|=0}^1 \left\| x_j \Big(\frac{D_x}{\langle D_x\rangle}\Big)^\mu v_\pm(t,\cdot)\right\|_{H^{1,\infty}}\le C(A+B)\varepsilon t^{\frac{\delta_2}{2}}, \label{norm_Linfty_xv-}
\end{gather}
\end{subequations}
and
\begin{equation} \label{norm_xu-}
\sum_{|\mu|=0}^1\left\|x_jD_x\Big(\frac{D_x}{|D_x|}\Big)^\mu u_\pm(t,\cdot)\right\|_{L^2}\le CB\varepsilon t^{1+\frac{\delta}{2}}.
\end{equation}
\proof
We warn the reader that, throughout the proof, $C$ will denote a positive constant that may change line after line.
As $v_+=-\overline{v_{-}}$ (resp. $u_+=-\overline{u_{-}}$), it is enough to prove the statement for $v_{-}$ (resp. for $u_{-}$).

Since $v_{-}$ is solution to equation \eqref{half KG} with $f=\Nlkg$, from \eqref{xjw_Zjw} it immediately follows that, for any $|\mu|\le 1$,
\begin{subequations}
\begin{equation} \label{xjv_H1}
\left\|x_j\Big(\frac{D_x}{\langle D_x\rangle}\Big)^\mu v_{-}(t,\cdot)\right\|_{H^1}\lesssim \|Z_jv_{-}(t,\cdot)\|_{L^2}+ t\|v_{-}(t,\cdot)\|_{H^1}+ \|x_j\textit{NL}_{kg}(t,\cdot)\|_{L^2(dx)}
\end{equation}
along with
\begin{equation} \label{xjv_Linfty}
\left\|x_j\Big(\frac{D_x}{\langle D_x\rangle}\Big)^\mu v_{-}(t,\cdot)\right\|_{H^{1,\infty}}\le \|Z_jv_{-}(t,\cdot)\|_{H^2}+ t\|v_{-}(t,\cdot)\|_{H^{2,\infty}}+ \|x_j\textit{NL}_{kg}(t,\cdot)\|_{L^\infty(dx)},
\end{equation}
\end{subequations}
derived by using the classical Sobolev injection. 
Observe that
\begin{subequations}
\begin{equation} \label{norm_Linfty_xNLkg}
\left\| x_j\textit{NL}_{kg}(t,\cdot)\right\|_{L^\infty}  \lesssim \left(\|x_jv_{-}(t,\cdot)\|_{L^\infty}+ \left\|x_j\frac{D_x}{\langle D_x\rangle}v_{-}(t,\cdot) \right\|_{L^\infty} \right) \sum_{\mu=0}^1\|\mathrm{R}^\mu_1U(t,\cdot)\|_{H^{2,\infty}},
\end{equation}
but also
\begin{equation} \label{norm_L2_xNLkg}
\|x_j\textit{NL}_{kg}(t,\cdot)\|_{L^2}\lesssim \sum_{|\mu|=0}^1 \left\| x_j \Big(\frac{D_x}{\langle D_x\rangle}\Big)^\mu v_\pm(t,\cdot)\right\|_{L^2}\left(\|U(t,\cdot)\|_{H^{2,\infty}}+\|\mathrm{R}_1U(t,\cdot)\|_{H^{2,\infty}}\right).
\end{equation}
\end{subequations}
Thus, if $\varepsilon_0>0$ is assumed sufficiently small to verify $\varepsilon_0<(2A)^{-1}$, 
by injecting \eqref{norm_L2_xNLkg} (resp. \eqref{norm_Linfty_xNLkg}) into \eqref{xjv_H1} (resp. in \eqref{xjv_Linfty}), and using a-priori estimates \eqref{est: bootstrap argument a-priori est}, we obtain that
\begin{equation*}
\begin{split}
\sum_{|\mu|=0}^1 \left\| x_j \Big(\frac{D_x}{\langle D_x\rangle}\Big)^\mu v_\pm(t,\cdot)\right\|_{H^1}&\le C\left[E^2_3(t;W)^\frac{1}{2}+ tE_3(t;W)^\frac{1}{2}\right] + \|\mathrm{R}_1U(t,\cdot)\|_{H^{2,\infty}}E_0(t;W)^\frac{1}{2}\\
&\le CB\varepsilon t^{1+\frac{\delta}{2}}
\end{split}
\end{equation*}
\begin{equation*}
\Big(\text{resp. }\sum_{|\mu|=0}^1 \left\| x_j \Big(\frac{D_x}{\langle D_x\rangle}\Big)^\mu v_\pm(t,\cdot)\right\|_{H^{1,\infty}}
\le C E^2_3(t;W)^\frac{1}{2}+t\|V(t,\cdot)\|_{H^{2,\infty}}
\le C(A+B)\varepsilon t^{\frac{\delta_2}{2}}\Big),
\end{equation*}\normalsize
and the conclusion of the proof of \eqref{norms_H1_Linfty_xv-}.

Analogously, from \eqref{equality_xDkw} with $w=u_{-}$ and $f=\Nlw$,
\begin{align*}
\sum_{|\mu|=0}^1\left\|x_j D_k \Big(\frac{D_x}{|D_x|}\Big)^\mu u_{-}(t,\cdot)\right\|_{L^2}&\lesssim \|Z_ju_\pm(t,\cdot)\|_{L^2}+ t\|u_\pm(t,\cdot)\|_{L^2}+\|x_j\Nlw(t,\cdot)\|_{L^2(dx)}\\
&\le CB\varepsilon t^{1+\frac{\delta}{2}},
\end{align*}
as follows \eqref{est: bootstrap Enn}, \eqref{est: bootstrap E02}, \eqref{norm_Linfty_xv-} and the fact that
\begin{equation}\label{xj_Nlw_Linfty}
\|x_j \Nlw(t,\cdot)\|_{L^2}\lesssim \sum_{|\mu|=0}^1\left\|x_j \Big(\frac{D_x}{\langle D_x\rangle}\Big)^\mu v_\pm(t,\cdot)\right\|_{L^\infty}\|v_\pm(t,\cdot)\|_{H^1}.
\end{equation}
\endproof
\end{lem}

\begin{cor}
There exists a constant $C>0$ such that, for every $j=1,2$, $t\in [1,T]$, 
\begin{subequations} \label{est:L2,Linfty xNLkg}
\begin{align}
\|x_j \textit{NL}_{kg}(t,\cdot)\|_{L^2}&\le C(A+B)B\varepsilon^2 t^\frac{\delta+\delta_2}{2}, \label{xj_NLkg}\\
\|x_j \textit{NL}_{kg}(t,\cdot)\|_{L^\infty}&\le C(A+B)B\varepsilon^2 t^{-\frac{1}{2}+\frac{\delta_2}{2}},
\end{align}
\end{subequations}
and
\begin{subequations}
\begin{align}
\|x_j \textit{NL}_w(t,\cdot)\|_{L^2}&\le C(A+B)B\varepsilon^2 t^\frac{\delta+\delta_2}{2}, \label{est_L2_xNLw}\\
\|x_j \textit{NL}_w(t,\cdot)\|_{L^\infty}&\le C(A+B)B\varepsilon^2 t^{-1+\frac{\delta_2}{2}}.
\end{align}
\end{subequations}
\proof
From
\begin{equation*}
\|x_j \textit{NL}_{kg}(t,\cdot)\|_{L^2}\lesssim \sum_{\mu=0}^1 \left\| x_j (D_x\langle D_x\rangle^{-1})^\mu v_\pm(t,\cdot)\right\|_{L^\infty}\|u_\pm(t,\cdot)\|_{H^1},
\end{equation*}
and \eqref{norm_Linfty_xNLkg}, together with \eqref{xj_Nlw_Linfty} and
\begin{gather*}
\|x_j\textit{NL}_w(t,\cdot)\|_{L^\infty}\lesssim  \sum_{\mu=0}^1 \left\| x_j (D_x\langle D_x\rangle^{-1})^\mu v_\pm(t,\cdot)\right\|_{L^\infty}\|v_\pm(t,\cdot)\|_{H^{2,\infty}},
\end{gather*}
we immediately derive the estimates of the statement using \eqref{norm_Linfty_xv-} and a-priori estimates. 
\endproof
\end{cor}

\begin{lem}
There exists a positive constant $C>0$ such that, for any multi-index $I$ of length $k$, with $1\le k\le 2$, any $j=1,2$, $t\in [1,T]$,
\begin{equation}\label{norm_L2_xj-GammaIv-}
\sum_{|\mu|=0}^1\left\|x_j\Big(\frac{D_x}{\langle D_x\rangle}\Big)^\mu (\Gamma^I v)_\pm(t,\cdot)\right\|_{H^1}+ \left\| x_j D_x\Big(\frac{D_x}{|D_x|}\Big)^\mu (\Gamma^I u)_\pm(t,\cdot)\right\|_{L^2} \le CB\varepsilon t^{1+\frac{\delta_{3-k}}{2}}.
\end{equation}
\proof
We warn the reader that, throughout the proof, $C$ will denote a positive constant that may change line after line. As $\Gamma^I w_+=-\overline{\Gamma^I w_{-}}$, for any $I$ and $w\in\{v,u\}$, it is enough to prove the statement for $\Gamma^I v_{-}$, $\Gamma^I u_{-}$.

From equalities \eqref{equalities_app_xw} together with the fact that, for any multi-index $I$, $(\Gamma^I v)_{-}$, $(\Gamma^I u)_{-}$ are solution to
\begin{subequations}
\begin{equation}\label{half KG Gammav}
[D_t  + \langle D_x\rangle] (\Gamma^I v)_{-}(t,x) = \Gamma^I\textit{NL}_{kg}
\end{equation}
and 
\begin{equation}
[D_t  + \langle D_x\rangle] (\Gamma^I u)_{-}(t,x) = \Gamma^I\textit{NL}_{w}
\end{equation}
\end{subequations}
respectively, we derive that, for any $j,k\in \{1,2\}$,
\begin{subequations}
\begin{equation} \label{preliminary_norm_xGammav-}
\sum_{|\mu|=0}^1\left\|x_j\Big(\frac{D_x}{\langle D_x\rangle}\Big)^\mu (\Gamma^I v)_\pm(t,\cdot)\right\|_{H^1}\le \|Z_j(\Gamma^I v)_{-}(t,\cdot)\|_{L^2}+t\|(\Gamma^I v)_{-}(t,\cdot)\|_{L^2}+\|x_j\Gamma^I\textit{NL}_{kg}(t,\cdot)\|_{L^2}
\end{equation}
together with
\begin{equation} \label{preliminary_norm_xGammau}
\sum_{|\mu|=0}^1 \left\| x_j D_x\Big(\frac{D_x}{|D_x|}\Big)^\mu (\Gamma^I u)_\pm(t,\cdot)\right\|_{L^2}\le \|Z_j(\Gamma^I u)_{-}(t,\cdot)\|_{L^2}+t\|(\Gamma^I u)_{-}(t,\cdot)\|_{L^2}+\|x_j\Gamma^I\textit{NL}_w(t,\cdot)\|_{L^2}.
\end{equation}
\end{subequations}
The first two quantities in above right hand sides are bounded by $CB\varepsilon t^{1+\delta_{3-k}/2}$ after \eqref{est: bootstrap E02}, so the quantities that need to be estimated in order to prove the statement are the $L^2$ norms of $x_j\Gamma^I\textit{NL}_{kg}$, $x_j\Gamma^I\textit{NL}_w$, for $1\le |I|\le 2$. 

We first prove \eqref{norm_L2_xj-GammaIv-} for $|I|=1$ and $\Gamma^I=\Gamma$, reminding that from \eqref{Gamma_nonlinearity},
\begin{subequations}
\begin{equation}\label{notation: NLGamma}
\Gamma\textit{NL}_{kg} = Q^\mathrm{kg}_0\big((\Gamma v)_\pm, D_1u_\pm\big) + Q^\mathrm{kg}_0\big(v_\pm, D_1(\Gamma u)_\pm\big) + G^\mathrm{kg}_1\big(v_\pm, Du_\pm\big)
\end{equation}
and
\begin{equation}\label{Gamma_Nlw}
\Gamma\textit{NL}_w = Q^\mathrm{w}_0\big((\Gamma v)_\pm, D_1v_\pm\big) + Q^\mathrm{w}_0\big(v_\pm, D_1(\Gamma v)_\pm\big) + G^\mathrm{w}_1\big(v_\pm, Dv_\pm\big),
\end{equation}
\end{subequations}
with $G^\mathrm{kg}_1\big(v_\pm, Du_\pm\big) = G_1(v, \partial u)$, $G^\mathrm{w}_1\big(v_\pm, Dv_\pm\big) = G_1(v, \partial v)$, and $G_1$ given by \eqref{def_G1}.

By multiplying $x_j$ against the Klein-Gordon component in each product of $\Gamma\textit{NL}_{kg}$ we find that
\begin{multline} \label{norm_xmZn_NLkg}
\|x_j\Gamma\textit{NL}_{kg}(t,\cdot)\|_{L^2}\lesssim \sum_{|\mu|=0}^1 \left\|x_j\Big(\frac{D_x}{\langle D_x\rangle}\Big)^\mu(\Gamma v)_{-}(t,\cdot)\right\|_{L^2}\left(\|U(t,\cdot)\|_{H^{2,\infty}}+\|\mathrm{R}_1U(t,\cdot)\|_{H^{2,\infty}}\right) \\
 +\sum_{|\mu|=0}^1\left\| x_j\Big(\frac{D_x}{\langle D_x\rangle}\Big)^\mu v_\pm(t,\cdot)\right\|_{L^\infty}\left(\|(\Gamma u)_\pm(t,\cdot)\|_{H^1}+ \|u_\pm(t,\cdot)\|_{H^1}+ \|D_t u_\pm(t,\cdot)\|_{L^2}\right),
\end{multline}
which injected into \eqref{preliminary_norm_xGammav-} with $\Gamma^I=\Gamma$, together with \eqref{Hs_norm_DtU} with $s=0$, \eqref{norm_Linfty_xv-}, and a-priori estimates \eqref{est: bootstrap argument a-priori est}, gives that 
\begin{equation*}
\sum_{|\mu|=0}^1\left\|x_j\Big(\frac{D_x}{\langle D_x\rangle}\Big)^\mu (\Gamma^I v)_\pm(t,\cdot)\right\|_{H^1}\le CB\varepsilon t^{1+\frac{\delta_2}{2}}.
\end{equation*}
Similarly, using the above estimate together with \eqref{Hs norm DtV} with $s=0$, \eqref{norm_Linfty_xv-} and a-priori estimates, we derive that
\begin{equation}\label{est_xGammaNLw_proof}
\begin{split}
&\|x_j\Gamma\textit{NL}_w(t,\cdot)\|_{L^2}\lesssim  \sum_{|\mu|=0}^1 \left\|x_j\Big(\frac{D_x}{\langle D_x\rangle}\Big)^\mu(\Gamma v)_{-}(t,\cdot)\right\|_{L^2}\|v_\pm(t,\cdot)\|_{H^{2,\infty}}\\
&+  \sum_{ |\mu|=0}^1 \left\| x_j\Big(\frac{D_x}{\langle D_x\rangle}\Big)^\mu v_\pm(t,\cdot)\right\|_{L^\infty}  \left(\| (\Gamma v)_\pm(t,\cdot)\|_{H^1}+\|v_\pm(t,\cdot)\|_{H^1}+\|D_t v_\pm(t,\cdot)\|_{L^2}\right) \\
&\le C(A+B)B\varepsilon^2 t^{\delta_2}.
\end{split}
\end{equation}
Plugging the above inequality in \eqref{preliminary_norm_xGammau} for $\Gamma^I=\Gamma$ and using again a-priori estimates we deduce that
\begin{equation*}
\sum_{|\mu|=0}^1\left\|x_j D_k\Big(\frac{D_x}{|D_x|}\Big)^\mu(\Gamma u)_{-}(t,\cdot)\right\|_{L^2}\le CB\varepsilon t^{1+\frac{\delta_2}{2}},
\end{equation*}
and conclude the proof of \eqref{norm_L2_xj-GammaIv-} when $|I|=1$.

When $|I|=2$ we observe that, from \eqref{Gamma_I_nonlinearity},
\begin{multline} \label{GammaI_NLkg}
\Gamma^I \textit{NL}_{kg} = Q^\mathrm{kg}_0(v^I_\pm, D_1u_\pm)+Q^\mathrm{kg}_0(v_\pm, D_1u^I_\pm)+ \sum_{\substack{(I_1,I_2)\in\mathcal{I}(I)\\ |I_1|=|I_2|=1}}Q^\mathrm{kg}_0(v^{I_1}_\pm, D_1u^{I_2}_\pm)\\ + \sum_{\substack{(I_1,I_2)\in\mathcal{I}(I)\\ |I_1|+|I_2|\le 1}}c_{I_1,I_2}Q^\mathrm{kg}_0(v^{I_1}_\pm, Du^{I_2}_\pm),
\end{multline}
with $c_{I_1,I_2}\in \{-1,0,1\}$.
Since the $L^2$ norm of terms indexed in $I_1,I_2$ with $|I_1|=|I_2|=1$ can be estimated using the Sobolev injection as follows:
\begin{equation}\label{xj_Q(vI1,uI2)_appendix}
\left\|x_j Q^\mathrm{kg}_0(v^{I_1}_\pm, D_1u^{I_2}_\pm)\right\|_{L^2} \lesssim \sum_{ |\mu|=0}^1\|v^{I_1}_\pm(t,\cdot)\|_{H^2}\left\|x_jD_1\Big(\frac{D_x}{|D_x|}\Big)^\mu u^{I_2}_\pm(t,\cdot)\right\|_{L^2},
\end{equation}
from \eqref{GammaI_NLkg} we derive that
\begin{align*}
&\|x_j\Gamma^I \textit{NL}_{kg}\|_{L^2}\lesssim \sum_{\substack{|J|\le 2\\ |\mu|,\nu=0}}^1 \left\| x_j\Big(\frac{D_x}{\langle D_x\rangle}\Big)^\mu(\Gamma^J v)_{-}(t,\cdot)\right\|_{L^2}\Big(\|\mathrm{R}_1^\nu u_\pm(t,\cdot)\|_{H^{2,\infty}}+ \|D_t\mathrm{R}_1^\nu u_\pm(t,\cdot)\|_{H^{1,\infty}}\Big) \\
& + \sum_{|\mu|=0}^1\left\|x_j\Big(\frac{D_x}{\langle D_x\rangle}\Big)^\mu v_\pm(t,\cdot)\right\|_{L^\infty}\Big[\|u^I_\pm(t,\cdot)\|_{H^1} + \sum_{|J|\le 1}\big(\|u^J_\pm(t,\cdot)\|_{H^1}+\|D_tu^J_\pm(t,\cdot)\|_{L^2}\big) \Big]\\
& + \sum_{\substack{|I_1|=|I_2|=1\\ |\mu|=0,1}}\|v^{I_1}_\pm(t,\cdot)\|_{H^2}\left\|x_jD_1\Big(\frac{D_x}{|D_x|}\Big)^\mu u^{I_2}_\pm(t,\cdot)\right\|_{L^2}.
\end{align*}
As before, injecting the above inequality into \eqref{preliminary_norm_xGammav-}, using a-priori estimates \eqref{est: bootstrap argument a-priori est} and the fact that $\varepsilon_0<(2A)^{-1}$, together with \eqref{Hs_norm_DtU} with $s=0$, \eqref{Hsinfty_norm_DtU}, \eqref{Hsinfty norm Dt R1U-new} with $s=1$, \eqref{DtUI}, \eqref{norm_Linfty_xv-}, and \eqref{norm_L2_xj-GammaIv-} with $k=1$, we obtain that
\begin{equation}\label{xGammaIv_I=2}
\sum_{|\mu|=0}^1\left\|x_j\Big(\frac{D_x}{\langle D_x\rangle}\Big)^\mu(\Gamma^I v)_{-}(t,\cdot)\right\|_{H^1}\le CB\varepsilon t^{1+\frac{\delta_1}{2}}.
\end{equation}
Analogously, since
\begin{multline*}
\Gamma^I \textit{NL}_{w} = Q^\mathrm{w}_0(v^I_\pm, D_1v_\pm)+Q^\mathrm{w}_0(v_\pm, D_1v^I_\pm)+ \sum_{\substack{(I_1,I_2)\in\mathcal{I}(I)\\ |I_1|=|I_2|=1}}Q^\mathrm{w}_0(v^{I_1}_\pm, D_1v^{I_2}_\pm)\\ + \sum_{\substack{(I_1,I_2)\in\mathcal{I}(I)\\ |I_1|+|I_2|<2}}c_{I_1,I_2}Q^\mathrm{w}_0(v^{I_1}_\pm, Dv^{I_2}_\pm),
\end{multline*}
we have that
\begin{equation*}
\begin{split}
&\|x_j\Gamma^I\textit{NL}_w\|_{L^2} \lesssim \sum_{\substack{|J|\le 2\\|\mu|=0,1}}\left\|x_j\Big(\frac{D_x}{\langle D_x\rangle}\Big)^\mu (\Gamma^J v)_\pm(t,\cdot)\right\|_{L^2}\Big(\|v_\pm(t,\cdot)\|_{H^{2,\infty}} + \|D_tv_\pm(t,\cdot)\|_{H^{1,\infty}}\Big)\\
&+ \sum_{ |\mu|=0}^1\left\| x_j \Big(\frac{D_x}{\langle D_x\rangle}\Big)^\mu v_\pm(t,\cdot)\right\|_{H^{1,\infty}}\left(\sum_{|J|\le 2}\|(\Gamma^J v)_\pm(t,\cdot)\|_{H^1}+\sum_{|J|\le 1}\|D_tv^J_\pm(t,\cdot)\|_{L^2}\right)\\
& + \sum_{\substack{|I_1|=|I_2|=1\\|\mu|=0,1 }}\|(\Gamma^{I_1}v)_\pm(t,\cdot)\|_{H^2}\left\|x_j\Big(\frac{D_x}{\langle D_x\rangle}\Big)^\mu (\Gamma^{I_2} v)_\pm(t,\cdot)\right\|_{L^2},
\end{split}
\end{equation*}
so from \eqref{Hs norm DtV} with $s=0$, \eqref{est: Hsinfty Dt V} with $s=1$, \eqref{DtVI}, \eqref{norm_Linfty_xv-}, \eqref{norm_L2_xj-GammaIv-} with $|I|=1$, \eqref{xGammaIv_I=2} and a-priori estimates \eqref{est: bootstrap argument a-priori est}, we deduce
\begin{equation*}
\sum_{|\mu|=0}^1\left\|x_j D_k\Big(\frac{D_x}{|D_x|}\Big)^\mu (\Gamma^I u)_{-}(t,\cdot)\right\|_{L^2}\lesssim CB\varepsilon t^{1+\frac{\delta_1}{2}},
\end{equation*}
and hence conclude the proof of inequality \eqref{norm_L2_xj-GammaIv-} also for the case $|I|=2$.
\endproof
\end{lem}

\begin{cor}
There exists a positive constant $C>0$ such that, for any $\Gamma\in\mathcal{Z}$, $j=1,2$, and every $t\in [1,T]$,
\begin{subequations}
\begin{align}
\|x_j\Gamma \textit{NL}_{kg}(t,\cdot)\|_{L^2} &\le C(A+B)B\varepsilon^2 t^{\frac{1}{2}+\frac{\delta_2}{2}}, \label{xjGamma_Nlkg}\\
\|x_j\Gamma \textit{NL}_w(t,\cdot)\|_{L^2} &\le C(A+B)B\varepsilon^2 t^{\delta_2}. \label{xjGamma_NLw}
\end{align}
\end{subequations}
\proof
Estimate \eqref{xjGamma_Nlkg} follows straightly from \eqref{norm_xmZn_NLkg}, \eqref{Hs_norm_DtU} with $s=0$, and estimates \eqref{est: bootstrap argument a-priori est}, \eqref{norm_Linfty_xv-}, and \eqref{norm_L2_xj-GammaIv-} with $k=1$, while
\eqref{xjGamma_NLw} has already been proved in \eqref{est_xGammaNLw_proof}.
\endproof
\end{cor}

\begin{lem}
There exists a constant $C>0$ such that, for every $i,j=1,2$, every $t\in [1,T]$,
\begin{subequations}\label{est_x2v-}
\begin{gather}
\sum_{|\mu|=0}^1\left\| x_jx_k \Big(\frac{D_x}{\langle D_x\rangle}\Big)^\mu v_\pm(t,\cdot)\right\|_{L^2}\le CB\varepsilon t^{2+\frac{\delta_2}{2}},  \label{norm_L2_xixjv-}\\
\sum_{|\mu|=0}^1\left\| x_jx_k \Big(\frac{D_x}{\langle D_x\rangle}\Big)^\mu v_\pm(t,\cdot)\right\|_{L^\infty}\le C(A+B)\varepsilon t^{1+\frac{\delta_2}{2}}. \label{norm_Linfty_xixjv-}
\end{gather}
\end{subequations}
Moreover, for any $\Gamma\in \mathcal{Z}$, 
\begin{equation} \label{est:xixjGamma v-}
\sum_{|\mu|=0}^1\left\|x_i x_j\Big(\frac{D_x}{\langle D_x\rangle}\Big)^\mu (\Gamma v)_\pm(t,\cdot)\right\|_{L^2}\le CB\varepsilon t^{2+\frac{\delta_2}{2}}.
\end{equation}
\proof
The proof of the statement follows from the fact that, by multiplying \eqref{xjw_Zjw} by $x_i$ and using that
\begin{equation*}
\|x_i x_j\textit{NL}_{kg}(t,\cdot)\|_{L^2}\lesssim \sum_{|\mu|=0}^1 \left\|x_i x_j \Big(\frac{D_x}{\langle D_x\rangle}\Big)^\mu v_{-}(t,\cdot)\right\|_{L^2} \left(\|u_\pm(t,\cdot)\|_{H^{2,\infty}}+\|\mathrm{R}_1u_\pm(t,\cdot)\|_{H^{2,\infty}}\right)
\end{equation*}
together with
\begin{equation*}
\|x_i x_j \textit{NL}_{kg}(t,\cdot)\|_{L^\infty} \lesssim\sum_{|\mu|=0}^1\left\| x_jx_k \Big(\frac{D_x}{\langle D_x\rangle}\Big)^\mu v_{-}(t,\cdot)\right\|_{L^\infty}\left(\|u_\pm (t,\cdot)\|_{H^{2,\infty}}+\|\mathrm{R}_1u_\pm(t,\cdot)\|_{H^{2,\infty}}\right),
\end{equation*}
we derive that
\begin{multline*}
\sum_{|\mu|=0}^1\left\| x_jx_k \Big(\frac{D_x}{\langle D_x\rangle}\Big)^\mu v_\pm(t,\cdot)\right\|_{L^2} \lesssim \sum_{|\mu|=0}^1\left(\|x^\mu_i (Z_jv)_{-}(t,\cdot)\|_{L^2}+ t\|x_i^\mu v_{-}(t,\cdot)\|_{L^2}\right)\\
+ \sum_{\mu=0}^1 \left\|x_i x_j \Big(\frac{D_x}{\langle D_x\rangle}\Big)^\mu v_{-}(t,\cdot)\right\|_{L^2} \left(\|u_\pm(t,\cdot)\|_{H^{2,\infty}}+\|\mathrm{R}_1u_\pm(t,\cdot)\|_{H^{2,\infty}}\right),
\end{multline*}
and using that operator $\langle D_x\rangle^{-1}$ is bounded from $H^1$ to $L^\infty$
\begin{multline*}
\sum_{|\mu|=0}^1\left\|x_i x_j\Big(\frac{D_x}{\langle D_x\rangle}\Big)^\mu v_{-}(t,\cdot)\right\|_{L^\infty} \lesssim \sum_{k=0}^1\left(\|x_i^k (Z_j v)_{-}(t,\cdot)\|_{H^1} + t \|x^k_i v_{-}(t,\cdot)\|_{H^{1,\infty}}\right)\\
 +\sum_{k, |\mu|=0}^1 \left\| x^k_i  x_j\Big(\frac{D_x}{\langle D_x\rangle}\Big)^\mu v_{-}(t,\cdot)\right\|_{L^\infty}\left(\|u_\pm(t,\cdot)\|_{H^{2,\infty}}+\|\mathrm{R}_1u_\pm(t,\cdot)\|_{H^{2,\infty}}\right).
\end{multline*}
As $\varepsilon_0>0$ verifies that $\varepsilon_0<(2A)^{-1}$, inequality \eqref{norm_xv-}, \eqref{norm_L2_xj-GammaIv-} with $k=1$, and a-priori estimates \eqref{est: bootstrap argument a-priori est} imply that
\begin{equation*}
\sum_{ |\mu|=0}^1 \left\|x_i x_j \Big(\frac{D_x}{\langle D_x\rangle}\Big)^\mu v_{-}(t,\cdot)\right\|_{L^2}\lesssim CB\varepsilon t^{2+\frac{\delta_2}{2}},
\end{equation*}
while from \eqref{norm_Linfty_xv-}, \eqref{norm_L2_xj-GammaIv-} with $k=1$ and a-priori estimates, 
\begin{equation*}
\sum_{ |\mu|=0}^1 \left\|x_i x_j \Big(\frac{D_x}{\langle D_x\rangle}\Big)^\mu v_{-}(t,\cdot)\right\|_{L^\infty}\le C(A+B)\varepsilon t^{1+\frac{\delta_2}{2}}.
\end{equation*}
As $v_+=-\overline{v_{-}}$, that implies the first part of the statement.

Analogously, using \eqref{xjw_Zjw} with $w= (\Gamma v)_{-}$ and multiplying that relation by $x_i$ we find that
\begin{multline} \label{xixjGammav_preliminary_1}
\sum_{|\mu|=0}^1\left\| x_i x_j\Big(\frac{D_x}{\langle D_x\rangle}\Big)^\mu (\Gamma v)_{-}(t,\cdot)\right\|_{L^2}\\
\lesssim \sum_{\mu=0}^1 \big( \|x_i^\mu Z_j (\Gamma v)_{-}(t,\cdot)\|_{L^2} + t\|x^\mu_i (\Gamma v)_{-}(t,\cdot)\|_{L^2} + \|x_i^\mu x_j\Gamma\textit{NL}_{kg}(t,\cdot)\|_{L^2}\big),
\end{multline}
and after \eqref{norm_L2_xj-GammaIv-}, \eqref{xjGamma_Nlkg} and a-priori estimates,
\begin{equation} \label{xixjGammav_preliminary_2}
\sum_{\mu=0}^1\big( \|x_i^\mu Z_j (\Gamma v)_{-}(t,\cdot)\|_{L^2} + t\|x^\mu_i (\Gamma v)_{-}(t,\cdot)\|_{L^2}\big) + \| x_j\Gamma\textit{NL}_{kg}(t,\cdot)\|_{L^2} \le CB\varepsilon t^{2+\frac{\delta_2}{2}}.
\end{equation}
By multiplying both $x_i, x_j$ against each Klein-Gordon factor in $\Gamma\textit{NL}_{kg}$ (see equality \eqref{notation: NLGamma}) we derive that
\begin{multline*}
\|x_ix_j \Gamma\textit{NL}_{kg}(t,\cdot)\|_{L^2}\lesssim \sum_{|\mu|\nu=0}^1 \left\|x_i x_j \Big(\frac{D_x}{\langle D_x\rangle}\Big)^\mu (\Gamma v)_{-}(t,\cdot)\right\|_{L^2}\|\mathrm{R}_1^\nu u_\pm(t,\cdot)\|_{H^{2,\infty}}\\
+ \sum_{|\mu|=0}^1 \left\|x_i x_j\Big(\frac{D_x}{\langle D_x\rangle}\Big)^\mu v_\pm(t,\cdot)\right\|_{L^\infty}\left(\|(\Gamma u)_\pm(t,\cdot)\|_{H^1}+\|u_\pm(t,\cdot)\|_{H^1}+\|D_tu_\pm(t,\cdot)\|_{L^2}\right),
\end{multline*}
so by \eqref{Hs_norm_DtU} with $s=0$, \eqref{norm_Linfty_xixjv-}, a-priori estimates and the fact that $\varepsilon_0<(2A)^{-1}$,
\begin{equation*}
\|x_ix_j \Gamma\textit{NL}_{kg}(t,\cdot)\|_{L^2} \le \frac{1}{2} \|x_i  x_j(\Gamma v)_{-}(t,\cdot)\|_{L^2}+ C(A+B)B\varepsilon^2 t^{1+\delta_2},
\end{equation*}
which injected in \eqref{xixjGammav_preliminary_1}, together with \eqref{xixjGammav_preliminary_2}, implies \eqref{est:xixjGamma v-}.
\endproof
\end{lem}

\begin{cor}
There exists a constant $C>0$ such that, for every $i,j=1,2$, every $t\in [1,T]$, \ref{Lem_appendix: Technical_estimates}, 
\begin{equation} \label{est_L2:xixj_NL}
\|x_i x_j\textit{NL}_{kg}(t,\cdot)\|_{L^2} + \|x_i x_j\textit{NL}_w(t,\cdot)\|_{L^2}\le C(A+B)B\varepsilon^2 t^{1+\frac{\delta+\delta_2}{2}}.
\end{equation}
\proof
Straightforward after \eqref{est: bootstrap enhanced Enn}, \eqref{norm_Linfty_xixjv-} and the following inequality
\begin{multline*}
\|x_i x_j\textit{NL}_{kg}(t,\cdot)\|_{L^2} + \|x_i x_j\textit{NL}_w(t,\cdot)\|_{L^2} \\
\lesssim \sum_{|\mu|=0}^1\left\|x_ix_j\Big(\frac{D_x}{\langle D_x\rangle}\Big)^\mu v_\pm(t,\cdot)\right\|_{L^\infty}\left(\|u_\pm(t,\cdot)\|_{H^1}+\|v_\pm(t,\cdot)\|_{H^1}\right).
\end{multline*}
\endproof
\end{cor}

\begin{lem}
There exists a constant $C>0$ such that, for any $i,j,k=1,2$, every $t\in [1,T]$, \ref{Lem_appendix: Technical_estimates}, 
\begin{equation}\label{est:x3_v-}
\sum_{|\mu|=0}^1\left\|x_i x_j x_k \Big(\frac{D_x}{\langle D_x\rangle}\Big)^\mu v_\pm(t,\cdot)\right\|_{L^2}\le CB\varepsilon t^{3+\frac{\delta_2}{2}}.
\end{equation}
\proof
Using equality \eqref{xjw_Zjw} we derive that
\begin{multline*}
\|x_i x_j x_k v_{-}(t,\cdot)\|_{L^2} \lesssim \sum_{\mu_1, \mu_2=0}^1\left(\|x^{\mu_1}_i x_j^{\mu_2} (Z_kv)_{-}(t,\cdot)\|_{L^2}+ t\|x^{\mu_1}_i x_j^{\mu_2} v_{-}(t,\cdot)\|_{L^2}\right)\\
+ \sum_{\mu_1, \mu_2, |\mu|=0}^1 \left\|x_i^{\mu_1} x_j^{\mu_2}x_k\Big(\frac{D_x}{\langle D_x\rangle}\Big)^\mu v_{-}(t,\cdot)\right\|_{L^2} \left(\|u_\pm(t,\cdot)\|_{H^{2,\infty}}+\|\mathrm{R}_1u_\pm(t,\cdot)\|_{H^{2,\infty}}\right),
\end{multline*}
so the result of the statement is a straight consequence of \eqref{norm_xv-}, \eqref{norm_L2_xj-GammaIv-}, \eqref{norm_L2_xixjv-}, \eqref{est:xixjGamma v-}, a-priori estimates, and the fact that $\varepsilon_0$ is smaller than $(2A)^{-1}$.
\endproof
\end{lem}

\section{First range of estimates} \label{sec_appB: first range of estimates}

The aim of this section is to show that, if a-priori estimates \eqref{est: bootstrap argument a-priori est} are satisfied for every $t\in [1,T]$, for some fixed $T>1$, then in the same interval the semi-classical Sobolev norms of the semi-classical functions $\ut, \vt$ introduced in \eqref{def utilde vtilde} grow in time at a moderate rate $t^\beta$, for some small $\beta>0$.
More precisely, in lemma \ref{Lem: from energy to norms in sc coordinates-WAVE} we prove that this is the case for the $H^s_h(\mathbb{R}^2)$ norm of $\ut$, $\ut^{\Sigma_j,k}$ (see definition \eqref{def utilde-Sigma,k}) for any $s\le n-15$, and for the $L^2(\mathbb{R}^2)$ norm of those functions when operators $\Omega_h$ and $\Mcal$, introduced in \eqref{def_Omega_h} and \eqref{def Mj} respectively, are acting on them and frequencies are less or equal than $h^{-\sigma}$, for some small $\sigma>0$.
Lemma \ref{Lem: from energy to norms in sc coordinates-KG} shows that this moderate growth is also enjoyed by the $H^s_h(\mathbb{R}^2)$ norm of $\vt$, again for $s\le n=15$, and by the $L^2(\mathbb{R}^2)$ norm of $\Lcal\vt$ (see \eqref{def Lj}) when restricted to frequencies $|\xi|\lesssim h^{-\sigma}$.
The proof of this latter lemma will require some intermediate results, among which lemma \ref{Lem_appendix: preliminary est VJ} that provides us with a first non-sharp estimate of the $L^\infty(\mathbb{R}^2)$ norm of Klein-Gordon functions $v_\pm$ when one Klainerman vector field is acting on them (and again frequencies are localized for $|\xi|\lesssim t^\sigma$).
This estimate will successively improved to the sharpest one \eqref{sharp_est_VI} in lemma \ref{Lem_appendix: sharp_est_VJ} of section \ref{Sec_App_4}.

As said at the beginning of this chapter, we prove the below results under the hypothesis that a-priori estimates \eqref{est: bootstrap argument a-priori est} are satisfied in some fixed $[1,T]$, with $\varepsilon_0<(2A+B)^{-1}$.
We remind here that, if $\chi\in C^\infty_0(\mathbb{R}^2)$ and $\sigma>0$, $\chi(t^{-\sigma}D_x)$ is a bounded operator from $H^s$ to $L^2$ with norm $O(t^{\sigma s})$, and on $L^\infty$ uniformly in time.

\begin{lem} \label{Lem: from energy to norms in sc coordinates-WAVE}
Let $\widetilde{u}, \widetilde{u}^{\Sigma_j,k}$ be defined, respectively, in \eqref{def utilde vtilde} and \eqref{def utilde-Sigma,k}, and $s\le n-15$.
There exists a constant $C>0$ such that, for any $\theta_0, \chi \in C^\infty_0(\mathbb{R}^2)$ and every $t\in [1,T]$,
\begin{subequations}\label{inequality: from_energy_to_norm_insc}
\begin{gather} 
\|\widetilde{u}(t,\cdot)\|_{H^s_h}+ \|\widetilde{u}^{\Sigma_j,k}(t,\cdot)\|_{H^s_h} \le CB\varepsilon t^{\frac{\delta}{2}+\kappa}, \label{est:utilde-Hs}\\
\|\Omega_h\widetilde{u}^{\Sigma_j,k}(t,\cdot)\|_{L^2}\le CB\varepsilon t^{\frac{\delta_2}{2}+\kappa},\label{est:Omega-utilde}\\
\sum_{|\mu|=1}\left(\|\oph(\chi(h^\sigma\xi))\mathcal{M}^\mu \widetilde{u}(t,\cdot)\|_{L^2}+ \|\mathcal{M}^\mu\widetilde{u}^{\Sigma_j,k}(t,\cdot)\|_{L^2}\right)  \le  C(A+B)\varepsilon t^{\frac{\delta_2}{2}+\kappa},\label{est:Mutilde}\\
\sum_{|\mu|=1}\|\theta_0(x)\Omega_h\mathcal{M}^\mu\widetilde{u}^{\Sigma_j,k}(t,\cdot)\|_{L^2} \le CB\varepsilon  t^{\frac{\delta_1}{2}+\kappa},\label{est:Omega-M-utilde}
\end{gather}
\end{subequations}
with $\kappa=\sigma\rho$ if $\rho\ge 0$, 0 otherwise.
\proof
We warn the reader that, throughout the proof, $C$ and $\beta$ will denote positive constants that may change line after line, with $\beta\rightarrow 0$ as $\sigma\rightarrow 0$. We will also use the following concise notation
\begin{gather*}
\phi^j_k(\xi):= \Sigma(\xi)(1-\chi_0)(h^{-1}\xi)\varphi(2^{-k}\xi)\chi_0(h^\sigma\xi), 
\end{gather*}
reminding that
\begin{equation} \label{Op(Sigma varphi)}
\left\|  \oph(\phi^j_k(\xi))\right\|_{\mathcal{L}(L^2)} =O(h^{-\kappa}),
\end{equation}
with $\kappa=\sigma\rho$ if $\rho\ge 0$, 0 otherwise.

Inequality \eqref{est:utilde-Hs} is straightforward after \eqref{Op(Sigma varphi)}, definitions \eqref{def utilde vtilde} and \eqref{def uNF}, inequality \eqref{Hs norm uNF- u-}, and a-priori estimate \eqref{est: boostrap vpm}. 
By commutating $\oph(\phi^j_k(\xi))$ with $\mathcal{M}$ (the commutator with $\Omega_h$ being zero if $\varphi, \chi_0$ are supposed to be radial) and using \eqref{Op(Sigma varphi)} we observe that there is some $\chi\in C^\infty_0(\mathbb{R}^2)$ such that
\begin{equation*}
\begin{gathered}
\|\Omega_h \widetilde{u}^{\Sigma_j,k}(t,\cdot)\|_{L^2}\lesssim h^{-\kappa}\|\oph(\chi_0(h^\sigma\xi))\Omega_h \widetilde{u}(t,\cdot)\|_{L^2},\\
\|\mathcal{M} \widetilde{u}^{\Sigma_j,k}(t,\cdot)\|_{L^2}\lesssim h^{-\kappa}\sum_{|\nu|=0}^1\|\oph(\chi(h^\sigma\xi))\mathcal{M}^\nu \widetilde{u}(t,\cdot)\|_{L^2}, 
\end{gathered}
\end{equation*}
\small
\begin{equation*}
\|\theta_0(x)\Omega_h \mathcal{M}\widetilde{u}^{\Sigma_j,k}(t,\cdot)\|_{L^2} 
\lesssim \|\theta_0(x)\oph(\phi^j_k(\xi))\Omega_h \mathcal{M}\widetilde{u}(t,\cdot)\|_{L^2} + h^{-\kappa}\sum_{\mu=0}^1\|\oph(\chi(h^\sigma\xi))\Omega^\mu_h\widetilde{u}(t,\cdot)\|_{L^2}.
\end{equation*}\normalsize
Therefore, as $h=t^{-1}$, in order to prove \eqref{est:Omega-utilde}-\eqref{est:Omega-M-utilde} it is enough to show that, for any $\chi\in C^\infty_0(\mathbb{R}^2)$,
\begin{subequations}
\begin{gather}
\|\oph(\chi(h^\sigma\xi))\Omega_h \widetilde{u}(t,\cdot)\|_{L^2}\le CB\varepsilon t^{\frac{\delta_2}{2}},\label{est_Omega-utilde_proof} \\
\|\oph(\chi(h^\sigma\xi))\mathcal{M} \widetilde{u}(t,\cdot)\|_{L^2} \le C(A+B)\varepsilon t^{\frac{\delta_2}{2}},\label{est_Mutilde-proof} \\
 \|\theta_0(x)\oph(\phi^j_k(\xi))\Omega_h \mathcal{M}\widetilde{u}(t,\cdot)\|_{L^2} \le CB\varepsilon t^{\frac{\delta_1}{2}}. \label{est_Omega-M-utilde-proof}
\end{gather}
\end{subequations}
Estimate \eqref{est_Omega-utilde_proof} follows from definitions \eqref{def utilde vtilde} and \eqref{def uNF}, inequality \eqref{est: L2 Omega integral D with cut-off} with $u=v=v_\pm$, and a-priori estimates \eqref{est: bootstrap argument a-priori est}, as
\begin{align*}
\|\oph(\chi(h^\sigma\xi))\Omega_h\widetilde{u}(t,\cdot)\|_{L^2} &\lesssim \|(\Omega u)_{-}(t,\cdot)\|_{L^2}+ \|\chi(t^{-\sigma}D_x)\Omega(u^{NF}-u_{-})(t,\cdot)\|_{L^2}\\
&\lesssim \|\Omega U(t,\cdot)\|_{L^2}+t^\beta \left(\|V(t,\cdot)\|_{L^2}+\|\Omega V(t,\cdot)\|_{L^2}\right)\|V(t,\cdot)\|_{H^{17,\infty}}\\
& \le C(1+A\varepsilon t^{-1+\beta})E^2_3(t;W)^\frac{1}{2}\le CB\varepsilon t^{\frac{\delta_2}{2}}.
\end{align*}
From equality \eqref{relation between Zju and Mj utilde-new} and definition \eqref{def uNF} of $u^{NF}$ we deduce that
\begin{equation} \label{L2_norm_Mutilde}
\begin{split}
&\|\oph(\chi(h^\sigma\xi))\mathcal{M}_n\widetilde{u}(t,\cdot) \|_{L^2}\lesssim \|Z_nU(t,\cdot)\|_{L^2}+ \|\chi(t^{-\sigma}D_x)Z_n(u^{NF}-u_{-})(t,\cdot)\|_{L^2}\\
& + \|\widetilde{u}(t,\cdot)\|_{L^2}+ \|\oph(\chi(h^\sigma\xi))[t(tx_j) [q_w+c_w](t,tx)]\|_{L^2(dx)}+ \|\chi(t^{-\sigma}D_x)(x_nr^{NF}_w)(t,\cdot)\|_{L^2},
\end{split}
\end{equation}
with $q_w$, $c_w$ and $r^{NF}_w$ given by \eqref{def_qw}, \eqref{def_cw} and \eqref{def rNF} respectively.
We first notice that, after inequality \eqref{est:L2 Z integral D with cut-off} with $u=v=v_\pm$, \eqref{Hs norm DtV} with $s=0$, a-priori estimates, and the fact that $A\varepsilon_0\le 1$, 
\begin{multline} \label{norm_Zn(uNf-u-)}
\|\chi(t^{-\sigma}D_x)Z_n(u^{NF}-u_{-})(t,\cdot)\|_{L^2}\\ \lesssim t^\beta \left(\|D_tV(t,\cdot)\|_{L^2}\|V(t,\cdot)\|_{H^{13}}+ \|V(t,\cdot)\|_{H^{15,\infty}}\|Z_nV(t,\cdot)\|_{L^2}\right)\le CB\varepsilon t^{\beta+\delta}.
\end{multline}
Let us also observe that from \eqref{def_qw}, \eqref{def_cw} we have that
\begin{equation} \label{qw+cw}
\begin{split}
q_w(t,x) + c_w(t,x) &= \frac{1}{2}\Im\left[\overline{v_{-}}\,  D_1 v_{-} - \overline{\frac{D_x}{\langle D_x\rangle}v_{-}}\cdot\frac{D_xD_1}{\langle D_x\rangle}v_{-}\right] (t,x)\\
& = \frac{h^2}{2} \Im \left[\overline{\widetilde{V}}\, \oph(\xi_1)\widetilde{V} - \overline{\oph\Big(\frac{\xi_1}{\langle \xi\rangle}\Big)\widetilde{V}}\cdot \oph\Big(\frac{\xi\xi_1}{\langle \xi\rangle}\Big)\widetilde{V}\right]\Big(t,\frac{x}{t}\Big),
\end{split}
\end{equation}
where $\widetilde{V}(t,x):=tv_{-}(t,tx)$ is such that, for every $s, \rho\ge 0$,
\begin{equation*}
\|\widetilde{V}(t,\cdot)\|_{H^s_h}=\| v_{-}(t,\cdot)\|_{H^s}, \quad \|\widetilde{V}(t,\cdot)\|_{H^{\rho,\infty}_h} = t\|v_{-}(t,\cdot)\|_{H^{\rho,\infty}}.
\end{equation*}
Moreover, by \eqref{relation_Zjw_Ljwidetilde(w)} with $w=v_{-}$ and $f=\textit{NL}_{kg}$
\begin{equation}\label{Lwidetilde(V)}
\begin{split}
\|\mathcal{L}_j\widetilde{V}(t,\cdot)\|_{H^1_h}&\lesssim \|Z_j v_{-}(t,\cdot)\|_{L^2}+\|v_{-}(t,\cdot)\|_{L^2}\\
&+ \left(\|x_jv_\pm(t,\cdot)\|_{L^\infty}+\left\|x_j\frac{D_x}{\langle D_x\rangle}v_\pm(t,\cdot)\right\|_{L^\infty}\right)\|U(t,\cdot)\|_{H^1}.
\end{split}
\end{equation}
Using \eqref{qw+cw} along with the definition of $\mathcal{L}_j$ in \eqref{def Lj} we derive that
\begin{equation} \label{x_jqw}
\begin{split}
& t(tx_j)[q_w+ c_w](t,tx) = \frac{1}{2}\Im\left[\overline{\widetilde{V}} \oph(\xi_1)(h\mathcal{L}_j\widetilde{V}) +\overline{\widetilde{V}}\oph\Big(\frac{\xi_1\xi_j}{\langle\xi\rangle}\Big)\widetilde{V} + \overline{\widetilde{V}}[x_j, \oph(\xi_1)]\widetilde{V} \right. \\
&- \overline{\oph\Big(\frac{\xi}{\langle\xi\rangle}\Big)\widetilde{V}}\cdot \oph\Big(\frac{\xi\xi_1}{\langle\xi\rangle}\Big)(h\mathcal{L}_j\widetilde{V}) - \overline{\oph\Big(\frac{\xi}{\langle\xi\rangle}\Big)\widetilde{V}}\cdot \oph\Big(\frac{\xi\xi_1\xi_j}{\langle\xi\rangle^2}\Big)\widetilde{V}\\
&\left. - \overline{\oph\Big(\frac{\xi}{\langle\xi\rangle}\Big)\widetilde{V}}\cdot \Big[x_j,\oph\Big(\frac{\xi\xi_1}{\langle\xi\rangle}\Big)\Big]\widetilde{V}\right](t,x),
\end{split}
\end{equation}
so after estimates \eqref{est: bootstrap argument a-priori est} and \eqref{norm_Linfty_xv-}
\begin{equation} \label{est_xj_qw+cw}
\begin{split}
\|t(tx_j)[q_w+c_w](t,t\cdot)\|_{L^2(dx)}& \lesssim \left[\|\widetilde{V}(t,\cdot)\|_{H^1_h}+ h\|\mathcal{L}_j\widetilde{V}(t,\cdot)\|_{H^1_h}\right]\|\widetilde{V}(t,\cdot)\|_{H^{1,\infty}_h} \\
&\le CA(A+B)\varepsilon^2 t^{\frac{\delta}{2}}.
\end{split}
\end{equation}
Moreover, from \eqref{def rNF}, the fact that $x_je^{ix\cdot\xi}=D_{\xi_j}e^{ix\cdot\xi}$, integration by parts, and inequalities  \eqref{est:L2 integral D with cut-off} with $\rho=2$ (after the first part of lemma \ref{Lem_Appendix: est on Dj1j2}), \eqref{est:L2 integral partialD with cut-off}, we get that
\begin{equation} \label{norm_xnrNF}
\begin{split}
\|\chi(t^{-\sigma}D_x)&(x_nr^{NF}_w)(t,\cdot)\|_{L^2}\\ 
&\lesssim t^\beta \left[\|x_nv_{-}(t,\cdot)\|_{L^\infty}\|\textit{NL}_{kg}(t,\cdot)\|_{H^{15}}+\|V(t,\cdot)\|_{H^{15}}\|x_n\textit{NL}_{kg}(t,\cdot)\|_{L^\infty} \right.\\
&\left. + \|\textit{NL}_{kg}(t,\cdot)\|_{L^2}\left(\|V(t,\cdot)\|_{H^{13}}+ \|V(t,\cdot)\|_{H^{13,\infty}}\right) + \|V(t,\cdot)\|_{H^{13}}\|\textit{NL}_{kg}(t,\cdot)\|_{L^\infty}\right]\\
&\le CB\varepsilon t^\frac{\delta_2}{2},
\end{split}
\end{equation}
where last estimate follows from \eqref{norm_Linfty_xv-}, \eqref{xj_NLkg}, inequalities \eqref{est L2 NLkg}, \eqref{est Linfty NLkg}, \eqref{est Hs NLkg-New} with $s=15$,
and a-priori estimates \eqref{est: bootstrap argument a-priori est}.
Consequently, from \eqref{L2_norm_Mutilde}, \eqref{norm_Zn(uNf-u-)}, \eqref{est_xj_qw+cw}, \eqref{norm_xnrNF}, \eqref{est:utilde-Hs} and a-priori estimate \eqref{est: bootstrap E02} with $k=2$, we obtain \eqref{est_Mutilde-proof}.

Let us now apply
$\theta_0\big(\frac{x}{t}\big)\phi^j_k(D_x)\Omega$ to both sides of \eqref{relation between Zju and Mj utilde-new} to deduce that \small
\begin{equation}\label{L2_norm_thetaM-utilde}
\begin{split}
&\left\| \theta_0(x) \oph(\phi^j_k(\xi))\Omega_h\mathcal{M}_n\widetilde{u}(t,\cdot)\right\|_{L^2} \lesssim \|\Omega Z_n U(t,\cdot)\|_{L^2}\\
&+ \left\|\theta_0\Big(\frac{x}{t}\Big)\phi^j_k(D_x)\Omega Z_n(u^{NF}-u_{-})(t,\cdot)\right\|_{L^2}+\sum_{\mu=0}^1 \|\oph(\chi_0(h^\sigma\xi))\Omega^\mu_h\widetilde{u}(t,\cdot)\|_{L^2}\\
&+ \left\|\theta_0(x)\oph(\phi^j_k(\xi))\Omega_h[t(tx_j)(q_w+c_w)(t,tx)]\right\|_{L^2(dx)}+ \left\|\theta_0(x)\oph(\phi^j_k(\xi))\Omega_h [t (tx_n)r^{NF}_w](t,tx)]\right\|_{L^2(dx)}.
\end{split}
\end{equation}\normalsize
In order to estimate the second addend in the above right hand side we first commute $Z_n$ to $\Omega$, reminding that
\[[\Omega, Z_1]=-Z_2 \quad \text{and}\quad [\Omega, Z_2]=Z_1,\]
and use that
\begin{equation*}
\theta_0\Big(\frac{x}{t}\Big)\phi^j_k(D_x) Z_j = \Big[t\theta^j_0\Big(\frac{x}{t}\Big)\phi^j_k(D_x) + \theta_0\Big(\frac{x}{t}\Big)[\phi^j_k(D_x),x_j]\Big]\partial_t + t \theta_0\Big(\frac{x}{t}\Big)\phi^j_k(D_x)\partial_j,
\end{equation*} 
with $\theta^j_0(z):=\theta_0(z)z_j$.
Observe that commutator $[\phi^j_k(D_x),x_j]$ is bounded on $L^2$ with norm $O(t)$, and that its symbol is still supported for moderate frequencies $|\xi|\lesssim t^{-\sigma}$. Therefore, for some new $\chi\in C^\infty_0(\mathbb{R}^2)$ we have that
\begin{equation*}
\begin{split}
\left\|\theta_0\Big(\frac{x}{t}\Big)\phi^j_k(D_x)\Omega Z_n(u^{NF}-u_{-})(t,\cdot)\right\|_{L^2} &\lesssim t  \left\|\chi(t^{-\sigma}D_x)\partial_{t,x}(u^{NF}-u_{-})(t,\cdot)\right\|_{L^2}\\
&+ t\left\| \chi(t^{-\sigma}D_x) \partial_{t,x}\Omega(u^{NF}-u_{-})(t,\cdot)\right\|_{L^2},
\end{split}
\end{equation*}
so using \eqref{est:L2 integral D with cut-off} with $\rho=2$ (because of first part of lemma \ref{Lem_Appendix: est on Dj1j2}) and \eqref{est: L2 Omega integral D with cut-off}, both considered with $u=\partial_{t,x}v_\pm, v=v_\pm$, and $u=v_\pm, v=\partial_{t,x}v_\pm$, we obtain that the above right hand side is estimated by \small
\begin{equation*} 
t^{1+\beta}\left[\left(\|\partial_{t,x}V(t,\cdot)\|_{L^2}+\|\Omega \partial_{t,x}V(t,\cdot)\|_{L^2}\right)\|V(t,\cdot)\|_{H^{17,\infty}}  + \left(\|V(t,\cdot)\|_{L^2}+\|\Omega V(t,\cdot)\|_{L^2} \right)\|\partial_{t,x}V(t,\cdot)\|_{H^{17,\infty}}\right]
\end{equation*} \normalsize
From \eqref{Hs norm DtV} and \eqref{est: Hsinfty Dt V} with $s=0$, \eqref{L2_norm_DtOmegaV} and a-priori estimates, we hence deduce that
\begin{equation}\label{OmegaZ(uNF-u)}
\left\|\theta_0\Big(\frac{x}{t}\Big)\phi^j_k(D_x)\Omega Z_n(u^{NF}-u_{-})(t,\cdot)\right\|_{L^2} \le CB\varepsilon t^{\beta+\frac{\delta_2}{2}}.
\end{equation}
As concerns, instead, the estimate of the fourth $L^2$ norm in the right hand side of \eqref{L2_norm_thetaM-utilde}, we observe that from equality \eqref{x_jqw}, Leibniz rule and \eqref{Op(Sigma varphi)}
\begin{multline}\label{Omegah(xjqw)}
\left\|\theta_0(x)\oph(\phi^j_k(\xi))\Omega_h[t(tx_j)[q_w+c_w](t,tx)]\right\|_{L^2} \lesssim \sum_{\mu=0}^1h^{-\kappa}\|\widetilde{V}(t,\cdot)\|_{H^{2,\infty}_h}\|\Omega^\mu_h\widetilde{V}(t,\cdot)\|_{H^1_h} \\
 +\sum_{\mu=0}^1 h^{1-\kappa}\|\widetilde{V}(t,\cdot)\|_{H^{1,\infty}_h}\|\Omega^\mu_h\mathcal{L}_j\widetilde{V}(t,\cdot)\|_{H^1}+ h^{1-\kappa}\|\Omega_h \widetilde{V}(t,\cdot)\|_{L^\infty}\|\mathcal{L}_j \widetilde{V}(t,\cdot)\|_{H^1_h},
\end{multline}
with $\kappa=\sigma\rho$ if $\rho\ge 0$, 0 otherwise.
Using the semi-classical Sobolev injection, \eqref{Lwidetilde(V)} and the fact that $\|\Omega_h\widetilde{V}(t,\cdot)\|_{H^s_h}=\|\Omega v_{-}(t,\cdot)\|_{H^s}$ for any $s\ge 0$, together with \eqref{norm_Linfty_xv-} and a-priori estimates, we see that
\begin{equation}\label{hOmegaV LjV}
 h\|\Omega_h \widetilde{V}(t,\cdot)\|_{L^\infty}\|\mathcal{L}_j \widetilde{V}(t,\cdot)\|_{H^1_h} \lesssim \|\Omega \widetilde{V}(t,\cdot)\|_{H^2_h}\|\mathcal{L}_j \widetilde{V}(t,\cdot)\|_{H^1_h} \le CB\varepsilon t^{\frac{3\delta_2}{2}}.
\end{equation}
Also, from \eqref{relation_Zjw_Ljwidetilde(w)} with $w=v_{-}$ and $f=\textit{NL}_{kg}$\small
\begin{equation*}
\|\Omega_h\mathcal{L}_j\widetilde{V}(t,\cdot)\|_{L^2} \lesssim \|\Omega Z_jv_{-}(t,\cdot)\|_{L^2}+\sum_{\mu=0}^1\|\Omega^\mu v_{-}(t,\cdot)\|_{L^2} + \left\|\Omega \left(x_j\textit{NL}_{kg}\right)(t,\cdot)\right\|_{L^2}\le C(A+B)B\varepsilon^2 t^{\frac{1}{2}+\frac{\delta_2}{2}},
\end{equation*}\normalsize
where last inequality is obtained using \eqref{est: bootstrap Enn}, \eqref{est: bootstrap E02} and estimates \eqref{xj_NLkg}, \eqref{xjGamma_Nlkg}.
Therefore
\begin{equation*}
 h\|\widetilde{V}(t,\cdot)\|_{H^{1,\infty}_h}\|\Omega_h\mathcal{L}_j\widetilde{V}(t,\cdot)\|_{L^2} \le CAB(A+B)\varepsilon^3 t^{-\frac{1}{2}+\frac{\delta_2}{2}},
\end{equation*}
which combined with \eqref{Omegah(xjqw)}, \eqref{hOmegaV LjV} and a-priori estimates gives that
\begin{equation} \label{est_Omegah(xqw)}
\left\|\theta_0(x)\oph(\phi^j_k(\xi))\Omega_h[t(tx_j)[q_w+c_w](t,tx)]\right\|_{L^2} \le CB\varepsilon t^{\frac{3\delta_2}{2}}.
\end{equation}
We estimate the latter $L^2$ norm in \eqref{L2_norm_thetaM-utilde} recalling definition \eqref{def rNF} of $r^{NF}_w$, commutating $\Omega$ and $x_n$, and using that
\begin{equation*}
\theta_0(x) \oph(\phi^j_k(\xi))x_n = \theta^n_0(x)\oph(\phi^j_k(\xi)) + \theta_0(x)[\oph(\phi^j_k(\xi)), x_n],
\end{equation*}
where 
\begin{equation*}
[\oph(\phi^j_k(\xi)), x_n] = -i h\oph(\partial_n \phi^j_k(\xi))
\end{equation*}
is uniformly bounded on $L^2$. 
After \eqref{est L2 rNF}, \eqref{est phi(D) Omega rNF-new} with $\theta\ll 1$ small, and a-priori estimates \eqref{est: bootstrap argument a-priori est} we derive that, for some $\chi\in C^\infty_0(\mathbb{R}^2)$,
\begin{equation*}
\begin{split}
& \left\|\theta_0(x)\oph(\phi^j_k(\xi))\Omega_h [t (tx_n)r^{NF}_w](t,tx)]\right\|_{L^2(dx)} \lesssim \sum_{\mu=0}^1 t\|\chi(t^{-\sigma}D_x)\Omega^\mu r^{NF}_w(t,\cdot)\|_{L^2}\le CB\varepsilon.
\end{split}
\end{equation*}
Combining \eqref{L2_norm_thetaM-utilde}, \eqref{OmegaZ(uNF-u)}, \eqref{est_Omegah(xqw)} and above estimate together with \eqref{est: bootstrap E02}, \eqref{est:utilde-Hs}, \eqref{est_Omega-utilde_proof}, and assuming $3\delta_2\le \delta_1$, we finally obtain \eqref{est_Omega-M-utilde-proof} and the conclusion of the proof.
\endproof
\end{lem}

In the following lemma we explain how we estimate the $L^2$ or the $L^\infty$ norm of products supported for moderate frequencies $|\xi|\lesssim t^\sigma$, when we have a control on high Sobolev norms of, at least, all factors but one. This type of estimate will be frequently used in most of the results that follow.

\begin{lem}\label{Lem_appendix:L_estimate of products}
Let $n\in\mathbb{N}$, $n\ge 2$, and $w_1,\dots, w_n$ such that $w_1\in L^2(\mathbb{R}^2)$, $w_2,\dots,w_n\in L^\infty(\mathbb{R}^2)\cap H^s(\mathbb{R}^2)$, for some large positive $s$.
Let also $\chi\in C^\infty_0(\mathbb{R}^2)$ and $\sigma>0$. There exists some $\chi_1\in C^\infty_0(\mathbb{R}^2)$, equal to 1 on the support of $\chi$, such that for $L=L^2$ or $L=L^\infty$
\begin{multline*}
\left\|\chi(t^{-\sigma}D_x)\left[w_1\dots w_n\right]\right\|_{L} \lesssim \left\|\left[\chi_1(t^{-\sigma}D_x)w_1\right] \prod_{j=2}^n \chi(t^{-\sigma}D_x)w_j\right\|_{L(dx)} \\
+ t^{-N(s)}\|w_1\|_{L^2(dx)}\sum_{j=2}^n\prod_{k\ne j}\|w_k\|_{L^\infty}\|w_j\|_{H^s(dx)},
\end{multline*}
with $N(s)$ as large as we want as long as $s>0$ is large.
\proof
The idea of the proof is to decompose each factor $w_j$, for $j=2,\dots,n$ into
\begin{equation}\label{dec_wj}
\chi(t^{-\sigma}D_x)w_j + (1-\chi)(t^{-\sigma}D_x)w_j,
\end{equation}
and to estimate the $L^2$ norm of product
\begin{equation}\label{prod_trunc}
\chi(t^{-\sigma}D_x)\left[w_1\prod_{\substack{k=2,\dots,n\\ k\ne j}}\widetilde{w}_k \left[ (1-\chi)(t^{-\sigma}D_x)w_j \right]\right],
\end{equation}
where $\widetilde{w}_k$ is either $w_k$ or $\chi(t^{-\sigma}D_x)w_k$, with the $L^2$ norm of $w_1$ times the $L^\infty$ norm of all remaining factors, reminding that $\chi(t^{-\sigma}D_x)$ is uniformly bounded on $L^\infty$ and that by Sobolev injection and \eqref{ineq:1-chi},
\begin{equation}
\left\|(1-\chi)(t^{-\sigma}D_x)w_j\right\|_{L^\infty(dx)}\lesssim t^{-N(s)}\|w_j\|_{H^s(dx)},
\end{equation}
with $N(s)$ as large as we want as long as $s>0$ is large.
The $L^\infty$ norm of \eqref{prod_trunc} is estimated in the same way, using firstly the $L^2-L^\infty$ continuity of operator $\chi(t^{-\sigma}D_x)$ acting on the entire product.

The end of the statement follows from the observation that, if $\text{supp}\chi\subset B_C(0)$ for some $C>0$, then
\begin{equation} \label{property_support}
\text{supp}\hat{w}_1\subset\{\xi : |\xi|\ge C_1>nC\} \quad \Rightarrow \quad\chi(t^{-\sigma}D_x)\Big[w_1\prod_{j=2}^n\chi(t^{-\sigma}D_x)w_j\Big]\equiv 0.
\end{equation}
\endproof
\end{lem}

\begin{remark}
Property \eqref{property_support} is more general, meaning that if $\chi,\chi_j\in C^\infty_0(\mathbb{R}^2)$ with $\text{supp}\chi\subset B_C(0)$, $\text{supp}\chi_j\subset B_{C_j}(0)$ for some $C,C_j>0$, for every $j=2,\dots,n$, then
\[\text{supp}\hat{w}_1\subset \Big\{\xi : |\xi|\ge C_1>C+\sum_{j=2}^n C_j\Big\} \quad \Rightarrow \quad\chi(t^{-\sigma}D_x)\Big[w_1\prod_{j=2}^n\chi_j(t^{-\sigma}D_x)w_j\Big]\equiv 0.\]
\end{remark}

We have seen at the beginning of section \ref{Sub: App_B1}, and already used in the previous lemma's proof, that, if $w\in H^s(\mathbb{R}^2)$ for some large $s>0$, the $L^2$ norm (resp. $L^\infty$ norm) of this function when restricted to large frequencies $|\xi|\gtrsim t^\sigma$ decays fast in time as $t^{-\sigma s}$ (resp. $t^{-\sigma(s-1)-1}$ after the semi-classical Sobolev injection).
The aim of the following lemma is to show that, even if we don't have a control on the $H^s(\mathbb{R}^2)$ norm of $(\Gamma u)_\pm$, $(\Gamma v)_\pm$, for $\Gamma\in\{\Omega,Z_m, m=1,2\}$ and $s$ larger than 2, the $L^2$ norm (resp. $L^\infty$ norm) of products as in \eqref{prod_app_w1wn} still have a good decay in time.

\begin{lem} \label{Lem_app:products_Gamma}
Let $w\in \{u,v\}$ and for any $\Gamma\in \{\Omega, Z_m , m=1,2\}$
\begin{equation*}
(\Gamma w)_\pm = 
\begin{cases}
(D_t\pm |D_x|)(\Gamma u), \quad & \text{if } w=u,\\
(D_t\pm \langle D_x\rangle)(\Gamma v), \quad & \text{if } w=v.
\end{cases}
\end{equation*}
Let also $n\in\mathbb{N}^*$, $w_1,\dots, w_n$ be such that $w_1, xw_1\in L^2(\mathbb{R}^2)\cap L^\infty(\mathbb{R}^2)$, $w_j\in L^\infty(\mathbb{R}^2)$ for $j=2,\dots,n$, $\chi\in C^\infty_0(\mathbb{R}^2)$, $\sigma>0$, and $a(D_x)= D_x^\alpha\big(\frac{D_x}{\langle D_x\rangle}\big)^\beta \big(\frac{D_x}{|D_x|}\big)^\gamma$ for any $\alpha,\beta,\gamma\in\mathbb{N}^2$ with $|\alpha|, |\beta|, |\gamma|\le 1$.
Then for $L=L^2$ or $L=L^\infty$ we have that
\begin{subequations}\label{lem_prod_Omega_Z}
\begin{equation} \label{lem_prod_Omegaw}
\begin{split}
\left\| a(D_x)(\Omega w)_\pm w_1\dots w_n \right\|_{L(dx)} &\lesssim \left\|\left[\chi(t^{-\sigma}D_x)a(D_x) (\Omega w)_\pm\right] \prod_{j=1}^n w_j\right\|_{L(dx)}\\
&+t^{-N(s)} \|w_\pm(t,\cdot)\|_{H^s}\Big(\sum_{|\mu|=0}^1\|x^\mu w_1\|_{L(dx)}\Big)\prod_{j=2}^n \|w_j\|_{L^\infty(dx)}
\end{split}
\end{equation}
and, for $m=1,2,$
\begin{multline}\label{lem_prod_Zmw}
\left\| a(D_x)(Z_m w)_\pm w_1\dots w_n \right\|_{L(dx)} \lesssim \left\|\left[\chi(t^{-\sigma}D_x)a(D_x) (Z_m w)_\pm\right] \prod_{j=1}^n w_j\right\|_{L(dx)}\\
+t^{-N(s)}\left( \|w_\pm(t,\cdot)\|_{H^s}+\|D_tw_\pm(t,\cdot)\|_{H^s}\right)\Big(\sum_{\mu=0}^1\|x^\mu_m w_1\|_{L(dx)}+ t\|w_1\|_{L(dx)}\Big)\prod_{j=2}^n \|w_j\|_{L^\infty(dx)},
\end{multline}
\end{subequations}
with $N(s)$ as large as we want as long as $s>0$ is large.
\proof
Let us remind definition \eqref{Omega, Zj} of Klainerman vector fields $\Omega, Z_m$, for $m=1,2$, and decompose factor $a(D_x)(\Gamma w)_\pm$ in frequencies by means of operator $\chi(t^{-\sigma}D_x)$.
When dealing with product
\begin{equation}\label{prod_app_w1wn}
\big[(1-\chi)(t^{-\sigma}D_x)a(D_x)(\Gamma w)_\pm\big]w_1\cdots w_n
\end{equation}
the idea is to discharge on $w_1$ factors $x$ and/or $t$ defining $\Gamma$, after a previous commutation between $D_t\pm |D_x|$ if $w=u$ (resp. $D_t\pm \langle D_x\rangle$ if $w=v$) and $\Gamma$, and between $(1-\chi)(t^{-\sigma}D_x)a(D_x)$ and the mentioned factors $x, t$.
For instance, if $w=u$ and $\Gamma=Z_1$
\begin{multline}
\big[(1-\chi)(t^{-\sigma}D_x)a(D_x)(Z_1 u)_\pm\big] w_1  = \big[ (1-\chi)(t^{-\sigma}D_x)a(D_x)(\partial_t u)_\pm\big](x_1 w_1)\\
 + \big[(1-\chi)(t^{-\sigma}D_x)a(D_x)(\partial_1 u)_\pm\big] (tw_1)+ \Big[(1-\chi)(t^{-\sigma}D_x)a(D_x)\frac{D_1}{|D_x|} u_\pm\Big] w_1 \\
 + \Big[\big[(1-\chi)(t^{-\sigma}D_x)a(D_x), x_1\big]D_t u_\pm\Big] w_1,
\end{multline}\normalsize
from which we deduce, using the Sobolev injection together with \eqref{ineq:1-chi}, that \small
\[
\left\|\big[ (1-\chi)(t^{-\sigma}D_x)a(D_x)(Z_1 u)_\pm\big] w_1\right\|_{L}
\lesssim t^{-N(s)}\left(\|u_\pm(t,\cdot)\|_{H^s}+\|D_t u(t,\cdot)\|_{H^s}\right)\Big(\sum_{\mu=0}^1 \|x_1^\mu w_1\|_{L} + t\|w_1\|_{L}\Big),
\]\normalsize
with $N(s)$ large as long as $s$ is large.
Analogous inequalities can be obtained for $\Gamma=\Omega, Z_2$ and/or $w=v$. 
This concludes the proof of the statement since the $L$ norm of \eqref{prod_app_w1wn} is bounded by the $L$ norm of $\big[(1-\chi)(t^{-\sigma}D_x)a(D_x)(Z_1 u)_\pm\big] w_1$ times the $L^\infty$ norm of the remaining factors.
\endproof
\end{lem}

\begin{cor}
If the hypothesis of lemma \ref{Lem_app:products_Gamma} are satisfied and in addition $w_1,\dots,w_n \in H^s(\mathbb{R}^2)$, we have that
\begin{subequations} \label{cor_app_est_1}\small
\begin{equation}
\begin{split} 
\left\| a(D_x)(\Omega w)_\pm w_1\cdots w_n \right\|_{L} &\lesssim \left\|\left[\chi(t^{-\sigma}D_x)a(D_x) (\Omega w)_\pm\right] \prod_{j=1}^n \chi(t^{-\sigma}D_x)w_j\right\|_{L}\\
&+t^{-N(s)} \|w_\pm(t,\cdot)\|_{H^s(dx)}\Big(\sum_{|\mu|=0}^1\|x^\mu w_1\|_{L}\Big)\prod_{j=2}^n \|w_j\|_{L^\infty}\\
&+ t^{-N(s)} \|(\Omega w)_\pm(t,\cdot)\|_{L^2}\sum_{j=1}^n \prod_{k\ne j}\|w_k\|_{L^\infty}\|w_j\|_{H^s}
\end{split}
\end{equation}\normalsize
and, for $m=1,2$,
\begin{equation}
\begin{split}
&\left\| a(D_x)(Z_m w)_\pm w_1\cdots w_n \right\|_{L} \lesssim \left\|\left[\chi(t^{-\sigma}D_x)a(D_x) (Z_m w)_\pm\right] \prod_{j=1}^n \chi(t^{-\sigma}D_x)w_j\right\|_{L}\\
& +t^{-N(s)}\left( \|w_\pm(t,\cdot)\|_{H^s}+\|D_tw_\pm(t,\cdot)\|_{H^s}\right)\Big(\sum_{\mu=0}^1\|x^\mu_m w_1\|_{L}+ t\|w_1\|_{L}\Big)\prod_{j=2}^n \|w_j\|_{L^\infty}\\
&+ t^{-N(s)} \|(Z_m w)_\pm(t,\cdot)\|_{L^2}\sum_{j=1}^n \prod_{k\ne j}\|w_k\|_{L^\infty}\|w_j\|_{H^s},
\end{split}
\end{equation}
\end{subequations}
with $N(s)$ as large as we want as long as $s>0$ is large.
Moreover, there exists $\chi_1\in C^\infty_0(\mathbb{R}^2)$ such that, for any fixed $j_0\in \{1,\dots,n\}$,
\begin{subequations}\label{cor_app_est_2}
\begin{equation}\label{cor_estOmega_2}
\begin{split}
&\left\|\chi(t^{-\sigma}D_x)\big[a(D_x)(\Omega w)_\pm w_1\cdots w_n\big] \right\|_{L} \\
&\lesssim \left\|\big[\chi(t^{-\sigma}D_x)a(D_x)(\Omega w)_\pm \big]\big[\chi_1(t^{-\sigma}D_x)w_{j_0}\big]\prod_{\substack{j=1,\dots,n\\j\ne j_0}}\chi(t^{-\sigma}D_x)w_j \right\|_{L} \\
&+t^{-N(s)} \|w_\pm(t,\cdot)\|_{H^s}\Big(\sum_{|\mu|=0}^1\|x^\mu w_1\|_{L}\Big)\prod_{j=2}^n \|w_j\|_{L^\infty}\\
&+ t^{-N(s)} \|(\Omega w)_\pm(t,\cdot)\|_{L^2}\sum_{\substack{j=1,\dots,n \\ j\ne j_0}} \prod_{k\ne j}\|w_k\|_{L^\infty}\|w_j\|_{H^s}
\end{split}
\end{equation}
and, for $m=1,2$,
\begin{equation} \label{cor_estZm_2}
\begin{split}
&\left\|\chi(t^{-\sigma}D_x)\big[a(D_x)(Z_m w)_\pm w_1\cdots w_n\big] \right\|_{L}\\
& \lesssim \left\|\big[\chi(t^{-\sigma}D_x)a(D_x) (Z_m w)_\pm \big]\big[\chi_1(t^{-\sigma}D_x)w_{j_0}\big]\prod_{\substack{j=1,\dots,n\\j\ne j_0}}\chi(t^{-\sigma}D_x)w_j \right\|_{L} \\
& +t^{-N(s)}\left( \|w_\pm(t,\cdot)\|_{H^s}+\|D_tw_\pm(t,\cdot)\|_{H^s}\right)\Big(\sum_{\mu=0}^1\|x^\mu_m w_1\|_{L}+ t\|w_1\|_{L}\Big)\prod_{j=2}^n \|w_j\|_{L^\infty}\\
&+ t^{-N(s)} \|(Z_m w)_\pm(t,\cdot)\|_{L^2}\sum_{\substack{j=1,\dots,n \\ j\ne j_0}} \prod_{k\ne j}\|w_k\|_{L^\infty}\|w_j\|_{H^s}.
\end{split}
\end{equation}
\end{subequations}
\proof
The inequalities of the statement mainly follows from \eqref{lem_prod_Omega_Z}. In fact, by decomposing each factor $w_j$ appearing in the first norm in the right hand sides of \eqref{lem_prod_Omega_Z} as in \eqref{dec_wj}, and then 
using the following inequality, for $\Gamma\in \{\Omega, Z_m, m=1,2\}$ and $\widetilde{w}_k$ either equal to $w_k$ or to $\chi(t^{-\sigma}D_x)w_k$,
\begin{multline*}
\left\| [\chi(t^{-\sigma}D_x)a(D_x)(\Gamma w)_\pm] \prod_{\substack{k=1,\dots,n\\ k\ne j}}\widetilde{w}_k \left[ (1-\chi)(t^{-\sigma}D_x)w_j \right]\right\|_{L} \\
\lesssim t^{-N(s)}\|(\Gamma w)_\pm(t,\cdot)\|_{L^2}\prod_{\substack{k=1,\dots,n\\ k\ne j}}\|w_k\|_{L^\infty}\|w_j\|_{H^s},
\end{multline*}
with $N(s)$ as large as we want as long as $s>0$,
which is obtained from \eqref{ineq:1-chi} together with the $L^2-L^\infty$ and $L^\infty-L^\infty$ continuity of operator $\chi(t^{-\sigma}D_x)$, we obtain \eqref{cor_app_est_1}.

On the other hand, if the product in the left hand side of \eqref{lem_prod_Omega_Z} is localized in frequencies by means of operator $\chi(t^{-\sigma}D_x)$, so it is for the product in the first norm of the same inequalities. 
Inequalities \eqref{cor_app_est_2} are then derived by bounding these $L$ norms by means lemma \ref{Lem_appendix:L_estimate of products}, where the role of $w_1$ is here played by $w_{j_0}$, for some fixed $j_0\in \{1,\dots,n\}$.
\endproof
\end{cor}

%

The following two lemmas are stated and proved in view of lemma \ref{Lem_appendix: preliminary est VJ}, in which we recover a first non-sharp estimate on the $L^\infty$ norm of the Klein-Gordon component when one Klainerman vector field is acting on it and its frequencies are less or equal than $t^\sigma$, for some small $\sigma>0$.
This estimate will be successively refined in lemma \ref{Lem_appendix: sharp_est_VJ}.

\begin{lem} \label{Lem_appendix:Linfty_bound_chi_w}
Let $\chi\in C^\infty_0(\mathbb{R}^2)$, $\sigma>0$ small, and $w=w(t,x)$ such that,
if $\widetilde{w}(t,x):=tw(t,tx)$, $\oph(\chi(h^\sigma\xi))\mathcal{L}^\mu \widetilde{w}(t,\cdot)\in L^2(\mathbb{R}^2)$ for any $|\mu|\le 1$. Then
\begin{equation} \label{ineq:norm_Linfty_chi-w}
\left\| \chi(t^{-\sigma}D_x)w(t,\cdot)\right\|_{L^\infty} \lesssim t^{-1+\beta} \sum_{|\mu|=0}^1\left\|\oph(\chi(h^\sigma\xi))\mathcal{L}^\mu \widetilde{w}(t,\cdot) \right\|_{L^2},
\end{equation}
with $\beta>0$ small, $\beta\rightarrow 0$ as $\sigma\rightarrow 0$.
\proof
Since
\begin{equation*}
\chi(t^{-\sigma}D_x)w(t,y) = t^{-1}\oph(\chi(h^\sigma\xi))\widetilde{w}(t,x)|_{x=\frac{y}{t}},
\end{equation*}
the goal is to prove that
\begin{equation} \label{norm_Linfty_vtildeGamma}
\left\|\oph(\chi(h^\sigma\xi))\widetilde{w}(t,\cdot) \right\|_{L^\infty}\lesssim  h^{-\beta}\sum_{|\mu|=0}^1\left\|\oph(\chi(h^\sigma\xi))\mathcal{L}^\mu \widetilde{w}(t,\cdot) \right\|_{L^2},
\end{equation}
for a small $\beta>0$, $\beta\rightarrow 0$ as $\sigma\rightarrow 0$. 
So let $w^\chi:=\oph(\chi(h^{\sigma}\xi)) \widetilde{w}$ and take $\chi_1\in C^\infty_0(\mathbb{R}^2)$ equal to 1 on the support of $\chi$, so that
\begin{equation*}
\oph(\chi(h^{\sigma}\xi)) \widetilde{w} = \oph(\chi_1(h^{\sigma}\xi)) \widetilde{w}^\chi.
\end{equation*}
For a $\gamma\in C^\infty_0(\mathbb{R}^2)$, equal to 1 in a neighbourhood of the origin and with sufficiently small support, we consider the following decomposition
\begin{equation*}
\oph\Big(\gamma\Big(\frac{x-p'(\xi)}{\sqrt{h}}\Big)\chi_1(h^\sigma\xi)\Big)\widetilde{w}^\chi + \oph\Big((1-\gamma)\Big(\frac{x-p'(\xi)}{\sqrt{h}}\Big)\chi_1(h^\sigma\xi)\Big)\widetilde{w}^\chi
\end{equation*}
and immediately observe that, from inequality \eqref{est_1L-Linfty},
\begin{equation*}
\left\| \oph\Big((1-\gamma)\Big(\frac{x-p'(\xi)}{\sqrt{h}}\Big)\chi_1(h^\sigma\xi)\Big)\widetilde{w}^\chi(t,\cdot) \right\|_{L^\infty}\lesssim h^{-\beta}\sum_{|\mu|=0}^1\| \oph(\chi_1(h^\sigma\xi))\mathcal{L}^\mu \widetilde{w}^\chi(t,\cdot)\|_{L^2}.
\end{equation*}
After lemma \ref{Lem:family_thetah} there exists a family of smooth cut-off functions $\theta_h(x)$ such that equality \eqref{cut-off-thetah} holds. Then, if we also consider a new cut-off function $\chi_2$ equal to 1 on the support of $\chi_1$ and a small $\sigma_1>\sigma$, by symbolic calculus and remark \ref{Remark:symbols_with_null_support_intersection}
we derive that for any $N\in\mathbb{N}$
\begin{equation*}
\begin{split}
\oph\Big(\gamma\left(\frac{x-p'(\xi)}{\sqrt{h}}\right)\chi_1(h^\sigma\xi)\Big)\widetilde{w}^\chi & = \theta_h(x)\oph(\chi_2(h^\sigma\xi)) \oph\Big(\gamma\left(\frac{x-p'(\xi)}{\sqrt{h}}\right)\chi_1(h^\sigma\xi)\Big)\widetilde{w}^\chi\\
& + \oph(r_\infty(x,\xi))\widetilde{w}^\chi + \theta_h(x) \oph(r^1_\infty(x,\xi))\widetilde{w}^\chi
\end{split}
\end{equation*}
with $r_\infty, r^1_\infty\in h^N S_{\frac{1}{2},\sigma}(\langle \frac{x-p'(\xi)}{\sqrt{h}}\rangle^{-1})$. It is enough to take $N=1$ to have, by proposition \ref{Prop : Continuity from $L^2$ to L^inf}, that 
\[\|\oph(r_\infty)\widetilde{w}^\chi(t,\cdot)\|_{L^\infty} + \|\theta_h(x)\oph(r^1_\infty)\widetilde{w}^\chi(t,\cdot)\|_{L^\infty}\le h^{-\beta}\|\widetilde{w}^\chi(t,\cdot)\|_{L^2}.\]
As function $\phi(x):=\sqrt{1-|x|^2}$ is well defined on the support of $\theta_h$ we are allowed to to write the following:
\begin{align*}
& \left\| \theta_h(x) \oph(\chi_2(h^{\sigma_1}\xi))\oph\Big(\gamma\left(\frac{x-p'(\xi)}{\sqrt{h}}\right)\chi_1(h^\sigma\xi)\Big)\widetilde{w}^\chi(t,\cdot)\right\|_{L^\infty} \\
& = \left\| e^{\frac{i}{h}\phi}\theta_h(x) \oph(\chi_2(h^{\sigma_1}\xi)) \oph\Big(\gamma\left(\frac{x-p'(\xi)}{\sqrt{h}}\right)\chi_1(h^\sigma\xi)\Big)\widetilde{w}^\chi(t,\cdot)\right\|_{L^\infty} \\
& \lesssim \left\| \oph(\chi_2(h^{\sigma_1}\xi))\left[ e^{\frac{i}{h}\phi}\theta_h(x)  \oph\Big(\gamma\left(\frac{x-p'(\xi)}{\sqrt{h}}\right)\chi_1(h^\sigma\xi)\Big)\widetilde{w}^\chi\right]\right\|_{L^\infty} + \|\oph(r_\infty)\widetilde{w}^\chi(t,\cdot)\|_{L^\infty},
\end{align*}
for a new $r_\infty\in h^NS_{\frac{1}{2},\sigma}\big(\langle\frac{x-p'(\xi)}{\sqrt{h}}\rangle^{-1}\big)$. This latter $r_\infty$ comes out from the commutation between $e^{\frac{i}{h}\phi}\theta_h(x)$ and $\oph(\chi_2(h^{\sigma_1}\xi))$, whose symbol is computed using \eqref{a sharp b asymptotic formula} until a large enough order $M$. We notice that we gain a factor $h^{|\alpha|(\sigma_1-\sigma)}$ at each order of the mentioned asymptotic development as $\sigma_1>\sigma$.
Moreover, those terms write in terms of the derivatives of $\chi_2$ and hence vanish on the support of $\chi_1$. By proposition \ref{Prop: a sharp b} and remark \ref{Remark:symbols_with_null_support_intersection} we then deduce that the composition of the mentioned commutator with $\oph\big(\gamma\big(\frac{x-p'(\xi)}{\sqrt{h}}\big)\chi_1(h^\sigma\xi)\big)$ is an operator of symbol $r_\infty$, with $N$ as large as we want.

Using the classical Sobolev injection, symbolic calculus and lemma \ref{Lem: (xi+dphi)Op(gamma)} we find that\small
\begin{align*}
& \left\| \oph(\chi_2(h^{\sigma_1}\xi))\left[ e^{\frac{i}{h}\phi}\theta_h(x)  \oph\Big(\gamma\left(\frac{x-p'(\xi)}{\sqrt{h}}\right)\chi_1(h^\sigma\xi)\Big)\widetilde{w}^\chi\right]\right\|_{L^\infty}\\
& \lesssim |\log h| \left[ \|\widetilde{w}^\chi(t,\cdot)\|_{L^2} + \sum_{j=1}^2\left\| D_j\left[ e^{\frac{i}{h}\phi}\theta_h(x)  \oph\Big(\gamma\left(\frac{x-p'(\xi)}{\sqrt{h}}\right)\chi_1(h^\sigma\xi)\Big)\widetilde{w}^\chi\right]\right\|_{L^2} \right] \\
&\lesssim |\log h| \left[ \|\widetilde{w}^\chi(t,\cdot)\|_{L^2} + \sum_{j=1}^2 h^{-1}\left\| \oph\big((\xi_j+d_j\phi(x))\theta_h(x)\big)\oph\Big(\gamma\left(\frac{x-p'(\xi)}{\sqrt{h}}\right)\chi_1(h^\sigma\xi)\Big)\widetilde{w}^\chi\right\|_{L^2}\right] \\
&\lesssim |\log h| \left[ \|\widetilde{w}^\chi(t,\cdot)\|_{L^2} + h^{-\beta}\sum_{|\mu|=0}^1 \left\|\oph(\chi_1(h^\sigma\xi))\mathcal{L}^\mu \widetilde{w}^\chi(t,\cdot)\right\|_{L^2}\right].
\end{align*}\normalsize
Finally, commutating $\mathcal{L}$ with $\oph(\chi(h^\sigma\xi))$ defining $\widetilde{w}^\chi$, and reminding that $\chi_1\equiv 1$ on the support of $\chi$, we obtain
\begin{equation*}
\|\oph(\chi(h^\sigma\xi))\widetilde{w}^\chi(t,\cdot)\|_{L^\infty}\lesssim h^{-\beta}\sum_{|\mu|=0}^1\|\oph(\chi(h^\sigma\xi))\mathcal{L}^\mu \widetilde{w}(t,\cdot)\|_{L^2},
\end{equation*}
for every $t\in[1,T]$, and hence \eqref{norm_Linfty_vtildeGamma}.
\endproof
\end{lem}

\begin{lem} \label{Lem_appendix:Intro_of_vNFGamma}
Let $I$ be a multi-index of length $j$, with $j=1,2$, and \index{vNFI@$\vNFGamma$, normal form function defined from $(\Gamma^Iv)_{-}$}
\begin{equation} \label{def_vNF-Gamma}
\vNFGamma(t,x): = (\Gamma^I v)_{-}(t,x)-\frac{i}{4(2\pi)^2}\sum_{j_1,j_2\in \{+,-\}}\int e^{ix\cdot\xi}B^1_{(j_1,j_2,+)}(\xi,\eta) \widehat{v^I_{j_1}}(\xi-\eta) \hat{u}_{j_2}(\eta) d\xi d\eta,
\end{equation}
with $B^1_{(j_1,j_2,+)}$ given by \eqref{def of B(i1,i2,i3)} with $j_3=+$ and $k=1$.
Then there exists a constant $C>0$ such that, for any $\chi\in C^\infty_0(\mathbb{R}^2)$, $\sigma>0$ small, and every $t\in [1,T]$,
\begin{equation} \label{ineq:Linfty_ vINF - vI}
\left\|\chi(t^{-\sigma}D_x)\left(\vNFGamma - (\Gamma^I v)_{-}\right)(t,\cdot) \right\|_{L^\infty}\le \frac{1}{2} \left\|\chi(t^{-\sigma}D_x)(\Gamma^I v)_{-}(t,\cdot) \right\|_{L^\infty} + CB\varepsilon t^{-1}.
\end{equation}
Moreover, for every $m=1,2,$ $t\in[1,T]$,
\begin{equation} \label{ineq:L2 Zm vINF-vI}
\left\|\chi(t^{-\sigma}D_x)Z_m\left(\vNFGamma - (\Gamma^I v)_{-}\right)(t,\cdot) \right\|_{L^2} \le C(A+B)B\varepsilon^2 t^{2\sigma+\frac{\delta_{3-j}+\delta_2}{2}}.
\end{equation}
\proof
We first notice that, after definition \eqref{def_vNF-Gamma} and inequalities \eqref{explicit integral B}, \eqref{def uIpm vIpm}, we have the following explicit expression:
\begin{equation} \label{explicit vNfGamma-vJ-}
\vNFGamma-(\Gamma^I v)_{-} = -\frac{i}{2}\left[(D_t\Gamma^I v)(D_1u) - (D_1\Gamma^I v)(D_tu) + D_1[(\Gamma^I v) D_tu] - \langle D_x\rangle [(\Gamma^I v) D_1u] \right].
\end{equation}
From the above equality together with lemma \ref{Lem_appendix:L_estimate of products} with $L=L^\infty$ and $w_1=(\Gamma^I v)_\pm$, and \eqref{def u+- v+-}, \eqref{def uIpm vIpm}, we deduce that there exists some $\chi_1\in C^\infty_0(\mathbb{R}^2)$, equal to 1 on the support of $\chi$, and $s>0$ sufficiently large such that
\begin{equation}\label{est_1_vNFGamma-vJ}
\begin{split}
&\left\|\chi(t^{-\sigma}D_x)(\vNFGamma - v^I_{-})(t,\cdot) \right\|_{L^\infty}\\
&\lesssim t^\sigma \sum_{\mu=0}^1\left\|[\chi_1(t^{-\sigma}D_x)(\Gamma^I v)_\pm(t,\cdot)] [\chi(t^{-\sigma}D_x)\mathrm{R}^\mu_1u_\pm](t,\cdot)\right\|_{L^\infty}\+ t^{-2}\|(\Gamma^I v)_\pm(t,\cdot)\|_{L^2}\|u_\pm(t,\cdot)\|_{H^s}\\
&\lesssim t^\sigma \sum_{\mu=0}^1\left\|[\chi_1(t^{-\sigma}D_x)(\Gamma^I v)_\pm(t,\cdot)] [\chi(t^{-\sigma}D_x)\mathrm{R}^\mu_1u_\pm](t,\cdot)\right\|_{L^\infty} + B^2\varepsilon^2 t^{-3/2},
\end{split}
\end{equation}
where the latter inequality follows from after a-priori energy estimates \eqref{est: bootstrap Enn}, \eqref{est: bootstrap E02}.
Our aim is to truncate factor $(\Gamma^I v)_\pm$ in the above right hand side rather with the same operator $\chi(t^{-\sigma}D_x)$ appearing on the left hand side. 
We hence proceed by picking some $\kappa\ge 1$ and decomposing $\chi(t^{-\sigma}D_x)\mathrm{R}^\mu_1 u_\pm$ as
\begin{equation}\label{dec_small_frequencies_Ru}
\chi(t^{-\sigma}D_x)\mathrm{R}^\mu_1 u_\pm = \chi(t^\kappa D_x) \mathrm{R}^\mu_1 u_\pm + (1-\chi)(t^\kappa D_x)\chi(t^{-\sigma}D_x)\mathrm{R}^\mu_1 u_\pm,
\end{equation}
noticing that, as $\chi(t^\kappa \xi)$ is supported for very small frequencies $|\xi|\lesssim t^{-\kappa}$, by Sobolev injection we have that
\begin{equation*}
\left\| \chi(t^\kappa D_x)\mathrm{R}^\mu_1 u_\pm(t,\cdot)\right\|_{L^\infty}\lesssim t^{-\kappa}\|u_\pm(t,\cdot)\|_{L^2}.
\end{equation*}
Consequently, using the $L^2-L^\infty$ continuity of $\chi_1(t^{-\sigma}D_x)$ along with a-priori estimates \eqref{est: bootstrap Enn}, \eqref{est: bootstrap E02}, we have that for any for $\mu=0,1$
\begin{align*}
\left\|[\chi_1(t^{-\sigma}D_x)(\Gamma^I v)_\pm(t,\cdot)] \chi(t^\kappa D_x)\mathrm{R}^\mu_1u_\pm(t,\cdot)\right\|_{L^\infty}&\lesssim t^{\sigma-\kappa} \|(\Gamma^I v)_\pm (t,\cdot)\|_{L^2}\|u_\pm(t,\cdot)\|_{L^2}\\
&\le CB\varepsilon t^{-\kappa+\sigma + \frac{\delta_{3-j}+\delta}{2}}.
\end{align*}
Choosing $\kappa=1+\sigma +\frac{\delta+\delta_1}{2}$ we deduce from \eqref{dec_small_frequencies_Ru} and the above inequality that
\begin{multline}\label{cut_2}
\left\|[\chi_1(t^{-\sigma}D_x)(\Gamma^I v)_\pm(t,\cdot)] \chi(t^{-\sigma}D_x)\mathrm{R}^\mu_1u_\pm(t,\cdot)\right\|_{L^\infty} \\
\lesssim \left\|[\chi_1(t^{-\sigma}D_x)(\Gamma^I v)_\pm(t,\cdot)] (1-\chi)(t^\kappa D_x)\chi(t^{-\sigma}D_x)\mathrm{R}^\mu_1u_\pm(t,\cdot)\right\|_{L^\infty} +CB\varepsilon t^{-1}.
\end{multline}
We then decompose $(\Gamma^Iv)_\pm$ in frequencies using the wished operator $\chi(t^{-\sigma}D_x)$.
In order to estimate the $L^\infty$ norm of 
\begin{equation*}
[(1-\chi)(t^{-\sigma} D_x)\chi_1(t^{-\sigma}D_x)(\Gamma^I v)_\pm] (1-\chi)(t^\kappa D_x)\chi(t^{-\sigma}D_x)\mathrm{R}^\mu_1u_\pm
\end{equation*}
we first commute $\Gamma^I$ to operator $D_t\pm \langle D_x\rangle$ (see \eqref{commutator_Z_Dt-<D>}) and successively look at it as a linear combination of derivations of the form $x^\alpha t^a\partial^\alpha_x\partial^b_t$, with $1\le |\alpha| + a\le 2$, $1\le |\beta| + b\le 2$.
By commutating $x^\alpha$ to $(1-\chi)(t^{-\sigma} D_x)\chi_1(t^{-\sigma}D_x)$, multiplying it against the wave factor, and successively combining the classical Sobolev injection with inequality \eqref{ineq:1-chi}, we find that
\begin{multline}\label{cut_1}
\left\|[(1-\chi)(t^{-\sigma} D_x)\chi_1(t^{-\sigma}D_x)(\Gamma^I v)_\pm(t,\cdot)] (1-\chi)(t^\kappa D_x)\chi(t^{-\sigma}D_x)\mathrm{R}^\mu_1u_\pm(t,\cdot)\right\|_{L^\infty} \\
\lesssim t^{-N(s)}\left(\|v_\pm(t,\cdot)\|_{H^s}+\|D_tv_\pm(t,\cdot)\|_{H^s} + \|D^2_tv_\pm(t,\cdot)\|_{H^s}\right) \\
\times \sum_{\substack{1\le |\alpha|+ a \le 2\\ |\mu|=0,1}}\left\|x^\alpha t^a (1-\chi)(t^\kappa D_x)\chi(t^{-\sigma}D_x)\mathrm{R}^\mu_1u_\pm\right\|_{L^\infty}.
\end{multline}
Using system \eqref{system for uI+-, vI+-} with $|I|=0$, \eqref{est Hs NLkg-New} with $s=1$, \eqref{Hs norm DtV} and a-priori estimates, it is straightforward to check that
\begin{equation} \label{est:A}
\|v_\pm(t,\cdot)\|_{H^s}+\|D_tv_\pm(t,\cdot)\|_{H^s} + \|D^2_tv_\pm(t,\cdot)\|_{H^s} \le CB\varepsilon t^\frac{\delta}{2}.
\end{equation}
Also,
\begin{equation}\label{est:B}
t^a\left\| (1-\chi)(t^\kappa D_x)\chi(t^{-\sigma}D_x)\mathrm{R}^\mu_1u_\pm\right\|_{L^\infty}\lesssim t^{a+\sigma}\|u_\pm(t,\cdot)\|_{L^2}\le CB\varepsilon t^{a+\sigma+\frac{\delta}{2}},
\end{equation}
and for $|\alpha|\in \{1,2\}$ we have that
\begin{equation}\label{ineq:x_alpha_R1u}
\left\|x^\alpha (1-\chi)(t^\kappa D_x)\chi(t^{-\sigma}D_x)\mathrm{R}^\mu_1u_\pm\right\|_{L^\infty}\le CB\varepsilon t^{|\alpha|+|\alpha|\kappa+\frac{\delta}{2}}.
\end{equation}
In fact, when $|\alpha|=1$ this can be proved by commutating $x^\alpha$ with $(1-\chi)(t^\kappa D_x)\chi(t^{-\sigma}D_x)$, using that
\begin{equation*}
[x_n, (1-\chi)(t^\kappa D_x)\chi(t^{-\sigma}D_x)] =-it^\kappa (\partial_n\chi)(t^\kappa D_x) +i t^{-\sigma}(\partial_n\chi)(t^{-\sigma}D_x), \quad n=1,2,
\end{equation*}
is bounded from $L^2$ to $L^\infty$ uniformly in $t$, and together with estimates \eqref{est: bootstrap E02}, \eqref{est_L2_xNLw}, and the following inequality
\begin{multline*}
\left\| (1-\chi)(t^\kappa D_x)\chi(t^{-\sigma}D_x)\big[x^\alpha \mathrm{R}^\mu_1u_\pm\big](t,\cdot)\right\|_{L^\infty} \\
\lesssim t^\kappa \left[\sum_{|\mu|=1}\|Z^\mu u_\pm(t,\cdot)\|_{L^2} + t\|u_\pm(t,\cdot)\|_{H^1} + \|x\textit{NL}_w(t,\cdot)\|_{L^2}\right],
\end{multline*}
which is obtained by writing 
\begin{align*}
&(1-\chi)(t^\kappa D_x)\chi(t^{-\sigma}D_x) x_n \mathrm{R}^\mu_1\\
&= t^\kappa \widetilde{\chi}_1(t^\kappa D_x)\chi(t^{-\sigma}D_x) x_n|D_x|\mathrm{R}^\mu_1 + t^\kappa \widetilde{\chi}_1(t^\kappa D_x)\chi(t^{-\sigma}D_x)[|D_x|,x]\mathrm{R}^\mu_1 \\
&= t^\kappa\widetilde{\chi}_1(t^\kappa D_x)\chi(t^{-\sigma}D_x)\mathrm{R}^\mu_1 \left[x_n|D_x| -tD_n + \frac{1}{2i}\frac{D_n}{|D_x|}\right] \\
&+ t^\kappa\widetilde{\chi}_1(t^\kappa D_x)\chi(t^{-\sigma}D_x)\mathrm{R}^\mu_1\left[tD_n - \frac{1}{2i}\frac{D_n}{|D_x|}\right]+\delta_{\mu=1} it^\kappa\widetilde{\chi}_1(t^\kappa D_x)\chi(t^{-\sigma}D_x)Op(|\xi|\partial_n(\xi_1|\xi|^{-1}))\\
&-it^\kappa \widetilde{\chi}_1(t^\kappa D_x)\chi(t^{-\sigma}D_x)\mathrm{R}_n\mathrm{R}^\mu_1
\end{align*}
with $\widetilde{\chi}(\xi):=(1-\chi)(\xi)|\xi|^{-1}$ and using relation \eqref{relation_w_wave_Zjw-new} with $w=u_\pm$.
The proof for $|\alpha|=2$ is analogous. 
It is based on the commutation of $x^\alpha$ with $(1-\chi)(t^\kappa D_x)\chi(t^{-\sigma}D_x)$ (the commutator is here a $L^2-L^\infty$ bounded operator with norm $O(t^{\kappa})$), on the fact that we can rewrite $(1-\chi)(t^\kappa D_x)\chi(t^{-\sigma}D_x) x^\alpha \mathrm{R}^\mu_1$ making appear $(x|D_x| -tD_x + \frac{1}{2i}\frac{D_x}{|D_x|})^\alpha$ by considering $\widetilde{\chi}_2(\xi) := (1-\chi)(\xi)|\xi|^{-2}$ instead of previous $\widetilde{\chi}_1$, and use relation \eqref{relation_w_wave_Zjw-new}.
Doing so we derive the following inequality\small 
\begin{multline*}
\left\| (1-\chi)(t^\kappa D_x)\chi(t^{-\sigma}D_x)\big[x^\alpha \mathrm{R}^\mu_1u_\pm\big](t,\cdot)\right\|_{L^\infty} \\
\lesssim t^{2\kappa} \left[\sum_{|\mu|=2}\|Z^\mu u_\pm(t,\cdot)\|_{L^2} + \sum_{|\mu|\le 1} t^{2-|\mu|}\|Z^\mu u_\pm(t,\cdot)\|_{H^1} + \sum_{|\mu|=1}^2\|x^\mu\textit{NL}_w(t,\cdot)\|_{L^2}\right] 
\le CB\varepsilon t^{2+2\kappa +\frac{\delta}{2}},
\end{multline*}\normalsize
last estimate following from a-priori estimates, \eqref{est_L2_xNLw} and \eqref{est_L2:xixj_NL}.
Summing up \eqref{est:A}, \eqref{est:B}, \eqref{ineq:x_alpha_R1u}, together with the previous choice of $\kappa$ and the fact that in \eqref{cut_1} $N(s)\ge 6$ if $s>0$ is sufficiently large, we deduce that
\begin{equation} \label{est_3_vNFGamma-vJ}
\left\|[(1-\chi)(t^{-\sigma} D_x)\chi_1(t^{-\sigma}D_x)(\Gamma^I v)_\pm(t,\cdot)] (1-\chi)(t^\kappa D_x)\chi(t^{-\sigma}D_x)\mathrm{R}^\mu_1u_\pm(t,\cdot)\right\|_{L^\infty}\le CB\varepsilon t^{-\frac{3}{2}}.
\end{equation}
Therefore, from \eqref{est_1_vNFGamma-vJ}, \eqref{cut_2},\eqref{est_3_vNFGamma-vJ}, and the uniform continuity on $L^\infty$ of $\chi(t^{-\sigma} D_x)$, we find that
\begin{equation*}
\left\| \chi(t^{-\sigma}D_x)\left(\vNFGamma-v^I_{-} \right)(t,\cdot)\right\|_{L^\infty} \lesssim \sum_{\mu=0}^1 t^\sigma \|\chi(t^{-\sigma}D_x)(\Gamma^I v)_\pm (t,\cdot)\|_{L^\infty}\|\mathrm{R}^\mu_1u_\pm(t,\cdot)\|_{L^\infty} 
+ CB\varepsilon t^{-1},
\end{equation*}
and as $\sigma$ is small and $\varepsilon_0<(2A)^{-1}$, from \eqref{est: bootstrap upm} we obtain \eqref{ineq:Linfty_ vINF - vI}.

\smallskip

In order to prove \eqref{ineq:L2 Zm vINF-vI} we apply $Z_m$ to equality \eqref{explicit vNfGamma-vJ-} and apply the Leibniz rule. As
\begin{equation} \label{commutator_Z_DtD1}
\begin{split}
[Z_m, D_t]= - D_m, \quad [Z_m, D_1]=-\delta_{m1} D_t, \quad [Z_m, \langle D_x\rangle]= - D_m\langle D_x\rangle^{-1} D_t,
\end{split}
\end{equation}
with $\delta_{m1}$ the Kronecker delta,
we find that \small
\begin{equation}\label{Zm(vgamma-vJ)}
\begin{split}
&2i\chi(t^{-\sigma}D_x) Z_m (\vNFGamma - v^I_{-}) \\
&=\chi(t^{-\sigma}D_x) \Big[ (D_t Z_m\Gamma^I v)(D_1u) - (D_1Z_m\Gamma^I v)(D_tu)+ D_1[(Z_m\Gamma^I v)(D_tu)] - \langle D_x\rangle [(Z_m \Gamma^I v)(D_1u)]\\
&\hspace{15pt}+ (D_t \Gamma^I v)(D_1Z_m u) - (D_1\Gamma^I v)(D_t Z_mu)+ D_1[(\Gamma^I v)(D_t Z_mu)] - \langle D_x\rangle [ (\Gamma^I v)(D_1Z_m u)]\\
&\hspace{15pt}- (D_m\Gamma^I v)(D_1u) + \delta_{m1}(D_t\Gamma^I v)(D_tu) - \delta_{m1} D_t[(\Gamma^I v)(D_t u)] + \frac{D_m}{\langle D_x\rangle} D_t[(\Gamma^I v)(D_1 u)] \\
& \hspace{15pt}-\delta_{m1} (D_t\Gamma^I v)(D_tu) + (D_1\Gamma^I v)(D_m u)-\delta_{m1} D_1[(\Gamma^I v)(D_t u)] + \delta_{m1}\langle D_x\rangle [(\Gamma^I v)(D_tu)]\Big].
\end{split}
\end{equation}\normalsize
The $L^2$ norm of all products in the above second, fourth and fifth line, i.e. those in which $Z_m$ is not acting on the wave component $u$,
is estimated by
\begin{multline}\label{first_terms_Zm(vgamma-vj)}
\sum_{\mu=0}^1 t^\sigma \left(\|(Z_m\Gamma^I v)_\pm (t,\cdot)\|_{L^2}+ \|(\Gamma^I v)_\pm (t,\cdot)\|_{L^2}\right)\left(\|\mathrm{R}^\mu_1 u_\pm (t,\cdot)\|_{L^\infty}+\|D_tu_\pm (t,\cdot)\|_{L^\infty}\right)\\
\le CAB\varepsilon^2 t^{-\frac{1}{2} + \frac{\delta_0}{2}+\sigma},
\end{multline}
after inequality \eqref{Hsinfty_norm_DtU} with $s=0$ and a-priori estimates.
The $L^2$ norm of products appearing in the second line are, instead, estimated by using \eqref{def uIpm vIpm} and \eqref{cor_estZm_2} with $L=L^2$, $\Gamma w= Z_mu$, $s>0$ sufficiently large so that $N(s)\ge 2$.
It is hence bounded by
\begin{equation*}
\begin{split}
& t^\sigma \left\|\chi(t^{-\sigma}D_x)(\Gamma^I v)_\pm(t,\cdot)\right\|_{L^\infty}\| (Z_mu)_\pm(t,\cdot)\|_{L^2}\\ 
&+ t^{-2}\Big(\sum_{|\mu|=0}^1 \|x^\mu (\Gamma^Iv)_\pm(t,\cdot)\|_{L^2} + t\|(\Gamma^I v)_\pm(t,\cdot)\|_{L^2}\Big)\left(\|u_\pm(t,\cdot)\|_{H^s}+ \|D_tu_\pm(t,\cdot)\|_{H^s}\right) \\
& \le CB^2\varepsilon^2 t^{2\sigma+ \frac{\delta_{3-j}+\delta_2}{2}},
\end{split}
\end{equation*}
where the latter estimate is obtained using the fact that $\chi(t^{-\sigma}D_x)$ is a bounded operator from $L^2$ to $L^\infty$ with norm $O(t^\sigma)$, together with \eqref{Hs_norm_DtU}, \eqref{norm_L2_xj-GammaIv-} and a-priori estimates.
That concludes, together with \eqref{first_terms_Zm(vgamma-vj)}, the proof of \eqref{ineq:L2 Zm vINF-vI} and of the statement.
\endproof
\end{lem}

\begin{lem}\label{Lem_appendix: preliminary est VJ}
There exists a constant $C>0$ such that, for any $\rho\in\mathbb{N}$, $\chi\in C^\infty_0(\mathbb{R}^2)$, equal to 1 in a neighbourhood of the origin, $\sigma>0$ small, and every $t\in [1,T]$,
\begin{equation} \label{Linfty_est_VJ}
\sum_{|I|=1}\|\chi(t^{-\sigma}D_x)V^I(t,\cdot)\|_{H^{\rho,\infty}} \le CB\varepsilon t^{-1+\beta+\frac{\delta_1}{2}},
\end{equation}
with $\beta>0$ small such that $\beta\rightarrow 0$ as $\sigma\rightarrow 0$.
\proof
Since $\chi(t^{-\sigma}D_x)$ is a bounded operator from $L^\infty$ to $H^{\rho,\infty}$ with norm $O(t^{\sigma\rho})$, for any $\rho\in\mathbb{N}$, it is enough to prove that the $L^\infty$ norm of $\chi(t^{-\sigma}D_x) V^I(t,\cdot)$ is bounded by the right hand side of \eqref{Linfty_est_VJ}. Moreover, as this latter inequality is automatically satisfied when $\Gamma$ is a spatial derivative after a-priori estimate \eqref{est: boostrap vpm} and the fact that operator $\chi(t^{-\sigma}D_x)$ is uniformly bounded on $L^\infty$, for the rest of the proof we will assume that $\Gamma\in \{\Omega, Z_j, j=1,2\}$ is a Klainerman vector field.
We also warn the reader that, throughout the proof, $C$ and $\beta$ will denote some positive constants that may change line after line, with $\beta\rightarrow 0$ as $\sigma\rightarrow 0$.

Instead of proving the result of the statement directly on $\chi(t^{-\sigma}D_x)v^I_\pm$ we do it for $\chi(t^{-\sigma}D_x)\vNFGamma$, where $\vNFGamma$ has been introduced in \eqref{def_vNF-Gamma} and is considered here for $|I|=1$ and $\Gamma^I=\Gamma$. In fact, by \eqref{ineq:Linfty_ vINF - vI}
\begin{equation}\label{vI-_bounded_by_vNFGamma}
\left\| \chi(t^{-\sigma}D_x)v^I_{-}(t,\cdot)\right\|_{L^\infty}\le 2 \left\| \chi(t^{-\sigma}D_x)\vNFGamma(t,\cdot)\right\|_{L^\infty}+CB\varepsilon t^{-1}.
\end{equation}
The advantage of dealing with this new function is related to the fact that it is solution to a half Klein-Gordon equation with a more suitable non-linearity (see \eqref{KG_vNF-Gamma}) than the equation satisfied by $v^I_{-}$.
In fact, it is a computation to show that from definition \eqref{def_vNF-Gamma}
\begin{equation}\label{KG_vNF-Gamma}
[D_t+\langle D_x\rangle]\vNFGamma(t,x) = \NLNF,
\end{equation}
where
\begin{equation} \label{def_NLNF}
\NLNF= r^{I,\textit{NF}}_{kg}(t,x) + Q^\mathrm{kg}_0(v_\pm, D_1 u^I_\pm) + G^\mathrm{kg}_1(v_\pm, Du_\pm),
\end{equation}
$G^\mathrm{kg}_1(v_\pm, Du_\pm) = G_1(v,\partial u)$ with $G_1$ given by \eqref{def_G1}, and
\begin{multline}\label{def:rNF-Gamma-kg}
r^{I,\textit{NF}}_{kg}(t,x)
 = -\frac{i}{4(2\pi)^2}\\
 \times \sum_{j_1,j_2\in \{+,-\}}\int e^{ix\cdot\xi} B^1_{(j_1,j_2,+)}(\xi,\eta) \left[\reallywidehat{\Gamma^I\textit{NL}_{kg}}(\xi-\eta)\hat{u}_{j_2}(\eta) - \widehat{v^I_{j_1}}(\xi-\eta) \widehat{\textit{NL}_w}(\eta)\right] d\xi d\eta,
\end{multline}
with $B^1_{(j_1,j_2,+)}$ given by \eqref{def of B(i1,i2,i3)} when $j_3=+$ and $k=1$. After \eqref{def_vNF-Gamma} and \eqref{explicit integral B} it appears that $r^{I,\textit{NF}}_{kg}$ has the following nice explicit expression
\begin{equation} \label{explicit rNFkg-Gamma}
r^{I,\textit{NF}}_{kg} =-\frac{i}{2}\left[(\Gamma^I\textit{NL}_{kg}) D_1u - (D_1\Gamma^I v)\textit{NL}_w + D_1[(\Gamma^I v)\textit{NL}_w]\right].
\end{equation}
Using lemma \ref{Lem_appendix:Linfty_bound_chi_w} and relation \eqref{relation_Zjw_Ljwidetilde(w)} with $w=\vNFGamma$, and reminding that $\|t w(t,t\cdot)\|_{L^2}=\|w(t,\cdot)\|_{L^2}$, we find the following 
\begin{multline}\label{prelimary_ineq_vNFGamma}
\left\| \chi(t^{-\sigma}D_x)\vNFGamma(t,\cdot)\right\|_{L^\infty} \\
\lesssim t^{-1+\beta} \sum_{|\mu|=0}^1\left\|\chi(t^{-\sigma}D_x)Z^\mu \vNFGamma(t,\cdot)\right\|_{L^2}+ \sum_{j=1}^2 t^{-1+\beta}\left\|\chi(t^{-\sigma}D_x) \big[x_j\NLNF\big](t,\cdot)\right\|_{L^2}.
\end{multline}
From equality \eqref{explicit vNfGamma-vJ-}, along with \eqref{def u+- v+-}, \eqref{def uIpm vIpm}, and a-priori estimates \eqref{est: bootstrap upm}, \eqref{est: bootstrap E02}, we immediately see that
\begin{equation} \label{est_preliminary_(vNFGamma - vJ-)}
\begin{split}
\left\| \chi(t^{-\sigma}D_x)(\vNFGamma -v^I_{-})(t,\cdot)\right\|_{L^2}&\lesssim t^\sigma\|v^I_\pm(t,\cdot)\|_{L^2}\left(\|u_\pm(t,\cdot)\|_{L^\infty}+\|\mathrm{R}_1u_\pm(t,\cdot)\|_{L^\infty}\right)\\
& \le CAB\varepsilon^2 t^{-\frac{1}{2}+\frac{\delta_2}{2}+\sigma},
\end{split}
\end{equation}
and as $\sigma, \delta_2\ll 1$ are small
\begin{equation}\label{est:vNFGamma}
\begin{split}
\left\|\chi(t^{-\sigma}D_x) \vNFGamma(t,\cdot)\right\|_{L^2} &\le \left\|\chi(t^{-\sigma}D_x) v^I_{-}(t,\cdot)\right\|_{L^2}+ \left\| \chi(t^{-\sigma}D_x)(\vNFGamma -v^I_{-})(t,\cdot)\right\|_{L^2} \\
&\le CB\varepsilon t^{\frac{\delta_2}{2}}.
\end{split} 
\end{equation}
Moreover, from \eqref{ineq:L2 Zm vINF-vI} and a-priori estimate \eqref{est: bootstrap E02} we have that, for every $m=1,2$, $t\in [1,T]$,
\begin{equation}\label{est_L2_ZmvNFGamma}
\left\|\chi(t^{-\sigma}D_x)Z_m \vNFGamma(t,\cdot) \right\|_{L^2}\le CB\varepsilon t^\frac{\delta_1}{2}.
\end{equation}
Finally, from \eqref{explicit rNFkg-Gamma}, \eqref{def u+- v+-}, \eqref{def uIpm vIpm}, \eqref{norm_Linfty_xv-}, \eqref{xjGamma_Nlkg} and a-priori estimates, we derive that
\begin{equation}\label{est:L2_xj_rINF_kg}
\begin{split}
&\left\|\chi(t^{-\sigma}D_x) \left[x_j r^{I,\textit{NF}}_{kg}\right](t,\cdot)\right\|_{L^2} \lesssim \|x_j\Gamma^I\textit{NL}_{kg}(t,\cdot)\|_{L^2}\left(\|u_\pm(t,\cdot)\|_{L^\infty}+ \|\mathrm{R}_1u_\pm(t,\cdot)\|_{L^\infty}\right)\\
& + \sum_{\mu=0}^1 t^\sigma\left(\|x_j^\mu v_\pm(t,\cdot)\|_{L^\infty}+\left\| x_j^\mu \frac{D_x}{\langle D_x\rangle}v_\pm(t,\cdot)\right\|_{L^\infty}\right)\|v^I_\pm(t,\cdot)\|_{L^2}\|v_\pm(t,\cdot)\|_{H^{2,\infty}}\\
& \le C(A+B)B\varepsilon^2 t^\frac{\delta_2}{2},
\end{split}
\end{equation}
while from \eqref{Hs_norm_DtU} with $s=0$, \eqref{norm_Linfty_xv-} and a-priori estimates
\begin{equation*}
\begin{split}
&\left\| \chi(t^{-\sigma}D_x)\left[x_j Q^\mathrm{kg}_0\left(v_\pm, D_1 u^I_\pm\right)\right](t,\cdot)\right\|_{L^2}+ \left\| \chi(t^{-\sigma}D_x)\left[x_j G^\mathrm{kg}_1\left(v_\pm, D u_\pm\right)\right](t,\cdot)\right\|_{L^2}\\
&\lesssim \left(\|x_j v_\pm(t,\cdot)\|_{L^\infty}+\left\| x_j \frac{D_x}{\langle D_x\rangle}v_\pm(t,\cdot)\right\|_{L^\infty}\right)\left(\|u^I_\pm(t,\cdot)\|_{H^1}+ \|D_tu_\pm(t,\cdot)\|_{L^2}\right)\\
& \le C(A+B)B\varepsilon t^{\delta_2}.
\end{split}
\end{equation*}
Therefore, from \eqref{def_NLNF} we deduce that
\begin{equation}\label{est_xj_NL-kg-NF-Gamma}
\|\chi(t^{-\sigma}D_x)\big[ x_j\NLNF\big](t,\cdot)\|_{L^2} \le C(A+B)B\varepsilon^2 t^{\delta_2},
\end{equation}
so injecting \eqref{est:vNFGamma}, \eqref{est_L2_ZmvNFGamma}, \eqref{est_xj_NL-kg-NF-Gamma} into \eqref{prelimary_ineq_vNFGamma}, and summing it up with \eqref{vI-_bounded_by_vNFGamma}, we obtain the result of the statement.
\endproof
\end{lem}

As done for the Klein-Gordon component in the above lemma, we also derive an estimate for the uniform norm of the wave component when a Klainerman vector field acts on it and its frequencies less or equal than $t^\sigma$ (see lemma \ref{Lem_appendix: est UJ}). We first need the following result.

\begin{lem}\label{Lem_appendix: L^2 estimates uJ}
Let $\Gamma\in \mathcal{Z}$, index $J$ be such that $\Gamma^J=\Gamma$, and $\widetilde{u}^J (t,x):= t(\Gamma u)_{-}(t,tx)$.
There exists a constant $C>0$ such that, for any $\theta_0, \chi\in C^\infty_0(\mathbb{R}^2)$, $\sigma>0$, and every $t\in [1,T]$,
\begin{subequations}
\begin{gather}
\|\widetilde{u}^J(t,\cdot)\|_{L^2}\le CB\varepsilon t^\frac{\delta_2}{2}, \label{utildeJ_L2}\\
\|\Omega_h \widetilde{u}^J(t,\cdot)\|_{L^2}\le CB\varepsilon t^{\frac{\delta_1}{2}}\label{Omega_utildeJ_L2} ,\\
\left\|\mathcal{M} \widetilde{u}^J(t,\cdot)\right\|_{L^2} \le CB\varepsilon t^{\frac{\delta_1}{2}},  \label{MuJ}\\
\left\| \theta_0(x) \oph(\chi(h^\sigma\xi))\Omega_h\mathcal{M}\widetilde{u}^J(t,\cdot)\right\|_{L^2}\le CB\varepsilon t^\frac{\delta_0}{2}.
\end{gather}
\end{subequations}
\proof
We warn the reader that, throughout the proof, $C$ will denote a positive constant that may change line after line. We also recall that
\begin{equation*}
[D_t + \langle D_x\rangle] (\Gamma u)_{-}(t,x)=\Gamma \textit{NL}_w(t,x).
\end{equation*}
Estimates \eqref{utildeJ_L2} and \eqref{Omega_utildeJ_L2} are straightforward after \eqref{est: bootstrap E02} and the fact that 
\[\|\widetilde{u}^J(t,\cdot)\|_{L^2}=\|(\Gamma u)_{-}(t,\cdot)\|_{L^2}, \quad \|\Omega_h\widetilde{u}^J(t,\cdot)\|_{L^2}=\|(\Omega\Gamma u)_{-}(t,\cdot)\|_{L^2}.\]
From \eqref{relation_Zjw_Mjwidetilde(w)} with $w=(\Gamma u)_{-}$ and $f=\Gamma \textit{NL}_w$, estimates \eqref{est: bootstrap E02}, \eqref{xjGamma_NLw}, along with the fact that $\delta_2\ll \delta_1$ (e.g. $2\delta_2\le\delta_1$), and $(A+B)\varepsilon_0<1$, we obtain \eqref{MuJ}.
By \eqref{relation_Zjw_Mjwidetilde(w)} we also derive that, for any $n=1,2,$
\begin{equation} \label{preliminary_Omega-M-uJ}
\begin{split}
\left\| \theta_0(x) \oph(\chi(h^\sigma\xi)) \Omega_h\mathcal{M}_n\widetilde{u}^J(t,\cdot)\right\|_{L^2} &\lesssim \|\Omega Z_n (\Gamma u)_{-}(t,\cdot)\|_{L^2} + \sum_{\mu=0}^1\|\Omega^\mu (\Gamma u)_{-}(t,\cdot)\|_{L^2} \\
& + \left\|\theta_0\Big(\frac{x}{t}\Big) \chi(t^{-\sigma}D_x)\Omega [x_n \Gamma \textit{NL}_w](t,\cdot) \right\|_{L^2}.
\end{split}
\end{equation}
The first two norms in the above right hand side are controlled by $E^0_3(t;W)^{1/2}$ and are hence bounded by $CB\varepsilon t^\frac{\delta_0}{2}$.
By commutating $x_n$ with $\chi(t^{-\sigma}D_x)\Omega$, and using that $\theta_0\big(\frac{x}{t}\big)x_n=t \theta^n_0\big(\frac{x}{t}\big)$, with $\theta^n_0(z):=\theta_0(z)z_n$, we deduce that
\begin{equation*}
\begin{split}
\left\|\theta_0\Big(\frac{x}{t}\Big)\chi(t^{-\sigma}D_x) \Omega [x_n \Gamma\textit{NL}_w](t,\cdot) \right\|_{L^2}\lesssim t \sum_{\mu=0}^1 \|\chi_1(t^{-\sigma}D_x)\Omega^\mu \Gamma  \textit{NL}_w\|_{L^2},
\end{split}
\end{equation*}
for some new $\chi_1\in C^\infty_0(\mathbb{R}^2)$.
On the one hand, using \eqref{Gamma_Nlw}, \eqref{Hs norm DtV} with $s=0$ and a-priori estimates we derive that
\begin{equation} \label{est_t_GammaNLw}
\begin{split}
t\|\Gamma \textit{NL}_w\|_{L^2} 
&\lesssim t \|v_\pm(t,\cdot)\|_{H^{2,\infty}}\left(\|(\Gamma v)_\pm(t,\cdot)\|_{H^1}+ \|v_\pm(t,\cdot)\|_{H^1}+\|D_tv_\pm(t,\cdot)\|_{L^2}\right)\lesssim CB\varepsilon t^\frac{\delta_2}{2}.
\end{split}
\end{equation}
On the other hand, when we compute $\Omega \Gamma \textit{NL}_w$ we find among the out-coming quadratic terms the following ones
\begin{equation*}
Q^\mathrm{w}_0((\Omega v)_\pm, D_1(\Gamma v)_\pm) \quad \text{and} \quad Q^\mathrm{w}_0((\Gamma v)_\pm, D_1( \Omega v)_\pm),
\end{equation*} 
which we estimate in the $L^2$ norm (when truncated for frequencies less or equal than $t^\sigma$) by means of \eqref{cor_estOmega_2} with $L=L^2$, $\Gamma w= \Omega v$, and $s>0$ large enough to have $N(s)\ge 3$. From \eqref{norm_L2_xj-GammaIv-}, \eqref{Linfty_est_VJ} and a-priori estimates, we obtain that
\begin{equation*}
\begin{split}
&\left\| \chi(t^{-\sigma}D_x)Q^\mathrm{w}_0((\Omega v)_\pm, D_1(\Gamma v)_\pm)\right\|_{L^2}+ \left\|\chi(t^{-\sigma}D_x)Q^\mathrm{w}_0((\Gamma v)_\pm, D_1( \Omega v)_\pm)\right\|_{L^2}\\
& \lesssim t^\sigma \|\chi(t^{-\sigma}D_x)(\Omega v)_\pm (t,\cdot)\|_{L^\infty}\|(\Gamma v)_\pm(t,\cdot)\|_{H^1} +\sum_{|\mu|=0}^1  t^{-3}\|v_\pm(t,\cdot)\|_{H^s}\|x^\mu (\Gamma v)_\pm (t,\cdot)\|_{H^1} \\
& \le CB^2\varepsilon^2 t^{-1+\beta+\frac{\delta_1+\delta_2}{2}},
\end{split}
\end{equation*}
with $\beta>0$ small such that $\beta\rightarrow 0$ as $\sigma\rightarrow 0$.
All remaining quadratic contributions to $\Omega \Gamma \textit{NL}_w$ are estimated with
\begin{multline*}
\|(\Omega \Gamma v)_\pm (t,\cdot)\|_{H^1}\|v_\pm(t,\cdot)\|_{H^{2,\infty}} + \|(\Omega v)_\pm(t,\cdot)\|_{L^2}\left(\|v_\pm(t,\cdot)\|_{H^{1,\infty}}+\|D_tv_\pm(t,\cdot)\|_{L^\infty}\right) \\
+ \|v_\pm(t,\cdot)\|_{H^{1,\infty}}\left(\|(\Omega v)_\pm(t,\cdot)\|_{H^1} + \|D_t(\Omega v)_\pm (t,\cdot)\|_{L^2}\right),
\end{multline*}
and are hence bounded by $C(A+B)B\varepsilon^2 t^{-1+\frac{\delta_1}{2}}$ after \eqref{est: Hsinfty Dt V}, \eqref{L2_norm_DtOmegaV} and the a-priori estimates.
This finally implies that
\begin{equation*}
t\left\|\chi(t^{-\sigma}D_x)\Omega \Gamma \textit{NL}_w(t,\cdot) \right\|_{L^2} \le C(A+B)B\varepsilon^2 t^{\beta+\frac{\delta_1+\delta_2}{2}},
\end{equation*}
which, together with \eqref{est_t_GammaNLw} and the fact that $\beta+\frac{\delta_1+\delta_2}{2}\le \frac{\delta_0}{2}$, as $\delta_2\ll \delta_1 \ll \delta_0$ and $\beta>0$ is as small as we want provided that $\sigma$ is small, gives
\begin{equation*}
\left\|\theta_0\Big(\frac{x}{t}\Big)\chi(t^{-\sigma}D_x) \Omega [x_n \Gamma\textit{NL}_w](t,\cdot) \right\|_{L^2} \le CB\varepsilon t^\frac{\delta_0}{2}.
\end{equation*}
\endproof
\end{lem}

\begin{lem}\label{Lem_appendix: est UJ}
There exists a constant $C>0$ such that, for any $\rho\in\mathbb{N}$, $\chi\in C^\infty_0(\mathbb{R}^2)$ equal to 1 in a neighbourhood of the origin, $\sigma>0$ small, and every $t\in [1,T]$,
\begin{equation} \label{Linfty_est_UJ}
\sum_{|J|=1}\sum_{|\mu|=0}^1\|\chi(t^{-\sigma}D_x)\mathrm{R}^\mu U^J(t,\cdot)\|_{H^{\rho,\infty}} \le C(A+B)\varepsilon t^{-\frac{1}{2}+\beta+\frac{\delta_1}{2}},
\end{equation}
for a small $\beta>0$, $\beta\rightarrow 0$ as $\sigma\rightarrow 0$.
\proof
We warn the reader that, throughout the proof, $C$ and $\beta$ will denote two positive constants that may change line after line, with $\beta\rightarrow 0$ as $\sigma\rightarrow 0$.
Moreover, since $\chi(t^{-\sigma}D_x)$ is a bounded operator from $L^\infty$ to $H^{\rho,\infty}$ with norm $O(t^{\sigma\rho})$, for any $\rho\in\mathbb{N}$, we can reduce to prove that the $L^\infty$ norm of $\chi(t^{-\sigma}D_x)\mathrm{R}^\mu U^J(t,\cdot)$ is bounded by the right hand side of \eqref{Linfty_est_UJ}.
We observe that this estimate is automatically satisfied when $J$ is such that $\Gamma^J$ is a spatial derivative, as a consequence of a-priori estimate \eqref{est: bootstrap upm}.
We therefore assume that $\Gamma^J$ is one of the Klainerman vector fields $\Omega, Z_m$, for $m\in\{1,2\}$.

Introducing $\widetilde{u}^J(t,x):=t u^J_{-}(t,tx)$, passing to the semiclassical setting ($t\mapsto t$, $x\mapsto \frac{x}{t}$, and $h:=1/t$), and reminding that $u^J_+ = -\overline{u^J_{-}}$, inequality \eqref{Linfty_est_UJ} becomes
\begin{equation}\label{est_Linfty_uJ_semiclassical}
\sum_{|\mu|=0}^1\left\|\oph\Big(\chi(h^\sigma\xi)(\xi|\xi|^{-1})^\mu\Big)\widetilde{u}^J_{-}(t,\cdot) \right\|_{L^\infty}\le C(A+B)\varepsilon h^{-\frac{1}{2}-\beta-\frac{\delta_1}{2}}.
\end{equation}
We consider a Littlewood-Paley decomposition such that
\begin{equation}\label{dec_UJ}
\chi(h^\sigma\xi)= \widetilde{\chi}(h^{-1}\xi)+\sum_k(1-\widetilde{\chi})(h^{-1}\xi)\psi(2^{-k}\xi)\chi(h^\sigma\xi),
\end{equation}
for some suitably supported $\widetilde{\chi}\in C^\infty_0(\mathbb{R}^2)$, $\psi\in C^\infty_0(\mathbb{R}^2\setminus \{0\})$, and immediately observe that the above sum is restricted to indices $k$ such that $h\lesssim 2^k\lesssim h^{-\sigma}$.
By the classical Sobolev injection, the uniform continuity of $\oph(\xi|\xi|^{-1})$ on $L^2$, and a-priori estimate \eqref{est: bootstrap E02}, we derive that for any $|\mu|\le 1$, every $t\in [1,T]$,
\begin{equation} \label{est_utildeJ_xi<h}
\begin{split}
\left\|\oph\big(\widetilde{\chi}(h^{-1}\xi)(\xi|\xi|^{-1})^\mu\big)\widetilde{u}^J(t,\cdot)\right\|_{L^\infty} &= \|\chi( D_x) \oph((\xi|\xi|^{-1})^\mu)\widetilde{u}^J(t,\cdot)\|_{L^\infty}\\
&\lesssim \| u^J_{-}(t,\cdot)\|_{L^2} \le CB\varepsilon t^{\frac{\delta_2}{2}}.
\end{split}
\end{equation}
If we concisely denote by $\phi_k(\xi)$ the $k$-th addend in decomposition \eqref{dec_UJ} and introduce two smooth cut-off functions $\chi_0$, $\gamma$, with $\chi_0$ radial and equal to 1 on the support of $\phi_k$, $\gamma$ with sufficiently small support, we can write
\begin{align*}
\oph\left(\phi_k(\xi)(\xi|\xi|^{-1})^\mu\right)\widetilde{u}^J &= \oph\Big(\gamma\Big(\frac{x|\xi|-\xi}{h^{1/2-\sigma}}\Big)\phi_k(\xi) (\xi|\xi|^{-1})^\mu\Big) \oph(\chi_0(h^\sigma\xi))\widetilde{u}^J \\
&+  \oph\Big((1-\gamma)\Big(\frac{x|\xi|-\xi}{h^{1/2-\sigma}}\Big)\phi_k(\xi) (\xi|\xi|^{-1})^\mu\Big)\oph(\chi_0(h^\sigma\xi))\widetilde{u}^J.
\end{align*}
On the one hand, after proposition \ref{Prop : continuity of Op(gamma1):X to L2}, the fact that $2^k\lesssim h^{-\sigma}$, a-priori estimate \eqref{est: bootstrap E02}, and the uniform $L^2$ continuity of $\oph(\chi_0(h^\sigma\xi))$, we have that for any $|\mu|\le 1$
\begin{multline}\label{gamma_uJ}
\left\| \oph\Big(\gamma\Big(\frac{x|\xi|-\xi}{h^{1/2-\sigma}}\Big)\phi_k(\xi) (\xi|\xi|^{-1})^\mu\Big)\oph(\chi_0(h^\sigma\xi))\widetilde{u}^J(t,\cdot)\right\|_{L^\infty} \\
\lesssim h^{-\frac{1}{2}-\beta}\left(\|\oph(\chi_0(h^\sigma\xi))\widetilde{u}^J(t,\cdot)\|_{L^2}+ \|\theta_0(x)\Omega_h\oph(\chi_0(h^\sigma\xi))\widetilde{u}^J(t,\cdot)\|_{L^2}\right)\\
\lesssim h^{-\frac{1}{2}-\beta}\left(\|u^J_{-}(t,\cdot)\|_{L^2}+ \|\Omega u^J_{-}(t,\cdot)\|_{L^2}\right) \le CB\varepsilon h^{-\frac{1}{2}-\beta-\frac{\delta_1}{2}}.
\end{multline}
On the other hand, using that $(1-\gamma)(z)=\gamma_1^j(z)z_j$, where $\gamma_1^j(z):=(1-\gamma)(z)z_j|z|^{-2}$ is such that $|\partial^\alpha_z\gamma_1^j(z)|\le \langle z\rangle^{-1-|\alpha|}$, we derive from \eqref{est: Linfty of gamma1 with L}, the commutation between $\mathcal{M}$ with $\oph(\chi_0(h^\sigma\xi))$, and lemma \ref{Lem_appendix: L^2 estimates uJ}, that
\begin{multline*}
\left\| \oph\Big((1-\gamma)\Big(\frac{x|\xi|-\xi}{h^{1/2-\sigma}}\Big)\phi_k(\xi) (\xi|\xi|^{-1})^\mu\Big)\oph(\chi_0(h^\sigma\xi))\widetilde{u}^J\right\|_{L^\infty}\\
\lesssim h^{-\beta}\sum_{\gamma,|\nu|=0}^1\|(\theta_0(x)\Omega_h)^\gamma \mathcal{M}^\nu \oph(\chi_0(h^\sigma\xi)) \widetilde{u}^J(t,\cdot)\|_{L^2}
\le CB\varepsilon t^{\beta+\frac{\delta_0}{2}}.
\end{multline*}
Combining this estimate with \eqref{gamma_uJ} we deduce that
\begin{equation*}
\|\oph\left(\phi_k(\xi)(\xi|\xi|^{-1})^\mu\right)\widetilde{u}^J(t,\cdot)\|_{L^\infty}\le C(A+B)\varepsilon h^{-\frac{1}{2}-\beta-\frac{\delta_1}{2}},
\end{equation*}
for any $|\mu|\le 1$, and hence \eqref{est_Linfty_uJ_semiclassical} 
after \eqref{dec_UJ}, \eqref{est_utildeJ_xi<h}, up to a further loss $|\log h|$, as a consequence of the fact that the sum in \eqref{dec_UJ} is finite and taken over indices $k$ such that $\log h \lesssim k\lesssim \log h^{-1}$.
\endproof
\end{lem}

\begin{lem} \label{Lem_appendix:xGammav_Linfty}
There exists a positive constant $C>0$ such that, for any $\chi\in C^\infty_0(\mathbb{R}^2)$ equal to 1 in a neighbourhood of the origin, $\sigma>0$, and every $t\in [1,T]$,
\begin{equation} \label{norm_Linfty_xjGammav-}
\sum_{|\mu|=0}^1\left\| \chi(t^{-\sigma}D_x)\left[x_j\Big(\frac{D_x}{\langle D_x\rangle}\Big)^\mu (\Gamma v)_\pm(t,\cdot)\right]\right\|_{L^\infty}\le CB\varepsilon t^{\beta+\frac{\delta_1}{2}},
\end{equation}
with $\beta>0$ small, $\beta\rightarrow 0$ as $\sigma\rightarrow 0$.
\proof
We warn the reader that, throughout the proof, $C$ and $\beta$ will denote two positive constants that may change line after line, with $\beta\rightarrow 0$ as $\sigma\rightarrow 0$. As $\Gamma v_+=-\overline{\Gamma v_{-}}$, it is enough to prove the statement for $\Gamma v_{-}$.

If $\Gamma$ is a spatial derivative, estimate \eqref{norm_Linfty_xjGammav-} is just consequence of the uniform continuity of $\chi(t^{-\sigma}D_x)$ on $L^\infty$ and of \eqref{norm_Linfty_xv-}. We then assume that $\Gamma\in \{\Omega, Z_m, m=1,2\}$ is a Klainerman vector field.
First of all, we observe that by \eqref{xjw_Zjw} with $w= (\Gamma v)_{-}$ and $f=\Gamma\textit{NL}_{kg}$, along with the classical Sobolev injection,
\begin{equation}\label{Linfty_xj_Gammav_preliminary}
\sum_{|\mu|=0}^1\left\| x_j\Big(\frac{D_x}{\langle D_x\rangle}\Big)^\mu (\Gamma v)_{-}(t,\cdot)\right\|_{L^\infty}\lesssim  \|Z_j(\Gamma v)_{-}(t,\cdot)\|_{H^1}+ t\|(\Gamma v)_{-}(t,\cdot)\|_{H^2}+ \sum_{\mu=0}^1\|x^\mu_j\Gamma\textit{NL}_{kg}(t,\cdot)\|_{L^\infty}.
\end{equation}
From equality \eqref{notation: NLGamma} and lemma \ref{Lem_app:products_Gamma} with $L=L^\infty$ and $s>0$ large enough so that $N(s)\ge3$, together with estimates \eqref{est: bootstrap argument a-priori est}, \eqref{Hs_norm_DtU}, \eqref{Hsinfty_norm_DtU}, \eqref{Hsinfty norm Dt R1U-new}, \eqref{Hs norm DtV}, \eqref{norm_xv-}, \eqref{norm_xu-}, \eqref{Linfty_est_VJ},  and \eqref{Linfty_est_UJ}, we get that \small
\begin{equation}\label{Gamma_NLkg_Linfty}
\begin{split}
&\|\Gamma\Nlkg(t,\cdot)\|_{L^\infty}\\
&\lesssim \sum_{\mu=0}^1\left(\|\chi(t^{-\sigma}D_x)(\Gamma v)_\pm(t,\cdot)\|_{H^{1,\infty}}\|\mathrm{R}^\mu_1u_\pm(t,\cdot)\|_{H^{2,\infty}}\right) + \|v_\pm(t,\cdot)\|_{H^{1,\infty}}\left\|\chi(t^{-\sigma}D_x)(\Gamma u)_\pm(t,\cdot)\right\|_{H^{2,\infty}} \\
&+ \|v_\pm(t,\cdot)\|_{H^{1,\infty}}\times \sum_{|\mu|=0}^1\left(\|\mathrm{R}^\mu u_\pm(t,\cdot)\|_{H^{2,\infty}}+ \|D_t\mathrm{R}^\mu u_\pm(t,\cdot)\|_{H^{1,\infty}}\right)\\
& + t^{-3}\left(\|v_\pm(t,\cdot)\|_{H^s}+\|D_tv_\pm(t,\cdot)\|_{H^s}\right)\Big(\sum_{|\mu|, |\nu|=0}^1\left\|x^\mu D_1\Big(\frac{D_x}{|D_x|}\Big)^\nu u_\pm(t,\cdot)\right\|_{L^2}+ t \|u_\pm(t,\cdot)\|_{L^2}\Big) \\
& + t^{-3}\Big(\sum_{|\mu|=0}^1 \|x^\mu v_\pm(t,\cdot)\|_{L^2}+t\|v_\pm(t,\cdot)\|_{L^2}\Big)\left(\|u_\pm(t,\cdot)\|_{H^s}+\|D_tu_\pm(t,\cdot)\|_{H^s}\right) \\
&\le CAB\varepsilon^2 t^{-\frac{3}{2}+\beta + \frac{\delta_1}{2}}.
\end{split}
\end{equation}\normalsize
Moreover, as
\begin{equation} \label{xj_QGammav_1}
\left\|x_j Q^\mathrm{kg}_0\big((\Gamma v)_\pm, D_1u_\pm\big)\right\|_{L^\infty} \lesssim \sum_{|\mu|, \nu=0}^1\left\|x_j \Big(\frac{D_x}{\langle D_x\rangle}\Big)^\mu(\Gamma v)_{-}(t,\cdot)\right\|_{L^\infty}\|\mathrm{R}^\nu_1u_\pm(t,\cdot)\|_{H^{2,\infty}}, 
\end{equation}
\begin{equation} \label{xjG1}
\left\| x_j G^\mathrm{kg}_1\big(v_\pm, D u_\pm\big)\right\|_{L^\infty}\lesssim \sum_{|\mu|=0}^1\left\|x_j \Big(\frac{D_x}{\langle D_x\rangle}\Big)^\mu v_\pm(t,\cdot)\right\|_{L^\infty} \left(\|u_\pm(t,\cdot)\|_{H^{2,\infty}}+\|D_tu_\pm(t,\cdot)\|_{H^{1,\infty}}\right),
\end{equation}
and by lemma \ref{Lem_app:products_Gamma} with $L=L^\infty$, $w= u$, and $s>0$ large enough so that $N(s)\ge 3$,
\begin{equation} \label{xj_QGammu}
\begin{split}
&\left\|x_j Q^\mathrm{kg}_0(v_\pm, D_1 (\Gamma u)_\pm)\right\|_{L^\infty}\lesssim \sum_{|\mu|=0}^1\left\|x_j \Big(\frac{D_x}{\langle D_x\rangle}\Big)^\mu v_\pm(t,\cdot)\right\|_{L^\infty} \|\chi(t^{-\sigma}D_x)(\Gamma u)_\pm(t,\cdot)\|_{H^{2,\infty}}\\
&  + t^{-3}\sum_{|\mu|,\nu=0}^1\left( \left\| x^\mu x_j^\nu v_\pm(t,\cdot)\right\|_{L^2} + t \|x^\nu_j v_\pm(t,\cdot)\|_{L^2}\right)\left(\|u_\pm(t,\cdot)\|_{H^s}+\|D_tu_\pm(t,\cdot)\|_{H^s}\right),
\end{split}
\end{equation}
we derive that
\begin{equation} \label{xjGamma_NL_preliminary}
\left\|x_j\Gamma\Nlkg(t,\cdot)\right\|_{L^\infty}\le CA\varepsilon t^{-\frac{1}{2}} \sum_{|\mu|, \nu=0}^1\left\|x_j \Big(\frac{D_x}{\langle D_x\rangle}\Big)^\mu(\Gamma v)_{-}(t,\cdot)\right\|_{L^\infty} + C(A+B)B\varepsilon^2 t^{-\frac{1}{2}+\beta+\frac{\delta_1+\delta_2}{2}},
\end{equation}
as follows using \eqref{Hs_norm_DtU}, \eqref{Hsinfty_norm_DtU} with $s=1$, \eqref{norm_xv-}, \eqref{norm_Linfty_xv-},\eqref{norm_L2_xixjv-}, \eqref{Linfty_est_UJ} and a-priori estimates.
By injecting the above inequality into \eqref{Linfty_xj_Gammav_preliminary} and using the fact that $\varepsilon_0<(2CA)^{-1}$, we initially obtain that
\begin{equation}\label{est_xj Gammav_preliminary}
\| x_j (\Gamma v)_{-}(t,\cdot)\|_{L^\infty}  + \left\|x_j\frac{D_x}{\langle D_x\rangle} (\Gamma v)_{-}(t,\cdot) \right\|_{L^\infty}\le CB\varepsilon t^{1+\frac{\delta_2}{2}}.
\end{equation}
If we take any smooth cut-off function $\chi$ and use equality \eqref{xjw_Zjw}, instead of \eqref{Linfty_xj_Gammav_preliminary} we find that
\begin{multline}\label{ineq_chi_xGammav-}
\sum_{|\mu|=0}^1\left\|\chi(t^{-\sigma} D_x)\Big[ x_j \Big(\frac{D_x}{\langle D_x\rangle}\Big)^\mu (\Gamma v)_{-}(t,\cdot)\Big]\right\|_{L^\infty}\lesssim  \|Z_j(\Gamma v)_{-}(t,\cdot)\|_{H^1}+ t\|\chi(t^{-\sigma}D_x)(\Gamma v)_{-}(t,\cdot)\|_{L^\infty}\\
+ \sum_{\mu=0}^1\left\|\chi(t^{-\sigma}D_x)\big[ x^\mu_j\Gamma\textit{NL}_{kg}(t,\cdot)\big]\right\|_{L^\infty},
\end{multline}
where now
\begin{equation*}
\left\|\chi(t^{-\sigma}D_x)\big[x_j\Gamma\textit{NL}_{kg}(t,\cdot)\big]\right\|_{L^\infty}\lesssim \left\|x_j\Gamma\textit{NL}_{kg}(t,\cdot)\right\|_{L^\infty} \le C(A+B)B\varepsilon^2 t^{\frac{1}{2}+\frac{\delta_2}{2}},
\end{equation*}
as follows injecting \eqref{est_xj Gammav_preliminary} into \eqref{xjGamma_NL_preliminary}. 
Therefore, from \eqref{Gamma_NLkg_Linfty}, lemma \ref{Lem_appendix: preliminary est VJ} and a-priori estimate \eqref{est: bootstrap E02} with $k=2$, we find that
\begin{equation} \label{xj_Gamma v-_preliminary2}
\left\|\chi(t^{-\sigma} D_x) \big[x_j (\Gamma v)_{-}(t,\cdot)\big]\right\|_{L^\infty} + \left\|\chi(t^{-\sigma} D_x) \left[x_j \frac{D_x}{\langle D_x\rangle} (\Gamma v)_{-}(t,\cdot)\right]\right\|_{L^\infty}  \le CB\varepsilon t^{\frac{1}{2}+\frac{\delta_2}{2}}.
\end{equation}
Finally, by means of lemma \ref{Lem_appendix:L_estimate of products} with $L=L^\infty$, $w_1 = x(\Gamma v)_\pm$, and $s>0$ such that $N(s)\ge 2$, we derive that for any $\chi\in C^\infty_0(\mathbb{R}^2)$ there is some $\chi_1\in C^\infty_0(\mathbb{R}^2)$ such that
\begin{multline*}
\left\|\chi(t^{-\sigma}D_x) x_jQ^\mathrm{kg}_0((\Gamma v)_\pm, D_1u_\pm) \right\|_{L^\infty} \\
\lesssim \sum_{|\mu|, \nu=0}^1\left\| \chi_1(t^{-\sigma}D_x)\left[x_j \Big(\frac{D_x}{\langle D_x\rangle}\Big)^\mu(\Gamma v)_{-}(t,\cdot)\right] \right\|_{L^\infty} \|\chi(t^{-\sigma}D_x)\mathrm{R}^\nu_1u_\pm(t,\cdot)\|_{H^{2,\infty}} \\
+ \sum_{\mu=0}^1 t^{-2}\left\|x^\mu_j (\Gamma v)_\pm(t,\cdot)\right\|_{L^2} \|u_\pm(t,\cdot)\|_{H^s}.
\end{multline*}
Then, combining such inequality with \eqref{xjG1}, \eqref{xj_QGammu}, together with \eqref{norm_L2_xj-GammaIv-}, \eqref{xj_Gamma v-_preliminary2}, and all the other inequalities to which we already referred before, from \eqref{notation: NLGamma} we find that
\begin{equation*}
\left\| \chi(t^{-\sigma}D_x) \big[x_j \Gamma\textit{NL}_{kg}(t,\cdot)\big]\right\|_{L^\infty}\le C(A+B)\varepsilon^2 t^{\delta_2},
\end{equation*}
which injected into \eqref{ineq_chi_xGammav-} finally implies, together with \eqref{est: bootstrap E02} with $k=2$, lemma \ref{Lem_appendix: preliminary est VJ}, and \eqref{Gamma_NLkg_Linfty}, the wished estimate \eqref{norm_Linfty_xjGammav-}.
\endproof
\end{lem}

Making use of lemmas \ref{Lem_appendix: preliminary est VJ} and \ref{Lem_appendix:xGammav_Linfty}, estimate \eqref{est:L2_xj_rINF_kg} can be improved of a factor $t^{-\frac{1}{2}}$. This improvement, that will be useful to derive  \eqref{first_estimate_xNLT}, is showed in the following lemma.

\begin{lem}
Let $I$ be a multi-index of length 1 and $r^{I,\textit{NF}}_{kg}$ be given by \eqref{def:rNF-Gamma-kg}.
There exists a constant $C>0$ such that, for any $\rho\in\mathbb{N}$, $\chi\in C^\infty_0(\mathbb{R}^2)$, equal to 1 in a neighbourhood of the origin, $\sigma>0$ small, $j=1,2$, and every $t\in [1,T]$,
\begin{equation} \label{enhanced_xjrINF_kg}
\left\|\chi(t^{-\sigma}D_x)\left[x_j r^{I,\textit{NF}}_{kg}\right](t,\cdot)\right\|_{L^2}\le C(A+B)AB\varepsilon^3 t^{-\frac{1}{2}+\beta+\frac{\delta+\delta_1}{2}},
\end{equation}
with $\beta>0$ small, $\beta\rightarrow 0$ as $\sigma\rightarrow 0$.
\proof
Let us remind the explicit expression \eqref{explicit rNFkg-Gamma} of $r^{I,\textit{NF}}_{kg}$ and consider the cubic term $x_j \Gamma^I\textit{NL}_{kg}(D_1u)$. Reminding \eqref{def u+- v+-} and applying lemma \ref{Lem_appendix:L_estimate of products} with $L=L^2$ and $s>0$ sufficiently large so that $N(s)\ge 2$, together with \eqref{xjGamma_Nlkg} and a-priori estimates, we derive that there is some $\chi_1\in C^\infty_0(\mathbb{R}^2)$ such that
\begin{equation}\label{est_xNLI_D1u_prel}
\begin{split}
& \left\|\chi(t^{-\sigma}D_x)\left[x_j \Gamma^I\textit{NL}_{kg}(D_1u)\right](t,\cdot)\right\|_{L^2}\\
&\lesssim \left\|\chi_1(t^{-\sigma}D_x)\left[x_j \Gamma^I\textit{NL}_{kg}\right](t,\cdot)\right\|_{L^2}\|\mathrm{R}_1u_\pm(t,\cdot)\|_{L^\infty} + t^{-2}\|x_j\Gamma^I\textit{NL}_{kg}(t,\cdot)\|_{L^2}\|u_\pm(t,\cdot)\|_{H^s}\\
& \le CA\varepsilon t^{-\frac{1}{2}}  \left\|\chi_1(t^{-\sigma}D_x)\left[x_j \textit{NL}^I_{kg}\right](t,\cdot)\right\|_{L^2} + C(A+B)B\varepsilon^2 t^{-1}.
\end{split}
\end{equation}
Then, recalling \eqref{notation: NLGamma} and using again lemma \ref{Lem_appendix:L_estimate of products} with $L=L^2$, $w_1 = (\Gamma v)_\pm$, and $s$ as before, in order to estimate the contribution coming from the first quadratic term in the right hand side of \eqref{notation: NLGamma}, we find that there is a new $\chi_2\in C^\infty_0(\mathbb{R}^2)$ such that
\begin{equation*}
\begin{split}
&  \left\|\chi_1(t^{-\sigma}D_x)\left[x_j \textit{NL}^I_{kg}\right](t,\cdot)\right\|_{L^2} \\
& \lesssim \sum_{|\mu|=0}^1\left\|\chi_2(t^{-\sigma}D_x)\Big[x_j \Big(\frac{D_x}{\langle D_x\rangle}\Big)^\mu (\Gamma v)_\pm\Big] (t,\cdot)\right\|_{L^\infty}\|u_\pm(t,\cdot)\|_{H^1}+ t^{-2}\|x_j(\Gamma v)_\pm(t,\cdot)\|_{L^2}\|u_\pm(t,\cdot)\|_{H^s} \\
&  +\sum_{|\mu|=0}^1 \left\| x_j \Big(\frac{D_x}{\langle D_x\rangle}\Big)^\mu v_\pm(t,\cdot)\right\|_{L^\infty}\left(\|(\Gamma u)_\pm(t,\cdot)\|_{H^1}+\|u_\pm(t,\cdot)\|_{H^1}+ \|D_t u_\pm(t,\cdot)\|_{L^2}\right)\\
&\le C(A+B)B\varepsilon^2 t^{\beta+\frac{\delta+\delta_1}{2}},
\end{split}
\end{equation*}
where the latter estimate is obtained from \eqref{Hs_norm_DtU} with $s=0$, \eqref{norm_Linfty_xv-}, \eqref{norm_L2_xj-GammaIv-} with $k=1$, \eqref{norm_Linfty_xjGammav-} and a-priori estimates.
This implies, combined with \eqref{est_xNLI_D1u_prel}, that
\begin{equation*}
 \left\|\chi(t^{-\sigma}D_x)\left[x_j \textit{NL}^I_{kg}(D_1u)\right](t,\cdot)\right\|_{L^2} \le C(A+B)AB\varepsilon^3 t^{-\frac{1}{2}+\beta+\frac{\delta+\delta_1}{2}},
\end{equation*}
and from \eqref{explicit rNFkg-Gamma}, \eqref{norm_Linfty_xv-} and a-priori estimates,
\begin{equation*}
\begin{split}
&\left\|\chi(t^{-\sigma}D_x) \left[x_j r^{I,\textit{NF}}_{kg}\right](t,\cdot)\right\|_{L^2} \lesssim  \left\|\chi(t^{-\sigma}D_x)\left[x_j \textit{NL}^I_{kg}(D_1u)\right](t,\cdot)\right\|_{L^2}\\
& + \sum_{\mu=0}^1 t^\sigma\left(\|x_j^\mu v_\pm(t,\cdot)\|_{L^\infty}+\left\| x_j^\mu \frac{D_x}{\langle D_x\rangle}v_\pm(t,\cdot)\right\|_{L^\infty}\right)\|v^I_\pm(t,\cdot)\|_{L^2}\|v_\pm(t,\cdot)\|_{H^{2,\infty}}\\
& \le  C(A+B)AB\varepsilon^3 t^{-\frac{1}{2}+\beta+\frac{\delta+\delta_1}{2}},
\end{split}
\end{equation*}
which concludes the proof of the statement.
\endproof
\end{lem}

\begin{lem}
Let $I$ be a multi-index of length 2.
There exists a constant $C>0$ such that, for every $j=1,2$, $t\in [1,T]$,
\begin{equation}\label{xj_GammaI_NLkg_I=2}
\left\| x_j \Gamma^I\textit{NL}_{kg}(t,\cdot)\right\|_{L^2}\le C(A+B)B\varepsilon^2 t^{\frac{1}{2}+\beta + \frac{\delta_1+\delta_2}{2}},
\end{equation}
with $\beta>0$ small, $\beta\rightarrow 0$ as $\sigma\rightarrow 0$.
\proof
We remind the reader about \eqref{GammaI_NLkg}. Instead of using \eqref{xj_Q(vI1,uI2)_appendix}, which was obtained by Sobolev injection, we apply lemma \ref{Lem_app:products_Gamma} with $L=L^2$, $\Gamma w =\Gamma^{I_2}u$, $s>0$ sufficiently large so that $N(s)\ge 3$, and exploit the fact that we have an estimate of the $H^{\rho,\infty}$ norm of $D_1u^{I_2}$ when truncated for frequencies less or equal than $t^\sigma$ (see lemma \ref{Lem_appendix: est UJ}). Therefore, for $(I_1,I_2)\in \mathcal{I}(I)$ such $|I_1|=|I_2|=1$ we obtain that
\begin{equation*}
\begin{split}
& \left\| x_j Q^\mathrm{kg}_0\left(v^{I_1}_\pm, D_1 u^{I_2}_\pm\right)(t,\cdot)\right\|_{L^2}\lesssim \sum_{\mu=0}^1\left\| x^\mu_j v^{I_1}_\pm(t,\cdot)\right\|_{L^2}\left\|\chi(t^{-\sigma}D_x) u^{I_2}_\pm (t,\cdot)\right\|_{H^{2,\infty}} \\
& + t^{-3}\left(\|u_\pm(t,\cdot)\|_{H^s}+\|D_tu_\pm(t,\cdot)\|_{H^s}\right) \left[\sum_{|\mu|=0}^2 \|x^\mu v^{I_1}_\pm(t,\cdot)\|_{L^2} + \sum_{|\mu|=0}^1 t\|x^\mu v^{I_1}_\pm(t,\cdot)\|_{L^2}\right]\\
& \le C(A+B)B\varepsilon^2 t^{\frac{1}{2}+\beta+\frac{\delta_1+\delta_2}{2}},
\end{split}
\end{equation*}
last estimate following from
 lemma \ref{Lem_appendix: est UJ} together with \eqref{Hs_norm_DtU}, \eqref{norm_L2_xj-GammaIv-} with $k=1$, \eqref{est:xixjGamma v-}, a-priori estimates, and the fact that $\delta_1,\delta_2\ll 1$ are small.
Consequently, from the following inequality
\begin{align*}
&\|x_j\Gamma^I \textit{NL}_{kg}(t,\cdot)\|_{L^2}\lesssim \sum_{\mu=0}^1\|\mathrm{R}_1^\mu u_\pm(t,\cdot)\|_{H^{2,\infty}}\sum_{\substack{|J|\le 2\\ \mu=0,1}}\| x^\mu_j(\Gamma^J v)_{-}(t,\cdot)\|_{L^2} \\
& + \sum_{|\mu|=0}^1\left\|x_j\Big(\frac{D_x}{\langle D_x\rangle}\Big)^\mu v_\pm(t,\cdot)\right\|_{L^\infty}\Big[\|u^I_\pm(t,\cdot)\|_{H^1} + \sum_{|J|<2}\big(\|u^J_\pm(t,\cdot)\|_{H^1}+\|D_tu^J_\pm(t,\cdot)\|_{L^2}\big) \Big]\\
& + \sum_{|I_1|=|I_2|=1}\left\|x_j Q^\mathrm{kg}_0\left(v^{I_1}_\pm, D_1 u^{I_2}_\pm\right)(t,\cdot)\right\|_{L^2},
\end{align*}
together with \eqref{norm_Linfty_xv-}, \eqref{Hs_norm_DtU} with $s=0$, \eqref{DtUI}, and \eqref{norm_L2_xj-GammaIv-} with $k=1$, we finally derive \eqref{xj_GammaI_NLkg_I=2}.
\endproof
\end{lem}

\begin{lem} \label{Lem: from energy to norms in sc coordinates-KG}
Let us fix $s\in\mathbb{N}$. There exists a constant $C>0$ such that, if we assume that a-priori estimates \eqref{est: bootstrap argument a-priori est} are satisfied in some interval $[1,T]$, for a fixed $T>1$, with $n\ge s+2$, then we have, for any $\chi \in C^\infty_0(\mathbb{R}^2)$ and $\sigma>0$ small,
\begin{subequations}
\begin{gather}
\|\widetilde{v}(t,\cdot)\|_{H^s_h} \le CB\varepsilon t^\frac{\delta}{2},\label{est:Hs_vtilde}\\
\sum_{|\mu|=1}\|\oph(\chi(h^\sigma\xi))\mathcal{L}^\mu\widetilde{v}(t,\cdot)\|_{L^2}  \le  CB\varepsilon t^\frac{\delta_2}{2}, \label{est:Lvtilde}
\end{gather}
\end{subequations}
for every $t\in [1,T]$.
\proof
We warn the reader that, throughout the proof, $C$ and $\beta$ will denote two positive constants that may change line after line, with $\beta>0$ is small as long as $\sigma$ is small.

It is straightforward to check that the $H^s_h$ norm of $\widetilde{v}$ is bounded by energy $E_n(t;W)^\frac{1}{2}$ whenever $n\ge s+2$, after definitions \eqref{def utilde vtilde}, \eqref{def vNF}, inequality \eqref{est_Hs_vnf-v}, and a-priori estimates \eqref{est: bootstrap upm}, \eqref{est: boostrap vpm}.

In order to prove \eqref{est:Lvtilde} we first use relation \eqref{relation between Zjv and Lj vtilde} and definition \eqref{def vNF} to derive that
\begin{equation} \label{norm_L2_Lnwidetilde(v)}
\begin{split}
\|\oph(\chi(h^\sigma\xi))\mathcal{L}_m\widetilde{v}(t,\cdot)\|_{L^2}&\lesssim \|Z_mV(t,\cdot)\|_{L^2}+\|\chi(t^{-\sigma}D_x)Z_m(v^{NF}-v_{-})(t,\cdot)\|_{L^2}\\
& + \|\widetilde{v}(t,\cdot)\|_{L^2}+ \|\chi(t^{-\sigma}D_x)[x_m r^{NF}_{kg}](t,\cdot)\|_{L^2},
\end{split}
\end{equation}
with $r^{NF}_{kg}$ given by \eqref{def rNF-kg}.
Using \eqref{def u+- v+-} we can rewrite \eqref{explicit_vNF-v_chapter5} and \eqref{explicit_rnfkg_chapter5} similarly to \eqref{explicit vNfGamma-vJ-}, \eqref{explicit rNFkg-Gamma}, as:
\begin{equation} \label{explicit vNf-v-}
v^{NF}-v_{-} = -\frac{i}{2}\left[(D_tv)(D_1u) - (D_1v)(D_tu) + D_1[v D_tu] - \langle D_x\rangle [v D_1u] \right]
\end{equation}
and
\begin{equation} \label{explicit rNFkg}
r^{NF}_{kg} =-\frac{i}{2}\left[\textit{NL}_{kg} D_1u - (D_1v)\textit{NL}_w + D_1(v\textit{NL}_w)\right].
\end{equation}
From \eqref{explicit rNFkg} and \eqref{def u+- v+-}, together with estimates \eqref{est: bootstrap argument a-priori est} and \eqref{norm_Linfty_xv-},
\begin{equation} \label{chi_xm_rNFkg}
\begin{split}
\|\chi(t^{-\sigma}D_x)&(x_m r^{NF}_{kg})(t,\cdot) \|_{L^2} \lesssim t^{\sigma}\left( \|x_n v_{-}(t,\cdot)\|_{L^\infty}+ \left\| x_n\frac{D_x}{\langle D_x\rangle} v_{-}(t,\cdot)\right\|_{L^\infty}\right) \\
&\times \Big[\left(\|U(t,\cdot)\|_{H^{2,\infty}} + \|\mathrm{R}_1U(t,\cdot)\|_{H^{2,\infty}}\right) \|U(t,\cdot)\|_{L^2} + \|V(t,\cdot)\|_{H^{2,\infty}}\|V(t,\cdot)\|_{L^2}\Big]\\
& + \|V(t,\cdot)\|^2_{H^{1,\infty}}\|V(t,\cdot)\|_{H^1} \le C(A+B)AB\varepsilon^3 t^{-\frac{1}{2}+\sigma+\frac{(\delta+\delta_2)}{2}}.
\end{split}
\end{equation}
Similarly to \eqref{Zm(vgamma-vJ)}, 
\begin{equation}\label{Zm(vNF-v-)}
\begin{split}
&2i\chi(t^{-\sigma}D_x) Z_m (v^{NF}-v_{-}) \\
& = \chi(t^{-\sigma}D_x) \Big[(D_t Z_m v)(D_1u) - (D_1Z_m v)(D_tu)+ D_1[(Z_m v)(D_tu)] - \langle D_x\rangle [(Z_m v)(D_1u)]\\
&\hspace{15pt}+ (D_t  v)(D_1Z_m u) - (D_1 v)(D_t Z_mu)+ D_1[v (D_t Z_mu)] - \langle D_x\rangle [ v(D_1Z_m u)]\\
&\hspace{15pt} - (D_m v)(D_1u) + \delta_{m1}(D_t v)(D_tu) - \delta_{m1} D_t[v (D_t u)] + \frac{D_m}{\langle D_x\rangle} D_t[v (D_1 u)] \\
& \hspace{15pt} -\delta_{m1} (D_t  v)(D_tu) + (D_1  v)(D_m u)-\delta_{m1} D_1[v(D_t u)] + \delta_{m1}\langle D_x\rangle [v(D_tu)]\Big].
\end{split}
\end{equation}
We bound the $L^2$ norm of all products in the first line of the above equality by means of lemma \ref{Lem_appendix:L_estimate of products}, and all the others by the $L^\infty$ norm of the Klein-Gordon factor times the $L^2$ norm of the wave one. In this way we get that, for some new $\chi_1\in C^\infty_0(\mathbb{R}^2)$ and $s>0$ sufficiently large, we derive that
\begin{multline*}
\left\|\chi(t^{-\sigma}D_x)Z_m (v^{NF}-v_{-})(t,\cdot) \right\|_{L^2}\\
\lesssim t^\sigma \left\|\chi_1(t^{-\sigma}D_x)(Z_m v)_\pm(t,\cdot)\right\|_{L^\infty}\|u_\pm(t,\cdot)\|_{L^2}
 + t^{-1}\|(Z_mv)_\pm(t,\cdot)\|_{L^2}\|u_\pm(t,\cdot)\|_{H^s}\\
+ t^\sigma \|v_\pm(t,\cdot)\|_{H^{1,\infty}}\left(\|(Z_mu)_\pm(t,\cdot)\|_{L^2}+\|u_\pm(t,\cdot)\|_{L^2}+\|D_tu_\pm(t,\cdot)\|_{L^2}\right).
\end{multline*}
Consequently, using  estimates \eqref{est: bootstrap argument a-priori est}, \eqref{Hs_norm_DtU} with $s=0$, and \eqref{Linfty_est_VJ}, we obtain that
\begin{equation}\label{est_Zm(vNF-v-)}
\|\chi(t^{-\sigma}D_x)Z_m(v^{NF}-v_{-})(t,\cdot)\|_{L^2}\le C(A+B)B\varepsilon^2 t^{-1+\beta+\frac{\delta+\delta_1}{2}},
\end{equation}
which plugged into \eqref{norm_L2_Lnwidetilde(v)}, along with \eqref{chi_xm_rNFkg}, \eqref{est:Hs_vtilde} and \eqref{est: bootstrap E02}, gives \eqref{est:Lvtilde}.
\endproof
\end{lem}

\section{Last range of estimates} \label{sec_appB: second range of estimates} 

The aim of this section is to show that a-priori estimates \eqref{est: bootstrap argument a-priori est} also infer a moderate growth in time of the $L^2(\mathbb{R}^2)$ norm of $\mathcal{L}^\mu\vt$, for $|\mu|=2$, when this function is restricted to frequencies less or equal than $h^{-\sigma}$, for $\sigma>0$ small.
This is proved in lemma  \ref{Lem_appendix: estimate L2vtilde}. Lemmas from \ref{Lem_appendix:Lm(Zv - vNF-Z)} to \ref{Lem_appendix: L xnrNFkg} are intermediate technical results.

\begin{lem} \label{Lem_appendix:Lm(Zv - vNF-Z)}
Let us consider $v^{NF}$ introduced in \eqref{def vNF} and $\vNFGamma$ as in \eqref{def_vNF-Gamma} with $|I|=1$ and $\Gamma^I=Z_n$, for $n\in \{1,2\}$.
There exists a constant $C>0$ such that, for any $\chi \in C^\infty_0(\mathbb{R}^2)$, $\sigma>0$ small, and every $t\in [1,T]$,
\begin{equation} \label{est:t(Zmv - vNFZn)}
\left\| \chi(t^{-\sigma}D_x)\big[ (Z_n v)_{-} - \vNFGamma\big](t,\cdot)\right\|_{L^2}\le C(A+B)B\varepsilon^2 t^{-1+\beta + \frac{\delta+\delta_1}{2}},
\end{equation}
\begin{multline} \label{est:xm Zn(vNF -v-)}
\left\| \chi(t^{-\sigma}D_x)\big[x_m Z_n(v^{NF}-v_{-})(t,\cdot)\big] \right\|_{L^2} + \left\| \chi(t^{-\sigma}D_x)\big[x_m \big((Z_n v)_{-} - \vNFGamma \big)\big](t,\cdot)\right\|_{L^2}\\
\le C(A+B)B \varepsilon^2 t^{\beta+ \frac{\delta_1+\delta_2}{2}}.
\end{multline}
The same estimates hold true when $Z_n$ is replaced with $\Omega$.
\proof
By comparing equality \eqref{explicit vNfGamma-vJ-}, with $|I|=1$ and $\Gamma^I=Z_n$, with \eqref{Zm(vNF-v-)} we see that $\chi(t^{-\sigma}D_x)(\vNFGamma - (Z_nv)_{-})$ corresponds to the first line in the right hand side of \eqref{Zm(vNF-v-)}.
Therefore, inequality \eqref{est:t(Zmv - vNFZn)} is automatically satisfied after \eqref{est_Zm(vNF-v-)}, which was obtained by estimating the right hand side of \eqref{Zm(vNF-v-)} term by term.
In order to prove \eqref{est:xm Zn(vNF -v-)}, let us consider equality \eqref{Zm(vNF-v-)} but with $\chi(t^{-\sigma}D_x)$ replaced with $\chi(t^{-\sigma}D_x)x_m$.
The $L^2$ norm of each product in the second to fourth line is bounded by
\begin{equation*}
t^\sigma\sum_{\mu,\nu=0}^1 \left\|x^\mu_m \Big(\frac{D_x}{\langle D_x\rangle}\Big)^\nu  v_\pm(t,\cdot)\right\|_{L^\infty}\left(\|(Z_m u)_\pm(t,\cdot)\|_{L^2}+\|u_\pm(t,\cdot)\|_{L^2}+\|D_tu_\pm(t,\cdot)\|_{L^2}\right),
\end{equation*}
and then by the right hand side of \eqref{est:xm Zn(vNF -v-)} after \eqref{est: bootstrap argument a-priori est}, \eqref{Hs_norm_DtU} with $s=0$, and \eqref{norm_Linfty_xv-}.
Using lemma \ref{Lem_appendix:L_estimate of products} with $L=L^2$ and $s>0$ large enough to have $N(s)\ge 2$, we obtain that the $L^2$ norm of products in the first line of (the modified) \eqref{Zm(vNF-v-)} is bounded by
\begin{multline*}
 \sum_{\mu,\nu=0}^1 \left\|\chi_1(t^{-\sigma}D_x)\Big[x^\mu_m \Big(\frac{D_x}{\langle D_x\rangle}\Big)^\nu (Z_m v)_\pm(t,\cdot)\Big]\right\|_{L^\infty}\|u_\pm(t,\cdot)\|_{L^2} \\
+ \sum_{\mu=0}^1 t^{-N(s)}\|x^\mu_m (Z_m v)_\pm(t,\cdot)\|_{L^2}\|u_\pm(t,\cdot)\|_{H^s},
\end{multline*}
for some smooth cut-off $\chi_1$, and hence by the right hand side of \eqref{est:xm Zn(vNF -v-)} after \eqref{est: bootstrap argument a-priori est}, \eqref{norm_L2_xj-GammaIv-} and \eqref{norm_Linfty_xjGammav-} with $\Gamma=Z_m$.
This concludes the proof of \eqref{est:xm Zn(vNF -v-)}.

When $Z_n$ is replaced with $\Omega$, instead of referring to \eqref{Zm(vNF-v-)} one uses that
\begin{equation*}
\begin{split}
2i\Omega \left(v^{NF}-v_{-}\right) & = (D_t \Omega v)(D_1u) - (D_1\Omega v)(D_tu)+ D_1[(\Omega v)(D_tu)] - \langle D_x\rangle [(\Omega v)(D_1u)]\\
&+ (D_t  v)(D_1\Omega u) - (D_1 v)(D_t \Omega u)+ D_1[v (D_t \Omega u)] - \langle D_x\rangle [ v(D_1\Omega u)]\\
& - (D_t  v)(D_2 u) + (D_2 v)(D_t  u)- D_2[v (D_t u)] + \langle D_x\rangle [ v(D_2 u)]
\end{split}
\end{equation*}
and applies the same argument as above to recover the wished estimates.
\endproof
\end{lem}

\begin{lem}\label{Lem_appendix: ZnvNF_LmZnvNF}
Let $v^{\textit{NF}}$ be defined as in \eqref{def vNF}. There exists a constant $C>0$ such that, for any $\chi\in C^\infty_0(\mathbb{R}^2)$, $\sigma>0$ small, $m=1,2$,
\begin{subequations}
\begin{equation} \label{est:ZnvNF}
\left\|\oph(\chi(h^\sigma\xi)) [tZ_nv^{NF}(t,tx)]\right\|_{L^2(dx)}\le CB\varepsilon t^{\frac{\delta_2}{2}},
\end{equation}
\begin{equation}\label{est:Lm_ZnvNF}
\left\|\oph(\chi(h^\sigma\xi))\mathcal{L}_m [tZ_nv^{NF}(t,tx)]\right\|_{L^2(dx)}\le CB\varepsilon t^{\frac{\delta_1}{2}},
\end{equation}
\end{subequations}
for every $t\in [1, T]$.
\proof
Let us write $Z_nv^{NF}$ as follows:
\begin{equation}\label{Zn_vNF}
Z_nv^{NF} = Z_n(v^{NF}-v_{-}) + \left[(Z_n v)_{-} - \vNFGamma\right] + \vNFGamma + \frac{D_n}{\langle D_x\rangle}v^{NF}+\frac{D_n}{\langle D_x\rangle}(v_{-}-v^{NF}),
\end{equation}
with $\vNFGamma$ given by \eqref{def_vNF-Gamma} with $|I|=1$ and $\Gamma^I=Z_n$.
From the fact that $\|tw(t,t\cdot)\|_{L^2}= \|w(t,\cdot)\|_{L^2}$ and estimates \eqref{est: bootstrap argument a-priori est}, \eqref{est_chi_vnf-v}, \eqref{est:vNFGamma},\eqref{est_Zm(vNF-v-)}, \eqref{est:t(Zmv - vNFZn)}, along with the following one
\begin{multline*}
\|\chi(t^{-\sigma}D_x)D_n\langle D_x\rangle^{-1} v^{NF}(t,\cdot)\|_{L^2}\\\le \|\chi(t^{-\sigma}D_x)D_n\langle D_x\rangle^{-1} v_{-}(t,\cdot)\|_{L^2} + \|\chi(t^{-\sigma}D_x)D_n\langle D_x\rangle^{-1}(v_{-} - v^{NF})(t,\cdot)\|_{L^2} 
\le CB\varepsilon  t^\frac{\delta}{2},
\end{multline*}
we immediately obtain \eqref{est:ZnvNF}.

From \eqref{Zn_vNF} we also derive that
\begin{equation}\label{ineq: L-Zn-vNF}
\begin{split}
&\left\|\oph(\chi(h^\sigma\xi))\mathcal{L}_m [t Z_nv^{NF}(t,tx)]\right\|_{L^2(dx)} \lesssim \left\|\oph(\chi(h^\sigma\xi))\mathcal{L}_m \big[tZ_n(v^{NF}-v_{-})(t,tx)\big] \right\|_{L^2(dx)}\\
& + \left\|\oph(\chi(h^\sigma\xi))\mathcal{L}_m \big[t\big((Z_n v)_{-}-\vNFGamma\big)(t,tx)\big]\right\|_{L^2(dx)}+ \left\|\oph(\chi(h^\sigma\xi))\mathcal{L}_m \big[t\vNFGamma(t,tx)\big]\right\|_{L^2(dx)}\\
&+ \left\|\oph(\chi(h^\sigma\xi))\mathcal{L}_m \big[tD_n\langle D_x\rangle^{-1}v^{NF}(t,tx)\big]\right\|_{L^2(dx)}\\
& + \left\|\oph(\chi(h^\sigma\xi))\mathcal{L}_m \big[tD_n\langle D_x\rangle^{-1}(v_{-} - v^{NF})(t,tx)\big]\right\|_{L^2(dx)}.
\end{split}
\end{equation}
By relation \eqref{relation_Zjw_Ljwidetilde(w)} with $w=\vNFGamma$ and estimates \eqref{est:vNFGamma}, \eqref{est_L2_ZmvNFGamma}, \eqref{est_xj_NL-kg-NF-Gamma}, it follows that
\begin{equation*}
\left\|\oph(\chi(h^\sigma\xi))\mathcal{L}_m \big[t\vNFGamma(t,tx)\big]\right\|_{L^2} \le CB\varepsilon t^\frac{\delta_1}{2},
\end{equation*}
while from \eqref{def utilde vtilde} and \eqref{est:Lvtilde} we have that
\begin{equation*}
 \left\|\oph(\chi(h^\sigma\xi))\mathcal{L}_m \big[tD_n\langle D_x\rangle^{-1}v^{NF}(t,tx)\big]\right\|_{L^2}\le CB\varepsilon t^\frac{\delta_2}{2}.
\end{equation*}
The remaining $L^2$ norms in the right hand side of \eqref{ineq: L-Zn-vNF} are estimated reminding definition \eqref{def Lj} of $\mathcal{L}_m$ and using the fact that
\begin{equation}\label{dec_norm_Lm}
\|\oph(\chi(h^\sigma\xi))\mathcal{L}_m [tw(t,tx)]\|_{L^2(dx)}\lesssim \|\chi(t^{-\sigma}D_x)[x_m w(t,\cdot)]\|_{L^2}+ t\|\chi(t^{-\sigma}D_x) w(t,\cdot)\|_{L^2}.
\end{equation}
Therefore, by \eqref{est_Zm(vNF-v-)} and lemma \ref{Lem_appendix:Lm(Zv - vNF-Z)} we derive that
\begin{multline*}
\left\|\oph(\chi(h^\sigma\xi))\mathcal{L}_m \big[tZ_n(v^{NF}-v_{-})(t,tx)\big] \right\|_{L^2(dx)}\\+  \left\|\oph(\chi(h^\sigma\xi))\mathcal{L}_m \big[t\big((Z_n v)_{-}-\vNFGamma\big)(t,tx)\big]\right\|_{L^2(dx)} 
\le C(A+B)B\varepsilon^2 t^{\beta + \frac{\delta+\delta_1}{2}},
\end{multline*}
while from \eqref{est_chi_vnf-v}, a-priori estimates, together with the following inequality
\begin{multline*}
\|\chi(t^{-\sigma}D_x)[x_m (v_{-} - v^{NF})](t,\cdot)\|_{L^2} \lesssim \sum_{\mu,\nu=0}^1t^\sigma \left\|x_m^\mu \Big(\frac{D_x}{\langle D_x\rangle}\Big)^\nu v_{-}(t,\cdot) \right\|_{L^\infty}\|u_\pm(t,\cdot)\|_{L^2}\\
\le C(A+B)B\varepsilon^2 t^{\sigma+ \frac{\delta+\delta_2}{2}},
\end{multline*}
which follows by \eqref{explicit vNf-v-}, \eqref{def u+- v+-}, \eqref{norm_Linfty_xv-}, \eqref{est: boostrap vpm}, \eqref{est: bootstrap Enn}, we derive
\begin{equation*}
\left\|\oph(\chi(h^\sigma\xi))\mathcal{L}_m \big[tD_n\langle D_x\rangle^{-1}(v_{-} - v^{NF})(t,tx)\big]\right\|_{L^2} \le C(A+B)B\varepsilon^2 t^{\sigma + \frac{(\delta+\delta_2)}{2}}.
\end{equation*}
\endproof
\end{lem}

In the following lemma we are going to prove that the product of the semiclassical wave function $\ut$ with the Klein-Gordon one $\vt$ enjoys a better $L^2$ (resp. $L^\infty$) estimate than the one roughly obtained by taking the $L^2$ (resp. $L^\infty$) norm of the former times the $L^\infty$ norm of the latter.
Estimates
\begin{align*}
\left\| \vt\ut (t,\cdot)\right\|_{L^2}& \lesssim \|\vt(t,\cdot)\|_{L^\infty}\|\ut(t,\cdot)\|_{L^2}\le CAB\varepsilon^2 h^{-\frac{\delta}{2}},\\
\left\| \vt\ut (t,\cdot)\right\|_{L^\infty}& \lesssim \|\vt(t,\cdot)\|_{L^\infty}\|\ut(t,\cdot)\|_{L^\infty}\le CA^2\varepsilon^2 h^{-\frac{1}{2}-\frac{\delta}{2}},
\end{align*}
which follows from \eqref{est:utilde-Hs}, \eqref{Hrho_infty_utilde_appendix}, \eqref{Hrho_infty_vtilde_appendix}, can be in fact improved of a factor $h^{1/2}$ (see \eqref{est_L2Linfty_Vtilde_Utilde}).
This comes from the fact that the main contribution to $\ut$ is localized around manifold $\Lw$ introduced in \eqref{def_Lw}, whereas $\vt$ concentrates around $\Lkg$ defined in \eqref{def_Lkg}, and these two manifolds are disjoint.

\begin{figure}[h]
\begin{center}
\begin{tikzpicture}[scale=1.9]

\draw[->] (0,-1.7) -- (0,1.7);
\draw[->] (-2,0) -- (2,0);
\node[below] at (1.9,0) {\small $x$};
\node[left] at (0,1.6) {\small $\xi$};

\draw[blue] (-1,-1.7) -- (-1,1.7);
\draw[blue] (1,-1.7) -- (1,1.7);
\node[below left] at (-1,0) {$-1$};
\node[below right] at (1,0) {$1$};
\node[right,blue] at (1,1) {$\Lw$};

\draw[red,thick] [domain=0:0.87] plot(\x, {\x/((1-\x^2)^(1/2))});
\draw[red,thick,red] [domain=0:0.87] plot(-\x, {-\x/((1-\x^2)^(1/2))});
\node[right,red] at (-0.8,-1.6) {$\Lkg$};
\end{tikzpicture}
\end{center}

\caption{Manifolds $\Lkg$ and $\Lw$.}
\label{picture: Lkg_Lw}
\end{figure}

\begin{lem} \label{Lem_appendix:product_Vtilde_Utilde}
Let $h=t^{-1}$, $\widetilde{u}, \widetilde{v}$ be defined in \eqref{def utilde vtilde}, $a_0(\xi)\in S_{0,0}(1)$, and $b_1(\xi)=\xi_j$ or $b_1(\xi)=\xi_j\xi_k|\xi|^{-1}$, with $j,k\in \{1,2\}$. 
There exists a constant $C>0$ such that, for any $\chi, \chi_1\in C^\infty_0(\mathbb{R}^2)$, $\sigma>0$, and every $t\in [1,T]$, we have that
\begin{subequations} \label{est_L2Linfty_Vtilde_Utilde}
\begin{equation} \label{est_L2_Vtilde_Utilde}
\big\|[\oph(\chi(h^\sigma\xi)a_0(\xi))\widetilde{v}(t,\cdot)][\oph(\chi_1(h^\sigma\xi)b_1(\xi))\widetilde{u}(t,\cdot)]\big\|_{L^2} \le C(A+B)B\varepsilon^2 h^{\frac{1}{2}-\beta-\frac{\delta+\delta_1}{2}},
\end{equation}
\begin{equation} \label{est_Linfty_Vtilde Utilde}
\big\|[\oph(\chi(h^\sigma\xi)a_0(\xi))\widetilde{v}(t,\cdot)][\oph(\chi_1(h^\sigma\xi)b_1(\xi))\widetilde{u}(t,\cdot)]\big\|_{L^\infty} \le C(A+B)B\varepsilon^2 h^{-\beta-\frac{\delta+\delta_1}{2}},
\end{equation}
\end{subequations}
with $\beta>0$ small, $\beta\rightarrow 0$ as $\sigma\rightarrow 0$.
\proof
Before entering in the details of the proof, we warn the reader that $C$ and $\beta$ denote two positive constants that may change line after line, with $\beta\rightarrow 0$ as $\sigma\rightarrow 0$. Also, we will denote by $R(t,x)$ any contribution, in what follows, that satisfies inequalities \eqref{est_L2Linfty_Vtilde_Utilde}, and by $\chi_2$ a smooth cut-off function, identically equal to 1 on the support of $\chi_1$, so that
\[\oph(\chi_1(h^\sigma\xi))\widetilde{u} = \oph(\chi_1(h^\sigma\xi))\oph(\chi_2(h^\sigma\xi))\widetilde{u},\] 
assuming that at any time $\widetilde{u}$ can be replaced with $\oph(\chi_2(h^\sigma\xi))\widetilde{u}$. Finally, it is useful to remind that from \eqref{def utilde vtilde}, \eqref{def uNF}, \eqref{Hrho-infty norm uNF- u-}, \eqref{Hrho-infty norm R(uNF-u-)}, and a-priori estimates,
\begin{equation} \label{Hrho_infty_utilde_appendix}
\|\widetilde{u}(t,\cdot)\|_{H^{\rho+1,\infty}_h}+\sum_{|\mu|=1}\| \oph\big((\xi|\xi|^{-1})^\mu\big)\widetilde{u}(t,\cdot)\|_{H^{\rho+1,\infty}_h} \le CA\varepsilon h^{-\frac{1}{2}},
\end{equation}
while by \eqref{def vNF}, \eqref{est_Hsinfty_vnf-v} (with $\theta \ll 1$ small enough) and a-priori estimates,
\begin{equation}\label{Hrho_infty_vtilde_appendix}
\|\widetilde{v}(t,\cdot)\|_{H^{\rho,\infty}_h}\le CA\varepsilon,
\end{equation}
for every $t\in [1,T]$.

First of all, we take $\gamma\in C^\infty_0(\mathbb{R}^2)$ equal to 1 in a neighbourhood of the origin and with sufficiently small support, and define
\begin{gather*}
\widetilde{v}_{\Lambda_{kg}}(t,x):= \oph\Big(\gamma\Big(\frac{x-p'(\xi)}{\sqrt{h}}\Big)\chi(h^\sigma\xi)a_0(\xi)\Big)\widetilde{v}(t,x), \\
 \widetilde{v}_{\Lambda^c_{kg}}(t,x):=\oph\Big((1-\gamma)\Big(\frac{x-p'(\xi)}{\sqrt{h}}\Big)\chi(h^\sigma\xi)a_0(\xi)\Big)\widetilde{v}(t,x),
\end{gather*}
with $p(\xi):=\langle \xi\rangle$, so that 
\begin{equation}\label{dec_vt_app}
\oph(\chi(h^\sigma\xi)a_0(\xi))\widetilde{v}=\widetilde{v}_{\Lambda_{kg}}+\widetilde{v}_{\Lambda^c_{kg}}.
\end{equation}
The following estimates hold:
\begin{subequations}
\begin{gather}
\|\widetilde{v}_{\Lambda_{kg}}(t,\cdot)\|_{L^\infty}\le CA\varepsilon h^{-\beta}, \label{est:vLambda_appendix} \\
\|\widetilde{v}_{\Lambda^c_{kg}}(t,\cdot)\|_{L^\infty} \le CB\varepsilon h^{\frac{1}{2}-\beta-\frac{\delta_1}{2}}. \label{est:vLambdac_appendix}
\end{gather}
\end{subequations}
The former one is a straight consequence of proposition \ref{Prop:Continuity Lp-Lp} with $p=+\infty$ and \eqref{Hrho_infty_vtilde_appendix}.
On the other hand, if we write
\begin{equation*}
(1-\gamma)\Big(\frac{x-p'(\xi)}{\sqrt{h}}\Big)\chi(h^\sigma\xi)a_0(\xi)  =\sum_{j=1}^2 \gamma_1^j\Big(\frac{x-p'(\xi)}{\sqrt{h}}\Big)\chi(h^\sigma\xi)a_0(\xi) \Big(\frac{x_j-p'_j(\xi)}{\sqrt{h}}\Big)
\end{equation*}
with $\gamma_1^j(z):=(1-\gamma)(z)z_j|z|^{-2}$ such that $|\partial^\alpha_z\gamma_1^j(z)|\lesssim \langle z\rangle^{-1-|\alpha|}$, and use \eqref{sharp gammatilde and its argument} with $c(x,\xi)=\chi(h^\sigma\xi)a_0(\xi)$, we obtain that
\begin{equation} \label{est:vLmabdakgc_appendix_preliminary}
\begin{split}
\|\widetilde{v}_{\Lambda^c_{kg}}(t,\cdot)\|_{L^\infty}&\lesssim \sum_{j=1}^2 \sqrt{h} \left\|\oph\Big(\gamma_1^j\Big(\frac{x-p'(\xi)}{\sqrt{h}}\Big)\chi(h^\sigma\xi)a_0(\xi)\Big)\mathcal{L}_j\widetilde{v}(t,\cdot) \right\|_{L^\infty}\\
& + \sum_{j=1}^2 \sqrt{h} \left\|\oph\Big(\gamma_1^j\Big(\frac{x-p'(\xi)}{\sqrt{h}}\Big)\partial_j\big(\chi(h^\sigma\xi)a_0(\xi)\big)\Big)\widetilde{v}(t,\cdot) \right\|_{L^\infty} \\
& + \sum_{j=1}^2 \sum_{|\alpha|=2}\sqrt{h} \left\|\oph\Big((\partial^\alpha\gamma_1^j)\Big(\frac{x-p'(\xi)}{\sqrt{h}}\Big)\chi(h^\sigma\xi)a_0(\xi)(\partial^\alpha_\xi p')(\xi)\Big)\widetilde{v}(t,\cdot) \right\|_{L^\infty} \\
&+ \|\oph(r(x,\xi))\widetilde{v}(t,\cdot)\|_{L^\infty},
\end{split}
\end{equation}
with $r\in h^{1-\beta}S_{\frac{1}{2},\sigma}(\langle\frac{x-p'(\xi)}{\sqrt{h}}\rangle^{-1})$.
Since $\gamma_1^j$ vanishes in a neighbourhood of the origin, we derive from inequality \eqref{est_1L-Linfty}, equation \eqref{KG equation vNF} and relation \eqref{relation_Zjw_Ljwidetilde(w)} with $w=v^{NF}$, lemmas \ref{Lem: from energy to norms in sc coordinates-KG}, \ref{Lem_appendix: ZnvNF_LmZnvNF}, and estimate \eqref{chi_xm_rNFkg}, that the first sum in the above right hand side is bounded by the right hand side of \eqref{est_Linfty_Vtilde Utilde}.
The same is true for the above second and third sums after \eqref{est_1L-Linfty} and lemma \ref{Lem: from energy to norms in sc coordinates-KG}, and for the above latter $L^\infty$ norm because of proposition \ref{Prop : Continuity from $L^2$ to L^inf} and estimate \eqref{est:Hs_vtilde}.

After decomposition \eqref{dec_vt_app} and estimates \eqref{est:utilde-Hs}, \eqref{Hrho_infty_utilde_appendix}, and \eqref{est:vLambdac_appendix}, we see that
\begin{equation*}
\oph(\chi(h^\sigma\xi)a_0(\xi))\widetilde{v} \, \oph(\chi(h^\sigma\xi)b_1(\xi))\widetilde{u} = \widetilde{v}_{\Lambda_{kg}}\oph(\chi(h^\sigma\xi)b_1(\xi))\widetilde{u} + R(t,x).
\end{equation*}
For some suitably supported $\chi_0\in C^\infty_0(\mathbb{R}^2)$, $\varphi\in C^\infty_0(\mathbb{R}^2\setminus\{0\})$, we also consider the following decomposition
\begin{equation*}
\oph(\chi_1(h^\sigma\xi)b_1(\xi))\widetilde{u} = \oph(\chi_0(h^{-1}\xi)b_1(\xi))\widetilde{u} + \sum_k \oph\big((1-\chi_0)(h^{-1}\xi)\varphi(2^{-k}\xi)\chi_1(h^\sigma\xi)b_1(\xi)\big)\widetilde{u},
\end{equation*}
and observe that, from proposition \ref{Prop : Continuity on H^s} and the classical Sobolev injection, 
\begin{gather*}
\left\|\oph(\chi_0(h^{-1}\xi)b_1(\xi))\widetilde{u}(t,\cdot)\right\|_{L^2}+ \left\|\oph(\chi_0(h^{-1}\xi)b_1(\xi))\widetilde{u}(t,\cdot)\right\|_{L^\infty}\lesssim h\|\widetilde{u}(t,\cdot)\|_{L^2}.
\end{gather*}
Combining the above decomposition and estimate with \eqref{est:vLambda_appendix} and \eqref{est:utilde-Hs} we derive that
\begin{equation}  \label{dec_k_app}
\widetilde{v}_{\Lambda_{kg}}\oph(\chi(h^\sigma\xi)b_1(\xi))\widetilde{u} = \sum_k \widetilde{v}_{\Lambda_{kg}}\oph(\phi_k(\xi)b_1(\xi))\widetilde{u} + R(t,x),
\end{equation}
where $\phi_k(\xi):=(1-\chi_0)(h^{-1}\xi)\varphi(2^{-k}\xi)\chi(h^\sigma\xi)$. 
We can further decompose $\oph(\phi_k(\xi)b_1(\xi))\widetilde{u}$ by defining
\begin{equation*}
\begin{gathered}
\widetilde{u}^k_{\Lambda_w}(t,x):=\oph\Big(\gamma\Big(\frac{x|\xi|-\xi}{h^{1/2-\sigma}}\Big)\phi_k(\xi)b_1(\xi)\Big)\widetilde{u}(t,x),\\
\widetilde{u}^k_{\Lambda^c_w}(t,x):=\oph\Big((1-\gamma)\Big(\frac{x|\xi|-\xi}{h^{1/2-\sigma}}\Big)\phi_k(\xi)b_1(\xi)\Big)\widetilde{u}(t,x),
\end{gathered}
\end{equation*}
and observe that
\begin{multline*}
\left\|\widetilde{u}^k_{\Lambda^c_w}(t,\cdot)\right\|_{L^2} \lesssim h^{\frac{1}{2}-\beta}\left[\|\widetilde{u}(t,\cdot)\|_{L^2} +\sum_{\mu, |\nu|=0}^1\|(\theta_0(x) \Omega_h)^\mu \mathcal{M}^\nu \oph(\chi_2(h^\sigma\xi))\widetilde{u}(t,\cdot)\|_{L^2}\right]\\\le CB\varepsilon h^{\frac{1}{2}-\beta-\frac{\delta_1}{2}},
\end{multline*}
and
\begin{multline*}
\left\|\widetilde{u}^k_{\Lambda^c_w}(t,\cdot)\right\|_{L^\infty} \lesssim h^{-\beta}\left[\|\widetilde{u}(t,\cdot)\|_{L^2} +\sum_{\mu, |\nu|=0}^1\|(\theta_0(x) \Omega_h)^\mu \mathcal{M}^\nu \oph(\chi_2(h^\sigma\xi))\widetilde{u}(t,\cdot)\|_{L^2}\right]\\\le CB\varepsilon h^{-\beta-\frac{\delta_1}{2}},
\end{multline*}
as follows by using the following equality
\begin{equation*}
(1-\gamma)\Big(\frac{x|\xi|-\xi}{h^{1/2-\sigma}}\Big)\phi_k(\xi)b_1(\xi) = \sum_{j=1}^2\gamma_1^j\Big(\frac{x|\xi|-\xi}{h^{1/2-\sigma}}\Big)\Big(\frac{x_j|\xi|-\xi_j}{h^{1/2-\sigma}}\Big)\phi_k(\xi)b_1(\xi),
\end{equation*}
with $\gamma_1^j(z):=(1-\gamma)(z)z_j|z|^{-2}$, together with \eqref{est: L2 Linfty with L} with $a\equiv 1$, $p=1$, and lemma \ref{Lem: from energy to norms in sc coordinates-WAVE}.
Then, as the sum over $k$ in the right hand side of \eqref{dec_k_app} is actually restricted to indices $k$ such that $h\lesssim 2^k\lesssim h^{-\sigma}$, the above estimates and \eqref{est:vLambda_appendix} imply that
\[\sum_k \widetilde{v}_{\Lambda_{kg}}\oph(\phi_k(\xi)b_1(\xi))\widetilde{u} = \sum_k \widetilde{v}_{\Lambda_{kg}}\ut^k_{\Lw} + R(t,x).\]
Moreover, using lemma \ref{Lem:family_thetah}, symbolic calculus and remark \ref{Remark:symbols_with_null_support_intersection}, each $\vt_{\Lkg}\ut^k_{\Lw}$ in the above right hand side can be replaced with
\[\frac{\theta_h(x)}{|x|^2-1}\vt_{\Lkg} \ (|x|^2-1)\ut^k_{\Lw}\]
up to a new remainder $R(t,x)$.
Since $|\theta_h(x) (|x|^2-1)^{-1}|\lesssim h^{-2\sigma}$ on the support of $\theta_h(x)$, from proposition \ref{Prop : Continuity on H^s} and estimates \eqref{est:Hs_vtilde}, \eqref{est:vLambda_appendix}, we get that
\begin{subequations} \label{est_uLambda_vLambda}
\begin{gather}
\left\|\theta_h(x)\widetilde{v}_{\Lambda_{kg}}(t,\cdot) \widetilde{u}^k_{\Lambda_w} (t,\cdot)\right\|_{L^2}\le CB\varepsilon h^{-\frac{\delta}{2}-\beta}\|\theta_h(x)(|x|^2-1)\widetilde{u}^k_{\Lambda_w}(t,\cdot)\|_{L^\infty},\\
\left\|\theta_h(x)\widetilde{v}_{\Lambda_{kg}} (t,\cdot)\widetilde{u}^k_{\Lambda_w} (t,\cdot)\right\|_{L^\infty}\le CA\varepsilon h^{-\beta}\|\theta_h(x)(|x|^2-1)\widetilde{u}^k_{\Lambda_w}(t,\cdot)\|_{L^\infty}.
\end{gather}
\end{subequations}
Then the end of the proof relies on the fact that $\theta_h(x)(|x|^2-1)\widetilde{u}^k_{\Lambda_w}$ can be expressed in terms of $h\Mcal \ut$. In fact, for a fixed $N\in\mathbb{N}$ and up to some negligible multiplicative constants, we have from proposition \ref{Prop: a sharp b} that
\begin{equation} \label{symb_prod_(x2-1)}
\begin{split}
\left[\theta_h(x)(|x|^2-1)\right]\sharp &\left[\gamma\Big(\frac{x|\xi|-\xi}{h^{1/2-\sigma}}\Big)\phi_k(\xi)b_1(\xi)\right] = \gamma\Big(\frac{x|\xi|-\xi}{h^{1/2-\sigma}}\Big)\phi_k(\xi)b_1(\xi)\theta_h(x) (|x|^2-1)\\
&+ h\left\{\theta_h(x) (|x|^2-1), \gamma\Big(\frac{x|\xi|-\xi}{h^{1/2-\sigma}}\Big)\phi_k(\xi)b_1(\xi) \right\}\\
& + \sum_{ |\alpha|=2}^{N-1} h^{|\alpha|} \partial^{\alpha}_x\left[\theta_h(x)(|x|^2-1)\right] \partial^\alpha_\xi\left[\gamma\Big(\frac{x|\xi|-\xi}{h^{1/2-\sigma}}\Big)\phi_k(\xi)b_1(\xi)\right] +r_N(x,\xi),
\end{split}
\end{equation}
with
\begin{multline}\label{rN_lemma_appendixB}
r_N(x,\xi) = \frac{h^N}{(\pi h)^4}\sum_{|\alpha|=N}\int e^{\frac{2i}{h}(\eta\cdot z-y\cdot\zeta)} \int_0^1 \partial^\alpha_x[\theta_h(x)(|x|^2-1)]|_{(x+tz)}(1-t)^{N-1}dt \\ 
\times \partial^\alpha_\xi\left[\gamma\Big(\frac{x|\xi|-\xi}{h^{1/2-\sigma}}\Big)\phi_k(\xi)b_1(\xi)\right]|_{(x+y, \xi+\eta)} dy dz d\eta d\zeta.
\end{multline}
As
\begin{equation*}
|x|^2-1 = x\cdot x - \frac{\xi\cdot\xi}{|\xi|^2} = (x|\xi|-\xi)\cdot \frac{x}{|\xi|} + (x|\xi|-\xi)\cdot\frac{\xi}{|\xi|^2},
\end{equation*}
the first term in the right hand side of \eqref{symb_prod_(x2-1)} appears to be linear combination of products of the form $\gamma\big(\frac{x|\xi|-\xi}{h^{1/2-\sigma}}\big)\phi_k(\xi)a(x)b_0(\xi)(x_j|\xi|-\xi_j)$, for some smooth compactly supported function $a(x)$, and $b_0(\xi)$ such that $|\partial^\alpha b_0(\xi)|\lesssim |\xi|^{-|\alpha|}$. From \eqref{est: Linfty of gamma1 with L} and lemma \ref{Lem: from energy to norms in sc coordinates-WAVE}, we hence deduce that
\begin{subequations}
\begin{gather}\label{est_first_term}
\left\| \oph\Big(\gamma\Big(\frac{x|\xi|-\xi}{h^{1/2-\sigma}}\Big)\phi_k(\xi)b_1(\xi)\theta_h(x) (|x|^2-1)\Big)\widetilde{u}(t,\cdot)\right\|_{L^\infty} \le CB\varepsilon h^{\frac{1}{2}-\beta-\frac{\delta_1}{2}}.
\end{gather}
\end{subequations}
An explicit computation shows that
\begin{multline*}
 h \left\{\theta_h(x) (|x|^2-1), \gamma\Big(\frac{x|\xi|-\xi}{h^{1/2-\sigma}}\Big)\phi_k(\xi)b_1(\xi) \right\}\\
 =\sum_i h^{\frac{1}{2}+\sigma}\partial_{x_i}[\theta_h(x)(|x|^2-1)]\sum_j(\partial_j\gamma)\Big(\frac{x|\xi|-\xi}{h^{1/2-\sigma}}\Big)\Big(x_j\frac{\xi_i}{|\xi|}-\delta_{ij}\Big)\phi_k(\xi)b_1(\xi) \\
+h  \gamma\Big(\frac{x|\xi|-\xi}{h^{1/2-\sigma}}\Big) \partial_x[\theta_h(x)(|x|^2-1)]\partial_\xi[\phi_k(\xi)b_1(\xi)],
\end{multline*}
with $\delta_{ij}=1$ being the Kronecker delta.
One the one hand, since the first contribution to the above right hand side is still supported for $|x|<1-ch^{2\sigma}$,
we can multiply and divide it by $|x|^2-1$ so that it writes as linear combination of terms of the form $h^{\frac{1}{2}-\sigma}\gamma_1\big(\frac{x|\xi|-\xi}{h^{1/2-\sigma}}\big)\phi_k(\xi)a(x)b_0(\xi)(x_j|\xi|-\xi_j)$, for a new $\gamma_1\in C^\infty_0(\mathbb{R}^2)$, and some new $a(x), b_0(\xi)$ with the same properties as the ones we considered before.
On the other hand, as $\partial_\xi[\phi_k(\xi)b_1(\xi)]$ is uniformly bounded and supported for frequencies of size $2^k$, the second term in the above right hand side writes as linear combination of products of the form $h\gamma\big(\frac{x|\xi|-\xi}{h^{1/2-\sigma}}\big)\phi^1_k(\xi)a(x)b_0(\xi)$, for some new $\phi^1_k\in C^\infty_0(\mathbb{R}^2\setminus\{0\})$.
Therefore, inequality \eqref{est: Linfty of gamma1 with L}, proposition \ref{Prop : continuity of Op(gamma1):X to L2}, and lemma \ref{Lem: from energy to norms in sc coordinates-WAVE}, give that
\begin{gather}\label{est_second_term}
\left\|  h \oph\Big(\Big\{\theta_h(x) (|x|^2-1), \gamma\Big(\frac{x|\xi|-\xi}{h^{1/2-\sigma}}\Big)\phi_k(\xi)b_1(\xi) \Big\}\Big)\widetilde{u}(t,\cdot)\right\|_{L^\infty} \le CB\varepsilon h^{\frac{1}{2}-\beta-\frac{\delta_1}{2}}.
\end{gather}
As concerns $|\alpha|$-order terms, for each fixed $2\le |\alpha|\le N-1$, we find using \eqref{derivatives of gamma_n 2} that they are given by
\begin{multline*}
h^{|\alpha|}\gamma\Big(\frac{x|\xi|-\xi}{h^{1/2-\sigma}}\Big) \partial^\alpha_x[\theta_h(x)(|x|^2-1)]\partial^\alpha_\xi (\phi_k(\xi)b_1(\xi))\\
 + \sum_{\substack{|\beta_1|+|\beta_2|=|\alpha|\\|\beta_1|\ge 1}}\sum_{j=1}^{|\beta_1|}h^{|\alpha| - j(\frac{1}{2}-\sigma)} \gamma_j\Big(\frac{x|\xi|-\xi}{h^{1/2-\sigma}}\Big)\widetilde{\theta}_j(x)b_{j-|\beta_1|}(\xi)\partial^{\beta_2}_\xi(\phi_k(\xi)b_1(\xi)),
\end{multline*}
for some $\gamma_j, \widetilde{\theta}_j\in C^\infty_0(\mathbb{R}^2)$. Since $|\alpha|\ge 2$ and $|\partial^\mu_\xi(\phi_k(\xi)b_1(\xi))|\lesssim 2^{-k(|\mu|-1)}$, for any $\mu\in\mathbb{N}^2$, by proposition
\ref{Prop : continuity of Op(gamma1):X to L2} and lemma \ref{Lem: from energy to norms in sc coordinates-WAVE} we obtain that the action of their quantization on $\widetilde{u}$ is estimated in the uniform norm by
\begin{multline} \label{est_third_term}
\left[h^{|\alpha|-\frac{1}{2}-\beta}2^{-k(|\alpha|-1)} + \sum_{1\le j\le |\alpha|} h^{|\alpha|-j(\frac{1}{2}-\sigma)}2^{k(j+1-|\alpha|)}h^{-\frac{1}{2}-\beta}\right]\\
\times \left[\|\widetilde{u}(t,\cdot)\|_{L^2} +\sum_{\mu, |\nu|=0}^1\|(\theta \Omega_h)^\mu \oph(\chi_1(h^\sigma\xi)\widetilde{u}(t,\cdot)\|_{L^2}\right]\le CB\varepsilon h^{\frac{1}{2}-\beta-\frac{\delta_1}{2}}.
\end{multline}
Finally, by integrating in $dyd\zeta$ and using \eqref{derivatives of gamma_n 1} in \eqref{rN_lemma_appendixB} we find that $r_N(x,\xi)$ can be written as
\begin{multline*}
\sum_{j\le N} h^{N-j(\frac{1}{2}-\sigma)}\frac{1}{(\pi h)^2}\int e^{\frac{2i}{h}\eta\cdot z}\int_0^1 \theta_N(x+tz)(1-t)^{N-1}dt \\
\times \gamma_j\Big(\frac{x|\xi+\eta|-(\xi+\eta)}{h^{1/2-\sigma}}\Big)\phi_k^j(\xi+\eta)b_{j+1-N}(\xi+\eta) dzd\eta,
\end{multline*}
for some new smooth compactly supported $\theta_N, \gamma_j, \phi_k^j$.From the last part of proposition \ref{Prop : Linfty est of integral remainders} then follows that the quantization of the above integral is a bounded operator from $L^2$ to $L^\infty$, with norm controlled by
\begin{equation*}\sum_{\substack{j\le N\\ i\le 6}}h^{N-j(\frac{1}{2}-\sigma)}2^{k(1+j-N)}(h^{-\frac{1}{2}+\sigma}2^k)^i (h^{-1}2^k)\lesssim h 
\end{equation*}
if $N$ is sufficiently large (e.g. $N\ge 10$), and consequently that
\begin{equation}
 \|\oph(r_N(x,\xi))\widetilde{u}(t,\cdot)\|_{L^\infty}\lesssim h\|\widetilde{u}^k(t,\cdot)\|_{L^2} \le CB\varepsilon h^{1-\frac{\delta}{2}}.
\end{equation}
Finally, summing up the above estimates with formulas from \eqref{symb_prod_(x2-1)} to \eqref{est_third_term} we obtain that
\begin{gather*}
\|\theta_h(x)(|x|^2-1)\widetilde{u}_{\Lambda_w}(t,\cdot)\|_{L^\infty}\lesssim CB\varepsilon h^{\frac{1}{2}-\beta-\frac{\delta_1}{2}},
\end{gather*}
which injected in \eqref{est_uLambda_vLambda} gives that $\theta_h(x)\widetilde{v}_{\Lambda_{kg}}\widetilde{u}^k_{\Lambda_w}$ is a remainder $R(t,x)$.
That concludes the proof of the statement.
\endproof
\end{lem}

A similar result to the one proved in lemma \ref{Lem_appendix:product_Vtilde_Utilde} holds true when $\widetilde{u}$ in the left hand side of \eqref{est_L2Linfty_Vtilde_Utilde} is replaced with 
\begin{equation} \label{def_utildeJ}
\widetilde{u}^J(t,x):= t (\Gamma u)_{-}(t,tx)
\end{equation}
with $\Gamma \in \{\Omega, Z_m, m=1,2\}$ being a Klainerman vector field, as briefly shown in the following:

\begin{lem} \label{Lem_appendix: product vtilde_utildeJ}
Let $h=t^{-1}$, $\widetilde{v}$ be defined in \eqref{def utilde vtilde}, $\widetilde{u}^J$ as in \eqref{def_utildeJ},  $a_0(\xi)\in S_{0,0}(1)$, and $b_1(\xi)=\xi_j$ or $b_1(\xi)=\xi_j\xi_k|\xi|^{-1}$, with $j,k\in \{1,2\}$. 
There exists a constant $C>0$ such that, for any $\chi, \chi_1\in C^\infty_0(\mathbb{R}^2)$, $\sigma>0$, and every $t\in [1,T]$, we have that
\begin{subequations} 
\begin{equation} \label{est_vtilde_utildeJ_L2}
\big\|[\oph(\chi(h^\sigma\xi)a_0(\xi))\widetilde{v}(t,\cdot)][\oph(\chi_1(h^\sigma\xi)b_1(\xi))\widetilde{u}^J(t,\cdot)]\big\|_{L^2} \le C(A+B)B\varepsilon^2 h^{\frac{1}{2}-\beta'},
\end{equation}
\begin{equation} \label{est_Linfty_vtilde_utildeJ}
\big\|[\oph(\chi(h^\sigma\xi)a_0(\xi))\widetilde{v}(t,\cdot)][\oph(\chi_1(h^\sigma\xi)b_1(\xi))\widetilde{u}^J(t,\cdot)]\big\|_{L^\infty} \le C(A+B)B\varepsilon^2 h^{-\beta'},
\end{equation}
\end{subequations}
with $\beta'>0$ small, $\beta\rightarrow 0$ as $\sigma,\delta_0\rightarrow 0$.
\proof
The proof of this result is analogous to that of lemma \ref{Lem_appendix:product_Vtilde_Utilde} except that, instead of referring to \eqref{Hrho_infty_utilde_appendix}, we should use that
\begin{equation}\label{est:Hrho_utildeJ}
\|\oph(\chi(h^\sigma\xi))\widetilde{u}^J(t,\cdot)\|_{H^{\rho+1,\infty}_h}+\sum_{|\mu|=1}\| \oph\big(\chi(h^\sigma\xi))(\xi|\xi|^{-1})^\mu\big)\widetilde{u}^J(t,\cdot)\|_{H^{\rho+1,\infty}_j} \le CA\varepsilon h^{-\frac{1}{2}-\beta-\frac{\delta_1}{2}},
\end{equation}
which is the semiclassical translation of \eqref{Linfty_est_UJ}, and to lemma \ref{Lem_appendix: L^2 estimates uJ} instead of lemma \ref{Lem: from energy to norms in sc coordinates-WAVE}.
\endproof
\end{lem}

\begin{lem}\label{Lem_appendix:Linfty_est_rNFkg}Let $a_0(\xi)\in S_{0,0}(1)$, $b_1(\xi)\in \{\xi_j, \xi_j\xi_k|\xi|^{-1}, |\xi|, j,k=1,2\}$, $b_0(\xi)\in \{1,\xi_j|\xi|^{-1}, j=1,2\}$. There exists a constant $C>0$ such that, for any $\chi\in C^\infty_0(\mathbb{R}^2)$, $\sigma>0$ small, and every $t\in [1,T]$,
\begin{equation}\label{est_Linfty_prod_a0v_b1u_R1u}
\left\| \chi(t^{-\sigma}D_x) \big[[a_0(D_x)v_{-}] [b_1(D_x)u_{-}] b_0(D_x)u_{-}\big](t,\cdot)\right\|_{L^\infty} \le C(A+B)AB\varepsilon^3 t^{-\frac{5}{2}+\beta + \frac{\delta+\delta_1}{2}},
\end{equation}
with $\beta>0$ small, $\beta\rightarrow 0$ as $\sigma\rightarrow0$. Consequently
\begin{equation}\label{est_Linfty_chi_rnfkg}
\left\|\chi(t^{-\sigma}D_x)r^{NF}_{kg}(t,\cdot) \right\|_{L^\infty}\le C(A+B)AB\varepsilon^3 t^{-\frac{5}{2}+\beta + \frac{\delta+\delta_1}{2}},
\end{equation}
where $r^{NF}_{kg}$ is given by \eqref{explicit rNFkg}.
\proof
We warn the reader that we denote by $C$ and $\beta$ two positive constants that may change line after line during this proof, with $\beta\rightarrow 0$ as $\sigma\rightarrow 0$.
Moreover, we are going to denote generically by $R(t,x)$ each term satisfying 
\[\|R(t,\cdot)\|_{L^\infty}\le C(A+B)AB\varepsilon^3 t^{-\frac{5}{2}+\beta + \frac{\delta+\delta_1}{2}}.\]

From lemma \ref{Lem_appendix:L_estimate of products} with $L=L^\infty$ and $s>0$ large enough to have $N(s)\ge 3$, and a-priori estimates \eqref{est: bootstrap argument a-priori est}, we can reduce ourselves to estimate the $L^\infty$ norm of the product in the left hand side of \eqref{est_Linfty_prod_a0v_b1u_R1u} when all its factors are supported for moderate frequencies less or equal than $t^\sigma$, up to remainders $R(t,x)$.
Moreover, since
\begin{subequations} \label{est:vNF-v-_uBF-u-_appendix}
\begin{equation} \label{a0_vNF-v-}
\left\|\chi(t^{-\sigma}D_x)a_0(D_x)[v^{NF}-v_{-}](t,\cdot) \right\|_{L^\infty}\le CA^2\varepsilon^2 t^{-\frac{3}{2}+\sigma}
\end{equation}
and 
\begin{equation} \label{b1_uNF-u-}
\left\|\chi(t^{-\sigma}D_x)b_1(D_x)[u^{NF}-u_{-}](t,\cdot) \right\|_{L^\infty}\le CA^2\varepsilon^2 t^{-2+\beta},
\end{equation}
\end{subequations}
as follows by \eqref{explicit vNf-v-} and \eqref{def uNF}, \eqref{est: Linfty integral D with cut-off} with $\rho=2$ (as consequence of lemma \ref{Lem_Appendix: est on Dj1j2}), together with a-priori estimates, we can also suppose $v_{-}$ (resp. $u_{-}$) be replaced with $v^{NF}$ (resp. $u^{NF}$), up to some new $R(t,x)$.
This reduces us to prove that
\begin{multline*}
\left\|[\chi(t^{-\sigma}D_x)a_0(D_x)v^{NF}][\chi(t^{-\sigma}D_x)b_1(D_x)u^{NF}][\chi(t^{-\sigma}D_x)b_0(D_x) u_{-}](t,\cdot)\right\|_{L^\infty}\\
\le C(A+B)AB\varepsilon^3 t^{-\frac{5}{2}+\beta+\frac{\delta+\delta_1}{2}},
\end{multline*}
or rather, reminding \eqref{est: bootstrap upm}, to show that
\begin{equation*}
\left\|[\chi(t^{-\sigma}D_x)a_0(D_x)v^{NF}][\chi(t^{-\sigma}D_x)b_1(D_x)u^{NF}](t,\cdot)\right\|_{L^\infty}\le C(A+B)B\varepsilon^2 t^{-2+\beta+\frac{\delta+\delta_1}{2}}.
\end{equation*}
But after writing the above product in the semi-classical setting and reminding definition \eqref{def utilde vtilde}, one can immediately check that this estimate is satisfied thanks to \eqref{est_Linfty_Vtilde Utilde}, which concludes the proof of \eqref{est_Linfty_prod_a0v_b1u_R1u}.

The last part of the statement follows from \eqref{explicit_rnfkg_chapter5}, the fact that
\begin{equation*}
\left\|\chi(t^{-\sigma}D_x)\left[ - \frac{D_x}{\langle D_x\rangle}(v_+ - v_{-})\, \Nlw + D_1\big[\langle D_x\rangle^{-1}(v_+-v_{-})\, \Nlw\right](t,\cdot) \right\|_{L^\infty}\le CA^3\varepsilon^3 t^{-3+\sigma}
\end{equation*}
for every $t\in [1,T]$, which is consequence of \eqref{est Linfty NLw} and a-priori estimate \eqref{est: boostrap vpm},
and from the observation that the remaining contributions to $\rnfkg$ are products of the form
\begin{equation*}
[a_0(D_x)v_{-}] [b_1(D_x)u_{-}] \mathrm{R}_1u_{-},
\end{equation*}
with $a_0(\xi)$ equal to 1 or to $\xi_j\langle\xi\rangle^{-1}$, and $b_1(\xi)$ equal to $\xi_1$ or to $\xi_j \xi_1|\xi|^{-1}$, for $j=1,2$.
\endproof
\end{lem}

\begin{lem} \label{Lem_appendix: L xnrNFkg}
Under the same assumptions as in lemma \ref{Lem_appendix:Linfty_est_rNFkg},
\begin{subequations}\label{est_xx_prod_a0v_b1u_R1u}
\begin{align}
\left\|\chi(t^{-\sigma}D_x) \big[x_n [a_0(D_x)v_{-}] [b_1(D_x)u_{-}] b_0(D_x) u_{-}\big](t,\cdot)\right\|_{L^2(dx)}&\le C(A+B)^2B\varepsilon^3 t^{-1+\beta+\frac{\delta_1}{2}}, \\
\left\|\chi(t^{-\sigma}D_x) \big[x_m x_n [a_0(D_x)v_{-}][ b_1(D_x)u_{-}]b_0(D_x)u_{-}\big](t,\cdot)\right\|_{L^2(dx)}&\le  C(A+B)^2B\varepsilon^3 t^{\beta+\frac{\delta_1}{2}},
\end{align}
\end{subequations}
for every $t\in [1,T]$, $m,n=1,2$, with $\beta>0$ small, $\beta\rightarrow 0$ as $\sigma\rightarrow 0$.
Moreover,
\begin{subequations} \label{est_xxrnfkg}
\begin{align}
\left\|\chi(t^{-\sigma}D_x) \big[x_n r^{NF}_{kg}(t,\cdot)\big]\right\|_{L^2(dx)} &\le C(A+B)^2B\varepsilon^3 t^{-1+\beta+\frac{\delta+\delta_1}{2}}, \label{est:t_xn_rNFkg}\\
\left\|\chi(t^{-\sigma}D_x) \big[x_m x_n r^{NF}_{kg}(t,\cdot)\big]\right\|_{L^2(dx)} &\le C(A+B)^2B\varepsilon^3 t^{\beta+\frac{\delta+\delta_1}{2}}. \label{est:xm_xn_rNFkg}
\end{align}
\end{subequations}
\proof
We warn the reader that we will denote by $C$ and $\beta$ two positive constants that may change line after line, with $\beta\rightarrow 0$ as $\sigma\rightarrow 0$. We also denote by $R(t,x)$ any contribution verifying
\begin{subequations} \label{est_R_lemma_rnfkg}
\begin{align}
\left\|\chi(t^{-\sigma}D_x) \big[x_nR(t,\cdot)\big]\right\|_{L^2(dx)} &\le C(A+B)^2B\varepsilon^3 t^{-1+\beta+\frac{\delta+\delta_1}{2}}, \label{xn_R} \\
\left\|\chi(t^{-\sigma}D_x) \big[x_m x_n R(t,\cdot)\big]\right\|_{L^2(dx)} &\le C(A+B)^2B\varepsilon^3 t^{\beta+\frac{\delta+\delta_1}{2}}.\label{xnxm_R}
\end{align}
\end{subequations}
Let us first notice that, after \eqref{est Linfty NLw}, \eqref{norm_xv-}, \eqref{norm_L2_xixjv-} and a-priori estimates, we have that
\begin{multline*}
\left\|\chi(t^{-\sigma}D_x)\left[ - x_n\frac{D_x}{\langle D_x\rangle}(v_+ - v_{-})\, \Nlw +x_n D_1\big[\langle D_x\rangle^{-1}(v_+-v_{-})\, \Nlw\right](t,\cdot) \right\|_{L^2} \\
\lesssim t^\sigma \sum_{\mu=0}^1\left\|x_n^\mu v_\pm(t,\cdot) \right\|_{L^2}\|\textit{NL}_w(t,\cdot)\|_{L^\infty}\le CA^2B\varepsilon^3 t^{-1+\sigma+\frac{\delta}{2}}
\end{multline*}
and
\begin{multline*}
\left\|\chi(t^{-\sigma}D_x)\left[ - x_mx_n\frac{D_x}{\langle D_x\rangle}(v_+ - v_{-})\, \Nlw + x_mx_nD_1\big[\langle D_x\rangle^{-1}(v_+-v_{-})\, \Nlw\right](t,\cdot) \right\|_{L^2} \\
\lesssim t^\sigma \sum_{\mu_1, \mu_2=0}^1\left\| x_m^{\mu_1}x_n^{\mu_2}  v_\pm(t,\cdot) \right\|_{L^2}\|\textit{NL}_w(t,\cdot)\|_{L^\infty}\le C(A+B)AB\varepsilon^3 t^{\sigma+\frac{\delta}{2}}.
\end{multline*}
Therefore, since from \eqref{explicit_rnfkg_chapter5} and \eqref{def_app_Nlkg} the remaining contributions to $\rnfkg$ are of the form
\begin{equation*}
[a_0(D_x)v_{-}] [b_1(D_x)u_{-}] \mathrm{R}_1u_{-}
\end{equation*}
with $a_0(\xi)$ equal to 1 or to $\xi_j\langle\xi\rangle^{-1}$, and $b_1(\xi)$ equal to $\xi_1$ or to $\xi_j \xi_1|\xi|^{-1}$, for $j=1,2$, estimates \eqref{est_xxrnfkg} will follow from \eqref{est_xx_prod_a0v_b1u_R1u}.
Our aim is hence to prove that the above product is a remainder $R(t,x)$.

Applying lemma \ref{Lem_appendix:L_estimate of products} with $L=L^2$, $w_1=x_n a_0(D_x) v_{-}$ (resp. $w_1=x_m x_n a_0(D_x) v_{-}$), $s>0$ sufficiently large so that $N(s)> 2$, and using estimates \eqref{norm_xv-} (resp. \eqref{norm_L2_xixjv-}), \eqref{est: bootstrap upm}, \eqref{est: bootstrap Enn}, we can suppose all above factors truncated for moderate frequencies less or equal than $t^\sigma$, up to remainders $R(t,x)$.
Let us also observe that, from \eqref{norm_Linfty_xv-}, \eqref{b1_uNF-u-} and \eqref{est: bootstrap Enn},\small
\begin{align*}
\left\|\chi(t^{-\sigma}D_x) \Big[ [\chi_1(t^{-\sigma}D_x)[x_n a_0(D_x)v_{-}]][ \chi(t^{-\sigma}D_x) b_1(D_x)(u^{NF}- u_{-})] [\chi(t^{-\sigma}D_x)b_0(D_x)u_{-}]\Big](t,\cdot)(t,\cdot) \right\|_{L^2} \\
\lesssim \sum_{|\mu|=0}^1 \left\|x_n\Big(\frac{D_x}{\langle D_x\rangle}\Big)^\mu v_\pm(t,\cdot)\right\|_{L^\infty}\| \chi(t^{-\sigma}D_x) b_1(D_x)(u^{NF}- u_{-})\|_{L^\infty}\| u_\pm(t,\cdot)\|_{L^2}\\
 \le C(A+B)A^2B\varepsilon^4 t^{-2+\beta+\frac{\delta+\delta_2}{2}},
\end{align*}\normalsize
and that, using additionally estimate \eqref{norm_Linfty_xixjv-}, \small
\begin{align*}
 \left\|\chi(t^{-\sigma}D_x) \Big[ [\chi_1(t^{-\sigma}D_x)[x_m x_n a_0(D_x)v_{-}][\chi(t^{-\sigma}D_x) b_1(D_x)(u^{NF}- u_{-})] [\chi(t^{-\sigma}D_x)b_0(D_x)u_{-}]\Big](t,\cdot)(t,\cdot) \right\|_{L^2} \\
\lesssim \sum_{|\mu|, |\nu| =0}^1 \left\|x_m x_n\Big(\frac{D_x}{\langle D_x\rangle}\Big)^\mu v_\pm(t,\cdot)\right\|_{L^\infty}\| \chi(t^{-\sigma}D_x) b_1(D_x)(u^{NF}- u_{-})\|_{L^2}\|\mathrm{R}^\nu u_\pm(t,\cdot)\|_{L^\infty}\\
 \le C(A+B)A^2B\varepsilon^4 t^{-1+\beta+\frac{\delta+\delta_2}{2}}.
\end{align*}\normalsize
This means that we can actually replace $u_{-}$ by $u^{NF}$ up to some new $R(t,x)$.
Furthermore, we can also substitute $\chi_1(t^{-\sigma}D_x)[x^k_m x_na_0(D_x)v_{-}]$ with $\chi(t^{-\sigma}D_x)[x^k_m x_n a_0(D_x)v^{NF}_{-}]$, for any $k\in \{0,1\}$, up to a new remainder $R(t,x)$ in consequence of a-priori estimate \eqref{est: bootstrap upm}, the fact that
\begin{equation} \label{est_L2_uNF}
\|u^{NF}(t,\cdot)\|_{L^2}\le CB\varepsilon t^\frac{\delta}{2},
\end{equation}
(see \eqref{est:utilde-Hs} in semi-classical coordinates), and the following inequalities
\begin{subequations}\label{xx_(vnf-v-)}
\begin{multline} \label{x_(vnf-v)}
\left\|\chi_1(t^{-\sigma}D_x)\left[x_na_0(D_x) (v^{NF}-v_{-})\right](t,\cdot)\right\|_{L^\infty}\\
\lesssim \sum_{\mu,\nu,\kappa=0}^1t^\sigma\left\|x^\mu_n \Big(\frac{D_x}{\langle D_x\rangle}\Big)^\nu v_\pm(t,\cdot) \right\|_{L^\infty} \|\mathrm{R}^\kappa_1u_\pm(t,\cdot)\|_{L^\infty} 
\le C(A+B)A\varepsilon^2 t^{-\frac{1}{2}+\sigma+\frac{\delta_2}{2}}
\end{multline}
and
\begin{multline} 
\left\| \chi_1(t^{-\sigma}D_x)\left[x_m x_n a_0(D_x) (v^{NF}-v_{-})\right](t,\cdot)\right\|_{L^\infty}\\ 
\lesssim \sum_{\mu_1,\mu_2,\nu,\kappa=0}^1\left\|x^{\mu_1}_m x^{\mu_2}_n \Big(\frac{D_x}{\langle D_x\rangle}\Big)^\nu v_\pm(t,\cdot) \right\|_{L^\infty} \|\mathrm{R}^\kappa_1u_\pm(t,\cdot)\|_{L^\infty}\le C(A+B)A\varepsilon^2 t^{\frac{1}{2}+\frac{\delta_2}{2}},
\end{multline}
\end{subequations}
derived from \eqref{explicit_vNF-v_chapter5}, \eqref{norm_Linfty_xv-}, \eqref{norm_Linfty_xixjv-}, \eqref{est: bootstrap upm} and \eqref{est: boostrap vpm}.
This reduces us to prove that, for $k=0,1$,
\begin{multline*}
\left\|\big[\chi_1(t^{-\sigma}D_x)[x_m^k x_n a_0(D_x)v^{NF}_{-}]\big] \big[\chi(t^{-\sigma}D_x)b_1(D_x)u^{NF}\big]\big[\chi(t^{-\sigma}D_x)b_0(D_x)u_{-}\big](t,\cdot)\right\|_{L^2(dx)}\\
\le C(A+B)^2B\varepsilon^3 t^{-1+k+\beta+\frac{\delta_1}{2}},
\end{multline*}
or rather, after \eqref{est: bootstrap upm}, that \small
\begin{equation*}
\left\| \left[\chi_1(t^{-\sigma}D_x)[x_m^k x_n a_0(D_x) v^{NF}_{-}]\right] [\chi(t^{-\sigma}D_x)b_1(D_x) u^{NF}](t,\cdot) \right\|_{L^2(dx)}\le C(A+B)B\varepsilon^2 t^{-\frac{1}{2}+k +\beta+\frac{\delta+\delta_1}{2}}.
\end{equation*}\normalsize
Passing to the semi-classical setting, this corresponds to prove that\small
\begin{equation} \label{est_sc_setting_appendix}
\sum_{k=0}^1\Big\| \left[\oph(\chi_1(h^\sigma\xi))[x_m^k x_n \oph(a_0(\xi))\widetilde{v}\right] [\oph(\chi(h^\sigma\xi)b_1(\xi))\widetilde{u}] (t,\cdot)\Big\|_{L^2(dx)}\le C(A+B)B\varepsilon^2 h^{\frac{1}{2}-\beta-\frac{\delta+ \delta_1}{2}}.
\end{equation}\normalsize
First of all let us notice that, from the commutation of $x_n$ with $\oph(a_0(\xi))$ and definition \eqref{def Lj} of $\mathcal{L}_n$, 
\begin{equation} \label{eq:xn_a0_vtilde}
\begin{split}
x_n \oph(a_0(\xi))\widetilde{v} & = h \oph(a_0(\xi))\mathcal{L}_n\widetilde{v} + \oph\Big(a_0(\xi)\frac{\xi_n}{\langle\xi\rangle}\Big)\widetilde{v} - \frac{h}{2i}\oph\big(\partial_{\xi_n}a_0(\xi)\big)\widetilde{v},
\end{split}
\end{equation}
while from the commutation of $x_m$ with $\oph(\chi(h^\sigma\xi)b_1(\xi))$, definition \eqref{def Mj} of $\mathcal{M}_m$, and symbolic calculus,
\begin{equation} \label{eq:xm_b1_utilde}
\begin{split}
x_m \oph(\chi(h^\sigma\xi)b_1(\xi))\widetilde{u}& = h \oph(\chi(h^\sigma\xi)b_1(\xi)|\xi|^{-1})\mathcal{M}_m\widetilde{u} - \frac{h}{2i}\oph\big(\partial_{\xi_m}(\chi(h^\sigma\xi)b_1(\xi)|\xi|^{-1})|\xi|\big)\widetilde{u}\\
& + \oph(\chi(h^\sigma\xi)b_1(\xi)\xi_m|\xi|^{-1})\widetilde{u} -\frac{h}{2i}\oph\big(\partial_{\xi_m}(\chi(h^\sigma\xi)b_1(\xi))\big)\widetilde{u}.
\end{split}
\end{equation}
On the one hand, using equality \eqref{eq:xn_a0_vtilde}, lemma \ref{Lem: from energy to norms in sc coordinates-KG}, and estimates \eqref{est_L2_Vtilde_Utilde}, \eqref{Hrho_infty_utilde_appendix}, we deduce that
\begin{equation} \label{ineq:xn_utilde_vtilde}
\begin{split}
&\Big\|\left[\oph(\chi_1(h^\sigma\xi)[x_n \oph(a_0(\xi))\widetilde{v}]\right] [\oph(\chi(h^\sigma\xi) b_1(\xi))\widetilde{u}](t,\cdot) \Big\|_{L^2} \\
&\le  \left\| \Big[\oph\Big(\chi_1(h^\sigma\xi)a_0(\xi)\frac{\xi_n}{\langle\xi\rangle}\Big)\widetilde{v}(t,\cdot)\Big] [\oph(\chi(h^\sigma\xi)b_1(\xi))\widetilde{u}(t,\cdot)] \right\|_{L^2} + CAB\varepsilon^2 h^{\frac{1}{2}-\beta-\frac{\delta_2}{2}}\\
& \le C(A+B)B\varepsilon^2 h^{\frac{1}{2}-\beta-\frac{\delta+\delta_1}{2}}.
\end{split}
\end{equation}
On the other hand, when we deal with the $L^2$ norm in the left hand side of \eqref{est_sc_setting_appendix} corresponding to $k=1$ we first commute $x_m$ with $\oph(\chi_1(h^\sigma\xi))$ and see, using symbolic calculus, that
\begin{align*}
&\Big\| \left[\oph(\chi_1(h^\sigma\xi))[x_m x_n \oph(a_0(\xi))\widetilde{v}\right] [\oph(\chi(h^\sigma\xi)b_1(\xi))\widetilde{u}] (t,\cdot)\Big\|_{L^2(dx)} \\
&\le \Big\| \left[h^\sigma \oph((\partial\chi_1)(h^\sigma\xi))[ x_n \oph(a_0(\xi))\widetilde{v}\right] [\oph(\chi(h^\sigma\xi)b_1(\xi))\widetilde{u}] (t,\cdot)\Big\|_{L^2(dx)} \\
& + \Big\| \left[\oph(\chi_1(h^\sigma\xi))[ x_n \oph(a_0(\xi))\widetilde{v}\right] [x_m\oph(\chi(h^\sigma\xi)b_1(\xi))\widetilde{u}](t,\cdot) \Big\|_{L^2(dx)}.
\end{align*}
The first norm in the above right hand side satisfies an inequality analogous to \eqref{ineq:xn_utilde_vtilde}.
In order to derive an estimate for the latter one, we first use equality \eqref{eq:xn_a0_vtilde} and observe the following: from the semi-classical Sobolev injection and estimates \eqref{est:Mutilde}, \eqref{est:Lvtilde}, we have that
\begin{multline}\label{est_Lvt_Mut}
\left\| h^2\big[\oph(\chi_1(h^\sigma\xi)a_0(\xi))\mathcal{L}_n\widetilde{v}\big]\big[ \oph(\chi(h^\sigma\xi)b_1(\xi)\xi_m|\xi|^{-1})\mathcal{M}_m\widetilde{u}\big](t,\cdot)\right\|_{L^2}\\
\lesssim h \left\| \oph(\chi_1(h^\sigma\xi)a_0(\xi))\mathcal{L}_n\widetilde{v}(t,\cdot)\right\|_{L^2}\left\| \oph(\chi(h^\sigma\xi)b_1(\xi)\xi_m|\xi|^{-1})\mathcal{M}_m\widetilde{u}(t,\cdot)\right\|_{L^2}\\
\le C(A+B)B\varepsilon^2 h^{1-\delta_2-\beta};
\end{multline}
a similar chain of inequalities as in \eqref{injection_R1w}, together with \eqref{Hs norm uNF- u-}, \eqref{Hrho-infty norm uNF- u-} and \eqref{est: bootstrap argument a-priori est}, gives that for any $\theta\in ]0,1[$ 
\begin{multline}\label{RRutillde}
\left\| \oph(b_1(\xi)\xi_m|\xi|^{-1})\widetilde{u}(t,\cdot)\right\|_{L^\infty}=t \left\|b_1(D_x)D_m|D_x|^{-1}u^{NF}(t,\cdot)\right\|_{L^\infty}\\ \lesssim t \|u^{NF}(t,\cdot)\|^{1-\theta}_{H^{3,\infty}}\|u^{NF}(t,\cdot)\|_{H^2} 
\le CA^{1-\theta}B^\theta\varepsilon t^{\frac{1}{2}+\frac{(1+\delta)}{2}\theta}.
\end{multline}
Therefore, from equality \eqref{eq:xm_b1_utilde} and estimates \eqref{est:Lvtilde}, \eqref{Hrho_infty_utilde_appendix}, \eqref{RRutillde}, \eqref{est_Lvt_Mut}, we find that
\begin{equation} \label{est_Lv_Mu}
h\left\|\oph(\chi_1(h^\sigma\xi)a_0(\xi))\mathcal{L}_n\widetilde{v}\, \big[x_m \oph(\chi(h^\sigma\xi)b_1(\xi))\widetilde{u}\big](t,\cdot)\right\|_{L^2}\le C(A+B)B\varepsilon^2 h^{\frac{1}{2}-\frac{\delta_2}{2} -\frac{(1+\delta)\theta}{2}}.
\end{equation}
Moreover, using again \eqref{eq:xm_b1_utilde} along with\eqref{est:utilde-Hs}, \eqref{est:Mutilde}, \eqref{est_L2_Vtilde_Utilde} and \eqref{Hrho_infty_vtilde_appendix},
\begin{align*}
&\left\|\Big[\oph\Big(\chi_1(h^\sigma\xi)a_0(\xi)\frac{\xi_n}{\langle\xi\rangle}\Big)\vt + h\oph\big(\chi_1(h^\sigma\xi) \partial_{\xi_n}a_0(\xi)\big)\vt\Big] \, \big[x_m \oph(\chi(h^\sigma\xi)b_1(\xi))\widetilde{u}\big](t,\cdot)\right\|_{L^2(dx)}\\
& \le \left\|\Big[\oph\Big(\chi_1(h^\sigma\xi)a_0(\xi)\frac{\xi_n}{\langle\xi\rangle}\Big)\vt\Big] \Big[\oph\Big(\chi(h^\sigma\xi)b_1(\xi)\frac{\xi_m}{|\xi|}\Big)\ut\Big](t,\cdot)\right\|_{L^2(dx)} + C(A+B)B\varepsilon^2 h^{1-\beta-\frac{\delta_2}{2}}\\
& \le C(A+B)B \varepsilon^2 h^{\frac{1}{2}-\beta-\frac{\delta+\delta_1}{2}}.
\end{align*}
Choosing $\theta\ll 1$ small enough, this concludes that
\begin{equation} \label{ineq:xmxn_vtilde_utilde}
\Big\| \left[\oph(\chi_1(h^\sigma\xi))[x_m x_n \oph(a_0(\xi))\widetilde{v}\right] [\oph(\chi(h^\sigma\xi)b_1(\xi))\widetilde{u}] (t,\cdot)\Big\|_{L^2(dx)} \le C(A+B)B\varepsilon^2 h^{\frac{1}{2}-\beta-\frac{\delta+\delta_1}{2}}
\end{equation}
and, together with \eqref{ineq:xn_utilde_vtilde}, the proof of \eqref{est_sc_setting_appendix}.
\endproof
\end{lem}

We can finally prove the following:
\begin{lem} \label{Lem_appendix: estimate L2vtilde}
There exists a constant $C>0$ such that, for any $\chi \in C^\infty_0(\mathbb{R}^2)$, $\sigma>0$ small, and every $t\in [1,T]$,
\begin{equation} \label{est:L2vtilde}
\sum_{|\mu|=2}\|\oph(\chi(h^\sigma\xi))\mathcal{L}^\mu\widetilde{v}(t,\cdot)\|_{L^2} \le CB\varepsilon t^{\beta+\frac{\delta+\delta_1}{2}},
\end{equation}
with $\beta>0$ small, $\beta\rightarrow 0$ as $\sigma\rightarrow 0$.
\proof
From relation \eqref{relation between Zjv and Lj vtilde} and the commutation between $\mathcal{L}_m$ and $\oph(\langle\xi\rangle)$ we deduce that
\begin{equation} \label{ineq_L2vtilde}
\begin{split}
&\left\|\oph(\chi(h^\sigma\xi))\mathcal{L}_m\mathcal{L}_n \widetilde{v}(t,\cdot)\right\|_{H^1_h}\lesssim \sum_{\mu=0}^1\Big[\left\|\oph(\chi(h^\sigma\xi))\mathcal{L}^\mu_m \big[tZ_n v^{NF}(t,tx)\big] \right\|_{L^2(dx)}\\
& + \left\|\oph(\chi(h^\sigma\xi))\mathcal{L}^\mu_m \oph\Big(\frac{\xi_n}{\langle\xi\rangle}\Big)\widetilde{v}(t,\cdot)\big]\right\|_{L^2(dx)}+ \left\|\oph(\chi(h^\sigma\xi))\mathcal{L}^\mu_m \big[t(tx_n)r^{NF}_{kg}(t,tx)\big]\right\|_{L^2(dx)}\Big],
\end{split}
\end{equation}
so the result of the statement follows from lemmas \ref{Lem: from energy to norms in sc coordinates-KG}, \ref{Lem_appendix: ZnvNF_LmZnvNF}, and inequalities \eqref{dec_norm_Lm}, \eqref{est_xxrnfkg}.
\endproof
\end{lem}

\section{The sharp decay estimate of the Klein-Gordon solution with a Klainerman vector field}\label{Sec_App_4}

This last section is devoted to prove that, for any admissible vector field $\Gamma$, the $L^\infty(\mathbb{R}^2)$ norm of functions $(\Gamma v)_\pm$, when restricted to moderate frequencies less or equal than $t^\sigma$, for some small $\sigma>0$, decays in time at the same sharp rate $t^{-1}$ of the two-dimensional linear Klein-Gordon solution.
This result is proved in lemma \ref{Lem_appendix: sharp_est_VJ} under the hypothesis that a-priori estimates \eqref{est: bootstrap argument a-priori est} are satisfied in some fixed interval $[1,T]$, with $\varepsilon_0<(2A+B)^{-1}$ and $0<\delta\ll \delta_2\ll\delta_1\ll \delta_0\ll 1$ sufficiently small, and is fundamental when proving lemmas \ref{Lem: L2 est nonlinearities} and \ref{Lem:L2 est nonlinearity Dt}.
All the other lemmas of this section are to be meant as preparatory intermediate results.

\begin{lem}
With the convention that $D=D_1$ whenever $|I_1|+|I_2|=2$, $D\in \{D_j, D_t, j=1,2\}$ otherwise, there exists a positive constant $C>0$ such that, for any $\chi\in C^\infty_0(\mathbb{R}^2)$, $\sigma>0$ small, $n=1,2$, and every $t\in [1,T]$,
\begin{equation} \label{xn_Qkg0_vI1uI2}
\sum_{\substack{|I_1|+|I_2|\le 2 \\ |I_1|<2}}\left\| \chi(t^{-\sigma}D_x)\big[x_n Q^\mathrm{kg}_0(v^{I_1}_\pm, D u^{I_2}_\pm)\big](t,\cdot)\right\|_{L^2(dx)}\le C(A+B)B\varepsilon^2 t^{\beta+\frac{\delta+\delta_2}{2}},
\end{equation}
with $\beta>0$ small such that $\beta\rightarrow 0$ as $\sigma\rightarrow 0$.
\proof
We estimate the $L^2$ norms in the left hand side of \eqref{xn_Qkg0_vI1uI2} separately.

$\bullet$ When $|I_1|=0$, $|I_2|=2$, we derive from \eqref{norm_Linfty_xv-} and \eqref{est: bootstrap E02} that
\begin{equation*}
\begin{split}
\left\| \chi(t^{-\sigma}D_x)\big[x_n Q^\mathrm{kg}_0(v_\pm, D_1 u^{I_2}_\pm)\big](t,\cdot)\right\|_{L^2(dx)} &\lesssim \sum_{|\mu|=0}^1\left\| x_n \Big(\frac{D_x}{\langle D_x\rangle}\Big)^\mu v_\pm(t,\cdot)\right\|_{L^\infty}\|u^{I_2}_\pm(t,\cdot)\|_{H^1}\\
&\le C(A+B)B\varepsilon t^\frac{\delta_1+\delta_2}{2};
\end{split}
\end{equation*}

$\bullet$ When $|I_1|=|I_2|=1$ and $\Gamma^{I_2}\in \{\Omega,Z_m, m=1,2\}$ is a Klainerman vector field we use inequalities \eqref{cor_app_est_2} with $L=L^2$, $w_{j_0} =x_n(D_x\langle D_x\rangle^{-1})^\mu v^{I_1}_\pm$ with $|\mu|=0,1,$ and $s$ large enough so that $N(s)\ge 2$, to derive that
\begin{equation*}
\begin{split}
& \left\|\chi(t^{-\sigma}D_x)\left[x_n Q^\mathrm{kg}_0(v^{I_1}_\pm, D_1 u^{I_2}_\pm)\right](t,\cdot)\right\|_{L^2(dx)}\\
& \lesssim \sum_{|\mu|=0}^1\left\|\chi_1(t^{-\sigma}D_x)\Big[ x_n \Big(\frac{D_x}{\langle D_x\rangle}\Big)^\mu v^{I_1}_\pm\Big](t,\cdot)\right\|_{L^\infty}\|u^{I_2}_\pm(t,\cdot)\|_{H^1}\\
& +\sum_{\substack{|\mu|=0,1,2 \\|\nu|=0,1 }}t^{-2}\left(\|x^\mu v^{I_1}_\pm(t,\cdot)\|_{L^2} + t\|x^\nu v^{I_1}_\pm(t,\cdot)\|_{L^2}\right) \left(\|u_\pm(t,\cdot)\|_{H^s} + \|D_tu_\pm(t,\cdot)\|_{H^s}\right)  \\
&\le CB^2\varepsilon^2 t^\frac{\delta_1+\delta_2}{2},
\end{split}
\end{equation*}
where last estimate is deduced using \eqref{Hs_norm_DtU}, \eqref{norm_L2_xj-GammaIv-}, \eqref{est:xixjGamma v-}, \eqref{norm_Linfty_xjGammav-} and \eqref{est: bootstrap E02};

$\bullet$ When $|I_1|=|I_2|=1$ and $\Gamma^{I_2}$ is a spatial derivative we use lemma \ref{Lem_appendix:L_estimate of products} with $L=L^2$, $w_1 =x_n(D_x\langle D_x\rangle^{-1})^\mu v^{I_1}_\pm$ with $|\mu|=0,1$, $s$ large enough so that $N(s)\ge 1$, and again estimates \eqref{norm_L2_xj-GammaIv-}, \eqref{norm_Linfty_xjGammav-} and \eqref{est: bootstrap E02}. We obtain that
\begin{equation*}
\begin{split}
& \left\|\chi(t^{-\sigma}D_x)\left[x_n Q^\mathrm{kg}_0(v^{I_1}_\pm, D_1 u^{I_2}_\pm)\right](t,\cdot)\right\|_{L^2(dx)}\\
& \lesssim \sum_{|\mu|=0}^1\left\|\chi_1(t^{-\sigma}D_x)\Big[ x_n \Big(\frac{D_x}{\langle D_x\rangle}\Big)^\mu v^{I_1}_\pm\Big](t,\cdot)\right\|_{L^\infty}\|u_\pm(t,\cdot)\|_{H^2}\\
& +\sum_{|\mu|=0}^1t^{-1}\|x^\mu v^{I_1}_\pm(t,\cdot)\|_{L^2} \|u_\pm(t,\cdot)\|_{H^s} \le CB^2\varepsilon^2 t^\frac{\delta+\delta_1}{2};
\end{split}
\end{equation*}

$\bullet$ When $|I_1|+|I_2|\le 1$, by the assumption derivative $D$ can be equal to $D_x$ or to $D_t$. Then

\hspace{6pt} - If $|I_1|=0$, after \eqref{Hs_norm_DtU}, \eqref{DtUI}, \eqref{norm_Linfty_xv-} and \eqref{est: bootstrap argument a-priori est}
\begin{equation*}
\begin{split}
&\left\|\chi(t^{-\sigma}D_x)\big[x_n Q^\mathrm{kg}_0(v_\pm,  D u^{I_2}_\pm)\big](t,\cdot) \right\|_{L^2(dx)}\\
& \lesssim \sum_{|\mu|=0}^1 \left\| x_n\Big(\frac{D_x}{\langle D_x\rangle}\Big)^\mu v_\pm(t,\cdot)\right\|_{L^\infty}\left(\|u^{I_2}_\pm(t,\cdot)\|_{H^1} + \|D_tu^{I_2}_\pm(t,\cdot)\|_{L^2}\right)  \le C(A+B)B\varepsilon^2 t^{\delta_2};
\end{split}
\end{equation*}

\hspace{6pt} - If $|I_1|=1$, $|I_2|=0$, using lemma \ref{Lem_appendix:L_estimate of products} as done above, together with \eqref{Hs_norm_DtU}, \eqref{norms_H1_Linfty_xv-}, \eqref{norm_L2_xj-GammaIv-}, \eqref{norm_Linfty_xjGammav-} and a-priori estimates, we derive that
\begin{equation}
\begin{split}
&\left\|\chi(t^{-\sigma}D_x)\big[x_n Q^\mathrm{kg}_0(v^{I_1}_\pm,  D u_\pm)\big] (t,\cdot)\right\|_{L^2(dx)} \\
&\lesssim \sum_{|\mu|=0}^1 \left\|\chi_1(t^{-\sigma}D_x)\Big[ x_n\Big(\frac{D_x}{\langle D_x\rangle}\Big)^\mu v^{I_1}_\pm(t,\cdot)\Big]\right\|_{L^\infty}\left(\|u_\pm(t,\cdot)\|_{H^1} + \|D_tu_\pm(t,\cdot)\|_{L^2}\right)\\
& + \sum_{|\mu|=0}^1t^{-1} \left\| x_n\Big(\frac{D_x}{\langle D_x\rangle}\Big)^\mu v^{I_1}_\pm(t,\cdot)\right\|_{L^2}\left(\|u_\pm(t,\cdot)\|_{H^s} + \|D_tu_\pm(t,\cdot)\|_{H^s}\right) \\
& \le CB^2\varepsilon^2 t^{\beta+\frac{\delta+\delta_1}{2}}.
\end{split}
\end{equation}
\endproof
\end{lem}

\begin{lem} \label{Lem_appendix:est vI I=2}
There exists a positive constant $C>0$ such that, for any $\chi\in C^\infty_0(\mathbb{R}^2)$, $\sigma>0$ small, $\rho\in\mathbb{N}$, and every $t\in [1,T]$,
\begin{equation} \label{est:Linfty_vI_I=2-new}
\sum_{|I|=2}\left\|\chi(t^{-\sigma}D_x)V^I(t,\cdot)\right\|_{H^{\rho,\infty}}\le C B\varepsilon t^{-1+\beta+\frac{\delta_0}{2}},
\end{equation}
with $\beta>0$ small such that $\beta\rightarrow 0$ as $\sigma\rightarrow 0$.
\proof
Estimate \eqref{est:Linfty_vI_I=2-new} is evidently satisfied when $\Gamma^I$ contains at least one spatial derivative thanks to lemma \ref{Lem_appendix: preliminary est VJ}.
We then focus on the case when $\Gamma^I$ is the product of two Klainerman vector fields.
As $v^I_+ = -\overline{v^I_{-}}$, we prove the statement for $\chi(t^{-\sigma}D_x)v^I_{-}$. Moreover, from the $L^\infty-H^{\rho,\infty}$ continuity of $\chi(t^{-\sigma}D_x)$ with norm $O(t^{\sigma\rho})$, we can assume the $H^{\rho,\infty}$ norm in \eqref{est:Linfty_vI_I=2-new} replaced with the $L^\infty$ one, up to a loss $t^{\sigma\rho}$.

As done in lemma \ref{Lem_appendix: preliminary est VJ}, instead of proving the statement directly on $\chi(t^{-\sigma}D_x)v^I_{-}$ we do it for $\chi(t^{-\sigma}D_x)\vNFGamma$, with $\vNFGamma$ introduced in \eqref{def_vNF-Gamma} and considered here with $|I|=2$. This is justified by inequality \eqref{vI-_bounded_by_vNFGamma}.
From definition \eqref{def_vNF-Gamma} of $\vNFGamma$, equation \eqref{half KG Gammav} and equality \eqref{GammaI_NLkg} one can check that
\begin{equation} \label{def_NLNF_I=2}
\begin{gathered}[]
[D_t + \langle D_x\rangle] \vNFGamma = \NLNF\\
\text{where } \NLNF = r^{I,\textit{NF}}_{kg}(t,x) + \sum_{\substack{(I_1,I_2)\in\mathcal{I}(I)\\ |I_1|<2}} c_{I_1,I_2}Q^\mathrm{kg}_0(v^{I_1}_\pm, D u^{I_2}_\pm),
\end{gathered}
\end{equation}
with $r^{I,\textit{NF}}_{kg}$ given by the same integral expression as in \eqref{def:rNF-Gamma-kg} but with $|I|=2$ (and hence having the explicit expression \eqref{explicit rNFkg-Gamma}), and $c_{I_1,I_2}\in \{-1,0,1\}$, $c_{I_1,I_2}=1$ when $|I_1|+|I_2|=2$ (in which case derivative $D$ corresponds to $D_1$).
It is straightforward to show that inequalities \eqref{prelimary_ineq_vNFGamma}, \eqref{est_preliminary_(vNFGamma - vJ-)}, \eqref{est:vNFGamma} hold even when $|I|=2$, up to replacing $\delta_2$ with $\delta_1$.
Therefore, using those latter ones together with
\[\sum_{j=1}^2 \left\|\chi(t^{-\sigma}D_x)Z_j\vNFGamma(t,\cdot)\right\|_{L^2}\le CB\varepsilon t^{\frac{\delta_0}{2}},\]
which is consequence of \eqref{est: bootstrap E02} with $k=0$ and of 
\eqref{ineq:L2 Zm vINF-vI} with $j=2$, we derive that
\[\left\|\chi(t^{-\sigma}D_x)\vNFGamma(t,\cdot)\right\|_{L^\infty}\le CB\varepsilon t^\frac{\delta_0}{2} + \sum_{j=1}^2 Ct^{-1+\beta}\left\|\chi(t^{-\sigma}D_x)\left[ x_j \NLNF\right](t,\cdot)\right\|_{L^2(dx)}.\]
The only thing we need to show in order to prove the statement is hence that
\begin{equation} \label{xj_NLNF_I=2}
\left\|\chi(t^{-\sigma}D_x)\left[ x_j \NLNF\right](t,\cdot)\right\|_{L^2(dx)}\le C(A+B)B\varepsilon^2 t^{\beta+\frac{\delta_1+\delta_2}{2}}.
\end{equation}
But from \eqref{def_NLNF_I=2} and \eqref{explicit rNFkg-Gamma} with $|I|=2$ we have that
\begin{equation} 
\begin{split}
&\left\|\chi(t^{-\sigma}D_x)\left[ x_j \NLNF\right](t,\cdot)\right\|_{L^2(dx)} \lesssim \|x_j\textit{NL}^I_{kg}(t,\cdot)\|_{L^2}\left(\|u_\pm(t,\cdot)\|_{L^\infty}+ \|\mathrm{R}_1u_\pm(t,\cdot)\|_{L^\infty}\right)\\
& + \sum_{\mu=0}^1 t^\sigma \left(\|x_j^\mu v_\pm(t,\cdot)\|_{L^\infty}+\left\| x_j^\mu \frac{D_x}{\langle D_x\rangle}v_\pm(t,\cdot)\right\|_{L^\infty}\right)\|v^I_\pm(t,\cdot)\|_{L^2}\|v_\pm(t,\cdot)\|_{H^{2,\infty}} \\
& + \sum_{\substack{(I_1,I_2)\in\mathcal{I}(I)\\ |I_1|<2}}\left\| \chi(t^{-\sigma}D_x)\left[x_j Q^\mathrm{kg}_0(v^{I_1}_\pm, D u^{I_2}_\pm)\right](t,\cdot)\right\|_{L^2}
\end{split}
\end{equation}
so \eqref{xj_NLNF_I=2} follows from a-priori estimates, \eqref{norm_Linfty_xv-}, \eqref{xj_GammaI_NLkg_I=2} and \eqref{xn_Qkg0_vI1uI2}. As $\delta_2\ll \delta_1\ll \delta_0$, that concludes that
\begin{equation}\label{est_Linfty_vINF_2}
\left\|\chi(t^{-\sigma}D_x)\vNFGamma(t,\cdot)\right\|_{L^\infty} \le C B\varepsilon t^{-1+\beta+\frac{\delta_0}{2}}.
\end{equation}
\endproof
\end{lem}

\begin{lem}
There exists a positive constant $C>0$ such that, for any multi-index $I$ of length 2, any $\chi\in C^\infty_0(\mathbb{R}^2)$, $\sigma>0$ small, $j=1,2$, and every $t\in [1,T]$
\begin{equation} \label{est:Linfty_xjGammaI_v}
\left\|\chi(t^{-\sigma}D_x)\left[x_j (\Gamma^Iv)_\pm\right](t,\cdot)\right\|_{L^\infty}\le CB\varepsilon t^{\beta + \frac{\delta_0}{2}},
\end{equation}
with $\beta>0$ small, $\beta\rightarrow 0$ as $\sigma\rightarrow 0$.
\proof
If $\Gamma^I$ contains at least one spatial derivative \eqref{est:Linfty_xjGammaI_v} is satisfied after \eqref{est: boostrap vpm}, \eqref{norm_Linfty_xv-} and \eqref{norm_Linfty_xjGammav-}.
Let us then assume that $\Gamma^I$ is product of two Klainerman vector fields.

From equation \eqref{half KG Gammav}, equality \eqref{xjw_Zjw} with $w=(\Gamma^I v)_{-}$, the $L^2-L^\infty$ continuity of operator $\chi(t^{-\sigma}D_x)\langle D_x\rangle^{-1}$ with norm $O(t^\sigma)$, and the $L^\infty$ continuity of $\chi(t^{-\sigma}D_x)D_x\langle D_x\rangle^{-1}$ with norm $O(t^{\sigma})$, we derive that
\begin{multline} \label{xj_GammaIv_preliminary}
\left\|\chi(t^{-\sigma}D_x)\left[x_j (\Gamma^Iv)_{-}\right](t,\cdot)\right\|_{L^\infty(dx)} \lesssim t^\sigma \|Z_j (\Gamma^I v)_{-}(t,\cdot)\|_{L^2}+ t\left\|\chi(t^{-\sigma}D_x)(\Gamma^I v)_{-}(t,\cdot)\right\|_{L^\infty}\\+t^{\sigma}\left\|\chi(t^{-\sigma}D_x)\left[x_j\Gamma^I \textit{NL}_{kg}\right](t,\cdot)\right\|_{L^\infty}.
\end{multline}
Reminding \eqref{GammaI_NLkg} and applying lemma \ref{Lem_appendix:L_estimate of products} with $L=L^\infty$ and $w_1 = (D_x\langle D_x\rangle^{-1})^\mu v^I_\pm$, for $|\mu|=0,1$, to the contribution coming from the first quadratic term in the right hand side of \eqref{GammaI_NLkg}, we find that there is some $\chi_1\in C^\infty_0(\mathbb{R}^2)$ such that 
\begin{equation} \label{xj_GammaI_NLkg_Linfty}
\begin{split}
\left\|\chi(t^{-\sigma}D_x)\left[x_j\Gamma^I \textit{NL}_{kg}\right](t,\cdot)\right\|_{L^\infty} &\lesssim \sum_{\mu, \nu=0}^1  \left\|\chi_1(t^{-\sigma}D_x)\left[x^\mu_j (\Gamma^I v)_\pm\right](t,\cdot)\right\|_{L^\infty}\|\mathrm{R}_1^\nu u_\pm(t,\cdot)\|_{H^{2,\infty}} \\
&+ t^{-N(s)}\sum_{\mu=0}^1 \left\|x^\mu_j(\Gamma^I v)_\pm(t,\cdot)\right\|_{L^2}\|u_\pm(t,\cdot)\|_{H^s}\\& + \sum_{\substack{(I_1,I_2)\in\mathcal{I}(I)\\ |I_1|<2}}\left\|\chi(t^{-\sigma}D_x)\left[x_jQ^\mathrm{kg}_0\left(v^{I_1}_\pm, D u^{I_2}_\pm\right)\right](t,\cdot)\right\|_{L^\infty}.
\end{split}
\end{equation}
Therefore, picking $s>0$ large so that $N(s)>1$ and using the $L^2-L^\infty$ continuity of $\chi_1(t^{-\sigma}D_x)$ with norm $O(t^{\sigma})$, together with the estimates \eqref{est: bootstrap argument a-priori est}, \eqref{norm_L2_xj-GammaIv-} with $k=2$, along with \eqref{xn_Qkg0_vI1uI2}, we find at first that
\begin{equation*}
\left\|\chi(t^{-\sigma}D_x)\left[x_j\Gamma^I \textit{NL}_{kg}\right](t,\cdot)\right\|_{L^\infty} \le CAB\varepsilon^2 t^{\frac{1}{2}+\sigma+\frac{\delta_1}{2}}.
\end{equation*}
Injecting the above estimate, together with \eqref{est: bootstrap E02} and \eqref{est:Linfty_vI_I=2-new}, into \eqref{xj_GammaIv_preliminary} we derive that
\begin{equation*}
\left\|\chi(t^{-\sigma}D_x)\left[x_j (\Gamma^Iv)_{-}\right](t,\cdot)\right\|_{L^\infty} \le CB\varepsilon t^{\frac{1}{2}+\sigma+\frac{\delta_1}{2}}.
\end{equation*}
The above inequality holds for any $\chi\in C^\infty_0(\mathbb{R}^2)$, so injecting it into \eqref{xj_GammaI_NLkg_Linfty} and using again a-priori estimates, \eqref{norm_L2_xj-GammaIv-}, \eqref{xn_Qkg0_vI1uI2}, together with the fact that $\beta+ (\delta + \delta_2)/2\le \delta_1/2$ as $\beta$ is as small as we want as long as $\sigma$ is small and $\delta, \delta_2\ll \delta_1$, we find the following enhanced estimate
\begin{equation*}
\left\|\chi(t^{-\sigma}D_x)\left[x_j\Gamma^I \textit{NL}_{kg}\right](t,\cdot)\right\|_{L^\infty}\le C(A+B)B\varepsilon^2 t^{\sigma +\frac{\delta_1}{2}}.
\end{equation*}
Consequently, summing up this latter one with \eqref{est: bootstrap E02} and \eqref{est:Linfty_vI_I=2-new}, we end up with \eqref{est:Linfty_xjGammaI_v}.
\endproof
\end{lem}

\begin{lem}
Let $\Gamma\in \mathcal{Z}$ be an admissible vector field. There exists a positive constant $C$ such that, for any $\chi\in C^\infty_0(\mathbb{R}^2)$, $\sigma>0$ small, $i,j=1,2$, and every $t\in [1,T]$,
\begin{equation} \label{est:Linfty_xixjGammav-new}
\left\|\chi(t^{-\sigma}D_x)\left[x_i x_j (\Gamma v)_\pm(t,\cdot)\right]\right\|_{L^\infty(dx)}\le CB\varepsilon t^{1+\beta+\frac{\delta_1}{2}},
\end{equation}
with $\beta>0$ such that $\beta\rightarrow 0$ as $\sigma\rightarrow 0$.
\proof
Since $(\Gamma v)_+ = -\overline{(\Gamma v)_+}$ we reduce to prove that inequality \eqref{est:Linfty_xixjGammav-new} holds true for $(\Gamma v)_{-}$.
Moreover, we only focus on the case where $\Gamma\in \{\Omega, Z_m, m=1,2\}$ is a Klainerman vector field, as \eqref{est:Linfty_xixjGammav-new} with $\Gamma$ being a spatial derivative is simply a consequence of \eqref{norm_Linfty_xixjv-}.

We remind that $(\Gamma v)_{-}$ is solution to non-linear Klein-Gordon equation \eqref{half KG Gammav} with $\Gamma^I=\Gamma$, and that the non-linearity $\Gamma\Nlkg$ is given by \eqref{notation: NLGamma}. Hence, multiplying $x_i$ to relation \eqref{xjw_Zjw} with $w=(\Gamma v)_{-}$ and making some commutations we find that
\begin{equation} \label{est_1_xixjGammav}
\begin{split}
& \left\|\chi(t^{-\sigma}D_x)\left[x_i x_j (\Gamma v)_\pm(t,\cdot)\right]\right\|_{L^\infty(dx)} \\
&\lesssim \sum_{\mu=0}^1 \left[\left\| \chi(t^{-\sigma}D_x)\left[x_i^\mu Z_j(\Gamma v)_{-}\right](t,\cdot)\right\|_{L^\infty(dx)} + t \left\| \chi(t^{-\sigma}D_x)\left[x_i^\mu (\Gamma v)_{-}\right](t,\cdot)\right\|_{L^\infty(dx)}\right] \\
& + \sum_{\mu=0}^1\left\|\chi(t^{-\sigma}D_x)\left[x_i^{\mu} x_j \Gamma \textit{NL}_{kg}\right](t,\cdot)\right\|_{L^\infty(dx)}.
\end{split}
\end{equation}
At first, we estimate the latter contribution in the above right hand side using that $\chi(t^{-\sigma}D_x)$ is a continuous $L^2-L^\infty$ operator with norm $O(t^{\sigma})$ together with estimates \eqref{est: bootstrap argument a-priori est}, \eqref{Hs_norm_DtU} with $s=0$, \eqref{norm_Linfty_xv-}, \eqref{norm_L2_xj-GammaIv-}, \eqref{norm_Linfty_xixjv-}, \eqref{est:xixjGamma v-}:
\begin{equation}\label{first_estimate_xixjGammaNL}
\begin{split}
& \sum_{\mu=0}^1\left\|\chi(t^{-\sigma}D_x)\left[x_i^{\mu} x_j \Gamma \textit{NL}_{kg}\right](t,\cdot)\right\|_{L^\infty(dx)} \lesssim \sum_{\mu=0}^1t^{\sigma}\left\|\chi(t^{-\sigma}D_x)\left[x_i^{\mu} x_j \Gamma \textit{NL}_{kg}\right](t,\cdot)\right\|_{L^2(dx)}\\
& \lesssim  \sum_{\mu_1, \mu_2, \nu=0}^1 t^{\sigma} \|x_i^{\mu_1} x_j^{\mu_2} (\Gamma v)_\pm(t,\cdot)\|_{L^2(dx)}\|\mathrm{R}^\nu_1 u_\pm(t,\cdot)\|_{H^{2,\infty}}\\
& + \sum_{\mu, |\nu|=0}^1 t^{\sigma} \left\|x_i^{\mu} x_j\Big(\frac{D_x}{\langle D_x\rangle}\Big)^\nu v_\pm(t,\cdot)\right\|_{L^\infty(dx)}\left(\| (\Gamma u)_\pm(t,\cdot)\|_{H^1}+\|u_\pm(t,\cdot)\|_{H^1}+\|D_tu_\pm(t,\cdot)\|_{L^2}\right)\\
& \le C(A+B)B \varepsilon^2 t^{\frac{3}{2}+\sigma+\frac{\delta_2}{2}}.
\end{split}
\end{equation}\normalsize
Injecting this estimate, along with \eqref{Linfty_est_VJ}, \eqref{norm_Linfty_xjGammav-}, \eqref{est:Linfty_vI_I=2-new} and \eqref{est:Linfty_xjGammaI_v}, into \eqref{est_1_xixjGammav} we deduce that for any smooth cut-off function $\chi$
\begin{equation} \label{first_est_xixjGammav}
 \left\|\chi(t^{-\sigma}D_x)\left[x_i x_j (\Gamma v)_\pm(t,\cdot)\right]\right\|_{L^\infty(dx)} \le CB\varepsilon t^{\frac{3}{2}+\sigma+\frac{\delta_2}{2}}.
\end{equation}
Now, if we change the approach of bounding the $L^\infty(dx)$ norm of $x_i^{\mu}x_j Q^\mathrm{kg}_0((\Gamma v)_{-}, D_1u_\pm)$, which is one of the contributions to $x_i^{\mu}x_j\Gamma\textit{NL}_{kg}$ after \eqref{notation: NLGamma}, and make use of lemma \ref{Lem_appendix:L_estimate of products} with $L=L^\infty$ instead of \eqref{first_estimate_xixjGammaNL}, we see that
\begin{equation*}
\begin{split}
&\sum_{\mu=0}^1 \left\|\chi(t^{-\sigma}D_x)\left[x_i^{\mu} x_j \Gamma \textit{NL}_{kg}\right](t,\cdot)\right\|_{L^\infty(dx)}\\
& \lesssim \sum_{\mu_1, \mu_2, \nu=0}^1\left\| \chi_1(t^{-\sigma}D_x)\left[x_i^{\mu_1}x_j^{\mu_2} (\Gamma v)_\pm\right](t,\cdot)\right\|_{L^\infty(dx)}\|\mathrm{R}^\nu_1u_\pm(t,\cdot)\|_{H^{2,\infty}}\\
& + \sum_{\mu_1, \mu_2=0}^1 t^{-N(s)}\|x_i^{\mu_1}x_j^{\mu_2} (\Gamma v)_\pm(t,\cdot)\|_{L^2(dx)}\|u_\pm(t,\cdot)\|_{H^s}\\
& + \sum_{\mu, |\nu|=0}^1 t^{\sigma} \left\|x_i^{\mu} x_j\Big(\frac{D_x}{\langle D_x\rangle}\Big)^\nu v_\pm(t,\cdot)\right\|_{L^\infty(dx)}\left(\| (\Gamma u)_\pm(t,\cdot)\|_{H^1}+\|u_\pm(t,\cdot)\|_{H^1}+\|D_tu_\pm(t,\cdot)\|_{L^2}\right).
\end{split}
\end{equation*}
Then, choosing $s>0$ sufficiently large so that $N(s)\ge 3$ and using again \eqref{est: bootstrap argument a-priori est}, \eqref{norm_Linfty_xv-}, \eqref{norm_L2_xj-GammaIv-} with $k=1$, \eqref{norm_Linfty_xixjv-}, \eqref{est:xixjGamma v-}, \eqref{norm_Linfty_xjGammav-}, together with \eqref{first_est_xixjGammav}, we obtain that
\begin{equation*}
\sum_{\mu=0}^1 \left\|\chi(t^{-\sigma}D_x)\left[x_i^{\mu} x_j \Gamma \textit{NL}_{kg}\right](t,\cdot)\right\|_{L^\infty(dx)} \le C(A+B)B\varepsilon^2 t^{1+\sigma+\frac{\delta_2}{2}},
\end{equation*}
which enhances \eqref{first_estimate_xixjGammaNL} of a factor $t^{1/2}$.
Combining the above estimate with \eqref{Linfty_est_VJ}, \eqref{norm_Linfty_xjGammav-}, \eqref{est:Linfty_vI_I=2-new} and \eqref{est:Linfty_xjGammaI_v}, we finally end up with the result of the statement.
\endproof
\end{lem}

\begin{lem}
Let $\Gamma\in \{\Omega, Z_m, m=1,2\}$ be a Klainerman vector field, $\vNFGamma$ the function defined in \eqref{def_vNF-Gamma} with $|I|=1$ and $\Gamma^I=\Gamma$, and $B^k_{(j_1,j_2,j_3)}(\xi,\eta)$ the multiplier introduced in \eqref{def of B(i1,i2,i3)} (resp. in \eqref{multiplier_B3}) for any $k=1,2$ (resp. $k=3$), any $j_i\in \{+,-\}$ for $i=1,2,3$. 
Let us define\index{VNFGamma@$\VNF$, normal form function defined from $(\Gamma v)_{-}$}
\begin{equation} \label{def_VNF}
\begin{split}
\VNF(t,x):= \vNFGamma(t,x) &-\frac{i}{4(2\pi)^2}\sum_{j_1,j_2\in \{+,-\}}\int e^{ix\cdot\xi}B^1_{(j_1,j_2,+)}(\xi,\eta) \hat{v}_{j_1}(\xi-\eta) \widehat{(\Gamma u)_{j_2}}(\eta) d\xi d\eta\\
& +\delta_{\Omega}\frac{i}{4(2\pi)^2}\sum_{j_1,j_2\in \{+,-\}}\int e^{ix\cdot\xi}B^2_{(j_1,j_2,+)}(\xi,\eta) \hat{v}_{j_1}(\xi-\eta) \hat{u}_{j_2}(\eta) d\xi d\eta \\
& + \delta_{Z_1}\frac{i}{4(2\pi)^2}\sum_{j_1,j_2\in \{+,-\}}\int e^{ix\cdot\xi}B^3_{(j_1,j_2,+)}(\xi,\eta) \hat{v}_{j_1}(\xi-\eta) \hat{u}_{j_2}(\eta) d\xi d\eta,
\end{split}
\end{equation}
where $\delta_\Omega$ (resp. $\delta_{Z_1}$) is equal to 1 if $\Gamma=\Omega$ (resp. if $\Gamma=Z_1$), 0 otherwise. There exists a constant $C>0$ such that, for any $\chi\in C^\infty_0(\mathbb{R}^2)$, $\sigma>0$, and every $t\in [1,T]$,
\begin{equation}\label{est:VNF-(Gammav)}
\left\|\chi(t^{-\sigma}D_x)\big(\VNF - (\Gamma v)_{-}\big)(t,\cdot) \right\|_{L^\infty} \le C(A+B)A\varepsilon^2 t^{-\frac{5}{4}}.
\end{equation}
Moreover, for every $m=1,2$ and $t\in [1,T]$
\begin{equation} \label{est:L2_Zm(VNF-Gammav)}
\left\||\chi(t^{-\sigma}D_x)Z_m\big(\VNF - (\Gamma v)_{-}\big)(t,\cdot) \right\|_{L^2} \le C(A+B)B\varepsilon^2 t^{3\sigma +\delta_2}.
\end{equation}
\proof
From definition \eqref{def_VNF} of $\VNF$ and equalities \eqref{explicit integral B}, \eqref{explicit_integral_B3}, we find that
\begin{equation} \label{explicit VNF - VJ}
\begin{split}
\VNF - (\Gamma v)_{-} &=  \vNFGamma - (\Gamma v)_{-}\\
& -\frac{i}{2}\left[(D_t v)(D_1 \Gamma u) - (D_1 v)(D_t \Gamma u) + D_1[v (D_t \Gamma u)] - \langle D_x\rangle [v (D_1 \Gamma u)] \right] \\
& +\delta_{\Omega} \frac{i}{2}\left[(D_tv)(D_2u) - (D_2v)(D_tu) + D_2[v D_tu] - \langle D_x\rangle [v D_2u] \right] \\
& +\delta_{Z_1} \frac{i}{2}\left[(D_t v)(D_t u) + v( |D_x|^2u )- \langle D_x\rangle[v (D_tu)]\right],
\end{split}
\end{equation}
where $\vNFGamma - (\Gamma v)_{-}$ has the explicit expression \eqref{explicit vNfGamma-vJ-}.
We use \eqref{def u+- v+-}, \eqref{def uIpm vIpm} and lemma \ref{Lem_appendix:L_estimate of products} with $L=L^\infty$, $w_1=\mathrm{R}_1^\mu u_\pm$ (resp. $w_1=\mathrm{R}_1^\mu(\Gamma u)_\pm$) for $\mu=0,1$, and $s>0$ large enough to have $N(s)\ge 2$, in order to estimate the $L^\infty$ norm of products appearing in \eqref{explicit vNfGamma-vJ-} (resp. in the second line in the above right hand side). For some new $\chi_1\in C^\infty_0(\mathbb{R}^2)$ we have that
\begin{equation*}
\begin{split}
&\left\|\chi(t^{-\sigma}D_x)\big(\VNF - (\Gamma v)_{-}\big)(t,\cdot) \right\|_{L^\infty}  \lesssim \left\|\chi(t^{-\sigma}D_x)\big(\vNFGamma - (\Gamma v)_{-}\big)(t,\cdot) \right\|_{L^\infty}  \\
& +\sum_{|\mu|=0}^1 t^\sigma \|v_\pm(t,\cdot)\|_{H^{1,\infty}}\left\|\chi_1(t^{-\sigma}D_x) \mathrm{R}^\mu_1 (\Gamma u)_\pm(t,\cdot)\right\|_{L^\infty} + t^{-2}\|v_\pm(t,\cdot)\|_{H^s} \|(\Gamma u)_\pm(t,\cdot)\|_{L^2}\\
& + \sum_{|\mu|=0}^1 t^\sigma \|v_\pm(t,\cdot)\|_{H^{1,\infty}} \|\mathrm{R}^\mu_1 u_\pm(t,\cdot)\|_{H^{2,\infty}},
\end{split}
\end{equation*}
with
\begin{equation*}
\begin{split}
\left\|\chi(t^{-\sigma}D_x)\big(\vNFGamma - (\Gamma v)_{-}\big)(t,\cdot) \right\|_{L^\infty}  &\lesssim \sum_{\mu=0}^1 t^\sigma \left\|\chi_1(t^{-\sigma}D_x)(\Gamma v)_\pm(t,\cdot)\right\|_{L^\infty}\|\mathrm{R}^\mu_1 u_\pm(t,\cdot)\|_{L^\infty}\\
& + t^{-2}\|(\Gamma v)_\pm(t,\cdot)\|_{L^2}\|u_\pm(t,\cdot)\|_{H^s}.
\end{split}
\end{equation*}
Estimate \eqref{est:VNF-(Gammav)} follows then from \eqref{est: bootstrap argument a-priori est}, \eqref{Linfty_est_VJ} and \eqref{Linfty_est_UJ}.

In order to derive \eqref{est:L2_Zm(VNF-Gammav)} we apply $Z_m$ to \eqref{explicit VNF - VJ} and use the Leibniz rule, reminding formulas \eqref{commutator_Z_DtD1}.
Among the quadratic terms coming out from the action of $Z_m$ on the second line in \eqref{explicit VNF - VJ} we see appear products where $Z_m$ is acting on $v$ and $\Gamma$ on $u$.
We estimate the $L^2$ norm of those ones, when truncated by operator $\chi(t^{-\sigma}D_x)$, using inequalities \eqref{cor_app_est_2} with $L=L^2$, $w=u$, $w_{j_0}=(D_x\langle D_x\rangle^{-1})^\mu Z_m v$ for $|\mu|=0,1$, and $s>0$ large enough to have $N(s)>1$.
We bound instead the $L^2$ norm of all other remaining products with the $L^\infty$ norm of factor that does not contain any vector field times the $L^2$ norm of the remaining one. Hence
\begin{equation}
\begin{split}
&\left\|\chi(t^{-\sigma}D_x)Z_m\big(\VNF - (\Gamma v)_{-}\big)(t,\cdot)\right\|_{L^2} \lesssim \left\|\chi(t^{-\sigma}D_x)Z_m\big(\vNFGamma - (\Gamma v)_{-}\big)(t,\cdot)\right\|_{L^2} \\
& + t^\sigma \left\|\chi_1(t^{-\sigma}D_x)(Z_m v)_\pm(t,\cdot)\right\|_{L^\infty}\|(\Gamma u)_\pm(t,\cdot)\|_{L^2}\\
& + t^{-N(s)}\Big(\sum_{|\mu|=0}^1\|x^\mu (Z_mv)_\pm(t,\cdot)\|_{L^2}+t\|(Z_m v)_\pm(t,\cdot)\|_{L^2}\Big)\left(\|u_\pm(t,\cdot)\|_{H^s}+\|D_tu_\pm(t,\cdot)\|_{H^s}\right)\\
&+t^\sigma \| v_\pm(t,\cdot)\|_{H^{1,\infty}}\left(\|(Z_m \Gamma u)_\pm(t,\cdot)\|_{L^2} + \|(\Gamma u)_\pm(t,\cdot)\|_{L^2}+ \|D_t(\Gamma u)_\pm(t,\cdot)\|_{L^2}\right)\\
& + \sum_{|\mu|=0}^1t^\sigma\|(\Gamma v)_\pm(t,\cdot)\|_{L^2}\|\mathrm{R}^\mu u_\pm(t,\cdot)\|_{H^{2,\infty}} \\
&+ t^\sigma\|v_\pm(t,\cdot)\|_{H^{1,\infty}}\left(\|(Z_m u)_\pm(t,\cdot)\|_{H^1}+ \|u_\pm(t,\cdot)\|_{H^1}+ \|D_tu_\pm(t,\cdot)\|_{L^2}\right),
\end{split}
\end{equation}
and estimate \eqref{est:L2_Zm(VNF-Gammav)} is obtained from \eqref{est: bootstrap argument a-priori est}, \eqref{Hs_norm_DtU},  \eqref{DtUI}, \eqref{norm_L2_xj-GammaIv-}, \eqref{ineq:L2 Zm vINF-vI} with $j=1$,  and \eqref{Linfty_est_VJ}.
\endproof
\end{lem}

\begin{lem} \label{Lem_appendix: LvtildeGamma}
Let $\Gamma\in \{\Omega, Z_m, m=1,2\}$ be a Klainerman vector field, $\VNF$ the function defined in \eqref{def_VNF} and  \index{VtildeGamma@$\widetilde{V}^\Gamma$, function $\VNF$ in semi-classical setting}
\begin{equation}\label{def_VtildeGamma}
\widetilde{V}^\Gamma(t,x) :=t\VNF(t,tx).
\end{equation}
There exists a positive constant $C>0$ such that, for any $\chi \in C^\infty_0(\mathbb{R}^2)$, $\sigma>0$ small, and every $t\in [1,T]$,
\begin{subequations}
\begin{gather}
\left\| \widetilde{V}^\Gamma(t,\cdot)\right\|_{L^2}\le CB\varepsilon t^\frac{\delta_2}{2}, \label{est:VtildeGamma}\\
\sum_{|\mu|=1}\left\|\oph(\chi(h^\sigma\xi))\mathcal{L}^\mu \widetilde{V}^\Gamma(t,\cdot)\right\|_{L^2}\le CB\varepsilon t^\frac{\delta_1}{2}.\label{est:Lcal_VtildeGamma}
\end{gather}
\end{subequations}
\proof
Let us recall equalities \eqref{explicit vNfGamma-vJ-} with $\Gamma^I=\Gamma$ and \eqref{explicit VNF - VJ}.
From a-priori estimates we immediately derive that, for every $t\in [1,T]$, 
\begin{equation*} 
\|[\VNF -(\Gamma v)_{-}](t,\cdot)\|_{L^2} \le CAB\varepsilon t^{-\frac{1}{2}+\frac{\delta_2}{2}+\sigma},
\end{equation*}
and consequently that
\begin{equation} \label{est_L2_VNF}
\|\widetilde{V}^\Gamma(t,\cdot)\|_{L^2}=\|\VNF(t,\cdot)\|_{L^2}\le CB\varepsilon t^{\frac{\delta_2}{2}}.
\end{equation}

Using definition \eqref{def_VNF} one can check that $\VNF$ is solution to
\begin{equation} \label{half KG VNF}
[D_t + \langle D_x\rangle]\VNF(t,x) = \NLT(t,x) - \delta_{Z_1} Q^\mathrm{kg}_0\big(v_\pm, Q^\mathrm{w}_0(v_\pm, D_1 v_\pm)\big)
\end{equation}
with
\begin{equation} \label{def_NLT}
\begin{split}
\NLT(t,x)& = r^{I,\textit{NF}}_{kg}(t,x)\\
&-  \frac{i}{4(2\pi)^2}\int e^{ix\cdot\xi}B^1_{(j_1,j_2,+)}(\xi,\eta) \left[\widehat{\textit{NL}_{kg}}(\xi-\eta) \widehat{(\Gamma u)_{j_2}}(\eta) - \hat{v}_{j_1}(\xi-\eta)\reallywidehat{\Gamma\textit{NL}_{w}}(\eta)\right] d\xi d\eta \\
& + \delta_{\Omega} \frac{i}{4(2\pi)^2}\int e^{ix\cdot\xi}B^2_{(j_1,j_2,+)}(\xi,\eta) \left[\widehat{\textit{NL}_{kg}}(\xi-\eta) \hat{u}_{j_2}(\eta) - \hat{v}_{j_1}(\xi-\eta)\reallywidehat{\textit{NL}_{w}}(\eta)\right] d\xi d\eta\\
&+ \delta_{Z_1} \frac{i}{4(2\pi)^2}\int e^{ix\cdot\xi}B^3_{(j_1,j_2,+)}(\xi,\eta) \left[\widehat{\textit{NL}_{kg}}(\xi-\eta) \hat{u}_{j_2}(\eta) - \hat{v}_{j_1}(\xi-\eta)\reallywidehat{\textit{NL}_{w}}(\eta)\right] d\xi d\eta,
\end{split}
\end{equation}
and $r^{I,\textit{NF}}_{kg}$ given by \eqref{def:rNF-Gamma-kg} (or, explicitly, by \eqref{explicit rNFkg-Gamma}) with $|I|=1$. Superscript \emph{c} in $\NLT$ stands for \emph{cubic} and wants to stress out the fact that, passing from function $(\Gamma v)_{-}$ to $\VNF$, we have replaced all quadratic terms in the right hand side of \eqref{half KG Gammav} (when $|I|=1$ and $\Gamma^I=\Gamma$) with cubic ones. Hence, from relation \eqref{relation_Zjw_Ljwidetilde(w)} with $w=\VNF$ and equation \eqref{half KG VNF} we get that
\begin{equation} \label{ineq:Lcal_VTildeGamma_proof}
\begin{split}
&\left\|\oph(\chi(h^\sigma\xi))\mathcal{L}_m \widetilde{V}^\Gamma(t,\cdot)\right\|_{H^1} \lesssim \|\chi(t^{-\sigma}D_x)Z_m\VNF(t,\cdot)\|_{L^2} + \left\|\oph(\chi(h^\sigma\xi)\xi_m\langle\xi\rangle^{-1})\widetilde{V}^\Gamma(t,\cdot)\right\|_{L^2}\\
& + \left\|\chi(t^{-\sigma}D_x)\left[x_m\NLT\right](t,\cdot) \right\|_{L^2(dx)} + \delta_{Z_1}\left\|\chi(t^{-\sigma}D_x)\left[x_mQ^\mathrm{kg}_0\big(v_\pm, Q^\mathrm{w}_0(v_\pm, D_1 v_\pm)\big)\right](t,\cdot) \right\|_{L^2(dx)}
\end{split}
\end{equation}
After \eqref{est: bootstrap E02} with $k=1$,\eqref{est:L2_Zm(VNF-Gammav)}, and the fact that $\sigma$ can be chosen sufficiently small so that $3\sigma+\delta_2\le \delta_1/2$, as $\delta_2\ll \delta_1$, it is straightforward to see that
\begin{equation}\label{est:L2_ZmVNF}
\|\chi(t^{-\sigma}D_x)Z_m\VNF(t,\cdot)\|_{L^2} \le CB\varepsilon t^\frac{\delta_1}{2}.
\end{equation}
Moreover, from \eqref{est L2 NLw}, \eqref{norm_Linfty_xv-} and a-priori estimates,
\begin{multline}\label{est:L2_xj_cubic}
\left\|\chi(t^{-\sigma}D_x)\left[x_mQ^\mathrm{kg}_0\big(v_\pm, Q^\mathrm{w}_0(v_\pm, D_1 v_\pm)\big)\right](t,\cdot) \right\|_{L^2(dx)} \\\lesssim \sum_{|\mu|=0}^1 \left\|x_n\Big(\frac{D_x}{\langle D_x\rangle}\Big)^\mu v_\pm(t,\cdot) \right\|_{L^\infty(dx)} \|\textit{NL}_w(t,\cdot)\|_{L^2} \le C(A+B)AB \varepsilon^3 t^{-1+\frac{\delta+\delta_2}{2}}.
\end{multline}
Using instead equalities \eqref{explicit integral B} and \eqref{explicit_integral_B3} we derive the following explicit expression for $\NLT$:
\begin{equation} \label{explicit NLT}
\begin{split}
\NLT(t,x)= r^{I,\textit{NF}}_{kg}(t,x) &- \frac{i}{2}\left[\textit{NL}_{kg}(D_1 \Gamma u) - (D_1v)\Gamma\textit{NL}_w + D_1[v \Gamma\textit{NL}_w]\right] \\
& + \delta_\Omega \frac{i}{2}\left[\textit{NL}_{kg}(D_2 u) - (D_2v)\textit{NL}_w+ D_2[v \textit{NL}_w]\right] \\
&+\delta_{Z_1}\left[\textit{NL}_{kg}(D_tu) + (D_tv)\textit{NL}_w - \langle D_x\rangle[v\textit{NL}_w]\right].
\end{split}
\end{equation}
Hence, reminding estimates \eqref{est: bootstrap argument a-priori est}, \eqref{est L2 NLw}, \eqref{Hs_norm_DtU} with $s=0$, \eqref{norm_Linfty_xv-}, \eqref{enhanced_xjrINF_kg}, and equality \eqref{Gamma_Nlw} from which follows that
\begin{equation}\label{L2_norm_NLIw}
\|\Gamma\textit{NL}_w(t,\cdot)\|_{L^2}\lesssim \|v_\pm(t,\cdot)\|_{H^{1,\infty}}\left(\|v^I_\pm(t,\cdot)\|_{H^1}+\|v_\pm(t,\cdot)\|_{H^1}+\|D_tv_\pm(t,\cdot)\|_{L^2}\right),
\end{equation}
we find that
\begin{equation}\label{first_estimate_xNLT}
\begin{split}
&\left\|\chi(t^{-\sigma}D_x)\left[x_j \NLT\right](t,\cdot)\right\|_{L^2(dx)} \lesssim \left\|\chi(t^{-\sigma}D_x)\left[x_j r^{I,\textit{NF}}_{kg}(t,\cdot)\right](t,\cdot)\right\|_{L^2(dx)} \\& + \sum_{|\mu|, \nu=0}^1 \left\|x_j\Big(\frac{D_x}{\langle D_x\rangle}\Big)^\mu v_\pm(t,\cdot) \right\|_{L^\infty(dx)} \|\mathrm{R}^\nu u_\pm(t,\cdot)\|_{H^{2,\infty}}\left(\|(\Gamma u)_\pm(t,\cdot)\|_{L^2}+ \|u_\pm(t,\cdot)\|_{L^2}\right) \\
& +  \sum_{k,|\mu|=0}^1 \left\|x_j^k\Big(\frac{D_x}{\langle D_x\rangle}\Big)^\mu v_\pm(t,\cdot) \right\|_{L^\infty(dx)} \left(\|\Gamma\textit{NL}_w(t,\cdot)\|_{L^2}+\|\textit{NL}_w(t,\cdot)\|_{L^2}\right)\\
& \le C(A+B)AB \varepsilon^2 t^{-\frac{1}{2}+\beta+\frac{\delta+\delta_1}{2}}.
\end{split}
\end{equation}
By injecting the above estimate, together with \eqref{est:VtildeGamma}, \eqref{est:L2_ZmVNF}, \eqref{est:L2_xj_cubic}, into \eqref{ineq:Lcal_VTildeGamma_proof} we finally deduce \eqref{est:Lcal_VtildeGamma} and conclude the proof of the statement.
\endproof
\end{lem}

\begin{lem} \label{Lem_appendix:est_Zm(vINF-Gammav)}
Let $\Gamma\in \{\Omega, Z_m, m=1,2\}$ be a Klainerman vector field and $I_1, I_2$ two multi-indices such that $\Gamma^{I_1}=\Gamma$ and $\Gamma^{I_2}=Z_m\Gamma$, with $m\in\{1,2\}$. Let also $v^{I,\textit{NF}}$ be the function defined in \eqref{def_vNF-Gamma} for a generic multi-index $I$ of length equal to 1 or 2.
There exists a constant $C>0$ such that, for any $\chi\in C^\infty_0(\mathbb{R}^2)$, $\sigma>0$ small, $m,n=1,2$, every $t\in [1,T]$,
\begin{subequations}\label{est:L2_Zn(vINF-Gammav)_xmZN(vINF-Gammav)}
\begin{multline}\label{est:L2_Zn(vINF-Gammav)}
\left\|\chi(t^{-\sigma}D_x)\left[Z_m \left(v^{I_1,\textit{NF}}  - (\Gamma v)_{-}\right)\right](t,\cdot) \right\|_{L^2} +  \left\|\chi(t^{-\sigma}D_x)\left(v^{I_2,\textit{NF}}  - (Z_m\Gamma v)_{-}\right)(t,\cdot) \right\|_{L^2}\\
 \le C(A+B)B\varepsilon^2 t^{-1+\beta + \frac{\delta+\delta_1+\delta_2}{2}}
\end{multline}
and
\begin{multline}\label{est:L2_xmZn(vINF-Gammav)}
\left\|\chi(t^{-\sigma}D_x)\left[x_nZ_m \left(v^{I_1,\textit{NF}} - (\Gamma v)_{-}\right)\right](t,\cdot) \right\|_{L^2}+\left\|\chi(t^{-\sigma}D_x)\left[x_n\left(v^{I_2,\textit{NF}}  - (Z_m\Gamma v)_{-}\right)\right](t,\cdot) \right\|_{L^2} \\
 \le C(A+B)B\varepsilon^2 t^{\beta + \frac{\delta+\delta_1+\delta_2}{2}},
\end{multline}
\end{subequations}
with $\beta>0$ small such that $\beta\rightarrow 0$ as $\sigma\rightarrow 0$.
Moreover, if $\VNF$ is the function defined in \eqref{def_VNF}, then for every $t\in [1,T]$
\begin{subequations} \label{est:L2_Zm(VNF-Gammav)_xnZm(VNF-Gammav)}
\begin{align}
\left\|\chi(t^{-\sigma}D_x)\left[Z_m \left(\VNF  - (\Gamma v)_{-}\right)\right](t,\cdot) \right\|_{L^2} &\le C(A+B)B\varepsilon^2 t^{-1+\beta + \frac{\delta+\delta_1+\delta_2}{2}}, \label{est:L2_Zm(VNF-Gammav)_enhanced} \\
\left\|\chi(t^{-\sigma}D_x)\left[x_nZ_m \left(\VNF - (\Gamma v)_{-}\right)\right](t,\cdot) \right\|_{L^2}
 &\le C(A+B)B\varepsilon^2 t^{\beta + \frac{\delta+\delta_1+\delta_2}{2}}. \label{est:L2_xnZm(VNF-Gammav)}
\end{align}
\end{subequations}
\proof
We warn the reader that throughout the proof we denote by $C$ and $\beta$ two positive constants that may change line after line, with $\beta\rightarrow 0$ as $\sigma\rightarrow 0$.

We refer to equality \eqref{Zm(vgamma-vJ)} with $I=I_1$ and bound the $L^2$ norm of each product in the first, third and fifth line of its right hand side by means of lemma \ref{Lem_appendix:L_estimate of products} with $L=L^2$.
The $L^2$ norm of the remaining products in the second line of the mentioned equality is instead estimated using inequalities \eqref{cor_app_est_2} with $L=L^2$ and $w_{j_0}=(D_x\langle D_x\rangle^{-1})^\mu (\Gamma^{I_1}v)_\pm$.
In this way we obtain that there is some $\chi_1\in C^\infty_0(\mathbb{R}^2)$ such that
\begin{equation*}
\begin{split}
&\left\|\chi(t^{-\sigma}D_x)\left[Z_m \left(v^{I_1,\textit{NF}}  - (\Gamma v)_{-}\right)\right](t,\cdot) \right\|_{L^2} \lesssim t^\sigma \|\chi_1(t^{-\sigma}D_x)(Z_m \Gamma v)_\pm(t,\cdot)\|_{L^\infty}\|u_\pm(t,\cdot)\|_{L^2} \\
&+ t^{-N(s)}\|(Z_m\Gamma v)_\pm(t,\cdot)\|_{L^2}\|u_\pm(t,\cdot)\|_{H^s} + t^\sigma \|\chi_1(t^{-\sigma}D_x)(\Gamma v)_\pm(t,\cdot)\|_{L^\infty} \|(Z_mu)_\pm(t,\cdot)\|_{L^2}\\
& +t^{-N(s)} \left(\sum_{\mu=0}^1\|x_m^\mu (\Gamma v)_\pm(t,\cdot)\|_{L^2} + t\|(\Gamma v)_\pm(t,\cdot)\|_{L^2}\right) \left(\|u_\pm(t,\cdot)\|_{H^s}+ \|D_tu_\pm(t,\cdot)\|_{H^s}\right)\\
&+ t^\sigma\|\chi_1(t^{-\sigma}D_x)(\Gamma v)_\pm(t,\cdot)\|_{L^\infty}\left(\|u_\pm(t,\cdot)\|_{L^2} + \|D_t u_\pm(t,\cdot)\|_{L^2}\right) \\
& +t^{-N(s)}\|(\Gamma v)_\pm(t,\cdot)\|_{L^2}\left(\|u_\pm(t,\cdot)\|_{H^s}+ \|D_t u_\pm(t,\cdot)\|_{H^s}\right).
\end{split}
\end{equation*} 
Choosing $s>0$ large so that $N(s)>1$ and using estimates \eqref{est: bootstrap argument a-priori est}, \eqref{Hs_norm_DtU}, \eqref{norm_L2_xj-GammaIv-}, together with lemmas \ref{Lem_appendix: preliminary est VJ} and \ref{Lem_appendix:est vI I=2}, we hence find that
\[\left\|\chi(t^{-\sigma}D_x)\left[Z_m \left(v^{I_1,\textit{NF}}  - (\Gamma v)_{-}\right)\right](t,\cdot) \right\|_{L^2}  \le C(A+B)B\varepsilon^2 t^{-1+\beta + \frac{\delta+\delta_1+\delta_2}{2}}.\]
Analogously,\small
\begin{equation*}
\begin{split}
&\left\|\chi(t^{-\sigma}D_x)\left[x_n Z_m \left(v^{I_1,\textit{NF}}  - (\Gamma v)_{-}\right)\right](t,\cdot) \right\|_{L^2} \lesssim t^\sigma \|\chi_1(t^{-\sigma}D_x)\left[x_n(Z_m \Gamma v)_\pm\right](t,\cdot)\|_{L^\infty} \|u_\pm(t,\cdot)\|_{L^2}\\
&+ t^{-N(s)}\|x_n(Z_m\Gamma v)_\pm(t,\cdot)\|_{L^2}\|u_\pm(t,\cdot)\|_{H^s} + t^\sigma \|\chi_1(t^{-\sigma}D_x)\left[x_n(\Gamma v)_\pm\right](t,\cdot)\|_{L^\infty} \|(Z_mu)_\pm(t,\cdot)\|_{L^2}\\
& +t^{-N(s)} \left(\sum_{\mu=0}^1\|x_m^\mu x_n (\Gamma v)_\pm(t,\cdot)\|_{L^2} + t\|x_n(\Gamma v)_\pm(t,\cdot)\|_{L^2}\right) \left(\|u_\pm(t,\cdot)\|_{H^s}+ \|D_tu_\pm(t,\cdot)\|_{H^s}\right)\\
& + t^\sigma \|\chi_1(t^{-\sigma}D_x)\left[x_n(\Gamma v)_\pm\right](t,\cdot)\|_{L^\infty}\left(\|u_\pm(t,\cdot)\|_{L^2} + \|D_t u_\pm(t,\cdot)\|_{L^2}\right)\\
&+ t^{-N(s)}\|x_n(\Gamma v)_\pm(t,\cdot)\|_{L^2}\left(\|u_\pm(t,\cdot)\|_{H^s}+ \|D_t u_\pm(t,\cdot)\|_{H^s}\right),\\
\end{split}
\end{equation*} \normalsize
so from \eqref{est: bootstrap argument a-priori est}, \eqref{Hs_norm_DtU}, \eqref{norm_L2_xj-GammaIv-}, \eqref{est:xixjGamma v-}, \eqref{norm_Linfty_xjGammav-} and \eqref{est:Linfty_xjGammaI_v}, we derive that
\begin{equation*}
\left\|\chi(t^{-\sigma}D_x)\left[x_n Z_m \left(v^{I_1,\textit{NF}}  - (\Gamma v)_{-}\right)\right](t,\cdot) \right\|_{L^2} \le C(A+B)B\varepsilon^2 t^{\beta+\frac{\delta+\delta_1+\delta_2}{2}}.
\end{equation*}
Inequalities \eqref{est:L2_Zn(vINF-Gammav)_xmZN(vINF-Gammav)} follows then just by the observation that, after the hypothesis on multi-indices $I_1,I_2$ and the comparison between \eqref{explicit vNfGamma-vJ-} with $I=I_2$ and \eqref{Zm(vgamma-vJ)} with $I=I_1$, $2i\chi(t^{-\sigma}D_x)\left(v^{I_2,\textit{NF}}  - (Z_m\Gamma v)_{-}\right)$ corresponds to the first line in the right hand side of \eqref{Zm(vgamma-vJ)}.

In order to derive estimate \eqref{est:L2_Zm(VNF-Gammav)_enhanced} we apply $Z_m$ to both sides of equality \eqref{explicit VNF - VJ}, use \eqref{est:L2_Zn(vINF-Gammav)_xmZN(vINF-Gammav)}, formulas \eqref{commutator_Z_DtD1}, and successively proceed as follows: products in which $Z_m$ acts on $v$ and $\Gamma$ on $u$, that arise from the action of $Z_m$ on the second line of \eqref{explicit VNF - VJ}, are estimated using inequalities \eqref{cor_app_est_2} with $L=L^2$ and $w=u$; products in which $Z_m$ is acting on $v$ and there are no Klainerman vector fields acting on $u$ are estimated applying lemma \ref{Lem_appendix:L_estimate of products} with $L=L^2$; the $L^2$ norm of the remaining ones are controlled by the $L^\infty$ norm of the Klein-Gordon factor times the $L^2$ norm of the wave one. In this way we get that
\begin{equation*}
\begin{split}
& \left\|\chi(t^{-\sigma}D_x)\left[Z_m \left(\VNF  - (\Gamma v)_{-}\right)\right](t,\cdot) \right\|_{L^2} \lesssim \left\|\chi(t^{-\sigma}D_x)\left[Z_m \left(v^{I_1,\textit{NF}}  - (\Gamma v)_{-}\right)\right](t,\cdot) \right\|_{L^2} \\
& +t^\sigma \left\|\chi_1(t^{-\sigma}D_x)(Z_m v)_\pm(t,\cdot)\right\|_{L^\infty}\left( \|(\Gamma u)_\pm(t,\cdot)\|_{L^2} + \|u_\pm(t,\cdot)\|_{L^2}\right)\\
& + t^{-N(s)}\left(\sum_{|\mu|=0}^1 \|x^\mu (Z_m v)_\pm(t,\cdot)\|_{L^2}+ t\|(Z_mv)_\pm(t,\cdot)\|_{L^2}\right)\left(\|u_\pm(t,\cdot)\|_{H^s}+\|D_tu_\pm(t,\cdot)\|_{H^s}\right) \\
& + \|v_\pm(t,\cdot)\|_{H^{1,\infty}}\left(\|(Z_m\Gamma u)_\pm(t,\cdot)\|_{L^2}+ \|(\Gamma u)_\pm(t,\cdot)\|_{L^2} + \|D_t(\Gamma u)_\pm(t,\cdot)\|_{L^2}\right.\\
&\hspace{9cm}\left. + \|u_\pm(t,\cdot)\|_{L^2} +\|D_tu_\pm(t,\cdot)\|_{L^2} \right).
\end{split}
\end{equation*}
Choosing $s>0$ large so that $N(s)>2$ and using \eqref{est: bootstrap argument a-priori est} \eqref{Hs_norm_DtU}, \eqref{DtUI}, \eqref{norm_L2_xj-GammaIv-} with $k=1$, \eqref{Linfty_est_VJ}, \eqref{est:L2_Zn(vINF-Gammav)}, we hence recover \eqref{est:L2_Zm(VNF-Gammav)_enhanced}. An analogous procedure leads us to the following inequality
\begin{equation*}
\begin{split}
& \left\|\chi(t^{-\sigma}D_x)\left[x_n Z_m \left(\VNF  - (\Gamma v)_{-}\right)\right](t,\cdot) \right\|_{L^2} \lesssim \left\|\chi(t^{-\sigma}D_x)\left[x_n Z_m \left(v^{I_1,\textit{NF}}  - (\Gamma v)_{-}\right)\right](t,\cdot) \right\|_{L^2} \\
& +t^\sigma \left\|\chi_1(t^{-\sigma}D_x)\left[x_n(Z_m v)_\pm\right](t,\cdot)\right\|_{L^\infty}\left( \|(\Gamma u)_\pm(t,\cdot)\|_{L^2} + \|u_\pm(t,\cdot)\|_{L^2}\right)\\
& + t^{-N(s)}\left(\sum_{|\mu|=0}^1 \|x^\mu x_n(Z_m v)_\pm(t,\cdot)\|_{L^2}+ t\|x_n (Z_mv)_\pm(t,\cdot)\|_{L^2}\right)\left(\|u_\pm(t,\cdot)\|_{H^s}+\|D_tu_\pm(t,\cdot)\|_{H^s}\right) \\
& +\sum_{|\mu|=0}^1 \left\|x_n \Big(\frac{D_x}{\langle D_x\rangle}\Big)^\mu v_\pm(t,\cdot)\right\|_{L^\infty}\left(\|(Z_m\Gamma u)_\pm(t,\cdot)\|_{L^2}+ \|(\Gamma u)_\pm(t,\cdot)\|_{L^2} + \|D_t(\Gamma u)_\pm(t,\cdot)\|_{L^2}\right.\\
&\hspace{10cm}\left. + \|u_\pm(t,\cdot)\|_{L^2} +\|D_tu_\pm(t,\cdot)\|_{L^2} \right) ,
\end{split}
\end{equation*}
and estimate \eqref{est:L2_xnZm(VNF-Gammav)} is obtained by choosing $s>0$ large so that $N(s)>1$ and using \eqref{Hs_norm_DtU}, \eqref{DtUI}, \eqref{norm_Linfty_xv-}, \eqref{norm_L2_xj-GammaIv-} with $k=1$, \eqref{est:xixjGamma v-}, \eqref{norm_Linfty_xjGammav-}, \eqref{est:L2_xmZn(vINF-Gammav)} and a-priori estimates.
\endproof
\end{lem}

\begin{lem} \label{Lem_appendix:Lm_ZnVNF}
Let $\Gamma\in \{\Omega, Z_m, m=1,2\}$ be a Klainerman vector field and $\VNF$ be the function defined in \eqref{def_VNF}. There exists a constant $C>0$ such that, for any $\chi\in C^\infty_0(\mathbb{R}^2)$, $\sigma>0$ small, $m,n=1,2$, and every $t\in [1,T]$,
\begin{equation} \label{est_Lm_ZnVNF}
\left\|\oph(\chi(h^\sigma\xi))\mathcal{L}_m\left[tZ_n\VNF(t,tx)\right]\right\|_{L^2(dx)}\le CB\varepsilon t^\frac{\delta_0}{2}.
\end{equation}
\proof
We warn the reader that, throughout the proof, we denote by $C$ and $\beta$ two positive constants that may change line after line, with $\beta\rightarrow 0$ as $\sigma\rightarrow 0$.

Let $\vNFGamma$ be the function defined in \eqref{def_vNF-Gamma} for a generic multi-index $I$ of length 1 or 2, and $I_1,I_2$ two multi-indices such that $\Gamma^{I_1}=\Gamma$, $\Gamma^{I_2}=Z_n\Gamma$.
Using \eqref{commutator_Z_Dt-<D>} we rewrite $Z_n\VNF$ as follows:
\begin{equation*}
Z_n\VNF = Z_n\left(\VNF - (\Gamma v)_{-}\right) + \left[(Z_n\Gamma v)_{-} - v^{I_2,\textit{NF}}\right] + v^{I_2,\textit{NF}} + \frac{D_n}{\langle D_x\rangle}v^{I_1,\textit{NF}} + \frac{D_n}{\langle D_x\rangle}\left[(\Gamma v)_{-} - v^{I_1,\textit{NF}}\right]
\end{equation*}
so that
\begin{equation} \label{ineq:Lm_ZnVNF}
\begin{split}
&\left\|\oph(\chi(h^\sigma\xi))\mathcal{L}_m\left[tZ_n\VNF(t,tx)\right]\right\|_{L^2(dx)} \\
&\lesssim \left\|\oph(\chi(h^\sigma\xi))\mathcal{L}_m\left[tZ_n\left(\VNF - (\Gamma v)_{-}\right)(t,tx)\right]\right\|_{L^2(dx)} \\
& + \left\|\oph(\chi(h^\sigma\xi))\mathcal{L}_m\left[t\left[(Z_n\Gamma v)_{-} - v^{I_2,\textit{NF}}\right](t,tx)\right]\right\|_{L^2(dx)}\\
&  + \left\|\oph(\chi(h^\sigma\xi))\mathcal{L}_m\left[t v^{I_2,\textit{NF}} (t,tx)\right]\right\|_{L^2(dx)}+ \left\|\oph(\chi(h^\sigma\xi))\mathcal{L}_m\left[t\frac{D_n}{\langle D_x\rangle}v^{I_1,\textit{NF}}(t,tx)\right]\right\|_{L^2(dx)}\\
& + \left\|\oph(\chi(h^\sigma\xi))\mathcal{L}_m\left[t\frac{D_n}{\langle D_x\rangle}\left[(\Gamma v)_{-} - v^{I_1,\textit{NF}}\right](t,tx)\right]\right\|_{L^2(dx)}.
\end{split}
\end{equation}
Since $v^{I_2,\textit{NF}}$ satisfies \eqref{def_NLNF_I=2} with $I=I_2$, we derive from relation \eqref{relation_Zjw_Ljwidetilde(w)} with $w=v^{I_2,\textit{NF}}$ that
\begin{equation*}
\begin{split}
 &\left\|\oph(\chi(h^\sigma\xi))\mathcal{L}_m\left[t v^{I_2,\textit{NF}} (t,tx)\right]\right\|_{L^2(dx)} \lesssim \left\|\chi(t^{-\sigma}D_x)Z_m(\Gamma^{I_2}v)_{-}(t,\cdot) \right\|_{L^2} \\
 &+\left\|\chi(t^{-\sigma}D_x)Z_m\left[v^{I_2, \textit{NF}}- (\Gamma^{I_2}v)_{-}\right](t,\cdot) \right\|_{L^2} + \left\|\chi(t^{-\sigma}D_x)v^{I_2,\textit{NF}}(t,\cdot) \right\|_{L^2} \\
 &+ \left\| \chi(t^{-\sigma}D_x)\left[x_m \textit{NL}^{I_2,\textit{NF}}_{kg}\right](t,\cdot)\right\|_{L^2}.
\end{split}
\end{equation*}
A-priori estimate \eqref{est: bootstrap E02} with $k=0$, \eqref{ineq:L2 Zm vINF-vI} with $I=I_2$, \eqref{xj_NLNF_I=2}, \eqref{est_Linfty_vINF_2}, the fact that $\delta\ll \delta_2\ll \delta_1\ll \delta_0$ and that $\beta$ is small as long as $\sigma$ is small, imply
\begin{equation*}
\left\|\oph(\chi(h^\sigma\xi))\mathcal{L}_m\left[t v^{I_2,\textit{NF}} (t,tx)\right]\right\|_{L^2(dx)} \le CB\varepsilon t^\frac{\delta_0}{2}.
\end{equation*}
Analogously, commutating $\mathcal{L}_m$ with $\oph(\xi_n \langle \xi\rangle^{-1})$, using \eqref{relation_Zjw_Ljwidetilde(w)} with $w=v^{I_1,\textit{NF}}$ and the fact that $v^{I_1,\textit{NF}}$ is solution to \eqref{KG_vNF-Gamma} with non-linear term given by \eqref{def_NLNF}, together with inequalities \eqref{est:vNFGamma}, \eqref{est_L2_ZmvNFGamma}, \eqref{est_xj_NL-kg-NF-Gamma}, we derive that
\begin{equation*}
\left\|\oph(\chi(h^\sigma\xi))\mathcal{L}_m\left[t\frac{D_n}{\langle D_x\rangle}v^{I_1,\textit{NF}}(t,tx)\right]\right\|_{L^2(dx)} \le CB\varepsilon t^\frac{\delta_1}{2}.
\end{equation*}
Finally, the remaining norms in the right hand side of \eqref{ineq:Lm_ZnVNF} are estimated by the right hand side of \eqref{est_Lm_ZnVNF} after \eqref{dec_norm_Lm} and lemma \ref{Lem_appendix:est_Zm(vINF-Gammav)}. 
\endproof
\end{lem}

Lemmas \ref{Lem_appendix: preliminary est VJ}, \ref{Lem_appendix: LvtildeGamma} and \ref{Lem_appendix:Lm_ZnVNF} allow us to prove an analogous result to that of lemma \ref{Lem_appendix:product_Vtilde_Utilde}, where $\widetilde{v}$ is replaced with $\widetilde{V}^\Gamma$ introduced in \eqref{def_VtildeGamma}.

\begin{lem} \label{Lem_appendix: product_VtildeGamma_utilde}
Let $h=t^{-1}$, $\ut, \widetilde{V}^\Gamma$ be respectively defined in \eqref{def utilde vtilde} and\eqref{def_VtildeGamma}, $a_0(\xi)\in S_{0,0}(1)$, and $b_1(\xi)=\xi_j$ or $b_1(\xi)=\xi_j\xi_k|\xi|^{-1}$, with $j,k\in \{1,2\}$. 
There exists a constant $C>0$ such that, for any $\chi, \chi_1\in C^\infty_0(\mathbb{R}^2)$, $\sigma>0$, and every $t\in [1,T]$, we have that
\begin{subequations} 
\begin{equation}\label{est_VtildeGamma_utilde}
\big\|[\oph(\chi(h^\sigma\xi)a_0(\xi))\widetilde{V}^\Gamma(t,\cdot)][\oph(\chi_1(h^\sigma\xi)b_1(\xi))\widetilde{u}(t,\cdot)]\big\|_{L^2} \le C(A+B)B\varepsilon^2 h^{\frac{1}{2}-\beta'},
\end{equation}
\begin{equation} \label{est_Linfty_VtildeGamma_utilde}
\big\|[\oph(\chi(h^\sigma\xi)a_0(\xi))\widetilde{V}^\Gamma(t,\cdot)][\oph(\chi_1(h^\sigma\xi)b_1(\xi))\widetilde{u}(t,\cdot)]\big\|_{L^\infty} \le C(A+B)B\varepsilon^2 h^{-\beta'},
\end{equation}
\end{subequations}
with $\beta'>0$ small, $\beta\rightarrow 0$ as $\sigma,\delta_0\rightarrow 0$.
\proof
The proof of this result has the same structure as that of lemma \ref{Lem_appendix:product_Vtilde_Utilde}. Only few differences occur due to to the fact that we are replacing $\widetilde{v}$ with $\widetilde{V}^\Gamma$. We limit here to indicate them.

Instead of referring to estimate \eqref{Hrho_infty_vtilde_appendix} we use the fact that, after \eqref{Linfty_est_VJ} in classical coordinates, there exists a constant $C>0$ such that for any $\rho\in\N$
\begin{equation} \label{est:Hrho_VtildeGamma}
\left\|\oph(\chi(h^\sigma\xi))\widetilde{V}^\Gamma(t,\cdot)\right\|_{H^{\rho,\infty}}\le CB\varepsilon h^{-\beta-\frac{\delta_1}{2}},
\end{equation}
with $\beta>0$ small such that $\beta\rightarrow 0$ as $\sigma\rightarrow 0$.
We successively decompose $\widetilde{V}^\Gamma$ into $\widetilde{V}^\Gamma_{\Lambda_{kg}}+\widetilde{V}^\Gamma_{\Lambda^c_{kg}}$, with
\begin{gather*}
\widetilde{V}^\Gamma_{\Lambda_{kg}}(t,x):= \oph\Big(\gamma\Big(\frac{x-p'(\xi)}{\sqrt{h}}\Big)\chi(h^\sigma\xi)a_0(\xi)\Big)\widetilde{V}^\Gamma(t,x), \\
 \widetilde{V}^\Gamma_{\Lambda^c_{kg}}(t,x):=\oph\Big((1-\gamma)\Big(\frac{x-p'(\xi)}{\sqrt{h}}\Big)\chi(h^\sigma\xi)a_0(\xi)\Big)\widetilde{V}^\Gamma(t,x).
\end{gather*}
On the one hand, from the fact that above operators are supported for frequencies $|\xi|\lesssim h^\sigma$,
together with proposition \ref{Prop:Continuity Lp-Lp} with $p=+\infty$ and \eqref{est:Hrho_VtildeGamma}, we have that
\begin{equation*}
\left\|\widetilde{V}^\Gamma_{\Lambda_{kg}}(t,\cdot)\right\|_{L^\infty}\le CB\varepsilon h^{-\beta-\frac{\delta_1}{2}}.
\end{equation*}
On the other hand, combining the analogous of \eqref{est:vLmabdakgc_appendix_preliminary} with lemma \ref{Lem_appendix: LvtildeGamma} (instead of \ref{Lem: from energy to norms in sc coordinates-KG}), estimates \eqref{est:L2_ZmVNF}, \eqref{est_Lm_ZnVNF} (instead of lemma \ref{Lem_appendix: ZnvNF_LmZnvNF}) and \eqref{first_estimate_xNLT} (instead of \eqref{chi_xm_rNFkg}),
\begin{equation*}
\left\|\widetilde{V}^\Gamma_{\Lambda^c_{kg}}(t,\cdot)\right\|_{L^\infty}\le CB\varepsilon h^{\frac{1}{2}-\beta-\frac{\delta_1}{2}}.
\end{equation*}
\end{lem}

\begin{lem} \label{Lem_appendix:xx(VNF-Gammav)}
Let $\Gamma\in \{\Omega, Z_m, m=1,2\}$ be a Klainerman vector field and $\VNF$ be the function defined in \eqref{def_VNF}. There exists a constant $C>0$ such that, for any $\chi\in C^\infty_0(\mathbb{R}^2)$, $\sigma>0$ small, $m,n=1,2$, and every $t\in [1,T]$,
\begin{subequations} \label{est:x_x_(VNF-Gammav)}
\begin{align}
\left\|\chi(t^{-\sigma}D_x)\left[x_m\left(\VNF - (\Gamma v)_{-}\right)\right](t,\cdot)\right\|_{L^\infty}& \le C(A+B)^2\varepsilon^2 t^{-\frac{1}{2}+\beta +\frac{\delta_1+\delta_2}{2}},\label{est:xn(VNF-Gammav)} \\
\left\|\chi(t^{-\sigma}D_x)\left[x_n x_m\left(\VNF - (\Gamma v)_{-}\right)\right](t,\cdot)\right\|_{L^\infty} &\le C(A+B)^2\varepsilon^2 t^{\frac{1}{2}+\beta +\frac{\delta_1+\delta_2}{2}}, \label{est:xmxn(VNF-Gammav)}
\end{align}
\end{subequations}
with $\beta>0$ small such that $\beta\rightarrow 0$ as $\sigma\rightarrow 0$.
\proof
We remind the reader about explicit expression \eqref{explicit VNF - VJ} of the difference $\VNF - (\Gamma v)_{-}$, and \eqref{explicit vNfGamma-vJ-}, here considered with $|I|=1$ such that $\Gamma^I=\Gamma$.

We first use equalities \eqref{def u+- v+-}, \eqref{def uIpm vIpm}, and, after some commutations, multiply $x_m$ (together with $x_n$ when proving \eqref{est:xmxn(VNF-Gammav)}) against each Klein-Gordon factor. Successively, we estimate the contribution coming from $\vNFGamma - (\Gamma v)_{-}$ using lemma \ref{Lem_appendix:L_estimate of products} with $L=L^\infty$, and all products coming from the second line of \eqref{explicit VNF - VJ} by means of inequalities \eqref{cor_app_est_2} with $L=L^\infty$, $w=u$ and $w_{j_0}=(D_x\langle D_x\rangle^{-1})^\mu Z_m v$ for $|\mu|=0,1$. On the one hand, we obtain that
\begin{equation*}
\begin{split}
& \left\|\chi(t^{-\sigma}D_x)\left[x_m\left(\VNF - (\Gamma v)_{-}\right)\right](t,\cdot)\right\|_{L^\infty} \\
&\lesssim \sum_{\mu=0}^1 t^\sigma \left\|\chi_1(t^{-\sigma}D_x)\left[x_m (\Gamma v)_\pm(t,\cdot)\right]\right\|_{L^\infty}\|\mathrm{R}^\mu_1 u_\pm(t,\cdot)\|_{L^\infty} + t^{-N(s)}\|x_m (\Gamma v)_\pm(t,\cdot)\|_{L^2}\|u_\pm(t,\cdot)\|_{H^s}\\
& + \sum_{|\mu|=0}^1 t^\sigma\left\|x_m \Big(\frac{D_x}{\langle D_x\rangle}\Big)^\mu v_\pm\right\|_{L^\infty}\|\chi_1(t^{-\sigma}D_x)(\Gamma u)_{-}(t,\cdot)\|_{L^\infty} \\
&+ \sum_{|\mu|=0}^2 t^{-N(s)}\|x^\mu v_\pm(t,\cdot)\|_{L^2}\left(\|u_\pm(t,\cdot)\|_{H^s}+\|D_tu_\pm(t,\cdot)\|_{H^s}\right) \\
& + \sum_{|\mu|, |\nu|=0}^1 t^\sigma\left\|x_m \Big(\frac{D_x}{\langle D_x\rangle}\Big)^\mu v_\pm\right\|_{L^\infty} \|\mathrm{R}^\mu u_\pm(t,\cdot)\|_{H^{2,\infty}}
\end{split}
\end{equation*}
and estimate \eqref{est:xn(VNF-Gammav)} follows choosing $s>0$ large enough to have $N(s)\ge 2$ and using \eqref{est: bootstrap argument a-priori est}, \eqref{Hs_norm_DtU}, \eqref{norms_H1_Linfty_xv-}, \eqref{norm_L2_xj-GammaIv-} with $k=1$, \eqref{norm_Linfty_xjGammav-}, \eqref{Linfty_est_UJ}.
On the other hand,
\begin{equation*}
\begin{split}
& \left\|\chi(t^{-\sigma}D_x)\left[x_n x_m\left(\VNF - (\Gamma v)_{-}\right)\right](t,\cdot)\right\|_{L^\infty} \\
&\lesssim \sum_{\mu=0}^1 t^\sigma \left\|\chi_1(t^{-\sigma}D_x)\left[x_n x_m (\Gamma v)_\pm(t,\cdot)\right]\right\|_{L^\infty}\|\mathrm{R}^\mu_1 u_\pm(t,\cdot)\|_{L^\infty} \\
&+ t^{-N(s)}\|x_nx_m (\Gamma v)_\pm(t,\cdot)\|_{L^2}\|u_\pm(t,\cdot)\|_{H^s}\\
& + \sum_{|\mu|=0}^1 t^\sigma\left\|x_nx_m \Big(\frac{D_x}{\langle D_x\rangle}\Big)^\mu v_\pm\right\|_{L^\infty}\|\chi_1(t^{-\sigma}D_x)(\Gamma u)_{-}(t,\cdot)\|_{L^\infty} \\
&+ \sum_{|\mu|=0}^3 t^{-N(s)}\|x^\mu v_\pm(t,\cdot)\|_{L^2}\left(\|u_\pm(t,\cdot)\|_{H^s}+\|D_tu_\pm(t,\cdot)\|_{H^s}\right) \\
& + \sum_{|\mu|, |\nu|=0}^1 t^\sigma\left\|x_n x_m \Big(\frac{D_x}{\langle D_x\rangle}\Big)^\mu v_\pm\right\|_{L^\infty} \|\mathrm{R}^\mu u_\pm(t,\cdot)\|_{H^{2,\infty}}
\end{split}
\end{equation*}
so picking the same $s$ as before and using \eqref{Hs_norm_DtU}, \eqref{norm_xv-}, \eqref{est_x2v-}, \eqref{est:xixjGamma v-}, \eqref{est:x3_v-}, \eqref{Linfty_est_UJ} and \eqref{est:Linfty_xixjGammav-new}, together with a-priori estimates, we derive \eqref{est:xmxn(VNF-Gammav)}. 
\endproof
\end{lem}

\begin{lem} \label{Lem_appendix:xxNLT}
Let $\Gamma\in \{\Omega, Z_m, m=1,2\}$ be a Klainerman vector field and $\NLT$ be given by \eqref{explicit NLT}. There exists a constant $C>0$ such that for any $\chi\in C^\infty_0(\mathbb{R}^2)$, $\sigma>0$ small, $m,n=1,2$, and every $t\in [1,T]$,
\begin{subequations}
\begin{align}
\left\|\chi(t^{-\sigma}D_x)\left[x_n \NLT\right](t,\cdot)\right\|_{L^2} &\le C(A+B)^2B\varepsilon^3 t^{-1+\beta'}, \\
\left\|\chi(t^{-\sigma}D_x)\left[x_m x_n \NLT\right](t,\cdot)\right\|_{L^2} &\le C(A+B)^2B\varepsilon^3 t^{\beta'},
\end{align}
\end{subequations}
with $\beta'>0$ small such that $\beta'\rightarrow 0$ as $\sigma, \delta_0\rightarrow 0$.
Moreover, in the same time interval
\begin{equation} \label{est_Linfty_NLT}
\left\|\chi(t^{-\sigma}D_x) \NLT(t,\cdot)\right\|_{L^\infty} \le C(A+B)^2B\varepsilon^3 t^{-\frac{5}{2}+\beta'}.
\end{equation}
\proof
We warn the reader that we will denote by $C,\beta, \beta'$ some positive constants that may change line after line, with $\beta\rightarrow 0$ (resp. $\beta'\rightarrow 0$) as $\sigma\rightarrow 0$ (resp. as $\sigma,\delta_0\rightarrow 0$).
For a seek of compactness we also denote by $R(t,x)$ any contribution verifying 
\begin{equation}\label{est_L2_R}
\begin{split}
\left\|\chi(t^{-\sigma}D_x)\left[x_nR(t,\cdot)\right]\right\|_{L^2}& \le C(A+B)^2B\varepsilon^3 t^{-1+\beta'}, \\
\left\|\chi(t^{-\sigma}D_x)\left[x_m x_n R(t,\cdot)\right]\right\|_{L^2}& \le C(A+B)^2B\varepsilon^3 t^{\beta'},
\end{split}
\end{equation}
together with 
\begin{equation}\label{est_Linfty_R}
\left\|\chi(t^{-\sigma}D_x)R(t,\cdot)\right\|_{L^\infty}\le C(A+B)^2B\varepsilon^3 t^{-\frac{5}{2}+\beta'}.
\end{equation}
Let us introduce $\textit{NL}^{\emph{cub}}_v$ as follows
\begin{equation} \label{NL_cub_v}
\begin{split}
\textit{NL}^{\emph{cub}}_v :=& -\frac{i}{2}\left[-(D_1\Gamma v)\Nlw +D_1\left[(\Gamma v)\Nlw\right]\right]
\\
& -\frac{i}{2}\left[-(D_1v)\Gamma\textit{NL}_w + D_1[v \Gamma\textit{NL}_w]\right]  + \frac{i}{2}\delta_\Omega \left[-(D_2v)\textit{NL}_w + D_2[v \textit{NL}_w]\right] \\
& + \delta_{Z_1}\left[(D_tv)\textit{NL}^I_w - \langle D_x\rangle[v \textit{NL}_w]\right],
\end{split}
\end{equation}
so that from \eqref{def_NLT}
\begin{equation} \label{dec_NLT}
\begin{split}
\NLT = \frac{i}{2}\left[(\Gamma\Nlkg) (D_1u) + \Nlkg(D_1\Gamma u)\right] + \delta_\Omega \frac{i}{2} \Nlkg (D_2u) + \delta_{Z_1}\Nlkg (D_tu) + \textit{NL}^\emph{cub}_v,
\end{split}
\end{equation}
with $\delta_\Omega$ (resp. $\delta_{Z_1}$) equal to 1 when $\Gamma = \Omega$ (resp. $\Gamma =Z_1$), 0 otherwise.
After \eqref{def u+- v+-}, \eqref{def uIpm vIpm}, and estimates \eqref{est: bootstrap argument a-priori est}, \eqref{est L2 NLw}, \eqref{Hs norm DtV} with $s=0$, \eqref{norm_Linfty_xv-}, \eqref{L2_norm_NLIw}, $\textit{NL}^\emph{cub}_v$ verifies the following:
 \small
\begin{equation*} 
\begin{split}
& \left\| \chi(t^{-\sigma}D_x)\left[x_n \textit{NL}^\emph{cub}_v\right](t,\cdot)\right\|_{L^2}\\
&\lesssim \sum_{\mu, |\nu|=0}^1 t^\sigma \left\|x^\mu_n\Big(\frac{D_x}{\langle D_x\rangle}\Big)^\nu v_\pm(t,\cdot)\right\|_{L^\infty} \left[\|\Gamma\textit{NL}_w(t,\cdot)\|_{L^2}+\|\textit{NL}_w(t,\cdot)\|_{L^2} + \|v_\pm(t,\cdot)\|_{H^{2,\infty}}\|(\Gamma^I v)_\pm(t,\cdot)\|_{L^2}\right] \\
& \le C(A+B)AB\varepsilon^3 t^{-1+\sigma+\delta_2}.
\end{split}
\end{equation*}\normalsize
From the mentioned inequalities and the additional \eqref{norm_Linfty_xixjv-}, it also satisfies
\begin{equation*} 
\begin{split}
& \left\| \chi(t^{-\sigma}D_x)\left[x_m x_n \textit{NL}^\emph{cub}_v\right](t,\cdot)\right\|_{L^2}\\
&\lesssim \sum_{\mu_1,\mu_2,|\nu|=0}^1 t^\sigma \left\|x^{\mu_1}_mx^{\mu_2}_n\Big(\frac{D_x}{\langle D_x\rangle}\Big)^\nu v_\pm(t,\cdot)\right\|_{L^\infty} \left[\|\Gamma\textit{NL}_w(t,\cdot)\|_{L^2}+\|\textit{NL}_w(t,\cdot)\|_{L^2}\right.\\
&\hspace{7.5cm}\left. + \|v_\pm(t,\cdot)\|_{H^{2,\infty}}\|(\Gamma^I v)_\pm(t,\cdot)\|_{L^2}\right]\\
&  \le C(A+B)AB\varepsilon^3 t^{\sigma+\delta_2}.
\end{split}
\end{equation*}
Moreover, applying twice lemma \ref{Lem_appendix:L_estimate of products} with $L=L^\infty$ and $s>0$ large enough to have $N(s)\ge 2$, the first time to estimate products involving $\Gamma v$ and $\textit{NL}_w$ in \eqref{NL_cub_v}, the second one to estimate the first two quadratic contributions to $\Gamma\textit{NL}_w$ (see \eqref{Gamma_Nlw}), we derive that there are two smooth cut-off functions $\chi_1, \chi_2$ such that
\begin{equation*}
\begin{split}
&\left\| \chi(t^{-\sigma}D_x)\textit{NL}^\emph{cub}_v(t,\cdot)\right\|_{L^\infty} \lesssim t^\sigma\left\|\chi_1(t^{-\sigma}D_x)(\Gamma v)_\pm(t,\cdot)\right\|_{L^\infty}\|\Nlw(t,\cdot)\|_{L^\infty} \\
&+ t^{-2}\|(\Gamma v)_\pm(t,\cdot)\|_{L^2}\|\Nlw(t,\cdot)\|_{H^s}+ t^\sigma \|\chi_1(t^{-\sigma}D_x)\Gamma\textit{NL}_w(t,\cdot)\|_{L^\infty}\|v_\pm(t,\cdot)\|_{H^{1,\infty}} \\
& + t^{-2}\|\textit{NL}^I_w(t,\cdot)\|_{L^2}\|v_\pm(t,\cdot)\|_{H^s}   + t^\sigma\|v_\pm(t,\cdot)\|_{H^{1,\infty}}\|\textit{NL}_w(t,\cdot)\|_{L^\infty}
\end{split}
\end{equation*} 
and
\begin{equation*}
\begin{split}
& \|\chi_1(t^{-\sigma}D_x)\Gamma\textit{NL}_w(t,\cdot)\|_{L^\infty} \lesssim \|\chi_2(t^{-\sigma}D_x)(\Gamma v)_\pm(t,\cdot)\|_{H^{2,\infty}}\|v_\pm(t,\cdot)\|_{H^{2,\infty}} \\
&+ t^{-2}\|(\Gamma v)_\pm(t,\cdot)\|_{H^1}\|v_\pm(t,\cdot)\|_{H^s} + \|v_\pm(t,\cdot)\|_{H^{1,\infty}}\left(\|v_\pm(t,\cdot)\|_{H^{2,\infty}}+ \|D_tv_\pm(t,\cdot)\|_{H^{1,\infty}}\right).
\end{split}
\end{equation*}
From a-priori estimates, \eqref{est Linfty NLw}, \eqref{est Hs NLw-New}, \eqref{Hs_norm_DtU} with $s=0$, \eqref{est: Hsinfty Dt V} with $s=1$ and $\theta\ll 1$ small, \eqref{Linfty_est_VJ}, \eqref{L2_norm_NLIw}, we then recover
\begin{equation*}
\left\| \chi(t^{-\sigma}D_x)\textit{NL}^\emph{cub}_v(t,\cdot)\right\|_{L^\infty} \le CA^2B\varepsilon^3 t^{-3+\beta'}.
\end{equation*}
Those inequalities make $\textit{NL}^\emph{cub}_v$ a contribution of the form $R(t,x)$, so from \eqref{dec_NLT} we are left to prove that the same is true for $\Gamma\textit{NL}_{kg}(D_1u)$, $\textit{NL}_{kg}(D_1\Gamma u)$, $\textit{NL}_{kg}(D_2 u)$ and $\textit{NL}_{kg}(D_tu)$.

We immediately observe, from \eqref{def_app_Nlkg} and \eqref{def u+- v+-}, that the cubic contributions to $\textit{NL}_{kg}(D_2u)$ and $\textit{NL}_{kg} (D_t u)$ are of the form
\begin{equation} \label{prod_a0v_b1u_b0u}
[a_0(D_x)v_{-}] [b_1(D_x)u_{-}] b_0(D_x)u_{-},
\end{equation}
with $a_0(\xi)\in \{1,\xi_j\langle \xi\rangle^{-1}, j=1,2\}$, $b_1(\xi)\in \{\xi_1, \xi_j\xi_1|\xi|^{-1}, j=1,2\}$, $b_0(\xi)\in \{1,\xi_2|\xi|^{-1}\}$. Therefore, lemmas \ref{Lem_appendix:Linfty_est_rNFkg}, \ref{Lem_appendix: L xnrNFkg} imply that $\textit{NL}_{kg}(D_2 u)$ and $\textit{NL}_{kg}(D_tu)$ are remainders $R(t,x)$.
Furthermore, from \eqref{notation: NLGamma}, \eqref{def_G1} and the equation satisfied by $u_\pm$ in \eqref{system for uI+-, vI+-} with $|I|=0$, 
\begin{multline*}
\Gamma\textit{NL}_{kg}= Q^\mathrm{kg}_0((\Gamma v)_\pm, D_1 u_\pm) + Q^\mathrm{kg}_0(v_\pm, D_1 (\Gamma u)_\pm) \\
-\delta_{\Omega}Q^\mathrm{kg}_0(v_\pm, D_2 u_\pm) - \delta_{Z_1} \Big[Q^\mathrm{kg}_0(v_\pm, |D_x| u_\pm) + Q^\mathrm{kg}_0\big(v_\pm, Q^\mathrm{w}_0(v_\pm, D_1 v_\pm)\big)\Big],
\end{multline*}
with $\delta_\Omega$ (resp. $\delta_{Z_1}$) equal to 1 if $\Gamma =\Omega$ (resp. $\Gamma =Z_1$), 0 otherwise. Estimates \eqref{est: bootstrap upm} and \eqref{est:L2_xj_cubic} imply that
\begin{equation*}
\left\|\chi(t^{-\sigma}D_x)\left[x_n Q^\mathrm{kg}_0\big(v_\pm, Q^\mathrm{w}_0(v_\pm, D_1 v_\pm)\big) (D_1u) \right](t,\cdot)\right\|_{L^2(dx)} \le C(A+B)A^2B\varepsilon^4 t^{-\frac{3}{2}+\frac{\delta+\delta_2}{2}},
\end{equation*}
while after \eqref{est: bootstrap upm}, \eqref{est L2 NLw}, \eqref{norm_Linfty_xixjv-},
\begin{multline*}
\left\|\chi(t^{-\sigma}D_x)\left[x_mx_n Q^\mathrm{kg}_0\big(v_\pm, Q^\mathrm{w}_0(v_\pm, D_1 v_\pm)\big) (D_1u) \right](t,\cdot)\right\|_{L^2(dx)}\\
\lesssim \sum_{|\mu|=0}^1\left\| x_mx_n\Big(\frac{D_x}{\langle D_x\rangle}\Big)^\mu v_\pm(t,\cdot)\right\|_{L^\infty(dx)}\|\textit{NL}_w(t,\cdot)\|_{L^2(dx)}\|\mathrm{R}_1 u_\pm(t,\cdot)\|_{L^\infty}\\
\le C(A+B)A^2B\varepsilon t^{-\frac{1}{2}+\frac{\delta+\delta_2}{2}}
\end{multline*}
Also, for any $\theta\in ]0,1[$,
\begin{multline*}
\left\|\chi(t^{-\sigma}D_x)\left[ Q^\mathrm{kg}_0\big(v_\pm, Q^\mathrm{w}_0(v_\pm, D_1 v_\pm)\big) (D_1u) \right](t,\cdot)\right\|_{L^\infty(dx)}\\
 \lesssim \|v_\pm(t,\cdot)\|_{H^{1,\infty}}\|\textit{NL}_w(t,\cdot)\|_{H^{1,\infty}}\|\mathrm{R}_1u_\pm(t,\cdot)\|_{L^\infty}
\le CA^{4-\theta}B^\theta \varepsilon^4 t^{-\frac{7}{2}+\theta(1+\frac{\delta}{2})},
\end{multline*}
as follows from \eqref{est Hsinfty for NLw-New} with $s=1$ and a-priori estimates.
Thus $Q^\mathrm{kg}_0\big(v_\pm, Q^\mathrm{w}_0(v_\pm, D_1 v_\pm)\big) (D_1u)$ is a remainder $R(t,x)$. The same holds true for
\[\left[-\delta_{\Omega}Q^\mathrm{kg}_0(v_\pm, D_2 u_\pm) -\delta_{Z_1}Q^\mathrm{kg}_0(v_\pm, |D_x| u_\pm)\right](D_1u)\]
thanks to lemmas \ref{Lem_appendix:Linfty_est_rNFkg} and \ref{Lem_appendix: L xnrNFkg}, since the above term is linear combination of products of the form
\[[a_0(D_x)v_{-}]\, [b_1(D_x)u_{-}] \, \mathrm{R}_1 u_{-},\]
with the same $a_0(\xi)$ as before and $b_1(\xi)\in \{\xi_2, \xi_2\xi_j|\xi|^{-1}, |\xi|, j=1,2\}$, as one can check using \eqref{Q0_pm} and \eqref{def u+- v+-}.

Summing up, the very contributions for which we have to prove estimates \eqref{est_L2_R} and \eqref{est_Linfty_R} are the following:
\begin{subequations}\label{fst_products}
\begin{gather}
[a_0(D_x)(\Gamma v)_{-}]\, [b_1(D_x)u_{-}] \, \mathrm{R}_1 u_{-} \label{first_product}\\
[a_0(D_x)v_{-}]\, [b_1(D_x)(\Gamma u)_{-}] \, \mathrm{R}_1 u_{-},\label{second_product} 
\end{gather}
which are the remaining types of products in $(\Gamma\textit{NL}_{kg})(D_1u)$, and
\begin{equation}
[a_0(D_x)v_{-}]\, [b_1(D_x)u_{-}] \, \mathrm{R}_1 (\Gamma u)_{-}, \label{third_product}
\end{equation}
\end{subequations}
which are the products appearing in $\Nlkg (D_1\Gamma u$), with $a_0$ being the same as above and $b_1(\xi)$ equal to $\xi_1$ or to $\xi_j\xi_1|\xi|^{-1}$, with $j=1,2$. 
All the manipulations we are going to make in what follows are aimed at showing that these estimates follow from lemmas \ref{Lem_appendix:product_Vtilde_Utilde}, \ref{Lem_appendix: product vtilde_utildeJ} and \ref{Lem_appendix: product_VtildeGamma_utilde}.

Firstly, we can assume that all factors in \eqref{fst_products} are truncated for moderate frequencies less or equal than $t^\sigma$, up to $R(t,x)$ contributions.
As regards \eqref{first_product}, this comes out from the application of lemma \ref{Lem_appendix:L_estimate of products}.
In fact, taking $L=L^2$, $w_1=x_m^k x_n a_0(D_x) (\Gamma v)_{-}$for $k\in \{0,1\}$, $s>0$ large enough to have $N(s)>2$, and using a-priori estimates and \eqref{norm_L2_xj-GammaIv-}, \eqref{est:xixjGamma v-}, we find that there is some $\chi_1\in C^\infty_0(\mathbb{R}^2)$ such that, for $k=0,1$,
\begin{equation*}
\begin{split}
&\left\|\chi(t^{-\sigma}D_x)\big[x_m^k x_n [a_0(D_x)(\Gamma v)_{-}]\, [b_1(D_x)u_{-}] \, \mathrm{R}_1 u_{-} \big]\right\|_{L^2(dx)} \\
&\lesssim \left\|\left[\chi_1(t^{-\sigma}D_x)[x_m^k x_na_0(D_x)(\Gamma v)_{-}]\right] \left[\chi(t^{-\sigma}D_x)b_1(D_x)u_{-}\right] \left[\chi(t^{-\sigma}D_x)\mathrm{R}_1 u_{-}\right]\right\|_{L^2(dx)} \\
&+ t^{-2}\sum_{\mu_1,\mu_2,|\nu|=0}^1 \|x_m^{\delta_{k1}\mu_1}x^{\mu_2}_n(\Gamma v)_{-}(t,\cdot)\|_{L^2(dx)} \|\mathrm{R}^\nu u_{-}(t,\cdot)\|_{H^{2,\infty}}\|u_{-}(t,\cdot)\|_{H^s}\\
&\lesssim \left\|\left[\chi_1(t^{-\sigma}D_x)[x_m^kx_na_0(D_x)(\Gamma v)_{-}]\right] \left[\chi(t^{-\sigma}D_x)b_1(D_x)u_{-}\right] \left[\chi(t^{-\sigma}D_x)\mathrm{R}_1 u_{-}\right]\right\|_{L^2(dx)} \\
&+ CAB^2\varepsilon^3 t^{-\frac{3}{2}(1-k)+\frac{\delta+\delta_2}{2}},
\end{split}
\end{equation*}
where $\delta_{k1}$ is the Kronecker delta. 
Taking instead $L=L^\infty$, from a-priori estimates we derive that
\begin{equation*}
\begin{split}
&\left\|\chi(t^{-\sigma}D_x)\big[ [a_0(D_x)(\Gamma v)_{-}]\, [b_1(D_x)u_{-}] \, \mathrm{R}_1 u_{-} \big]\right\|_{L^\infty(dx)} \\
& \lesssim \left\|\left[\chi_1(t^{-\sigma}D_x)[a_0(D_x)(\Gamma v)_{-}]\right] \left[\chi(t^{-\sigma}D_x)b_1(D_x)u_{-}\right] \left[\chi(t^{-\sigma}D_x)\mathrm{R}_1 u_{-}\right]\right\|_{L^\infty(dx)}\\
& + t^{-2}\sum_{|\mu|=0}^1 \|(\Gamma v)_{-}(t,\cdot)\|_{L^2} \|\mathrm{R}^\mu u_{-}(t,\cdot)\|_{H^{2,\infty}}\|u_{-}(t,\cdot)\|_{H^s}\\
& \lesssim \left\|\left[\chi_1(t^{-\sigma}D_x)[a_0(D_x)(\Gamma v)_{-}]\right] \left[\chi(t^{-\sigma}D_x)b_1(D_x)u_{-}\right] \left[\chi(t^{-\sigma}D_x)\mathrm{R}_1 u_{-}\right]\right\|_{L^\infty(dx)}\\
& + CAB^2\varepsilon^3 t^{-\frac{5}{2}+\frac{\delta+\delta_2}{2}}.
\end{split}
\end{equation*}
As concerns instead products \eqref{second_product} and \eqref{third_product}, this follows applying inequalities \eqref{cor_app_est_2} with $w=u$, $w_{j_0}=x_m^kx_na_0(D_x)v_{-}$ for $k=0,1$, and $s>0$ such that $N(s)\ge$.
In fact, for $L=L^2$ we use estimates \eqref{est: bootstrap argument a-priori est}, \eqref{norms_H1_Linfty_xv-}, \eqref{est_x2v-}, together with \eqref{est:x3_v-}, to derive that for $k\in \{0,1\}$
\begin{equation*}
\begin{split}
& \left\|\chi(t^{-\sigma}D_x)\left[x_m^kx_n [a_0(D_x)v_{-}]\, [b_1(D_x)(\Gamma u)_{-}] \, \mathrm{R}_1 u_{-} \right]\right\|_{L^2(dx)}\\
& \lesssim \left\|\left[\chi_1(t^{-\sigma}D_x)[x_m^kx_n a_0(D_x) v_{-}]\right] \left[\chi(t^{-\sigma}D_x)b_1(D_x)(\Gamma u)_{-}\right] \chi(t^{-\sigma}D_x)\mathrm{R}_1 u_{-} \right\|_{L^2(dx)} \\
& + t^{-3}\sum_{|\mu_1|, \mu_2, \mu_3=0}^1 \left(\left\|x^{\mu_1}x_m^{\delta_{k1}\mu_2}x^{\mu_3}_n v_{-}(t,\cdot) \right\|_{L^2(dx)} + t \|x^{\mu_2}_n v_{-}(t,\cdot)\|_{L^2}\right) \|u_\pm(t,\cdot)\|_{H^s}\|\mathrm{R}_1u_{-}(t,\cdot)\|_{L^\infty} \\
&+ t^{-3}\left\|x^k_mx_na_0(D_x)v_{-}(t,\cdot)\right\|_{L^\infty(dx)}\|(\Gamma u)_{-}(t,\cdot)\|_{L^2}\|u_{-}(t,\cdot)\|_{H^s}\\
& \lesssim \left\| \left[\chi_1(t^{-\sigma}D_x)[x_m^kx_n a_0(D_x) v_{-}]\right] \left[\chi(t^{-\sigma}D_x)b_1(D_x)(\Gamma u)_{-}\right] \chi(t^{-\sigma}D_x)\mathrm{R}_1 u_{-} \right\|_{L^2(dx)}\\& + CAB^2\varepsilon^2 t^{-(1-k) + \frac{\delta+\delta_2}{2}}.
\end{split}
\end{equation*}
Using instead \eqref{cor_app_est_2} with $L=L^\infty$ along with \eqref{est: bootstrap argument a-priori est} and \eqref{norm_Linfty_xv-},
\begin{equation*}
\begin{split}
& \left\|\chi(t^{-\sigma}D_x)\left[ [a_0(D_x)v_{-}]\, [b_1(D_x)(\Gamma u)_{-}] \, \mathrm{R}_1 u_{-} \right]\right\|_{L^\infty} \\
& \lesssim \left\| \left[\chi_1(t^{-\sigma}D_x)a_0(D_x) v_{-}\right] \left[\chi(t^{-\sigma}D_x)b_1(D_x)(\Gamma u)_{-}\right] \chi(t^{-\sigma}D_x)\mathrm{R}_1 u_{-} \right\|_{L^\infty} + CAB^2\varepsilon^2 t^{-\frac{5}{2}+ \frac{\delta+\delta_2}{2}}.
\end{split}
\end{equation*} 

Secondly, we can assume that in \eqref{first_product} (resp. in \eqref{third_product}) $b_1(D_x) u_{-}$ is replaced with $b_1(D_x)\unf$ (with $\unf$ introduced in \eqref{def uNF}). This is justified up to some $R(t,x)$ terms that satisfy \eqref{est_L2_R} as consequence of \eqref{est: bootstrap upm}, \eqref{norm_L2_xj-GammaIv-}, \eqref{est:xixjGamma v-} (resp. \eqref{norm_xv-}, \eqref{norm_L2_xixjv-}), \eqref{b1_uNF-u-}, and also \eqref{est_Linfty_R} because of \eqref{est: bootstrap upm}, \eqref{Linfty_est_VJ} (resp. \eqref{est: boostrap vpm}) and \eqref{b1_uNF-u-}.
Hence we are led to estimate the $L^2$ norm of
\begin{subequations}
\begin{align}
& \left[\chi_1(t^{-\sigma}D_x)[x_m^kx_n^l a_0(D_x) (\Gamma v)_{-}]\right] \left[\chi(t^{-\sigma}D_x)b_1(D_x) \unf\right] \chi(t^{-\sigma}D_x)\mathrm{R}_1 u_{-}  \label{first}\\
& \left[\chi_1(t^{-\sigma}D_x)[x_m^kx_n^l a_0(D_x) v_{-}]\right] \left[\chi(t^{-\sigma}D_x)b_1(D_x)(\Gamma u)_{-}\right] \chi(t^{-\sigma}D_x)\mathrm{R}_1 u_{-} \label{second}\\
& \left[\chi_1(t^{-\sigma}D_x)[x_m^kx_n^l a_0(D_x) v_{-}\right] \left[\chi(t^{-\sigma}D_x)b_1(D_x) \unf\right] \chi(t^{-\sigma}D_x)\mathrm{R}_1 (\Gamma u)_{-}\label{third}
\end{align}
\end{subequations}
for $k=0,1$, $l=1$, and the $L^\infty$ norm of above products when $k=l=0$.

Thirdly, we can think of $a_0(D_x)(\Gamma v)_{-}$ in \eqref{first} and of $a_0(D_x)v_{-}$ in \eqref{second}, \eqref{third} as replaced with $a_0(D_x)\VNF$ and $a_0(D_x)\vnf$ respectively, where $\VNF$ has been introduced in \eqref{def_VNF} and $\vnf$ in \eqref{def vNF}.
For \eqref{first} (resp. \eqref{third}) this substitution is justified up to some $R(t,x)$ terms that satisfy \eqref{est_L2_R} and \eqref{est_Linfty_R}, the former because of a-priori estimate \eqref{est: bootstrap upm}, \eqref{est_L2_uNF} and \eqref{est:x_x_(VNF-Gammav)} (resp. \eqref{Linfty_est_UJ}, \eqref{est_L2_uNF} and \eqref{xx_(vnf-v-)}), the latter after \eqref{est: bootstrap upm}, \eqref{est:VNF-(Gammav)} (resp. \eqref{a0_vNF-v-}, \eqref{Linfty_est_UJ}) and the classical translation of the semi-classical \eqref{Hrho_infty_utilde_appendix}
\begin{equation*} 
\|\unf(t,\cdot)\|_{H^{\rho,\infty}}+ \|\mathrm{R}\unf(t,\cdot)\|_{H^{\rho,\infty}}\le CB\varepsilon t^{-\frac{1}{2}}.
\end{equation*}
Therefore, in order to conclude the proof we must prove that, for some $\chi, \chi_1\in C^\infty_0(\R^2)$ and $k\in \{0,1\}$,
\begin{align*}
& \left\|\left[\chi_1(t^{-\sigma}D_x)[x_m^kx_n a_0(D_x) (\Gamma v)_{-}]\right] \left[\chi(t^{-\sigma}D_x)b_1(D_x) \unf\right] \chi(t^{-\sigma}D_x)\mathrm{R}_1 u_{-}  \right\|_{L^2(dx)} \\
& + \left\|\left[\chi_1(t^{-\sigma}D_x)[x_m^kx_n a_0(D_x) v_{-}]\right] \left[\chi(t^{-\sigma}D_x)b_1(D_x)(\Gamma u)_{-}\right] \chi(t^{-\sigma}D_x)\mathrm{R}_1 u_{-} \right\|_{L^2(dx)} \\
& +\left\| \left[\chi_1(t^{-\sigma}D_x) [x_m^k x_n^la_0(D_x) v_{-}]\right] \left[\chi(t^{-\sigma}D_x)b_1(D_x) \unf\right] \chi(t^{-\sigma}D_x)\mathrm{R}_1 (\Gamma u)_{-}\right\|_{L^2(dx)} \\
&\le C(A+B)^2B\varepsilon^3 t^{-1+k+\beta'}
\end{align*}
and
\begin{align*}
& \left\|\left[\chi_1(t^{-\sigma}D_x) a_0(D_x) (\Gamma v)_{-}\right] \left[\chi(t^{-\sigma}D_x)b_1(D_x) \unf\right] \chi(t^{-\sigma}D_x)\mathrm{R}_1 u_{-}  \right\|_{L^\infty(dx)} \\
& + \left\|\left[\chi_1(t^{-\sigma}D_x) a_0(D_x) v_{-}\right] \left[\chi(t^{-\sigma}D_x)b_1(D_x)(\Gamma u)_{-}\right] \chi(t^{-\sigma}D_x)\mathrm{R}_1 u_{-} \right\|_{L^\infty(dx)} \\
& + \left\| \left[\chi_1(t^{-\sigma}D_x) a_0(D_x) v_{-}\right] \left[\chi(t^{-\sigma}D_x)b_1(D_x) \unf\right] \chi(t^{-\sigma}D_x)\mathrm{R}_1 (\Gamma u)_{-}\right\|_{L^\infty(dx)}\\
& \le C(A+B)^2B\varepsilon^3 t^{-\frac{5}{2}+\beta'}.
\end{align*}
Actually, using \eqref{est: bootstrap upm}, \eqref{Linfty_est_UJ}, and passing to the semi-classical framework and unknowns with $\widetilde{V}^\Gamma$ introduced in \eqref{def_VtildeGamma}, $\ut, \vt$ in \eqref{def utilde vtilde}, and $\ut^I(t,x)=t^{-1}(\Gamma u)_{-}(t, t^{-1}x)$, above inequalities will follow respectively from
\begin{equation}\label{est_L2_prod_sem} 
\begin{split}
\sum_{k=0}^1&\Big[\left\| \left[\oph(\chi_1(h^\sigma\xi))[x^k_m x_n\oph(a_0(\xi))\widetilde{V}^\Gamma]\right]\big[ \oph(\chi(h^\sigma\xi)b_1(\xi))\widetilde{u}\big](t,\cdot) \right\|_{L^2(dx)}\\
 &+\left\| \left[\oph(\chi_1(h^\sigma\xi))[x^k_m x_n\oph(a_0(\xi))\vt]\right]\big[ \oph(\chi(h^\sigma\xi)b_1(\xi))\ut^I\big](t,\cdot) \right\|_{L^2(dx)}\\
&+\left\| \left[\oph(\chi_1(h^\sigma\xi))[x^k_m x_n\oph(a_0(\xi))\vt]\right]\big[ \oph(\chi(h^\sigma\xi)b_1(\xi))\ut\big](t,\cdot) \right\|_{L^2(dx)}\Big]\\
& \le C(A+B)B\varepsilon^3 h^{-\frac{1}{2}-\beta'}
\end{split}
\end{equation}
and
\begin{equation}\label{est_Linfty_prod_sem}
\begin{split}
&\left\| \big[\oph(\chi_1(h^\sigma\xi)a_0(\xi))\widetilde{V}^\Gamma\big]\big[ \oph(\chi(h^\sigma\xi)b_1(\xi))\widetilde{u}\big](t,\cdot) \right\|_{L^\infty(dx)}\\
&+\left\| \big[\oph(\chi_1(h^\sigma\xi)a_0(\xi))\vt\big]\big[ \oph(\chi(h^\sigma\xi)b_1(\xi))\ut^I\big](t,\cdot) \right\|_{L^\infty(dx)} \\
&+\left\| \big[\oph(\chi_1(h^\sigma\xi)a_0(\xi))\vt\big]\big[ \oph(\chi(h^\sigma\xi)b_1(\xi))\widetilde{u}\big](t,\cdot) \right\|_{L^\infty(dx)}\le C(A+B)B\varepsilon^3 h^{-\beta'}.
\end{split}
\end{equation}
We immediately obtain from inequalities \eqref{ineq:xn_utilde_vtilde} and \eqref{ineq:xmxn_vtilde_utilde} that\small
\[\sum_{k=0}^1\left\| \left[\oph(\chi_1(h^\sigma\xi))[x^k_m x_n\oph(a_0(\xi))\vt]\right]\big[ \oph(\chi(h^\sigma\xi)b_1(\xi))\ut\big](t,\cdot) \right\|_{L^2(dx)}\le C(A+B)B\varepsilon^3 h^{-\frac{1}{2}-\beta'}.\]\normalsize
Moreover, one can check that
\begin{multline*} 
\left\| \left[\oph(\chi_1(h^\sigma\xi))[ x_n\oph(a_0(\xi))\widetilde{V}^\Gamma]\right]\big[ \oph(\chi(h^\sigma\xi)b_1(\xi))\widetilde{u}\big](t,\cdot) \right\|_{L^2(dx)} \\
\le  \left\| \Big[\oph\Big(\chi_1(h^\sigma\xi)a_0(\xi)\frac{\xi_n}{\langle\xi\rangle}\Big)\widetilde{V}^\Gamma\Big] [\oph(\chi(h^\sigma\xi)b_1(\xi))\widetilde{u}](t,\cdot) \right\|_{L^2} + CAB\varepsilon^2 h^{\frac{1}{2}-\beta'},
\end{multline*}
\begin{equation*}
\begin{split}
&\left\| \left[\oph(\chi_1(h^\sigma\xi))[ x_m x_n\oph(a_0(\xi))\widetilde{V}^\Gamma]\right]\big[ \oph(\chi(h^\sigma\xi)b_1(\xi))\widetilde{u}\big](t,\cdot) \right\|_{L^2(dx)} \\
&\lesssim  \left\| \Big[\oph\Big(\chi_1(h^\sigma\xi)a_0(\xi)\frac{\xi_n}{\langle\xi\rangle}\Big)\widetilde{V}^\Gamma\Big]\Big[\oph\Big(\chi(h^\sigma\xi)b_1(\xi)\frac{\xi_m}{|\xi|}\Big)\widetilde{u}\Big](t,\cdot) \right\|_{L^2(dx)}+ C(A+B)B\varepsilon^2 h^{\frac{1}{2}-\beta'},
\end{split}
\end{equation*}
and
\begin{multline*} 
\left\| \left[\oph(\chi_1(h^\sigma\xi))\big[ x_n \oph(a_0(\xi))\vt\big]\right]\big[ \oph(\chi(h^\sigma\xi)b_1(\xi))\widetilde{u}^I\big](t,\cdot) \right\|_{L^2(dx} \\
\le  \left\| \Big[\oph\Big(\chi_1(h^\sigma\xi)a_0(\xi)\frac{\xi_n}{\langle\xi\rangle}\Big)\vt\Big] [\oph(\chi(h^\sigma\xi)b_1(\xi))\widetilde{u}^I](t,\cdot) \right\|_{L^2} + CAB\varepsilon^2 h^{\frac{1}{2}-\beta'},
\end{multline*}
\begin{equation*}
\begin{split}
& \left\| \left[\oph(\chi_1(h^\sigma\xi))\big[x_m x_n \oph(a_0(\xi))\vt\big]\right]\big[ \oph(\chi(h^\sigma\xi)b_1(\xi))\widetilde{u}^J\big](t,\cdot) \right\|_{L^2} \\
& \lesssim \left\| \Big[\oph\Big(\chi_1(h^\sigma\xi)a_0(\xi)\frac{\xi_n}{\langle\xi\rangle}\Big)\vt\Big]\Big[\oph\Big(\chi(h^\sigma\xi)b_1(\xi)\frac{\xi_m}{|\xi|}\Big)\widetilde{u}^I\Big](t,\cdot) \right\|_{L^2(dx)} + C(A+B)B\varepsilon^2 h^{\frac{1}{2}-\beta'}.
\end{split}
\end{equation*}
This can be done using a similar argument to the one that led us to \eqref{ineq:xn_utilde_vtilde} and \eqref{ineq:xmxn_vtilde_utilde}, up to replacing $\widetilde{v}$ with $\widetilde{V}^\Gamma$ in \eqref{eq:xn_a0_vtilde}, referring to lemma \ref{Lem_appendix: LvtildeGamma} instead of \ref{Lem: from energy to norms in sc coordinates-KG}, and to estimate \eqref{est:Hrho_VtildeGamma} instead of \eqref{Hrho_infty_vtilde_appendix}, in order to derive the former two inequalities;
up to replacing $\widetilde{u}$ with $\widetilde{u}^I$ in \eqref{eq:xm_b1_utilde}, using lemma \ref{Lem_appendix: L^2 estimates uJ} instead of \eqref{est:utilde-Hs}, \eqref{est:Mutilde}, estimate \eqref{est:Hrho_utildeJ} instead of \eqref{Hrho_infty_utilde_appendix}, and the fact that for any $\theta\in]0,1[$
\begin{multline*}
\left\| \oph(\chi(h^\sigma\xi)b_1(\xi)\xi_m|\xi|^{-1})\widetilde{u}^I(t,\cdot)\right\|_{L^\infty}=t \left\|\chi(t^{-\sigma}D_x)b_1(D_x)D_m|D_x|^{-1}(\Gamma u)_{-}(t,\cdot)\right\|_{L^\infty}\\ \lesssim t \|\chi(t^{-\sigma}D_x) (\Gamma u)_{-}(t,\cdot)\|^{1-\theta}_{H^{3,\infty}}\|(\Gamma u)_{-}(t,\cdot)\|_{H^2}^\theta 
\le C(A+B)^{1-\theta}B^\theta\varepsilon t^{\frac{1}{2}+\beta+\frac{\delta_1}{2}+\frac{(1+\delta_1+\delta_2)}{2}\theta},
\end{multline*}
which is the analogous of \eqref{RRutillde} (last estimate deduced using \eqref{Linfty_est_UJ} and \eqref{est: bootstrap E02} with $k=1$), to demonstrate the latter two ones.
Therefore, above inequalities and \eqref{est_vtilde_utildeJ_L2}, \eqref{est_VtildeGamma_utilde} imply  \eqref{est_L2_prod_sem}.
Finally, \eqref{est_Linfty_prod_sem} is consequence of \eqref{est_Linfty_Vtilde Utilde}, \eqref{est_Linfty_vtilde_utildeJ} and \eqref{est_Linfty_VtildeGamma_utilde}.
That concludes the proof of the statement.
\endproof
\end{lem}

\begin{lem} \label{Lem_appendix:Lm xn NLT}
Let $\NLT$ be given by \eqref{explicit NLT}.
There exists a constant $C>0$ such that, for any $\chi\in C^\infty_0(\mathbb{R}^2)$, $\sigma>0$ small, $m,n=1,2$, and every $t\in [1,T]$,
\begin{gather*}
\left\| \oph(\chi(h^\sigma\xi))\mathcal{L}_m \left[t(tx_n)\NLT(t,tx)\right]\right\|_{L^2(dx)}\le C(A+B)^2B\varepsilon^3 t^{\beta'}, \\
\left\| \oph(\chi(h^\sigma\xi))\mathcal{L}_m \left[t(tx_n) Q^\mathrm{kg}_0\big(v_\pm, Q^\mathrm{w}_0(v_\pm, D_1 v_\pm)\big)(t,tx)\right]\right\|_{L^2(dx)}\le C(A+B)AB\varepsilon^3 t^{\frac{\delta+\delta_2}{2}},
\end{gather*}
with $\beta'>0$ such that $\beta'\rightarrow 0$ as $\sigma,\delta_0\rightarrow 0$.
\proof
Straightforward after \eqref{dec_norm_Lm}, lemma \ref{Lem_appendix:xxNLT}, estimate \eqref{est:L2_xj_cubic} and the following inequality
\begin{multline*}
\left\|\chi(t^{-\sigma}D_x)\left[x_mx_n Q^\mathrm{kg}_0\big(v_\pm, Q^\mathrm{w}_0(v_\pm, D_1 v_\pm)\big)\right](t,\cdot)\right\|_{L^2} \\\lesssim \sum_{|\mu|=0}^1 \left\|x_mx_n \Big(\frac{D_x}{\langle D_x\rangle}\Big)^\mu v_\pm(t,\cdot)\right\|_{L^\infty} \|\textit{NL}_w(t,\cdot)\|_{L^2}
 \le C(A+B)AB\varepsilon^3 t^{\frac{\delta+\delta_2}{2}},
\end{multline*}
deduced from \eqref{est: bootstrap argument a-priori est}, \eqref{est L2 NLw} and \eqref{norm_Linfty_xixjv-}.
\endproof
\end{lem}

\begin{lem}\label{Lem_appendix:Lcal2_VtildeGamma}
Let $\widetilde{V}^\Gamma$ be the function defined in \eqref{def_VtildeGamma}.
There exists some positive constant $C$ such that, for any $\chi\in C^\infty_0(\mathbb{R}^2)$, $\sigma>0$ small, and every $t\in [1,T]$,
\begin{equation}
\sum_{|\mu|=2}\left\|\oph(\chi(h^\sigma\xi)\mathcal{L}^\mu\widetilde{V}^\Gamma(t,\cdot)\right\|_{L^2}\le CB\varepsilon t^{\beta'},
\end{equation}
with $\beta'>0$ small, $\beta'\rightarrow 0$ as $\sigma,\delta_0\rightarrow 0$.
\proof
First of all we remind that $\VNF$ is solution to \eqref{half KG VNF}.
From relation \eqref{relation between Zjv and Lj vtilde} and the commutation between $\mathcal{L}_m$ and $\oph(\langle\xi\rangle)$ we deduce that, for any $m,n=1,2$,
\begin{equation*}
\begin{split}
&\left\|\oph(\chi(h^\sigma\xi))\mathcal{L}_m\mathcal{L}_n \widetilde{V}^\Gamma(t,\cdot)\right\|_{H^1_h}\lesssim \sum_{\mu=0}^1\Big[\left\|\oph(\chi(h^\sigma\xi))\mathcal{L}^\mu_m \big[tZ_n \VNF(t,tx)\big] \right\|_{L^2(dx)}\\
& + \left\|\oph(\chi(h^\sigma\xi))\mathcal{L}^\mu_m \oph\Big(\frac{\xi_n}{\langle\xi\rangle}\Big)\widetilde{V}^\Gamma(t,\cdot)\big]\right\|_{L^2}+ \left\|\oph(\chi(h^\sigma\xi))\mathcal{L}^\mu_m \big[t(tx_n)\NLT(t,tx)\big]\right\|_{L^2(dx)}\\
&+  \left\| \oph(\chi(h^\sigma\xi))\mathcal{L}^\mu_m \left[t(tx_n) Q^\mathrm{kg}_0\big(v_\pm, Q^\mathrm{w}_0(v_\pm, D_1 v_\pm)\big)(t,tx)\right]\right\|_{L^2(dx)}\Big].
\end{split}
\end{equation*}
The result of the statement follows then from \eqref{est:L2_ZmVNF}, \eqref{est:L2_xj_cubic}, \eqref{first_estimate_xNLT}, \eqref{est_Lm_ZnVNF}, and lemmas \ref{Lem_appendix: LvtildeGamma}, \ref{Lem_appendix:Lm xn NLT}.
\endproof
\end{lem}

\begin{lem}\label{Lem_appendix: sharp_est_VJ}
There exists a constant $C>0$ such that, for any $\chi\in C^\infty_0(\mathbb{R}^2)$ equal to 1 in a neighbourhood of the origin, $\sigma>0$ small, and every $t\in [1,T]$,
\begin{equation} \label{sharp_est_VI}
\sum_{|I|=1}\|\chi(t^{-\sigma}D_x)V^I(t,\cdot)\|_{L^\infty} \le CB\varepsilon t^{-1}.
\end{equation}
\proof
As this estimate is evidently satisfied when $I$ is such that $\Gamma^I$ is a spatial derivative after a-priori estimate \eqref{est: boostrap vpm}, we focus on proving the statement for $\Gamma^I\in \{\Omega, Z_m, m=1,2\}$ being a Klainerman vector field. For simplicity, we refer to $\Gamma^I$ simply by $\Gamma$.

Instead of proving the result of the statement directly on $(\Gamma v)_\pm$ we show that
\begin{equation} \label{sharp_est_VNF}
\left\|\VNF(t,\cdot)\right\|_{L^\infty}\le CB\varepsilon t^{-1},
\end{equation}
where $\VNF$ has been introduced in \eqref{def_VNF}. After \eqref{est:VNF-(Gammav)}, the above inequality evidently implies the statement.
The main idea to derive the sharp decay estimate in \eqref{sharp_est_VNF} is to use the same argument that, in subsection \ref{Subsection : The Derivation of the ODE Equation}, led us to the propagation of a-priori estimate \eqref{est: boostrap vpm}, i.e. to move to the semi-classical setting and deduce an ODE from equation \eqref{half KG VNF} satisfied by $\VNF$. The most important feature that will provide us with \eqref{sharp_est_VNF} is that the uniform norm of all involved non-linear terms is integrable in time.
Before going into the details, we also remind the reader our choice to denote by $C, \beta$ and $\beta'$ some positive constants that may change line after line, with $\beta\rightarrow 0$ (resp. $\beta'\rightarrow 0$) as $\sigma\rightarrow 0$ (resp. as $\sigma,\delta_0\rightarrow 0$).

Let us consider $\widetilde{V}^\Gamma(t,x):=t \VNF(t,tx)$, operator $\Gamma^{kg}$ as follows
\begin{equation*}
\Gamma^{kg}:=\oph\Big(\gamma\Big(\frac{x-p'(\xi)}{\sqrt{h}}\Big)\chi_1(h^\sigma\xi)\Big),
\end{equation*}
with $\gamma, \chi_1\in C^\infty_0(\mathbb{R}^2)$ such that $\gamma\equiv 1$ close to the origin, $\chi_1\equiv 1$ on the support of $\chi$, $p(\xi):=\langle \xi\rangle$, and
\begin{align*}
\widetilde{V}^\Gamma_{\Lambda_{kg}}(t,x)&:= \Gamma^{kg}\oph(\chi(h^\sigma\xi))\widetilde{V}^\Gamma(t,x), \\ \widetilde{V}^\Gamma_{\Lambda^c_{kg}}(t,x)&:=\oph\Big((1-\gamma)\Big(\frac{x-p'(\xi)}{\sqrt{h}}\Big)\chi_1(h^\sigma\xi)\Big)\oph(\chi(h^\sigma\xi))\widetilde{V}^\Gamma(t,x),
\end{align*}
so that
\begin{equation*}
\oph(\chi(h^\sigma\xi))\widetilde{V}^\Gamma(t,\cdot) =\widetilde{V}^\Gamma_{\Lambda_{kg}}+ \widetilde{V}^\Gamma_{\Lambda^c_{kg}}.
\end{equation*}
It immediately follows from inequality \eqref{est_2L-Linfty} and lemmas \ref{Lem_appendix: LvtildeGamma}, \ref{Lem_appendix:Lcal2_VtildeGamma}, that 
\begin{equation} \label{est_widetildeV_Lambda-kg-c}
\left\|\widetilde{V}^\Gamma_{\Lambda^c_{kg}}(t,\cdot)\right\|_{L^\infty}\lesssim \sum_{|\mu|=0}^2 h^{\frac{1}{2}-\beta}\left\|\oph(\chi(h^\sigma\xi))\mathcal{L}^\mu \widetilde{V}^\Gamma(t,\cdot) \right\|_{L^2}\le CB\varepsilon t^{-\frac{1}{2}+\beta'}.
\end{equation}
On the other hand, as $\VNF$ is solution to \eqref{half KG VNF} an explicit computation shows that $\widetilde{V}^\Gamma$ satisfies the following semi-classical pseudo-differential equation:
\begin{equation*}
\left[D_t - \oph(x\cdot\xi - \langle\xi\rangle)\right]\widetilde{V}^\Gamma(t,x)=h^{-1}\NLT(t,tx) - \delta_{Z_1}h^{-1} Q^\mathrm{kg}_0\big(v_\pm, Q^\mathrm{w}_0(v_\pm, D_1 v_\pm)\big)(t,tx),
\end{equation*}
with $\NLT$ given explicitly by \eqref{explicit NLT}.
Applying successively operators $\oph(\chi(h^\sigma\xi))$ and $\Gamma^{kg}$ to the above equation we find, from symbolic calculus and the first part of lemma \ref{Lem: Commutator Gamma-kg}, that $\widetilde{V}^\Gamma_{\Lambda_{kg}}$ satisfies
\begin{multline} \label{KG_VtildeGamma_Lambda}
\left[D_t - \oph(x\cdot\xi - \langle\xi\rangle)\right]\widetilde{V}^\Gamma_{\Lambda_{kg}}(t,x) = h^{-1}\Gamma^{kg}\oph(\chi(h^\sigma\xi))\left[\NLT(t,tx)\right]\\ - \delta_{Z_1} h^{-1}\Gamma^{kg}\oph(\chi(h^\sigma\xi))\big[  Q^\mathrm{kg}_0\big(v_\pm, Q^\mathrm{w}_0(v_\pm, D_1 v_\pm)\big)(t,tx)\big]
 - \oph(b(x,\xi))\oph(\chi(h^\sigma\xi))\widetilde{V}^\Gamma(t,x)\\
+ i\sigma h^{1+\sigma}\Gamma^{kg}\oph\big((\partial\chi)(h^\sigma\xi)\cdot(h^\sigma\xi)\big)\widetilde{V}^\Gamma,
\end{multline}
with symbol $b$ given by \eqref{symbol of commutator}.
Since $\gamma$'s derivatives vanish in a neighbourhood of the origin and $\partial\chi_1\equiv 0$ on the support of $\chi$, from  symbolic calculus of lemma \ref{Lem : a sharp b} and remark \ref{Remark:symbols_with_null_support_intersection}, \ref{Lem : composition gamma-1 and its argument}, together with inequalities \eqref{est_1L-Linfty}, \eqref{est_2L-Linfty}, that
\begin{multline*}
\left\| \oph(b(x,\xi))\oph(\chi(h^\sigma\xi))\widetilde{V}^\Gamma(t,\cdot)\right\|_{L^\infty}\\
\lesssim h^{\frac{3}{2}-\beta} \sum_{|\mu|=0}^2 \|\oph(\chi(h^\sigma\xi))\mathcal{L}^\mu \widetilde{V}^\Gamma(t,\cdot)\|_{L^2} + h^2\left\|\widetilde{V}^\Gamma(t,\cdot)\right\|_{L^2} \le CB\varepsilon t^{-\frac{3}{2}+\beta'},
\end{multline*}
where last estimate is obtained using lemmas \ref{Lem_appendix: LvtildeGamma}, \ref{Lem_appendix:Lcal2_VtildeGamma}.
Moreover, reminding lemma \ref{Lem:family_thetah} and using symbolic calculus we see that, for any $N\in \N$ as large as we want,
\begin{multline}\label{Gamma_Op(partialchi)_VtildeGamma}
h^{1+\sigma} \left\| \Gamma^{kg}\oph((\partial\chi(h^\sigma\xi)\cdot(h^\sigma\xi))\widetilde{V}^\Gamma(t,\cdot)\right\|_{L^\infty}\\
\le h^{1+\sigma} \left\|\Gamma^{kg}\theta_h(x)\oph((\partial\chi(h^\sigma\xi)\cdot(h^\sigma\xi))\widetilde{V}^\Gamma(t,\cdot) \right\|_{L^\infty} 
+ h^N\|\widetilde{V}^\Gamma(t,\cdot)\|_{L^2},
\end{multline}
where $\theta_h(x)$ is a smooth cut-off function supported in closed ball $\overline{B_{1-ch^{2\sigma}}(0)}$, with $c>0$ small.
Denoting $(\partial\chi)(\xi)\cdot \xi$ concisely by $\widetilde{\chi}(\xi)$, we observe from proposition \ref{Prop:Continuity Lp-Lp} with $p=+\infty$, together with the uniform continuity on $L^\infty$ of operator $\widetilde{\chi}(t^{-\sigma}D_x)$, the definition of $\widetilde{V}^\Gamma$ in terms of $\VNF$, and \eqref{est:VNF-(Gammav)}, that
\begin{multline*}
h^{1+\sigma}\left\| \Gamma^{kg} \theta_h(x)\oph(\widetilde{\chi}(h^\sigma\xi)) \widetilde{V}^\Gamma(t,\cdot)\right\|_{L^\infty} \lesssim h^{1-\beta}\left\|\theta_h(x)\oph(\widetilde{\chi}(h^\sigma\xi)) \widetilde{V}^\Gamma(t,\cdot)\right\|_{L^\infty}
\\
\le t^{\beta} \left\|\theta_h\Big(\frac{\cdot}{t}\Big)\widetilde{\chi}(t^{-\sigma}D_x) (\Gamma v)_{-}(t,\cdot) \right\|_{L^\infty}+C(A+B)B\varepsilon^2 t^{-\frac{5}{4}+\beta}.
\end{multline*}
Using the fact that, for $\theta_h^j(z):=\theta_h(z)z_j$, 
\begin{equation*}
\theta_h\Big(\frac{x}{t}\Big) (\Omega v)_{-} = t\Big[\theta_h^1\Big(\frac{x}{t}\Big) \partial_2v_{-} - \theta_h^2\Big(\frac{x}{t}\Big) \partial_1 v_{-}\Big]
\end{equation*}
and
\begin{equation*}
\theta_h\Big(\frac{x}{t}\Big) (Z_m v)_{-} = t\Big[\theta_h^m\Big(\frac{x}{t}\Big) \partial_t v_{-} + \theta_h\Big(\frac{x}{t}\Big)\partial_m v_{-}\Big] + \theta_h\Big(\frac{x}{t}\Big) \frac{D_m}{\langle D_x\rangle}v_{-}, \quad m=1,2,
\end{equation*}
and making some commutations, we can express $(\Gamma v)_{-}$ in terms of $v_{-}$ and its derivatives up to a loss in $t$. Thus, from the classical Sobolev injection combined inequality \eqref{ineq:1-chi}, we obtain that
\begin{equation*}
\begin{split}
t^{-\beta}\left\|\widetilde{\chi}(t^{-\sigma}D_x)\theta_h\Big(\frac{\cdot}{t}\Big)(\Gamma v)_{-}(t,\cdot) \right\|_{L^\infty}& \lesssim t^{-N(s)+1+\beta}\left(\|D_tv_\pm(t,\cdot)\|_{H^s}+\|v_\pm(t,\cdot)\|_{H^s}\right)\\
& \le CB\varepsilon t^{-\frac{3}{2}},
\end{split}
\end{equation*}
last estimate following by taking $s>0$ large enough to have $N(s)\ge 3$ and using a-priori estimates along with \eqref{Hs norm DtV} with $s=0$.
From \eqref{est:VtildeGamma} and \eqref{Gamma_Op(partialchi)_VtildeGamma} we hence derive that
\begin{equation*}
h^{1+\sigma} \left\| \Gamma^{kg}\oph((\partial\chi(h^\sigma\xi)\cdot(h^\sigma\xi))\widetilde{V}^\Gamma(t,\cdot)\right\|_{L^\infty} \le CB\varepsilon t^{-\frac{3}{2}},
\end{equation*}
so the last two terms in the right hand side of equation \eqref{KG_VtildeGamma_Lambda} are remainders $R(t,x)$ such that
\begin{equation}\label{remainder R(V)}
\|R(t,\cdot)\|_{L^\infty}\le CB\varepsilon t^{-\frac{5}{4}},
\end{equation}
for every $t\in [1,T]$.

After proposition \ref{Prop:Continuity Lp-Lp} with $p=+\infty$, estimate \eqref{est_Linfty_NLT}, and the fact that for any $\theta\in ]0,1[$,
\begin{equation*}
\left\|Q^{\mathrm{kg}}_0\left(v_\pm, Q^\mathrm{w}_0(v_\pm , D_1v_\pm)\right) (t,\cdot)\right\|_{L^\infty(dx)}\le CA^{3-\theta}B^\theta\varepsilon^3 t^{-3+\theta(1+\frac{\delta}{2})},
\end{equation*}
as follows by \eqref{est Hs NLw-New} with $s=1$ and a-priori estimates, we deduce (up to taking $\theta\ll 1$ small in the above inequality) that also the first two non-linear terms in the right hand side of \eqref{KG_VtildeGamma_Lambda} satisfy \eqref{remainder R(V)} and can be included into $R(t,x)$.
Therefore, $\widetilde{V}^\Gamma_{\Lambda_{kg}}$ satisfies
\begin{equation*}
\left[D_t - \oph(x\cdot\xi - \langle\xi\rangle)\right]\widetilde{V}^\Gamma_{\Lambda_{kg}}(t,x) = R(t,x),
\end{equation*}
and using \eqref{Op(a)v development} along with inequality \eqref{est Linfty of R2(widetildev)}, together with lemmas \ref{Lem_appendix: LvtildeGamma}, \ref{Lem_appendix:Lcal2_VtildeGamma}, we deduce that, for the same family of cut-off functions $\theta_h$ introduced above, $\widetilde{V}^\Gamma_{\Lambda_{kg}}$ is solution to the following ODE:
\begin{equation} \label{ODE_Vtilde}
D_t \widetilde{V}^\Gamma_{\Lambda_{kg}}(t,x) = -\theta_h(x)\phi(x)\widetilde{V}^\Gamma_{\Lambda_{kg}}(t,x)+ R(t,x),
\end{equation}
with $\phi(x)=\sqrt{1-|x|^2}$.
Since the inhomogeneous term $R(t,x)$ decays, in the uniform norm, at a rate which is integrable in time, we get that
\begin{equation*}
\|\widetilde{V}^\Gamma_{\Lambda_{kg}}(t,\cdot)\|_{L^\infty}\lesssim \|\widetilde{V}^\Gamma_{\Lambda_{kg}}(1,\cdot)\|_{L^\infty}+CB\varepsilon \le CB\varepsilon,
\end{equation*}
which summed up with \eqref{est_widetildeV_Lambda-kg-c} implies \eqref{sharp_est_VNF}, and hence the conclusion of the proof.
\endproof
\end{lem}

\printbibliography

\printindex
\end{document}